\documentclass[a4paper, BCOR=0.0mm, DIV=calc]{scrartcl}
\usepackage[T1]{fontenc}
\usepackage[utf8]{inputenc}
\usepackage[ngerman]{babel}
\usepackage[osf,sc]{mathpazo}
\usepackage{textcomp}
\usepackage{microtype}
\usepackage{graphicx, color}
\usepackage{amsmath, amssymb, amsfonts}
\usepackage{booktabs}
\usepackage{nicefrac}
\usepackage{wasysym,pifont}
\usepackage{tikz}
\usepackage{graphicx}
\usepackage{float}
\usepackage{wrapfig}
\usetikzlibrary{calc,intersections}
\KOMAoptions{twoside=false, twocolumn=false, headinclude=false, footinclude=false, mpinclude=false, pagesize=auto}
\recalctypearea
\setkomafont{section}{\mdseries\scshape\Large}
\addtokomafont{subsection}{\mdseries\scshape}
\setkomafont{title}{\mdseries\scshape}

\begin{document}
\title{On the Interpolation of series
\footnote{
Original title: "'De Interpolatione Serierum"', first published as part of the book \textit{Institutiones calculi differentialis cum eius usu in analysi finitorum ac doctrina serierum}, 1755, reprinted in in \textit{Opera Omnia}: Series 1, Volume 10, pp. 619 - 647, Eneström-Number E212, translated by: Alexander Aycock for the Euler-Kreis Mainz}}
\author{Leonhard Euler}
\date{ }
\maketitle

\paragraph*{§389}

A series is said to be interpolated, if its terms corresponding to fractional or even surdic indices are assigned. Therefore, if the general term of a series was known, the interpolation will be of no difficulty, since, whatever number is substituted for the index $x$, this expression yields the corresponding term. But if the series was of such a nature that its general term can not be exhibited in any way,  the interpolation of series of such a kind is highly difficult and in most cases the terms corresponding to non-integer indices can only be defined by means of infinite series. Since in the preceding chapter we determined the values corresponding to any indices of expressions of this kind which can not be expressed finitely by means of the customary techniques, this treatment will have a very great use for interpolations. Therefore, we will demonstrate the applicability following from the preceding chapter for this task more diligently.

\paragraph*{§390}

Therefore, let any series be propounded

\begin{alignat*}{9}
&1    && 2  && 3   && 4 \quad \quad \cdots  &&x \\
&A~+~ &&B~+~&& C~+~&& D~+~ \cdots ~+~ &&X,
\end{alignat*}
whose general term $X$ we assume to be known, but the summatory term $S$ is not. From this form another series whose general term is equal to the summatory term of that series, and this new series will be

\begin{equation*}
\renewcommand{\arraystretch}{1,2}
\setlength{\arraycolsep}{0.0mm}
\begin{array}{ccccccccccc}
1 & 2 & 3 & 4 & 5 \\
A, \quad & (A+B), \quad & (A+B+C), \quad & (A+B+C+D), \quad & (A+B+C+D+E) \quad & \text{etc.}
\end{array}
\end{equation*} 
and its general term or the one corresponding to the indefinite index $x$ will be

\begin{equation*}
=A+B+C+D+ \cdots +X=S;
\end{equation*}
since it is not explicitly known, the interpolation of this new series will lead to the same difficulties we mentioned before. Therefore, in order to interpolate this series, one needs the values of $S$  it takes on, if  any non-integer numbers are substituted for $x$. For, if $x$ would be an integer number,  the corresponding value of $S$ would be found without difficulty, of course by means of addition of as many terms of the series $A+B+C+D+\text{etc.}$ as $x$ contains units.

\paragraph*{§391}

Therefore, to use the results treated in the preceding chapter, let us put that $x$ is an integer number such that the value corresponding to it, $S=A+B+C+\cdots+X$, is known, and let us find the value $\Sigma$  $S$ is transformed into, if  one writes $x+\omega$ instead of $x$ while $\omega$ is any fraction; and $\Sigma$ will be the  term corresponding to the index $x+ \omega$ of the propounded series to be interpolated; therefore, having found the latter, the interpolation of this series will be done. Let $Z$ be the term corresponding to the index $x+ \omega$ of the series $A$, $B$, $C$, $D$, $E$ etc., and let $Z^{\prime}$, $Z^{\prime \prime}$, $Z^{\prime \prime \prime}$ be  consecutive terms corresponding to the indices $x+\omega +1$, $x+ \omega +2$, $x+\omega +3$ etc. And at first let us put that the infinitesimal terms of the series $A$, $B$, $C$, $D$ etc. vanish. Therefore, having constituted all this, the series

\begin{equation*}
\renewcommand{\arraystretch}{1,2}
\setlength{\arraycolsep}{0.0mm}
\begin{array}{ccccccccccc}
1 & 2 & 3 & 4  \\
A, \quad & (A+B), \quad & (A+B+C), \quad & (A+B+C+D) \quad & \text{etc.},
\end{array}
\end{equation*} 
whose term corresponding to the index $x$ is $S=A+B+C+\cdots +X$, will be interpolated by finding the term $\Sigma$  corresponding to the fractional index $x+\omega$. But, as we already found, it will be

\begin{alignat*}{9}
&\Sigma = S &&+X^{\prime} &&+X^{\prime \prime} &&+X^{\prime \prime \prime} &&+X^{\prime \prime \prime \prime} &&+ \text{etc.} \\
&           &&-Z^{\prime} &&-Z^{\prime \prime} &&-Z^{\prime \prime \prime} &&-Z^{\prime \prime \prime \prime}&&-\text{etc.}
\end{alignat*}
and hence one will have an infinite series equal to the term $\Sigma$ in question, which, because of

\begin{equation*}
Z=X+\frac{\omega dX}{dx}+\frac{\omega^2 ddX}{1 \cdot 2 dx^2}+\frac{\omega^3 d^3 X}{1 \cdot 2 \cdot 3 dx^3}+\text{etc.},
\end{equation*}
is transformed into 

\begin{alignat*}{9}
&\Sigma = S &&-~\frac{\omega}{dx} &&d. &&(X^{\prime}+X^{\prime \prime}+X^{\prime \prime \prime}+X^{\prime \prime \prime \prime}+\text{etc.}) \\[3mm]
&  &&-\frac{\omega^2}{2dx^2} &&dd. &&(X^{\prime}+X^{\prime \prime}+X^{\prime \prime \prime}+X^{\prime \prime \prime \prime}+\text{etc.}) \\[3mm]
&  &&-\frac{\omega^3}{6dx^3} &&d^3. &&(X^{\prime}+X^{\prime \prime}+X^{\prime \prime \prime}+X^{\prime \prime \prime \prime}+\text{etc.}) \\
& && &&\text{etc.}
\end{alignat*}
of which formulas the  more convenient more in each case can be applied.

\paragraph*{§392}

Let us  take the arbitrary harmonic series for $A$, $B$, $C$, $D$ etc. and consider

\begin{equation*}
\frac{1}{a}+\frac{1}{a+b}+\frac{1}{a+2b}+\frac{1}{a+3b}+\text{etc.},
\end{equation*}
whose general term or the one corresponding to the index $x$ is $=\frac{1}{a+(x-1)b}=X$. Therefore, let this series be formed

\begin{equation*}
\renewcommand{\arraystretch}{1,7}
\setlength{\arraycolsep}{0.0mm}
\begin{array}{ccccccccccc}
1 & 2 & 3 & 4  \\
\dfrac{1}{a}, \quad & \left(\dfrac{1}{a}+\dfrac{1}{a+b}\right), \quad & \left(\dfrac{1}{a}+\dfrac{1}{a+b}+\dfrac{1}{a+2b}\right), \quad & \left(\dfrac{1}{a}+\dfrac{1}{a+b}+\dfrac{1}{a+2b}+\dfrac{1}{a+3b}\right) \quad & \text{etc.},
\end{array}
\end{equation*} 
whose term corresponding to the index $x$ will therefore be

\begin{equation*}
S=\frac{1}{a}+\frac{1}{a+b}+\frac{1}{a+2b}+\frac{1}{a+3b}+\cdots + \frac{1}{a+(x-1)b}.
\end{equation*}
If now $\Sigma$ denotes the term corresponding to the index $x+\omega$ of this series, because of $Z=\frac{1}{a+(x+\omega -1)b}$, it will be

\begin{alignat*}{9}
&X^{\prime} &&= \frac{1}{a+bx}, \quad && Z^{\prime} &&=\frac{1}{a+bx+b \omega} \\[3mm]
&X^{\prime \prime} &&= \frac{ 1 }{a+b+bx}, \quad && Z^{\prime \prime} &&=\frac{1}{a+b+bx+b \omega} \\[3mm]
&X^{\prime \prime \prime} &&= \frac{1}{a+2b+bx}, \quad && Z^{\prime \prime \prime} &&=\frac{1}{a+2b+bx+b \omega} \\
&           &&\text{etc.}            &&   &&\text{etc.}
\end{alignat*}
and hence it will result

\begin{alignat*}{9}
&\Sigma = S &&+\frac{\quad \quad ~ 1 \quad \quad ~ }{a+bx}&&+\frac{\quad \quad \quad ~ 1 \quad \quad \quad  ~}{b+b+bx}&&+\frac{\quad \quad \quad  ~~ 1 \quad \quad \quad ~}{a+2b+bx}&&+\text{etc.} \\[3mm]
&           &&-\frac{1}{a+bx+b \omega}&&-\frac{1}{a+b+bx+b \omega}&&-\frac{1}{a+2b+bx +b \omega}&&-\text{etc.}
\end{alignat*}
but the other expression will be of this kind

\begin{alignat*}{9}
& \Sigma = S &&+ b\omega && \left(\frac{1}{(a+bx)^2}+\frac{1}{(a+b+bx)^2}+\frac{1}{(a+2b+bx)^2}+\text{etc.}\right) \\[3mm]
&   &&- b^2\omega^2 && \left(\frac{1}{(a+bx)^3}+\frac{1}{(a+b+bx)^3}+\frac{1}{(a+2b+bx)^3}+\text{etc.}\right) \\[3mm]
&  &&+ b^3\omega^3 && \left(\frac{1}{(a+bx)^4}+\frac{1}{(a+b+bx)^4}+\frac{1}{(a+2b+bx)^4}+\text{etc.}\right) \\
& && && \text{etc.}
\end{alignat*}

\subsection*{Example 1}

\textit{Let this series be propounded}

\begin{equation*}
\renewcommand{\arraystretch}{1,7}
\setlength{\arraycolsep}{0.0mm}
\begin{array}{ccccccccccc}
1 & 2 & 3 & 4  \\
1, \quad & \left(1+\dfrac{1}{2}\right), \quad & \left(1+\dfrac{1}{2}+\dfrac{1}{3}\right), \quad & \left(1+\dfrac{1}{2}+\dfrac{1}{3}+\dfrac{1}{4}\right) \quad & \text{etc.},
\end{array}
\end{equation*} 
\textit{whose terms corresponding to fractional indices are to be found.}\\[2mm]
Therefore, it will be $a=1$ and $b=1$; hence, if the term corresponding to the integer index $x$ is put

\begin{equation*}
S=1 +\frac{1}{2}+\frac{1}{3}+\cdots + \frac{1}{x}
\end{equation*}
and  the term corresponding to the fractional index $x+\omega$ is called $=\Sigma$, it will be

\begin{alignat*}{9}
& \Sigma &&=S &&+ \quad \frac{1}{1+x}&&+\quad \frac{1}{2+x}&&+\quad\frac{1}{3+x}&&+ \quad \frac{1}{4+x}&&+\quad \frac{1}{5+x}&&+\text{etc.} \\[3mm]
&   &&  &&- \frac{1}{1+x+\omega}&&+\frac{1}{2+x+\omega}&&+\frac{1}{3+x+\omega}&&+\frac{1}{4+x+\omega}&&+\frac{1}{5+x+\omega}&&+\text{etc.} \\
\end{alignat*}
But it is to be noted, if the term corresponding to the fractional index $\omega$ was found, which we want to put $=T$, that from it the term of the index $x+ \omega$ can easily be found; for, if $T^{\prime}$, $T^{\prime \prime}$, $T^{\prime \prime \prime}$ etc. denote the terms corresponding to the indices $1+\omega$, $2+\omega$, $3+\omega$ etc., it will be

\begin{alignat*}{9}
&T^{\prime} &&=T +\frac{1}{1+\omega} \\[3mm]
&T^{\prime \prime} &&=T +\frac{1}{1+\omega}+\frac{1}{2+\omega} \\[3mm]
&T^{\prime \prime \prime} &&=T +\frac{1}{1+\omega}+\frac{1}{2+\omega} +\frac{1}{3+ \omega}\\
& &&\text{etc.},
\end{alignat*}
whence it suffices to have investigated only the terms corresponding to the indices $\omega$ smaller than $1$. For this purpose, let us put $x=0$; it will also be $S=0$ and the term $T$ of the series corresponding to the fractional index $\omega$ will be expressed this way

\begin{alignat*}{9}
& T &&~=~ &&+ \quad \frac{1}{1}&&+ \quad \frac{1}{2}&&+\quad \frac{1}{3}&&+ \quad \frac{1}{4}&&+ \text{etc.} \\[3mm]
&   &&    &&-\frac{1}{1+\omega}&&-\frac{1}{2+\omega}&&-\frac{1}{3+\omega}&&-\frac{1}{4+\omega}&& -\text{etc.}
\end{alignat*}
or, after having converted these fractions into infinite series, another expression will result

\begin{alignat*}{9}
& T &&~=~ &&+\omega && \left(1+\frac{1}{2^2}+\frac{1}{3^2}+\frac{1}{4^2}+\frac{1}{5^2}+\text{etc.}\right) \\[3mm]
&  && &&-\omega^2 && \left(1+\frac{1}{2^3}+\frac{1}{3^3}+\frac{1}{4^3}+\frac{1}{5^3}+\text{etc.}\right) \\[3mm]
&  && &&+\omega^3 && \left(1+\frac{1}{2^4}+\frac{1}{3^4}+\frac{1}{4^4}+\frac{1}{5^4}+\text{etc.}\right) \\[3mm]
&  && &&-\omega^4 && \left(1+\frac{1}{2^5}+\frac{1}{3^5}+\frac{1}{4^5}+\frac{1}{5^5}+\text{etc.}\right) \\
& &&  && &&\text{etc.},
\end{alignat*}
which value is more than apt to find the value of $T$ approximately.\\[2mm]
Therefore, let the term corresponding to the index $\frac{1}{2}$ of the propounded series be in question; if it is put $=T$, it will be

\begin{equation*}
T=1-\frac{2}{3}+\frac{1}{2}-\frac{2}{5}+\frac{1}{3}-\frac{2}{7}+\frac{1}{4}-\frac{2}{9}+\text{etc.}
\end{equation*}
or

\begin{equation*}
T=2 \left(\frac{1}{2}-\frac{1}{3}+\frac{1}{4}-\frac{1}{5}+\frac{1}{6}-\text{etc.}\right),
\end{equation*}
the value of which series is $=2 -2 \log 2$, and so the term corresponding to the index $=\frac{1}{2}$ can be expressed finitely. Therefore, the following terms whose indices are $\frac{1}{2}$, $\frac{3}{2}$, $\frac{5}{2}$, $\frac{7}{2}$ etc. will be expressed this way

\begin{equation*}
\renewcommand{\arraystretch}{1,7}
\setlength{\arraycolsep}{0.0mm}
\begin{array}{lccccccccccc}
\text{Ind.} \quad &1 & 2 & 3 & 4  \\
\text{Term} \quad & 2-2 \log 2, \quad & 2+\dfrac{2}{3}-2 \log 2, \quad & 2+\dfrac{2}{3}+\dfrac{2}{5}-2 \log 2, \quad & 2+\dfrac{2}{3}+\dfrac{2}{5}+\dfrac{2}{7}-2 \log 2 \quad & \text{etc.},
\end{array}
\end{equation*} 

\subsection*{Example 2}

\textit{Let this series be propounded}

\begin{equation*}
\renewcommand{\arraystretch}{1,7}
\setlength{\arraycolsep}{0.0mm}
\begin{array}{ccccccccccc}
1 & 2 & 3 & 4  \\
1, \quad & \left(1+\dfrac{1}{3}\right), \quad & \left(1+\dfrac{1}{3}+\dfrac{1}{5}\right), \quad & \left(1+\dfrac{1}{3}+\dfrac{1}{5}+\dfrac{1}{7} \right) \quad & \text{etc.},
\end{array}
\end{equation*} 
\textit{whose general terms corresponding to fractional indices are to be expressed.}\\[2mm]
Therefore, it will be $a=1$ and $b=2$; hence, if the term corresponding to the integer index $x$ is put

\begin{equation*}
S=1+\frac{1}{3}+\frac{1}{5}+\cdots + \frac{1}{2x-1}
\end{equation*}
and the term corresponding to the fractional index $x+\omega$ is called $=\Sigma$, it will be

\begin{alignat*}{9}
&\Sigma =S &&+ \quad \frac{1}{1+2x}&&+\quad ~\frac{1}{3+2x}&&+~\quad\frac{1}{5+2x}&&+~\quad \frac{1}{7+2x}&&+\text{etc.} \\[3mm]
& &&-\frac{1}{1+2(x+\omega)}&&-\frac{1}{3+2(x+\omega)}&&-\frac{1}{5+2(x+\omega)}&&-\frac{1}{7+2(x+\omega)}&&-\text{etc.}
\end{alignat*}
Since  therefore it suffices to have assigned the terms corresponding to indices smaller than $1$, let $x=0$ and $S=0$; therefore, if the term corresponding to the index $\omega$ is put $=T$, it will be

\begin{alignat*}{9}
& T ~=~ &&+~~~~\frac{1}{1}&&+~~~~\frac{1}{3}&&+~~~~\frac{1}{5}&&+~~~~\frac{1}{7}&&+~~~~\frac{1}{9}&&+\text{etc.} \\[3mm]
& &&-\frac{1}{1+2 \omega}&&-\frac{1}{3+ 2 \omega}&&-\frac{1}{5+2 \omega}&&-\frac{1}{7+2 \omega}&&-\frac{1}{9+ 2 \omega}&&-\text{etc.}
\end{alignat*}
and if $\omega$ is put to denote any number, since $T$ is the term corresponding to the index $\omega$, $T$ will be the general term of the  propounded series which will also be expressed this way

\begin{equation*}
T=\frac{2 \omega}{1(1+2 \omega)}+\frac{2 \omega}{3(3+2 \omega)}+\frac{2 \omega}{5(5+2 \omega)}+\frac{2 \omega}{7(7+2 \omega)}+\text{etc.}
\end{equation*}
or this way

\begin{alignat*}{9}
&T &&= 2 \omega && \left(1+\frac{1}{3^2}+\frac{1}{5^2}+\frac{1}{7^2}+\frac{1}{9^2}+\text{etc.}\right) \\[3mm]
& &&-4 \omega^2 && \left(1+\frac{1}{3^3}+\frac{1}{5^3}+\frac{1}{7^3}+\frac{1}{9^3}+\text{etc.}\right) \\[3mm]
& &&+8 \omega^3 && \left(1+\frac{1}{3^4}+\frac{1}{5^4}+\frac{1}{7^4}+\frac{1}{9^4}+\text{etc.}\right) \\[3mm]
& &&-16 \omega^4 && \left(1+\frac{1}{3^5}+\frac{1}{5^5}+\frac{1}{7^5}+\frac{1}{9^5}+\text{etc.}\right) \\
&   && &&\text{etc.}
\end{alignat*}
Let us put that $\omega= \frac{1}{2}$; therefore, the term corresponding to this index  will be

\begin{equation*}
T=1-\frac{1}{2}+\frac{1}{3}-\frac{1}{4}+\frac{1}{5}-\text{etc.}=\log 2
\end{equation*}
and we will have

\begin{equation*}
\renewcommand{\arraystretch}{2,0}
\setlength{\arraycolsep}{0.0mm}
\begin{array}{lccccccccccc}
\text{Indices} \quad &\dfrac{1}{2} & \dfrac{3}{2} & \dfrac{5}{2} & \dfrac{7}{2}  \\
\text{Terms} \quad & \log 2, \quad & \dfrac{1}{2}+ \log 2, \quad & \dfrac{1}{2}+\dfrac{1}{4}+ \log 2, \quad & \dfrac{1}{2}+\dfrac{1}{4}+\dfrac{1}{6}+ \log 2 \quad & \text{etc.},
\end{array}
\end{equation*} 
If $\omega =\frac{1}{4}$, it will be

\begin{alignat*}{9}
& T &&= 1 &&+\frac{1}{3}&&+\frac{1}{5}&&+\frac{1}{7}&&+\text{etc.} \\[3mm]
&   &&-\frac{2}{3} &&-\frac{2}{7}&&-\frac{2}{11}&&-\frac{2}{15}&&-\text{etc.}
\end{alignat*}
or

\begin{equation*}
T=1-\frac{1}{3}+\frac{1}{5}-\frac{1}{7}+\text{etc.}-\frac{1}{2}\log 2 =\frac{\pi}{4}-\log 2.
\end{equation*}

\paragraph*{§393}

If the general term corresponding to the index $=\frac{1}{2}$ of this general series

\begin{equation*}
\frac{1}{a}, \quad \left(\frac{1}{a}+\frac{1}{a+b}\right), \quad \left(\frac{1}{a}+\frac{1}{a+b}+\frac{1}{a+2b}\right) \quad \text{etc.}
\end{equation*}
is in question,  put $x=0$ in the expressions of the preceding paragraph and $\omega =\frac{1}{2}$ and it will be $S=0$ and the term corresponding to the index $\frac{1}{2}$ and in question will be

\begin{equation*}
\Sigma = \frac{1}{a}-\frac{2}{2a+b}+\frac{1}{a+b}-\frac{2}{2a+3b}+\frac{1}{a+2b}-\frac{2}{2a+5b}+\text{etc.}
\end{equation*}
or, having rendered the terms more uniform, it will be

\begin{equation*}
\frac{1}{2}\Sigma =\frac{1}{2a}-\frac{1}{2a+b}+\frac{1}{2a+2b}-\frac{1}{2a+3b}+\frac{1}{2a+4b}-\text{etc.};
\end{equation*}
since in this series the signs $+$ and $-$ alternate, by taking the continued differences applying the method explained above [§ 8], the value of $\frac{1}{2}\Sigma$ will be expressed by a series converging more rapidly.\\[2mm]
The series of the differences will be

\begin{equation*}
\frac{-b}{2a(2a+b)}, \quad \frac{-b}{(2a+b)(2a+2b)}, \quad \frac{-b}{(2a+2b)(2a+3b)} \quad \text{etc.}
\end{equation*}
\begin{equation*}
\frac{2bb}{2a(2a+b)(2a+2b)}, \quad \frac{2bb}{(2a+b)(2a+2b)(2a+3b)} \quad \text{etc.}
\end{equation*}
\begin{equation*}
\frac{-6b^3}{2a(2a+b)(2b+2a)(2a+3b)} \quad \text{etc.}
\end{equation*}
\begin{equation*}
\text{etc.}
\end{equation*}
From these it is concluded that 

\begin{equation*}
\frac{1}{2}\Sigma =\frac{1}{4a}+\frac{1b}{8a(2a+b)}+\frac{1\cdot 2 bb}{16a (2a+b)(2a+2b)}+\frac{1 \cdot 2 \cdot 3 b^3}{32 a(2a+b)(2a+2b)(2a+3b)}+\text{etc.}
\end{equation*}
And hence one will have

\begin{equation*}
\Sigma =\frac{1}{2a}+\frac{\frac{1}{2}b}{2a(2a+b)}+\frac{\frac{1}{2}\cdot \frac{2}{2}bb}{2a(2a+b)(2a+2b)}+\frac{\frac{1}{2}\cdot \frac{2}{2}\cdot \frac{3}{2}b^3}{2a(2a+b)(2a+2b)(2a+3b)}+\text{etc.},
\end{equation*}
which series is most convergent and exhibits the value of the term $\Sigma$ without much work.

\paragraph*{§394}

But if in general the infinitesimal terms of the series $A$, $B$, $C$, $D$, $E$ etc. vanish and the term corresponding to the index $\omega$ was $=Z$ and the following ones  corresponding to the indices $\omega +1$, $\omega +2$, $\omega +3$ etc. are $Z^{\prime}$, $Z^{\prime \prime}$, $Z^{\prime \prime \prime}$, $Z^{\prime \prime \prime \prime}$ etc., if in the above expressions (§ 391) one puts $x=0$ that  $S=0$ and $X^{\prime}=A$, $X^{\prime \prime}=B$, $X^{\prime \prime \prime}=C$ etc., it will follow, if a series of this kind is given

\begin{equation*}
\renewcommand{\arraystretch}{2,0}
\setlength{\arraycolsep}{0.0mm}
\begin{array}{ccccccccccc}
1  & 2 \quad & 3  & 4  &   \\
\quad A, \quad  & \quad (A+B), \quad  & \quad (A+B+C), \quad & \quad (A+B+C+D) \quad  & \text{etc.}
\end{array}
\end{equation*}
and its term corresponding to the index $\omega$ is put $=\Sigma$, that it will be

\begin{equation*}
\Sigma = (A-Z^{\prime})+(B-Z^{\prime \prime})+(C-Z^{\prime \prime \prime})+(D-Z^{\prime \prime \prime \prime})+\text{etc.},
\end{equation*}
from which expression any intermediate terms can be defined. But it will suffice to have investigated the terms  corresponding to the indices $\omega$ smaller than $1$ for the interpolation. If the term $\Sigma$ corresponding to an arbitrary index $\omega$ of this kind was found and those corresponding to the indices $\omega +1$, $\omega +2$, $\omega +3$ etc. are put $\Sigma^{\prime}$, $\Sigma^{\prime \prime}$, $\Sigma^{\prime \prime \prime}$ etc., it will be

\begin{alignat*}{9}
&\Sigma^{\prime} &&= \Sigma +Z^{\prime} \\
&\Sigma^{\prime \prime} &&= \Sigma +Z^{\prime}+Z^{\prime \prime} \\
&\Sigma^{\prime \prime \prime} &&= \Sigma +Z^{\prime}+Z^{\prime \prime}+Z^{\prime \prime \prime} \\
& && \text{etc.}
\end{alignat*}

\subsection*{Example 1}

\textit{To interpolate this series}

\begin{equation*}
\renewcommand{\arraystretch}{2,0}
\setlength{\arraycolsep}{0.0mm}
\begin{array}{ccccccccccc}
1  & 2 \quad & 3  & 4  &   \\
\quad 1, \quad  & \quad \left(1+\dfrac{1}{4}\right), \quad  & \quad \left(1+\dfrac{1}{4}+\dfrac{1}{9}\right), \quad & \quad \left(1+\dfrac{1}{4}+\dfrac{1}{9}+\dfrac{1}{16}\right) \quad  & \text{etc.}
\end{array}
\end{equation*}
Let $\Sigma$ be the term of this series corresponding to the index $\omega$, and because this series is formed from the summation of this one

\begin{equation*}
1+\frac{1}{4}+\frac{1}{9}+\frac{1}{16}+\frac{1}{25}+\text{etc.},
\end{equation*}
whose term corresponding to the index $\omega$ is $=\frac{1}{\omega^2}$, it will be

\begin{alignat*}{9}
&\Sigma &&=\quad ~~1 &&+\quad ~~ \frac{1}{4}&&+\quad ~~\frac{1}{9}&&+\quad ~~\frac{1}{16}&&+\text{etc.} \\[3mm]
&     &&-\frac{1}{(1+\omega)^2}&&-\frac{1}{(2+\omega)^2}&&-\frac{1}{(3+\omega)^2}&&-\frac{1}{(4+\omega)^2}&&-\text{etc.} \\
\end{alignat*}
If the term corresponding to the index $\frac{1}{2}$ of the propounded series is in question, one will have to put $\omega= \frac{1}{2}$ and it will be

\begin{equation*}
\Sigma =1-\frac{4}{9}+\frac{1}{4}-\frac{4}{25}+\frac{1}{9}-\frac{4}{49}+\text{etc.}
\end{equation*}
or

\begin{equation*}
\Sigma = 4 \left(\frac{1}{4}-\frac{1}{9}+\frac{1}{16}-\frac{1}{25}+\frac{1}{36}-\frac{1}{49}+\text{etc.}\right).
\end{equation*}
Since 

\begin{equation*}
1-\frac{1}{4}+\frac{1}{9}-\frac{1}{25}+\text{etc.}=\frac{\pi^2}{12},
\end{equation*}
it will be

\begin{equation*}
\Sigma = 4 \left(1-\frac{\pi\pi}{12}\right)=4- \frac{1}{3}\pi^2,
\end{equation*}
which is the term corresponding to the index $\frac{1}{2}$. Therefore, we will have the correspondences

\begin{equation*}
\renewcommand{\arraystretch}{2,0}
\setlength{\arraycolsep}{0.0mm}
\begin{array}{lccccccccccc}
\text{Indices} \quad &\dfrac{1}{2} & \dfrac{3}{2} & \dfrac{5}{2}  \\
\text{Terms} \quad & 4-\dfrac{1}{3}\pi^2, \quad & \dfrac{4}{1}+\dfrac{4}{9}-\dfrac{1}{3}\pi^2, \quad & \dfrac{4}{1}+\dfrac{4}{9}+\dfrac{4}{25}-\dfrac{1}{3}\pi^2 \quad  & \text{etc.}.
\end{array}
\end{equation*} 

\subsection*{Example 2}

\textit{To interpolate this series}

\begin{equation*}
\renewcommand{\arraystretch}{2,0}
\setlength{\arraycolsep}{0.0mm}
\begin{array}{ccccccccccc}
1  & 2 \quad & 3  & 4  &   \\
\quad 1, \quad  & \quad \left(1+\dfrac{1}{9}\right), \quad  & \quad \left(1+\dfrac{1}{9}+\dfrac{1}{25}\right), \quad & \quad \left(1+\dfrac{1}{9}+\dfrac{1}{25}+\dfrac{1}{49}\right) \quad  & \text{etc.}
\end{array}
\end{equation*}
Let $\Sigma$ be the term corresponding to the index $\omega$, and because this series is formed from the summation of this one

\begin{equation*}
1+\frac{1}{9}+\frac{1}{25}+\frac{1}{49}+\frac{1}{81}+\text{etc.},
\end{equation*}
from which the term corresponding to the index $\omega$ becomes $Z=\frac{1}{(2 \omega -1)}$, it will be

\begin{equation*}
Z^{\prime}=\frac{1}{(2 \omega +1)^2}, \quad Z^{\prime \prime}=\frac{1}{(2 \omega +3)^2}, \quad Z^{\prime \prime \prime}=\frac{1}{(2 \omega +5)^2} \quad \text{etc.}
\end{equation*}
Therefore, one will have

\begin{alignat*}{9}
& \Sigma &&= \quad ~~ 1 &&+ \quad ~~ \frac{1}{9}&&+\quad ~~ \frac{1}{25}&&+\quad ~~ \frac{1}{49} &&+\text{etc.} \\[3mm]
 &        &&-\frac{1}{(1+2 \omega)^2}&&-\frac{1}{(3- 2\omega)^2}&&-\frac{1}{(5+2 \omega)^2}&&-\frac{1}{(7+2 \omega)^2}&&-\text{etc.}
\end{alignat*}
Let us put $\omega =\frac{1}{2}$ to find the term of the propounded series corresponding to the index $=\frac{1}{2}$ which will be

\begin{equation*}
\Sigma = 1- \frac{1}{4}+\frac{1}{9}-\frac{1}{16}+\frac{1}{25}-\frac{1}{36}+\text{etc.}=\frac{\pi \pi }{12},
\end{equation*}
from which the terms falling in the middle between each two given ones will be expressed as follows. We will have

\begin{equation*}
\renewcommand{\arraystretch}{2,0}
\setlength{\arraycolsep}{0.0mm}
\begin{array}{lccccccccccc}
\text{Indices} \quad & \dfrac{1}{2} & \dfrac{3}{2} & \dfrac{5}{2}  & \dfrac{7}{2} & \quad \text{etc.}\\
\text{Terms} \quad & \quad \dfrac{\pi\pi}{12}, \quad & \dfrac{1}{4}+\dfrac{\pi \pi}{12}, \quad & \dfrac{1}{4}+\dfrac{1}{16}+\dfrac{\pi \pi}{12},  \quad  & \dfrac{1}{4}+\dfrac{1}{16}+\dfrac{1}{36}+\dfrac{\pi \pi}{12} & \quad \text{etc.}.
\end{array}
\end{equation*} 

\subsection*{Example 3}

\textit{To interpolate this series}

\begin{equation*}
\renewcommand{\arraystretch}{2,0}
\setlength{\arraycolsep}{0.0mm}
\begin{array}{ccccccccccc}
1  & 2 \quad & 3  & 4  &   \\
\quad 1, \quad  & \quad \left(1+\dfrac{1}{2^n}\right), \quad  & \quad \left(1+\dfrac{1}{2^n}+\dfrac{1}{3^n}\right), \quad & \quad \left(1+\dfrac{1}{2^n}+\dfrac{1}{3^n}+\dfrac{1}{4^n}\right) \quad  & \text{etc.}
\end{array}
\end{equation*}
Let, as before, $\Sigma$  be the term corresponding to the index $\omega$; it will be $Z=\frac{1}{\omega^n}$ and

\begin{equation*}
Z^{\prime}=\frac{1}{(1+\omega)^n}, \quad Z^{\prime \prime}=\frac{1}{(2+\omega)^n}, \quad Z^{\prime \prime \prime}=\frac{1}{(3+\omega)^n} \quad \text{etc.}
\end{equation*}
And one will have

\begin{alignat*}{9}
&\Sigma &&= \quad ~~ 1&& +\quad ~~\frac{1}{2^n}&&+\quad ~~\frac{1}{3^n}&&+ ~~\quad \frac{1}{4^n}&&+\text{etc.} \\[3mm]
& &&-\frac{1}{(1+\omega)^n}&&-\frac{1}{(2+\omega)^n}&&-\frac{1}{(3+\omega)^n}&&-\frac{1}{(4+\omega)^n}&&-\text{etc.}
\end{alignat*}
If the term corresponding to the index $\frac{1}{2}$ is desired, it will be

\begin{equation*}
=1-\frac{2^n}{3^n}+\frac{1}{2^n}-\frac{2^n}{5^n}+\frac{1}{3^n}-\frac{2^n}{7^n}+\text{etc.}
\end{equation*}
or

\begin{equation*}
=2^n \left(\frac{1}{2^n}-\frac{1}{3^n}+\frac{1}{4^n}-\frac{1}{5^n}+\frac{1}{6^n}-\frac{1}{7^n}+\text{etc.}\right).
\end{equation*}
If one puts

\begin{equation*}
\mathfrak{N}=1-\frac{1}{2^n}+\frac{1}{3^n}-\frac{1}{4^n}+\frac{1}{5^n}-\frac{1}{6^n}+\text{etc.},
\end{equation*}
the term of the propounded series corresponding to the index $\frac{1}{2}$ will be $=2^n(1-\mathfrak{N})$; and hence we will have the correspondences

\begin{equation*}
\renewcommand{\arraystretch}{2,0}
\setlength{\arraycolsep}{0.0mm}
\begin{array}{lccccccccccc}
\text{Indices} \quad & \dfrac{1}{2} & \dfrac{3}{2} & \dfrac{5}{2}  & \quad \text{etc.}\\
\text{Terms} \quad & \quad 2^n-2^n \mathfrak{N}, \quad & 2^n+\dfrac{2^n}{3^n}-2^n \mathfrak{N}, \quad & 2^n+\dfrac{2^n}{3^n}+\dfrac{2^n}{5^n}-2^n \mathfrak{N} & \quad \text{etc.}.
\end{array}
\end{equation*} 

\subsection*{Example 3}

\textit{To interpolate this series}

\begin{equation*}
\renewcommand{\arraystretch}{2,0}
\setlength{\arraycolsep}{0.0mm}
\begin{array}{ccccccccccc}
1  & 2 \quad & 3  & 4  &   \\
\quad 1, \quad  & \quad \left(1+\dfrac{1}{3^n}\right), \quad  & \quad \left(1+\dfrac{1}{3^n}+\dfrac{1}{5^n}\right), \quad & \quad \left(1+\dfrac{1}{3^n}+\dfrac{1}{5^n}+\dfrac{1}{7^n}\right) \quad  & \text{etc.}
\end{array}
\end{equation*}
Let $\Sigma$ be the term corresponding to the index $\omega$ and since $Z=\frac{1}{(2\omega -1)^n}$, it will be

\begin{equation*}
Z^{\prime}=\frac{1}{(2 \omega +1)^n}, \quad Z^{\prime \prime}=\frac{1}{(2 \omega +3)^n}, \quad Z^{\prime \prime \prime}=\frac{1}{(2 \omega +5)^n} \quad \text{etc.}
\end{equation*}
and

\begin{alignat*}{9}
& \Sigma &&= \quad ~~~ 1 &&+ \quad ~~~\frac{1}{3^n} &&+\quad ~~~ \frac{1}{5^n}&&+\quad ~~~ \frac{1}{7^n}&&+\text{etc.} \\[3mm]
& &&-\frac{1}{(1+2 \omega)^n}&&-\frac{1}{(3+2 \omega)^n}&&-\frac{1}{(5+2 \omega)^n}&&-\frac{1}{(7+2 \omega)^n}&&-\text{etc.}
\end{alignat*}
Put $\omega=\frac{1}{2}$ and the term corresponding to the index $\frac{1}{2}$ will result as

\begin{equation*}
=1-\frac{1}{2^n}+\frac{1}{3^n}-\frac{1}{4^n}+\frac{1}{5^n}-\frac{1}{6^n}+\text{etc.}=\mathfrak{N},
\end{equation*}
from which further the remaining terms in the middle between two given ones will be

\begin{equation*}
\renewcommand{\arraystretch}{2,0}
\setlength{\arraycolsep}{0.0mm}
\begin{array}{lccccccccccc}
\text{Indices} \quad & \dfrac{1}{2} & \dfrac{3}{2} & \dfrac{5}{2}  & \quad \text{etc.}\\
\text{Terms} \quad & \quad  \mathfrak{N}, \quad & \dfrac{1}{2^n} +\mathfrak{N}, \quad & \dfrac{1}{2^n}+\dfrac{1}{4^n}+ \mathfrak{N} & \quad \text{etc.}.
\end{array}
\end{equation*} 

\paragraph*{§395}

Now let us put that the infinitesimal terms of the series $A$, $B$, $C$, $D$, $E$ etc., from which the series to be interpolated is formed, do not vanish, but are of such a nature that their differences vanish and let $X$ be the term of this series corresponding to the index $x$ and $Z$ the term corresponding to the exponent $x+\omega$; then, let $X^{\prime}$, $X^{\prime \prime}$, $X^{\prime \prime \prime}$, $X^{\prime \prime \prime \prime}$ etc. be the terms following $X$ and $Z^{\prime}$, $Z^{\prime \prime}$, $Z^{\prime \prime \prime}$, $Z^{\prime \prime \prime \prime}$ etc. the terms following $Z$. Having constituted all this, let this series  be propounded to be interpolated

\begin{equation*}
\renewcommand{\arraystretch}{2,0}
\setlength{\arraycolsep}{0.0mm}
\begin{array}{cccccccccccc}
1 \quad & 2 & 3 & 4  & \quad\\
\quad A \quad & \quad  (A+B), \quad & \quad (A+B+C), \quad & \quad (A+B+C+D) \quad & \quad \text{etc.}.
\end{array}
\end{equation*} 
whose term corresponding to the index $x$ we want to be $=S$, but the term corresponding to the index $x+ \omega$ we want to be $=\Sigma$, and from the results discussed in the preceding chapter it follows

\begin{alignat*}{9}
&\Sigma = S &&+X^{\prime} &&+X^{\prime \prime}&&+X^{\prime \prime \prime} &&+\text{etc.} \\[3mm]
&          &&-Z^{\prime} &&-Z^{\prime \prime} &&-Z^{\prime \prime \prime} &&-\text{etc.}
\end{alignat*}
\begin{equation*}
+\omega X^{\prime}+\omega \left\lbrace
\renewcommand{\arraystretch}{1,3}
\setlength{\arraycolsep}{0.0mm}
\begin{array}{lllllllllll}
& & X^{\prime \prime} & ~+~ & X^{\prime \prime \prime} & ~+~  & X^{\prime \prime \prime \prime} &~+~& \text{etc.}\\
&~-~ & X^{\prime} & ~+~ & X^{\prime \prime } & ~+~  & X^{\prime \prime \prime } &~-~& \text{etc.}\\
\end{array}
\right\rbrace
\end{equation*}
But because, as before, it suffices to have investigated the terms corresponding to indices smaller than $1$, let us put $x=0$ that  $S=0$, $X^{\prime}=A$, $X^{\prime \prime}=B$ etc., and the term corresponding to the index $\omega$ will be

\begin{equation*}
\Sigma = (A-Z^{\prime})+(B-Z^{\prime \prime})+(C-Z^{\prime \prime \prime})+(D-Z^{\prime \prime \prime \prime}) \quad \text{etc.}
\end{equation*} 
\begin{equation*}
\omega A+ \omega ((B-A)+(C-B)+(D-C)+(E-D)+\text{etc.})
\end{equation*}
But if we want to express these differences in the way we did above, i.e. as $\Delta A =B-A$, $\Delta B = C-B$ etc., one will have

\begin{equation*}
\Sigma = (A-Z^{\prime})+(B-Z^{\prime \prime})+(C-Z^{\prime \prime \prime})+(D-Z^{\prime \prime \prime \prime}) \quad \text{etc.}
\end{equation*} 
\begin{equation*}
+\omega (A+\Delta A+\Delta B+\Delta C+\Delta D+\text{etc.})
\end{equation*}

\paragraph*{§396}

But if the infinitesimal terms of the series $A$, $B$, $C$, $D$, $E$ etc., from whose summation the series to be interpolated is formed, neither vanish nor have vanishing first differences, then more series have to added to express the value of $\Sigma$, namely, until one eventually reaches vanishing differences of the infinitesimal terms. For, as before, let the term corresponding to the index $x$ of the series $A$, $B$, $C$, $D$, $E$ etc.  be $=X$ and the following ones $X^{\prime}$, $X^{\prime \prime}$, $X^{\prime \prime \prime}$ etc., but let the term $Z$ correspond to the index $x+\omega$, which term we want to be followed by $Z^{\prime}$, $Z^{\prime \prime}$, $Z^{\prime \prime \prime}$ etc. and let this series be propounded

\begin{equation*}
\renewcommand{\arraystretch}{2,0}
\setlength{\arraycolsep}{0.0mm}
\begin{array}{cccccccccccc}
1 \quad & 2 & 3 & 4  & \quad\\
\quad A \quad & \quad  (A+B), \quad & \quad (A+B+C), \quad & \quad (A+B+C+D) \quad & \quad \text{etc.}
\end{array}
\end{equation*} 
whose term corresponding to the index $x$ we want to be

\begin{equation*}
S=A+B+C+D+\cdots +X,
\end{equation*}
but let the term $\Sigma$ corresponding to the index $x+ \omega$ such that 

\begin{equation*}
\renewcommand{\arraystretch}{2,0}
\setlength{\arraycolsep}{0.0mm}
\begin{array}{l|lllllllllllllll}
\text{to the indices} \quad & \quad \text{these terms correspond}\\
x +\omega +1 \quad & \quad \Sigma' = \Sigma+ Z' \\
x +\omega +2 \quad & \quad \Sigma' = \Sigma+ Z' +Z'' \\
x +\omega +3 \quad & \quad \Sigma' = \Sigma+ Z' +Z'' +Z''' \\
   \text{etc.} & \quad \text{etc.},
\end{array}
\end{equation*}
If now the differences are expressed in such a way that 

\begin{equation*}
\Delta X^{\prime} = X^{\prime \prime}- X^{\prime}, \quad \Delta X^{\prime \prime}=X^{\prime \prime \prime}-X^{\prime \prime}, \quad \Delta X^{\prime \prime \prime}=X^{\prime \prime \prime \prime}- X^{\prime \prime \prime} \quad \text{etc.}
\end{equation*}
\begin{equation*}
\Delta^2 X^{\prime} =\Delta X^{\prime \prime}- \Delta X^{\prime}, \quad \Delta^2 X^{\prime \prime}= \Delta X^{\prime \prime \prime}- \Delta X^{\prime \prime}, \quad \Delta^2 X^{\prime \prime \prime}= \Delta X^{\prime \prime \prime \prime}-\Delta X^{\prime \prime \prime} \quad \text{etc.} 
\end{equation*}
\begin{equation*}
\Delta^3 X^{\prime}=\Delta^2 X^{\prime \prime}-\Delta X^{\prime}, \quad \Delta^3 X^{\prime \prime}= \Delta^2 X^{\prime \prime \prime}-\Delta X^{\prime \prime} \quad \text{etc.},
\end{equation*}
from § 377 the term $\Sigma$ will be expressed the following way:

 \begin{alignat*}{9}
&\Sigma = S &&+X^{\prime} &&+X^{\prime \prime}&&+X^{\prime \prime \prime} &&+\text{etc.} \\
&          &&-Z^{\prime} &&-Z^{\prime \prime} &&-Z^{\prime \prime \prime} &&-\text{etc.}
\end{alignat*}
\begin{equation*}
+\omega (X^{\prime}+\Delta X^{\prime}+\Delta X^{\prime \prime}+\Delta X^{\prime \prime \prime}+\Delta X^{\prime \prime \prime \prime}+\text{etc.})
\end{equation*}
\begin{equation*}
+\frac{\omega (\omega -1)}{1 \cdot 2} (\Delta X^{\prime}+\Delta^2 X^{\prime}+\Delta^2 X^{\prime \prime}+\Delta^2 X^{\prime \prime \prime}+\Delta^2 X^{\prime \prime \prime \prime}+\text{etc.})
\end{equation*}
\begin{equation*}
+\frac{\omega (\omega -1)(\omega -2)}{1 \cdot 2 \cdot 3} (\Delta^2 X^{\prime}+\Delta^3 X^{\prime}+\Delta^3 X^{\prime \prime}+\Delta^3 X^{\prime \prime \prime}+\Delta^3 X^{\prime \prime \prime \prime}+\text{etc.})
\end{equation*}
\begin{equation*}
\text{etc.}
\end{equation*}

\paragraph*{§397}

As we already noted, it suffices to have added so many terms until one gets to vanishing differences of the infinitesimal terms; for, if we also want to continue these series to infinity or at least until the differences of the finite terms vanish, then, because of

\begin{equation*}
Z^{\prime}=X^{\prime}+ \omega \Delta X^{\prime}+\frac{\omega (\omega-1)}{1 \cdot 2}\Delta^2 X^{\prime}+\frac{\omega (\omega -1)(\omega -2)}{1 \cdot 2 \cdot 3}\Delta^3 X^{\prime} +\text{etc.},
\end{equation*} 
the whole found expression will be contracted into

\begin{equation*}
\Sigma =S + \omega X^{\prime}+\frac{\omega(\omega -1)}{1 \cdot 2}\Delta X^{\prime}+\frac{\omega (\omega -1)(\omega -2)}{1 \cdot 2 \cdot 3}\Delta^2 X^{\prime}+\text{etc.},
\end{equation*}
which involves the summatory term of the series $A+B+C+D+\text{etc.}$; but, if this term would be known, the interpolation would have no difficulty. Nevertheless, also this formula can be used, which, if it terminates, exhibits the term to  be interpolated expressed finitely and algebraically; but if it continues to infinity, it is in the most cases more convenient to apply the first formula in which the infinitesimal terms are taken into account. But this formula, if one puts $x=0$ that $\Sigma$ denotes the term corresponding to the index $\omega$, because of $S=0$, takes on this form

 \begin{alignat*}{9}
&\Sigma &&= A &&+B&&+C &&+D &&+\text{etc.} \\
&          &&-Z^{\prime} &&-Z^{\prime \prime} &&-Z^{\prime \prime \prime} &&-Z^{\prime \prime \prime \prime}&&-\text{etc.}
\end{alignat*}
\begin{equation*}
+\omega (A+\Delta A+ \Delta B+ \Delta C+ \Delta D+\text{etc.})
\end{equation*}
\begin{equation*}
+\frac{\omega (\omega -1)}{1 \cdot 2}(\Delta A+ \Delta^2 A +\Delta^2 B+ \Delta^2 C+\Delta^2 D+\text{etc.})
\end{equation*}
\begin{equation*}
+\frac{\omega (\omega -1)(\omega -2	)}{1 \cdot 2 \cdot 3}(\Delta^2 A+\Delta^3 A+ \Delta^3 B+\Delta^3 C +\Delta^3 D +\text{etc.})
\end{equation*}
\begin{equation*}
\text{etc.}
\end{equation*}
Or if, for the sake of brevity, one puts

\begin{equation*}
\omega = \alpha, \quad \frac{\omega(\omega -1)}{1 \cdot 2}=\beta, \quad \frac{\omega(\omega-1)(\omega -2)}{1\cdot 2 \cdot 3}=\gamma \quad \text{etc.},
\end{equation*}
it will be

\begin{alignat*}{9}
&\Sigma &&= \alpha A && +\beta \Delta A &&+ \gamma \Delta^2 A &&+ \delta \Delta^3 A &&+ \text{etc.} \\[3mm]
& &&+A && +\alpha \Delta A &&+\beta \Delta^2 A &&+ \gamma \Delta^3 A &&+\text{etc.} -Z^{\prime}\\[3mm]
& &&+B && +\alpha \Delta B &&+\beta \Delta^2 B &&+ \gamma \Delta^3 B &&+\text{etc.} -Z^{\prime \prime} \\[3mm]
& &&+C && +\alpha \Delta C &&+\beta \Delta^2 C &&+ \gamma \Delta^3 C &&+\text{etc.} -Z^{\prime \prime \prime}\\
& && && &&\text{etc.},
\end{alignat*}
the number of which horizontal series might seem to be infinite, but each of them consists of a finite number of terms.

\subsection*{Example}

\textit{To interpolate this series}

\begin{equation*}
\renewcommand{\arraystretch}{2,0}
\setlength{\arraycolsep}{0.0mm}
\begin{array}{ccccccccccc}
1  & 2 \quad & 3  & 4  &   \\
\quad \dfrac{1}{2}, \quad  & \quad \left(\dfrac{1}{2}+\dfrac{2}{3}\right), \quad  & \quad \left(\dfrac{1}{2}+\dfrac{2}{3}+\dfrac{3}{4}\right), \quad & \quad \left(\dfrac{1}{2}+\dfrac{2}{3}+\dfrac{3}{4}+\dfrac{4}{5}\right) \quad  & \text{etc.}
\end{array}
\end{equation*}
Let the term of this series corresponding to the index $\omega$ be $=\Sigma$, and since it results from the summation of this series $\frac{1}{2}$, $\frac{2}{3}$, $\frac{3}{4}$, $\frac{4}{5}$ etc., it will be $Z=\frac{\omega}{\omega +1}$, and since the infinitesimal terms already have vanishing first differences, only the first differences are to be taken, which, because of

\begin{equation*}
A=\frac{1}{2}, \quad  B=\frac{2}{3}, \quad C=\frac{3}{4}, \quad D=\frac{4}{5} \quad \text{etc.},
\end{equation*} 
will be

\begin{equation*}
\Delta A = \frac{1}{2 \cdot 3}, \quad \Delta B =\frac{1}{3 \cdot 4}, \quad \Delta C =\frac{1}{4 \cdot 5} \quad \text{etc.}
\end{equation*}
Therefore, one will have 

\begin{alignat*}{9}
&\Sigma =\frac{\omega}{2}&&+ \quad \frac{1}{2}&&+\quad \frac{2}{3} &&+\quad\frac{3}{4}&& +\quad\frac{4}{5}&& +\text{etc.} \\[3mm]
&                        &&+~~\frac{\omega}{2 \cdot 3}&&+~~\frac{\omega}{3 \cdot 4} &&+~~\frac{\omega}{4 \cdot 5}&&+~~ \frac{\omega}{5 \cdot 6}&&+\text{etc.} \\[3mm]
& &&-\frac{\omega +1}{\omega +2}&&-\frac{\omega +2}{\omega +3}&&-\frac{\omega +3}{\omega +4}&&-\frac{\omega +4}{\omega +5}&&- \text{etc.}
\end{alignat*}
or, because of

\begin{equation*}
\frac{\omega}{2}+\frac{\omega}{2 \cdot 3}+\frac{\omega}{3 \cdot 4}+\frac{\omega}{4 \cdot 5}+\text{etc.}=\omega,
\end{equation*}
it will be

\begin{alignat*}{9}
&\Sigma =\omega &&+ \quad \frac{1}{2}&&+\quad \frac{2}{3} &&+\quad\frac{3}{4}&& +\quad\frac{4}{5}&& +\text{etc.} \\[3mm]
& &&-\frac{\omega +1}{\omega +2}&&-\frac{\omega +2}{\omega +3}&&-\frac{\omega +3}{\omega +4}&&-\frac{\omega +4}{\omega +5}&&- \text{etc.}
\end{alignat*}
If  the term corresponding to the index $\frac{1}{2}$ is in question, it will be

\begin{equation*}
\Sigma = \frac{1}{2}+\frac{1}{2}-\frac{3}{5}+\frac{2}{3}-\frac{5}{7}+\frac{3}{4}-\frac{7}{9}+\frac{4}{5}-\frac{9}{11}+\text{etc.}
\end{equation*}
or

\begin{equation*}
\Sigma =\frac{1}{2}-\frac{1}{2 \cdot 5}-\frac{1}{3 \cdot 7}-\frac{1}{4 \cdot 9}-\frac{1}{5 \cdot 11}-\frac{1}{6 \cdot 13}-\text{etc.}
\end{equation*}
and hence

\begin{equation*}
\frac{1}{2}\Sigma =\frac{1}{4}-\frac{1}{4 \cdot 5}-\frac{1}{6 \cdot 7}-\frac{1}{8 \cdot 9}-\frac{1}{10 \cdot 11}-\frac{1}{12 \cdot 13}-\text{etc.}
\end{equation*}
or

\begin{alignat*}{9}
&\frac{1}{2}\Sigma =\frac{1}{4}&&-\frac{1}{4}&&-\frac{1}{6}&&-\frac{1}{8}&&-\frac{1}{10}&&-\frac{1}{12}&&-\text{etc.} \\[3mm]
&                              &&+\frac{1}{5}&&+\frac{1}{7}&&+\frac{1}{9}&&+\frac{1}{11}&&+\frac{1}{13}&&+\text{etc.} \\
\end{alignat*}
Therefore, because 

\begin{equation*}
1-\frac{1}{2}+\frac{1}{3}-\frac{1}{4}+\frac{1}{5}-\frac{1}{6}+\text{etc.}= \log 2,
\end{equation*}
it will be

\begin{equation*}
\frac{1}{2}\Sigma = \log 2 -1+\frac{1}{2}-\frac{1}{3}+\frac{1}{4}= \log 2- \frac{7}{12}
\end{equation*}
and hence

\begin{equation*}
\Sigma =2\log 2- \frac{7}{6}.
\end{equation*}

\paragraph*{§398}

Now, let us proceed to the interpolation of series whose terms are conflated of factors, and let this most general series be propounded

\begin{equation*}
\renewcommand{\arraystretch}{1,2}
\setlength{\arraycolsep}{0.0mm}
\begin{array}{ccccccccccc}
1  & 2  & 3  & 4  &   5  \\
\quad A, \quad  & \quad AB, \quad  & \quad ABC, \quad & \quad ABCD \quad  & \quad ABCDE \quad &\text{etc.}
\end{array}
\end{equation*}
whose term corresponding to the index $\omega$ we want to be $=\Sigma$. Therefore, $\log \Sigma$ will be the term corresponding to index $\omega$ in this series

\begin{small}
\begin{equation*}
\renewcommand{\arraystretch}{1,2}
\setlength{\arraycolsep}{0.0mm}
\begin{array}{ccccccccccc}
1  & 2 \quad & 3  & 4  &   \\
\quad \log A, \quad  & \quad \left(\log A+ \log B\right), \quad  & \quad \left(\log A+\log B +\log C\right), \quad & \quad \left(\log A+ \log B+ \log C+\log D\right) \quad  & \text{etc.}
\end{array}
\end{equation*}
\end{small}If we put that the infinitesimal terms of this series vanish and the terms of the series $A$, $B$, $C$, $D$, $E$ etc. corresponding to the index $\omega$ is $Z$ and the following ones corresponding to the indices $\omega +1$, $\omega +2$, $\omega +3$, $\omega +4$ etc.  are $Z^{\prime}$, $Z^{\prime \prime}$, $Z^{\prime \prime \prime}$, $Z^{\prime \prime \prime \prime}$ etc.,  applying the results demonstrated above, it will be

\begin{alignat*}{9}
&\log \Sigma &&= \log A &&+ \log B &&+ \log C &&+ \log D &&+\text{etc.} \\
&           &&- \log Z^{\prime} &&- \log Z^{\prime \prime} &&- \log Z^{\prime \prime \prime} &&- \log Z^{\prime \prime \prime \prime} && -\text{etc.}
\end{alignat*}
Therefore,  going back to numbers, one will have

\begin{equation*}
\Sigma =\frac{A}{Z^{\prime}}\cdot \frac{B}{Z^{\prime \prime}}\cdot \frac{C}{Z^{\prime \prime \prime}}\cdot \frac{D}{Z^{\prime \prime \prime \prime}}\cdot \text{etc.}
\end{equation*}

\paragraph*{§399}

But if the logarithms of the infinitesimal terms of the series $A$, $B$, $C$, $D$ etc. do not vanish, but have vanishing differences, it will, as we saw, be

\begin{alignat*}{9}
&\log \Sigma &&= \log A &&+ \log B &&+ \log C  &&+\text{etc.} \\
&           &&- \log Z^{\prime} &&- \log Z^{\prime \prime} &&- \log Z^{\prime \prime \prime}  && -\text{etc.}
\end{alignat*}
\begin{equation*}
+ \omega \log A +\omega \left(\log \frac{B}{A}+\log \frac{C}{B}+\log \frac{D}{C}+\text{etc.}\right)
\end{equation*}
and hence, going from logarithms to numbers, it will be

\begin{equation*}
\Sigma = A^{\omega} \cdot \frac{A^{1-\omega}B^{\omega}}{Z^{\prime}} \cdot  \frac{B^{1-\omega}C^{\omega}}{Z^{\prime \prime}}\cdot  \frac{C^{1-\omega}D^{\omega}}{Z^{\prime \prime \prime}} \cdot  \frac{D^{1-\omega}E^{\omega}}{Z^{\prime \prime \prime \prime}}\cdot \text{etc.}
\end{equation*}
But if just the second differences of those infinitesimal logarithms vanish, it will be

\begin{equation*}
\log \Sigma = \log A+ \log B+ \log C+ \log D+ \text{etc.}
\end{equation*}
\begin{equation*}
- \log Z^{\prime}- \log Z^{\prime \prime} -\log Z^{\prime \prime \prime}- \log Z^{\prime \prime \prime \prime}-\text{etc.}
\end{equation*}
\begin{equation*}
+\omega \left(\log A +\log \frac{B}{A}+\log \frac{C}{B}+\log \frac{D}{C}+\log \frac{E}{D}+\text{etc.}\right)
\end{equation*}
\begin{equation*}
+\frac{\omega (\omega-1)}{1 \cdot 2}\left(\log \frac{B}{A}+\log \frac{AC}{B^2}+\log \frac{BD}{C^2}+\log \frac{CE}{D^2}+\log \frac{DF}{E^2}+\text{etc.}\right)
\end{equation*}
From these one will therefore obtain

\begin{equation*}
\Sigma = A^{\frac{\omega(3- \omega)}{2}}\cdot B^{\frac{\omega(\omega-1)}{1 \cdot 2}} \cdot \frac{A^{\frac{(\omega -1)(\omega -2)}{1 \cdot 2}} B^{\omega(2-\omega)}C^{\frac{\omega(\omega -1)}{1 \cdot 2}}}{Z^{\prime}} \cdot \frac{B^{\frac{(\omega -1)(\omega -2)}{1 \cdot 2}} C^{\omega(2-\omega)}D^{\frac{\omega(\omega -1)}{1 \cdot 2}}}{Z^{\prime \prime}}\cdot \text{etc.}
\end{equation*}
which, if $\omega <1$, is expressed more conveniently this way

\begin{equation*}
\Sigma = \frac{A^{\frac{\omega(3-\omega)}{1 \cdot 2}}}{B^{\frac{\omega(1-\omega)}{1 \cdot 2}}} \cdot \frac{A^{\frac{(1-\omega)(2-\omega)}{1 \cdot 2}}B^{\omega(2-\omega)}}{C^{\frac{\omega(1-\omega)}{1 \cdot 2}}Z^{\prime}}\cdot \frac{B^{\frac{(1-\omega)(2-\omega)}{1 \cdot 2}}C^{\omega(2-\omega)}}{D^{\frac{\omega(1-\omega)}{1 \cdot 2}}Z^{\prime \prime}}\cdot \text{etc.}
\end{equation*}

\paragraph*{§400}

Let us apply this interpolation to this series

\begin{equation*}
\renewcommand{\arraystretch}{2,0}
\setlength{\arraycolsep}{0.0mm}
\begin{array}{ccccccccccc}
1  & 2 \quad & 3  & 4  &   \\
\quad \dfrac{a}{b}, \quad  & \quad \dfrac{a(a+c)}{b(b+c)}, \quad  & \quad \dfrac{a(a+c)(a+2c)}{b(b+c)(b+2c)}, \quad & \quad \dfrac{a(a+c)(a+2c)(a+3c)}{b(b+c)(b+2c)(b+3c)} \quad  & \text{etc.}
\end{array}
\end{equation*}
whose factors are taken from this series 

\begin{equation*}
\renewcommand{\arraystretch}{2,0}
\setlength{\arraycolsep}{0.0mm}
\begin{array}{ccccccccccc}
1  & 2 \quad & 3  & 4  &   \\
\quad \dfrac{a}{b}, \quad  & \quad \dfrac{a+c}{b+c}, \quad  & \quad \dfrac{a+2c}{b+2c}, \quad & \quad \dfrac{a+3c}{b+3c} \quad  & \text{etc.}
\end{array}
\end{equation*}
whose logarithms of the infinitesimal terms are $=0$. Therefore, it will be

\begin{equation*}
Z=\frac{a-c+c \omega }{b-c+c \omega}, \quad Z^{\prime}=\frac{a+c \omega}{b+c \omega} \quad \text{etc.}
\end{equation*}
Hence, if the term of that series corresponding to the index $\omega$ is put $=\Sigma$,  from § 398 it will be

\begin{equation*}
\Sigma = \frac{a(b+c \omega)}{b(a +c \omega)}\cdot \frac{(a+c)(b+c+c \omega)}{(b+c)(a+c+c \omega)}\cdot \frac{(a+2c)(b+2c+c \omega)}{(b+2c)(a+2c+c \omega)}\cdot \text{etc.}
\end{equation*}
If the term corresponding to the index $\frac{1}{2}$ is desired, having put $\omega =\frac{1}{2}$, it will be

\begin{equation*}
\Sigma = \frac{a(2b+c)}{b(2a+c)}\cdot \frac{(a+c)(2b+3c)}{(b+c)(2a+3c)}\cdot \frac{(a+2c)(2b+5c)}{(b+2c)(2a+5c)}\cdot \text{etc.}
\end{equation*}

\subsection*{Example}

\textit{To interpolate this series}

\begin{equation*}
\renewcommand{\arraystretch}{2,0}
\setlength{\arraycolsep}{0.0mm}
\begin{array}{ccccccccccc}
1  & 2 \quad & 3  & 4  & 5   \\
\quad \dfrac{1}{2}, \quad  & \quad \dfrac{1 \cdot 3}{2 \cdot 4}, \quad  & \quad \dfrac{1 \cdot 3 \cdot 5}{2 \cdot 4 \cdot 6}, \quad & \quad \dfrac{1 \cdot 3 \cdot 5 \cdot 7}{2 \cdot 4 \cdot 6 \cdot 8} \quad  & \dfrac{1 \cdot 3 \cdot 5 \cdot 7 \cdot 9}{2 \cdot 4 \cdot 6 \cdot 8 \cdot 10}\quad &  \text{etc.}
\end{array}
\end{equation*}
Since here  $a=1$, $b=2$ and $c=2$, if the term corresponding to the index $\omega$ is put $=\Sigma$, it will be

\begin{equation*}
\Sigma =\frac{1(2+2 \omega)}{2(1+2\omega)}\cdot \frac{3(4+2 \omega)}{4(3 +2 \omega)}\cdot \frac{5(6+2\omega)}{6(5+2 \omega)}\cdot \frac{7(8 +2 \omega)}{8(7+2 \omega)}\cdot \text{etc.}
\end{equation*}
Hence, if the terms  corresponding to the indices $\omega +1$, $\omega +2$, $\omega +3$ etc. are put $\Sigma^{\prime}$, $\Sigma^{\prime \prime}$, $\Sigma^{\prime \prime \prime}$ etc., it will be

\begin{alignat*}{9}
&\Sigma^{\prime} &&= \frac{1+2 \omega}{2+2 \omega}\cdot \Sigma \\[3mm]
&\Sigma^{\prime \prime} &&= \frac{1+2 \omega}{2+2 \omega}\cdot \frac{3+2 \omega}{4+2 \omega} \cdot \Sigma \\[3mm]
&\Sigma^{\prime \prime \prime} &&= \frac{1+2 \omega}{2+2 \omega}\cdot \frac{3+2 \omega}{4+2 \omega} \cdot \frac{5+2 \omega}{6+2 \omega}\cdot \Sigma \\
& && \text{etc.}
\end{alignat*}
If the term corresponding to the index $\frac{1}{2}$ is desired, having put $\omega=\frac{1}{2}$, it will be

\begin{equation*}
\Sigma =\frac{1 \cdot 3 }{2 \cdot 2} \cdot \frac{3 \cdot 5}{4 \cdot 4} \cdot \frac{5 \cdot 7}{6 \cdot 6}\cdot \frac{7 \cdot 9}{8 \cdot 8}\cdot \frac{9 \cdot 11}{10 \cdot 10}\cdot \text{etc.}
\end{equation*}
But, having put $\pi =$ the half of the circumference of the circle whose radius is $=1$, we showed above that 

\begin{equation*}
\pi =2 \cdot \frac{2 \cdot 2}{1 \cdot 3}\cdot \frac{4 \cdot 4}{3 \cdot 5}\cdot \frac{6 \cdot 6}{5 \cdot 7}\cdot \frac{8 \cdot 8}{7 \cdot 9}\cdot \text{etc.}
\end{equation*}
Therefore, the intermediate terms corresponding to the indices $\frac{1}{2}$, $\frac{3}{2}$, $\frac{5}{2}$ etc. can be expressed by the circumference of the circle this way

\begin{equation*}
\renewcommand{\arraystretch}{2,5}
\setlength{\arraycolsep}{0.0mm}
\begin{array}{lccccccccccc}
\text{Indices:}\quad &  \dfrac{1}{2}  & \dfrac{3}{2} & \dfrac{5}{2}  & \dfrac{7}{2}     \\
\text{Terms:}  &  \quad \dfrac{2}{\pi}, \quad & \quad \dfrac{2}{3}\cdot \dfrac{2}{\pi} \quad  & \quad \dfrac{2 \cdot 4}{3 \cdot 5}\cdot \dfrac{2}{\pi}\quad & \quad \dfrac{2 \cdot 4 \cdot 6}{3 \cdot 5 \cdot 7}\cdot \dfrac{2}{\pi} \quad &  \text{etc.}
\end{array}
\end{equation*}
Wallis found this same interpolation in the \textit{Arithmetica infinitorum}.

\paragraph*{§401}

Now, let us consider this series

\begin{equation*}
\renewcommand{\arraystretch}{1,2}
\setlength{\arraycolsep}{0.0mm}
\begin{array}{ccccccccccc}
1  & 2  & 3  & 4  &   \\
\quad a, \quad  & \quad a(a+b), \quad  & \quad a(a+b)(a+2b), \quad & \quad a(a+b)(a+2b)(a+3b) \quad  &   \text{etc.}
\end{array}
\end{equation*}
whose factors constitute this arithmetic progression

\begin{equation*}
a, \quad (a+b), \quad (a+2b), \quad (a+3b), \quad (a+4b), \quad \text{etc.}
\end{equation*}
whose infinitesimal terms are of such a nature that the differences of their logarithms vanish. Since therefore

\begin{equation*}
Z=a-b+b \omega
\end{equation*}
and

\begin{equation*}
Z^{\prime}=a+b \omega, \quad Z^{\prime \prime}=a+b+b \omega, \quad Z^{\prime \prime \prime}=a+2b+b \omega \quad \text{etc.},
\end{equation*}
if $\Sigma$ denotes the term of the propounded series whose index is $=\omega$, it will be

\begin{equation*}
\Sigma =a^{\omega} \cdot \frac{a^{1-\omega}(a+b)^{\omega}}{a+b \omega}\cdot \frac{(a+b)^{1-\omega}(a+2b)^{\omega}}{a+b+b \omega}\cdot \frac{(a+2b)^{1-\omega}(a+3b)^{\omega}}{a+2b+b \omega}\cdot \text{etc.} 
\end{equation*}
And having found this value, if $\omega$ denotes a certain fractional number smaller than $1$, the following terms corresponding to the indices $1+\omega$, $2+\omega$, $3+\omega$ etc. will be determined in such a way that

\begin{alignat*}{9}
&\Sigma^{\prime} &&= (a+b \omega) \Sigma \\
&\Sigma^{\prime \prime} &&= (a+b \omega)(a+b+b\omega) \Sigma \\
&\Sigma^{\prime \prime \prime} &&= (a+b \omega)(a+b+b\omega)(a+2b+b \omega) \Sigma \\
& &&\text{etc.}
\end{alignat*} 
If the term corresponding to the index $\frac{1}{2}$ is desired, having put $\omega =\frac{1}{2}$, it will be

\begin{equation*}
\Sigma = a^{\frac{1}{2}}\cdot \frac{a^{\frac{1}{2}}(a+b)^{\frac{1}{2}}}{a+\frac{1}{2}b}\cdot \frac{(a+b)^{\frac{1}{2}}(a+2b)^{\frac{1}{2}}}{a+\frac{3}{2}b}\cdot \frac{(a+2b)^{\frac{1}{2}}(a+3b)^{\frac{1}{2}}}{a+\frac{5}{2}b}\cdot \text{etc.}
\end{equation*}
and hence, having taken the squares,

\begin{equation*}
\Sigma^2 = a \cdot \frac{a(a+b)}{(a+\frac{1}{2}b)(a+\frac{1}{2}b)}\cdot \frac{(a+b)(a+2b)}{(a+\frac{3}{2}b)(a+\frac{3}{2}b)}\cdot \frac{(a+2b)(a+3b)}{(a+\frac{5}{2}b)(a+\frac{5}{2}b)}\cdot \text{etc.}
\end{equation*}

\paragraph*{§402}
In the series we treated above [§ 400]
\begin{equation*}
\renewcommand{\arraystretch}{1,2}
\setlength{\arraycolsep}{0.0mm}
\begin{array}{ccccccccccc}
1  & 2  & 3  & 4  &   \\
\quad \dfrac{f}{g}, \quad  & \quad \dfrac{f(f+h)}{g(g+h)}, \quad  & \quad \dfrac{f(f+h)(f+2h)}{g(g+h)(g+2h)}, \quad & \quad \dfrac{f(f+h)(f+2h)(f+3h)}{g(g+h)(g+2h)(g+3h)} \quad  &   \text{etc.}
\end{array}
\end{equation*}
put the term corresponding to the index $\frac{1}{2}$ $=\Theta$; it will be

\begin{equation*}
\Theta =\frac{f(g+\frac{1}{2}h)}{g(f+\frac{1}{2}h)}\cdot \frac{(f+h)(g+\frac{3}{2}h)}{(g+h)(f+\frac{3}{2}h)}\cdot \frac{(f+2h)(g+\frac{5}{2}h)}{(g+2h)(f+\frac{5}{2}h)}\cdot \text{etc.};
\end{equation*}
now, set

\begin{equation*}
f=a, \quad g=a+\frac{1}{2}b \quad \text{and} \quad h=b;
\end{equation*}
it will be

\begin{equation*}
\Theta =\frac{a(a+b)}{(a+\frac{1}{2}b)(a+\frac{1}{2}b)}\cdot \frac{(a+b)(a+2b)}{(a+\frac{3}{2}b)(a+\frac{3}{2}b)}\cdot \text{etc.}
\end{equation*}
and hence it will be $\Sigma^2 = a \Theta $ and $\Sigma =\sqrt{a\Theta}$. Therefore, if the term of this series

\begin{equation*}
\renewcommand{\arraystretch}{2,0}
\setlength{\arraycolsep}{0.0mm}
\begin{array}{ccccccccccc}
1  & 2  & 3  & 4  &   \\
\quad a, \quad  & \quad a(a+b), \quad  & \quad a(a+b)(a+2b), \quad & \quad a(a+b)(a+2b)(a+3b) \quad  &   \text{etc.}
\end{array}
\end{equation*}
corresponding to the index $\frac{1}{2}$ is put $=\Theta$, it will be $\Sigma =\sqrt{a \Theta}$.\\[2mm]
Because here the term of the series of the numerators  corresponding to the index $\frac{1}{2}$ is $=\Sigma$, if in the series of denominators the term corresponding to the index $\frac{1}{2}$ is put $=\Lambda$, it will be $\Theta =\frac{\Sigma}{\Lambda}$; but  $\Theta =\frac{\Sigma^2}{a}$, whence it will be $\Sigma =\frac{a}{\Lambda}$ or $\Sigma \Lambda =a$  by means of which theorems the interpolation of series of this kind is  illustrated very well.

\subsection*{Example 1}

\textit{Let this series  be propounded to be interpolated}

\begin{equation*}
1, \quad 1 \cdot 2, \quad 1 \cdot 2 \cdot 3, \quad 1 \cdot 2 \cdot 3 \cdot 4 \quad \text{etc.}
\end{equation*}
Since here $a=1$ and $b=1$, if the term corresponding to the index $\omega$ is put $=\Sigma$, it will be

\begin{equation*}
\Sigma = \frac{1^{1-\omega}\cdot 2^{\omega}}{1+\omega}\cdot \frac{2^{1-\omega}\cdot 3^{\omega}}{2+\omega}\cdot \frac{3^{1-\omega}\cdot 4^{\omega}}{3+\omega}\cdot \frac{4^{1-\omega}\cdot 5^{\omega}}{4+\omega}\cdot \text{etc.}
\end{equation*}
Here,  one can always take a fraction smaller than $1$ for $\omega$; for, nevertheless the interpolation  extends to the whole series. For, if the terms corresponding $1+\omega$, $2+\omega$, $3+\omega$ etc. are put $\Sigma^{\prime}$, $\Sigma^{\prime \prime}$, $\Sigma^{\prime \prime \prime}$ etc., it will be

\begin{alignat*}{9}
&\Sigma^{\prime} &&= (1+\omega)\Sigma \\
&\Sigma^{\prime \prime} &&= (1+\omega)(2+ \omega)\Sigma \\
&\Sigma^{\prime \prime \prime} &&= (1+\omega)(2+ \omega)(3+\omega)\Sigma \\
& &&\text{etc.}
\end{alignat*}
Therefore, the term of the propounded series corresponding to the index $\frac{1}{2}$ will be

\begin{equation*}
\Sigma =\frac{1^{\frac{1}{2}}\cdot 2^{\frac{1}{2}}}{1^{\frac{1}{2}}}\cdot \frac{2^{\frac{1}{2}}\cdot 3^{\frac{1}{2}}}{2^{\frac{1}{2}}}\cdot \frac{3^{\frac{1}{2}}\cdot 4^{\frac{1}{2}}}{3^{\frac{1}{2}}}\cdot \text{etc.}
\end{equation*}
or

\begin{equation*}
\Sigma^2 = \frac{2 \cdot 4}{3 \cdot 3}\cdot \frac{4 \cdot 6}{5 \cdot 5}\cdot \frac{6 \cdot 8}{7 \cdot 7}\cdot \frac{8 \cdot 10}{9 \cdot 9}\cdot \text{etc.}
\end{equation*}
Therefore, because 

\begin{equation*}
\pi =2 \cdot \frac{2 \cdot 2}{1 \cdot 3}\cdot \frac{4 \cdot 4}{3 \cdot 5}\cdot \frac{6 \cdot 6}{5 \cdot 7}\cdot \text{etc.},
\end{equation*}
it will be

\begin{equation*}
\Sigma^2 =\frac{\pi}{4} \quad \text{and} \quad \Sigma =\frac{\sqrt{\pi}}{2}
\end{equation*}
and hence we have the correspondences

\begin{equation*}
\renewcommand{\arraystretch}{2,5}
\setlength{\arraycolsep}{0.0mm}
\begin{array}{lccccccccccc}
\text{Indices:}\quad &  \dfrac{1}{2}  & \dfrac{3}{2} & \dfrac{5}{2}  & \dfrac{7}{2}  & \quad \text{etc.}   \\
\text{Terms:}  &  \quad \dfrac{\sqrt{\pi}}{2}, \quad & \quad \dfrac{3}{2}\cdot \dfrac{\sqrt{\pi}}{2} \quad  & \quad \dfrac{3 \cdot 5}{2 \cdot 2}\cdot \dfrac{\sqrt{\pi}}{2}\quad & \quad \dfrac{3 \cdot 5 \cdot 7}{2 \cdot 2 \cdot 2}\cdot \dfrac{\sqrt{\pi}}{2} \quad & \quad \text{etc.}
\end{array}
\end{equation*}

\subsection*{Example 2}

\textit{Let this series be propounded to be interpolated}

\begin{equation*}
\renewcommand{\arraystretch}{1,2}
\setlength{\arraycolsep}{0.0mm}
\begin{array}{cccccccccccc}
  1  & 2 & 3  & 4  &   \\
 \quad 1, \quad & \quad 1 \cdot 3, \quad  & \quad 1 \cdot 3 \cdot 5, \quad & \quad 1 \cdot 3 \cdot 5 \cdot 7 \quad & \quad \text{etc.}
\end{array}
\end{equation*}
Since here  $a=1$, $b=2$, if the term corresponding to the index $\omega$ is put $=\Sigma$, it will be

\begin{equation*}
\Sigma = \frac{1^{1-\omega}\cdot 3^{\omega}}{1+2\omega}\cdot \frac{3^{1-\omega}\cdot 5^{\omega}}{3+2\omega}\cdot \frac{5^{1-\omega}\cdot 7^{\omega}}{5+2\omega}\cdot \text{etc.}
\end{equation*}
the term following in order will be of such a nature:

\begin{alignat*}{9}
&\Sigma^{\prime} &&= (1+2\omega) \Sigma \\
&\Sigma^{\prime \prime} &&= (1+2\omega)(3+ 2 \omega) \Sigma \\
&\Sigma^{\prime \prime \prime} &&= (1+2\omega)(3+ 2 \omega)(5+2 \omega) \Sigma \\
& &&\text{etc.}
\end{alignat*}
Therefore, if the term corresponding to the index $\frac{1}{2}$ of the propounded series  is desired and it is called $=\Sigma$, it will be

\begin{equation*}
\Sigma = \frac{\sqrt{1 \cdot 3}}{2}\cdot \frac{\sqrt{3 \cdot 5}}{4}\cdot \frac{\sqrt{5 \cdot 7}}{6}\cdot \frac{\sqrt{7 \cdot 9}}{8}\cdot \text{etc.},
\end{equation*}
therefore,

\begin{equation*}
\Sigma^2 = \frac{1 \cdot 3}{2 \cdot 2}\cdot \frac{3 \cdot 5}{4 \cdot 4}\cdot \frac{5 \cdot 7}{6 \cdot 6}\cdot \frac{7 \cdot 9}{8 \cdot 8}\cdot \text{etc.}=\frac{2 }{\pi}
\end{equation*}
and one will have $\Sigma =\sqrt{\frac{2}{\pi}}$. But we will have the correspondences

\begin{equation*}
\renewcommand{\arraystretch}{2,5}
\setlength{\arraycolsep}{0.0mm}
\begin{array}{lccccccccccc}
\text{Indices:}\quad &  \dfrac{1}{2}  & \dfrac{3}{2} & \dfrac{5}{2}  & \dfrac{7}{2}  & \quad \text{etc.}   \\
\text{Terms:}  &  \quad \sqrt{\dfrac{2}{\pi}}, \quad & \quad 2 \cdot \sqrt{\dfrac{2}{\pi}} \quad  & \quad 2 \cdot4 \sqrt{\dfrac{2}{\pi}}\quad & \quad 2 \cdot 4 \cdot 6 \sqrt{\dfrac{2}{\pi}} \quad & \quad \text{etc.}
\end{array}
\end{equation*}
If the first series and this one are multiplied by each other that one has this series

\begin{equation*}
\renewcommand{\arraystretch}{2,0}
\setlength{\arraycolsep}{0.0mm}
\begin{array}{ccccccccccc}
1  & 2  & 3  & 4  & 5   \\
\quad 1^2, \quad  & \quad 1^2\cdot 2 \cdot 3, \quad  & \quad 1^2 \cdot 2 \cdot 3^2 \cdot 5, \quad & \quad 1^2 \cdot 2 \cdot 3^2 \cdot 4 \cdot 5 \cdot 7 \quad  & 1^2 \cdot 2 \cdot 3^2 \cdot 4 \cdot 5^2 \cdot 7 \cdot 9 & \quad  \text{etc.},
\end{array}
\end{equation*}
the term corresponding to the index $\frac{1}{2}$ will be $=\frac{\sqrt{\pi}}{2}\cdot \sqrt{\frac{2}{\pi}}=\frac{1}{\sqrt{2}}$; this is easily seen, if this form is attributed to that series

\begin{equation*}
\renewcommand{\arraystretch}{2,0}
\setlength{\arraycolsep}{0.0mm}
\begin{array}{ccccccccccc}
1  & 2  & 3  & 4  &   \\
\quad \dfrac{1 \cdot 2}{2}, \quad  & \quad \dfrac{1 \cdot 2 \cdot 3 \cdot 4}{2^2}, \quad  & \quad \dfrac{1 \cdot 2 \cdot 3 \cdot 4 \cdot 5 \cdot 6}{2^3}, \quad & \quad \dfrac{1 \cdot 2 \cdot 3 \cdot 4 \cdot 5 \cdot 6 \cdot 7 \cdot 8}{2^4} \quad  &   \text{etc.},
\end{array}
\end{equation*}
whose term corresponding to the index $\frac{1}{2}$ obviously is $=\frac{1}{\sqrt{2}}$.

\subsection*{Example 3}

\textit{Let this series be propounded to be interpolated}

\begin{equation*}
\renewcommand{\arraystretch}{1,6}
\setlength{\arraycolsep}{0.0mm}
\begin{array}{ccccccccccc}
1  & 2  & 3  & 4  &   \\
\quad \dfrac{n}{1}, \quad  & \quad \dfrac{n(n-1)}{1 \cdot 2}, \quad  & \quad \dfrac{n(n-1)(n-2)}{1 \cdot 2 \cdot 3}, \quad & \quad \dfrac{n(n-1)(n-2)(n-3)}{1 \cdot 2 \cdot 3 \cdot 4} \quad  &   \text{etc.}
\end{array}
\end{equation*}
Consider the numerators and denominators of this series separately, and since the numerators are

\begin{equation*}
\renewcommand{\arraystretch}{1,2}
\setlength{\arraycolsep}{0.0mm}
\begin{array}{ccccccccccc}
1  & 2  & 3  & 4  &   \\
\quad n, \quad  & \quad n(n-1), \quad  & \quad n(n-1)(n-2), \quad & \quad n(n-1)(n-2)(n-3) \quad  &   \text{etc.},
\end{array}
\end{equation*} 
 after an application of the result before it will be $a=n$ and $b=-1$, whence the term corresponding to the index $\omega$ of this series will be

\begin{equation*}
=n^{\omega} \cdot \frac{n^{1- \omega}(n-1)^{\omega}}{n-\omega} \cdot \frac{(n-1)^{1- \omega}(n-2)^{\omega}}{n-1-\omega} \cdot \frac{(n-2)^{1- \omega}(n-3)^{\omega}}{n-2-\omega}\cdot \text{etc.},
\end{equation*}
which expression, because of the factors going over into negatives, does not show anything certain. Therefore,  for the sake of brevity having put $1 \cdot 2 \cdot 3 \cdots n = N$, transform the propounded series into this one

\begin{equation*}
\renewcommand{\arraystretch}{2,0}
\setlength{\arraycolsep}{0.0mm}
\begin{array}{ccccccccccc}
1  & 2  & 3  &    \\
\quad \dfrac{N}{1 \cdot 1 \cdot 2 \cdot 3 \cdots (n-1)}, \quad  & \quad \dfrac{N}{1 \cdot 2 \cdot 1 \cdot 2 \cdot 3 \cdots (n-2)}, \quad  & \quad \dfrac{N}{1 \cdot 2 \cdot 3 \cdot 1 \cdot 2 \cdot 3 \cdots (n-3)} \quad &   \text{etc.};
\end{array}
\end{equation*}
since its denominators consist of two factors, the one group will constitute this series
\begin{equation*}
\renewcommand{\arraystretch}{2,0}
\setlength{\arraycolsep}{0.0mm}
\begin{array}{ccccccccccc}
1  & 2  & 3  & 4  &   \\
\quad 1 \cdot 2 \cdot 3 \cdots (n-1), \quad  & \quad 1 \cdot 2 \cdot 3 \cdots (n-2), \quad  & \quad 1 \cdot 2 \cdot 3 \cdots (n-3)  &   \quad\text{etc.}
\end{array}
\end{equation*}
whose term corresponding to the index $\omega$ coincides with term of this series

\begin{equation*}
\renewcommand{\arraystretch}{1,5}
\setlength{\arraycolsep}{0.0mm}
\begin{array}{ccccccccccc}
1  & 2  & 3  & 4  & 5  \\
\quad 1, \quad  & \quad 1 \cdot 2, \quad  & \quad 1 \cdot 2 \cdot 3, \quad & \quad 1 \cdot 2 \cdot 3 \cdot 4 \quad  & \quad 1 \cdot 2 \cdot 3 \cdot 4 \cdot 5 \quad & \quad  \text{etc.}
\end{array}
\end{equation*}
corresponding to the index $n-\omega$ which is

\begin{equation*}
\frac{1^{1-n- \omega}\cdot 2^{n-\omega}}{1+n-\omega} \cdot \frac{2^{1-n- \omega}\cdot 3^{n-\omega}}{2+n-\omega} \cdot \frac{3^{1-n- \omega}\cdot 4^{n-\omega}}{3+n-\omega}\cdot \text{etc.}
\end{equation*}
But let the term of this series corresponding to the index $1- \omega$ be $=\Theta$; it will be

\begin{equation*}
\Theta =\frac{1^{\omega}\cdot 2^{1-\omega}}{2- \omega} \cdot \frac{2^{\omega}\cdot 3^{1-\omega}}{3- \omega}\cdot \frac{3^{\omega}\cdot  4^{1-\omega}}{4- \omega}\cdot \text{etc.},
\end{equation*}
and, because of the correspondences

\begin{equation*}
\renewcommand{\arraystretch}{2,0}
\setlength{\arraycolsep}{0.0mm}
\begin{array}{lccccccccccc}
\text{Indices:}\quad &  1- \omega  & 2- \omega & 3- \omega    & \quad \text{etc.}   \\
\text{Terms:}  &  \quad \Theta, \quad & \quad (2- \omega)\Theta \quad  & \quad (2- \omega)(3- \omega)\Theta \quad & \quad \text{etc.},
\end{array}
\end{equation*}
 the following term will correspond to the index $n- \omega$

\begin{equation*}
(2- \omega)(3- \omega)(4- \omega) \cdots (n -\omega)\Theta.
\end{equation*}
Further, the other factors of those denominators will constitute this series

\begin{equation*}
\renewcommand{\arraystretch}{1,5}
\setlength{\arraycolsep}{0.0mm}
\begin{array}{ccccccccccc}
1  & 2  & 3  & 4  & 5  \\
\quad 1, \quad  & \quad 1 \cdot 2, \quad  & \quad 1 \cdot 2 \cdot 3, \quad & \quad 1 \cdot 2 \cdot 3 \cdot 4 \quad  & \quad 1 \cdot 2 \cdot 3 \cdot 4 \cdot 5 \quad & \quad  \text{etc.};
\end{array}
\end{equation*}
if the term corresponding to the index $\omega$ is put $=\Lambda$, it will be

\begin{equation*}
\Lambda = \frac{1^{1-\omega} \cdot 2^{\omega}}{1+\omega}\cdot \frac{2^{1-\omega} \cdot 3^{\omega}}{2+\omega}\cdot \frac{3^{1-\omega} \cdot 4^{\omega}}{3+\omega}\cdot \text{etc.}  
\end{equation*}
Having found these, if the term of the propounded series itself

\begin{equation*}
\renewcommand{\arraystretch}{2,0}
\setlength{\arraycolsep}{0.0mm}
\begin{array}{ccccccccccc}
1  & 2  & 3  & 4  &   \\
\quad \dfrac{n}{1}, \quad  & \quad \dfrac{n(n-1)}{1 \cdot 2}, \quad  & \quad \dfrac{n(n-1)(n-2)}{1 \cdot 2 \cdot 3}, \quad & \quad \dfrac{n(n-1)(n-2)(n-3)}{1 \cdot 2 \cdot 3 \cdot 4} \quad  &   \text{etc.}
\end{array}
\end{equation*}
corresponding to the index $\omega$ is put $\Sigma$, it will be

\begin{equation*}
\Sigma =\frac{N}{\Lambda \cdot (2- \omega)(3- \omega)(4- \omega)\cdots (n- \omega)\Theta}.
\end{equation*}
But on the other hand

\begin{equation*}
\frac{N}{(2- \omega)(3- \omega)(4- \omega)\cdots (n- \omega)}=\frac{2}{2- \omega}\cdot \frac{3}{3- \omega}\cdot \frac{4}{4- \omega}\cdots \frac{n}{n- \omega}
\end{equation*}
and 

\begin{equation*}
\Lambda \Theta = \frac{1 \cdot 2}{(1+\omega)(2-\omega)}\cdot \frac{2 \cdot 3}{(2+ \omega)(3- \omega)}\cdot \frac{3 \cdot 4}{(3+\omega)(4-\omega)}\cdot \text{etc.}
\end{equation*}
From these the term in question corresponding to the index $\omega$ will be

\begin{equation*}
\Sigma = \frac{2}{2- \omega}\cdot \frac{3}{3- \omega}\cdot \frac{4}{4- \omega}\cdot \frac{5}{5-\omega}\cdots \frac{n}{n-\omega}
\end{equation*}
\begin{equation*}
\cdot \frac{(1+\omega)(2- \omega)}{1 \cdot 2}\cdot \frac{(2+ \omega)(3- \omega)}{2 \cdot 3}\cdot \frac{(3+ \omega)(4- \omega)}{3 \cdot 4}\cdot \text{etc. up to infinity.}
\end{equation*}
Therefore,  the following term will correspond to the index $\frac{1}{2}$

\begin{equation*}
\frac{4}{3}\cdot \frac{6}{5}\cdot \frac{8}{7}\cdot \frac{10}{9}\cdot \frac{12}{11} \cdots \frac{2n}{2n-1}\cdot \frac{3 \cdot 3}{2 \cdot 4}\cdot \frac{5 \cdot 5}{4 \cdot 6} \cdot \frac{7 \cdot 7}{6 \cdot 8} \cdot \frac{9 \cdot 9}{8 \cdot 10}\cdot \text{etc.},
\end{equation*}
which is reduced to

\begin{equation*}
\frac{4 \cdot 6 \cdot 8 \cdot 10 \cdots 2n}{3 \cdot 5 \cdot 7 \cdot 9 \cdots (2n-1)}\cdot \frac{4}{\pi},
\end{equation*}
or it will be

\begin{equation*}
=\frac{2 }{\pi} \cdot \frac{2 \cdot 4 \cdot 6 \cdot 8 \cdot 10 \cdots 2n}{1 \cdot 3 \cdot 5 \cdot 7 \cdot 9 \cdots (2n-1)}.
\end{equation*}
If $n=2$, this series to be interpolated will result

\begin{equation*}
\renewcommand{\arraystretch}{1,5}
\setlength{\arraycolsep}{0.0mm}
\begin{array}{ccccccccccc}
0 & 1  & 2  & 3  & 4  & 5 & 6   \\
\quad 1, \quad  & \quad 2, \quad  & \quad 1, \quad & \quad 0, \quad  & \quad 0, \quad & \quad 0, \quad & \quad 0, \quad  & \quad   \text{etc.}
\end{array}
\end{equation*}
whose term corresponding to the index $\frac{1}{2}$ therefore is $=\frac{16}{3 \pi}$.

\subsection*{Example 4}

\textit{Let the term corresponding to the index $=\frac{1}{2}$ in this series be in question}

\begin{equation*}
\renewcommand{\arraystretch}{2,0}
\setlength{\arraycolsep}{0.0mm}
\begin{array}{ccccccccccc}
0& 1  & 2  & 3  & 4  &   \\
\quad 1, \quad  & \quad +\dfrac{1}{2}, \quad  & \quad -\dfrac{1 \cdot 1}{2 \cdot 4}, \quad & \quad +\dfrac{1 \cdot 1 \cdot 3}{2 \cdot 4 \cdot 6} \quad  & \quad - \dfrac{1 \cdot 1 \cdot 3 \cdot 5}{2 \cdot 4 \cdot 6 \cdot 8} \quad & \quad   \text{etc.}
\end{array}
\end{equation*}
This series results from the preceding, if one puts $n=\frac{1}{2}$, and therefore the term in question, which we want to be $=\Sigma$, will be

\begin{equation*}
\Sigma = \frac{2}{\pi}\cdot \frac{2 \cdot 4 \cdot 6 \cdot 8 \cdots 2n}{1 \cdot 3 \cdot 5 \cdot 7 \cdots (2n-1)},
\end{equation*}
having put $n=\frac{1}{2}$. Put 

\begin{equation*}
\frac{2 \cdot 4 \cdot 6 \cdot 8 \cdots 2n}{1 \cdot 3 \cdot 5 \cdot 7 \cdots (2n-1)}= \Theta,
\end{equation*}
if $n=\frac{1}{2}$, and $\Theta $ will be the term corresponding to the index $\frac{1}{2}$ in this series

\begin{equation*}
\frac{2}{1}, \quad \frac{2 \cdot 4}{1 \cdot 3}, \quad \frac{2 \cdot 4 \cdot 6}{1 \cdot 3 \cdot 5}, \quad \frac{2 \cdot 4 \cdot 6 \cdot 8}{1 \cdot 3 \cdot 5 \cdot 7} \quad \text{etc.},
\end{equation*}
which results from the above ones as $=\frac{\pi}{2}$. Therefore, the term of the propounded series corresponding to the index $\frac{1}{2}$ which is in question will be $=1$. But since in this series, if the term corresponding to any arbitrary index $\omega$ is put $=\Sigma$, the one following it will be $\Sigma^{\prime}= \frac{1- 2 \omega}{2+2 \omega}\Sigma$, the propounded series will be interpolated  by terms falling in the respective middles this way:

\begin{equation*}
\renewcommand{\arraystretch}{2,0}
\setlength{\arraycolsep}{0.0mm}
\begin{array}{lccccccccccc}
\text{Indices:}\quad &  0  & \dfrac{1}{2} & 1    & \dfrac{3}{2} & 2 & \dfrac{5}{2} & 3 & \dfrac{7}{2} & \quad \text{etc.}   \\
\text{Terms:}  &  \quad 1, \quad & \quad 1, \quad  &\quad  \dfrac{1}{2} \quad &\quad 0, \quad & \quad \dfrac{-1 \cdot 1}{2 \cdot 4}, \quad & \quad 0, \quad & \dfrac{1 \cdot 1 \cdot 3}{2 \cdot 4 \cdot 6}, \quad & \quad 0 \quad & \quad \text{etc.},
\end{array}
\end{equation*}

\subsection*{Example 5}

\textit{If $n$ was an arbitrary fractional number, to find the term corresponding to the index $\omega$ in the series}

\begin{equation*}
\renewcommand{\arraystretch}{2,0}
\setlength{\arraycolsep}{0.0mm}
\begin{array}{ccccccccccc}
0& 1  & 2  & 3  & 4  &   \\
\quad 1, \quad & \quad \dfrac{n}{1}, \quad  & \quad \dfrac{n(n-1)}{1 \cdot 2}, \quad  & \quad \dfrac{n(n-1)(n-2)}{1 \cdot 2 \cdot 3}, \quad & \quad \dfrac{n(n-1)(n-2)(n-3)}{1 \cdot 2 \cdot 3 \cdot 4} \quad  &   \text{etc.}
\end{array}
\end{equation*}
If we compare the expression

\begin{equation*}
\frac{2}{2- \omega}\cdot \frac{3}{3- \omega}\cdot \frac{4}{4- \omega}\cdots \frac{n}{n- \omega}
\end{equation*}
to  § 400, we have $a=1$, $c=1$, $b=1- \omega$ and, having written $n$ instead of $\omega$ there, it will be

\begin{equation*}
\frac{1}{1- \omega}\cdot \frac{2}{2- \omega}\cdot \frac{3}{3- \omega}\cdots \frac{n}{n- \omega}= \frac{1(1- \omega +n)}{(1- \omega)(1+n)}\cdot \frac{2(2- \omega +n)}{(2- \omega)(2+n)}\cdot \text{etc.},
\end{equation*}
whence the term in question corresponding to the index $\omega$, if that term is put $=\Sigma$,  will be

\begin{equation*}
\Sigma = \frac{(1- \omega +n)\cdot 2}{(1+n)(2- \omega)}\cdot \frac{(2- \omega +n)3}{(2+n)(3- \omega)}\cdot \text{etc.} \cdot \frac{(1+\omega)(2- \omega)}{1 \cdot 2}\cdot \frac{(2+ \omega)(3- \omega)}{2 \cdot 3} \cdot \text{etc.}
\end{equation*}
and hence

\begin{equation*}
\Sigma =\frac{(1+\omega)(1+n- \omega)}{1(1+n)} \cdot \frac{(2+ \omega)(2+n- \omega)}{2(2+n)}\cdot \frac{(3+ \omega)(3+n- \omega)}{3(3+n)}\cdot \text{etc.}.;
\end{equation*}
therefore, if $n- \omega$ was an integer number, the value of $\Sigma$ can be expressed rationally.

\begin{alignat*}{9}
&\text{So if} \quad n=\omega, \quad && \text{it will be} \quad \Sigma =1; \\
&\text{if} \quad n=1+\omega, \quad && \text{it will be} \quad \Sigma =n; \\
&\text{if} \quad n=2+\omega, \quad && \text{it will be} \quad \Sigma =\frac{n(n-1)}{1\cdot 2}; \\[3mm]
&\text{if} \quad n=3+\omega, \quad && \text{it will be} \quad \Sigma =\frac{n(n-1)(n-2)}{1 \cdot 2 \cdot 3} \\
& && \text{etc.}
\end{alignat*}
But if $\omega -n$ was an positive integer, it will always be $\Sigma =0.$
\end{document}


\title{On the transformation of series
\footnote{
Original title: "' De Transformatione serierum"', first published as part of the book "'\textit{Institutiones calculi differentialis cum eius usu in analysi finitorum ac doctrina serierum}"', 1755, reprinted in \textit{Opera Omnia}: Series 1, Volume 10, pp. 217 - 234, Eneström-Number E212, translated by: Alexander Aycock for the project "`Euler-Kreis Mainz"'}}
\author{Leonhard Euler}
\date{ }
\maketitle

\paragraph*{§1}

Since it is propounded to us to show the use of differential calculus, both in the whole field of analysis and in the doctrine of series, several auxiliary tools from common algebra which are usually are not discussed will have to be covered here. Although we have already covered a huge part  in the \textit{Introductio}, some things were nevertheless left aside there, either on purpose, because it is convenient to explain them just then when they are actually needed, or because all the things which will be necessary could not have been foreseen at that point. This concerns the transformation of series we devote this chapter to and  by means of which a given series is transformed into innumerable others such that, if the sum of the propounded series is known, the resulting ones can all be summed at the same time. Indeed, having discussed this subject in advance, we will be able to develop the doctrine of series even further by means of differential and integral calculus.

\paragraph*{§2}

But we will mainly consider series whose terms are multiplied by successive powers of a certain variable quantity, since these extend further and are of greater utility.\\[2mm]
Therefore, let the following general series be propounded, whose sum, either known or not, we want to put $=S$, and let

\begin{equation*}
S= ax+bx^2+cx^3+dx^4+ex^5+\text{etc.}
\end{equation*} 
Now, put $x=\frac{y}{1+y}$ and because by an  infinite series 

\begin{equation*}
\renewcommand{\arraystretch}{1,5}
\setlength{\arraycolsep}{0.0mm}
\begin{array}{rlrlrrrrrrrrrrrrr}
&x   &~=~&y    &~-& y^2&~+~&  y^3&~-~ & y^4&~+~ &y^5&~-~&y^6&~+~&\text{etc.}\\
&x^2 &~=~&y^2  &~-~&2y^3&~+~&3 y^4&~-~&4 y^5&~+~&5y^6&~-~&6y^7&~+~&\text{etc.}\\
&x^3 &~=~&y^3  &~-~&3y^4&~+~&6 y^5&~-~&10y^6&~+~&15y^7&~-~&21y^8&~+~&\text{etc.}\\
&x^4 &~=~&y^4  &~-~&4y^4&~+~&10y^6&-~&20y^7&~+~&35y^8&~-~&35y^8&~+~&\text{etc.}\\
&    &&      &&     &&      && \text{etc.},
\end{array}
\end{equation*}
these values, having substituted them and having arranged the series according to powers of $y$, will give

\begin{equation*}
\renewcommand{\arraystretch}{1,0}
\setlength{\arraycolsep}{0.0mm}
\begin{array}{rrrrrrrrrrrrrrrrrrrrrrrrr}
&S &~=~&a&y &~-~&a&y^2 &~+~&a&y^3&~-~&a&y^4&~+~&a&y^5 & \quad \text{etc.} \\
&  &   & &  &~+~&b&    &~-~&2b&  &~+~&3b&  &~-~&4b&   & \\
&  &   & &  &   & &    &~+~&c&   &~-~&3c&  &~+~&6c&   & \\
&  &   & &  &   & &    &   & &   &~+~&d &  &~-~&3d&   & \\
&  &   & &  &   & &    &   & &   &   &  &  &~+~&e&   & \\
\end{array}
\end{equation*}

\paragraph*{§3}

Since we put $x=\frac{y}{1+y}$, it will be $y=\frac{x}{1-x}$; having substituted this value for $y$, the propounded series

\begin{equation*}
S=ax+bx^2+cx^3+dx^4+ex^5+\text{etc.}
\end{equation*}
will be transformed into this one:

\begin{equation*}
S=a \frac{x}{1-x}+(b-a)\frac{x^2}{(1-x)^2}+(c-2b+a)\frac{x^3}{(1-x)^3}+\text{etc.},
\end{equation*}
in which the coefficient of the second term $b-a$ is the first difference of $a$ from the series $a$, $b$, $c$, $d$, $e$ etc., which difference we denoted by $\Delta a$ above; the coefficient of the third term $c-2b+a$ is the second difference $\Delta^2 a$; the coefficient of the fourth term is the third difference of $\Delta^3 a$ etc. Therefore, using the iterated differences of $a$ which are formed  from the series $a$, $b$, $c$, $d$, $e$ etc., the transformed series will go over into this one

\begin{equation*}
S=\frac{x}{1-x}a+\frac{x^2}{(1-x)^2}\Delta a+\frac{x^3}{(1-x)^3}\Delta^2 a+\frac{x^4}{(1-x)^4}\Delta^3 a+\text{etc.},
\end{equation*} 
the sum of which series is therefore known, if the sum of the propounded series was known.

\paragraph*{§4}

Therefore, if the series $a$, $b$, $c$, $d$ etc. was of such a nature that it finally leads to constant differences, what happens, if its general term was a polynomial, the series $\frac{x}{1-x}a+\frac{x^2}{(1-x)^2}\Delta a+$etc.  will have terms vanishing eventually and hence its sum can be exhibited by a finite expression. Therefore, if the first differences of the series $a$, $b$, $c$, $d$ etc.  were already constant, the sum of this series $ax+bx^2+cx^3+dx^4+\text{etc.}$ will be

\begin{equation*}
=\frac{x}{1-x}a+\frac{x^2}{(1-x)^2}\Delta a.
\end{equation*}
But if just the second differences of the coefficients of that series become constant, the sum of the propounded series will be

\begin{equation*}
=\frac{x}{1-x}a+\frac{x^2}{(1-x)^2}\Delta a+\frac{x^3}{(1-x)^3}\Delta \Delta a.
\end{equation*}
Therefore, the sums of series of this kind are easily found from the differences of the coefficients.\\[2mm]

I. \textit{Let the sum of this series be in question}

\begin{equation*}
1x+3x^2+5x^3+7x^4+9x^5+\text{etc.},
\end{equation*}
\begin{equation*}
\text{Diff. I} \quad 2, \quad 2, \quad 2,  \quad 2 \quad \text{etc.}
\end{equation*}
Therefore, since the first differences are constant, because of $a=1$ and $\Delta a=2$, the sum of the propounded series will be

\begin{equation*}
=\frac{x}{1-x}+\frac{2xx}{(1-x)^2}=\frac{x+xx}{(1-x)^2}.
\end{equation*}

II. \textit{Let the sum of this series be in question}

\begin{equation*}
1x+4xx+9x^3+16x^4+25x^5+\text{etc.}
\end{equation*}
\begin{alignat*}{19}
&\text{Diff. I} \quad && 3, \quad && \quad &&5, \quad && \quad &&7, \quad && \quad && 9, \quad && \text{etc.} \\
&\text{Diff. II} \quad &&  \quad && 2 \quad && \quad && 2 \quad && \quad && 2 \quad && \text{etc.} \\
\end{alignat*}
Therefore, since $a=1$, $\Delta a=3$, $\Delta^2 a=2$, the sum of the propounded series will be

\begin{equation*}
=\frac{x}{1-x}+\frac{3xx}{(1-x)^2}+\frac{2x^3}{(1-x)^3}=\frac{x+xx}{(1-x)^3}.
\end{equation*}

III. \textit{Let the sum of this series be in question}

\begin{equation*}
S=4x+15x^2+40x^3+85x^4+156x^5+259x^6+\text{etc.}
\end{equation*}
\begin{alignat*}{19}
&\text{Diff. I} \quad && 11, \quad && \quad &&25, \quad && \quad &&45, \quad && \quad && 71, \quad && \quad && 103 \quad && \quad && \text{etc.} \\ 
&\text{Diff. II} \quad &&  \quad &&14, \quad && \quad && 20, \quad && \quad && 26, \quad &&  \quad && 32, \quad && \quad &&    \text{etc.} \\ 
&\text{Diff. III} \quad &&  \quad && \quad &&6, \quad && \quad &&6, \quad && \quad &&6,  \quad && \quad  &&   \text{etc.} \\ 
\end{alignat*}
Because $a=4$, $\Delta a=11$, $\Delta^2 a=14$, $\Delta^3 a=6$, the sum will be

\begin{equation*}
S=\frac{4x}{1-x}+\frac{11xx}{(1-x)^2}+\frac{14x^3}{(1-x)^3}+\frac{6x^4}{(1-x)^4}
\end{equation*}
or

\begin{equation*}
S=\frac{4x-xx+4x^3-x^4}{(1-x)^4}=\frac{x(1+xx)(4-x)}{(1-x)^4}.
\end{equation*}

\paragraph*{§5}

Although this way the sums of these infinite series are found, nevertheless using the same principles the finite counterparts of these series, i.e. series consisting of a finite number of terms, can be summed. For, let this series be propounded

\begin{equation*}
S=ax+bx^2+cx^3+dx^4+\cdots\cdots+ox^n,
\end{equation*}
and first let its sum be in question, if the series actually an infinite series; then the sum will be

\begin{equation*}
=\frac{x}{1-x}a+\frac{x^2}{(1-x)^2}\Delta a+\frac{x^3}{(1-x)^3}\Delta^2 a+\text{etc.}
\end{equation*}
Now consider the terms of the same series following after the last $ox^n$ which we want to put

\begin{equation*}
px^{n+1}+qx^{n+2}+rx^{n+3}+sx^{n+4}+\text{etc.};
\end{equation*}
the sum of this series, if  divided by $x^n$, can be found as before; this sum, multiplied by $x^n$ again, will be

\begin{equation*}
\frac{x^{n+1}}{1-x}p+\frac{x^{n+2}}{(1-x)^2}\Delta p+\frac{x^{n+3}}{(1-x)^3}\Delta^2 p+\text{etc.};
\end{equation*}
if the sum of this series is subtracted from the sum of the infinite series, the sum of the propounded portion in question will remain, i.e.

\begin{equation*}
S=\frac{x}{1-x}(a-x^p)+\frac{x^2}{(1-x)^2}(\Delta a-x^n \Delta p)+\frac{x^3}{(1-x)^3}(\Delta^2 a- x^n \Delta^2 p)+\text{etc.}
\end{equation*}

I. \textit{ Let the sum of this finite series be in question}

\begin{equation*}
S=1x+2x^2+3x^3+4x^4+\cdots + nx^n.
\end{equation*}
Find the differences so of these coefficients as of the ones following the last term\\[2mm]

\begin{tabular}{c|c}
$ \begin{aligned}
&1, ~~ & ~~ &2, ~~ & ~~ & 3, ~~ & ~~ &4, ~~ & ~~ & \text{etc.} \\ 
& ~~ & 1,~~  & ~~ &1, ~~ & ~~ &1,  ~~ &  ~~  & \text{etc.} \\ 
\end{aligned} \qquad $  &
$ \qquad \begin{aligned}
&n+1, ~~ & ~~ &n+2, ~~ & ~~ & n+3, ~~ & ~~ & \text{etc.} \\ 
& ~~ & 1,~~  & ~~ &1, ~~ & ~~ & \text{etc.} \\ 
\end{aligned} $
\end{tabular}\\[3mm]
and it will be $a=1$, $\Delta a=1$, $p=n+1$, $\Delta p=1$, whence the sum in question is

\begin{equation*}
s=\frac{x}{1-x}(1-(n+1)x^n)+\frac{x^2}{(1-x)^2}(1-x^n)
\end{equation*}
or

\begin{equation*}
S=\frac{x-(n+1)x^{n+1}+nx^{n+2}}{(1-x)^2}.
\end{equation*}

II.\textit{ Let the sum of this finite series be in question}

\begin{equation*}
S=1+x+4x²+9x^3+16x^4+\cdots \cdots+n^2x^n.
\end{equation*}
At first, investigate the differences this way\\[2mm]

\begin{center}
\begin{small}
\begin{tabular}{c|c}
$ \begin{aligned}
&1, ~~ & ~~ &4, ~~ & ~~ & 9, ~~ & ~~ &16, ~~ & ~~ & \text{etc.} \\ 
& ~~ & 3,~~  & ~~ &5, ~~ & ~~ &7,  ~~ &  ~~  & \text{etc.} \\ 
& ~~ & ~~  &2, ~~ & ~~ &2, ~~ &  ~~ & \text{etc.}  \\ 
\end{aligned} \qquad $  &
$ \begin{aligned}
&(n+1)^2, &  &(n+2)^2,  & & (n+3)^2, &  & \text{etc.} \\ 
&  & 2n+3,  &  &2n+5,  &    & \text{etc.} \\ 
&  &    &\qquad 2,  &  & \qquad\text{etc.}  \\ 
\end{aligned} \qquad $
\end{tabular}
\end{small}
\end{center}
having found which, the sum in question will be

\begin{equation*}
S=\frac{x}{1-x}(1-(n+1)^2x^n)+\frac{x^2}{(1-x)^2}(3-(2n+3)x^n)+\frac{x^3}{(1-x)^3}(2-2x^n)
\end{equation*}
or 

\begin{equation*}
S=\frac{x+xx-(n+1)^2 x^{n+1}+(2nn+2n-1)x^{n+2}-nnx^{n+3}}{(1-x)^3}.
\end{equation*}

\paragraph*{§6}

But if the propounded series does not have coefficients which are finally reduced to constant differences,  the transformation exhibited here is not of any use to determine its sum. Furthermore, the sum  can not even be approximated in a more convenient way applying said transformation than it is possible by actual addition of the terms of the propounded series itself. For, if in the series $ax+bx^2+cx^3+dx^4+\text{etc.}$ it was $x<1$, in which case only the summation, in the sense explained above, is actually possible, it will be $\frac{x}{1-x}>x$ and hence the new series converges less than the initial one. But if in the propounded series it was $x=1$,  all terms of the new series even become infinite, in which case this transformation will therefore be completely useless.

\paragraph*{§7}

Let us consider the series in which the signs $+$ and $-$ alternate and which will be deduced from the preceding by assuming $x$  to be negative. Therefore, if it was

\begin{equation*}
S=ax-bx^2+cx^3-dx^4+ex^5-\text{etc.},
\end{equation*}
the negative of which series results, if in the preceding series one takes a negative $x$, as before, let us take the differences $\Delta a$, $\Delta^2 a$, $\Delta^3 a$ etc. of the series of coefficients $a$, $b$, $c$, $d$, $e$ etc., having attributed the signs to the powers of $x$, and the propounded series will be transformed into this one

\begin{equation*}
S=\frac{x}{1+x}a-\frac{x^2}{(1+x)^2}\Delta a+\frac{x^3}{(1+x)^3}\Delta^2 a-\frac{x^4}{(1+x)^4}\Delta^3 a+\text{etc.},
\end{equation*}
whence it  is seen that the propounded series can be summed in the same cases as the preceding one, of course, if the series $a$, $b$, $c$, $d$ etc. has finally constant iterated differences.

\paragraph*{§8}

But in this case, this transformation yields a convenient approximation of the value of the propounded series $ax-bx^2+cx^3-dx^4+ex^5-fx^6+\text{etc.}$; for, no matter how large the number $x$ is, the fraction $\frac{x}{1+x}$, in powers of which the other series is expanded, becomes smaller than $1$; and if $x=1$, it will be $\frac{x}{1+x}=\frac{1}{2}$. But if  $x<1$, say $x=\frac{1}{n}$, it will be $\frac{x}{1+x}=\frac{1}{n+1}$ and hence the series found by means of the transformation will always converge more rapidly than the initial one. Let us especially consider the case, in which  $x=1$, which is especially useful for the summation of series, and let 

\begin{equation*}
S=a-b+c-d+e-f+\text{etc.},
\end{equation*}
and denote the first, second and following differences of $a$, which the progression $a$, $b$, $c$, $d$, $e$ etc. yields, by $\Delta a$, $\Delta^2 a $, $\Delta^3 a$ etc.; having found these, it will be

\begin{equation*}
S=\frac{1}{2}a-\frac{1}{4}\Delta a+\frac{1}{8}\Delta^2 a-\frac{1}{16}\Delta^3 a+\text{etc.},
\end{equation*}
which series, if it does not actually terminate, exhibits the approximate sum conveniently.

\paragraph*{§9}

Therefore, let us show the use of this last transformation, in which we took $x=1$, in some examples and at first certainly in examples in which  the true sum can be expressed finitely. Such series are divergent series, in which the numbers $a$, $b$, $c$, $d$ etc. finally lead to constant differences; since the sums of these  series  can not be exhibited in the usual sense of the word sum, we understand the word sum here in this sense we gave  it above [§ 111 of the first part], such that the word sum means the value of the finite expression, from whose the expansion  the propounded series results.\\[2mm]

I. \textit{Therefore, let this series due to Leibniz be propounded}

\begin{equation*}
S=1-1+1-1+1-1+\text{etc.};
\end{equation*}
because in this series all terms are equal, all differences will become $=0$ and hence, because of $a=1$, it will be $S=\frac{1}{2}$.\\[2mm]

II. \textit{Let this series be propounded, i.e.}

\begin{equation*}
S=1-2+3-4+5-6+\text{etc.}
\end{equation*}
\begin{equation*}
\text{Diff. I} \quad 1, \quad 1, \quad 1, \quad 1, \quad 1 \quad \text{etc.}
\end{equation*}
Therefore, because $a=1$, $\Delta a=1$, it will be $S=\frac{1}{2}-\frac{1}{4}=\frac{1}{4}$.\\[2mm]

III. \textit{Let this series be propounded}

\begin{equation*}
S=1-3+5-7+9-\text{etc.}
\end{equation*}
\begin{equation*}
\text{Diff. II} \quad 2, \quad 2, \quad 2, \quad 2, \quad \text{etc.}
\end{equation*}
Because of $a=1$ and $\Delta a=2$, we have $S=\frac{1}{2}-\frac{2}{4}=0$.\\[2mm]

IV. \textit{Let this series of the triangular numbers be propounded, i.e.}

\begin{equation*}
1-3+6-10+15-21+\text{etc.}
\end{equation*}
\begin{alignat*}{19}
&\text{Diff. I} \quad &&2, \quad && \quad && 3, \quad && \quad && 4, \quad && \quad && 5, \quad && \quad && 6, \quad && \quad && \text{etc.} \\
&\text{Diff. II} \quad && \quad &&1 \quad &&  \quad &&1, \quad && \quad &&1, \quad  && \quad &&1 \quad && \quad && \text{etc.} \\
\end{alignat*}
Here, because of $a=1$, $\Delta a=2$ and $\Delta \Delta a=1$, it will be $S=\frac{1}{2}-\frac{2}{4}+\frac{1}{8}=\frac{1}{8}$.\\[2mm]

V. \textit{Let the series of the squares be propounded, i.e.}

\begin{equation*}
S=1-4+9-16+25-36+\text{etc.}
\end{equation*}
\begin{alignat*}{19}
&\text{Diff. I} \quad &&3, \quad && \quad && 5, \quad && \quad && 7, \quad && \quad && 9, \quad && \quad && 11, \quad && \quad && \text{etc.} \\
&\text{Diff. II} \quad && \quad &&2 \quad &&  \quad &&2, \quad && \quad &&2, \quad  && \quad &&2 \quad && \quad && \text{etc.} \\
\end{alignat*}
Because of $a=1$, $\Delta a=3$, $\Delta \Delta a =2$ it will be $S=\frac{1}{2}-\frac{3}{4}+\frac{2}{8}=0$.\\[2mm]

VI. \textit{Let this series of the fourth powers be propounded, i.e.}

\begin{equation*}
S=1-16+81-256+625-1296+\text{etc.}
\end{equation*}
\begin{alignat*}{19}
&\text{Diff. I} \quad &&15, \quad && \quad && 65, \quad && \quad && 175, \quad && \quad && 175, \quad && \quad && 369, \quad && \quad && 671, \quad && \quad && \text{etc.} \\
&\text{Diff. II} \quad && \quad &&50, \quad && \quad &&110, \quad &&  \quad && 194, \quad && \quad && 302, \quad && \quad && \text{etc.}  \\
&\text{Diff. III} \quad && \quad && \quad && 60, \quad && \quad &&   84, \quad && \quad && 108, \quad && \quad && \text{etc.}   \\
&\text{Diff. IV} \quad && \quad && \quad && \quad &&24, \quad &&  \quad &&24, \quad && \quad && \text{etc.} 
\end{alignat*}
Therefore, it will be $S=\frac{1}{2}-\frac{15}{4}+\frac{50}{8}-\frac{60}{16}+\frac{24}{32}=0$.

\paragraph*{§10}

If the series diverges more rapidly, as the geometric series and other similar ones, applying the transformation it is immediately transformed into a  series converging more, which, if it did not already converge quickly enough, in like manner will be converted into another  converging more rapidly.\\[2mm]

I. \textit{Let this geometric series be propounded}

\begin{equation*}
s=1-2+4-8+16-32+\text{etc.}
\end{equation*}
\begin{alignat*}{19}
&\text{Diff. I} \quad &&1, \quad && \quad && 2, \quad && \quad && 4, \quad  && \quad && 8, \quad && \quad && 16, \quad && \quad &&  \quad && \text{etc.} \\
&\text{Diff. II} \quad && \quad &&1, \quad && \quad &&2, \quad &&  \quad && 4, \quad && \quad && 8, \quad && \quad && \text{etc.}  \\
&\text{Diff. III} \quad && \quad && \quad && 1, \quad && \quad &&   2, \quad && \quad && 4, \quad && \quad && \text{etc.}   
\end{alignat*}
Therefore, because in all these differences the first term is $=1$, the sum of the series will be expressed this way

\begin{equation*}
S=\frac{1}{2}-\frac{1}{4}+\frac{1}{8}-\frac{1}{16}+\frac{1}{32}-\frac{1}{64}+\text{etc.},
\end{equation*}
the sum of which series is $=\frac{1}{3}$; for, it results from the expansion of the fraction $\frac{1}{2+1}$, whereas the propounded series results from $\frac{1}{1+2}$.\\[2mm]

II. \textit{Let this recurring series be propounded}
\begin{equation*}
S=1-2+5-12+29-70+169-\text{etc.}
\end{equation*}
\begin{alignat*}{19}
&\text{Diff. I} \quad && 1, \quad && \quad &&3, \quad && \quad &&7, \quad && \quad && 17, \quad && \quad && 41, \quad && \quad && 99, \quad && \quad && \text{etc.}\\
&\text{Diff. II} \quad &&  \quad && 2, \quad && \quad && 4, \quad && \quad && 10, \quad &&  \quad && 24, \quad &&  \quad && 58, \quad && \quad && \text{etc.} \\
&\text{Diff. III} \quad &&  \quad && \quad && 2, \quad && \quad && 6,\quad &&  \quad && 14,  \quad &&  \quad && 34,  \quad &&  \quad && \text{etc.}  \\
&\text{Diff. IV} \quad &&  \quad && \quad &&  \quad && 4, \quad && \quad && 8,  \quad &&  \quad && 20,  \quad &&  \quad && \text{etc.}  \\
&\text{Diff. V} \quad &&  \quad && \quad &&  \quad  && \quad && 4, \quad && \quad && 12,  \quad &&  \quad && \text{etc.}  \\
&\text{Diff. VI} \quad &&  \quad && \quad &&  \quad  && \quad && \quad && 8, \quad && \quad && \text{etc.}  \\
& &&  \quad && \quad &&  \quad  && \quad && \quad && \quad && \text{etc.}
\end{alignat*}
Therefore, the first terms of the continued differences constitute this double geometric series $1$, $1$, $2$, $2$, $4$, $4$, $8$, $8$, $16$, $16$ etc., whence it will be

\begin{equation*}
S=\frac{1}{2}-\frac{1}{4}+\frac{2}{8}-\frac{2}{16}+\frac{4}{32}-\frac{4}{64}+\frac{8}{128}-\text{etc.};
\end{equation*}
therefore, because, except for the first, each two terms  cancel each other, it will be $S=\frac{1}{2}$. Indeed, the propounded series results from the expansion of the fraction $\frac{1}{1+2-1}=\frac{1}{2}$, as we showed in the discussion of the nature of recurring series.\\[2mm]

III. \textit{Let the hypergeometric series be propounded, i.e.}

\begin{equation*}
S=1-2+6-24+120-720+5040-\text{etc.},
\end{equation*}
whose continued differences we will investigate more conveniently this way:
\begin{small}
\begin{center}
\begin{tabular}{r|r|r|r|r}
 & Diff. I & Diff. II & Diff. III & \\ 
 $1$ & $1$ & $3$ & $11$ & \\
 $2$ & $4$ & $14$ & $64$ & \\
 $6$ & $18$ & $78$ & $426$ & \\
 $24$ & $96$ & $504$ & $3216$ & \\
 $120$ & $600$ & $3720$ & $27240$ & etc.\\
 $720$ & $4320$ & $30960$ & $256320$ & \\
 $5040$ & $35280$ & $287280$ & $2656080$ & \\
 $40320$ & $322560$ & $2943360$ &  & \\
 $362880$ & $3265920$ &  &  & \\
 $3628800$ &  &  &  & \\
\end{tabular}
\end{center}
\end{small}
Having  continued these differences further, it will be

\begin{equation*}
S=\frac{1}{2}-\frac{1}{4}+\frac{3}{8}-\frac{11}{16}+\frac{53}{32}-\frac{309}{64}+\frac{2119}{128}-\frac{16687}{256}+\frac{148329}{512}-\frac{1468457}{1024}
\end{equation*}
\begin{equation*}
+\frac{16019531}{2048}-\frac{190899411}{4096}+\text{etc.}
\end{equation*}
Collect the two initial terms and it will be $S=\frac{1}{4}+A$, where

\begin{equation*}
A=\frac{3}{8}-\frac{11}{16}+\frac{53}{32}-\frac{309}{64}+\frac{2119}{128}-\text{etc.}
\end{equation*}
If now in the same way the differences are taken, it will be

\begin{equation*}
A=\frac{3}{2^4}-\frac{5}{2^6}+\frac{21}{2^8}-\frac{99}{2^{10}}+\frac{615}{2^{12}}-\frac{4401}{2^{14}}+\frac{36585}{2^{16}}-\frac{342207}{2^{18}}
\end{equation*}
\begin{equation*}
+\frac{3565323}{2^{20}}-\frac{40866525}{2^{22}}+\text{etc.}
\end{equation*}
Collect the two initial terms, because they converge, and it will be $A=\frac{7}{2^6}+B$ while $B=\frac{21}{2^8}-\frac{99}{2^{10}}+\text{etc.}$; taking the differences of this series again, it will be

\begin{equation*}
B=\frac{21}{2^9}-\frac{15}{2^{12}}+\frac{159}{2^{15}}-\frac{429}{2^{18}}+\frac{5241}{2^{21}}-\frac{26283}{2^{24}}+\frac{338835}{2^{27}}-\frac{2771097}{2^{30}}+\text{etc.}
\end{equation*}
Collect the four initial terms into one and put $B=\frac{153}{2^{12}}+\frac{843}{2^{20}}+C$ while

\begin{equation*}
C=\frac{5241}{2^{21}}-\frac{26283}{2^{24}}+\text{etc.}
\end{equation*}
and actually summing  some terms it will approximately be $C=\frac{15645}{2^{24}}-\frac{60417}{2^{30}}$. Therefore, from these   the sum of the series will finally be concluded to be $S=0.40082055$, which can nevertheless only be considered to be accurate hardly further than to three or four digits because of the divergence of the series; it is nevertheless certainly smaller than the correct value. For, in another paper I found this sum to be $=0.4036524077$, where not even the last digit deviates from the true value.

\paragraph*{§11}

But  this transformation is especially  useful for the transformation of an already but slowly converging series into others which converge a lot more rapidly. Since indeed the following terms are smaller than the preceding ones, the first differences become negative; therefore, in the following the nature of the sign is carefully to be taken into account.\\[2mm]

I. \textit{Let this series be propounded}

\begin{equation*}
S=1-\frac{1}{2}+\frac{1}{3}-\frac{1}{4}+\frac{1}{5}-\frac{1}{6}+\text{etc.}
\end{equation*}

\begin{alignat*}{9}
&\text{Diff. I} \quad &&-\frac{1}{2}, \quad  -\frac{1}{2 \cdot 3}, \quad  -\frac{1}{3 \cdot 4}, \quad  -\frac{1}{4 \cdot 5}, \quad  \frac{1}{5 \cdot 6 } \quad \text{etc.}\\[3mm]
&\text{Diff. II} \quad &&+\frac{1}{3}, \quad  \frac{2}{2 \cdot 3 \cdot 4}, \quad  \frac{2}{3 \cdot 4 \cdot 5}, \quad  \frac{2}{4 \cdot 5 \cdot 6}  \quad \text{etc.}\\[3mm]
&\text{Diff. III} \quad &&-\frac{1}{4}, \quad  -\frac{2 \cdot 3}{2 \cdot 3 \cdot 4 \cdot 5}, \quad  -\frac{2 \cdot 3}{3 \cdot 4 \cdot 5 \cdot 6} \quad \text{etc.} \\[3mm]
&\text{Diff. IV} \quad &&+\frac{1}{5} \quad  \text{etc.}\\
&  && \text{etc.}
\end{alignat*} 
Therefore, it will be

\begin{equation*}
S=\frac{1}{2}+\frac{1}{2 \cdot 4}+\frac{1}{3 \cdot 8}+\frac{1}{4 \cdot 16}+\frac{1}{5 \cdot 32}+\text{etc.};
\end{equation*}
but we already showed in the \textit{Introductio} that both series exhibit the hyperbolic logarithm of two.\\[2mm]

II. \textit{Let this series for the circle be propounded, i.e.}

\begin{equation*}
S=1-\dfrac{1}{3}+\dfrac{1}{5}-\dfrac{1}{7}+\dfrac{1}{9}-\dfrac{1}{11}+\text{etc.}
\end{equation*}
\begin{alignat*}{9}
&\text{Diff. I} \quad && -\frac{2}{1 \cdot 3}, \quad -\frac{2}{3 \cdot 5}, \quad -\frac{2}{5\cdot 7}, \quad  -\frac{2}{7 \cdot 9}, \quad  -\frac{2}{9 \cdot 11} \quad \text{etc.} \\[3mm]
&\text{Diff. II} \quad && +\frac{2 \cdot 4}{1 \cdot 3 \cdot 5}, \quad \frac{2 \cdot 4}{3 \cdot 5 \cdot 7}, \quad \frac{2\cdot 4}{5\cdot 7 \cdot 9}, \quad  \frac{2 \cdot 4}{7 \cdot 9 \cdot 11} \quad    \text{etc.} \\[3mm]
&\text{Diff. III} \quad && -\frac{2 \cdot 4 \cdot 6}{1 \cdot 3 \cdot 5 \cdot 7}, \quad \frac{2 \cdot 4 \cdot 6}{3 \cdot 5 \cdot 7 \cdot 9} \quad \text{etc.} \\
\end{alignat*}
Therefore,  the sum of the series will also be

\begin{equation*}
S=\frac{1}{2}+\frac{1}{3 \cdot 2}+\frac{1 \cdot 2}{3 \cdot 5 \cdot 2}+\frac{1 \cdot 2 \cdot 3}{3 \cdot 5 \cdot 7 \cdot 2}+\text{etc.}
\end{equation*}
or

\begin{equation*}
2S= 1+\frac{1}{3}+\frac{1 \cdot 2}{3\cdot 5}+\frac{1 \cdot 2 \cdot 3}{3 \cdot 5 \cdot 7}+\frac{1 \cdot 2 \cdot 3 \cdot 4}{3 \cdot 5 \cdot 7 \cdot 9}+\text{etc.}
\end{equation*}

III. \textit{Let the value of this infinite series be in question}

\begin{equation*}
S= \log 2- \log 3 +\log 4 -\log 5 + \log 6 - \log 7+\log 8- \log 9+\text{etc.}
\end{equation*}
Since the differences  become too irregular at the beginning, let us actually sum the terms up to $\log 10$ using tables, the value of which sum will be found to be $=-0.3911005$, and it will be

\begin{equation*}
S=- 0.3911005+\log 10- \log 11+ \log 12- \log 13 +\log 14- \log 15+\text{etc.}
\end{equation*}
Take those logarithms from tables and look for their differences this way:

\begin{center}
\begin{tabular}{c|c|c|c|c|c}
 & Diff. I & Diff. II& Diff. III & Diff. IV & Diff. V \\
$\log 10 =1.0000000$ & $+$       & $-$     & $+$ & $-$ & $+$ \\
                    &            &        &     &     &    \\  
$\log 11 =1.0413927$ & $413927$  &         &     &  \\
                    &           & $36042$ & &\\
$\log 12 =1.0791812$ & $377885$  &         & $5779$ & &\\
                    &           & $30263$ & & $1292$ & \\
$\log 13 = 1.1139434$& $347622$  & & $4487$ & & $368$ \\
                    &           & $25776$ & & $924$ & \\
$\log 14 = 1.1461280$& $321846$  & & $3563$& &   \\
                    & & $22213$ & &  &  \\
$\log 15= 1.1760913$ & $299633$ & & & &                                                     
\end{tabular}
\end{center}
From these one finds

\begin{equation*}
\log 10 -\log 11 +\log 12 - \log 13 +\text{etc.}
\end{equation*}
\begin{equation*}
=\frac{1.0000000}{2}-\frac{413927}{4}-\frac{36042}{8}-\frac{5779}{16}-\frac{1292}{32}-\frac{368}{64}=0.4891606.
\end{equation*}
Therefore, the value of the propounded series will be
 
\begin{equation*}
 S= \log 2 - \log 3 +\log 4 - \log 5+\text{etc.}=0.0980601,
\end{equation*}
to which logarithm the number $1.253315$ corresponds.
 
\paragraph*{§12} 

As we obtained these transformations by substituting the fraction $\frac{y}{1\pm y}$ for $x$ in the series, so innumerable other transformations will result, if  other functions of $y$ are substituted for $x$. Let again this series be propounded

\begin{equation*}
S=ax+bx^2+cx^3+dx^4+ex^5+fx^6+\text{etc.}
\end{equation*}
and put $x=y(1-y)$, having done which the following series will result

\begin{equation*}
\renewcommand{\arraystretch}{1,3}
\setlength{\arraycolsep}{0.0mm}
\begin{array}{rrrrrrrrrrrrrrrrrrrrrrrrr}
&S &~=~& ay&~-~&ayy \\
&  &   &   &~+~&byy&~-~&2by^3&~+~& by^4  \\
&  &   &   &   &   &~+~&cy^3 &~-~& 3cy^4&~+~& 3cy^5 &~-~& cy^6  \\
&  &   &   &   &   &   &     &~+~& dy^4 &~-~& 4dy^5 &~+~& 6dy^6  \\
&  &   &   &   &   &   &     &   &      &~+~& ey^5  &~-~& 5ey^6  \\
&  &   &   &   &   &   &     &   &      &   &       &~+~& fy^6 & \quad  \text{etc.}  \\
\end{array}
\end{equation*}
Therefore, if the one of these series was summable, at the same time the sum of the other will be known. Hence, if one puts

\begin{equation*}
S=x+x^2+x^3+x^4+x^5+\text{etc.}=\frac{x}{1-x},
\end{equation*}
it will be

\begin{equation*}
S=y-y^3-y^4+y^6+y^7-y^9-y^{10}+\text{etc.}
\end{equation*}
The sum of this series will therefore  also be $=\frac{y-yy}{1-y+yy}$.

\paragraph*{§13}

If the one series terminates anywhere, the sum of the other can be exhibited explicitly. Let us put $a=1$ and that in the found series all terms after the first  vanish, so that $S=y$;  hence, because of $x=y-yy$, the sum of the first will be $=\frac{1}{2}-\sqrt{\frac{1}{4}-x}$. But, because of $a=1$, it will be

\begin{equation*}
\renewcommand{\arraystretch}{2,5}
\setlength{\arraycolsep}{0.0mm}
\begin{array}{rrrrrl}
&b &~=~ & 1&~=~& \dfrac{1}{4}\cdot 2^2 \\
&c &~=~ & 2&~=~& \dfrac{1 \cdot 3}{4 \cdot 6} \cdot 2^4 \\
&d &~=~ &5 &~=~& \dfrac{1 \cdot 3 \cdot 5}{4 \cdot 6 \cdot 8}\cdot 2^6 \\
&e &~=~ &14&~=~& \dfrac{1 \cdot 3 \cdot 5 \cdot 7}{4 \cdot 6 \cdot 8 \cdot 10} \cdot 2^8 \\
&f&~=~&42 &~=~&\dfrac{1 \cdot 3 \cdot 5 \cdot 7 \cdot 9}{4 \cdot 6 \cdot 8 \cdot 10 \cdot 12}\cdot 2^{10} \\
&g &~=~&132 &~=~& \dfrac{1 \cdot 3 \cdot 5 \cdot 7 \cdot 9 \cdot 11}{4 \cdot 6 \cdot 8 \cdot 10 \cdot 12 \cdot 14} \cdot 2^{12} \\
& && && \quad \text{etc.},
\end{array}
\end{equation*}
whence the first series will go over into this one

\begin{equation*}
S=\frac{1}{2}-\sqrt{\frac{1}{4}-x}=x+x^2+2x^3+5x^4+14x^5+42x^6+132x^7+\text{etc.},
\end{equation*}
which same series is found, if the surdic quantity $\sqrt{\frac{1}{4}-x}$ is expanded into a series and is subtracted from $\frac{1}{2}$.

\paragraph*{§14}

For  the transformation to extend further, let us put $x=y(1+ny)^{\nu}$ and the propounded series

\begin{equation*}
S=ax+bx^2+cx^3+dx^4+ex^5+\text{etc.}
\end{equation*}
will be transformed into the following

\begin{small}
\begin{equation*}
\renewcommand{\arraystretch}{2,5}
\setlength{\arraycolsep}{0.0mm}
\begin{array}{rrrrrrrrrrrrrrrrrrrrrrrrrrr}
&S&~=~ay &~+~&\dfrac{\nu}{1}nay^2 &~+~&\dfrac{\nu (\nu -1)}{1 \cdot 2}n^2 ay^3 &~+~&\dfrac{\nu (\nu -1)(\nu -2)}{1 \cdot 2 \cdot 3}n^3 ay^4&~+~&\dfrac{\nu (\nu-1)(\nu-2)(\nu-3) }{1 \cdot 2 \cdot 3 \cdot 4}n^4 a y^5\\
& & &  ~+~&by^2 &~+~&\dfrac{2 \nu}{1}n b y^3 &~+~&\dfrac{2 \nu (2 \nu -1)}{1 \cdot 2}n^2 b y^4 &~+~& \dfrac{2 \nu (2 \nu -1)(2 \nu -2)}{1 \cdot 2 \cdot 3}n^3 by^5 \\
& & & & &~+~&cy^3 &~+~&\dfrac{3 \nu}{1}n cy^4 &~+~&\dfrac{3 \nu (3 \nu -1)}{1 \cdot 2}n^2 cy^5 \\
&  & & & & & & ~+~&dy^4 &~+~&\dfrac{4 \nu}{1}n dy^5 \\
&  & & & & & &\text{etc.}        
\end{array}
\end{equation*}
\end{small}Therefore, if the sum of this series was known, one will at the same time also have the sum of the first one and vice versa. Since $n$ and $\nu$ can be taken arbitrarily, from one summable series innumerable other summable ones can be found.

\paragraph*{§15}

One can also do transformations of such a kind that the sum of the found series becomes irrational in the following way.\\[2mm]
Let this series be propounded

\begin{equation*}
S=ax+bx^3+cx^5+dx^7+ex^9+fx^{11}+\text{etc.};
\end{equation*}
it will be

\begin{equation*}
Sx= ax^2+bx^4+cx^6+dx^8+ex^{10}+fx^{12}+\text{etc.}
\end{equation*}
Now, put

\begin{equation*}
x=\frac{y}{\sqrt{1-nyy}};
\end{equation*}
it will be $xx=\frac{y^2}{1-nyy}$ and the propounded series will be transformed into this one

\begin{equation*}
\renewcommand{\arraystretch}{2,0}
\setlength{\arraycolsep}{0.0mm}
\begin{array}{rrrrrrrrrrrrrrrrrrrrr}
&\dfrac{Sy}{\sqrt{1-nyy}}&~=~& ay^2&~+~&nay^4&~+~&n^2 ay^6&~+~&n^3 ay^8&~+~&n^4 a y^{10}&~+~&\text{etc.}\\
& & & &~+~&by^4&~+~&2n by^6&~+~&3n^2 by^8&~+~&4n^3b  y^{10}&~+~&\text{etc.}\\
&  & & &   &    &~+~&cy^6&~+~&3n cy^8&~+~&6n^2c y^{10}&~+~&\text{etc.}\\
&  & & &   &    &   &    &~+~&dy^8&~+~&4nd  y^{10}&~+~&\text{etc.}\\
& & & &   &    &   &    &   &    &~+~&e  y^{10}&~+~&\text{etc.}\\
& & & & & & &\text{etc.}
\end{array}
\end{equation*}
Therefore, if the sum $S$ was known from the first series, one will at the same time have the sum of the following series

\begin{equation*}
\frac{S}{\sqrt{1-nyy}}=ay+(na+b)y^3+(n^2a+2nb+c)y^5+(n^3a +3n^2b+3nc+d)y^7+\text{etc.}
\end{equation*}

\paragraph*{§16}

If one takes $n=-1$, the coefficients of this series will be the iterated differences of $a$ from the series $a$, $b$, $c$, $d$ etc.; but if in the propounded series the signs  alternate, the coefficients will be these differences for $n=1$. Therefore, let $\Delta a$, $\Delta^2 a$, $\Delta^3 a$, $\Delta^4 a$ etc.  denote the first, second, third etc. differences of $a$ of the series of the numbers $a$, $b$, $c$, $d$, $e$, $f$ etc. And if it was

\begin{equation*}
S=ax+bx^3+cx^5+dx^7+ex^9+\text{etc.},
\end{equation*}
having put $x=\frac{y}{\sqrt{1+yy}}$, it will be

\begin{equation*}
\frac{S}{\sqrt{1+yy}}=ay+\Delta a \cdot y^3+\Delta^2 a \cdot y^5+\Delta^3 a \cdot y^7+\text{etc.}
\end{equation*}
But if it was

\begin{equation*}
S=ax-bx^3+cx^5-dx^7+ex^9-\text{etc.}
\end{equation*}
and one puts $x=\frac{y}{\sqrt{1-yy}}$, it will be

\begin{equation*}
\frac{S}{\sqrt{1-yy}}=ay-\Delta a \cdot y^3+\Delta^2 a \cdot y^5- \Delta^3 a \cdot y^7+\text{etc.}
\end{equation*}
Therefore, if the series $a$, $b$, $c$, $d$, $e$ etc. finally leads to constant differences, both series can be summed explicitly; but this summation also follows from the preceding paragraphs.

\paragraph*{§17}

Let us put that the coefficients $a$, $b$, $c$, $d$ etc. constitute this series

\begin{equation*}
1, \quad \frac{1}{3}, \quad \frac{1}{5}, \quad \frac{1}{7}, \quad \frac{1}{9} \quad \text{etc.}
\end{equation*}
and, as we already saw above [§ 11, II.], it will be

\begin{equation*}
a=1, \quad \Delta a=-\frac{2}{3}, \quad \Delta^2 a=\frac{2 \cdot 4}{3 \cdot 5}, \quad \Delta^3 a=-\frac{2 \cdot 4 \cdot 6}{3 \cdot 5 \cdot 7} \quad \text{etc.},
\end{equation*}
whence we will sum the following two series.\\[2mm]

I. Let $S=x+\frac{1}{3}x^3+\frac{1}{5}x^5+\frac{1}{7}x^7+\text{etc.}$; it will be $S=\frac{1}{2}\log \frac{1+x}{1-x}$. Now, having put $x=\frac{y}{\sqrt{1+yy}}$, it will be

\begin{equation*}
S=\frac{1}{2} \log \frac{\sqrt{1+yy}+y}{\sqrt{1+yy}-y}=\log (\sqrt{1+yy}+y),
\end{equation*}
whence it will be

\begin{equation*}
\frac{\log (\sqrt{1+yy}+y)}{\sqrt{1+yy}}=y-\frac{2}{3}y^3+\frac{2 \cdot 4}{3 \cdot 5}y^5-\frac{2 \cdot 4 \cdot 6}{3 \cdot 5 \cdot 7}y^7+\text{etc.}
\end{equation*}

II. Let $S=x-\frac{1}{3}x^3+\frac{1}{5}x^5-\frac{1}{7}x^7+\text{etc.}$; it will be $S=\arctan x$. Now, having put $x=\frac{y}{\sqrt{1-yy}}$, it will be

\begin{equation*}
S= \arctan \frac{y}{\sqrt{1-yy}}= \arcsin y =\arccos \sqrt{1-yy}.
\end{equation*}
Therefore, one will obtain this summation

\begin{equation*}
\frac{\arcsin y}{\sqrt{1-yy}}=y+\frac{2}{3}y^3+\frac{2 \cdot 4}{3 \cdot 5}y^5+\frac{2 \cdot 4 \cdot 6}{3 \cdot 5 \cdot 7}y^7+\text{etc.}
\end{equation*}

\paragraph*{§18}

One can also substitute transcendental functions of $y$ for $x$ and  can discover other  summations more difficult to find this way; but nevertheless, for the new series to not become too complex, one has to pick functions whose powers can easily be exhibited, as it is the case for the exponential quantities $e^y$. Therefore, having propounded this series

\begin{equation*}
S=ax+bx^2+cx^3+dx^4+ex^5+fx^6+\text{etc.},
\end{equation*}
put $x=e^{ny}y$, where $e$ denotes the number whose hyperbolic logarithm is $=1$; it will be $x^2=e^{2ny}y^2$, $x^3=e^{3ny}y^3$ etc. Indeed, in general, as it is known,

\begin{equation*}
e^z=1+z+\frac{z^2}{1 \cdot 2}+\frac{z^3}{1 \cdot 2 \cdot 3}+\frac{z^4}{1 \cdot 2 \cdot3 \cdot 4}+\text{etc.}
\end{equation*}
Therefore, the propounded series will be transformed into this one

\begin{equation*}
\renewcommand{\arraystretch}{2,5}
\setlength{\arraycolsep}{0.0mm}
\begin{array}{rrrrrrrrrrrrrrrrrrrrr}
&S&~=~ay&~+~&1nay^2&~+~&\dfrac{1}{2}&n^2ay^3&~+~&\dfrac{1}{6}&n^3 ay^4&~+~&\dfrac{1}{24}&n^4a y^5&~ +~&\text{etc.} \\
& &     &~+~&b y^2&~+~&\dfrac{2}{1}&nby^3&~+~&\dfrac{4}{2}&n^2by^4&~+~&\dfrac{8}{6}&n^3b y^5 &~ +~&\text{etc.} \\
& &     &   &        &~+~&  &cy^3&~+~&\dfrac{3}{1}&ncy^4&~+~&\dfrac{9}{2}&n^2 cy^5 &~+~&\text{etc.} \\
& & & & & & &  &~+~& &dy^4&~+~&\dfrac{4}{1}&n dy^5 &~ +~&\text{etc.} \\
& & & & & & &  &   & &    &~+~& &ey^5 &~+~&\text{etc.} \\
& && && && && \text{etc.}
\end{array}
\end{equation*}

I. Let the geometric series be propounded, i.e. $S=x+x^2+x^3+x^4+x^5+\text{etc.}$; it will be $S=\frac{x}{1-x}$. Now, put $n=-1$ so that $x=e^{-y}y$ and $S=\frac{e^{-y}y}{1-e^{-y}y}=\frac{y}{e^y-y}$; one will find this sum

\begin{equation*}
\frac{y}{e^y-y}=y-\frac{1}{2}y^3-\frac{1}{6}y^4+\frac{5}{24}y^5+\frac{19}{120}y^6-\text{etc.},
\end{equation*}
the law of which series is not recognized\footnote{By this Euler means that it is not possible to see an explicit formula for the general coefficient of this power series.}.\\[2mm]

II. In the other series, let all terms except for the first be $=0$; it will be

\begin{equation*}
b=-na, \quad c=\frac{3}{2}n^2a, \quad d=-\frac{8}{3}n^3 a, \quad e=\frac{125}{24}n^4a, \quad f=-\frac{54}{5}n^5a \quad \text{etc.}
\end{equation*}
Because therefore the sum is $S=ay$ and $x=ye^{ny}$, it will be

\begin{equation*}
y=x-nx^2+\frac{3}{2}n^2 x^3-\frac{8}{3}n^3 x^4+\frac{125}{24}n^4 x^5-\frac{54}{5}n^5 x^6+\text{etc.}
\end{equation*}
Since in these series the structure of the progression is not obvious, the summations deduced from this substitution have hardly any use. But  the transformations derived from the substitution $x=\frac{y}{1 \pm y}$, which not only yield extraordinary summations but also appropriate ways to approximate the sums of series, stand out especially. Therefore, having mentioned these things in advance without resorting to differential calculus, we want to proceed to show the use of this calculus in the doctrine of series.


\title{On the investigation of summable series
\footnote{
Original title: "'De Investigatione serierum summabilium"', first published as part of the book "`\textit{Institutiones calculi differentialis cum eius usu in analysi finitorum ac doctrina serierum}, 1755"', reprinted in \textit{Opera Omnia}: Series 1, Volume 10, pp. 235 - 255, Eneström-Number E212, translated by: Alexander Aycock for the "`Euler-Kreis Mainz"'}}
\author{Leonhard Euler}
\date{ }
\maketitle

\paragraph*{§19}

If the sum of a series, whose terms contain the variable quantity $x$, was known, and which series will therefore be a function of $x$, then, whatever value is attributed to $x$, one will always be able to assign the sum of the series. Therefore, if one writes $x+dx$ instead of $x$, the sum of the resulting series will be equal to the sum of the first series and the differential: Hence it follows that the differential of the sum will be $=$  the differential of the series. Because this way so the sum as each term will be multiplied by $dx$, if one divides by $dx$ everywhere, one will have a new series, whose sum will be known. In like manner, if this series is differentiated again and is divided by $dx$, a new series will result together with its sum and this way  new likewise summable series will be found from one summable series involving the variable quantity $x$, if that series is differentiated several times. 

\paragraph*{§20}

In order to understand all this more clearly, let the variable geometric progression be propounded, whose sum is known; for,

\begin{equation*}
\frac{1}{1-x}=1+x+x^2+x^3+x^4+x^5+x^6+\text{etc.}
\end{equation*}
If this equation is now differentiated with respect to $x$, it will be

\begin{equation*}
\frac{dx}{(1-x)^2}=dx+2xdx+3x^2dx+4x^3dx+5x^4dx+\text{etc.},
\end{equation*}
and having divided by $dx$, one will have

\begin{equation*}
\frac{1}{(1-x)^2}=1+2x+3x^2+4x^3+5x^4+\text{etc.}
\end{equation*}
If one differentiates the last equation once more and divides by $dx$, this equation will result

\begin{equation*}
\frac{2}{(1-x)^3}=2+2\cdot 3 x+ 3 \cdot 4 x^2+4 \cdot 5 x^3+5 \cdot 6 x^4+\text{etc.}
\end{equation*}
or

\begin{equation*}
\frac{1}{(1-x)^3}=1+3x+6x^2+10x^3+15x^4+21x^5+\text{etc.}
\end{equation*}
where the coefficients are the triangular numbers. If one differentiates once again and divides by $3dx$, one will obtain

\begin{equation*}
\frac{1}{(1-x)^4}=1+4x+10x^2+20x^3+35x^4+\text{etc.},
\end{equation*}
whose coefficients are the first pyramidal numbers. And, proceeding further this way, the same series result, which are known to originate from the expansion of the fraction $\frac{1}{(1-x)^n}$.

\paragraph*{§21}

This investigation will extend even further, if, before the differentiation is done, the series and the sum are multiplied by a certain power of $x$ or  even a function of $x$. Therefore, because

\begin{equation*}
\frac{1}{1-x}=1+x+x^2+x^3+x^4+x^5+\text{etc.},
\end{equation*}
multiply by $x^m$ on both sides and it will be

\begin{equation*}
\frac{x^m}{1-x}=x^m+x^{m+1}+x^{m+2}+x^{m+3}+x^{m+4}+\text{etc.}
\end{equation*}
Now differentiate this series and, having divided  the result by $dx$, it will be

\begin{equation*}
\frac{mx^{m-1}-(m-1)x^m}{(1-x)^2}=mx^{m-1}+(m+1)x^m+(m+2)x^{m+1}+(m+3)x^{m+2}+\text{etc.}
\end{equation*}
Now divide by $x^{m-1}$; one will have

\begin{equation*}
\frac{m-(m-1)x}{(1-x)^2}=\frac{m}{1-x}+\frac{x}{(1-x)^2}=m+(m+1)x+(m+2)x^2+\text{etc.}
\end{equation*}
Before differentiating again multiply this equation by $x^n$, so that 
\begin{equation*}
\frac{mx^n}{1-x}+\frac{x^{n+1}}{(1-x)^2}=mx^n+(m+1)x^{n+1}+(m+2)x^{n+2}+\text{etc.}
\end{equation*}
Now, let us differentiate, and having divided by $dx$, it will be

\begin{equation*}
\frac{mnx^{n-1}}{1-x}+\frac{(m+n+1)x^n}{(1-x)^2}+\frac{2x^{n+1}}{(1-x)^3}
\end{equation*}
\begin{equation*}
=mnx^{n-1}+(m+1)(n+1)x^n+(m+2)(n+2)x^{n+1}+\text{etc.}
\end{equation*}
But, having divided by $x^{n-1}$, it will be

\begin{equation*}
\frac{mn}{1-x}+\frac{(m+n+1)x}{(1-x)^2}+\frac{2xx}{(1-x)^3}
\end{equation*}
\begin{equation*}
=mn+(m+1)(n+1)x+(m+2)(n+2)x^2+\text{etc.}
\end{equation*}
and it will be possible to proceed further this way; indeed one will always find the same series which result from the expansions of the fractions constituting the sum.

\paragraph*{§22}

Since the sum of the geometric progression assumed at first can be assigned up to any given term, this way also series consisting of a  finite  number of terms will be summed. Because 

\begin{equation*}
\frac{1-x^{n+1}}{1-x}=1+x+x^2+x^3+x^4+\cdots+x^n,
\end{equation*}
after the differentiation and having divided the result by $dx$ it will be

\begin{equation*}
\frac{1}{(1-x)^2}-\frac{(n+1)x^n-nx^{n+1}}{(1-x)^2}=1+2x+3x^2+4x^3+\cdots +nx^{n-1}.
\end{equation*}
Therefore, the sum of the powers of natural numbers up to a certain term can be found. For, multiply this series by $x$, so that

\begin{equation*}
\frac{x-(n+1)x^{n+1}+nx^{n+2}}{(1-x)^2}=x+2x^2+3x^3+\cdots+nx^n,
\end{equation*}
which, having differentiated it again and divided  by $dx$, will give

\begin{equation*}
\frac{1+x-(n+1)x^n+(2nn+2n-1)x^{n+1}-nnx^{n+2}}{(1-x)^3}=1+4x+9x^2+\cdots+n^2x^{n-1};
\end{equation*}
this equation multiplied by $x$ will give

\begin{equation*}
\frac{x+x^2-(n+1)^2 x^{n+1}+(2nn+2n-1)x^{n+2}-nnx^{n+3}}{(1-x)^3}=x+4x^2+9x^3+\cdots+n^2x^n,
\end{equation*}
which equality, if differentiated, divided by $dx$ and multiplied by $x$, will lead to this series

\begin{equation*}
x+8x^2+27x^3+\cdots+n^2x^n,
\end{equation*}
whose sum can therefore be assigned. And from this in like manner it is possible to find the indefinite sum of the fourth powers and higher powers.

\paragraph*{§23}

Therefore, this method can be applied to all series  containing a variable quantity and whose sum is known, of course. Because, aside from the geometric series, all recurring series enjoy the same property that they can be summed not only up to infinity but also to any given term, one will be able to find innumerable other summable series from these by the same method. Because  a lot of work would be necessary, if we wanted to study this in more detail, let us consider only one single example.\\[2mm]
Let this series be propounded

\begin{equation*}
\frac{x}{1-x-xx}=x+x^2+2x^3+3x^4+5x^5+8x^6+13x^7+\text{etc.},
\end{equation*}
which equation, if  differentiated and divided by $dx$, will give

\begin{equation*}
\frac{1+xx}{(1-x-xx)^2}=1+2x+6x^2+12x^3+25x^4+48x^5+91x^6+\text{etc.}
\end{equation*}
But obviously all series resulting this way will also be recurring, whose sums can therefore even be found more naturally.

\paragraph*{§24}

Therefore, in general, if the sum of a certain series of this form

\begin{equation*}
ax+bx^2+cx^3+dx^4+\text{etc.}
\end{equation*}
was known, which sum we want to put $=S$, one will be able to find the sum of the same series, if each term is multiplied by terms of an arithmetic progression. For, let

\begin{equation*}
S=ax+bx^2+cx^3+dx^4+ex^5+\text{etc.};
\end{equation*}
multiply by $x^m$; it will be

\begin{equation*}
Sx^m= ax^{m+1}+bx^{m+2}+cx^{m+3}+dx^{m+4}+\text{etc.};
\end{equation*}
differentiate this equation and divide by $dx$ to find

\begin{equation*}
mSx^{m-1}+x^m\frac{dS}{dx}=(m+1)ax^m+(m+2)bx^{m+1}+(m+3)cx^{m+2}+\text{etc.};
\end{equation*}
divide by $x^{m-1}$ and it will be

\begin{equation*}
mS+\frac{xdS}{dx}=(m+1)ax+(m+2)bx^2+(m+3)cx^3+\text{etc.}
\end{equation*}
Therefore, if one wants to find the sum of the following series

\begin{equation*}
\alpha a x+(\alpha +\beta)bx^2+(\alpha +2 \beta)cx^3+(\alpha +3 \beta)dx^4+\text{etc.},
\end{equation*}
multiply the above series by $\beta$ and put $m \beta +\beta =\alpha$, so that $M=\frac{\alpha -\beta}{\beta}$, and the sum of this series will be

\begin{equation*}
=(\alpha -\beta)S+\frac{\beta xdS}{dx}.
\end{equation*}

\paragraph*{§25}

One will also be able to find the sum of this propounded series, if  each term is multiplied by a term of a progression of second order, whose second differences are just constant, of course. For, because we already found

\begin{equation*}
mS+\frac{xdS}{dx}=(m+1)ax+(m+2)bx^2+(m+3)cx^3+\text{etc.},
\end{equation*}
multiply this equation by $x^n$ that

\begin{equation*}
mSx^n+\frac{x^{n+1}dS}{dx}=(m+1)ax^{n+1}+(m+2)bx^{n+2}+\text{etc.};
\end{equation*}
differentiate this equation, where $dx$ is assumed to be constant, and divide by $dx$; this gives

\begin{equation*}
mnSx^{n-1}+\frac{(m+n+1)x^nS}{dx}+\frac{x^{n+1}ddS}{dx^2}
\end{equation*}
\begin{equation*}
=(m+1)(n+1)ax^n+(m+2)(n+2)bx^{n+1}+\text{etc.}
\end{equation*}
Divide by $x^{n-1}$ and multiply by $k$, so that

\begin{equation*}
mnkS+\frac{(m+n+1)kxdS}{dx}+\frac{kx^2ddD}{dx^2}
\end{equation*}
\begin{equation*}
=(m+1)(n+1)kax+(m+2)(n+2)kbx^2+(m+3)(n+3)kcx^3+\text{etc.}
\end{equation*}
Now, compare this series to the initial one; it will be

\begin{center}
\begin{tabular}{ccc}
 & Diff. I & Diff.  II \\
$ 
\renewcommand{\arraystretch}{2,0}
\setlength{\arraycolsep}{0.0mm}
\begin{array}{rrrrrrrrrrrrrrrrr}
kmn &~+~&1km &~+~&1kn &~ +~&1k &~=~&\alpha \\
knm &~+~&2km &~+~&2kn &~ +~&4k &~=~&\alpha &~+~&1\beta  \\
lnm &~+~&3km &~+~&3kn &~ +~&9k &~=~&\alpha &~+~&2 \beta &~ +~&\gamma
\end{array}$ &
$ \renewcommand{\arraystretch}{2,0}
\setlength{\arraycolsep}{0.0mm}
\begin{array}{rrrrrrrrrrrr}
km &~+~&kn &~+~&3k &~=~&\beta \\
km &~+~&kn &~+~&5k &~=~&\beta &~+~& \gamma
\end{array} $ &
$ \renewcommand{\arraystretch}{2,0}
\setlength{\arraycolsep}{0.0mm}
\begin{array}{rrrrrr}
2k &~=~& \gamma
\end{array}$ 
\end{tabular}
\end{center}
Therefore, $k=\frac{1}{2}\gamma$ and $m+n=\frac{2 \beta}{\gamma}-3$ and

\begin{equation*}
mn=\frac{\alpha}{k}-m-n-1=\frac{2 \alpha}{\gamma}-\frac{2 \beta}{\gamma}+2= \frac{2 (\alpha - \beta +\gamma)}{\gamma}.
\end{equation*}
Therefore, the sum of the series in question will be

\begin{equation*}
(\alpha -\beta +\gamma)S+\frac{(\beta- \gamma)xdS}{dx}+\frac{\gamma x^2 ddS}{2dx^2}.
\end{equation*}

\paragraph*{§26}

In like manner, one will be able to find the sum of this series

\begin{equation*}
Aa+Bbx+Ccx^2+Ddx^3+Eex^4+\text{etc.},
\end{equation*}
if the sum $S$ of this series was known, of course, i.e.

\begin{equation*}
S=a+bx+cx^2+dx^3+dx^3+ex^4+fx^5+\text{etc.}
\end{equation*}
and $A$, $B$, $C$, $D$ etc. constitute a series eventually reduced to constant differences. For, since its general form is concluded from the preceding, assume this sum

\begin{equation*}
\alpha S+ \frac{\beta x dS}{dx}+\frac{\gamma x^2 ddS}{2dx^2}+\frac{\delta x^3 d^3 S}{6 dx^3}+\frac{\varepsilon x^4 d^4 S}{24dx^4}+\text{etc.}
\end{equation*}
Now to find the letters $\alpha$, $\beta$, $\gamma$, $\delta$ etc., expand each series and it will be

\begin{equation*}
\renewcommand{\arraystretch}{2,5}
\setlength{\arraycolsep}{0.0mm}
\begin{array}{rrrrrrrrrrrrrrrr}
\quad \alpha S &~=~& \alpha a &~+~&\alpha b x&~+~& \alpha cx^2 &~+~&\alpha dx^3 &~+~& \alpha e x^4 &~+~&\text{etc.} \\
\dfrac{\beta x dS}{dx}  &~=~&  &~+~&\beta b x&~+~&2 \beta cx^2 &~+~&3\beta dx^3 &~+~&4 \beta e x^4 &~+~&\text{etc.} \\
\dfrac{\gamma x^2 ddS}{2dx^2}  &~=~&  & & &~+~& \gamma cx^2 &~+~&3 \gamma dx^3 &~+~&6 \gamma e x^4 &~+~&\text{etc.} \\
\dfrac{\delta x^3 d^3S}{6dx^3}  &~=~& & & & & &~+~& \delta dx^3 &~+~&4 \delta e x^4 &~+~&\text{etc.} \\
\dfrac{\varepsilon x^4 d^4S}{24dx^4}  &~=~&  && &&  && &~+~& \varepsilon e x^4 &~+~&\text{etc.} \\
 & & & & & \text{etc.}
\end{array}
\end{equation*}
compare this series, having arranged it according to the powers of $x$, to the propounded one, i.e.

\begin{equation*}
Z=Aa+Bbx+Ccx^2+Ddx^3+Eex^4+\text{etc.}
\end{equation*} 
and having made the comparison for each term, we find

\begin{alignat*}{9}
&\alpha &&=A \\
&\beta  &&=B- \alpha =B-A\\
&\gamma &&=C-2 \beta - \alpha  =C-2B+A\\
&\delta &&=D-3 \gamma -3 \beta - \alpha =D-3C+3B-A\\
&       &&  \qquad \qquad \text{etc.}
\end{alignat*}
Having found these values, the sum in question will therefore be

\begin{equation*}
Z= AS+\frac{(B-A)xdS}{1dx}+\frac{(C-2B+A)x^2ddS}{1 \cdot 2 dx^2}+\frac{(D-3C+3B-A)x^3d^3 S}{1 \cdot 2 \cdot 3 dx^3}+\text{etc.},
\end{equation*}
or, if the differences of the series $A$, $B$, $C$, $D$, $E$ etc. are denoted as usual, it will be

\begin{equation*}
Z=AS+\frac{\Delta A \cdot xdS}{1 dx}+\frac{\Delta^2 A \cdot x^2 d^2 S}{1 \cdot 2 dx^2}+\frac{\Delta^3 A \cdot x^3 d^3 S}{1 \cdot 2 \cdot 3 dx^3}+\text{etc.},
\end{equation*} 
if, as we assumed, it was

\begin{equation*}
S=a+bx+cx^2+dx^3+ex^4+fx^5+\text{etc.}
\end{equation*}
Therefore, if the series $A$, $B$, $C$, $D$ etc.  has eventually constant differences, one will be able to express the sum of the series $Z$ finitely.

\paragraph*{§27}

Since, having taken $e$ for the number whose hyperbolic logarithm is $=1$,  

\begin{equation*}
e^x =1+\frac{x}{1}+\frac{x^2}{1 \cdot 2}+\frac{x^3}{1 \cdot 2 \cdot 3}+\frac{x^4}{1 \cdot 2 \cdot 3 \cdot 4}+\frac{x^5}{1 \cdot 2 \cdot 3 \cdot 4 \cdot 5}+\text{etc.},
\end{equation*}
assume this series for the first, and because $S=e^x$, it will be $\frac{dS}{dx}=e^x$, $\frac{ddS}{dx^2}=e^x$ etc. Therefore, the sum of this series composed of that one and this one: $A$, $B$, $C$, $D$ etc., i.e. the series

\begin{equation*}
A+\frac{Bx}{1}+\frac{Cx^2}{1 \cdot 2}+\frac{Dx^3}{1 \cdot 2 \cdot 3}+\frac{Ex^4}{1 \cdot 2 \cdot 3 \cdot 4}+\text{etc.}
\end{equation*} 
will be expressed this way

\begin{equation*}
e^x \left(A+\frac{x \Delta A}{1}+\frac{xx \Delta^2 A}{1 \cdot 2}+\frac{x^3 \Delta^3 A}{1 \cdot 2 \cdot 3}+\frac{x^4 \Delta^4 A}{1 \cdot 2 \cdot 3 \cdot 4}+\text{etc.}\right).
\end{equation*}
Hence, if this series is propounded

\begin{equation*}
2+\frac{5x}{1}+\frac{10x^2}{1 \cdot 2}+\frac{17 x^2}{1 \cdot 2 \cdot 3}+\frac{26 x^4}{1 \cdot 2 \cdot 3 \cdot 4}+\frac{37 x^5}{1 \cdot 2 \cdot 3 \cdot 4 \cdot 5}+\text{etc.},
\end{equation*}
because of the series

\begin{equation*}
A, \quad B, \quad C, \quad D, \quad E \quad \text{etc.},
\end{equation*}
\begin{alignat*}{19}
&A &&=\quad && 2, \quad && \quad &&5, \quad && \quad &&10, \quad && \quad && 17, \quad && \quad && 26 \quad && \quad && \text{etc.}\\
&\Delta A&&=\quad && \quad && 3, \quad && \quad &&5, \quad && \quad &&7, \quad && \quad && 9 \quad && \quad && \text{etc.}\\
&\Delta \Delta A&&=\quad && \quad && \quad && 2, \quad && \quad &&2, \quad && \quad &&2 \quad && \quad  && \text{etc.},
\end{alignat*}
the sum of this series

\begin{equation*}
2+5x+\frac{10x^2}{2}+\frac{17x^3}{6}+\frac{26x^4}{24}+\text{etc.}
\end{equation*}
will be

\begin{equation*}
=e^x(2+3x+xx)=e^x(1+x)(2+x),
\end{equation*}
which is immediately clear. For, 

\begin{equation*}
\renewcommand{\arraystretch}{2,5}
\setlength{\arraycolsep}{0.0mm}
\begin{array}{rcccccccccccccccccc}
2e^x &~=~& 2&~+~&\dfrac{2x}{1}&~+~&\dfrac{2x^2}{2}&~+~&\dfrac{2x^3}{6}&~+~&\dfrac{2x^4}{24}&~+~&\text{etc.} \\
3xe^x &~=~& &~+~&3x&~+~&\dfrac{3x^2}{1}&~+~&\dfrac{3x^2}{2}&~+~&\dfrac{3x^4}{6}&~+~&\text{etc.} \\
xxe^x &~=~ & & & &~+~&xx&~+~&\dfrac{x^3}{1}&~+~&\dfrac{x^4}{2}&~+~&\text{etc.}\\ 
\text{and in total} \\ \hline
e^x(1+3x+xx) &~=~&2 &~+~&5x&~+~&\dfrac{10xx}{2}&~+~&\dfrac{17x^3}{6}&~+~&\dfrac{24	x^4}{24}&~+~&\text{etc.} \\
\end{array}
\end{equation*}

\paragraph*{§28}

The things discussed up to this point apply not only to infinite series, but also to sums of series consisting of a finite number of terms; for, the coefficients $a$, $b$, $c$, $d$ etc. can either be continued to infinity or can terminate at any arbitrary point. But because this does not require any further explanation, let us consider in more detail what follows from the results found up to this point. Therefore, having propounded any arbitrary series, whose  terms each consist of two factors, the one group of which factors  constitutes a series leading to constant differences, one will be able to assign the sum of this series, as long as, having omitted these factors, the sum was summable. Of course, if this series is propounded

\begin{equation*}
Z=Aa+Bbx+Ccx^2+Ddx^3+Eex^4+\text{etc.}
\end{equation*}
in which the quantities $A$, $B$, $C$, $D$, $E$ etc. constitute a series eventually reduced to constant differences, one will be able to exhibit the sum of this series, if  the sum $S$ of the following  series is known

\begin{equation*}
S=a+bx+cx^2+dx^3+ex^4+\text{etc.}
\end{equation*}
For, having calculated the iterated differences of the progression $A$, $B$, $C$, $D$, $E$ etc., as we showed at the beginning of this book,

\begin{alignat*}{19}
&A, && && B, && && C,&& &&D, && &&E,&& &&F,&& &&\text{etc.} \\
& &&\Delta A&& &&\Delta B,&& && \Delta C,&& && \Delta D, && && \Delta E && &&\text{etc.} \\
& && &&\Delta^2 A&& &&\Delta^2 B,&& && \Delta^2 C,&& && \Delta^2 D&& &&\text{etc.} \\
& && && &&\Delta^3 A&& &&\Delta^3 B,&& && \Delta^3 C,&& &&\text{etc.} \\
& && && && &&\Delta^4 A&& &&\Delta^4 B,&& &&\text{etc.} \\
& && && && && &&\Delta^5 A&& &&\text{etc.} \\
& && && && && && &&\text{etc.} \\
\end{alignat*}
the sum of the propounded series will be

\begin{equation*}
Z=SA+\frac{xdS}{1dx}\Delta A+\frac{x^2 ddS}{1 \cdot 2 dx^2}\Delta^2 A+\frac{x^3 d^3 S}{1 \cdot 2 \cdot 3 dx^3}\Delta^3 A+\text{etc.},
\end{equation*}
having put $dx$ to be constant in the higher powers of $S$, of course.

\paragraph*{§29}

Therefore, if the series $A$, $B$, $C$, $D$ etc. never leads to constant differences, the sum of the series $Z$ will be expressed by means of a new infinite series converging more rapidly than the initial one, and hence this series will be transformed into another one equal to it. To illustrate this, let this series be propounded

\begin{equation*}
Y=y+\frac{y^2}{2}+\frac{y^3}{3}+\frac{y^4}{4}+\frac{y^5}{5}+\frac{y^6}{6}+\text{etc.},
\end{equation*}
which is known to express $\log \frac{1}{1-y}$ such that $Y=-\log(1-y)$. Divide the series by $y$ and put $y=x$ and $Y=yZ$, so that

\begin{equation*}
Z=-\frac{1}{y}\log(1-y)=-\frac{1}{x}\log(1-x);
\end{equation*}
it will be

\begin{equation*}
Z=1+\frac{x}{2}+\frac{x^2}{3}+\frac{x^3}{4}+\frac{x^4}{5}+\frac{x^5}{6}+\text{etc.},
\end{equation*}
which series compared to

\begin{equation*}
S=1+x+x^2+x^3+x^4+x^5+x^6+\text{etc.}=\frac{1}{1-x}
\end{equation*}
will give these values for the series $A$, $B$, $C$, $D$, $E$ etc.

\begin{alignat*}{19}
&1, && &&\quad \frac{1}{2}, && && \quad \frac{1}{3}, && &&\quad \frac{1}{4}, && &&\quad \frac{1}{5} && &&\text{etc.} \\[3mm]
& &&-\frac{1}{1 \cdot 2}, && && \quad -\frac{1}{2 \cdot 3}, && &&\quad -\frac{1}{3 \cdot 4}, && &&-\frac{1}{4\cdot 5}  &&  &&\text{etc.} \\[3mm]
& && &&\frac{1\cdot 2}{1 \cdot 2 \cdot 3}, && &&\frac{1 \cdot 2}{2 \cdot 3 \cdot 4}, && && \frac{1\cdot 2}{3 \cdot 4 \cdot 5} &&  &&\text{etc.} \\[3mm]
& && && &&-\frac{1\cdot 2\cdot 3}{1 \cdot 2 \cdot 3 \cdot 4}, && &&-\frac{1 \cdot 2 \cdot 3}{2 \cdot 3 \cdot 4 \cdot 5}  &&  &&\text{etc.} \\
& && && && && \quad \text{etc.}
\end{alignat*}
Therefore, it will be

\begin{equation*}
A=1, \quad \Delta A=-\frac{1}{2}, \quad \Delta^2 A=\frac{1}{3}, \quad \Delta^3 A=-\frac{1}{4} \quad \text{etc.}
\end{equation*}
Further, because $S=\frac{1}{1-x}$, it will be

\begin{equation*}
\frac{dS}{dx}=\frac{1}{(1-x)^2}, \quad \frac{ddS}{1 \cdot 2 dx^2}=\frac{1}{(1-x)^3}, \quad \frac{d^3 S}{1 \cdot 2 \cdot 3 dx^3}=\frac{1}{(1-x)^4} \quad \text{etc.}
\end{equation*}
Having substituted these values, this sum will result

\begin{equation*}
Z=\frac{1}{1-x}-\frac{x}{2(1-x)^2}+\frac{x^2}{3(1-x)^3}-\frac{x^3}{4(1-x)^4}+\frac{x^4}{5(1-x)^5}-\text{etc.}
\end{equation*}
Therefore, because $x=y$ and $Y=-\log(1-y)=yZ$, it will be

\begin{equation*}
-\log(1-y)=\frac{y}{1-y}-\frac{y^2}{2(1-y)^2}+\frac{y^3}{3(1-y)^3}-\frac{y^4}{4(1-y)^4}+\text{etc.},
\end{equation*}
which series obviously expresses $\log \left(1+\frac{y}{1-y}\right)=\log \frac{1}{1-y}=-\log(1-y)$, the validity of which is even immediate considering  the results demonstrated before.

\paragraph*{§30}

To illustrate the application for a case in which only odd powers occur and the signs alternate, let this series be propounded

\begin{equation*}
Y=y-\frac{y^3}{3}+\frac{y^5}{5}-\frac{y^7}{7}+\frac{y^9}{9}-\frac{y^{11}}{11}+\text{etc.},
\end{equation*}
which  equation is equivalent to $Y=\arctan y$.\\[2mm]
Divide this series by $y$ and put $\frac{Y}{y}=Z$ and $yy=x$; it will be

\begin{equation*}
Z=1-\frac{x}{3}+\frac{xx}{5}-\frac{x^3}{7}+\frac{x^4}{9}-\frac{x^5}{11}+\text{etc.}
\end{equation*}
If it is compared to the following

\begin{equation*}
S=1-x+xx-x^3+x^4-x^5+\text{etc.},
\end{equation*}
it will be $S=\frac{1}{1+x}$ and the series of coefficients $A$, $B$, $C$, $D$ etc. will become

\begin{alignat*}{19}
&A&&= \quad &&1, \quad && \quad &&\frac{1}{3}, \quad && \quad && \frac{1}{5}, \quad && \quad  &&\frac{1}{7}, \quad && \quad && \frac{1}{9}, \quad && \quad && \text{etc.} \\[3mm]
&\Delta A&&= \quad && \quad && -\frac{2}{3}, \quad && \quad &&-\frac{2}{3 \cdot 5}, \quad && \quad && -\frac{2}{5 \cdot 7}, \quad && \quad  &&-\frac{2}{7 \cdot 9} \quad && \quad && \text{etc.} \\[3mm]
&\Delta^2 A&&= \quad && \quad && \quad && \frac{2\cdot 4}{3 \cdot 5}, \quad && \quad &&\frac{2 \cdot 4}{3 \cdot 5 \cdot 7}, \quad && \quad && \frac{2 \cdot 4}{5 \cdot 7 \cdot 9} \quad &&  \quad && \text{etc.} \\[3mm]
&\Delta^3 A&&= \quad && \quad && \quad && \quad && -\frac{2\cdot 4 \cdot 6}{3 \cdot 5 \cdot 7}, \quad && \quad &&-\frac{2 \cdot 4 \cdot 6}{3 \cdot 5 \cdot 7 \cdot 9} \quad &&  \quad && \text{etc.} \\[3mm]
&\Delta^4 A&&= \quad && \quad && \quad && \quad && \quad && \frac{2\cdot 4 \cdot 6 \cdot 8}{3 \cdot 5 \cdot 7 \cdot 9} \quad && \quad  \text{etc.} \\
& && \quad && \quad && \quad && \quad && \quad && \quad \text{etc.} 
\end{alignat*}
But, since $S=\frac{1}{1+x}$, it will be

\begin{equation*}
\frac{dS}{1 dx}=-\frac{1}{(1+x)^2}, \quad \frac{ddS}{1 \cdot 2 dx^2}=\frac{1}{(1+x)^3}, \quad \frac{d^3 S}{1 \cdot 2 \cdot 3 dx^3}=-\frac{1}{(1+x)^4} \quad \text{etc.}
\end{equation*}
Hence, having substituted these values, the form will become

\begin{equation*}
Z=\frac{1}{1+x}+\frac{2x}{3(1+x)^2}+\frac{2 \cdot 4 x^2}{3 \cdot 5(1+x)^3}+\frac{2 \cdot 4 \cdot 6 x^3}{3 \cdot 5 \cdot 7(1+x)^4}+\text{etc.};
\end{equation*}
having  substituted $x=yy$ again and multiplied by $y$, it will be

\begin{equation*}
Y= \arctan y = \frac{y}{1+yy}+\frac{2y^3}{3(1+yy)^2}+\frac{2 \cdot 4 y^5}{3 \cdot 5 (1+yy)^3}+\frac{2 \cdot 4 \cdot 6 y^7}{3 \cdot 5 \cdot 7(1+yy)^4}+\text{etc.}
\end{equation*}

\paragraph*{§31}

One can also transform the above series expressing the arc of a circle in another way by comparing it to the logarithmic series. \\[2mm]
For, let us consider the series

\begin{equation*}
Z=1-\frac{x}{3}+\frac{xx}{5}-\frac{x^3}{7}+\frac{x^4}{9}-\frac{x^{5}}{11}+\text{etc.},
\end{equation*}
which we want to compare to this one

\begin{equation*}
S=\frac{1}{0}-\frac{x}{2}+\frac{xx}{4}-\frac{x^3}{6}+\frac{x^4}{8}-\text{etc.}=\frac{1}{0}-\frac{1}{2}\log(1+x),
\end{equation*}
and the values of the letters $A$, $B$, $C$, $D$ etc. will be

\begin{alignat*}{19}
&A &&= \quad &&\frac{0}{1}, \quad && \quad && \frac{2}{3}, \quad && \quad && \frac{4}{5}, \quad && \quad && \frac{6}{7}, \quad && \quad && \frac{8}{9} \quad && \quad && \text{etc.}\\
&\Delta A &&= \quad && \quad &&\frac{2}{3}, \quad && \quad && \frac{+2}{3 \cdot 5}, \quad && \quad && \frac{+2}{5 \cdot 7}, \quad && \quad && \frac{+2}{7 \cdot 9} \quad && \quad && \text{etc.}\\
&\Delta^2 A &&= \quad && \quad && \quad &&\frac{-2 \cdot 4}{3 \cdot 5}, \quad && \quad && \frac{-2 \cdot 4}{3 \cdot 5 \cdot 7}, \quad && \quad && \frac{-2 \cdot 4}{5 \cdot 7 \cdot 9} \quad && \quad && \text{etc.}\\
&\Delta^3 A &&= \quad && \quad && \quad && \quad &&\frac{2 \cdot 4 \cdot 6}{3 \cdot 5 \cdot 7} \quad && \quad \text{etc.}\\
& && \quad && \quad && \quad && \quad  && \quad \text{etc.}\\
\end{alignat*}
Further, since $S=\frac{1}{0}-\frac{1}{2}\log(1+x)$, it will be

\begin{equation*}
\frac{dS}{1dx}=-\frac{1}{2(1+x)}, \quad \frac{ddS}{1 \cdot 2 dx^2}=\frac{1}{4(1+x)^2}, \quad 
\end{equation*}
\begin{equation*}
\frac{d^3 S}{1 \cdot 2 \cdot 3 dx^3}=-\frac{1}{6(1+x)^3}, \quad \frac{d^4 S}{1 \cdot 2 \cdot 3 \cdot 4 dx^4}=\frac{1}{8(1+x)^4} \quad \text{etc.}
\end{equation*}
Therefore, it will be $SA=S\frac{0}{1}=1$ and from the remaining terms it will be

\begin{equation*}
Z=1-\frac{x}{3(1+x)}-\frac{2xx}{3 \cdot 5(1+x)^2}-\frac{2 \cdot 4 x^3}{3 \cdot 5 \cdot 7(1+x)^3}-\text{etc.}
\end{equation*}
Now, let us put $x=yy$ and multiply by $y$; it will be

\begin{equation*}
Y= \arctan y= y-\frac{y^3}{3(1+yy)}-\frac{2 y^5}{3 \cdot 5(1+yy)^2}-\frac{2 \cdot 4 y^7}{3 \cdot 5 \cdot 7 (1+yy)^3}-\text{etc.}
\end{equation*}
This transformation will therefore not be obstructed by the infinite term $\frac{1}{0}$ entering the series $S$. But if  there remains any doubt, just expand each but the first term into power series in $y$ and one will discover that indeed the series propounded initially results.

\paragraph*{§32}

Up to this point we considered only series in which all powers  of the variable occurred. Now, we want proceed to other series which in each term contain the same power of the variable; the following series is of this kind

\begin{equation*}
S=\frac{1}{a+x}+\frac{1}{b+x}+\frac{1}{c+x}+\frac{1}{d+x}+\text{etc.}
\end{equation*}
For, if the sum $S$ of this series was known and is expressed by a certain function of $x$, by differentiating and by dividing by $-dx$, it will be

\begin{equation*}
\frac{-dS}{dx}=\frac{1}{(a+x)^2}+\frac{1}{(b+x)^2}+\frac{1}{(c+x)^2}+\frac{1}{(d+x)^2}+\text{etc.}
\end{equation*} 
If this series is differentiated again and divided by $-2dx$, one will recognize the series of the cubes

\begin{equation*}
\frac{ddS}{2dx^2}=\frac{1}{(a+x)^3}+\frac{1}{(b+x)^3}+\frac{1}{(c+x)^3}+\frac{1}{(d+x)^3}+\text{etc.}
\end{equation*} 
and this series, differentiated again and divided by $-3dx$, will give

\begin{equation*}
\frac{-d^3S}{dx^3}=\frac{1}{(a+x)^4}+\frac{1}{(b+x)^4}+\frac{1}{(c+x)^4}+\frac{1}{(d+x)^4}+\text{etc.}
\end{equation*} 
And in the same way, the sum of all following powers will be found, if the sum of the first series was known.

\paragraph*{§33}

But we found series of fractions of this kind involving a variable quantity in the \textit{Introductio}, where we showed, if the  half of the circumference of the circle, whose radius is $=1$, is set $=\pi$, that 

\begin{alignat*}{9}
&~\frac{\pi}{n \sin \frac{m}{n}\pi}&&=\frac{1}{m}&&+\frac{1}{n-m}&&-\frac{1}{n+m}&&-\frac{1}{2n-m}&&+\frac{1}{2n+m}&&+\frac{1}{3n-m}&&-\text{etc.}\\[3mm]
&\frac{\pi \cos \frac{m}{n}\pi}{n \sin \frac{m}{n}\pi}&&=\frac{1}{m}&&-\frac{1}{n-m}&&+\frac{1}{n+m}&&-\frac{1}{2n-m}&&+\frac{1}{2n+m}&&-\frac{1}{3n-m}&&+\text{etc.}
\end{alignat*}
Therefore, because it is possible to assume any arbitrary numbers for $m$ and $n$, let us set $n=1$ and $m=x$ so that we obtain a series similar to that one we  propounded in the preceding paragraph; having done this, it will be

\begin{alignat*}{9}
&~~\frac{\pi }{\sin \pi x}&&=\frac{1}{x}&&+\frac{1}{1-x}&&-\frac{1}{1+x}&&-\frac{1}{2-x}&&+\frac{1}{2+x}&&+\frac{1}{3-x}&&-\text{etc.}\\[3mm]
&\frac{\pi \cos \pi }{\sin \pi x}&&=\frac{1}{x}&&-\frac{1}{1-x}&&+\frac{1}{1+x}&&-\frac{1}{2-x}&&+\frac{1}{2+x}&&-\frac{1}{3-x}&&+\text{etc.}\\
\end{alignat*}
Therefore, one will be able to exhibit the sums of any powers of fractions resulting from these fractions by means of differentiations.

\paragraph*{§34}

Let us consider the first series and for the sake of brevity put $\frac{\pi}{\sin \pi x}=S$; take its higher differentials  alwas assuming $dx$ to be constant, and it will be

\begin{alignat*}{9}
&~~~~ \quad S&&=~\frac{1}{x}&&+~~\frac{1}{1-x}&&-~~\frac{1}{1+x}&&-~~\frac{1}{2-x}&&+~~\frac{1}{2+x}&&+~~\frac{1}{3-x}&&-\text{etc.}\\
&~~\frac{-dS}{dx}&&=\frac{1}{xx}&&-\frac{1}{(1-x)^2}&&-\frac{1}{(1+x)^2}&&+\frac{1}{(2-x)^2}&&+\frac{1}{(3+x)^2}&&-\frac{1}{(3-x)^2}&&-\text{etc.}\\
&~~~\frac{ddS}{2d^2x}&&=\frac{1}{x^3}&&+\frac{1}{(1-x)^3}&&-\frac{1}{(1+x)^3}&&-\frac{1}{(2-x)^3}&&+\frac{1}{(3+x)^3}&&+\frac{1}{(3-x)^3}&&-\text{etc.}\\
&~~\frac{-d^3S}{6d^3x}&&=\frac{1}{x^4}&&-\frac{1}{(1-x)^4}&&-\frac{1}{(1+x)^4}&&+\frac{1}{(2-x)^4}&&+\frac{1}{(3+x)^4}&&-\frac{1}{(3-x)^4}&&-\text{etc.}\\
&~~\frac{d^4S}{24d^4x}&&=\frac{1}{x^5}&&+\frac{1}{(1-x)^5}&&-\frac{1}{(1+x)^5}&&-\frac{1}{(2-x)^5}&&+\frac{1}{(3+x)^5}&&+\frac{1}{(3-x)^5}&&-\text{etc.}\\
&\frac{-d^5S}{120d^3x}&&=\frac{1}{x^6}&&-\frac{1}{(1-x)^6}&&-\frac{1}{(1+x)^6}&&+\frac{1}{(2-x)^6}&&+\frac{1}{(3+x)^6}&&-\frac{1}{(3-x)^6}&&-\text{etc.}\\
&  && && && &&\text{etc.}
\end{alignat*}
where it is to be noted that in the series of even powers the signs follow the same law and in like manner in the series of odd powers the structure of the signs is always the same. Therefore, the sums of all these series are found from the differentials of the expression $S=\frac{\pi}{\sin \pi x}$. 

\paragraph*{§35}

To express this differentials more conveniently, let us put

\begin{equation*}
\sin \pi =p \quad \text{and} \quad \cos \pi =q;
\end{equation*}
it will be

\begin{equation*}
dp=\pi dx \cos \pi x =\pi q dx \quad \text{and} \quad dq=- \pi p dx.
\end{equation*}
Therefore, because $S=\frac{\pi}{p}$, it will be

\begin{alignat*}{9}
&~\frac{-dS}{dx}&&=\frac{\pi^2 q}{pp} \\
&~~\frac{ddS}{dx^2}&&=\frac{\pi^3(pp+2qq)}{p^3}=\frac{\pi^3(qq+1)}{p^3} \quad \text{since} \quad pp+qq=1\\
&\frac{-d^3S}{dx^3}&&=\pi^4 \left(\frac{5q}{pp}+\frac{6q^3}{p^4}\right)=\frac{\pi^4(q^3+5q)}{p^4}\\
&~~\frac{d^4 S}{dx^4}&&=\pi^5 \left(\frac{24q^4}{p^5}+\frac{28q^2}{p^3}+\frac{5}{p}\right)=\frac{\pi^5(q^4+18q^2+5)}{p^5}\\
&\frac{-d^5 S}{dx^5}&&=\pi^6 \left(\frac{120 q^5}{p^6}+\frac{180q^3}{p^4}+\frac{61q}{pp}\right)=\frac{\pi^6(q^5+58q^3+61q)}{p^6}\\
&~~\frac{d^6 S}{dx^6}&&=\pi^7 \left(\frac{720 q^6}{p^7}+\frac{1320 q^4}{p^5}+\frac{662q^2}{p^3}+\frac{61}{p}\right)=\frac{\pi^7(q^6+179q^4+479q^2+61)}{p^7}\\
&\frac{-d^7 S}{dx^7}&&=\pi^8 \left(\frac{5040q^7}{p^8}+\frac{10920 q^5}{p^6}+\frac{7266 q^3}{p^4}+\frac{1385 q}{p^2}\right)
\intertext{or}
& &&=\frac{\pi^8}{p^8}(q^7+543q^5+3111q^3+1385q)\\
&~~\frac{d^8 S}{dx^8}&&=\pi^ 9 \left(\frac{40320q^8}{p^9}+\frac{100800 q^6}{p^7}+\frac{83664 q^4}{p^5}+\frac{24568 q^2}{p^3}+\frac{1385}{p}\right)
\intertext{or}
& &&=\frac{\pi^9}{p^9}(q^8+1636q^6+18270q^4+19028q^2+1385)\\
& && \qquad \qquad \text{etc.}
\end{alignat*}
These expressions are easily continued arbitrarily far; for, if it was

\begin{equation*}
\pm \frac{d^n S}{dx^n}= \pi^{n+1} \left(\frac{\alpha q^n}{p^{n+1}}+\frac{\beta q^{n-2}}{p^{n-1}}+\frac{\gamma q^{n-4}}{p^{n-3}}+\frac{\delta q^{n-6}}{p^{n-5}}+\text{etc.}\right),
\end{equation*}
then its differential, having changed the signs, will be

\begin{equation*}
\mp \frac{d^{n+1}S}{dx^{n+1}}
\left\lbrace
\begin{aligned}
&(n+1) \alpha \frac{q^{n+1}}{p^{n+2}}+(n\alpha +(n-1)\beta)\frac{q^{n-1}}{p^n}+((n-2)\beta +(n-3)\gamma)\frac{q^{n-3}}{p^{n-2}} \\
&  +((n-4)\gamma +(n-5)\delta)\frac{q^{n-5}}{p^{n-4}}+\text{etc.}
\end{aligned}
\right\rbrace
\end{equation*}

\paragraph*{§36}

Therefore, from these series one will obtain the following sums of the  series exhibited in § 34 

\begin{alignat*}{9}
 S &=\pi \cdot \frac{1}{p}\\
 \frac{-dS}{dx} &= \frac{\pi^2}{1} \cdot \frac{q}{p^2}\\ 
 \frac{ddS}{24dx^2}&= \frac{\pi^3}{2}\left(\frac{2q^2}{p^3}+\frac{1}{p}\right)\\
 \frac{-d^3 S}{6dx^3}&=\frac{\pi^4}{6}\left(\frac{6q^3}{p^4}+\frac{5q}{p^2}\right)\\
\frac{d^4 S}{24 dx^4}& =\frac{\pi^5}{24}\left(\frac{24 q^4}{p^5}+\frac{28q^2}{p^3}+\frac{5}{p}\right)\\
\frac{-d^5 S}{120 dx^5}& =\frac{\pi^6}{120}\left(\frac{120 q^5}{p^6}+\frac{180 q^3}{p^4}+\frac{61q}{p^2}\right)\\
\frac{d^6 S}{720 dx^6}&= \frac{\pi^7 }{720}\left(\frac{720 q^6}{p^7}+\frac{1320 q^4}{p^5}+\frac{662 q^2}{p^3}+\frac{61}{p}\right) \\
\frac{-d^7 S}{720 dx^6}&= \frac{\pi^8}{5040}\left(\frac{5040 q^7}{p^8}+\frac{10920 q^5}{p^6}+\frac{7266q^3}{p^4}+\frac{1385 q}{p^2}\right)\\
\frac{d^8 S}{40320 dx^8}&=\frac{\pi^9}{40320}\left(\frac{40320 q^8}{p^9}+\frac{100800q^6}{p^7}+\frac{83664q^4}{p^5}+\frac{24568 q^2}{p^3}+\frac{1385}{p}\right)\\
 &\text{etc.}
\end{alignat*}

\paragraph*{§37}

Let us treat the other series found above [§ 33] in the same way, i.e. the series

\begin{equation*}
\frac{\pi \cos \pi x}{\sin \pi x}=\frac{1}{x}-\frac{1}{1-x}+\frac{1}{1+x}-\frac{1}{2-x}+\frac{1}{2+x}-\frac{1}{3-x}+\text{etc.}
\end{equation*}
and, for the sake of brevity  having put $\frac{\pi \cos \pi x}{\sin \pi x}=T$,  the following summations will result

\begin{alignat*}{9}
&\quad ~~~~ T &&=~\frac{1}{x} &&-~~\frac{1}{1-x}&&+~~\frac{1}{1+x}&&-~~\frac{1}{2-x}&&+~~\frac{1}{2+x}&&-\text{etc.} \\
&~~~\frac{-dT}{dx}&&=\frac{1}{x^2}&&+\frac{1}{(1-x)^2}&&+\frac{1}{(1+x)^2}&&+\frac{1}{(2-x)^2}&&+\frac{1}{(2+x)^2}&&+\text{etc.}\\
&~~~\frac{ddT}{2dx^2}&&=\frac{1}{x^3}&&-\frac{1}{(1-x)^3}&&+\frac{1}{(1+x)^3}&&-\frac{1}{(2-x)^3}&&+\frac{1}{(2+x)^2}&&-\text{etc.}\\
&~~\frac{-d^3T}{6d^3x}&&=\frac{1}{x^4}&&+\frac{1}{(1-x)^4}&&+\frac{1}{(1+x)^4}&&+\frac{1}{(2-x)^4}&&+\frac{1}{(2+x)^4}&&+\text{etc.}\\
&~~\frac{d^4T}{24dx^4}&&=\frac{1}{x^5}&&-\frac{1}{(1-x)^5}&&+\frac{1}{(1+x)^5}&&-\frac{1}{(2-x)^5}&&+\frac{1}{(2+x)^5}&&-\text{etc.}\\
&\frac{-d^5T}{120d^5x}&&=\frac{1}{x^6}&&+\frac{1}{(1-x)^6}&&+\frac{1}{(1+x)^6}&&+\frac{1}{(2-x)^6}&&+\frac{1}{(2+x)^6}&&+\text{etc.}\\
& && && && \text{etc.},
\end{alignat*}
where in the series of even powers all terms are positive, but in the series of odd  powers the signs $+$ and $-$ alternate.

\paragraph*{§38}

To find the values of these differentials, as before, let us put

\begin{equation*}
\sin \pi x =p \quad \text{and} \quad dq=- \pi pdx
\end{equation*}
that $pp+qq=1$; it will be

\begin{equation*}
dp = \pi q dx \quad \text{and} \quad dq=-\pi p dx.
\end{equation*}
Having added these values, it will be

\begin{alignat*}{9}
T &= \pi \cdot \frac{ q}{p}\\
\frac{-dT}{dx}&= \pi^2 \left(\frac{qq}{pp}+1\right)=\frac{\pi^2}{pp} \\
\frac{ddT}{dx^2}&=\pi^3 \left(\frac{2q^3}{p^3}+\frac{2q}{p}\right)=\frac{2 \pi^3 q}{p^3}\\
\frac{-d^3 T}{dx^3}&=\pi^4 \left(\frac{6q^4}{p^4}+\frac{8qq}{pp}+2\right)=\pi^4 \left(\frac{6qq}{p^4}+\frac{2}{pp}\right)\\
\frac{d^4 T}{dx^4}&= \pi^5 \left(\frac{24q^3}{p^5}+\frac{16q}{p^3}\right)\\
\frac{-d^5 T}{dx^5}&=\pi^6 \left(\frac{120q^4}{p^6}+\frac{120 qq}{p^4}+\frac{16}{pp}\right)\\
\frac{d^6 T}{dx^6}&=\pi^7 \left(\frac{720 q^5}{p^7}+\frac{960q^3}{p^5}+\frac{272q}{p^3}\right)\\
\frac{-d^7 T}{dx^7}&=\pi^8 \left(\frac{5040q^6}{p^8}+\frac{8400 q^4}{p^6}+\frac{3696 q^2}{q^4}+\frac{272}{p^2}\right)\\
\frac{d^8 T}{dx^8}&= \pi^9 \left(\frac{40320 q^7}{p^9}+\frac{80640q^5}{p^7}+\frac{48384 q^3}{p^5}+\frac{7936q}{p^3}\right)\\
 &\text{etc.}
\end{alignat*}
These formulas can easily be continued arbitrarily far. For, if 

\begin{equation*}
\pm \frac{d^n T}{dx^n}=\pi^{n+1} \left(\frac{\alpha q^{n-1}}{p^{n+1}}+\frac{\beta q^{n-3}}{p^{n-1}}+\frac{\gamma q^{n-5}}{p^{n-3}}+\frac{\delta q^{n-7}}{p^{n-5}}+\text{etc.}\right),
\end{equation*}
the expression for the following differential will be

\begin{equation*}
\mp \frac{d^{n+1}T}{dx^{n+1}}=\pi^{n+2} \left(\frac{(n+1)\alpha q^n}{p^{n+2}}+\frac{(n-1)(\alpha +\beta)q^{n-2}}{p^n}+\frac{(n-3)(\beta +\gamma)q^{n-4}}{p^{n-2}}+\text{etc.}\right)
\end{equation*}

\paragraph*{§39}

Therefore,  having put $\sin \pi x =p$ and $\cos \pi  x =q$, the series of powers given in § 37 will have the following sums 

\begin{alignat*}{9}
T &= \pi \cdot \frac{q}{p} \\
\frac{-dT}{dx}&=\pi^2 \frac{1}{pp} \\
\frac{ddT}{2 dx^2}&=\pi^3 \frac{q}{p^3}\\
\frac{-d^3 T}{6dx^3}&= \pi^4 \left(\frac{qq}{\pi^4}+\frac{1}{3pp}\right) \\
\frac{-d^4 T}{24 dx^4}&=\pi^5 \left(\frac{q^3}{p^5}+\frac{2q}{3p^3}\right)\\
\frac{-d^5 T}{120 dx^5}&= \pi^6 \left(\frac{q^4}{p^6}+\frac{3qq}{p^4}+\frac{2}{15pp}\right)\\
\frac{d^6 T}{720 dx^6}&= \pi^7 \left(\frac{q^5}{p^7}+\frac{4 q^3}{3 p^5}+\frac{17q}{45p^3}\right)\\
\frac{-d^7 T}{5040 dx^7}&= \pi^8 \left(\frac{q^6}{p^8}+\frac{5q^4}{3p^6}+\frac{11 q^2}{15p^4}+\frac{17}{315pp}\right)\\
\frac{d^8 T}{40320 dx^8}&=\pi^9 \left(\frac{q^7}{p^9}+\frac{6q^5}{3p^7}+\frac{6q^3}{5p^5}+\frac{62q}{315p^3}\right)\\
& \text{etc.}
\end{alignat*}

\paragraph*{§40}

Aside from these series, we found several others in the \textit{Introductio} from which  others can be derived by means  of differentiation in the same way.\\[2mm]
For, we showed that

\begin{equation*}
\frac{1}{2x}-\frac{\pi \sqrt{x}}{2 x \tan \pi \sqrt{x}}=\frac{1}{1-x}+\frac{1}{4-x}+\frac{1}{9-x}+\frac{1}{16-x}+\frac{1}{25-x}+\text{etc.}
\end{equation*}
Let us put that the sum of this series is $=S$, so that 

\begin{equation*}
S=\frac{1}{2x}-\frac{\pi}{2 \sqrt{x}}\cdot \frac{\cos \pi \sqrt{x}}{\sin \pi \sqrt{x}};
\end{equation*}
it will be

\begin{equation*}
\frac{dS}{dx}=-\frac{1}{2xx}+\frac{\pi}{4x \sqrt{x}}\cdot \frac{\cos \pi \sqrt{x}}{\sin \pi \sqrt{x}}+\frac{\pi \pi}{4x(\sin \pi \sqrt{x})^2},
\end{equation*}
which expression therefore yields the sum of this series

\begin{equation*}
\frac{1}{(1-x)^2}+\frac{1}{(4-x)^2}+\frac{1}{(9-x)^2}+\frac{1}{(16-x)^2}+\frac{1}{(25-x)^2}+\text{etc.}
\end{equation*}
Further, we also showed that

\begin{equation*}
\frac{\pi}{2 \sqrt{x}}\cdot \frac{e^{2\pi \sqrt{x}}+1}{e^{2 \pi \sqrt{x}}-1}-\frac{1}{2x}=\frac{1}{1+x}+\frac{1}{4+x}+\frac{1}{9+x}+\frac{1}{16+x}+\text{etc.}
\end{equation*}
Therefore, if this sum is put $=S$, it will be

\begin{equation*}
\frac{-dS}{dx}=\frac{1}{(1+x)^2}+\frac{1}{(4+x)^2}+\frac{1}{(9+x)^2}+\frac{1}{(16+x)^2}+\text{etc.}
\end{equation*}
But

\begin{equation*}
\frac{dS}{dx}=\frac{- \pi}{4x\sqrt{x}}\cdot \frac{e^{2\pi \sqrt{x}}+1}{e^{2\pi \sqrt{x}}-1}-\frac{\pi \pi}{x}\cdot \frac{e^{2\pi \sqrt{x}}}{(e^{2 \pi \sqrt{x}}-1)^2}+\frac{1}{2xx}.
\end{equation*}
Therefore, the sum of this series will be

\begin{equation*}
\frac{-dS}{dx}=\frac{\pi}{4x \sqrt{x}}\cdot \frac{e^{2 \pi \sqrt{x}}+1}{e^{2 \pi \sqrt{x}}-1}+\frac{ \pi \pi }{x}\cdot \frac{e^{2 \pi \sqrt{x}}}{(e^{2 \pi \sqrt{x}}-1)^2}-\frac{1}{2xx}.
\end{equation*}
And in like manner  the sums of the following powers will be found by means of further differentiation.

\paragraph*{§41}

If the value of a certain product composed of factors involving the  variable letter $x$ was known, one will be able to find innumerable summable series from it by means of the same method. For, let the value of this product

\begin{equation*}
(1+\alpha x)(1+\beta x)(1+\gamma x)(1+\delta x)(1+\varepsilon x)\text{etc.}
\end{equation*}
be $=S$, a function of $x$, of course;  taking logarithms it will be

\begin{equation*}
\log S= \log (1+\alpha x)+\log (1+\beta x)+\log (1+\gamma x)+\log (1+\delta x)+\text{etc.}
\end{equation*}
Now take the differentials; after division by $dx$ it will be

\begin{equation*}
\frac{dS}{Sdx}=\frac{\alpha }{1+\alpha x}+\frac{\beta}{1 +\beta x}+\frac{\gamma}{1+ \gamma x}+\frac{\delta }{1+ \delta x}+\text{etc.},
\end{equation*}
from  further differentiation of which series the sums of any powers of these fractions will be found,  precisely as we explained it in more detail in the preceding examples.

\paragraph*{§42}

But, in the \textit{Introductio} we exhibited several expressions we want to apply this method to. If $\pi$ is the arc of $180^{\circ}$ of the circle whose radius is $=1$, we showed that 

\begin{equation*}
\sin \frac{m \pi}{2n}=\frac{m \pi}{2n}\cdot \frac{4nn-mm}{4nn}\cdot \frac{16nn -mm}{16nn}\cdot \frac{26nn-mm}{36nn}\cdot \text{etc.}
\end{equation*}
\begin{equation*}
\cos \frac{m \pi}{2n}= \frac{nn -mm}{nn}\cdot \frac{9nn-mm}{9nn}\cdot \frac{25nn-mm}{25nn}\cdot \frac{49nn-mm}{49nn}\cdot \text{etc.}
\end{equation*}
Let us put $n=1$ and $m=2x$, so that

\begin{equation*}
\sin \pi x= \pi x \cdot \frac{1-xx}{1}\cdot \frac{4-xx}{4}\cdot \frac{9-xx}{9}\cdot \frac{16-xx}{16}\cdot \text{etc.}
\end{equation*}
or

\begin{equation*}
\sin \pi x= \pi x \cdot \frac{1-x}{1}\cdot \frac{1+x}{1}\cdot \frac{2-x}{2}\cdot \frac{2+x}{2}\cdot \frac{3-x}{3}\cdot \frac{3+x}{3}\cdot \frac{4-x}{4}\cdot \text{etc.}
\end{equation*}
and

\begin{equation*}
\cos \pi x = \frac{1-4xx}{1}\cdot \frac{9-4xx}{9}\cdot \frac{25-4xx}{25}\cdot \frac{49-4xx}{49}\cdot \text{etc.}
\end{equation*}
or

\begin{equation*}
\cos \pi x= \frac{1-2x}{1}\cdot \frac{1+2x}{1}\cdot \frac{3-2x}{3}\cdot \frac{3+2x}{3}\cdot \frac{5-2x}{5}\cdot \frac{5+2x}{5}\cdot \text{etc.}
\end{equation*}
Therefore, from these expressions, if one takes logarithms, it will be

\begin{alignat*}{9}
& \log \sin \pi x &&= \log \pi x +\log \frac{1-x}{1}+\log \frac{1+x}{1}+\log \frac{2-x}{2}+\log \frac{2+x}{2}+\log \frac{3-x}{3}+\text{etc.} \\[3mm]
& \log \cos \pi x &&= \log \frac{1-2x}{1}+\log \frac{1+2x}{1}+\log \frac{3-2x}{3}+\log \frac{3+2x}{3}+\log \frac{5-2x}{5}+\text{etc.}
\end{alignat*}

\paragraph*{§43}

Now let us take the differentials of these series of logarithms and, having divided by $dx$ everywhere, the first series will give

\begin{equation*}
\frac{\pi \cos \pi x}{\sin \pi x}= \frac{1}{x}-\frac{1}{1-x}+\frac{1}{1+x}-\frac{1}{2-x}+\frac{1}{2+x}-\frac{1}{3-x}+\text{etc.}
\end{equation*} 
which is the series we discussed in § 37. The other series on the other hand will give

\begin{equation*}
\frac{- \pi \sin \pi x}{\cos \pi x}=-\frac{2}{1-2x}+\frac{2}{1+2x}-\frac{2}{3-2x}+\frac{2}{3+2x}-\frac{2}{5-2x}+\text{etc.}
\end{equation*}
Let us put $2x=z$, so that $x=\frac{z}{2}$, and divide by $-2$; it will be

\begin{equation*}
\frac{\pi \sin \frac{1}{2}\pi z}{2 \cos \frac{1}{2}\pi z}= \frac{1}{1-z}-\frac{1}{1+z}+\frac{1}{3-z}-\frac{1}{3+z}+\frac{1}{5-z}-\text{etc.}
\end{equation*}
But, since

\begin{equation*}
\sin \frac{1}{2}\pi z=\sqrt{\frac{1- \cos \pi z}{2}} \quad \text{and} \quad \cos \frac{1}{2}\pi z= \sqrt{\frac{1+\cos \pi z}{2}},
\end{equation*}
it will be

\begin{equation*}
\frac{\pi \sqrt{1- \cos \pi z}}{\sqrt{1+ \cos \pi z}}=\frac{2}{1-z}-\frac{2}{1+z}+\frac{2}{3-z}+\frac{2}{5-z}-\text{etc.}
\end{equation*}
or, writing $x$ instead of $z$,

\begin{equation*}
\frac{\pi \sqrt{1- \cos \pi x}}{\sqrt{1+\cos \pi x}}=\frac{2}{1-x}-\frac{2}{1+x}+\frac{2}{3-x}-\frac{2}{3+x}+\frac{2}{5-x}-\text{etc.}
\end{equation*}
Add this series to the one found first

\begin{equation*}
\frac{\pi \cos \pi x}{\sin \pi x}=\frac{1}{x}-\frac{1}{1-x}+\frac{1}{1+x}-\frac{1}{2-x}+\frac{1}{2+x}-\frac{1}{3-x}+\text{etc.}
\end{equation*}
and one will find the sum of this series

\begin{equation*}
\frac{1}{x}+\frac{1}{1-x}-\frac{1}{1+x}-\frac{1}{2-x}+\frac{1}{2+x}+\frac{1}{3-x}-\frac{1}{3+x}-\text{etc.}
\end{equation*}
to be $=\frac{\pi \sqrt{1- \cos \pi x}}{\sqrt{1+\cos \pi}}+\frac{\pi \cos \pi x}{\sin \pi x}$. But this fraction $\frac{\sqrt{1- \cos \pi x}}{\sqrt{1+\cos \pi x}}$, if the numerator and denominator are multiplied by $\sqrt{1- \cos \pi x}$, goes over into $\frac{1- \cos \pi x}{\sin \pi}$. Therefore, the sum of the series will be $=\frac{\pi }{\sin \pi x}$, which is the series we considered in § 34; therefore, we will not prosecute this any further.


\title{On Finding finite differences
\footnote{
Original title: "'De Inventione Differentiarum finitarum"', first published as part of the book "`\textit{Institutiones calculi differentialis cum eius usu in analysi finitorum ac doctrina serierum}, 1755"', reprinted in \textit{Opera Omnia}: Series 1, Volume 10, pp. 256 - 275, Eneström-Number E212, translated by: Alexander Aycock for the "`Euler-Kreis Mainz"'}}
\author{Leonhard Euler}
\date{ }
\maketitle

\paragraph*{§44}

At the beginning of this book we explained, how from the finite differences of the functions their differentials can easily be found, and  even derived the principle of differentials from this source. For, if the differences, if they were assumed to be finite, vanish and go over into zero, the differentials result; and because in this case many and often innumerable terms, which constitute the finite difference, are neglected, the differentials can be found a lot more easily and can be expressed both more conveniently and succinctly than the finite differences. And therefore, there seems to be no way to ascend from differentials to finite differences. Nevertheless, by the method we will use here, one will be able to determine the finite differences from the differentials of all orders of any function.

\paragraph*{§45}

Let $y$ be an arbitrary function of $x$; because this function, having written $x+dx$ instead of $x$, goes over into $y+dy$, if  one writes $x+dx$ instead of $x$ once more, the value $y+dy$ will be increased by its differential $dy+ddy$  and it will be $=y+2dy+ddy$ which value therefore corresponds to $x+2dx$  written  instead of $x$ in the function $y$.  In the same way, if we assume that $x$ is continuously increased by its differential $dx$ that it successively takes on the values $x+dx$, $x+2dx$, $x+3dx$, $x+4dx$ etc., the corresponding values of $y$ will be those seen in the following table.

\begin{tabular}{c|c}
Values of & Corresponding Values of the Function \\
 $
\renewcommand{\arraystretch}{1,4}
\setlength{\arraycolsep}{0.0mm}
\begin{array}{rrrr}
 x & \\
 x &~+~ &dx \\
 x &~+&2dx \\
 x &~+~&3dx \\
 x &~+~&4dx \\
 x &~+~&5dx \\
 x &~+~&6dx \\
   & \text{etc.}
\end{array}
 $
 &
 $ 
 \renewcommand{\arraystretch}{1,3}
\setlength{\arraycolsep}{0.0mm}
\begin{array}{rrrrrrrrrrrrrrr}
 & y \\
 & y &~+~&dy \\
 & y &~+~&2dy &~ +~&ddy \\
 & y &~+~&3dy &~ +~&3ddy  &~+~&d^3 y \\
 & y &~+~&4dy &~ +~&6ddy  &~+~&4d^3y &~+~&d^4y \\
 & y &~+~&5dy &~+~ &10ddy &~+~&10d^3y&~+~&5d^4 y &~+~&d^5y  \\
 & y &~+~&6dy &~+~ &15ddy &~+~&20d^3y&~+~&15 d^4y &~+~&6d^5 y &~+~&d^6 y\\
 &   &   &    &    & \text{etc.}  
 \end{array}$
\end{tabular}

\paragraph*{§46}

Therefore, if in general $x$ goes over into $x+ndx$, the function $y$ will take on this form

\begin{equation*}
y+\frac{n}{1}dy+\frac{n(n-1)}{1 \cdot 2}ddy+\frac{n(n-1)(n-2)}{1 \cdot 2 \cdot 3}d^3y +\frac{n(n-1)(n-2)(n-3)}{1 \cdot 2 \cdot 3 \cdot 4}d^4y+\text{etc.};
\end{equation*}
even though in this expression each term is infinite times smaller than its preceding term, we nevertheless did not omit any term, so that this formula  can be used for the present task. For, we will assume $n$ to be an infinitely large number, and since we know that it can happen that the product of an infinitely large and an infinitely small quantity is equal to a finite quantity, the second term can certainly become homogeneous to the first term, i.e. the quantity $ndy$ can represent a finite quantity. For the same reason, the third term $\frac{n(n-1)}{1 \cdot 2}ddy$, even though $ddy$ is infinite times smaller than $dy$, because the one factor $\frac{n(n-1)}{1 \cdot 2}$ is infinite times larger than $n$, can also express a finite quantity; and hence, having put $n$ to be an infinite number, it is not possible to neglect any term of that expression.

\paragraph*{§47}

But having put $n$ to be an infinite number, by whatever finite number it is either increased or decreased, the resulting number will have the ratio of $1$ to $n$ and hence one can  write the number $n$ for each of the factors $n-1$, $n-2$, $n-3$, $n-4$ etc. everywhere. For, because $\frac{n(n-1)}{1 \cdot 2}ddy= \frac{1}{2}nnddy-\frac{1}{2}nddy$, the first term $\frac{1}{2}nnddy$ will have the same ratio to the second $\frac{1}{2}nddy$ as  $n$ to $1$ and that second term will vanish with respect to the first; therefore,  one will be able to write $\frac{1}{2}nn$ instead of $\frac{n(n-1)}{1\cdot 2}$. In like manner, the coefficient of the fourth term $\frac{n(n-1)(n-3)}{1 \cdot 2 \cdot 3}$ can be contracted to $\frac{n^3}{6}$ and in like manner one can neglect the numbers subtracted from $n$ in the factors in the following. But having done this, the function $y$, if in  one writes $x+ndx$ instead of $x$, while $n$ is an infinite number,  will have the following value

\begin{equation*}
y+\frac{ndy}{1}+\frac{nnddy}{1 \cdot 2}+\frac{n^3 d^3y}{1 \cdot 2 \cdot 3}+\frac{n^4 d^4 y}{1 \cdot 2 \cdot 3 \cdot 4}+\frac{n^5 d^5y}{1 \cdot 2 \cdot 3 \cdot 4 \cdot 5}+\text{etc.}
\end{equation*}

\paragraph*{§48}

Therefore, because having assumed $n$ as infinitely large number, even though $dx$ is infinitely small, the product $ndx$ can express a finite quantity, let us put $ndx= \omega$, so that $n=\frac{\omega}{dx}$; $n$ will certainly be an infinite number being the quotient resulting from a division of the finite quantity $\omega$ by the infinitely small quantity $dx$. But having used this value  instead of $n$, if the variable quantity $x$ is increased by a certain quantity $\omega$ or if  one writes $x +\omega$ instead of $x$, a function $y$ of that variable $x$ will go over into this form

\begin{equation*}
y+\frac{\omega dy}{1 dx}+\frac{\omega^2 ddy}{1 \cdot 2 dx^2}+\frac{\omega^3 d^3 y}{1 \cdot 2 \cdot 3 dx^3}+\frac{\omega^4 d^4 y}{1 \cdot 2 \cdot 3 \cdot 4 dx^4}+\text{etc.},
\end{equation*}
each term of which expression can be found by iterated differentiation of $y$. For, because $y$ is a function $x$, we showed above that these functions $\frac{dy}{dx}$, $\frac{ddy}{dx^2}$, $\frac{d^3y}{dx^3}$ etc. all exhibit finite quantities.

\paragraph*{§49}

Therefore, because, while the variable quantity $x$ is increased by the finite quantity $\omega$, any function $y$ of that $x$ is increased by its first difference, which we indicated by $\Delta y$ above, where we had $\omega =\Delta x$, one will be able to find the difference of $y$ by iterated differentiation; for, it will be

\begin{equation*}
\Delta y= \frac{ \omega dy}{dx}+\frac{\omega^2 ddy}{2dx^2}+\frac{\omega^3 d^3y}{6dx^3}+\frac{\omega^4 d^4y}{24 dx^4}+\text{etc.}
\end{equation*}
or 
\begin{equation*}
\Delta y= \frac{\Delta x}{1}\cdot \frac{d y}{d x}+\frac{\Delta x^2}{2}\cdot \frac{ddy}{dx^2}+\frac{\Delta x^3 }{6} \cdot \frac{d^3 y}{dx^3}+\frac{\Delta x^4}{24}\cdot \frac{d^4y}{dx^4}+\text{etc.}
\end{equation*}
And hence the finite difference $\Delta y$ is expressed by a progression whose  terms proceed in powers of $\Delta x$. Therefore, vice versa it is obvious, if the quantity $x$ is increased only by an infinitely small quantity, i.e. $\Delta x$ goes over into its differential $dx$, that all terms  vanish with respect to the first term and that it will be $\Delta y= dy$; for, having set $\Delta x=dx$, the difference $\Delta y$, by definition, goes over into the differential $dy$.

\paragraph*{§50}

Because, if  one writes $x+\omega$ instead of $x$, any function $y$ of that variable $x$ then has the following value

\begin{equation*}
y+\frac{\omega dy}{dx}+\frac{\omega^2 ddy}{2dx^2}+\frac{\omega^3 d^3 y}{6 dx^3}+\frac{\omega^4 d^4 y}{24 dx^4}+\text{etc.};
\end{equation*}
the validity of this expression can checked in examples in which the higher differentials of $y$ finally vanish; for, in these cases the number of terms of the above expression will become finite.

\subsection*{Example 1} 

\textit{Let the value of the expression $xx-x$ be in question, if  one writes $x+1$ instead of $x$}.\\[2mm]
Put $y=xx-x$, and because we assume that $x$ goes over into $x+1$, it will be $\omega =1$; now, having taken the differentials, it will be

\begin{equation*}
\frac{dy}{dx}=2x-1, \quad \frac{ddy}{dx^2}=2, \quad \frac{d^3 y}{dx^3}=0 \quad \text{etc.}
\end{equation*}
Therefore, the function $y=xx-x$, having written $x+1$ instead of $x$, will go over into

\begin{equation*}
xx-x+1(2x-1)+\frac{1}{2}\cdot 2= xx+x.
\end{equation*}
But if in $xx-x$ one actually writes $x+1$ instead of $x$,

\begin{equation*}
\renewcommand{\arraystretch}{1,4}
\setlength{\arraycolsep}{0mm}
\begin{array}{rrrrrrrrrrrrrrrrrrrrrrr}
 &xx & \quad \text{goes over into} \quad & xx &~+~&2x&~+~&1 \\
 &x  & \quad \text{goes over into} \quad &    &   & x&~+~&1 \\
\text{Therefore} \quad \\ \cline{2-8}
 &x-xx& \quad \text{goes over into} \quad &xx &~+~&x &
\end{array}
\end{equation*}

\subsection*{Example 2} 

\textit{Let the value of the expression $x^3+xx+x$ be in question, if  one puts $x+2$ instead of $x$}.\\[2mm]
Put $y=x^3+xx+x$ and it will be $\omega =2$; now, since

\begin{equation*}
y=x^3+xx+x,
\end{equation*}
it will be

\begin{equation*}
\frac{dy}{dx}=3xx+2x+1, \quad \frac{ddy}{dx^2}=6x+2, \quad \frac{d^3y}{dx^3}=6, \quad \frac{d^4y}{dx^4}=0 \quad \text{etc.}
\end{equation*}
Hence the value of the function $y=x^3+xx+x$, if  one substitutes $x+2$ for $x$,  will be the following

\begin{equation*}
x^3+xx+x+2(3xx+2x+1)+\frac{4}{2}(6x+2)+\frac{8}{6}\cdot 6=x^3+7xx+17x+14,
\end{equation*}
which same function arises, if $x+2$ is actually substituted for $x$.

\subsection*{Example 3}

\textit{Let the value of the expression $xx+3x+1$ be in question,  if one writes $x-3$ instead of $x$.}\\[2mm]
Therefore, it will be $\omega=-3$, and having put

\begin{equation*}
y=xx+3x+1,
\end{equation*}
it will be

\begin{equation*}
\frac{dy}{dx}=2x+3, \quad \frac{ddy}{dx^2}=2, \quad \frac{d^3y}{dx^3}=0 \quad \text{etc.},
\end{equation*}
whence, written put $x-3$ instead of $x$, the function $x^2+3x+1$ will go over into

\begin{equation*}
x^2+3x+1-\frac{3}{1}(2x+3)+\frac{9}{2}\cdot 2= x^2-3x+1.
\end{equation*}

\paragraph*{§51}

If a negative number is taken for $\omega$, one will find the value, which a function of $x$ has, if the quantity $x$ is decreased by the given quantity $\omega$. Of course, if  one writes $x- \omega$ instead of $x$, an arbitrary function $y$ of $x$ will then have this value

\begin{equation*}
y-\frac{\omega dy}{dx}+\frac{\omega^2 ddy}{2dx^2}-\frac{\omega^3 d^3 y}{6dx^3}+\frac{\omega^4 d^4 y}{24 dx^4}-\text{etc.},
\end{equation*}
whence all variations the function $y$ can undergo, while the quantity $x$ is changed arbitrarily, can be found. But if $y$ was a polynomial function of $x$, since in that case one eventually gets to vanishing differentials, the varied value will be expressed finitely; but if $y$ was not a function of this kind, the varied value will be expressed by  an infinite series, whose sum, since, if the substitution is actually done, the varied value is easily assigned, can be exhibited by a finite expression.

\paragraph*{§52}

But as the first difference was found, so also the following differences can be exhibited by similar expressions. For, let $x$ successively take on the values $x+\omega$, $x+2 \omega$, $x+3 \omega$, $x+4 \omega$ etc. and denote the corresponding values of $y$ by $y^{\text{I}}$, $y^{\text{II}}$, $y^{\text{III}}$, $y^{\text{IV}}$ etc., as we did at the beginning of this book. Therefore, since $y^{\text{I}}$, $y^{\text{II}}$, $y^{\text{III}}$, $y^{\text{IV}}$ etc. are the values $y$ will have, if  one respectively writes $x+ \omega$, $x+2 \omega$, $x+3 \omega$, $x+4 \omega$ etc. instead of $y$, using the demonstrated method, the values of these $y$s will be expressed this way:

\begin{equation*}
\renewcommand{\arraystretch}{2,5}
\setlength{\arraycolsep}{0.0mm}
\begin{array}{lccccccccccccccccc}
y^{\text{I}}&~=~&y&~+~&\dfrac{\omega dy}{dx}&~+~&\dfrac{\omega^2 ddy}{2dx^2}&~+~&\dfrac{\omega^3 d^3 y}{6dx^3}&~+~&\dfrac{\omega^4 d^4y}{24 dx^4}&~+~&\text{etc.}\\
y^{\text{II}}&~=~&y&~+~&\dfrac{2\omega dy}{dx}&~+~&\dfrac{4\omega^2 ddy}{2dx^2}&~+~&\dfrac{8\omega^3 d^3 y}{6dx^3}&~+~&\dfrac{16\omega^4 d^4y}{24 dx^4}&~+~&\text{etc.}\\
y^{\text{III}}&~=~&y&~+~&\dfrac{3\omega dy}{dx}&~+~&\dfrac{9\omega^2 ddy}{2dx^2}&~+~&\dfrac{27\omega^3 d^3 y}{6dx^3}&~+~&\dfrac{81\omega^4 d^4y}{24 dx^4}&~+~&\text{etc.}\\
y^{\text{IV}}&~=~&y&~+~&\dfrac{4\omega dy}{dx}&~+~&\dfrac{16\omega^2 ddy}{2dx^2}&~+~&\dfrac{64\omega^3 d^3 y}{6dx^3}&~+~&\dfrac{256\omega^4 d^4y}{24 dx^4}&~+~&\text{etc.}\\
& && && && \text{etc.}
\end{array}
\end{equation*}

\paragraph*{§53}

Therefore, because, if $\Delta y$, $\Delta^2 y$, $\Delta^3 y$, $\Delta^4 y$ etc. denote the first, second, third, fourth etc. differences, 

\begin{alignat*}{9}
&\Delta y   &&= y^{\text{I}}&&-y \\
&\Delta^2 y &&= y^{\text{II}}&&-2y^{\text{I}}&&+y\\
&\Delta^3 y &&= y^{\text{III}}&&-3y^{\text{II}}&&+3y^{\text{I}}&&-y \\
&\Delta^4 y &&= y^{\text{IV}}&&-4y^{\text{III}}&&+6y^{\text{II}}&&-4y^{\text{I}}&&+y\\
& && && &&\text{etc.},	
\end{alignat*}
these differences  will be expressed this way by means of differentials:

\begin{alignat*}{9}
&\Delta y&&=\frac{\omega dy}{dx}+\frac{\omega^2 ddy}{dx^2}+\frac{\omega^3 d^3y}{6dx^3}+\frac{\omega^4 d^4y}{24 dx^4}+\text{etc.}\\[3mm]
&\Delta^2 y &&= \frac{(2^2-2\cdot 1)\omega^2 ddy}{2 dx^2}+\frac{(2^3-2 \cdot 1)\omega^3 d^3 y}{6dx^3}+\frac{(2^4-2 \cdot 1)\omega^4 d^4y}{24 dx^4}+\text{etc.}\\[3mm]
&\Delta^3 y &&=\frac{(3^3-3 \cdot 2^3+3 \cdot 1)\omega^3 d^3 y}{6dx^3}+\frac{(3^4-3 \cdot 2^4+3 \cdot 1)\omega^4 d^4 y}{24 dx^4}+\text{etc.} \\[3mm]
&\Delta^4 y &&=\frac{(4^4-4 \cdot 3^4+6 \cdot 2^4-4 \cdot 1)\omega^4 d^4y}{24dx^4}+\frac{(4^5-4 \cdot 3^5+6 \cdot 2^5-4 \cdot 1)\omega^5 d^5 y}{120 dx^5}+\text{etc.} \\
& &&  \text{etc.}
\end{alignat*}

\paragraph*{§54}

It is immediately clear, of how much use these expressions of differences are in the doctrine of series and progressions,  and we will explain  it in greater detail later parts of this chapter. Meanwhile,  we want to consider applications following  immediately follows from this for the understanding of series. Although the indices of the terms of a certain series are usually assumed to constitute an arithmetic progression whose difference is $1$, for the sake of  broader and easier applicability, let us set the difference $=\omega$, such that, if the general term or that corresponding to the index $x$, was $y$, the following terms correspond to the indices $x+\omega$, $x+2 \omega$, $x+3 \omega$ etc. Therefore, if the following terms correspond to these indices

\begin{alignat*}{9}
&x, \quad && x+\omega, \quad && x+2 \omega, \quad && x+3 \omega, \quad && x+4 \omega \quad && \text{etc.} \\
&y, \quad && ~~~P, \quad && ~~~Q, \quad && ~~~R, \quad && ~~~S, \quad && \text{etc.}
\end{alignat*}
each term will be defined from $y$ and its differentials this way:

\begin{equation*}
\renewcommand{\arraystretch}{2,5}
\setlength{\arraycolsep}{0.0mm}
\begin{array}{cccccccccccccccc}
P &~=~& y &~+~&\dfrac{\omega dy}{dx}&~+~&\dfrac{\omega^2 ddy}{2 dx^2}&~+~&\dfrac{\omega^3 d^3y}{6 dx^3}&~+~&\dfrac{\omega^4 d^4y}{24dx^4}&~+~&\text{etc.} \\
Q &~=~& y &~+~&\dfrac{2\omega dy}{dx}&~+~&\dfrac{4\omega^2 ddy}{2 dx^2}&~+~&\dfrac{9\omega^3 d^3y}{6 dx^3}&~+~&\dfrac{16\omega^4 d^4y}{24dx^4}&~+~&\text{etc.} \\
R &~=~& y &~+~&\dfrac{3\omega dy}{dx}&~+~&\dfrac{9\omega^2 ddy}{2 dx^2}&~+~&\dfrac{27\omega^3 d^3y}{6 dx^3}&~+~&\dfrac{81\omega^4 d^4y}{24dx^4}&~+~&\text{etc.} \\
S &~=~& y &~+~&\dfrac{4\omega dy}{dx}&~+~&\dfrac{16\omega^2 ddy}{2 dx^2}&~+~&\dfrac{64\omega^3 d^3y}{6 dx^3}&~+~&\dfrac{256\omega^4 d^4y}{24dx^4}&~+~&\text{etc.} \\
T &~=~& y &~+~&\dfrac{5\omega dy}{dx}&~+~&\dfrac{25\omega^2 ddy}{2 dx^2}&~+~&\dfrac{125\omega^3 d^3y}{6 dx^3}&~+~&\dfrac{625\omega^4 d^4y}{24dx^4}&~+~&\text{etc.} \\
& && && && \text{etc.}
\end{array}
\end{equation*}

\paragraph*{§55}

If these expressions are subtracted from each other, $y$ will no longer enter the differences and it will be

\begin{equation*}
\renewcommand{\arraystretch}{2,5}
\setlength{\arraycolsep}{0.0mm}
\begin{array}{cccccccccccccccc}
P&~-~&y &~=~&\dfrac{\omega dy}{dx}&~+~&\dfrac{\omega^2 ddy}{2dx^2}&~+~&\dfrac{\omega^3 d^3y}{6dx^3}&~+~&\dfrac{\omega^4 d^4y}{24 dx^4}&~+~&\text{etc.} \\
Q&~-~&P &~=~&\dfrac{\omega dy}{dx}&~+~&\dfrac{3\omega^2 ddy}{2dx^2}&~+~&\dfrac{7\omega^3 d^3y}{6dx^3}&~+~&\dfrac{15\omega^4 d^4y}{24 dx^4}&~+~&\text{etc.} \\
R&~-~&Q &~=~&\dfrac{\omega dy}{dx}&~+~&\dfrac{5\omega^2 ddy}{2dx^2}&~+~&\dfrac{19\omega^3 d^3y}{6dx^3}&~+~&\dfrac{65\omega^4 d^4y}{24 dx^4}&~+~&\text{etc.} \\
S&~-~&R &~=~&\dfrac{\omega dy}{dx}&~+~&\dfrac{7\omega^2 ddy}{2dx^2}&~+~&\dfrac{37\omega^3 d^3y}{6dx^3}&~+~&\dfrac{175\omega^4 d^4y}{24 dx^4}&~+~&\text{etc.} \\
T&~-~&S &~=~&\dfrac{\omega dy}{dx}&~+~&\dfrac{9\omega^2 ddy}{2dx^2}&~+~&\dfrac{61\omega^3 d^3y}{6dx^3}&~+~&\dfrac{369\omega^4 d^4y}{24 dx^4}&~+~&\text{etc.} \\
& && && &&\text{etc.}
\end{array}
\end{equation*}
If these expressions are again subtracted from each other, the first differentials will also cancel each other and it will be

\begin{equation*}
\renewcommand{\arraystretch}{2,5}
\setlength{\arraycolsep}{0.0mm}
\begin{array}{ccccccccccccc}
Q&~-~&2P &~+~&y &~=~& \dfrac{2 \omega^2 ddy}{2dx^2}&~+~&\dfrac{6 \omega^3 d^3 y}{6 dx^3}&~+~&\dfrac{14 \omega^4 d^4y}{24 dx^4}&~+~&\text{etc.}\\
R&~-~&2Q &~+~&P &~=~& \dfrac{2 \omega^2 ddy}{2dx^2}&~+~&\dfrac{12 \omega^3 d^3 y}{6 dx^3}&~+~&\dfrac{50 \omega^4 d^4y}{24 dx^4}&~+~&\text{etc.}\\
S&~-~&2R &~+~&Q &~=~& \dfrac{2 \omega^2 ddy}{2dx^2}&~+~&\dfrac{18 \omega^3 d^3 y}{6 dx^3}&~+~&\dfrac{110 \omega^4 d^4y}{24 dx^4}&~+~&\text{etc.}\\
T&~-~&2S &~+~&R &~=~& \dfrac{2 \omega^2 ddy}{2dx^2}&~+~&\dfrac{24 \omega ^3 d^3 y}{6 dx^3}&~+~&\dfrac{194 \omega^4 d^4y}{24 dx^4}&~+~&\text{etc.}\\
& && && &&  \text{etc.}
\end{array}
\end{equation*}
Having subtracted them from each other again, the second differentials will also go out of the computation:

\begin{equation*}
\renewcommand{\arraystretch}{2,5}
\setlength{\arraycolsep}{0.0mm}
\begin{array}{ccccccccccccccccccc}
R &~-~&3Q &~+~&3P &~-~&y&~=~& \dfrac{6 \omega^3 d^3 y}{6dx^3}&~+~&\dfrac{36 \omega^4 d^4 y}{24dx^4}&~+~&\text{etc.}\\
S &~-~&3R &~+~&3Q &~-~&P&~=~& \dfrac{6 \omega^3 d^3 y}{6dx^3}&~+~&\dfrac{60 \omega^4 d^4 y}{24dx^4}&~+~&\text{etc.}\\
T &~-~&3S &~+~&3R &~-~&Q&~=~& \dfrac{6 \omega^3 d^3 y}{6dx^3}&~+~&\dfrac{84 \omega^4 d^4 y}{24dx^4}&~+~&\text{etc.}\\
& && &&   \text{etc.}
\end{array}
\end{equation*}
By continuing the subtraction it will be

\begin{alignat*}{9}
&S &&-4R &&+6Q &&-4P &&+y &&=\frac{24 \omega^4 d^4 y}{24 dx^4}+\text{etc.}\\
&T &&-4S &&+6R &&-4Q &&+P &&=\frac{24 \omega^4 d^4 y}{24 dx^4}+\text{etc.}\\
& && && && && &&\text{etc.}
\end{alignat*}
and

\begin{equation*}
T-5S+10R-10 Q+5P-y =\frac{120 \omega^5 d^5 y}{120 dx^5}+\text{etc.}
\end{equation*}
\begin{equation*}
\text{etc.}
\end{equation*}

\paragraph*{§56}

Therefore, if $y$ was a polynomial function of $x$, since its higher differentials will vanish eventually, proceeding this way, one will reach  vanishing expressions. Therefore, because these expressions are differences of $y$, let us consider their forms and coefficients more diligently.

\begin{equation*}
\renewcommand{\arraystretch}{2,5}
\setlength{\arraycolsep}{0.0mm}
\begin{array}{rccccccccccccccccccccccccccccccccccc}
y&~=~&y \\
\Delta y &~=~&\dfrac{\omega dy}{dx}&~+~&\dfrac{\omega^2 ddy}{2dx^2}&~+~&\dfrac{\omega^3 d^3 y}{6dx^3}&~+~&\dfrac{\omega^4 d^4y}{24 dx^4}&~+~&\dfrac{\omega^5 d^5 y}{120 dx^5}&~+~&\text{etc.}\\
\Delta^2 y&~ =~& \dfrac{\omega^2 ddy}{dx^2}&~+~&\dfrac{3 \omega^3 d^3 y}{3 dx^3}&~+~&\dfrac{7 \omega^4 d^4 y}{3 \cdot 4 dx^4}&~+~&\dfrac{15 \omega^5 d^5 y}{3 \cdot 4 \cdot 5 dx^5}&~+~&\dfrac{31 \omega^6 d^6 y}{3 \cdot 4 \cdot 5 \cdot 6 dx^6}&~+~&\text{etc.}\\
 \Delta^3 y &~=~& \dfrac{\omega^3 d^3y}{dx^3}&~+~&\dfrac{6 \omega^4 d^4 y}{4 dx^4}&~+~&\dfrac{25 \omega^5 d^5 y}{4 \cdot 5 dx^5}&~+~&\dfrac{90 \omega^6 d^6 y}{4 \cdot 5 \cdot 6 dx^6}&~+~&\dfrac{301 \omega^7 d^7 y}{4 \cdot 5 \cdot 6 \cdot 7 dx^7}&~+~&\text{etc.}\\
\Delta^4 y&~=~& \dfrac{\omega^4 d^4 y}{dx^4}&~+~&\dfrac{10 \omega^5 d^5 y}{5 dx^5}&~+~&\dfrac{65 \omega^6 d^6 y}{5 \cdot 6 dx^6}&~+~&\dfrac{350 \omega^7 d^7y}{5 \cdot 6 \cdot 7 dx^7}&~+~&\dfrac{1701 \omega^8 d^8 y}{5 \cdot 6 \cdot 7 \cdot 8 dx^8}&~+~&\text{etc.}\\
\Delta^5 y&~=~& \dfrac{\omega^5 d^5 y}{dx^5}&~+~&\dfrac{15 \omega^6 d^6 y}{6 dx^6}&~+~&\dfrac{140 \omega^7 d^7 y}{6 \cdot 7 dx^7}&~+~&\dfrac{1050 \omega^8 d^8y}{6 \cdot 7 \cdot 8 dx^8}&~+~&\dfrac{6951 \omega^9 d^9 y}{6 \cdot 7 \cdot 8 \cdot 9 dx^9}&~+~&\text{etc.}\\
\Delta^6 y&~=~& \dfrac{\omega^6 d^6 y}{dx^6}&~+~&\dfrac{21 \omega^7 d^7 y}{7 dx^7}&~+~&\dfrac{266 \omega^8 d^8 y}{7 \cdot 8 dx^8}&~+~&\dfrac{2646 \omega^9 d^9 y}{7 \cdot 8 \cdot 9 dx^9}&~+~&\dfrac{22827 \omega^{10} d^{10} y}{7 \cdot 8 \cdot 9 \cdot 10 dx^{10}}&~+~&\text{etc.}\\
& && && && &&\text{etc.}
\end{array}
\end{equation*}

\paragraph*{§57}

Let us also consider the same series continued backwards at the same time, which contains the terms corresponding to the indices $x- \omega$, $x-2 \omega$, $x-3 \omega$ etc.

\begin{alignat*}{16}
&x-4 \omega \quad &&x-3\omega, \quad &&x- 2 \omega, \quad && x- \omega, \quad && x, \quad && x+ \omega, \quad && x+2 \omega, \quad && x+3 \omega, \quad && x+4 \omega \quad &&\text{etc.}\\
& \quad s, &&\quad r, && \quad q, && \quad p, &&y, && \quad P, && \quad Q, && \quad R, && \quad S &&  \text{etc.} 
\end{alignat*}
Therefore, since

\begin{equation*}
\renewcommand{\arraystretch}{2,5}
\setlength{\arraycolsep}{0.0mm}
\begin{array}{cccccccccccccccccccc}
p &~=~& y &~-~&\dfrac{\omega dy}{dx}&~+~&\dfrac{\omega^2 ddy}{2 dx^2}&~-~& \dfrac{\omega^3 d^3y}{6 dx^3}&~+~&\dfrac{\omega^4 d^4 y}{24 dx^4}&~-~&\text{etc.} \\
q &~=~& y &~-~& \dfrac{2\omega dy}{dx}&~+~&\dfrac{4\omega^2 ddy}{2 dx^2}&~-~& \dfrac{8\omega^3 d^3y}{6 dx^3}&~+~&\dfrac{16\omega^4 d^4 y}{24 dx^4}&~-~&\text{etc.} \\
r &~=~& y &~-~& \dfrac{3\omega dy}{dx}&~+~&\dfrac{9\omega^2 ddy}{2 dx^2}&~-~& \dfrac{27\omega^3 d^3y}{6 dx^3}&~+~&\dfrac{81\omega^4 d^4 y}{24 dx^4}&~-~&\text{etc.} \\
s &~=~& y &~-~& \dfrac{4\omega dy}{dx}&~+~&\dfrac{16\omega^2 ddy}{2 dx^2}&~-~& \dfrac{64\omega^3 d^3y}{6 dx^3}&~+~&\dfrac{256\omega^4 d^4 y}{24 dx^4}&~-~&\text{etc.} \\
& &&  && &&\text{etc.},
\end{array}
\end{equation*}
by subtracting these values from the above ones, i.e. from $P$, $Q$, $R$, $S$ etc., it will be

\begin{equation*}
\renewcommand{\arraystretch}{2,5}
\setlength{\arraycolsep}{0.0mm}
\begin{array}{cccccccccccccccc}
\dfrac{P-p}{2}&~=~&\dfrac{\omega dy}{dx}&~+~&\dfrac{\omega^3 d^3 y}{6 dx^3}&~+~&\dfrac{\omega^5 d^5 y}{120 dx^5}&~+~&\text{etc.} \\
\dfrac{Q-q}{2}&~=~&\dfrac{2\omega dy}{dx}&~+~&\dfrac{8\omega^3 d^3 y}{6 dx^3}&~+~&\dfrac{32\omega^5 d^5 y}{120 dx^5}&~+&\text{etc.} \\
\dfrac{R-r}{2}& ~=~&\dfrac{3\omega dy}{dx}&~+~&\dfrac{27\omega^3 d^3 y}{6 dx^3}&~+~&\dfrac{243\omega^5 d^5 y}{120 dx^5}&~+~&\text{etc.} \\
\dfrac{S-s}{2}&~=~&\dfrac{4\omega dy}{dx}&~+~&\dfrac{64\omega^3 d^3 y}{6 dx^3}&~+~&\dfrac{1024\omega^5 d^5 y}{120 dx^5}&~+~&\text{etc.} \\
& && &&   \text{etc.}
\end{array}
\end{equation*}
But if these terms  are added to the above ones, then, as  the differentials of even orders are missing here, the odd differentials will go out of the computation in this case. For, it will be

\begin{equation*}
\renewcommand{\arraystretch}{2,5}
\setlength{\arraycolsep}{0.0mm}
\begin{array}{ccccccccccccccccc}
\dfrac{P+p}{2}&~=~&y &~+~&\dfrac{\omega^2 ddy}{2dx^2}&~+~&\dfrac{\omega^4 d^4 y}{24 dx^4}&~+~~&\dfrac{\omega^6 d^6y}{720 dx^6}&~+~&\text{etc.}\\
 \dfrac{Q+q}{2}&~=~& y &~+~&\dfrac{4\omega^2 ddy}{2dx^2}&~+~&\dfrac{16\omega^4 d^4 y}{24 dx^4}&~+~&\dfrac{64\omega^6 d^6y}{720 dx^6}&~+~&\text{etc.}\\
\dfrac{R+r}{2}&~=~&y &~+~&\dfrac{9\omega^2 ddy}{2dx^2}&~+~&\dfrac{81\omega^4 d^4 y}{24 dx^4}&~+~&\dfrac{729\omega^6 d^6y}{720 dx^6}&~+~&\text{etc.}\\
\dfrac{S+s}{2}&~=~&y &~+~&\dfrac{16\omega^2 ddy}{2dx^2}&~+~&\dfrac{256\omega^4 d^4 y}{24 dx^4}&~+~&\dfrac{4096\omega^6 d^6y}{720 dx^6}&~+~&\text{etc.}\\
& && && &&\text{etc.}
\end{array}
\end{equation*}

\paragraph*{§58}

Since  the preceding terms can all be expressed, if they are collected into one sum, the summatory term of the propounded series will result. Let the first term correspond to the index $x-n \omega$ and the first term itself will be

\begin{equation*}
=y-\frac{n \omega dy}{dx}+\frac{n^2 \omega^2}{2dx^2}-\frac{n^3 \omega^3 d^3 y}{6 dx^3}+\frac{n^4 \omega^4 d^4 y}{24 dx^4}+\text{etc.}
\end{equation*}
Therefore, since the term corresponding to the index $x$ will be $=y$ and the number of all terms is $=n+1$, the sum of all, from the first up to $y$, i.e. the summatory term will be

\begin{alignat*}{9}
&=(n+1)y &&- ~~\frac{\omega dy}{dx}&&(1&&+2&&+3&&+\cdots &&+n)\\[3mm]
& &&+ ~\frac{\omega^2 ddy}{2dx^2}&&(1&&+2^2&&+3^2&&+\cdots &&+n^2)\\[3mm]
& &&- ~\frac{\omega^3 d^3y}{6dx^3}&&(1&&+2^3&&+3^3&&+\cdots &&+n^3)\\[3mm]
& &&+ ~\frac{\omega^4 d^4y}{24dx^4}&&(1&&+2^4&&+3^4&&+\cdots &&+n^4)\\[3mm]
& &&- \frac{\omega^5 d^5y}{120dx^5}&&(1&&+2^5&&+3^5&&+\cdots &&+n^5)\\
& &&+\text{etc.}
\end{alignat*}

\paragraph*{§59}

Above we exhibited the sums of each of these series [§ 62 of the first part]; if these are substituted here, the sum of the propounded series will be

\begin{alignat*}{9}
&=(n+1)y &&- ~~\frac{\omega dy}{dx}&&\left(\frac{1}{2}nn+\frac{1}{2}n\right)\\[3mm]
& &&+~\frac{\omega^2 ddy}{2dx^2}&&\left(\frac{1}{3}n^3+\frac{1}{2}nn+\frac{1}{6}n\right)\\[3mm]
& &&-~\frac{\omega^3 d^3y}{6dx^3}&&\left(\frac{1}{4}n^4+\frac{1}{2}n^3+\frac{1}{4}n^2\right)\\[3mm]
& &&+~\frac{\omega^4 d^4y}{24dx^4}&&\left(\frac{1}{5}n^5+\frac{1}{2}n^4+\frac{1}{3}n^3-\frac{1}{30}n\right)\\[3mm]
& &&-\frac{\omega^5 d^5y}{120dx^5}&&\left(\frac{1}{6}n^6+\frac{1}{2}n^5+\frac{5}{12}n^4-\frac{1}{12}n^2 \right)\\[3mm]
& && &&\text{etc.},
\end{alignat*}
where $n$ will be given from the index of the first term which is the initial term of the sum. If  one puts $\omega=1$ and the index of the first term is $=1$, the index of the second $=2$ and of the last $=x$ such that this series is propounded

\begin{alignat*}{9}
&1, \quad &&2, \quad && 3, \quad && 4,  &&\cdots \cdots &&x\\
&a, \quad &&b, \quad && c, \quad && d,  &&\cdots \cdots &&y,
\end{alignat*}
 because of $x-n=1$ and $n=x-1$, the sum of this series will be

\begin{alignat*}{9}
&=xy &&- ~~~~~\frac{dy}{dx}&&\bigg(\frac{1}{2}xx&&-\frac{1}{2}x\bigg)\\[3mm]
&    &&+~~~\frac{ddy}{2dx^2}&&\bigg(\frac{1}{3}x^3&&-\frac{1}{2}xx&&+\frac{1}{6}x\bigg)\\[3mm]  
&    &&-~~~\frac{d^3y}{6dx^3}&&\bigg(\frac{1}{4}x^4&&-\frac{1}{2}x^3&&+\frac{1}{4}xx\bigg)\\[3mm] 
&    &&+~~\frac{d^4y}{24dx^4}&&\bigg(\frac{1}{5}x^5&&-\frac{1}{2}x^4&&+\frac{1}{3}xx&&-\frac{1}{30}x\bigg)\\[3mm]  
&    &&+~\frac{d^5y}{120dx^4}&&\bigg(\frac{1}{6}x^6&&-\frac{1}{2}x^5&&+\frac{5}{12}x^4&&-\frac{1}{12}x^2\bigg)\\[3mm]
&    &&+~\frac{d^6y}{720dx^4}&&\bigg(\frac{1}{7}x^7&&-\frac{1}{2}x^6&&+\frac{1}{2}x^5&&-\frac{1}{6}x^3+\frac{1}{42}x\bigg)\\
& && && &&\text{etc.}  
\end{alignat*}

\paragraph*{§60}

From this expression, since the number of coefficients will increase immensely, if $x$ was a large number, hardly anything of use for the doctrine of series follows; nevertheless, it will be helpful to have mentioned other properties following from these considerations. Let the general term be $x^n$ and indicate the summatory term by $S.y$ or $S.x$. Having used this notation everywhere, it will be

\begin{alignat*}{9}
&\frac{1}{2}xx&&-\frac{1}{2}x=S.x-x\\[3mm]
&\frac{1}{3}x^3&&-\frac{1}{2}x^2+\frac{1}{6}x=S.x^2-x^2\\[3mm]
&\frac{1}{4}x^4&&-\frac{1}{2}x^3+\frac{1}{4}xx=S.x^3-x^3\\[3mm]
& && \qquad \qquad \text{etc.}
\end{alignat*}
Therefore, from the above expression one will obtain 

\begin{equation*}
S.x^n= x^{n+1}-nx^{n+1}S.x+nx^n
\end{equation*}
\begin{equation*}
+\frac{n(n-1)}{1 \cdot 2}x^{n-2}S.x^2-\frac{n(n-1)}{1 \cdot 2}x^n-\frac{n(n-1)(n-2)}{1\cdot 2 \cdot 3}x^{n-3}S.x^3+\frac{n(n-1)(n-2)}{1 \cdot 2 \cdot 3}x^n+\text{etc.}
\end{equation*}
But, because

\begin{equation*}
(1-1)^n=0=1-n+\frac{n(n-1)}{1\cdot 2}-\frac{n(n-1)(n-2)}{1 \cdot 2 \cdot 3}+\text{etc.},
\end{equation*}
it will be

\begin{equation*}
n-\frac{n(n-1)}{1\cdot 2}+\frac{n(n-1)(n-2)}{1 \cdot 2\cdot 3}-\text{etc.}=1
\end{equation*}
and hence, having excluded the case $n=0$ in which this expression becomes $=0$,

\begin{equation*}
S.x^n=x^{n+1}+x^n-nx^{n-1}S.x+\frac{n(n-1)}{1 \cdot 2}x^{n-2}S.x^2-\frac{n(n-1)(n-2)}{1 \cdot 2 \cdot 3}x^{n-3}S.x^3
\end{equation*}
\begin{equation*}
+\frac{n(n-1)(n-2)(n-3)}{1 \cdot 2 \cdot 3 \cdot 4}x^{n-4}S.x^4-\text{etc.}
\end{equation*}

\paragraph*{§61}

To see both the validity and the applicability of this formula more clearly, let us expand some special cases and at first let  $n=1$ and it will be $S.x= x^2+x-S.x$ and hence $S.x=\frac{xx+x}{2}$, which is well-known, of course. Therefore, let us put $n=2$ and it will be

\begin{equation*}
S.x^2=x^3+xx-2xS.x+S.x^2,
\end{equation*}
which equation, since the terms $S.x^2$ appearing on both sides cancel, give the same equation, i.e. $S.x=\frac{xx+x}{2}$. But if $n=3$, it will be

\begin{equation*}
S.x^3=\frac{3}{2}xS.x^2-\frac{3}{2}x^2S.x+\frac{1}{2}x^3(x+1);
\end{equation*}
if one puts $n=4$,  this expression will result

\begin{equation*}
S.x^4=x^5+x^4-4x^3S.x+6x^2S.x^2-4xS.x^3+S.x^4,
\end{equation*}
whence, because of the terms $S.x^4$ cancelling each other, it will be

\begin{equation*}
S.x^3=\frac{3}{2}xS.x^2-x^2S.x+\frac{1}{4}x^3(x+1);
\end{equation*} 
if from the triple of this sum the double of the preceding sum is subtracted, it will remain

\begin{equation*}
S.x^3=\frac{3}{2}xS.x^2-\frac{1}{4}x^3(x+1).
\end{equation*}
If one puts $n=5$, it will be

\begin{equation*}
S.x^5=x^5+x^5-5x^4S.x+10x^3S.x^2-10x^2S.x^3+5xS.x^4-S.x^5
\end{equation*}
or

\begin{equation*}
S.x^5=\frac{5}{2}xS.x^4-5x^2S.x^3+5.x^3S.x^2-\frac{5}{2}x^4S.x+\frac{1}{2}x^5(x+1)
\end{equation*}
and from $n=6$ it follows

\begin{equation*}
S.x^6=x^7+x^6-6x^5S.x+15x^4S.x^2-20x^3S.x^3+15x^2S.x^4-6xS.x^5+S.x^6
\end{equation*}
or

\begin{equation*}
S.x^5=\frac{5}{2}xS.x^4-\frac{10}{3}x^2S.x^3+\frac{5}{2}x^3S.x^2-x^4S.x+\frac{1}{6}x^5(x+1).
\end{equation*}

\paragraph*{§62}

From these results we therefore conclude, if $n=2m+1$, that in general it will be

\begin{equation*}
S.x^{2m+1}=\frac{2m+1}{2}xS.x^{2m}-\frac{(2m+1)2m}{ 2 \cdot 1 \cdot 2}x^2 S.x^{2m-1}
\end{equation*}
\begin{equation*}
+\frac{(2m+1)2m(2m-1)}{2 \cdot 1 \cdot 2 \cdot 3}x^3 S.x^{2m-2}-\cdots-\frac{2m+1}{2}x^{2m}S.x+\frac{1}{2}x^{2m+1}(x+1).
\end{equation*}
But if it was $n=2m+2$, since the terms $S.x^{2m+2}$ cancel each other, one will find

\begin{equation*}
S.x^{2m+1}=\frac{2m+1}{2}xS.x^{2m}-\frac{(2m+1)2m}{2 \cdot 3}x^2S.x^{2m-1}
\end{equation*}
\begin{equation*}
+\frac{(2m-1)2m(2m+1)}{2 \cdot 3 \cdot 4}x^3 S.x^{2m-2}-\cdots-x^{2m}S.x+\frac{1}{2m+2}x^{2m+1}(x+1).
\end{equation*}
Therefore, the sum of the odd powers can  be determined from the sums of the lower powers in two ways and from the various combinations of these formulas infinitely many others can be formed.

\paragraph*{§63}

But the sum of the odd powers can be determined a lot more easily from the preceding ones and for this it suffices to know only the sum of the preceding even power. For, from the sums of powers exhibited above [§ 62 of the first part] it is known that the number of terms constituting the sums is only increased in the sum of the even powers, such that the sum of the odd powers consists of as many terms as the sum of the preceding even power. So, if the sum of the even power $x^{2n}$ is

\begin{equation*}
S.x^{2n}=\alpha x^{2n+1}+\beta x^{2n}+\gamma x^{2n-1}-\delta x^{2n-3}+\varepsilon x^{2n-5}-\text{etc.}
\end{equation*}
(for, we saw that after the third term each second term is missing and at the same time the signs alternate), hence the sum of the following power $x^{2n+1}$ will be found, if each term is  respectively multiplied by these numbers

\begin{equation*}
\frac{2n+1}{2n+2}x, \quad \frac{2n+1}{2n+1}x, \quad \frac{2n+1}{2n}x, \quad \frac{2n+1}{2n-1}x, \quad \frac{2n+1}{2n-2}x \quad \text{etc.}
\end{equation*}
not omitting the missing terms; therefore, it will be

\begin{equation*}
S.x^{2n+1}=\frac{2n+1}{2n+2}\alpha x^{2n+2}+\frac{2n+1}{2n+1}\beta x^{2n+1}+\frac{2n+1}{2n}\gamma x^{2n}-\frac{2n+1}{2n-1}\delta x^{2n-2}
\end{equation*}
\begin{equation*}
+\frac{2n+1}{2n-4}\varepsilon x^{2n-4}-\frac{2n+1}{2n-6}\zeta x^{2n-6}+\text{etc.}
\end{equation*}
Therefore, if the sum of the power $x^{2n}$ is known, from it the sum of the following power $x^{2n+1}$ can be constructed conveniently.

\paragraph*{§64}

This investigation of the following sums is also extended to the even powers; but since the sums of these receive a new term, this term is not found by this method; nevertheless, it can always be found from the nature of the series, from which it is clear, if one puts $x=1$, that the sum has also to become $=1$. But vice versa from a sum of certain power one will always be able to find the sum of the preceding power. For, if it was

\begin{equation*}
S.x^n= \alpha x^{n+1}+\beta x^n+\gamma x^{n-1}-\delta x^{n-3} +\varepsilon x^{n-5}-\zeta x^{n-7}+\text{etc.},
\end{equation*}
for the preceding power it will be 

\begin{equation*}
S. x^{n-1}=\frac{n+1}{n}\alpha x^n+\frac{n}{n}\beta x^{n-1}+\frac{n-1}{n}\gamma x^{n-2}-\frac{n-3}{n}\delta x^{n-4}+\text{etc.}
\end{equation*}
and hence one can go backwards arbitrarily far. But it is to be noted that  always  $\alpha=\frac{1}{2}$ and $\beta=\frac{1}{2}$ as is it already clear from the formulas given above.

\paragraph*{§65}

The attentive reader will immediately see that the sum of $x^{n-1}$  results, if the sum of the powers $x^n$ is differentiated and its differential is divided by $ndx$; and it will be $d.S.x^n=ndx \cdot S.x^{n-1}$, and since $d.x^n=nx^{n-1}dx$, it will be

\begin{equation*}
d.S.x^n=S.nx^{n-1}dx=S.d.x^n;
\end{equation*}
hence the differential of the sum is equal to the sum of the differentials; so, if the general term of a certain series was $=y$ and $S.y$ was its summatory term, it will also be $S.dy=d.S.y$ in general, i.e. the sum of all differentials is equal to the differential of the sum of the terms. The validity of this equality is easily seen from those results derived above on the differentiation of series. For, because 

\begin{equation*}
S.x^n=x^n+(x-1)^n+(x-2)^n+(x-3)^n+(x-4)^n+\text{etc.},
\end{equation*}
it will be

\begin{equation*}
\frac{d.S.x^n}{ndx}=x^{n-1}+(x-1)^{n-1}+(x-2)^{n-2}+(x-3)^{n-1}+\text{etc.}=S.x^{n-1},
\end{equation*}
which proof extends to all other series.

\paragraph*{§66}

But let us return to our starting point, i.e. to the differences of functions, on which still several things are to be remarked. Because we saw, if $y$ was any function of $x$ and  one writes $x \pm \omega$ instead of $x$ everywhere, that the function $y$ will have the following value

\begin{equation*}
y \pm \frac{\omega dy}{1dx}+\frac{\omega^2 ddy}{1 \cdot 2 dx^2}\pm \frac{\omega^3 d^3 y}{1 \cdot 2 \cdot 3dx^3}+\frac{\omega^4 d^4 y}{1 \cdot 2 \cdot 3 \cdot 4 dx^4}\pm \frac{\omega^5 d^5 y}{1 \cdot 2 \cdot 3 \cdot 4 \cdot 5 dx^5}+\text{etc.},
\end{equation*}
this expression will be valid, no matter whether for $\omega$ any constant quantity or even a variable value depending on $x$ is taken. For, having found the values of $\frac{dy}{dx}$, $\frac{ddy}{dx^2}$, $\frac{d^3y}{dx^3}$ etc. by  differentiation, the variability in the factors $\omega$, $\omega^2$, $\omega^3$ etc. is not considered    and hence it does not matter, whether $\omega$ denotes a constant quantity or a variable quantity depending on $x$.

\paragraph*{§67}

Therefore, let us put that $\omega =x$ and in the function $y$ the value  $x-x=0$ is written instead of $x$. Therefore, if in any function of $x$  one writes $0$ instead of $x$ everywhere, the value of the function will be

\begin{equation*}
y-\frac{xdy}{1dx}+\frac{x^2 ddy}{1 \cdot 2 dx^2}-\frac{x^3d^3}{1 \cdot 2 \cdot 3 dx^3}+\frac{x^4 d^4 y}{1 \cdot 2 \cdot 3 \cdot 4 dx^4}-\text{etc.}
\end{equation*} 
Therefore, this expression always indicates the value the function $y$ has for $x=0$, the validity of which statement will be illustrated by the following examples.

\subsection*{Example 1}

\textit{Let $y=xx+ax+ab$, whose value, if one puts $x=0$, is in question, which value is known to be $=ab$, of course.}\\[2mm]
Because $y=xx+ax+ab$, it will be

\begin{equation*}
\frac{dy}{1dx}=2x+a, \quad \frac{ddy}{1 \cdot 2 dx^2}=1
\end{equation*}
and hence the value in question results as

\begin{equation*}
=xx+ax+ab-x(2x+a)+xx\cdot 1=ab.
\end{equation*}

\subsection*{Example 2}

\textit{Let $y=x^3-2x+3$, whose value, having put $x=0$, is in question, which value is known to be $=3$.}\\[2mm]
Because $y=x^3-2x+3$, it will be

\begin{equation*}
\frac{dy}{dx}=3xx-2, \quad \frac{ddy}{1 \cdot 2 dx^2}=3x, \quad \frac{d^3 y}{1 \cdot 2 \cdot 3 dx^3}=1;
\end{equation*}
the value in question will be found to be

\begin{equation*}
=x^3-2x+3-x(3xx-2)+xx \cdot 3x-x^3 \cdot 1=3.
\end{equation*}

\subsection*{Example 3}

\textit{Let $y=\frac{x}{1-x}$, whose value, having put $x=0$, is in question, which is known to be $=0$.}\\[2mm]
Because $y=\frac{x}{1-x}$, it will be

\begin{equation*}
\frac{dy}{dx}=\frac{1}{(1-x)^2}, \quad \frac{ddy}{1 \cdot 2 dx^2}=\frac{1}{(1-x)^3}, \quad \frac{d^3y}{1 \cdot 2 \cdot 3 dx^3}=\frac{1}{(1-x)^4}\quad \text{etc.}
\end{equation*}
Hence the value in question will be

\begin{equation*}
=\frac{x}{1-x}-\frac{x}{(1-x)^2}+\frac{xx}{(1-x)^3}-\frac{x^3}{(1-x)^4}+\frac{x^4}{(1-x)^5}-\text{etc.}
\end{equation*}
and therefore the value of this series is $=0$. This is also obvious from the fact that this series without the first term, i.e. $\frac{x}{(1-x)^2}-\frac{xx}{(1-x)^3}+\frac{x^3}{(1-x)^4}-\text{etc.}$ is a geometric series and its sum is $=\frac{x}{(1-x)^2+x(1-x)}=\frac{x}{1-x}$, whence the value found will be

\begin{equation*}
=\frac{x}{1-x}-\frac{x}{1-x}=0.
\end{equation*}

\subsection*{Example 4}

\textit{Let $y=e^x$, while $e$ denotes the number whose hyperbolic logarithm is $1$, and let the value of $y$ be in question, if one puts $x=0$, which value is known to be $=1$.}\\[2mm]
Because $y=e^x$, it will be

\begin{equation*}
\frac{dy}{dx}=e^x,\quad \frac{ddy}{dx^2}=e^x \quad \text{etc.}
\end{equation*}
and hence the value in question will be

\begin{equation*}
=e^x-\frac{e^x}{1}+\frac{e^x xx}{1 \cdot 2}-\frac{e^x x^3}{1 \cdot 2 \cdot 3}+\frac{e^x x^4}{1 \cdot 2 \cdot 3 \cdot 4}-\text{etc.}
\end{equation*}
\begin{equation*}
=e^x \left(1-\frac{x}{1}+\frac{xx}{1 \cdot 2}-\frac{x^3}{1 \cdot 2 \cdot 3}+\frac{x^4}{1 \cdot 2 \cdot 3 \cdot 4}-\text{etc.}\right).
\end{equation*}
But above we saw that the series

\begin{equation*}
1-\frac{x}{1}+\frac{xx}{1 \cdot 2}-\frac{x^3}{1 \cdot 2 \cdot 3}+\text{etc.}
\end{equation*}
expresses the value $e^{-x}$; therefore, the value in question will be $e^x \cdot e^{-x}=\frac{e^x}{e^x}=1$, of course.

\subsection*{Example 5}

\textit{Let $y=\sin x$ and having put $x=0$ it is obvious that it will be $y=0$, what also the general formula will indicate.}\\[2mm]
For, if $y=\sin x$, it will be

\begin{equation*}
\frac{dy}{dx}=\cos x, \quad \frac{ddy}{dx^2}=-\sin x, \quad \frac{d^3 y}{dx^3}=-\cos x, \quad \frac{d^4 y}{dx^4}=\sin x \quad \text{etc.}
\end{equation*}
Having put $x=0$, the value of $y$ will be

\begin{equation*}
\sin x -\frac{x}{1}\cos x- \frac{xx}{1 \cdot 2}\sin x+\frac{x^3}{1 \cdot 2 \cdot 3}\cos x+\frac{x^4}{1 \cdot 2 \cdot 3 \cdot 4}\sin x-\text{etc.}
\end{equation*}
which is

\begin{equation*}
=\sin x \left(1-\frac{xx}{1 \cdot 2}+\frac{x^4}{1 \cdot 2 \cdot 3 \cdot 4}-\frac{x^6}{1 \cdot 2 \cdot 3 \cdots 6}+\text{etc.}\right)
\end{equation*}
\begin{equation*}
-\cos x \left(\frac{x}{1}-\frac{x^3}{1 \cdot 2 \cdot 3}+\frac{x^5}{1 \cdot 2 \cdot 3 \cdot 4 \cdot 5}-\frac{x^7}{1 \cdot 2 \cdot 3 \cdots 7}+\text{etc.}\right)
\end{equation*}
But the upper of these series expresses $\cos x$, the lower expresses $\sin x$, whence the value in question will be

\begin{equation*}
= \sin x \cos x - \cos x \cdot \sin x=0.
\end{equation*}

\paragraph*{§68}

Hence, vice versa we can conclude, if $y$ was a function of $x$  vanishing for $x=0$, that it will be

\begin{equation*}
y- \frac{xdy}{1 dx}+\frac{xx ddy}{1 \cdot 2 dx^2}-\frac{x^3 d^3 y}{1 \cdot 2 \cdot 3 dx^3}+\frac{x^4 d^4 y}{1 \cdot 2 \cdot 3 \cdot 4 dx^4}-\text{etc.}=0.
\end{equation*}
Hence this is the general equation of completely all functions of $x$, which, if $x=0$, also become zero. And therefore, this equation is of such a nature, that, no matter which function of $x$, as long as it vanishes for vanishing $x$, is substituted for $y$, it is always satisfied. But if  $y$ was a function of $x$ which, having put $x=0$, has a given value $=A$, it will be

\begin{equation*}
=y-\frac{xdy}{1 dx}+\frac{x^2 ddy}{1\cdot 2dx^2}-\frac{x^3 d^3y}{1 \cdot 2 \cdot 3 dx^3}+\frac{x^4 d^4y}{1 \cdot 2 \cdot 3 \cdot 4 dx^4}-\text{etc.}=A,
\end{equation*}
which equation contains all functions of $x$  which go over into $A$ for $x=0$.

\paragraph*{§69}

If  one writes $2x$ or $x+x$ instead of $x$, any function of $x$, which we want to denote by $y$, will have this value

\begin{equation*}
y+\frac{xdy}{1 dx}+\frac{x^2 ddy}{1 \cdot 2 dx^2}+\frac{x^3 d^3y}{1 \cdot 2 \cdot 3 dx^3}+\frac{x^4 d^4 y}{1 \cdot 2 \cdot 3 \cdot 4 dx^4}+\text{etc.}
\end{equation*}
And if we write $nx$ instead of $x$, i.e. $x+(n-1)x$, the function $y$ will have the following value

\begin{equation*}
y+\frac{(n-1)xdy}{1 dx}+\frac{(n-1)^2xxddy}{1 \cdot 2 dx^2}+\frac{(n-1)^3 x^3 d^3 y}{1 \cdot 2 \cdot 3 dx^3}+\text{etc.}
\end{equation*}
but if in general we write $t$ for $x$, because of $t=x+t-x$, any function $y$ of $x$ will  be transformed into the following form

\begin{equation*}
y+\frac{(t-x)dy}{1dx}+\frac{(t-x)^2 ddy}{1 \cdot 2 dx^2}+\frac{(t-x)^3 d^3y}{1 \cdot 2 \cdot 3 dx^3}+\text{etc.}
\end{equation*} 
Therefore, if $v$ was such a function of $t$ as $y$ is of $x$, since $v$ results from $y$ by writing $t$ instead of $x$, it will be

\begin{equation*}
v=y+\frac{(t-x)dy}{1dx}+\frac{(t-x)^2ddy}{1 \cdot 2 dx^2}+\frac{(t-x)^3 d^3y}{1 \cdot 2 \cdot 3 dx^3}+\text{etc.},
\end{equation*}
the validity of which formula can be checked in any arbitrary example.

\subsection*{Example}

\textit{For, let $y=xx-x$;  having written $t$ instead of $x$,   it  will be obviously be $v=tt-t$, which same equation the expression  we found will also reveal.}\\[2mm]
For, because of $y=xx-x$, it will be

\begin{equation*}
\frac{dy}{dx}=2x-1 \quad \text{and} \quad \frac{ddy}{2dx^2}=1;
\end{equation*}
therefore, it will be

\begin{equation*}
v=xx-x+(t-x)(2x-1)+(t-x)^2
\end{equation*}
\begin{equation*}
=xx-x+2tx-2xx-t+x+tt-2tx+xx=tt-t.
\end{equation*}
Therefore, if $y$ was a function of $x$, which goes over into $A$ for $x=a$, because of $t=a$ and $v=A$, it will be

\begin{equation*}
A=y+\frac{(a-x)dy}{1dx}+\frac{(a-x)^2 ddy}{1 \cdot 2 dx^2}+\frac{(a-x)^3d^3 y}{1 \cdot 2 \cdot 3 dx^3}+\text{etc.}
\end{equation*}
and hence all functions of $x$, which go over into $A$ for $x=a$, satisfy this equation.


\title{On the Conversion of Functions into Series
\footnote{
Original title: "'De Conversione Functionum in Series"', first published as part of the book "`\textit{Institutiones calculi differentialis cum eius usu in analysi finitorum ac doctrina serierum}, 1755"', reprinted in \textit{Opera Omnia}: Series 1, Volume 10, pp. 276 - 308, Eneström-Number E212, translated by: Alexander Aycock for the "`Euler-Kreis Mainz"'}}
\author{Leonhard Euler}
\date{ }
\maketitle

\paragraph*{§70}

In the last chapter we have already partially shown the application the general expressions found there for the finite differences have for the investigation of series exhibiting the value of a certain function of $x$. For, if $y$ was a given function of $x$, the value it has for $x=0$, will be known; and if this value is put $=A$, it will be, as we found,

\begin{equation*}
y-\frac{xdy}{dx}+\frac{x^2 ddy}{1 \cdot 2 dx^2}-\frac{x^3 d^3 y}{1 \cdot 2 \cdot 3}+\frac{x^4 d^4 y}{1 \cdot 2 \cdot 3 \cdot 4 dx^4}-\text{etc.}=A.
\end{equation*}
Therefore,  we  not only have a, in most cases infinite, series, whose sum is equal to the constant quantity $A$, even though the variable quantity $x$ is contained in each term, but we will also be able to express the function $y$ by means of a series; for, it will be

\begin{equation*}
y= A+\frac{xdy}{dx}-\frac{xxddy}{1 \cdot 2 dx^2}+\frac{x^3 d^3 y}{1 \cdot 2 \cdot 3 dx^3}-\frac{x^4 dy^4}{1 \cdot 2 \cdot 3 \cdot 4 dx^4}+\text{etc.},
\end{equation*}
several examples of which  were already mentioned. 

\paragraph*{§71}

But for this investigation to extend further, let us put that the function $y$ goes over into $z$, if  one writes $x+ \omega$ instead of $x$ everywhere, so that $z$ is such a function of $x+\omega$ as $y$ is of $x$, and we showed [§ 48] that it will be

\begin{equation*}
z=y+\frac{\omega dy}{dx}+\frac{\omega^2 ddy}{1 \cdot 2 dx^2}+\frac{\omega^3 d^3 y}{1 \cdot 2 \cdot 3 dx^3}+\frac{\omega^4 d^4 y}{1 \cdot 2 \cdot 3 \cdot 4 dx^4}+\text{etc.}
\end{equation*}
Therefore, since each term of this series can be found by iterated differentiation of $y$, having put $dx$ to be constant, and at the same time the value of $z$ can actually be exhibited by substituting $x+\omega$ for $x$, this way one will always obtain a series equal to the value of $z$, which, if $\omega$ was a very small quantity, converges rapidly and,  by taking many terms, will yield an approximately true value of $z$. Hence the use of this formula  will be huge for approximations.

\paragraph*{§72}

Therefore,  to proceed in order in the demonstration of the vast applicability of this  formula, let us at first substitute algebraic functions of $x$ for $y$. Let $y=x^n$ and, if  one writes $x+\omega$ instead of $x$, it will be $z=(x+\omega)^n$. Therefore, since

\begin{equation*}
\frac{dy}{dx}=nx^{n-1}, \quad \frac{ddy}{dx^2}=n(n-1)x^{n-2}, \quad \frac{d^3y}{dx^3}=n(n-1)(n-2)x^{n-3},
\end{equation*} 
\begin{equation*}
\frac{d^4y}{dx^4}=n(n-1)(n-2)(n-3)x^{n-4} \quad \text{etc.},
\end{equation*}
having substituted these values, it will be

\begin{equation*}
(x+ \omega)^n =x^n+\frac{n}{1}x^{n-1}\omega +\frac{n(n-1)}{1 \cdot 2}x^{n-2} \omega^2+\frac{n(n-1)(n-2)}{1 \cdot 2 \cdot 3}x^{n-3} \omega^3+\text{etc.},
\end{equation*}
which is the well-known Newtonian expression, by which the power of the binomial $(x+\omega)^n$ is converted into a series. And the number of terms of this series is always finite, if $n$ was a positive integer.

\paragraph*{§73}

Hence we will also be able to find the progression  expressing the value of the power of the binomial in such a way that it terminates, if the  exponent of the power was a negative number. For, let us set

\begin{equation*}
\omega =\frac{-ux}{x+u};
\end{equation*}
it will be

\begin{equation*}
z=(x+\omega)^n=\left(\frac{xx}{x+u}\right)^n
\end{equation*}
and hence one will have

\begin{equation*}
\frac{x^{2n}}{(x+u)^n}=x^n-\frac{nx^nu}{1(x+u)}+\frac{n(n-1)x^n u^2}{1 \cdot 2 (x+u)^2}-\frac{n(n-1)(n-2)x^nu^3}{1 \cdot 2 \cdot 3 (x+u)^3}+\text{etc.}
\end{equation*}
 Divide by $x^{2n}$ everywhere and it will be
 
\begin{equation*}
(x+u)^{-n}=x^{-n}-\frac{nx^{-n}u}{1(x+u)}+\frac{n(n-1)x^{-n}u^2}{1 \cdot 2 (x+u)^2}-\frac{n(n-1)(n-2)x^{-n}u^3}{1 \cdot 2 \cdot 3 (x+u)^3}+\text{etc.} 
\end{equation*}
Now put $-n=m$ and this equation will result

\begin{equation*}
(x+u)^m=x^m+\frac{mx^mu}{1(x+u)}+\frac{m(m+1)x^mu^2}{1 \cdot 2 (x+u)^2}+\frac{m(m+1)(m+2)x^mu^3}{1 \cdot 2 \cdot 3 (x+u)^3}+\text{etc.},
\end{equation*}
which series, if $m$ is a negative integer, will consist of a finite number of terms. Therefore, this series is equal to the one found first, if  one writes $u$ and $m$ instead of $\omega$ and $n$; for, hence it will be

\begin{equation*}
(x+u)^m=x^m+\frac{mx^{m-1}u}{1}+\frac{m(m-1)x^{m-2}u^2}{1 \cdot 2}+\frac{m(m-1)(m-2) x^{m-3}u^3}{1 \cdot 2 \cdot 3}+\text{etc.}
\end{equation*}

\paragraph*{§74}

This same series can also be deduced from the expression given at the beginning of § 70. For, because, if $y$ goes over into $A$ for $x=0$,

\begin{equation*}
y-\frac{xdy}{dx}+\frac{xxddy}{1 \cdot 2 dx^2}-\frac{x^3 d^3y}{1 \cdot 2 \cdot 3 dx^3}+\frac{x^4 d^4y}{1 \cdot 2 \cdot 3 \cdot 4 dx^4}-\text{etc.}=A,
\end{equation*}
put $y=(x+a)^n$ and it will be $A=a^n$ and, because of 

\begin{equation*}
\frac{dy}{dx}=n(x+a)^{n-1}, \quad \frac{ddy}{dx^2}=n(n-1)(x+a)^{n-2},
\end{equation*}
\begin{equation*}
\frac{d^3y}{dx^3}=n(n-1)(n-2)(x+a)^{n-3}\quad \text{etc.},
\end{equation*}
it will be

\begin{equation*}
(x+a)^n-\frac{n}{1}x(x+a)^{n-1}+\frac{n(n-1)}{1 \cdot 2}x^2(x+a)^{n-2}-\text{etc.}=a^n;
\end{equation*}
divide by $a^n(x+a)^n$ and this equation will result

\begin{equation*}
(x+a)^{-n}=a^{-n}-\frac{na^{-n}x}{1(x+a)}+\frac{n(n-1)a^{-n}x^2}{1 \cdot 2 (x+a)^2}-\text{etc.},
\end{equation*}
which, having  substituted $u$, $x$ and $-m$ for $x$, $a$ and $n$ respectively, this expression will turn out to be the series found before.

\paragraph*{§75}

If one substitutes fractional numbers for $m$, both series will continue forever; nevertheless, if $u$  was a very small quantity with respect to $x$, the series will converge to the true value rapidly. Therefore, let $m=\frac{\mu}{\nu}$ and $x=a^{\nu}$;  from the series found first it will be

\begin{equation*}
(a^{\nu}+u)^{\frac{\mu}{\nu}}=a^{\mu} \left(1+\frac{\mu }{\nu}\cdot \frac{u}{a^{\nu}}+\frac{\mu(\mu- \nu)}{\nu \cdot 2 \nu}\cdot \frac{uu}{a^{2 \nu}}+\frac{\mu(\mu- \nu)(\mu-2 \nu)}{\nu \cdot 2 \nu \cdot 3 \nu} \cdot \frac{u^3}{a^{3 \nu}}+\text{etc.}\right).
\end{equation*}
But the series found later will give

\begin{equation*}
(a^ {\nu}+u)^{\frac{\mu}{\nu}}=a^{\mu}\left(1+\frac{\mu u}{\nu (a^{\nu}+u)}+\frac{\mu(\mu +\nu)u^2}{\nu \cdot 2 \nu (a^{\nu}+u)^2}+\frac{\mu(\mu +\nu)(\mu +2 \nu)u^3}{\nu \cdot 2 \nu \cdot 3 \nu (a^{\nu}+u)^3}+\text{etc.}\right).
\end{equation*}
But this last series converges more rapidly than the first, since its terms also decrease, if it was $u>a^{\nu}$, in which case the first series even diverges.\\[2mm]
Therefore, let $\mu=1$, $\nu=2$, it will be

\begin{equation*}
\sqrt{a^2+u}=a \left(1+\frac{1 u}{2(a^2+u)}+\frac{1 \cdot 3 u^2}{2 \cdot 4(a^2+u)^2}+\frac{1 \cdot 3 \cdot 5 u^3}{2 \cdot 4 \cdot 6 (a^2+u)^3+\text{etc.}}\right).
\end{equation*}
In like manner, by substituting the numbers $3$, $4$, $5$ etc. for $\nu$, while  still  $\mu=1$, it will be

\begin{alignat*}{9}
&\sqrt[3]{a^3+u}&&=a \bigg(1&&+\frac{1u}{3(a^3+u)}&&+~\frac{1 \cdot 4 u^2}{3 \cdot 6(a^3+u)^2}&&+~~\frac{1 \cdot 4 \cdot 7 u^3}{3 \cdot 6 \cdot 9 (a^3+u)^3}&&+\text{etc.}\bigg) \\[3mm]
&\sqrt[4]{a^4+u}&&=a \bigg(1&&+\frac{1u}{4(a^4+u)}&&+~\frac{1 \cdot 5 u^2}{4 \cdot 8(a^4+u)^2}&&+~\frac{1 \cdot 5 \cdot 9 u^3}{4 \cdot 8 \cdot 12 (a^4+u)^3}&&+\text{etc.}\bigg) \\[3mm]
&\sqrt[5]{a^5+u}&&=a \bigg(1&&+\frac{1u}{5(a^5+u)}&&+\frac{1 \cdot 6 u^2}{5 \cdot 10(a^5+u)^2}&&+\frac{1 \cdot 6 \cdot 11 u^3}{5 \cdot 10 \cdot 15 (a^5+u)^3}&&+\text{etc.}\bigg) \\
& && && &&\text{etc.}
\end{alignat*}

\paragraph*{§76}

From these formulas one can therefore easily find the root of  a certain power of any number. For, having propounded the number $c$, find the power closest to it, either larger or smaller; in the first case $u$ will become a negative number, in the second a positive number. But if the resulting series does not converge fast enough, multiply the number $c$ by any power, say by $f^{\nu}$, if the root of the power $\nu$ has to be extracted, and find the root of the number $f^{\nu}c$, which divided by $f$ will give the root in question of the number $c$.  The greater the number $f$ is assumed, the more the series will converge and  that especially, if a similar power $a^{\nu}$ does not deviate much from $f^{\nu}c$.

\subsection*{Example 1}

\textit{Let the square root of the number $2$ be in question.}\\[2mm]
If without any further preparation one puts $a=1$ and $u=1$, it will be

\begin{equation*}
\sqrt{2}=1+\frac{1}{2 \cdot 2}+\frac{1 \cdot 3}{2 \cdot 4 \cdot 2^2}+\frac{1 \cdot 3 \cdot 5}{2 \cdot 4 \cdot 6 \cdot 2^3}+\text{etc.};
\end{equation*}
even though this series already converges rapidly, it is nevertheless preferable to multiply the number $2$ by a square, as $25$, before, so that the product $50$ deviates from another square $49$ as less as possible. Therefore, find the square root of $50$, which divided by $5$ will give $\sqrt{2}$. But then it will be $a=7$ and $u=1$, whence it will be

\begin{equation*}
\sqrt{50}= 5\sqrt{2}=7 \left(1+\frac{1}{2 \cdot 50}+\frac{1 \cdot 3}{2 \cdot 4 \cdot 50^2}+\frac{1 \cdot 3 \cdot 5}{2 \cdot 4 \cdot 6 \cdot 50^3}+\text{etc.}\right)
\end{equation*}
or

\begin{equation*}
\sqrt{2}=\frac{7}{5}\left(1+\frac{1}{100}+\frac{1 \cdot 3}{100 \cdot 200}+\frac{1 \cdot 3 \cdot 5}{100 \cdot 200 \cdot 300}+\text{etc.}\right),
\end{equation*}
which is most appropriate for the calculation in decimal numbers. For, it will be
 
 \begin{center}
 \begin{tabular}{rrr}
$\dfrac{7}{5}$ & $=$ & $1.4000000000000$ \\[3mm]
$\dfrac{7}{5} \cdot \dfrac{1}{100}$ & $=$ & $140000000000$ \\ [3mm]
$\dfrac{7}{5} \cdot \dfrac{1}{100} \cdot \dfrac{3}{200}$ & $=$ & $2100000000$ \\ [3mm]
$\dfrac{7}{5} \cdot \dfrac{1}{100} \cdot \dfrac{3}{200} \cdot \dfrac{5}{300}$ & $=$ & $35000000$ \\ [3mm]
the preceding by $\dfrac{7}{400}$ & $=$ & $612500$ \\  [3mm]
the preceding by $\dfrac{9}{500}$ & $=$ & $11025$ \\  [3mm]
the preceding by $\dfrac{11}{600}$ & $=$ & $202$ \\  [3mm]
the preceding by $\dfrac{13}{400}$ & $=$ & $3$ \\ [3mm] \cline{3-3} 
Therefore, $\sqrt{2}$ &$=$ & $1.4142135623730$. 
 \end{tabular}
 \end{center}

\subsection*{Example 2}

\textit{Let the cube root of $3$ be in question.}\\[2mm]
Multiply $3$ by the cube $8$ and find the cube root of $24$; for, it will be $\sqrt[3]{24}=2\sqrt[3]{3}$. Therefore, put $a=3$ and $u=-3$ and it will be 

\begin{equation*}
\sqrt[3]{24}=3\left(1- \frac{1 \cdot 3}{3 \cdot 24}+\frac{1 \cdot 4 \cdot 3^2}{3 \cdot 6 \cdot 24^2}-\frac{1 \cdot 4 \cdot 7 \cdot 3^2}{3 \cdot 6 \cdot 9 \cdot 24^3}+\text{etc.}\right)
\end{equation*} 
and

\begin{equation*}
\sqrt[3]{3}=\frac{3}{2}\left(1-\frac{1}{3 \cdot 8}+\frac{1 \cdot 4}{3 \cdot 6 \cdot 8^2}-\frac{1 \cdot 4 \cdot 7}{3 \cdot 6 \cdot 9 \cdot 8^3}+\text{etc.}\right)
\end{equation*}
or

\begin{equation*}
\sqrt[3]{3}=\frac{3}{2}\left(1-\frac{1}{24}+\frac{1}{24}\cdot \frac{4}{48}-\frac{1}{24}\cdot \frac{4}{48}\cdot \frac{7}{72}+\text{etc.}\right),
\end{equation*}
which series already converges rapidly, since every term is more than eight times smaller than the preceding. But if $3$ is multiplied by the cube $729$, it will be $2187$ and $\sqrt[3]{2187}=\sqrt{13^3-10}=9\sqrt[3]{3}$. Therefore, because of $a=13$ and $u=-10$, it will be

\begin{equation*}
\sqrt[3]{3}=\frac{13}{9}\left(1-\frac{1 \cdot 10}{3 \cdot 2187}+\frac{1 \cdot 4 \cdot 10^2}{3 \cdot 6 \cdot 2187^2}-\frac{1 \cdot 4 \cdot 7 \cdot 10^3}{3 \cdot 6 \cdot 9 \cdot 2187^3}+\text{etc.}\right),
\end{equation*}
 every term of which series is more than two hundred times smaller than the preceding.

\paragraph*{§77}

The expansion of the power of the binomial extends so far that all algebraic functions are comprehended by it. For, if for the sake of an example the value of this function $\sqrt{a+2bx+cxx}$ expressed as  a series is in question, this question can be answered by means of the preceding formulas  considering two terms as one. Furthermore, this expansion can be done by means of the expression given first; for, if one puts $\sqrt{a+2bx+cxx}=y$, since, having put $x=0$,  $y=a$, it will be $A=\sqrt{a}$, and since the differentials of $y$ will be as follows

\begin{equation*}
\frac{dy}{dx}=\frac{b+cx}{\sqrt{a+2bx+cxx}}, \quad \frac{ddy}{dx^2}=\frac{ac-bb}{(a+2bx+cxx)^{\frac{3}{2}}}, \quad \frac{d^3y}{dx^3}=\frac{3(bb-ac)(b+cx)}{(a+2bx+cxx)^{\frac{5}{2}}},
\end{equation*}
\begin{equation*}
\frac{d^4y}{dx^4}=\frac{3(bb-ac)(ac-5bb-8bcx-4ccxx)}{(a+2bx+cxx)^{\frac{7}{2}}} \quad \text{etc.},
\end{equation*}
from these one will therefore obtain

\begin{equation*}
\sqrt{a+2bx+cxx}-\frac{(b+cx)x}{\sqrt{a+2bx+cxx}}-\frac{(bb-ac)xx}{2(a+2bx+cxx)^{\frac{3}{2}}}-\frac{(bb-ac)(b+cx)x^3}{2(a+2bx+cxx)^{\frac{5}{2}}}
\end{equation*}
\begin{equation*}
-\frac{(bb-ac)(5bb-ac+8bcx+4ccxx)x^4}{8(a+2bx+cxx)^{\frac{7}{2}}}-\text{etc.}=\sqrt{a}
\end{equation*}
Therefore, if one multiplies by $\sqrt{a+2bx+cxx}$ everywhere, the series will be rational and it will be

\begin{equation*}
\sqrt{a(a+2bx+cxx)}=a+2bx+cxx-(b+cx)x-\frac{(bb-ac)xx}{2(a+2bx+cxx)}
\end{equation*}
\begin{equation*}
-\frac{(bb-ac)(b+cx)x^3}{2(a+2bx+cxx)^2}-\frac{(bb-ac)(5bb-ac+8bcx+4ccxx)x^4}{8(a+2bx+cxx)^3}-\text{etc.}
\end{equation*}
or

\begin{equation*}
\sqrt{a+2bx+cxx}=\sqrt{a}+\frac{bx}{\sqrt{a}}-\frac{(bb-ac)xx}{2(a+2bx+cxx)\sqrt{a}}-\frac{(bb-ac)(b+cx)x^3}{2(a+2bx+cxx)^2\sqrt{a}}-\text{etc.}
\end{equation*}

\paragraph*{§78}

Hence let us go over to transcendental functions and let us  substitute them for $y$. Therefore, at first let $y= \log x$ and having put $x +\omega$ instead of $x$ it will be $z= \log(x+ \omega)$. But let these logarithms  have a ratio of $n:1$ to the hyperbolic logarithms; and for the hyperbolic logarithms it will be $n=1$ and for the tabulated logarithms it will be $n=0.4343944819032$. Hence the differentials of $y= \log x$ will be

\begin{equation*}
\frac{dy}{dx}=\frac{n}{x}, \quad \frac{ddy}{dx^2}=-\frac{n}{x^2}, \quad \frac{d^3y}{dx^3}=\frac{2n}{x^3} \quad \text{etc.},
\end{equation*}
from which one concludes

\begin{equation*}
\log (x+ \omega)= \log x +\frac{n \omega}{x}-\frac{n \omega^2}{2 x^2}+\frac{n \omega^3}{3 x^3}-\frac{n \omega^4}{4x^4}+\text{etc.}
\end{equation*}
In like manner, if $\omega$ is assumed to be negative, it will be\begin{equation*}
\log (x- \omega)=\log x -\frac{n \omega}{x}-\frac{n \omega^2}{2 x^2}-\frac{n \omega^3}{3 x^3}-\frac{n \omega^4}{4x^4}-\text{etc.}
\end{equation*}
Therefore, if this series is subtracted from the first, it will be

\begin{equation*}
\log \frac{x+ \omega}{x- \omega}=2n \left(\frac{\omega}{x}+\frac{\omega^3}{3x^3}+\frac{\omega^5}{5x^5}+\frac{\omega^7}{7x^7}+\text{etc.}\right).
\end{equation*}

\paragraph*{§79}

If in the series found first

\begin{equation*}
\log (x+ \omega)= \log x +\frac{n \omega}{x}-\frac{n \omega^2}{2x^2}+\frac{n \omega^3}{3 x^3}-\frac{n \omega^4}{4 x^4}+\text{etc.}
\end{equation*}
one puts

\begin{equation*}
\omega = \frac{xx}{u-x},
\end{equation*}
it will be $x+ \omega = \frac{ux}{u-x}$ and

\begin{equation*}
\log (x+\omega)= \log u+ \log x - \log (u-x)= \log x+\frac{nx}{u-x}-\frac{nxx}{2(u-x)^2}+\text{etc.}
\end{equation*}
and

\begin{equation*}
\log (u-x)= \log u -\frac{nx}{u-x}+\frac{nxx}{2(u-x)^2}-\frac{nx^3}{3(u-x)^3}+\text{etc.}
\end{equation*}
and, having taken a negative $x$,  one will have

\begin{equation*}
\log (u+x)= \log u +\frac{nx}{u+x}+\frac{nxx}{2(u+x)^2}+\frac{nx^3}{3(u+x)^3}+\frac{nx^4}{4(u+x)^4}+\text{etc.}
\end{equation*}
Therefore,  the logarithms can be found by means of these series in a convenient manner, if these series converge rapidly. Indeed, the following series, which are easily deduced from those already found, will be of this kind

\begin{alignat*}{9}
&\log (x+1) &&= \log x &&+n \bigg(\frac{1}{x}&&-\frac{1}{2xx}&&+\frac{1}{3x^3}&&-\frac{1}{4x^4}&&+\text{etc.}\bigg) \\[3mm]
&\log (x-1) &&= \log x &&-n \bigg(\frac{1}{x}&&+\frac{1}{2xx}&&+\frac{1}{3x^3}&&+\frac{1}{4x^4}&&+\text{etc.}\bigg);
\end{alignat*}
because these two series differ only in regard to the signs, if they are used for a calculation, given the logarithm of the number $x$, at the same time the logarithms of the two numbers $x-1$ and $x+1$ will be found. Furthermore, from the remaining series it will be

\begin{alignat*}{9}
& \log (x+1) &&= \log (x-1)+2n \left(\frac{1}{x}+\frac{1}{3x^3}+\frac{1}{5x^5}+\frac{1}{7x^7}+\text{etc.}\right)\\[3mm]
& \log (x-1) &&= \log x -n \left(\frac{1}{x-1}-\frac{1}{2(x-1)^2}+\frac{1}{3(x-1)^3}-\frac{1}{4(x-1)^4}+\text{etc.}\right)\\[3mm]
& \log (x+1) &&= \log x +n \left(\frac{1}{x+1}+\frac{1}{2(x+1)^2}+\frac{1}{3(x+1)^3}+\frac{1}{4(x+1)^4}+\text{etc.}\right).
\end{alignat*}

\paragraph*{§80}

Therefore, given the logarithm of the number $x$, the logarithms of the contiguous numbers $x+1$ and $x-1$ can easily be found; given the logarithm of the number $x-1$, even the logarithm of the number greater by two units and vice versa will be found. Although this was shown in much detail in the \textit{Introductio}, we will nevertheless add certain examples here.

\subsection*{Example 1}

\textit{Given hyperbolic logarithm of the number $10$, which is $2.3025850919940$, to find the hyperbolic logarithms of the numbers $11$ and $9$.}\\[2mm]
Since this question concerns hyperbolic logarithms, it will be $n=1$ and hence one will have these series

\begin{alignat*}{9}
& \log 11 &&= \log 10 &&+\frac{1}{10}&&-\frac{1}{2 \cdot 10^2}&&+\frac{1}{3 \cdot 10^3}&&-\frac{1}{4 \cdot 10^4}&&+\frac{1}{5 \cdot 10^5}&&-\text{etc.} \\[3mm]
& \log 9 &&= \log 10 &&+\frac{1}{10}&&+\frac{1}{2 \cdot 10^2}&&+\frac{1}{3 \cdot 10^3}&&+\frac{1}{4 \cdot 10^4}&&+\frac{1}{5 \cdot 10^5}&&-\text{etc.} \\
\end{alignat*}
To find the sums of these series, collect the even and odd terms separately and it will be

\begin{center}
\begin{tabular}{ccc|ccc}
$\dfrac{1}{10}$ & $=$ & $0.1000000000000$ & $\dfrac{1}{2 \cdot 10^2}$ & $=$ & $0.0050000000000$ \\[2mm]
$\dfrac{1}{3 \cdot 10^3}$ & $=$ & $0.0003333333333$ & $\dfrac{1}{4 \cdot 10^4}$ & $=$ & $0.0000250000000$ \\[2mm]
$\dfrac{1}{5 \cdot 10^5}$ & $=$ & $0.0000020000000$ & $\dfrac{1}{6 \cdot 10^6}$ & $=$ & $0.0000001666666$ \\[2mm]
$\dfrac{1}{7 \cdot 10^7}$ & $=$ & $0.0000000142857$ & $\dfrac{1}{8 \cdot 10^8}$ & $=$ & $0.0000000012500$ \\[2mm]
$\dfrac{1}{3 \cdot 10^9}$ & $=$ & $0.0000000001111$ & $\dfrac{1}{10 \cdot 10^{10}}$ & $=$ & $0.0000000000100$ \\[2mm]
$\dfrac{1}{11 \cdot 10^{11}}$ & $=$ & $0.0000000000009$ & $\dfrac{1}{12 \cdot 10^{12}}$ & $=$ & $0.0000000000001$ \\
 & & & & & \\\hline
Sum &$=$ & $0.1003353477310$ & Sum & $=$ & $0.05050251679267$  
\end{tabular}
\end{center}

\begin{center}
\begin{tabular}{llll}
The sum of both will be & & & $0.1053605156577$ \\
The difference of both will be & & & $0.0953101798043$ \\
Now  & $\log 10$ & $=$ & $2.3025850929940$ \\ \cline{4-4}
Therefore, it will be & $\log 11$ & $=$ & $2.397895272793$ \\
and & $\log 9$ & $=$ & $2.1972245773363$ \\ \cline{4-4} 
Hence further  & $\log 3$ & $=$ & $1.0986122886681$ \\
and & $\log 99$ & $=$ & $4.5951198501346$ \\
\end{tabular}
\end{center}

\subsection*{Example 2}

\textit{Using the hyperbolic logarithm of the number $99$ just found to find the logarithm of the number $101$.} \\[2mm]
For this, apply the series found above, i.e. 

\begin{equation*}
\log (x+1)=\log (x-1)+\frac{2}{x}+\frac{2}{3x^3}+\frac{2}{5x^5}+\frac{2}{7x^7}+\text{etc.},
\end{equation*}
in which one has to put $x=100$, and it will be

\begin{equation*}
\log 101 = \log 99 +\frac{2}{100}+\frac{2}{3 \cdot 100^3}+\frac{2}{5 \cdot 100^5}+\frac{2 }{7 \cdot 100^7}+\text{etc.},
\end{equation*}
the sum of which series  is calculated to be $=0.0200006667066$ from these four terms, which number added to $\log 99$ will give $\log 101 =4.6151205168412$.

\subsection*{Example 3}

\textit{Using the given tabulated logarithm of the number $10$, which is $=1$, to find the logarithm of the numbers $11$ and $9$.}\\[2mm]
Since here we look for  the common logarithm, it will be

\begin{equation*}
n=0.434244819032;
\end{equation*}
therefore, having put $x=10$, it will be

\begin{alignat*}{9}
&\log 11 &&= \log 10 &&+\frac{n}{10}&&-\frac{n}{2 \cdot 10^2}&&+\frac{n}{3 \cdot 10^3}&&-\frac{n}{4 \cdot 10^4}&&+\text{etc.} \\[3mm]
&\log 9 &&= \log 10 &&-\frac{n}{10}&&-\frac{n}{2 \cdot 10^2}&&-\frac{n}{3 \cdot 10^3}&&-\frac{n}{4 \cdot 10^4}&&-\text{etc.} \\
\end{alignat*}
Therefore, collect the even and odd terms separately

\begin{center}
\begin{tabular}{ccc|ccc}
$\dfrac{n}{10}$ & $=$ & $0.0434294481903$ & $\dfrac{n}{2 \cdot 10^2}$ & $=$ & $0.0021714724095$ \\[2mm]
$\dfrac{n}{3 \cdot 10^3}$ & $=$ & $0.0001447648273$ & $\dfrac{n}{4 \cdot 10^4}$ & $=$ & $0.0000108573620$ \\[2mm]
$\dfrac{n}{5 \cdot 10^5}$ & $=$ & $0.0000008685889$ & $\dfrac{n}{6 \cdot 10^6}$ & $=$ & $0.0000000723824$ \\[2mm]
$\dfrac{n}{7 \cdot 10^7}$ & $=$ & $0.0000000062042$ & $\dfrac{n}{8 \cdot 10^8}$ & $=$ & $0.0000000005428$ \\[2mm]
$\dfrac{n}{3 \cdot 10^9}$ & $=$ & $0.0000000000482$ & $\dfrac{n}{10 \cdot 10^{10}}$ & $=$ & $0.0000000000043$ \\[2mm]
$\dfrac{n}{11 \cdot 10^{11}}$ & $=$ & $0.0000000000004$ & $\dfrac{n}{12 \cdot 10^{12}}$ & $=$ & $0.0000000000000$ \\
 & & & & & \\\hline
Sum &$=$ & $0.0435750878593$ & Sum & $=$ & $0.0021824027010$  
\end{tabular}
\end{center}
\begin{tabular}{llll}
The aggregate of both is & &$=$ &$0.0457574905603$ \\ 
{Their difference is}    & &$=$ &$0.0413926851583$ \\
Therefore, because & $\log 10$ &$=$ &$ 1.0000000000000$ \\ 
it will be & $ \log 11$ & $=$ & $ 1.0413926851582$ \\ 
and & $\log 9$ &$=$ &  $0.9542425094397$ \\
hence & $\log 3$ & $=$ & $0.4771212547198$ \\ 
and & $\log 99$ &$=$ &$1.9956351945980$
\end{tabular}

\subsection*{Example 4}

\textit{Using the tabulated logarithm of the number $99$ found here to find the tabulated logarithm of the number $101$.}\\[2mm]
Here, by applying the same series we used in the second example, we will have

\begin{equation*}
\log 101 = \log 99 +2n \left(\frac{1}{100}+\frac{1}{3 \cdot 100^3}+\frac{1}{5 \cdot 100^5}+\text{etc.}\right),
\end{equation*}
the sum of which series, having substituted the corresponding value for $n$, will quickly be found to be

\begin{tabular}{llll}
  & &$=$ &$0.0086861791849$ \\
having added which to & $\log 99$ & $=$ & $1.9956351945980$ \\ \cline{4-4}
it   results & $\log 101$ & $=$ & $2.0043213737829$
\end{tabular}

\paragraph*{§81}

Now, in our general expression, let us  attribute  an exponential value to $y$ and let $y=a^x$; having substituted $x+\omega$ for $x$, it will be $z=a^{x+\omega}$, whose value, because of the differentials

\begin{equation*}
\frac{dy}{dx}=a^x \log a, \quad \frac{ddy}{dx^2}= a^x (\log a)^2, \quad \frac{d^3y}{dx^3}=a^x (\log a)^3 \quad \text{etc.},
\end{equation*} 
will be

\begin{equation*}
a^{x+\omega}=a^x \left(1+\frac{\omega \log a}{1}+\frac{\omega^2 (\log a)^2}{1 \cdot 2}+\frac{\omega^3 (\log a)^3}{1 \cdot 2 \cdot 3}+\text{etc.}\right);
\end{equation*}
if this equation is divided by $a^x$, the series expressing the values of an exponential quantity will result, which series we already found in the \textit{Introductio}, i.e.

\begin{equation*}
a^{\omega}=1+\frac{\omega \log a}{1}+\frac{\omega^2 (\log a)^2}{1 \cdot 2}+\frac{\omega^3 (\log a)^3}{1 \cdot 2 \cdot 3}+\frac{\omega^4 (\log a)^4}{1 \cdot 2 \cdot 3 \cdot 4}+\text{etc.}
\end{equation*}
In like manner, for negative $\omega$ it will be

\begin{equation*}
a^{-\omega}=1-\frac{\omega \log a}{1}+\frac{\omega^2 (\log a)^2}{1 \cdot 2}-\frac{\omega^3 (\log a)^3}{1 \cdot 2 \cdot 3}+\frac{\omega^4 (\log a)^4}{1 \cdot 2 \cdot 3 \cdot 4}-\text{etc.},
\end{equation*}
from whose combination these equations result

\begin{alignat*}{9}
&\frac{a^{\omega}+a^{-\omega}}{2}&&=1+\frac{\omega^2(\log a)^2}{1 \cdot 2}+\frac{\omega^4 (\log a)^4}{1 \cdot 2 \cdot 3 \cdot 4}+\frac{\omega^6 (\log a)^6}{1 \cdot 2 \cdot 3 \cdot 4 \cdot 5 \cdot 6}+\text{etc.} \\[3mm]
& \frac{a^{\omega}-a^{-\omega}}{2}&&=\frac{ \omega \log a}{1}+\frac{\omega^3(\log a)^3}{1 \cdot 2 \cdot 3}+\frac{\omega^5 (\log a)^5}{1 \cdot 2 \cdot 3 \cdot 4 \cdot 5}+\text{etc.},
\end{alignat*}
where it is to be noted that $\log a $ denotes the hyperbolic logarithm of the number $a$.

\paragraph*{§82}

By means of this formula, given a logarithm, one will be able to find the number corresponding to it. For, let any logarithm $u$ be propounded, for which the logarithm of the number $a$ is set $=1$. In the same base find the logarithm coming closest to $u$ and let $u=x+ \omega$, but  let the number corresponding to $x$ be $y=a^x$; the number corresponding to the logarithm $u=x+\omega$ will be $=a^{x+\omega}=z$ and it will be

\begin{equation*}
z=y \left(1+\frac{\omega \log a}{1}+\frac{\omega^2 (\log a)^2}{1 \cdot 2}+\frac{\omega^3 (\log a)^3}{1 \cdot 2 \cdot 3}+\frac{\omega^4 (\log a)^4}{1 \cdot 2 \cdot 3 \cdot 4}+\text{etc.}\right),
\end{equation*}
which series, because of the very small number $\omega$, converges rapidly; we will show its use   in the following example.

\subsection*{Example}

\textit{Let the number equal to this power of two $2^{2^{24}}$ be in question.}\\[2mm]
Since  $2^{24}=16777216$, it will be $2^{2^{24}}=2^{16777216}$ and by taking common logarithms the logarithm of this number will be $=16777216 \log 2$. But since

\begin{equation*}
\log 2 = 0.30102999566398119521373889,
\end{equation*}
the logarithm of the number in question will be

\begin{equation*}
5050445.259733675932039063,
\end{equation*}
whose characteristic indicates that the number in question has $5050446$ digits; since they cannot all be exhibited, it will suffice to have assigned the first few digits, which must be investigated from the mantissa

\begin{equation*}
,259733675932029063=u.
\end{equation*}
But from tables one concludes that the number, whose logarithm comes closest to this, will be $=1.818$, which number we want to  put $y$; its logarithm is

\begin{tabular}{lrll}
 & $x$ & $=$ & $0.259593878885948644$ \\
whence it will be & $\omega$ & $=$& $0.000139797046090419$ \\
Because now  & $a$ & $=$ & $10$ \\
it will be & $\log a$ & $=$& $2.3025850929940456840179914$ \\ \cline{2-4} 
and & $\omega \log a $ & $=$ & $0.000321894594372400$ \\ \cline{2-4}
Further, it will be & $y$ & $=$     & $1.818000000000000000$ \\[2mm]
& $\dfrac{\omega \log a}{1}y$ & $=$ & $0.000585204372569023$ \\[3mm]
& $\dfrac{\omega^2 (\log a)^2}{1 \cdot 2}y$ & $=$ & $0.000000094187062066$ \\ [3mm]
& $\dfrac{\omega^3 (\log a)^3}{1 \cdot 2 \cdot 3}y$ & $=$ & $0.000000000010106102$ \\ [3mm]
& $\dfrac{\omega^4 (\log a)^4}{1 \cdot 2 \cdot 3 \cdot 4}y$ & $=$ & $0.000000000000000813$ \\ [3mm] \cline{2-4}
& & & $1818585298569738004$
\end{tabular}
and these are the first few digits of the number in question, of which  all digits, except for the maybe the last, are correct.

\paragraph*{§83}

Let us consider quantities depending on the circle and let, as usual, the radius of the circle be $=1$ and let $y$ denote the arc of the circle, whose sine is $=x$, or let $y=\arcsin x$. Write $x+ \omega$ instead of $x$ and it will be $z=\arcsin (x+ \omega)$; to express this value, find the differentials of $y$ [§ 200 of the first part]

\begin{equation*}
\frac{dy}{dx}=\frac{1}{\sqrt{1-xx}}, \quad \frac{ddy}{dx^2}=\frac{x}{(1-xx)^{\frac{3}{2}}}, \quad \frac{d^3 y}{dx^3}=\frac{1+2xx}{(1-xx)^{\frac{5}{2}}}, \quad \frac{d^4y}{dx^4}=\frac{9x+6x^3}{(1-xx)^{\frac{7}{2}}},
\end{equation*}
\begin{equation*}
\frac{d^5y}{dx^5}=\frac{9+72x^2+24x^4}{(1-xx)^{\frac{9}{2}}}, \quad \frac{d^6y}{dx^6}=\frac{225x+600x^3+120x^5}{(1-xx)^{\frac{11}{2}}}
\end{equation*}
\begin{center}
etc.
\end{center}
Therefore,  from these one finds

\begin{equation*}
\arcsin (x+ \omega)=\arcsin x+ \frac{\omega}{\sqrt{1-xx}}+\frac{\omega^2 x}{2(1-xx)^{\frac{3}{2}}}+\frac{\omega^3 (1+2xx)}{6(1-xx)^{\frac{5}{2}}}
\end{equation*}
\begin{equation*}
+\frac{\omega^4(9x+6x^3)}{24(1-xx)^{\frac{7}{2}}}+\frac{\omega^5(9+72x^2+24x^4)}{120(1-xx)^{\frac{9}{2}}}+\text{etc.}
\end{equation*}

\paragraph*{§84}

Therefore, if the arc, whose sine is $=x$, was known, then, by means of this formula one will be able to find the arc, whose sine is $x+ \omega$, if $\omega$ was a very small quantity. But the series, whose sum must be added, will be expressed in parts of the radius, which will easily be reduced to an arc, as it will be understood from this example.

\subsection*{Example}

\textit{Let the arc of the circle, whose sine is $=\frac{1}{3}=0.3333333333$, be in question.}\\[2mm]
From tables find the arc, whose sine is approximately $\frac{1}{3}$ but a little bit smaller, which arc will be $19^{\circ}28^{\text{I}}$, whose sine is $=0.3332584$. Therefore, put $19^{\circ}28^{\text{I}}= \arcsin x=y$; it will be $x=0.3332584$ and $\omega = 0.000749$ and from tables $\sqrt{1-xx}=\cos y = 0.9428356$. Therefore, the arc in question $z$, whose sine   $=\frac{1}{3}$ is propounded, will be

\begin{equation*}
=19^{\circ}28^{\text{I}}+\frac{\omega}{\cos y}+\frac{\omega \omega \sin y}{2 \cos^3 y},
\end{equation*}
which expression already suffices; therefore, using logarithms in the calculation, it will be

\begin{tabular}{rrrrrr}
$\log \omega$ &$=$ & $5.8744818$ \\
$\log \cos y$ & $=$ & $9.9744359$ \\ \cline{3-3}
$ \log \dfrac{\omega}{\cos y}$ & $=$ & $5.9000459$ & $\dfrac{\omega}{\cos y}$ & $=$ & $0.0000794412$ \\[3mm]
$ \log \dfrac{\omega^2}{\cos^2 y}$ & $=$ & $1.8000918$ \\[3mm]
$ \log \dfrac{\omega^3}{\cos^3 y}$ & $=$ & $9.5483452$ \\ \cline{3-3} 
 &  & $1.3484370$ \\
 $ \log 2$ & $=$ & $0.3010300$ \\
 $ļog \dfrac{\omega^2 \sin y}{2 \cos^2 y}$ & $=$ & $1.0474070$ & $\dfrac{\omega^2 \sin y}{2 \cos^3 y}$ & $=$ & $0.0000000011$ \\ \cline{6-6}
 & & & sum & $=$ & $0.000794442,$     
\end{tabular}\\
which is the value of the arc to be added to $19^{\circ}28^{\text{I}}$, to express which  in minutes and seconds let us take its logarithm, which is
\begin{center}
\begin{tabular}{ll}
& $5.9000518$ \\
from which subtract & $4.6855749$ \\ \cline{2-2}
                    & $1.2144769$ \\ 
to which logarithm corresponds the number &$=16.38615$,
\end{tabular}
\end{center}
which is the number of minutes and seconds; but, expressing this fraction in thirds and quarters, the arc in question will be

\begin{equation*}
=19^{\circ} 28^{\text{I}}16^{\text{II}}23^{\text{III}}10^{\text{IV}}8^{\text{V}}14^{\text{VI}}.
\end{equation*}

\paragraph*{§85}

In like manner, the expression for the cosine will be found; for, having put $y= \arccos x$, since $dy=\frac{-dx}{\sqrt{1-xx}}$, the series we found  remains unchanged, as long as its signs are changed. Therefore, it will be

\begin{equation*}
\arccos (x+\omega)=\arccos x -\frac{\omega}{\sqrt{1-xx}}-\frac{\omega^2 x}{2(1-xx)^{\frac{3}{2}}}-\frac{\omega^3(1+2xx)}{6(1-xx)^{\frac{5}{2}}}
\end{equation*}
\begin{equation*}
-\frac{\omega^4(9x+6x^3)}{24(1-xx)^{\frac{7}{2}}}-\frac{\omega^5(9+72x^2+24x^4)}{120(1-xx)^{\frac{9}{2}}}-\text{etc.},
\end{equation*}
which series, as the preceding, will always converge rapidly, if, consulting tables of sines, angles close to the true one are chosen, such that in most cases the first term $\frac{\omega}{\sqrt{1-xx}}$ alone suffices. Nevertheless, if $x$ is approximately equal to $1$ or to a sinus totus, then, because of the very small denominators, the series will only converge slowly. Therefore, in these cases, in which $x$ does not deviate much from $1$, since the differences become very small, it will be more convenient to use only the usual method of interpolation.

\paragraph*{§86}

Therefore, let us also substitute the arc, whose tangent is given, for $y$ and let $y=\arctan x$ and $z=\arctan (x+ \omega)$, so that 

\begin{equation*}
z=y+\frac{\omega dy}{dx}+\frac{\omega^2 ddy}{2 dx^2}+\frac{\omega^3 d^3y}{6 dx^3}+\text{etc.}
\end{equation*}   
To investigate the terms, find each differential of $y$

\begin{equation*}
\frac{dy}{dx}=\frac{1}{1+xx}, \quad \frac{ddy}{dx^2}=\frac{-2x}{(1+xx)^2}, \quad \frac{d^3y}{dx^3}=\frac{-2+6xx}{(1+xx)^3}, \quad \frac{d^4y}{dx^4}=\frac{24x-24x^3}{(1+xx)^4},
\end{equation*}
\begin{equation*}
\frac{d^5y}{dx^5}=\frac{24.240x^2+120x^4}{(1+xx)^5}, \quad \frac{d^5y}{dx^5}=\frac{-720x+2400x^3-720x^5}{(1+xx)^5}
\end{equation*}
\begin{center}
etc.,
\end{center}
whence one concludes that

\begin{equation*}
\arctan(x+ \omega)=\arctan x
\end{equation*}
\begin{equation*}
+\frac{\omega}{1+xx}-\frac{\omega^2 x}{(1+xx)^2}+\frac{\omega^3}{(1+xx)^3}\left(xx-\frac{1}{3}\right)-\frac{\omega^4}{(1+xx)^4}(x^3-x)
\end{equation*}
\begin{equation*}
+\frac{\omega^5}{(1+xx)^5}\left(x^4-2x^2+\frac{1}{5}\right)-\frac{\omega^6}{(1+xx)^6}\left(x^5-\frac{10}{3}x^3+x\right)+\text{etc.}
\end{equation*}

\paragraph*{§87}

This series, whose law of progression is not that obvious, can be transformed into another form, whose structure is  immediately clear. For this aim, put $\arctan x= 90^{\circ}-u$, so that  $x= \cot u = \frac{\cos u}{\sin u}$; it will be $1+xx=\frac{1}{\sin^2 u}$, whence $\frac{dy}{dx}=\frac{1}{1+xx}=\sin^2 u$. Further, since $dx=\frac{-du}{\sin^2 u}$ or $du=-dx \sin^2 u$,  taking more differentials, it will be

\begin{equation*}
\frac{ddy}{dx}=2du \sin u \cdot \cos u = du \sin 2 u = -dx \sin^2 u \cdot \sin 2u
\end{equation*}
and hence

\begin{equation*}
\frac{ddy}{1dx^2}=-\sin^2 u \cdot \sin 2u,
\end{equation*}
\begin{equation*}
\frac{d^3y}{2dx^2}=-du \sin u \cdot \cos u \cdot \sin 2u -du \sin^2 u \cdot \cos 2u = - du \sin u \cdot \sin 3 u
\end{equation*}
\begin{equation*}
=dx \sin^3 u \cdot \sin 3u
\end{equation*}
and hence

\begin{equation*}
\frac{d^3y}{1 \cdot 2 dx^3}=+\sin^3 u \cdot \sin 3u,
\end{equation*}
\begin{equation*}
\frac{d^4 y}{1 \cdot 2 \cdot 3 dx^3}= du \sin^2 u (\cos u \cdot \sin 3u +\sin u \cdot \cos 3u)=du \sin^2 u \cdot \sin 4u
\end{equation*}
\begin{equation*}
=- dx \sin^4 \cdot \sin 4u
\end{equation*}
and hence

\begin{equation*}
\frac{d^4 y}{1 \cdot 2 \cdot 3 dx^4}=- \sin^4 \cdot \sin 4u,
\end{equation*}
\begin{equation*}
\frac{d^5 y}{1 \cdot 2 \cdot 3 \cdot 4 dx^4}=-du \sin^3 u(\cos u \cdot \sin 4u +\sin u \cdot \cos 4u)=-du \cdot \sin^3 u \cdot \sin 5u
\end{equation*}
\begin{equation*}
=+dx \sin^5 u \cdot 5u
\end{equation*}
and hence

\begin{equation*}
\frac{d^5y}{1 \cdot 2 \cdot 3 \cdot 4 dx^5}=+\sin^5 u \cdot \sin 5u
\end{equation*}
\begin{center}
etc.
\end{center}
From these one concludes that it will be

\begin{equation*}
\arctan (x+\omega)=\arctan x+ \frac{\omega}{1}\sin u \cdot \sin u -\frac{\omega^2}{2}\sin^2 u \cdot \sin 2u +\frac{\omega^3}{3}\sin^3 u \cdot \sin 3u
\end{equation*}
\begin{equation*}
-\frac{\omega^4}{4}\sin^4 u \cdot \sin 4u +\frac{\omega^5}{5}\sin^5 u \cdot \sin 5u -\frac{\omega^6}{6}\sin^6 \cdot \sin 6u +\text{etc.};
\end{equation*}
because here $\arctan x =y$ and $\arctan x= 90^{\circ}-u$, it will be $y=90^{\circ}-u$. 

\paragraph*{§88}

If one puts $\operatorname{arccot} x=y$ and $\operatorname{arccot}(x+\omega)=z$, it will be

\begin{equation*}
z=y+\frac{\omega dy}{dx}+\frac{\omega^2 ddy}{1 \cdot 2 dx^2}+\frac{\omega^3 d^3y}{1 \cdot 2 \cdot 3 dx^3}+\frac{\omega^4 d^4y}{1 \cdot2 \cdot 3 \cdot 4 dx^4}+\text{etc.}
\end{equation*}
But, since $dy=\frac{-dx}{1+xx}$, all terms with exception of the the first of this series agree with the ones found before, just with different signs. Hence, if, as before, one puts $\arctan x = 90^{\circ}-u$ or $\operatorname{arccot} x=u$, that $u=y$, it will be

\begin{equation*}
\operatorname{arccot}(x+\omega)= \operatorname{arccot}x -\frac{\omega}{1}\sin u \cdot \sin u +\frac{\omega^2}{2}\sin^2 u \cdot \sin 2 u -\frac{\omega^3}{3}\sin^3 u \cdot \sin 3u
\end{equation*}
\begin{equation*}
+\frac{\omega^4}{4}\sin^4 u \cdot \sin 4u -\frac{\omega^5}{5}\sin^5 u	\cdot \sin 5u +\text{etc.},
\end{equation*}
which expression follows from the preceding  immediately; for, since 

\begin{equation*}
\operatorname{arccot}(x+ \omega)= 90^{\circ}-\operatorname{arctan}(x+\omega) \quad \text{and} \quad \operatorname{arccot} x= 90^{\circ}-\arctan x,
\end{equation*}
it will be

\begin{equation*}
\operatorname{arccot}(x+\omega)-\operatorname{arccot}x=- \arctan (x+\omega)+\arctan x.
\end{equation*}

\paragraph*{§89}

From these expressions many extraordinary corollaries follow, depending on which values are substituted for the given $x$ and $\omega$. Therefore, first let $x=0$, and because  $u=90^{\circ}-\arctan x$, it will be $u=90^{\circ}$ and $\sin u=1$, $\sin 2u =0$, $\sin 3u =-1$, $\sin 4u =0$, $\sin 5u =1$, $\sin 6u =0$, $\sin 7u =-1$ etc., whence it will be

\begin{equation*}
\arctan \omega =\frac{\omega}{1}-\frac{\omega^3}{3}+\frac{\omega^5}{5}-\frac{\omega^7}{7}+\frac{\omega^9}{9}-\frac{\omega^{11}}{11}+\text{etc.},
\end{equation*}
which is the  well-known series expressing the arc whose tangent is $=\omega$.\\[2mm]
Let $x=1$; it will be $\arctan x=45^{\circ}$ and hence $u=45^{\circ}$, hence $\sin u=\frac{1}{\sqrt{2}}$, $\sin 2u =1$, $\sin 3u =\frac{1}{\sqrt{2}}$, $\sin 4u =0$, $\sin 5u =-\frac{1}{\sqrt{2}}$, $\sin 6u =-1$, $\sin 7u =-\frac{1}{\sqrt{2}}$, $\sin 8u =0$, $\sin 9u =\frac{1}{\sqrt{2}}$ etc.  Hence

\begin{equation*}
\arctan (1+\omega)=45^{\circ}+\frac{\omega}{2}-\frac{\omega^2}{2 \cdot 2}+\frac{\omega^3}{3 \cdot 4}-\frac{\omega^5}{5 \cdot 8}+\frac{\omega^6}{6 \cdot 8}-\frac{\omega^7}{7 \cdot 16}+\frac{\omega^9}{9 \cdot 32}-\frac{\omega^{10}}{10 \cdot 32}
\end{equation*}
\begin{equation*}
+\frac{\omega^{11}}{11 \cdot 64}-\frac{\omega^{13}}{13 \cdot 128}+\frac{\omega^{14}}{14 \cdot 128}-\text{etc.}
\end{equation*}
Therefore, if $\omega =-1$, because of $\arctan (1+\omega)=0$ and $45^{\circ}=\frac{\pi}{4}$, it will be

\begin{equation*}
\frac{\pi}{4}=\frac{1}{1 \cdot 2}+\frac{1}{2\cdot 2}+\frac{1}{3 \cdot 2^3}-\frac{1}{6 \cdot 2^3}-\frac{1}{7 \cdot 2^4}+\frac{1}{9 \cdot 2^5}+\frac{1}{10 \cdot 2^5}+\frac{1}{11 \cdot 2^6}-\text{etc.};
\end{equation*}
If this value  is substituted for $45^{\circ}$ in that expression, it will be
\begin{equation*}
\arctan(1+ \omega)=\frac{\omega +1}{1 \cdot 2}-\frac{\omega^2-1}{2 \cdot 2}+\frac{\omega^3+1}{3 \cdot 2^2}-\frac{\omega^5+1}{5\cdot 2^3}+\frac{\omega^6-1}{6 \cdot 2^3}-\frac{\omega^7+1}{7 \cdot 2^4}+\text{etc.}
\end{equation*}
But that series is most appropriate to find the value of $\frac{\pi}{4}$ approximately.

\paragraph*{§90}

Because

\begin{equation*}
\frac{\pi}{4}=\frac{1}{1 \cdot 2}+\frac{1}{2 \cdot 2}+\frac{1}{3 \cdot 2^2}-\frac{1}{5 \cdot 2^3}-\frac{1}{6 \cdot 2^3}-\frac{1}{7 \cdot 2^4}+\text{etc.},
\end{equation*}
but the terms containing $2$, $6$, $10$ etc. in the denominators, i.e.

\begin{equation*}
\frac{1}{2 \cdot 2}-\frac{1}{6 \cdot 2^3}+\frac{1}{10 \cdot 2^5}-\frac{1}{14 \cdot 2^7}+\text{etc.},
\end{equation*}
express  $\frac{1}{2}\arctan \frac{1}{2}$, it will be

\begin{equation*}
\frac{\pi}{4}=\frac{1}{2}\arctan \frac{1}{2}+\frac{1}{1 \cdot 2}+\frac{1}{3 \cdot 2^2}-\frac{1}{5 \cdot 2^3}-\frac{1}{7 \cdot 2^4}+\frac{1}{9 \cdot 2^5}+\frac{1}{11 \cdot 2^6}-\text{etc.}
\end{equation*}
But because in the other formula for negative $\omega$

\begin{alignat*}{9}
&\arctan (1- \omega) &&=\frac{1}{1 \cdot 2}&&+\frac{1}{2 \cdot 2}&&+\frac{1}{3 \cdot 2^2}&&-\frac{1}{5 \cdot 2^3}&&-\frac{1}{6 \cdot 2^3}&&-\frac{1}{7 \cdot 2^4}&&+\text{etc.} \\[3mm]
&  &&-\frac{\omega}{1 \cdot 2}&&-\frac{\omega^2}{2 \cdot 2}&&-\frac{\omega^3}{3 \cdot 2^2}&&+\frac{\omega^5}{5 \cdot 2^3}&&+\frac{\omega^6}{6 \cdot 2^3}&&+\frac{\omega^7}{7 \cdot 2^4}&&-\text{etc.}, \\
\end{alignat*}
if  $\omega= \frac{1}{2}$ , it will be 

\begin{alignat*}{9}
&\arctan (1- \omega) &&=\frac{1}{1 \cdot 2}&&+\frac{1}{2 \cdot 2}&&+\frac{1}{3 \cdot 2^2}&&-\frac{1}{5 \cdot 2^3}&&-\frac{1}{6 \cdot 2^3}&&-\frac{1}{7 \cdot 2^4}&&+\text{etc.} \\[3mm]
&  &&-\frac{1}{1 \cdot 2}&&-\frac{1}{2 \cdot 2}&&-\frac{1}{3 \cdot 2^2}&&+\frac{1}{5 \cdot 2^3}&&+\frac{1}{6 \cdot 2^3}&&+\frac{1}{7 \cdot 2^4}&&-\text{etc.} \\
\end{alignat*}
and, having taken the terms divided by $2$, $6$, $10$ etc., it will be

\begin{alignat*}{19}
&\arctan \frac{1}{2}&&= \frac{1}{2}\arctan \frac{1}{2}&&+\frac{1}{1 \cdot 2}&&+\frac{1}{3 \cdot 2^3}&&-\frac{1}{5 \cdot 2^3}&&-\frac{1}{7 \cdot 2^4}&&+\frac{1}{9 \cdot 2^5}&&+\text{etc.} \\[3mm]
& &&- \frac{1}{2}\arctan \frac{1}{2}&&-\frac{1}{1 \cdot 2^2}&&-\frac{1}{3 \cdot 2^5}&&+\frac{1}{5 \cdot 2^8}&&+\frac{1}{7 \cdot 2^{11}}&&-\frac{1}{9 \cdot 2^{14}}&&-\text{etc.} \\
\end{alignat*}
and hence

\begin{alignat*}{9}
&\frac{1}{2}\arctan \frac{1}{2}&&= &&+\frac{1}{1 \cdot 2}&&+\frac{1}{3 \cdot 2^2}&&-\frac{1}{5 \cdot 2^3}&&-\frac{1}{7 \cdot 2^4}&&+\text{etc.} \\[3mm]
& &&-\frac{1}{2}\arctan \frac{1}{8}&&-\frac{1}{1 \cdot 2^2}&&-\frac{1}{3 \cdot 2^5}&&+\frac{1}{5\cdot 5^8}&&+\frac{1}{7 \cdot 2^{11}}&&-\text{etc.};
\end{alignat*}
if this value is substituted in the above series and $\arctan \frac{1}{8}$ is converted into a series, one will find

\begin{equation*}
\frac{\pi}{4}= \left\lbrace
\renewcommand{\arraystretch}{2,5}
\setlength{\arraycolsep}{0.0mm}
\begin{array}{rcccccccccccccccccc}
1 &~+~& \dfrac{1}{3 \cdot 2^1}& ~-~& \dfrac{1}{5 \cdot 2^2} &~-~& \dfrac{1}{7 \cdot 2^3}& ~+~& \dfrac{1}{9 \cdot 2^4}&~+~& \text{etc.} \\ 
-~\dfrac{1}{1 \cdot 2^2} &~-~& \dfrac{1}{3 \cdot 2^5}& ~+~& \dfrac{1}{5 \cdot 2^8} &~+~& \dfrac{1}{7 \cdot 2^{11}}& ~-~& \dfrac{1}{9 \cdot 2^{14}}&~-~& \text{etc.} \\
-~\dfrac{1}{1 \cdot 2^4} &~-~& \dfrac{1}{3 \cdot 2^{10}}& ~-~& \dfrac{1}{5 \cdot 2^{16}} &~+~& \dfrac{1}{7 \cdot 2^{22}}& ~-~& \dfrac{1}{9 \cdot 2^{28}}&~+~& \text{etc.} \\
\end{array}
\right.
\end{equation*} 

\paragraph*{ § 90a}

These and many others follow, if one puts $x=1$; but if we put $x=\sqrt{3}$ that $\arctan x = 60^{\circ}$, it will be $u=30^{\circ}$ and $\sin u= \frac{1}{2}$, $\sin 2u =\frac{\sqrt{3}}{2}$, $\sin 3u = 1$, $\sin 4u =\frac{\sqrt{3}}{2}$, $\sin 5u =\frac{1}{2}$, $\sin 6u =0$, $\sin 7u =-\frac{1}{2}$ etc., whence it will be

\begin{equation*}
\arctan (\sqrt{3}+\omega)= 60^{\circ}+\frac{\omega}{1 \cdot 2^2}-\frac{\omega^2 \sqrt{3}}{2 \cdot 2^3}+\frac{\omega^3}{3 \cdot 2^3}-\frac{\omega^4 \sqrt{3}}{4 \cdot 2^5}+\frac{\omega^5}{5 \cdot 2^6}-\frac{\omega^7}{7 \cdot 2^8}
\end{equation*}
\begin{equation*}
+\frac{\omega^8 \sqrt{3}}{8 \cdot 2^9}-\frac{\omega^9}{9 \cdot 2^9}+\frac{\omega^{10}\sqrt{3}}{10 \cdot 2^{11}}-\frac{\omega^{11}}{11 \cdot 2^{12}}+\text{etc.}
\end{equation*}
But if one puts $x=\frac{1}{\sqrt{3}}$, so that $\arctan x= 30^{\circ}$, it will be $u=60^{\circ}$ and $\sin u =\frac{\sqrt{3}}{2}$, $\sin 2u =\frac{\sqrt{3}}{2}$, $\sin 3u =0$, $\sin 4u =-\frac{\sqrt{3}}{2}$, $\sin 6u =0$, $\sin 7u =\frac{\sqrt{3}}{2}$ etc., having substituted which values it will be

\begin{equation*}
\arctan \left(\frac{1}{\sqrt{3}}+\omega\right)= 30^{\circ}+\frac{3 \omega}{1 \cdot 2^2}-\frac{3 \omega^2 \sqrt{3}}{2 \cdot 2^3}+\frac{3^2 \omega^4 \sqrt{3}}{4 \cdot 2^5}-\frac{3^3 \omega^5}{5 \cdot 2^5}+\text{etc.};
\end{equation*}
therefore, if $\omega =-\frac{1}{\sqrt{3}}$, because of $30^{\circ}=\frac{\pi}{6}$, it will be

\begin{equation*}
\frac{\pi}{6 \sqrt{3}}=\frac{1}{1 \cdot 2^2}+\frac{1}{2 \cdot 2^3}-\frac{1}{4 \cdot 2^5}-\frac{1}{5\cdot 2^6}+\frac{1}{7 \cdot 2^8}+\frac{1}{8 \cdot 2^9}-\text{etc.}
\end{equation*}

\paragraph*{§91}

Let us return to the general expression we found, i.e.

\begin{equation*}
\arctan (x+ \omega)
\end{equation*}
\begin{equation*}
= \arctan x+ \frac{ \omega}{1}\sin u \cdot \sin u -\frac{\omega^2}{2}\sin^2 u \cdot \sin 2u +\frac{\omega^3}{3}\sin^3 u \cdot \sin 3u -\text{etc.}
\end{equation*}
and let us put $\omega =-x$, so that  $\arctan (x+ \omega)=0$ and it will be

\begin{equation*}
\arctan x= \frac{x}{1}\sin u \cdot \sin u +\frac{x^2}{2}\sin^2 u \cdot \sin 2u +\frac{x^3}{3}\sin^3 u \cdot \sin 3u +\text{etc.}
\end{equation*}
But because $\arctan x = 90^{\circ}-u=\frac{\pi}{2}-u$, it will be $x= \cot u= \frac{\cos u}{\sin u}$. Therefore, it will be

\begin{equation*}
\frac{\pi}{2}= u+\cos u \cdot \sin u +\frac{1}{2}\cos^2 u \cdot \sin 2u +\frac{1}{3}\cos^3 u \cdot \sin 3u +\frac{1}{4} \cos^4 u \cdot \sin 4u +\text{etc.},
\end{equation*}
which series is even more remarkable, since, whatever arc is taken for $u$, the value of the series always turns out to be  the same, $=\frac{\pi}{2}$.\\[2mm]
But if  $\omega =-2x$, because of $\arctan (-x)=-\arctan x$, it will be

\begin{equation*}
2 \arctan x =\frac{2x}{1}\sin u \cdot \sin u+\frac{4x^2}{2}\sin^2 u \cdot \sin 2u +\frac{8}{3}\sin^3 u \cdot \sin 3u +\text{etc.}
\end{equation*}
But because  $\arctan x=\frac{\pi}{2}-u$ and $x=\frac{\cos u}{\sin u}$, it will be

\begin{equation*}
\pi =2u +\frac{2}{1}\cos u \cdot \sin u+\frac{2^2}{2}\cos^2 u \cdot \sin 2u +\frac{2^3}{3}\cos^3 u \cdot \sin 3u +\text{etc.}
\end{equation*}
Let $u=45^{\circ}$; it will be $\cos u= \frac{1}{\sqrt{2}}$, $\sin u= \frac{1}{\sqrt{2}}$, $\sin 2u =1$, $\sin 3u =\frac{1}{\sqrt{2}}$, $\sin 4u =0$, $\sin 5u = \frac{-1}{\sqrt{2}}$, $\sin 6u =-1$, $\sin 7u =\frac{-1}{\sqrt{2}}$, $\sin 8u =0$, $\sin 9u =\frac{1}{\sqrt{2}}$ etc. and it will be

\begin{equation*}
\frac{\pi}{2}=\frac{1}{1}+\frac{2}{2}+\frac{2}{3}-\frac{2^5}{5}-\frac{2^3}{6}-\frac{2^3}{7}+\frac{2^4}{9}+\frac{2^5}{10}+\frac{2^{5}}{11}-\text{etc.},
\end{equation*}
which series, even though it diverges,  nevertheless is remarkable for its simplicity.

\paragraph*{§92}

In the general expression we found put

\begin{equation*}
\omega =-x-\frac{1}{x}=\frac{-1}{\sin u \cdot \cos u};
\end{equation*}
because of $x=\frac{\cos u}{\sin u}$, it will be

\begin{equation*}
\arctan (x+\omega)=\arctan \left(-\frac{1}{x}\right)=-\arctan \frac{1}{x}=-\frac{\pi}{2}+\arctan x.
\end{equation*}
Therefore, one will hence obtain the following expression

\begin{equation*}
\frac{\pi}{2}=\frac{\sin u}{1 \cos u}+\frac{\sin 2u}{2 \cos^2 u}+\frac{\sin 3u}{3 \cos^3 u}+\frac{\sin 4u}{4 \cos^4 u}+\frac{\sin 5u}{5 \cos^5 u}+\text{etc.},
\end{equation*}
which, having put $u=45^{\circ}$, gives the same series we found last. \\[2mm]
But if we put $\omega =-\sqrt{1+xx}$, because of $x=\frac{\cos u}{\sin u}$, it will be

\begin{equation*}
\omega =-\frac{1}{\sin u}
\end{equation*}
and

\begin{equation*}
\arctan (x-\sqrt{1+xx})=-\arctan (\sqrt{1+xx}-x)
\end{equation*}
\begin{equation*}
=-\frac{1}{2}\arctan \frac{1}{x}=-\frac{1}{2}\left(\frac{\pi}{2}-\arctan x\right)=-\frac{1}{2}u
\end{equation*}
and 

\begin{equation*}
\arctan x =\frac{\pi}{2}-u.
\end{equation*}
Therefore, it will be

\begin{equation*}
\frac{\pi}{2}=\frac{1}{2}u+\frac{1}{1}\sin u +\frac{1}{2}\sin 2u+\frac{1}{3}\sin 3u +\frac{1}{4}\sin 4u +\text{etc.}
\end{equation*}
Therefore, if this equation is differentiated, it will be

\begin{equation*}
0=\frac{1}{2}+\cos u +\cos 2u +\cos 3u +\cos 4u+ \cos 5u +\text{etc.},
\end{equation*}
whose correctness is seen from the nature of recurring series.

\paragraph*{§93}

In like manner, if the series found before are differentiated, new summable series will be found. First, from the series

\begin{equation*}
\arctan (1+\omega)=\frac{\pi}{4}+\frac{\omega}{2}-\frac{\omega^2}{2 \cdot 2}+\frac{\omega^3}{3 \cdot 4}-\frac{\omega^5}{5 \cdot 8}+\frac{\omega^6}{6 \cdot 8}-\text{etc.}
\end{equation*}
it follows

\begin{equation*}
\frac{1}{2+2\omega +\omega^2}=\frac{1}{2}-\frac{\omega}{2}+\frac{\omega^2}{4}-\frac{\omega^4}{8}+\frac{\omega^5}{8}-\frac{\omega^6}{16}+\frac{\omega^8}{32}-\text{etc.},
\end{equation*}
which results from the expansion of this fraction $\frac{2-2 \omega +\omega^2}{4+\omega^4}=\frac{1}{2+ 2 \omega +\omega^2}$.\\[2mm]
Further, this series

\begin{equation*}
\frac{\pi}{2}=u +\cos u \cdot \sin u + \frac{1}{2}\cos^2 u \cdot \sin 2u +\frac{1}{3}\cos^3 u \cdot \sin 3u+\frac{1}{4} \cos^4 \cdot \sin 4u +\text{etc.},
\end{equation*}
by means of differentiation, will give

\begin{equation*}
0=1 +\cos 2u + \cos \cdot \cos 3u +\cos^2 u \cdot \cos 4u +\cos^3 u \cdot \cos 5u +\text{etc.}
\end{equation*}
Finally, the series

\begin{equation*}
\frac{\pi}{2}= \frac{\sin u}{\cos u}+\frac{\sin 2u}{2 \cos^2 u}+ \frac{\sin 3u}{3 \cos^3 u}+\frac{\sin 4u}{4 \cos^4 u}+\text{etc.}
\end{equation*}
gives

\begin{equation*}
0=\frac{1}{\cos^2 u}+\frac{\cos u}{\cos^3 u}+\frac{\cos 2u}{\cos^4 u}+\frac{\cos 3u}{\cos^5 u}+\frac{\cos 4u}{\cos^6 u}+\text{etc.}
\end{equation*}
or 

\begin{equation*}
0=1+\frac{\cos u}{\cos u}+\frac{\cos 2u}{\cos^2 u}+\frac{\cos 3u}{\cos^3 u}+\frac{\cos 4u}{\cos^4 u}+\frac{\cos 5u}{\cos^5 u}+\text{etc.}
\end{equation*}

\paragraph*{§94}

But especially the expression we found

\begin{equation*}
\arctan (x+ \omega)
\end{equation*}
\begin{equation*}
=\arctan x+ \frac{\omega}{1}\sin u \cdot \sin u-\frac{\omega^2}{2}\sin^2 u \cdot \sin 2u +\frac{\omega^3}{3}\sin^3 u \cdot \sin 3u -\text{etc.},
\end{equation*}
while 

\begin{equation*}
x= \cot u \quad \text{or} \quad u= \operatorname{arccot}x= 90^{\circ}-\arctan x,
\end{equation*}
can be applied to the angle or the arc corresponding to a given certain tangent. For, let the tangent $=t$ be propounded and, consulting tables, find the tangent coming closet to this, call it $x$, to which the arc $=y$ corresponds, and it will be $u=90^{\circ}-y$. Then, put $x+ \omega =t$ or $\omega =t-x$ and the arc in question will be

\begin{equation*}
=y+\frac{\omega}{1}\sin u \cdot \sin u -\frac{\omega^2}{2}\cdot \sin 2u +\text{etc.},
\end{equation*}
which rule is especially useful, when the propounded tangent was very large and therefore the arc in question hardly deviates from $90^{\circ}$. For, in these cases, because of the rapidly increasing tangents, the usual method of interpolation leads too far away from the true value. Therefore, let this example be propounded.

\subsection*{Example}

\textit{Let the arc be in question, whose tangent is $=100$, having put the radius $=1$, of course.}\\[2mm]
The arc approximately equal to the one in question is $89^{\circ}25^{\text{I}}$, whose tangent is

\begin{tabular}{rrrr}
 & $x$ & $=$ & $98.217943$ \\
subtract this from & $t$ & $=$ & $100.000000$ \\ \cline{4-4}
it will remain & $\omega$ & $=$ & $1.782057$
\end{tabular}\\
Further, because  $y=89^{\circ}25^{\text{I}}$, it will be $u=35^{\text{I}}$, $2u=1^{\circ}10^{\text{I}}$, $3u =1^{\circ}45^{\text{I}}$ etc. Now investigate each term by means of logarithms.  \\[2mm]

\begin{center}
\begin{tabular}{lrrrr}
To & $\log \omega$ & $=$  & $0.2509125$ \\
add  & $\log \sin u $ & $=$ & $8.0077867$ \\
    & $\log \sin u $ & $=$ & $8.0077867$ \\ \cline{4-4}
subtract   & $\log \omega \sin u \cdot \sin u $ & $=$ & $6.2664949$ \\
           &  &  & $4.6855749$ \\  
            &  & $=$ & $1.5809200$ \\  
Therefore & $\omega \sin u \cdot \sin u$ & $=$ & $38.09956$ & seconds \\
To        & $\log \omega \sin^2 u$ & $=$ & $6.2664949$ \\
add       & $\log \omega$ & $=$ & $0.2509215$ \\
          & $\log 2u $    & $=$ & $8.3087941$ \\ \cline{4-4} 
subtract  &               &     & $4.8262105$ \\
          & $\log 2$      & $=$ & $0.3010300$ \\ \cline{4-4}
          & $\log \dfrac{1}{2}\omega^2 \sin^2 u \cdot 2u$ & $=$ & $4.5251805$ \\
subtract  &               &     & $4.6855749$ \\ \cline{4-4}
it remains &              &     & $9.8396056$ \\  
Therefore  & $\dfrac{1}{2}\omega^2 \sin^2 u \cdot \sin 2u$ & $=$ & $0.6912000$ & seconds \\ \cline{4-4} 
Further to & $\log \omega^3$ & $=$ & $0.7527645$ \\
           & $\log \sin^3 u$ & $=$ & $4.0233601$ \\
subtract           & $\log \sin 3u $ & $=$ & $3.2609725$ \\   
                   & $\log 3$ & $=$ & $0.4771213$ \\ \cline{4-4}
                   &          &     & $2.7838512$ \\
subtract           &          &     & $4.6855749$ \\ \cline{4-4}
                   &          &     & $8.0982763$ \\ \cline{4-4}                
Therefore          & $\dfrac{1}{3}\omega^3 \sin^3 u \cdot \sin 3u$ & $=$ & $0.0125400$ & seconds.    \\               
\end{tabular} \\[2mm]
\begin{tabular}{lrrrr}
Finally, to & $\log \omega^4$ & $=$& $1.0036860$ \\
add         & $\log \sin^4 u$ & $=$& $2.0311468$ \\
            & $\log \sin 4u$  & $=$& $8.6097341$ \\ \cline{4-4}
            &                 &    & $1.6445669$ \\
subtract    & $\log 4$        & $=$& $0.6020600$ \\ \cline{4-4}
            &                 &    & $1.0425069$ \\
subtract    &                 &    & $4.6855749$ \\ \cline{4-4}
            &                 &    & $6.3569320$                                    
\end{tabular} 
\end{center}
Therefore,

\begin{equation*}
\frac{1}{4}\omega^4 \sin^4 u =0.00023 \quad \text{seconds}
\end{equation*}
Hence

\begin{center}
\begin{tabular}{cc|c}
& Terms to be added & Terms to be subtracted \\
& $38.09956$ & $0.69120$ \\
& $0.01254$  & $0.00023$ \\ \cline{2-3}
subtract & $0.69143$ & 
\end{tabular}
\end{center}
Hence in total

\begin{equation*}37.4067  =37^{\text{II}}25^{\text{III}}14^{\text{IV}}24^{\text{V}}36^{\text{VI}}.
\end{equation*}
Therefore, the arc, whose tangent is hundred times  the radius, will be

\begin{equation*}
89^{\circ}25^{\text{I}}37^{\text{II}}25^{\text{III}}14^{\text{IV}}24^{\text{V}}36^{\text{VI}}
\end{equation*}
and the error does not affect the fourth, but can only occur the fifth, whence we will be able to confirm that this angle  is almost $=89^{\circ}25^{\text{I}}37^{\text{II}}25^{\text{III}}14^{\text{IV}}$. If an even greater tangent is propounded, even though $\omega$ might turn out to be larger,  because of the still small angle $u$, one will nevertheless be able to define the arc in a convenient way.

\paragraph*{§95}

Since here we substituted an arc of the circle for $y$, let us now substitute the inverse functions for $y$, i.e. $\sin x$, $\cos x$, $\tan x$, $\cot x$ etc. Therefore, let $y=\sin x$, and, having written $x+ \omega$ instead of $x$, it will be $z= \sin (x+ \omega)$ and the equation

\begin{equation*}
z = y+ \frac{\omega dy}{dx}+\frac{\omega^2 ddy}{2dx^2}+\frac{\omega^3 d^3 y}{6dx^3}+\frac{\omega^4 d^4 y}{24 dx^4}+\text{etc.},
\end{equation*}
because of

\begin{equation*}
\frac{dy}{dx}=\cos x, \quad \frac{ddy}{dx^2}=-\sin x, \quad \frac{d^3y}{dx^3}=-\cos x, \quad \frac{d^4y}{dx^4}=\sin x \quad \text{etc.},
\end{equation*}
will give 

\begin{equation*}
\sin (x+ \omega)=\sin x +\omega \cos x-\frac{1}{2}\omega^2 \sin x -\frac{1}{6}\omega^3 \cos x  +\frac{1}{24}\omega^4 \sin x+\text{etc.}
\end{equation*}
and, having taken a negative $\omega$, it will be

\begin{equation*}
\sin (x- \omega)=\sin x -\omega \cos x-\frac{1}{2}\omega^2 \sin x +\frac{1}{6}\omega^3 \cos x  +\frac{1}{24}\omega^4 \sin x-\text{etc.}
\end{equation*}
So, if one sets $y= \cos x$, because of

\begin{equation*}
\frac{dy}{dx}=-\sin x, \quad \frac{ddy}{dx^2}=-\cos x, \quad \frac{d^3y}{dx^3}= \sin x, \quad \frac{d^4 y}{dx^4}= \cos x \quad \text{etc.},
\end{equation*}
it will be

\begin{equation*}
\cos (x+ \omega)=\cos x - \omega \sin x -\frac{1}{2} \omega^2 \cos x +\frac{1}{6}\omega^3 \sin x +\frac{1}{24} \omega^4 \cos x -\text{etc.}
\end{equation*}
and for a negative $\omega$ it will be

\begin{equation*}
\cos (x- \omega)=\cos x + \omega \sin x -\frac{1}{2} \omega^2 \cos x -\frac{1}{6}\omega^3 \sin x +\frac{1}{24} \omega^4 \cos x +\text{etc.}
\end{equation*}

\paragraph*{§96}

The use of these formulas is immense both for the construction and interpolation of tables of sines and cosines. For, if  the sines and cosines of a certain arc $x$ were known, from them  the sines and cosines of the angles $x+ \omega$ and $x-\omega$ can be easily found, if the difference $\omega$ was sufficiently small; for, in this case the found series converge rapidly. For this it is necessary that the arc $\omega$ is expressed in parts of the radius; this, because the arc of $180^{\circ}$ is

\begin{equation*}
3.14159265358979323846,
\end{equation*}
is easily done; for, after a division by $180$ it will be

\begin{alignat*}{9}
& \text{arc of} ~1^{\circ} &&=0.017453292519943295769 \\
& \text{arc of} ~1^{\text{I}} &&=0.000290888208665721596 \\
& \text{arc of} ~1^{\text{II}} &&=0.000048481368110953599. \\
\end{alignat*}

\subsection*{Example 1}

\textit{To find the sine and the cosine of the angles $45^{\circ}1^{\text{I}}$ and $44^{\circ}59^{\text{I}}$, the given sine and cosine of the angle $45^{\circ}$, both of which are $=\frac{1}{\sqrt{2}}=0.707167811865$.} \\[2mm]
Therefore, since 

\begin{alignat*}{9}
& \sin x && ~=~ &&\cos x && =0.7071067811865 
\intertext{and}
&        &&  && \quad ~~ \omega   && =0.0002908882086,  
\end{alignat*}
in order to perform the multiplication more easily, note that it will be

\begin{alignat*}{9}
& 2 \omega && = 0.0005817764173 \\
& 3 \omega && = 0.0008726646259 \\
& 4 \omega && = 0.0011635528346 \\
& 5 \omega && = 0.0014544410433 \\
& 6 \omega && = 0.0017453292519 \\
& 7 \omega && = 0.0020362174606 \\
& 8 \omega && = 0.0023271056693 \\
& 9 \omega && = 0.0026179938779 \\
\end{alignat*}
Therefore, $\omega \sin x $ and $\omega \cos x$ will be found this way:
\begin{center}
\begin{tabular}{rrr}
$7 \quad$ &$\cdot$& $0.00020362174606$ \\
$0 \quad$        &$\cdot$& \\
$7 \quad$ &$\cdot$& $0.00000203621746$ \\
$1 \quad$ &$\cdot$& $2908882$ \\
$0 \quad$ &$\cdot$& \\
$6 \quad$ &$\cdot$& $174532$ \\ 
$7 \quad$ &$\cdot$& $20362$ \\
$8 \quad$ &$\cdot$& $2327$ \\
$1 \quad$ &$\cdot$& $29$ \\
$1 \quad$ &$\cdot$& $2$ \\
$8 \quad$ &$\cdot$& $2$ \\
$6 \quad$ &$\cdot$& $0$ \\ \hline
\end{tabular}
\end{center}
In total

\begin{equation*}
\omega \sin x = \omega \cos x = 0.00020568902488
\end{equation*}
Therefore,

\begin{tabular}{rrrr}
&$\dfrac{1}{2}\omega \cos x$ & $=$ & $0.00010284451244$ \\ \cline{4-4}
by $\omega$ & $1$ & $\cdot$ & $0.00000002908882$ \\
  & $0$ & $\cdot$ &  \\
  & $2$ & $\cdot$ & $58177$ \\
  & $8$ & $\cdot$ & $23271$ \\
  & $4$ & $\cdot$ & $1163$ \\
  & $4$ & $\cdot$ & $116$ \\
  & $5$ & $\cdot$ & $14$ \\ \cline{4-4}
  & $\dfrac{1}{2}\omega^2 \cos x$ & $=$ & $0.00000002991623$ \\[3mm]
  & $\dfrac{1}{6}\omega^3 \cos x$ & $=$ & $0.00000000997208$ \\[3mm] \cline{4-4}
 by $\omega$ & $9$ & $\cdot $ & $0.0000000000261$ \\ 
  & $9$ & $\cdot $ & $26$ \\ 
  & $7$ & $\cdot $ & $2$ \\ \cline{4-4}
  & $\dfrac{1}{6}\omega^3 \cos x$ & $=$ & $0.0000000000289$ 
\end{tabular} \\[2mm]
Therefore, to find $\sin 45^{\circ}1^{\text{I}}$ to 

\begin{tabular}{rrrrrr}
         & $\sin x $ & $=$ & $0.7071067811865$ \\
add  & $\omega \cos x $ & $=$ & $2056890249$ \\ \cline{4-4}
  & & & $0.7073124702114$ \\
  subtract & $\dfrac{1}{2}\omega^2 \sin x$ & $=$ & $299162$ \\ \cline{4-4}
   & & & $0.7073124402952$ \\ 
    subtract & $\dfrac{1}{6}\omega^3 \cos x$ & $=$ & $29$ \\ \cline{4-4}
      & $\sin 45^{\circ} 1^{\text{I}}$ & $=$ & $0.7073124402923$ & $=$ & $\cos 44^{\circ}59^{\text{I}}$
\end{tabular}\\[2mm]
But to find $\cos 45^{\circ}1^{\text{I}}$ from

\begin{tabular}{rrrrrr}
         & $\cos x $ & $=$ & $0.7071067811865$ \\
subtract & $\omega \sin x $ & $=$ & $2056890249$ \\ \cline{4-4}
  & & & $0.70769010921616$ \\
  subtract & $\dfrac{1}{2}\omega^2 \cos x$ & $=$ & $299162$ \\ \cline{4-4}
   & & & $0.7069010622454$ \\ 
    add & $\dfrac{1}{6}\omega^3 \cos x$ & $=$ & $29$ \\ \cline{4-4}
      & $\cos 45^{\circ} 1^{\text{I}}$ & $=$ & $0.7069010622483$ & $=$ & $\sin 44^{\circ}59^{\text{I}}$
\end{tabular}\\[2mm]

\subsection*{Example 2}

\textit{Given sine and cosine of the arc $67^{\circ}30^{\text{I}}$, to find the sine and the cosine of the arcs $67^{\circ}31^{\text{I}}$ and $67^{\circ}29^{\text{I}}$}. \\[2mm]
Let us perform this calculation in decimal fractions up to $7$ digits, as the common tables are usually constructed, and hence the task will easily be solved applying logarithms. Because $x=67^{\circ}30^{\text{I}}$ and $\omega = 0.000290888$, it will be

\begin{equation*}
\log \omega =6.4637259
\end{equation*}

\begin{tabular}{rrrrrrr}
and & $\log \sin x $ & $=$ & $9.9656153$ \quad \quad \quad & $\log \cos x$ & $=$ & $9.5828397$ \\
  & $\log \omega $ & $=$ & $6.6437259$ \quad \quad \quad & $\log \omega$ & $=$ & $6.6437259$ \\ \cline{2-4} \cline{5-7}
    & $\log \omega \sin x$ & $=$ & $6.4293412$ \quad \quad \quad & $\log \omega \cos x$ & $=$ & $6.0465656$ \\[3mm]
      & $\log \dfrac{1}{2}\omega$ & $=$ & $6.1626959$ \quad \quad \quad & $\log \dfrac{1}{2} \omega^2$ & $=$ & $6.1626959$ \\[3mm] \cline{4-4} \cline{7-7}
       & $\log \dfrac{1}{2}\omega^2 \sin x$ & $=$ & $2.59200371$ \quad \quad \quad & $\log \dfrac{1}{2} \omega^2 \cos x$ & $=$ & $2.2092615$ \\[4mm]
 Therefore      & $\omega \sin x$ & $=$ & $0.00026874$ \quad \quad \quad & $ \omega \cos x$ & $=$ & $0.00011132$ \\
 & $\dfrac{1}{2}\omega^2 \sin x$ & $=$ & $0.00000004$ \quad \quad \quad & $ \dfrac{1}{2}\omega^2 \cos x$ & $=$ & $0.00000001$ \\[3mm]
whence & $ \sin 67^{\circ}31^{\text{I}}$ & $=$ & $0.9239908$ \quad \quad \quad & $ \cos 67^{\circ}31^{\text{I}}$ & $=$ & $0.3824147$ \\
 & $ \sin 67^{\circ}29^{\text{I}}$ & $=$ & $0.9237681$ \quad \quad \quad & $ \cos 67^{\circ}29^{\text{I}}$ & $=$ & $0.3829522$ \\
\end{tabular}\\[2mm]
where not even the terms $\frac{1}{2}\omega^2 \sin x$ and $\frac{1}{2}\omega^2 \cos x$ were necessary.

\paragraph*{§97}

From the series we found above,

\begin{alignat*}{9}
&\sin (x+ \omega) &&= \sin x &&+\omega \cos x &&-\frac{1}{2}\omega^2 \sin x && -\frac{1}{6}\omega^3 \cos x && +\frac{1}{24}\omega^4 \sin x &&+\text{etc.} \\[3mm]
&\cos (x+ \omega) &&= \cos x &&-\omega \sin x &&-\frac{1}{2}\omega^2 \cos x && +\frac{1}{6}\omega^3 \sin x && +\frac{1}{24}\omega^4 \cos x &&-\text{etc.} \\[3mm]
&\sin (x- \omega) &&= \sin x &&-\omega \cos x &&-\frac{1}{2}\omega^2 \sin x && +\frac{1}{6}\omega^3 \cos x && +\frac{1}{24}\omega^4 \sin x &&-\text{etc.} \\[3mm]
&\cos (x- \omega) &&= \cos x &&+\omega \sin x &&-\frac{1}{2}\omega^2 \cos x && -\frac{1}{6}\omega^3 \sin x && +\frac{1}{24}\omega^4 \cos x &&+\text{etc.},
\end{alignat*}
combining them we will find

\begin{equation*}
\frac{\sin (x+\omega)+\sin (x- \omega)}{2}
\end{equation*}
\begin{equation*}
= \sin x -\frac{1}{2}\omega^2 \sin x +\frac{1}{24}\sin x -\frac{1}{720}\omega^6 \sin x +\text{etc.}= \sin x \cdot \cos \omega
\end{equation*}
and

\begin{equation*}
\frac{\sin (x+\omega)-\sin (x- \omega)}{2}
\end{equation*}
\begin{equation*}
= \omega \cos x -\frac{1}{6}\omega^3 \cos x +\frac{1}{120}\omega^5 \cos x -\frac{1}{5040}\omega^7 \cos x +\text{etc.}=\cos x \cdot \sin \omega,
\end{equation*}
whence the series found above for the sines and cosines result to be

\begin{equation*}
\cos \omega =1- \frac{1}{2}\omega^2 +\frac{1}{24}\omega^4 -\frac{1}{720}\omega^6 + \text{etc.}
\end{equation*}
\begin{equation*}
\sin \omega = \omega -\frac{1}{6}\omega^3 +\frac{1}{120}\omega^5 -\frac{1}{5040}\omega^7 +\text{etc.},
\end{equation*}
which  same series follow from the first for $x=0$; for, because  $\cos x=1$ and $\sin x =0$,  the first series will exhibit $\sin \omega$, the second on the other hand $\cos \omega$.

\paragraph*{§98}

Now let us also put $y= \tan x $, so that $z=\tan (x+ \omega)$; because of $y=\frac{\sin x}{\cos x}$ [§ 206 of the first part]

\begin{equation*}
\frac{dy}{dx}=\frac{1}{\cos^2 x}, \quad \frac{ddy}{2dx^2}=\frac{\sin x}{\cos^3 x}, \quad \frac{d^3y}{2dx^3}=\frac{1}{\cos^2 x}+\frac{3 \sin^2 x}{\cos^4 x}=\frac{3}{\cos^4 x}-\frac{2}{\cos^2 x},
\end{equation*} 
\begin{equation*}
\frac{d^4y}{2 \cdot 4 dx^4}=\frac{3 \sin x}{\cos^5 x}-\frac{\sin x}{\cos^3 x}, \quad \frac{d^5y}{2 \cdot 4}=\frac{15}{\cos^6 x}-\frac{15}{\cos^4 x}+\frac{2}{\cos^2 x},
\end{equation*}
whence it follows that

\begin{equation*}
\tan (x+ \omega)= \tan x +
\left\lbrace
\renewcommand{\arraystretch}{2,5}
\setlength{\arraycolsep}{0.0mm}
\begin{array}{rcccccccccccccccccc}
\dfrac{\omega}{\cos^2 x}&~+~& \dfrac{\omega^2 \sin x}{\cos^3 x}&~+~& \dfrac{\omega^3}{\cos^4 x}&~+~& \dfrac{\omega^4 \sin x}{\cos^5 x}&~+~& \text{etc.} \\
& & &~-~& \dfrac{2\omega^3}{3\cos^2 x}&~-~& \dfrac{\omega^4 \sin x}{3\cos^3 x}&~-~& \text{etc.} \\
\end{array}
\right.
\end{equation*}
by means of which formula, given tangent of any angle, one can find the tangents of  angles very close to it. Since the above series is a geometric one, having collected it into one sum, it will be

\begin{equation*}
\tan (x+ \omega)= \tan x+ \frac{\omega +\omega^2 \tan x}{\cos^2 x-\omega^2}-\frac{2 \omega^3}{3 \cos^2 x}-\frac{\omega^4 \sin x}{3 \cos^3 x}-\text{etc.}
\end{equation*}
or

\begin{equation*}
\tan (x+ \omega)=\frac{\sin x \cdot \cos x +\omega}{\cos^2 x- \omega^2}-\frac{2 \omega^3}{3 \cos^2 x}-\frac{\omega^4 \sin x}{3 \cos^3 x}-\text{etc.},
\end{equation*}
which formula is applied more conveniently for this aim.

\paragraph*{§99}

Similar expressions can also be found for the logarithms of sines, cosines and tangents. For, let $y=$ a logarithm of the sine of the angle $x$, which we want to express as $y= \log \sin x$, and $z=\log \sin (x+ \omega)$; because of $\frac{dy}{dx}=\frac{n \cos x}{\sin x}$, it will be $\frac{ddy}{dx^2}=\frac{-n}{\sin^2 x}$, $\frac{d^3y}{dx^3}=\frac{+ 2n \cos x}{\sin^3 x}$ etc., whence it will be

\begin{equation*}
z = \log \sin (x+ \omega)= \log \sin x +\frac{n \omega \cos x}{\sin x}-\frac{n \omega^2}{2 \sin^2 x}+\frac{n \omega^3 \cos x}{3 \sin^3 x}-\text{etc.},
\end{equation*}
where $n$ denotes the number, by which the hyperbolic must be multiplied that the propounded logarithms result. But if $y=\tan x$ and $z= \log \tan (x+ \omega)$, it will be $\frac{dy}{dx}= \frac{n}{\sin x \cdot \cos x}=\frac{2n}{ \sin 2x}$, $\frac{ddy}{2dx^2}=\frac{-2n \cos 2x}{(\sin 2x)^2}$ etc. and hence

\begin{equation*}
\log \tan (x+ \omega)= \log \tan x +\frac{2n \omega}{\sin 2x}-\frac{2n \omega^2 \cos 2x}{(\sin 2x)^2}+\text{etc.},
\end{equation*}
by means of which formulas the logarithms of sines and tangents can be interpolated.

\paragraph*{§100}

Let us put that $y$ denotes the arc, whose logarithm of the sine is $=x$, or that $y= \operatorname{A}. \log x $, and that $z$ is the arc, whose logarithm of the sine we want to put $=x+\omega$ or $z= \operatorname{A}.\log \sin (x+\omega)$; it will be $x= \log \sin y$ and 

\begin{equation*}
\frac{dx}{dy}=\frac{n \cos y}{\sin y}, \quad \text{whence} \quad \frac{dy}{dx}=\frac{\sin y}{n \cos y};
\end{equation*}
it will be

\begin{equation*}
\frac{ddy}{dx}=\frac{dy}{n \cos^2 y}=\frac{dx \sin y}{n^2 \cos^3 y}, \quad \text{therefore} \quad \frac{ddy}{dx^2}=\frac{\sin y}{n^2 \cos^3 y}.
\end{equation*}
As a logical consequence

\begin{equation*}
z=y +\frac{\omega \sin y}{n \cos y}+\frac{\omega^2 \sin y}{2n^2 \cos^3 y}+\text{etc.}
\end{equation*}
In like manner, given the logarithm of a cosine, the expression will be found. \\[2mm]
But if $y=\operatorname{A}.\log \tan x$ and $z=\operatorname{A}. \log (x+\omega)$, since $x= \log \tan y$, it will be

\begin{equation*}
\frac{dx}{dy}=\frac{n}{\sin y \cdot \cos y} \quad \text{and} \quad \frac{dy}{dx}=\frac{\sin y \cdot \cos y}{n}= \frac{\sin 2y}{2n},
\end{equation*}
whence

\begin{equation*}
\frac{ddy}{dx}=\frac{2dy \cos 2y}{2n}=\frac{dx \sin 2y \cdot \cos 2y}{2nn}
\end{equation*}
and

\begin{equation*}
\frac{ddy}{dx^2}=\frac{\sin 2y \cdot \cos 2y}{2nn}=\frac{\sin 4y}{4nn}, \quad \frac{d^3y}{dx^3}=\frac{\sin 2y \cdot 4y}{2n^3} \quad \text{etc.};
\end{equation*}
hence

\begin{equation*}
z=y+\frac{\omega \sin 2y}{2n}+\frac{\omega^2 \sin 2y \cdot \cos 2y}{4nn}+\frac{\omega^3 \sin 2y \cdot \cos 4y}{12n^3}+\text{etc.}
\end{equation*}

\paragraph*{§101}

Since the use of these expressions for the construction of tables of logarithms of sines and cosines can easily be seen from the preceding paragraphs, we will not treat this here any longer. Therefore, lastly let us consider the value $y=e^x \sin nx$ and let $z=e^{x+ \omega}\sin n(x+\omega)$; since

\begin{alignat*}{9}
&\frac{dy}{dx}&&=e^x (\sin nx+ n \cos nx) \\
&\frac{ddy}{dx^2}&&=e^x ((1-nn)\sin nx+2n \cos nx) \\
&\frac{d^3y}{dx^3}&&=e^x((1-3nn)\sin nx+ n(3-nn)\cos nx) \\
&\frac{d^4y}{dx^4}&&=e^x((1-6nn+n^4)\sin nx +n(4-4nn)\cos nx) \\
&\frac{d^5 y}{dx^5}&&= e^x((1-10nn+5n^4)\sin nx+n(5-10nn+n^4)\cos nx),
\end{alignat*}
\begin{center}
etc.;
\end{center}
having substituted these values and having divided by $e^x$, it will be

\begin{equation*}
e^{\omega}\sin n(x+\omega)=\sin nx
\end{equation*}
\begin{equation*}
+\omega \sin nx +\frac{1-nn}{2}\omega^2 \sin nx +\frac{1-3nn}{6}\omega^3 \sin nx +\frac{1-6nn+n^4}{24}\omega^4 \sin nx+ \text{etc.}
\end{equation*}
\begin{equation*}
+n \omega \cos nx +\frac{2n}{2}\omega^2 \cos nx +\frac{n(3-nn)}{6}\omega^3 \cos nx +\frac{n(4-4nn)}{24}\omega^4 \cos nx +\text{etc.}
\end{equation*}

\paragraph*{§102}

Hence many extraordinary corollaries can be deduced; but it suffices for us to have mentioned the following things here. If  it was $x=0$, it will be

\begin{equation*}
e^{\omega} \sin n \omega = n \omega +\frac{2n}{2}\omega^2 +\frac{n(3-nn)}{6}\omega^3 +\frac{n(4-nn)}{24}\omega^4 +\frac{n(5-10n^2+n^4)}{120}\omega^5+\text{etc.}
\end{equation*}
If $\omega =-x$, because of $\sin n(x+\omega)=0$, it will be

\begin{equation*}
\tan nx= \frac{nx -\frac{2n}{2}x^2+\frac{n(3-nn)}{6}x^3-\frac{n(4-4nn)}{24}x^4+\frac{n(5-10n^2+n^4)}{120}x^5-\text{etc.}}{1-x+\frac{1-nn}{2}x^2-\frac{1-3nn}{6}x^3+\frac{1-6nn+n^4}{24}x^4-\text{etc.}}
\end{equation*}
But in general, if $n=1$, one will have

\begin{equation*}
e^{\omega}\sin (x+\omega)= \sin x \left(1+\omega -\frac{1}{3}\omega^3-\frac{1}{6}\omega^4-\frac{1}{30}\omega^5+\frac{1}{630}\omega^7+\text{etc.}\right)
\end{equation*}
\begin{equation*}
+\omega \cos x \left(1+\omega +\frac{1}{3}\omega^2-\frac{1}{30}\omega^4-\frac{1}{90}\omega^5-\frac{1}{630}\omega^6+\text{etc.}\right)
\end{equation*}
But if it is $n=0$, because $\sin n(x+\omega)=n(x+\omega)$ and $\sin nx=nx$ and $\cos nx=1$, and if one divides by $n$ everywhere, we have

\begin{equation*}
e^{\omega}(x+\omega)=x+ \omega x+ \frac{1}{2}\omega^2 x+\frac{1}{6}\omega^3 x+ \frac{1}{24}\omega^4x+ \text{etc.}
\end{equation*}
\begin{equation*}
+\omega +\omega^2+ \frac{1}{2}\omega^3+\frac{1}{6}\omega^4+\frac{1}{24}\omega^5+\text{etc.},
\end{equation*}
the validity of which  equation is obvious.


\title{Investigation of the sum of series from the general Term
\footnote{
Original title: "'Investigatio summae serierum ex termino generali"', first published as part of the book"`\textit{Institutiones calculi differentialis cum eius usu in analysi finitorum ac doctrina serierum}, 1755"', reprinted in \textit{Opera Omnia}: Series 1, Volume 10, pp. 309 - 336, Eneström-Number E212, translated by: Alexander Aycock for the "`Euler-Kreis Mainz"'}}
\author{Leonhard Euler}
\date{ }
\maketitle

\paragraph*{§103}

Let the general term corresponding to the index $x$ of a certain series be $=y$, so that $y$ is a function of $x$. Further, let $Sy$ be the sum or the summatory  term of the series expressing the aggregate of all terms from the first or another fixed term up to $y$. Indeed, in the following we will compute the sum of the series starting from the first term, whence, if  $x=1$, $y$ will give the first term and $Sy$ will exhibit this first term $y$; but if one puts $x=0$, the summatory term $Sy$ has to go over into zero, because there are no terms to be summed. Therefore, the summatory term $Sy$ will be a function of $x$ vanishing for $x=0$.

\paragraph*{§104}

If the general term $y$ consists of several parts so that $y=p+q+r+\text{etc.}$, then one can consider the series  as conflated of several other series, whose general terms are $p$, $q$, $r$ etc. Hence, if the sums of these series are known, one will also be able to assign the sum of the propounded series; for, it will be the aggregate of all single series. Therefore, if $y=p+q+r+\text{etc.}$, it will be $Sy=Sp+Sq+Sq+\text{etc.}$ Therefore, because above we exhibited the sums of series, whose general terms are arbitrary positive powers of $x$ with positive integer coefficients, hence  one will be able to find the summatory term of any series, whose general term is $ax^{\alpha}+bx^{\beta}+cx^{\gamma}+\text{etc.}$, while $\alpha$, $\beta$, $\gamma$ etc. are positive integer numbers, or whose general term is a polynomial function of $x$.

\paragraph*{§105}

In this series whose general term or the term corresponding to the exponent $x$ is $=y$ let the term preceding this one or the term corresponding to the index $x-1$ be $=v$; since $v$ results from $y$, if  one writes $x-1$ instead of $x$, it will be

\begin{equation*}
v=y-\frac{dy}{dx}+\frac{ddy}{2dx^2}-\frac{d^3y}{6dx^3}+\frac{d^4y}{24dx^4}-\frac{d^5y}{120 dx^5}+\text{etc.}
\end{equation*}
Therefore, if $y$ was the general term of this series

\begin{alignat*}{9}
&1 \quad &&2 \quad && 3 \quad && 4 \cdots &&  \cdots ~ && x-1 \quad  && x\\
&a~+~ \quad &&b~+~ \quad &&c~+~ \quad &&d~+~ \quad && \cdots &&+~ v~+~  \quad &&y
\end{alignat*}
and the term corresponding to the index $0$ of this series was $=A$, $v$, since it is a function of $x$, will be the general term of this series

\begin{alignat*}{9}
&1 \quad &&2 \quad && 3 \quad &&4 \quad && 5 \cdots \cdots \cdot && x \\
&A~+~ \quad &&a~+~ \quad &&b ~+~ \quad &&c ~+~ \quad &&d~+~   \cdots ~+~&&v,
\end{alignat*}
whence, if $Sv$ denotes the sum of this series, it will be $Sv=Sy-y+A$. Therefore, having put $x=0$, since  $Sy=0$ and $y=A$, also $Sv$ will vanish.

\paragraph*{§106}

Therefore, since

\begin{equation*}
v=y -\frac{dy}{dx}+\frac{ddy}{2dx^2}-\frac{d^3y}{6dx^3}+\text{etc.},
\end{equation*}
 by means of the results demonstrated before, it will be

\begin{equation*}
Sv=Sy-S\frac{dy}{dx}+S\frac{ddy}{2dx^2}-S\dfrac{d^3y}{6dx^3}+S\frac{d^4 y}{24dx^4}-\text{etc.}
\end{equation*}
and, because of $Sv=Sy-y+A$, it will be

\begin{equation*}
y-A= S\frac{dy}{dx}-S\frac{ddy}{2dx^2}+S\frac{d^3 y}{6dx^3}-S\frac{d^4y}{24 dx^4}+\text{etc.}
\end{equation*}
and hence one will have

\begin{equation*}
S\frac{dy}{dx}=y-A+S\frac{ddy}{2dx^2}-S\frac{d^3y}{6dx^3}+S\frac{d^4y}{24dx^4}-\text{etc.}
\end{equation*}
Therefore, if one knows the summatory terms of the series, whose general terms are $\frac{ddy}{dx^2}$, $\frac{d^3y}{dx^3}$, $\frac{d^4y}{dx^4}$ etc.,  from them one will obtain the summatory term of the series, whose general term is $\frac{dy}{dx}$. The quantity $A$ has to be of such a nature that for $x=0$ the summatory $S\frac{dy}{dx}$ term vanishes, and using this condition, it is  determined more easily than if we would say that it is the term corresponding to the index $0$ in the series whose general term is $=y$.

\paragraph*{§107}

This ansatz is usually made to investigate the sums of the powers of natural numbers. For, let $y=x^{n+1}$; since 

\begin{equation*}
\frac{dy}{dx}=(n+1)x^n, \quad \frac{ddy}{2dx^2}=\frac{(n+1)n}{1 \cdot 2}x^{n-1}, \quad \frac{d^3y}{6dx^3}=\frac{(n+1)n(n-1)}{1 \cdot 2 \cdot 3}x^{n-2},
\end{equation*}
\begin{equation*}
\frac{d^4 y}{24 dx^4}=\frac{(n+1)n (n-1)(n-2)}{1 \cdot 2 \cdot 3 \cdot 4}x^{n-3} \quad \text{etc.},
\end{equation*}
having substituted these values, it will be

\begin{equation*}
(n+1)Sx^n=x^{n+1}-A+\frac{(n+1)n}{1 \cdot 2}Sx^{n-1}-\frac{(n+1)n(n-1)}{1 \cdot 2 \cdot 3}Sx^{n-2}+\text{etc.};
\end{equation*}
and if one divides by $n+1$ on both sides, it will be

\begin{equation*}
Sx^n=\frac{1}{n+1}x^{n+1}+\frac{n}{2}Sx^{n-1}-\frac{n(n-1)}{2 \cdot 3}Sx^{n-2}+\frac{n(n-1)(n-2)}{2 \cdot 3 \cdot 4}Sx^{n-3}-\text{etc.}-\text{Const.},
\end{equation*}
which constant has to be taken in such a way that for $x=0$ the summatory term vanishes. Therefore, using this formula, from the already known sums of lower powers, whose general terms are $x^{n-1}$, $x^{n-2}$ etc., one will be able to find the sum of the higher powers expressed by the general term $x^n$.

\paragraph*{§108}

If in this expression $n$ denotes a positive integer, the number of terms will be finite. And hence the sum of infinitely many  powers will be found explicitly; for, if $n=0$, it will be

\begin{equation*}
Sx^0=x.
\end{equation*}
And having known this one, it will be possible to proceed to the sums of higher powers; for, having put $n=1$, it will be

\begin{equation*}
Sx^1=\frac{1}{2}x^2+\frac{1}{2}Sx^0= \frac{1}{2}x^2+\frac{1}{2}x;
\end{equation*}
if one further sets $n=2$, this equation will result

\begin{equation*}
Sx^2=\frac{1}{3}x^3+Sx-\frac{1}{3}Sx^0= \frac{1}{3}x^3+\frac{1}{2}x^2+\frac{1}{6}x,
\end{equation*}
moreover,

\begin{equation*}
Sx^3= \frac{1}{4}x^4+\frac{3}{2}Sx^2-Sx+\frac{1}{4}Sx^0= \frac{1}{4}x^4+\frac{1}{2}x^3+\frac{1}{4}x^2,
\end{equation*} 
\begin{equation*}
Sx^4= \frac{1}{5}x^5+\frac{4}{2}Sx^3-\frac{4}{2}Sx^2+Sx-\frac{1}{5}Sx^0
\end{equation*}
or

\begin{equation*}
Sx^4= \frac{1}{5}x^5+\frac{1}{2}x^4+\frac{1}{3}x^3-\frac{1}{30}x.
\end{equation*}
And so forth, the successive sums of all higher powers are derived from the lower ones; but the same is achieved a lot easier as follows.

\paragraph*{§109}

Since we found above that

\begin{equation*}
S\frac{dy}{dx}= y+\frac{1}{2}S \frac{ddy}{dx^2}-\frac{1}{6}S\frac{d^3 y}{dx^3}+\frac{1}{24}S\frac{d^4 y}{dx^4}-\frac{1}{120}S\frac{d^5 y}{dx^5}+\text{etc.}
\end{equation*}
if we put $\frac{dy}{dx}=z$, it will be $\frac{ddy}{dx^2}=\frac{dz}{dx}$, $\frac{d^3y}{dx^3}=\frac{ddz}{dx^2}$ etc. But then, because of $dy=zdx$, $y$ will be the quantity whose differential is $=zdx$ which we denote by $y=\int zdx$. Although this way to find $y$ from  given $z$ depends on integral calculus, we will nevertheless be able to use this formula $\int zdx$ here, if  we substitute only  functions of $x$ of such a kind for $z$ that this function whose differential is $=zdx$ can be exhibited from the preceding ones. Therefore, having substituted these values, it will be

\begin{equation*}
Sz= \int zdx+\frac{1}{2}S\frac{dz}{dx}-\frac{1}{6}S\frac{ddz}{dx^2}+\frac{1}{24}S\frac{d^3 z}{dx^3}-\text{etc.},
\end{equation*} 
 adding such a constant that for $x=0$ the sum $Sz$ also vanishes.

\paragraph*{§110}

But by substituting the letter $z$ for $y$ in the above expression or, which is the same, by differentiating this equation it will be

\begin{equation*}
S \frac{dz}{dx}=z+\frac{1}{2}S\frac{ddz}{dx^2}-\frac{1}{6}S\frac{d^3 z}{dx^3}+\frac{1}{24}S\frac{d^4z}{dx^4}-\text{etc.};
\end{equation*}
but if one writes $\frac{dz}{dx}$ instead of $y$, it will be

\begin{equation*}
S\frac{ddz}{dx^2}=\frac{dz}{dx}+\frac{1}{2}S \frac{d^3z}{dx^3}-\frac{1}{6}S \frac{d^4 z}{dx^4}+\frac{1}{24}S\frac{d^5 z}{dx^5}-\text{etc.}
\end{equation*}
But if in like manner  one successively substitutes the values $\frac{ddz}{dx^2}$, $\frac{d^3z}{dx^3}$ etc. for $y$, one will find

\begin{alignat*}{9}
&S\frac{d^3 z}{dx^3}&&=\frac{ddz}{dx^2}&&+\frac{1}{2}S\frac{d^4z}{dx^4}&&-\frac{1}{6}S\frac{d^5 z}{dx^5}&&+\frac{1}{24}S\frac{d^6 z}{dx^6}&&-\text{etc.} \\[3mm]
&S\frac{d^4 z}{dx^4}&&=\frac{d^3z}{dx^3}&&+\frac{1}{2}S\frac{d^5z}{dx^5}&&-\frac{1}{6}S\frac{d^6 z}{dx^6}&&+\frac{1}{24}S\frac{d^7 z}{dx^7}&&-\text{etc.} 
\end{alignat*}
and so forth to infinity.

\paragraph*{§111}

If now these values are successively substituted for $S\frac{dz}{dx}$, $S\frac{ddz}{dx^2}$, $S\frac{d^3z}{dx^3}$ etc. in the expression

\begin{equation*}
Sz =\int zdx+ \frac{1}{2}S \frac{dz}{dx}-\frac{1}{6}S \frac{ddz}{dx^2}+\frac{1}{24}S\frac{d^3z}{dx^3}-\text{etc.},
\end{equation*}
one will find an expression for $Sz$ consisting of these terms $\int zdx$, $z$, $\frac{dz}{dx}$, $\frac{ddz}{dx^2}$, $\frac{d^3z}{dx^3}$ etc., whose coefficients are investigated more easily the following way. Put

\begin{equation*}
Sz = \int zdx +\alpha z +\frac{\beta dz}{dx}+\frac{\gamma ddz}{dx^2}+\frac{\delta d^3 z}{dx^3}+\frac{\varepsilon d^4 z}{dx^4}+\text{etc.}
\end{equation*}
and for these terms substitute their values  they obtain from the preceding series, from which

\begin{equation*}
\renewcommand{\arraystretch}{2,5}
\setlength{\arraycolsep}{0.0mm}
\begin{array}{rcccccccccccccccccccccccccccccccc}
 \int zdx &~=~&Sz &~-~&\dfrac{1}{2}& S\dfrac{dz}{dx}&~+~&\dfrac{1}{6}&S \dfrac{ddz}{dx^2}&~-~&\dfrac{1}{24}&S \dfrac{d^3z}{dx^3}&~+~&\dfrac{1}{120}&S \dfrac{d^4 z}{dx^4}&~-~&\text{etc.}\\
 \alpha z &~=~& &~+~&\alpha & S \dfrac{dz}{dx}&~-~&\dfrac{\alpha}{2}&S\dfrac{ddz}{dx^2}&~+~& \dfrac{\alpha}{6}& S \dfrac{d^3z}{dx^3}&~-~&\dfrac{\alpha}{24}& S\dfrac{d^4z}{dx^4}&~+~&\text{etc.} \\
\dfrac{\beta dz}{dx} &~=~& &&  &  &~+~&\beta &S\dfrac{ddz}{dx^2}&~-~&\dfrac{\beta}{2}& S \dfrac{d^3z}{dx^3}&~+~&\dfrac{\beta}{6}& S\dfrac{d^4z}{dx^4}&~+~&\text{etc.} \\
\dfrac{\gamma ddz}{dx^2} &~=~ &  &&  && && &~+~&\gamma &  S \dfrac{d^3z}{dx^3}&~-~&\dfrac{\gamma}{2}& S\dfrac{d^4z}{dx^4}&~+~&\text{etc.} \\
\dfrac{\delta d^3z}{dx^3} &~=~& &&  &&  && &&  &&~+~&\delta  & S\dfrac{d^4z}{dx^4}&~+~&\text{etc.} \\
& && &&  && &&\text{etc.}
\end{array}
\end{equation*}
Since, having added all these values, they have to produce $Sz$, the coefficients $\alpha$, $\beta$, $\gamma$, $\delta$ etc. will be defined from the following equations

\begin{equation*}
\alpha -\frac{1}{2}=0, \quad \beta -\frac{\alpha }{2}+\frac{1}{6}=0, \quad \gamma -\frac{\beta}{2}+\frac{\alpha}{6}-\frac{1}{24}=0,
\end{equation*}
\begin{equation*}
\delta -\frac{\gamma}{2}+\frac{\beta}{6}-\frac{\alpha}{24}+\frac{1}{120}=0, \quad \varepsilon -\frac{\delta}{2}+\frac{\beta}{24}+\frac{\alpha}{120}-\frac{1}{720}=0,
\end{equation*}
\begin{equation*}
\zeta -\frac{\varepsilon}{2}+\frac{\delta}{6}-\frac{\gamma}{24}+\frac{\beta}{120}-\frac{\alpha}{720}+\frac{1}{5040}=0 \quad \text{etc.}
\end{equation*}
and continuing this way, one will  find terms each second of which vanishes. Therefore, the third, fifth, seventh letter and in general all odd ones will be $=0$ except for the first, which seems to violate the law of continuity. Therefore, it is even more necessary to prove rigorously that all odd terms except for the first necessarily vanish.

\paragraph*{§113}

Since each letter is determined according to a constant law from the preceding ones, they will constitute a recurring series. To make this explicit, assume this series

\begin{equation*}
1+ \alpha u+ \beta u^2+ \gamma u^3+ \delta u^4+ \varepsilon u^5+ \zeta u^6+ \text{etc.},
\end{equation*}
whose value we want to put $=V$, and it is obvious that this recurring series results from the expansion of this fraction

\begin{equation*}
V=\frac{1}{1-\frac{1}{2}u+\frac{1}{6}u^2-\frac{1}{24}u^3+\frac{1}{240}u^4-\text{etc.}}
\end{equation*}
And if this fraction   can be resolved into a power series in $u$ in another way, it is necessary that always the same series 

\begin{equation*}
V=1+\alpha u+\beta u^2+\gamma u^3+\delta u^4+\varepsilon u^5+\text{etc.}
\end{equation*}
results; and this way another law, by which the same values $\alpha$, $\beta$, $\gamma$, $\delta$ etc. are determined, will be found.

\paragraph*{§114}

Since, if $e$ denotes the number, whose hyperbolic logarithm is equal to $1$, it will be

\begin{equation*}
e^{-u}=1-u+\frac{1}{2}u^2-\frac{1}{6}u^3+\frac{1}{24}u^4-\frac{1}{120}u^5+\text{etc.},
\end{equation*}
it will also be

\begin{equation*}
\frac{1-e^{-u}}{u}=1-\frac{1}{2}u+\frac{1}{6}u^2-\frac{1}{24}u^3+\frac{1}{120}u^4-\text{etc.}
\end{equation*}
and hence

\begin{equation*}
V=\frac{u}{1-e^{-u}}.
\end{equation*}
Now cancel the second term $\alpha u=\frac{1}{2}u$ from the series that

\begin{equation*}
V-\frac{1}{2}u=1+\beta u^2+\gamma u^3 +\delta u^4+ \varepsilon u^5 +\zeta u^6+\text{etc.};
\end{equation*}
it will be

\begin{equation*}
V-\frac{1}{2}u=\frac{\frac{1}{2}u(1+e^{-u})}{1-e^{-u}}.
\end{equation*}
Multiply the numerator and denominator by $e^{\frac{1}{2}u}$ and it will be

\begin{equation*}
V-\frac{1}{2}u=\frac{u\left(e^{\frac{1}{2}u}+e^{-\frac{1}{2}u}\right)}{\left(e^{\frac{1}{2}u}-e^{-\frac{1}{2}u}\right)}
\end{equation*}
and, having converted the quantities $e^{\frac{1}{2}u}$ and $e^{-\frac{1}{2}u}$ into series, it will be

\begin{equation*}
V-\frac{1}{2}u=\frac{1+\frac{u^2}{2 \cdot 4}+\frac{u^4}{2 \cdot 4 \cdot 6 \cdot 8}+\frac{u^6}{2 \cdot 4 \cdot 6 \cdot 8 \cdot 10 \cdot 12}+\text{etc.}}{2\left(\frac{1}{2}+\frac{u^2}{2 \cdot 4 \cdot 6}+\frac{u^4}{2 \cdot 4 \cdot 6 \cdot 8 \cdot 10}+\text{etc.}\right)}
\end{equation*}
or

\begin{equation*}
V-\frac{1}{2}u=\frac{1+\frac{u^2}{2 \cdot 4}+\frac{u^4}{2 \cdot 4 \cdot 6 \cdot 8}+\frac{u^8}{2 \cdot 4 \cdots 12}+\frac{u^8}{2 \cdot 4 \cdots 16}+\text{etc.}}{1+\frac{u^2}{4 \cdot 6}+\frac{u^4}{4 \cdot 6 \cdot 8 \cdot 10}+\frac{u^8}{4 \cdot 6 \cdots 14}+\frac{u^8}{4 \cdot 6 \cdots 18}+\text{etc.}}
\end{equation*}

\paragraph*{§115}

Therefore, because in this fraction the odd powers are completely missing, its power series will contain no odd powers at all; therefore, because $V-\frac{1}{2}u$ becomes equal to this series

\begin{equation*}
1+\beta u^2+\gamma u^3+\delta u^4+\varepsilon u^5+\zeta u^6+\text{etc.},
\end{equation*}
the coefficients of the odd powers, $\gamma$, $\varepsilon$, $\eta$, $\iota$ etc. will all vanish. And this is the reason why in the series $1+\alpha u+\beta u^2+\gamma u^3+\delta u^4+\text{etc.}$ each second terms except for the first are $=0$ and the law of continuity is nevertheless not violated. Therefore, it will be

\begin{equation*}
V=1+\frac{1}{2}u+\beta u^2+\delta u^4+\zeta u^6+\theta u^8+ \varkappa x^{10}+\text{etc.}
\end{equation*}
and, having determined the letters $\beta$, $\delta$, $\zeta$, $\theta$, $\varkappa$ etc. by expansion of the fraction above, we will obtain the summatory term $Sz$ of the series whose general term corresponding to the index $x$ is $=z$ expressed this way

\begin{equation*}
Sz= \int zdx+\frac{1}{2}z+\frac{\beta dz}{dx}+\frac{\delta d^3z}{dx^3}+\frac{\zeta d^5 z}{dx^5}+\frac{\theta d^7z }{dx^7}+\text{etc.}
\end{equation*} 

\paragraph*{§116}

Since the series $1+\beta u^2 +\delta u^4 +\zeta u^6+\theta u^8+\text{etc.}$ results from the expansion of this fraction

\begin{equation*}
\frac{1+\frac{u^2}{2 \cdot 4}+\frac{u^4}{2 \cdot 4 \cdot 6 \cdot 8}+\frac{u^6}{2 \cdot 4 \cdot 6 \cdot 8 \cdot 10 \cdot 12}+\text{etc.}}{1+\frac{u^2}{4 \cdot 6}+\frac{u^4}{4 \cdot 6 \cdot 8 \cdot 10}+\frac{u^8}{4 \cdot 6 \cdot 8 \cdot 10 \cdot 12 \cdot 14}+\text{etc.}},
\end{equation*}
the letters $\beta$, $\delta$, $\zeta$, $\theta$ etc. will obey the following law

\begin{alignat*}{9}
&\beta &&= \frac{\quad \quad ~~   1 \quad   \quad ~~}{ 2 \cdot 4}&&- \frac{1}{4 \cdot 6}\\[3mm]
&\gamma &&= \frac{~~~~~~~~~1~~~~~~~~~}{2 \cdot 4 \cdot 6 \cdot 8}&&-\frac{\beta}{4 \cdot 6}&&-\frac{1}{4 \cdot 6 \cdot 8 \cdot 10}\\[3mm]
&\delta &&= \frac{1}{2 \cdot 4 \cdot 6 \cdots 12}&&-\frac{\delta}{4 \cdot 6}&&- \frac{\beta}{4 \cdot 6 \cdot 8 \cdot 10}&&- \frac{1}{4 \cdot 6 \cdots 11}\\[3mm]
&\theta &&=\frac{1}{2 \cdot 4 \cdot 6 \cdots 16}&&-\frac{\zeta}{4 \cdot 6}&&-\frac{\delta}{4 \cdot 6 \cdot 8 \cdot 10}&&-\frac{\beta}{4 \cdot 6 \cdots 14}&&- \frac{1}{4 \cdot 6 \cdots 18}\\
& && && && \text{etc.}
\end{alignat*}
But these values are  positive and negative alternately.

\paragraph*{§117}

Therefore, if each second of these letters is assumed to be negative so that

\begin{equation*}
Sz =\int zdx +\frac{1}{2}z-\frac{\beta dz}{dx}+\frac{\delta d^3 z}{dx^3}-\frac{\zeta d^5 z}{dx^5}+\frac{\theta d^7 z}{dx^7}-\text{etc.},
\end{equation*}
the letters $\beta$, $\delta$, $\zeta$, $\theta$ can be defined from this fraction

\begin{equation*}
\frac{1-\frac{u^2}{2 \cdot 4}+\frac{u^4}{2 \cdot 4 \cdot 6 \cdot 8}-\frac{u^6}{2 \cdot 4 \cdot 6 \cdot 8 \cdot 10 \cdot 12}-\text{etc.}}{1-\frac{u^2}{4 \cdot 6}+\frac{u^4}{4 \cdot 6 \cdot 8 \cdot 10}-\frac{u^8}{4 \cdot 6 \cdot 8 \cdot 10 \cdot 12 \cdot 14}-\text{etc.}},
\end{equation*}
by expanding it into the series

\begin{equation*}
1+\beta u^2+\delta u^4+\zeta u^6+\theta u^8+\text{etc.};
\end{equation*}
therefore, it will be

\begin{alignat*}{9}
&\beta &&= \frac{1}{4 \cdot 6}&&-\frac{\quad  \quad 1 \quad \quad ~}{2 \cdot 4}\\[3mm]
&\delta &&= \frac{\beta }{4 \cdot 6}&&-\frac{1}{4 \cdot 6 \cdot 8 \cdot 10}&&+\frac{1}{2 \cdot 4 \cdot 6 \cdot 8}\\[3mm]
&\zeta &&= \frac{\delta}{4 \cdot 6}&&-\frac{\beta}{4 \cdot 6 \cdot 8 \cdot 10}&&+\frac{1}{4 \cdot 6 \cdots 14}&&-\frac{1}{2 \cdot 4 \cdots 12}\\
& && &&+\text{etc.};
\end{alignat*}
but now all terms will become negative.

\paragraph*{§118}

Therefore, let us put $\beta =-A$, $\delta =-B$, $\zeta =-C$ etc. such that 

\begin{equation*}
Sz =\int zdx +\frac{1}{2}z+\frac{Adz}{dx}-\frac{Bd^3z}{dx^3}+\frac{Cd^5z}{dx^5}-\frac{Dd^7 z}{dx^7}+\text{etc.},
\end{equation*}
and to define the letters $A$, $B$, $C$, $D$ etc., consider this series

\begin{equation*}
1-Au^2-Bu^4-Cu^6-Du^8-Eu^{10}-\text{etc.},
\end{equation*}
which results from the expansion of this fraction

\begin{equation*}
\frac{1-\frac{u^2}{2 \cdot 4}+\frac{u^4}{2 \cdot 4 \cdot 6 \cdot 8}-\frac{u^6}{2 \cdot 4 \cdots 12}+\frac{u^8}{2 \cdot 4 \cdots 16}-\text{etc.}}{1-\frac{u^2}{4 \cdot 6}+\frac{u^4}{4 \cdot 6 \cdot 8 \cdot 10}-\frac{u^6}{4 \cdot 6 \cdots 10}+\frac{u^8}{4 \cdot 6 \cdots 18}-\text{etc.}},
\end{equation*}
or consider this series

\begin{equation*}
\frac{1}{u}-Au-Bu^3 -Cu^5-Du^7-Eu^9-\text{etc.}=s,
\end{equation*}
which results from the expansion of this fraction

\begin{equation*}
s=\frac{1-\frac{u^2}{2 \cdot 4}+\frac{u^4}{2 \cdot 4 \cdot 6 \cdot 8}-\frac{u^6}{2 \cdot 4 \cdots 12}+\text{etc.}}{u-\frac{u^3}{4 \cdot 6}+\frac{u^5}{4 \cdot 6 \cdot 8 \cdot 10}-\frac{u^7}{4 \cdot 6 \cdots 14}+\text{etc.}}.
\end{equation*}
But since

\begin{equation*}
\cos \frac{1}{2}u=1-\frac{u^2}{2 \cdot 4}+\frac{u^4}{2 \cdot 4 \cdot 6 \cdot 8}-\frac{u^6}{2 \cdot 4 \cdots 12}+\text{etc.},
\end{equation*}
\begin{equation*}
\sin \frac{1}{2}u= \frac{u}{2}-\frac{u^3}{2 \cdot 4 \cdot 6}+\frac{u^5}{2 \cdot 4 \cdot 6 \cdot 8 \cdot 10}-\frac{u^7}{2 \cdot 4 \cdots 14}+\text{etc.},
\end{equation*}
it follows that

\begin{equation*}
s=\frac{\cos \frac{1}{2}u}{2 \sin \frac{1}{2}u}=\frac{1}{2}\cot \frac{1}{2}u.
\end{equation*}
Therefore, if the cotangent of the arc $\frac{1}{2}u$ is converted into a series whose terms are ascending  powers of $u$, from it one will find the values of the letters $A$, $B$, $C$, $D$, $E$ etc.

\paragraph*{§119}

Therefore, since $s=\frac{1}{2}\cot \frac{1}{2}u$, it will be $\frac{1}{2}u= \operatorname{arccot} 2s$ and, by differentiating this equation, it will be $\frac{1}{2}du=\frac{-2ds}{1+4ss}$ or $4ds+du+4ssdu=0$ or

\begin{equation*}
\frac{4ds}{du}+1+4ss=0.
\end{equation*}
But since

\begin{equation*}
s=\frac{1}{u}-Au-Bu^3-Cu^5-\text{etc.},
\end{equation*}
it will be

\begin{equation*}
\renewcommand{\arraystretch}{2,5}
\setlength{\arraycolsep}{0.0mm}
\begin{array}{rrrrrrrrrrrrrrrrrrrrrrrrrr}
\dfrac{4ds}{du}&~=~&-\dfrac{4}{uu}&~-~&4A &~-~&3 \cdot 4 Bu^2&~-~& 5 \cdot 4 Cu^4 &  ~-~& 7\cdot 4 Du^6 &~-~&\text{etc.} \\
1 &~=~&  &&1 \\
4ss &~ =~& \dfrac{4}{uu}&~-~&8A &~-~& 8Bu^2 &~-~& 8Cu^4 &~-~& 8Du^6 &~-~&\text{etc.} \\
&   && & & ~+~&4 A^2u^2 &~+~&8ABu^4 &~ +~&8 ACu^6 &~+~&\text{etc.} \\
& && && & & & &~+~&4 BBu^6 &~+~&\text{etc.}
\end{array}
\end{equation*}
Having set these homogeneous terms equal to zero, it will be

\begin{equation*}
A=\frac{1}{12}, \quad B=\frac{A^2}{5}, \quad C=\frac{2AB}{7}, \quad D=\frac{2AC+BB}{9}, \quad E=\frac{2AD+2BD}{11},
\end{equation*}
\begin{equation*}
F=\frac{2AE+2BD+CC}{13}, \quad G=\frac{2AF+2BE+2CD}{15}, \quad H=\frac{2AG+2BF+2CE+DD}{17},
\end{equation*}
\begin{equation*}
\text{etc.}
\end{equation*}
From these formulas it now  obviously follows that each value is positive.

\paragraph*{§120}

But since the denominators of these values become immensely large and impede the calculation quite a lot, instead of the letters $A$, $B$, $C$, $D$ etc. let us introduce these new letters

\begin{equation*}
A=\frac{\alpha}{1 \cdot 2 \cdot 3}, \quad B=\frac{\beta}{1 \cdot 2 \cdot 3 \cdot 4 \cdot 5}, \quad C=\frac{\gamma}{1 \cdot 2 \cdot 3 \cdots 7},
\end{equation*}
\begin{equation*}
D=\frac{\delta}{1 \cdot 2 \cdot 3 \cdots 9}, \quad E=\frac{\varepsilon}{1 \cdot 2 \cdot 3 \cdots 11} \quad \text{etc.}
\end{equation*}
And one will find 

\begin{equation*}
\alpha =\frac{1}{2}, \quad \beta=\frac{2}{3}\alpha^2, \quad \gamma =2 \cdot \frac{2}{3}\alpha \beta, \quad \delta =2 \cdot \frac{4}{3} \alpha \gamma +\frac{8 \cdot 7}{4 \cdot 5} \beta^2,
\end{equation*}
\begin{equation*}
\varepsilon=2 \cdot \frac{5}{3}\alpha \delta +2 \cdot \frac{10 \cdot 9 \cdot 8}{1 \cdot 2 \cdots 5}\beta \gamma, \quad \zeta= 2 \cdot \frac{12}{1 \cdot 2 \cdot 3}\alpha \varepsilon +2 \frac{12 \cdot 11 \cdot 10}{1 \cdot 2 \cdots 5}\beta \delta +\frac{12 \cdot 11 \cdot 10 \cdot 9 \cdot 8}{1 \cdot 2 \cdots 7}\gamma \gamma,
\end{equation*}
\begin{equation*}
\eta = 2 \cdot \frac{14}{1 \cdot 2 \cdot 3}\alpha \zeta +2 \cdot \frac{14 \cdot 13 \cdot 12}{1 \cdot 2 \cdots 5}\beta \varepsilon+2 \cdot \frac{14 \cdot 13 \cdot 12 \cdot 11 \cdot 10}{1 \cdot 2 \cdots 7}\gamma \delta
\end{equation*}
\begin{equation*}
\text{etc.}
\end{equation*}

\paragraph*{§121}
But we will  apply these formulas more conveniently

\begin{equation*}
\alpha =\frac{1}{2}, \quad \beta =\frac{4}{3}\cdot \frac{\alpha \alpha}{2}, \quad \gamma =\frac{6}{3}\cdot \alpha \beta, \quad \delta= \frac{8}{3}\cdot \alpha \gamma+ \frac{8 \cdot 7 \cdot 6}{3 \cdot 4 \cdot 5}\cdot \frac{\beta \beta }{2},
\end{equation*}
\begin{equation*}
\varepsilon =\frac{10}{3}\cdot \alpha \delta +\frac{10 \cdot 9 \cdot 8}{3 \cdot 4 \cdot 5}\cdot \beta \gamma, \quad \zeta =\frac{12}{3}\cdot \alpha \varepsilon +\frac{12 \cdot 11 \cdot 10 }{3 \cdot 4 \cdot 5}\cdot \beta \delta+\frac{12 \cdot 11 \cdot 10 \cdot 9 \cdot 8}{3 \cdot 4 \cdot 5 \cdot 6 \cdot 7}\cdot \frac{\gamma \gamma}{2},
\end{equation*}
\begin{equation*}
\eta = \frac{14}{3} \cdot \alpha \zeta +\frac{14 \cdot 13 \cdot 12}{3 \cdot 4 \cdot 5}\cdot \beta \varepsilon +\frac{14 \cdot 13 \cdot 12 \cdot 11 \cdot 10}{3 \cdot 4 \cdot 5 \cdot 6 \cdot 7}\cdot \gamma \delta,
\end{equation*}
\begin{equation*}
\theta = \frac{16}{3}\cdot \alpha \eta + \frac{16 \cdot 15 \cdot 14}{3 \cdot 4 \cdot 5}\cdot \beta \zeta +\frac{16 \cdot 15 \cdots 12}{3 \cdot 4 \cdots 7}\gamma \varepsilon +\frac{16 \cdot 15 \cdots 10}{3 \cdot 4 \cdots 9}\cdot \frac{\delta \delta}{2}
\end{equation*}
\begin{equation*}
\text{etc.}
\end{equation*}
Using this law, which simplifies the calculation a lot, if the values of the letters $\alpha$, $\beta$, $\gamma$, $\delta$ etc. were found, then the summatory term of any arbitrary series whose general term or the term corresponding to the index $x$ was $=z$, will be expressed as follows

\begin{equation*}
Sz= \int zdx+\frac{1}{2}z+\frac{\alpha dz}{1 \cdot 2 \cdot 3 \cdot dx}- \frac{\beta d^3 z}{1 \cdot 2 \cdot 3 \cdot 4 \cdot 5 dx^5}+\frac{\gamma d^5 z}{1 \cdot 2 \cdots 7 dx^5}
\end{equation*}
\begin{equation*}
-\frac{\delta d^7 z}{1 \cdot 2 \cdots 9 dx^7}+\frac{\varepsilon d^9 z}{1 \cdot 2 \cdots 11 dx^9}-\frac{\zeta d^{11}z}{1 \cdot 2 \cdots 13 dx^{11}}+\text{etc.}
\end{equation*}
But these letters $\alpha$, $\beta$, $\gamma$, $\delta$ etc. were found to have the following values:

\begin{alignat*}{18}
\alpha & =\frac{1}{2} \qquad \text{or} \qquad  & 1 \cdot 2 \alpha &=1  \\
\beta & =\frac{1}{6} \qquad & 1 \cdot 2 \cdot 3 \beta &=1 \\
\gamma &=\frac{1}{6} \qquad & 1 \cdot 2 \cdot 3 \cdot 4 \gamma &=4\\
\delta & =\frac{3}{10}  \qquad & 1 \cdot 2 \cdot 3 \cdots 5 \delta &=36 \\
\varepsilon &= \frac{5}{6} \qquad & 1 \cdot 2 \cdot 3 \cdots 6 \varepsilon &= 600 \\
\zeta & =\frac{691}{210} \qquad & 1 \cdot 2 \cdot 3 \cdots 7 \zeta &= 24 \cdot 691 \\
\eta & = \frac{35}{2} \qquad & 1 \cdot 2 \cdot 3 \cdots 8 \eta &= 20160 \cdot 35 \\
\theta &= \frac{3617}{30} \qquad & 1 \cdot 2 \cdot 3 \cdots 9 \theta &= 12096 \cdot 3617 \\
\iota &= \frac{43867}{42} \qquad & 1 \cdot 2 \cdot 3 \cdots 10 \iota &= 86400 \cdot 43867 \\
\varkappa &= \frac{1222277}{110} \qquad & 1 \cdot 2 \cdot 3 \cdots 11 \varkappa &= 362880 \cdot 1222277 \\
\lambda & =\frac{854513}{6} \qquad & 1 \cdot 2 \cdot 3 \cdots 12 \lambda &= 79833600 \cdot 854513 \\
\mu & =\frac{1181820455}{546} \qquad & 1 \cdot 2 \cdot 3 \cdots 13 \mu &= 11404800 \cdot 1181820455 \\
\nu &= \frac{76977927}{2} \qquad & 1 \cdot 2 \cdot 3 \cdots 14 \nu &= 43589145600 \cdot 76977927 \\
\xi &= \frac{23749461029}{30} \qquad & 1 \cdot 2 \cdot 3 \cdots 15 \xi &= 43589145600 \cdot 23749461029 \\
\pi &= \frac{8615841276005}{462} \qquad & 1 \cdot 2 \cdot 3 \cdots 16 \pi &= 45287424000 \cdot 8615841276005 \\ 
\end{alignat*}

\paragraph*{§122}

These numbers have the greatest use throughout the whole doctrine of series. For, first using these numbers one can  form the last terms in the sums of the even powers, on which we remarked above [§ 63 of the first part] that they cannot be found in the same way as the remaining terms from the sums of the preceding power. For, in the even powers the last terms containing $x$ of the sums are multiplied by certain numbers which numbers for the powers II, IV, VI, VII etc. are $\frac{1}{6}$, $\frac{1}{30}$, $\frac{1}{42}$, $\frac{1}{30}$ etc. affected with alternating signs. But these numbers result, if the values of the letters $\alpha$, $\beta$, $\gamma$, $\delta$ etc. found above are respectively divided by the odd numbers $3$, $5$, $7$, $9$ etc. whence these numbers, which  are usually called Bernoulli numbers after its discoverer Jacob Bernoulli, will be

\begin{alignat*}{9}
&\frac{\alpha}{3}&&=~~~\frac{1}{6}&&= \mathfrak{A}, \qquad \qquad &&\frac{\iota}{19}&&=~~~~~~~\frac{43867}{798}&& =\mathfrak{I}\\[3mm]
&\frac{\beta}{5}&&=~~\frac{1}{30}&&= \mathfrak{B}, \qquad \qquad &&\frac{\varkappa}{21}&&=~~~~~~\frac{174611}{330}&& =\mathfrak{K}=\frac{283 \cdot 617}{330}\\[3mm]
&\frac{\gamma}{7}&&=~~\frac{1}{42}&&= \mathfrak{C}, \qquad \qquad &&\frac{\lambda}{21}&&=~~~~~~\frac{854513}{138}&& =\mathfrak{L}= \frac{11 \cdot 131 \cdot 593}{2 \cdot 3 \cdot 23}\\[3mm]
&\frac{\alpha}{9}&&=~~\frac{1}{30}&&= \mathfrak{D}, \qquad \qquad &&\frac{\mu}{25}&&=~~~~\frac{236364091}{2730}&& =\mathfrak{M}\\[3mm]
&\frac{\varepsilon}{11}&&=~~\frac{5}{66}&&= \mathfrak{E}, \qquad \qquad &&\frac{\nu}{27}&&=~~~~~\frac{8553103}{6}&& =\mathfrak{N}= \frac{13 \cdot 657931}{6}\\[3mm]
&\frac{\zeta}{13}&&=\frac{691}{2730}&&= \mathfrak{F}, \qquad \qquad &&\frac{\xi}{27}&&=~~\frac{23749461029}{870}&& =\mathfrak{D}\\[3mm]
&\frac{\eta}{15}&&=~~~\frac{7}{6}&&= \mathfrak{G}, \qquad \qquad &&\frac{\pi}{31}&&=\frac{8615841276005}{14322}&& =\mathfrak{P}\\[3mm]
&\frac{\theta}{17}&&=\frac{3617}{510}&&= \mathfrak{H}, \qquad \qquad &&\text{etc.}
\end{alignat*} 

\paragraph*{§123}

Therefore, one will be able to find these Bernoulli numbers $\mathfrak{A}$, $\mathfrak{B}$, $\mathfrak{C}$ etc. immediately from the following equations

\begin{alignat*}{19}
&\mathfrak{A}&&= && ~~~\frac{1}{6}\\[3mm]
&\mathfrak{B}&&=~~\frac{4 \cdot 3}{1 \cdot 2}&&\cdot ~\frac{1}{5} \mathfrak{A}^2\\[3mm]
&\mathfrak{C}&&=~~\frac{6 \cdot 5}{1 \cdot 2}&&\cdot ~\frac{2}{7}\mathfrak{AB}\\[3mm]
&\mathfrak{D}&&= ~~\frac{8 \cdot 7}{1 \cdot 2}&&\cdot ~\frac{2}{9}\mathfrak{A C}&&+~~~\frac{8 \cdot 7 \cdot 6 \cdot 5}{1 \cdot 2 \cdot 3 \cdot 4}&& \cdot ~\frac{1}{9}\mathfrak{B}^2 \\[3mm]
&\mathfrak{E}&&=~\frac{10 \cdot 9}{1 \cdot 2}&&\cdot \frac{2}{11}\mathfrak{A D}&&+~\frac{12 \cdot 11 \cdot 10 \cdot 9}{1 \cdot 2 \cdot 3 \cdot 4}&&\cdot \frac{2}{13}\mathfrak{BD}&&+~~\frac{12 \cdot 11 \cdot 10 \cdot 9 \cdot 8 \cdot 7}{1 \cdot 2 \cdot 3 \cdot 4 \cdot 5 \cdot 6}&&\cdot \frac{1}{13}\mathfrak{C}^2\\[3mm]
&\mathfrak{G}&&= \frac{14 \cdot 13}{1 \cdot 2}&&\cdot \frac{2}{15}\mathfrak{AF}&&+\frac{14 \cdot 13 \cdot 12 \cdot 11}{1 \cdot 2 \cdot 3 \cdot 4}&&\cdot \frac{2}{15}\mathfrak{BE}&&+\frac{14 \cdot 13 \cdot 12 \cdot 12 \cdot 11 \cdot 9}{1 \cdot 2 \cdot 3  \cdot 4 \cdot 5 \cdot 6}&&\cdot \frac{2}{15}\mathfrak{CD}\\
& && && &&\text{etc.},
\end{alignat*}
the structure of which equations is clear per se, if one only notes, that where the square of a certain letter occurs, its coefficient is half as small as it seems to have to be according to the rule. But, the terms containing the products of different letters are to be considered to occur twice; for  example it will be

\begin{equation*}
13 \mathfrak{F}=\frac{12 \cdot 11}{1 \cdot 2}\mathfrak{A E}+\frac{12 \cdot 11 \cdot 10 \cdot 9}{1 \cdot 2 \cdot 3 \cdot 4}\mathfrak{BD}+\frac{12 \cdot 11 \cdot 10 \cdot 9 \cdot 8 \cdot 7}{1 \cdot 2 \cdot 3 \cdot 4 \cdot 5 \cdot 6}\mathfrak{CC}
\end{equation*} 
\begin{equation*}
+\frac{12 \cdot 11 \cdot 10 \cdots 5}{1 \cdot 2 \cdot 3 \cdots 8}\mathfrak{DB}+\frac{12 \cdot 11 \cdot 10 \cdots 3}{1 \cdot 2 \cdot 3 \cdots 10}\mathfrak{EA}.
\end{equation*}

\paragraph*{§124}

Further, the same numbers $\alpha$, $\beta$, $\gamma$, $\delta$ etc. enter the expressions of the sums of the series of fractions contained in this general form

\begin{equation*}
1+\frac{1}{2^n}+\frac{1}{3^n}+\frac{1}{4^n}+\frac{1}{5^n}+\frac{1}{6^n}+\text{etc.},
\end{equation*}
if $n$ is a positive even number. For, we gave these sums  expressed in terms of the half of the circumference of the circle $\pi$ whose radius is $=1$ in the \textit{Introductio} and  these numbers $\alpha$, $\beta$, $\gamma$, $\delta$ etc. are detected to enter  the coefficients of these powers. But to understand that this does not happen accidentally but has to happen, let us investigate the same sums in a special way so that the structure of those sums will become clear more easily. Since we found above (§ 43) that

\begin{equation*}
\frac{\pi}{n}\cot \frac{m}{n}\pi=\frac{1}{m}-\frac{1}{n-m}+\frac{1}{n+m}-\frac{1}{2n-m}+\frac{1}{2n+m}-\frac{1}{3n-m}+\text{etc.},
\end{equation*} 
 combining each two terms, we will have

\begin{equation*}
\frac{\pi}{n}\cot \frac{m}{n}\pi=\frac{1}{m}-\frac{2m}{nn-m^2}-\frac{2m}{4n^2-m^2}-\frac{2m}{9n^2-m^2}-\frac{2m}{16n^2-m^2}-\text{etc.},
\end{equation*}
whence we conclude

\begin{equation*}
\frac{1}{n^2-m^2}+\frac{1}{4n^2-m^2}+\frac{1}{9n^2-m^2}+\frac{1}{16n^2-m^2}+\text{etc.}=\frac{1}{2mm}-\frac{\pi}{2mn}\cot \frac{m}{n}\pi.
\end{equation*}
Now, let us set $n=1$ and for $m$ let us put $u$ that

\begin{equation*}
\frac{1}{1-u^2}+\frac{1}{4-u^2}+\frac{1}{9-u^2}+\frac{1}{16-u^2}+\text{etc.}=\frac{1}{2uu}-\frac{\pi }{2u}\cot \pi u.
\end{equation*} 
Resolve each of these fractions into series:

\begin{alignat*}{9}
&~\frac{1}{1-u^2}&&=~1~&&+u^2&&+u^4&&+u^6&&+u^8&&+\text{etc.}\\[3mm]
&~\frac{1}{4-u^2}&&=\frac{1}{2^2}&&+\frac{u^2}{2^4}&&+\frac{u^4}{2^6}&&+\frac{u^6}{2^8}&&+\frac{u^8}{2^8}&&+\text{etc.}\\[3mm]
&~\frac{1}{9-u^2}&&=\frac{1}{3^2}&&+\frac{u^2}{3^4}&&+\frac{u^4}{3^6}&&+\frac{u^6}{3^8}&&+\frac{u^8}{3^{10}}&&+\text{etc.}\\[3mm]
&\frac{1}{16-u^2}&&=\frac{1}{4^2}&&+\frac{u^2}{4^4}&&+\frac{u^6}{4^8}&&+\frac{u^8}{4^{10}}&&+\text{etc.} \\
& && && && \text{etc.}
\end{alignat*}

\paragraph*{§125}

Therefore, if one puts

\begin{alignat*}{19}
&1 &&+\frac{1}{2^2}&&+\frac{1}{3^2}&&+\frac{1}{4^2}&&+\text{etc.}&&=\mathfrak{a} \qquad \qquad &1 &&+\frac{1}{2^8}&&+\frac{1}{3^8}&&+\frac{1}{4^8}&&+\text{etc.}&&=\mathfrak{b} \\[3mm]
&1 &&+\frac{1}{2^4}&&+\frac{1}{3^4}&&+\frac{1}{4^4}&&+\text{etc.}&&=\mathfrak{c} \qquad \qquad &1 &&+\frac{1}{2^{10}}&&+\frac{1}{3^{10}}&&+\frac{1}{4^{10}}&&+\text{etc.}&&=\mathfrak{d} \\[3mm]
&1 &&+\frac{1}{2^6}&&+\frac{1}{3^6}&&+\frac{1}{4^6}&&+\text{etc.}&&=\mathfrak{c} \qquad \qquad &1 &&+\frac{1}{2^{12}}&&+\frac{1}{3^{12}}&&+\frac{1}{4^{12}}&&+\text{etc.}&&=\mathfrak{f} \\[3mm]
& && && && && && \text{etc.}
\end{alignat*}
the above series will be transformed into this one

\begin{equation*}
\mathfrak{a}+\mathfrak{b}u^2+\mathfrak{c}u^4+\mathfrak{d}u^6+\mathfrak{e}u^8+\mathfrak{f}u^{10}+\text{etc.}=\frac{1}{2uu}-\frac{\pi}{2u}\cot \pi u.
\end{equation*}
Therefore, because in § 118 the letters $A$, $B$, $C$, $D$ etc. were found to be of such a nature that, having put

\begin{equation*}
s=\frac{1}{u}-Au-Bu^3-Cu^5-Du^7-Eu^9-\text{etc.},
\end{equation*}
 $s=\frac{1}{2}\cot \frac{1}{2}u$, having written $\pi u$ instead of $\frac{1}{2}u$ or $2 \pi u$ instead of $u$, it will be

\begin{equation*}
\frac{1}{2}\cot \pi u=\frac{1}{2 \pi u}-2A \pi u-2^3 B \pi^3 u^3-2^5 C \pi^5 u^5 -2^7 D \pi^7 u^7 -\text{etc.},
\end{equation*}
whence, by multiplying by $\frac{\pi}{u}$, it will be

\begin{equation*}
\frac{\pi}{2 u}\cot \pi u=\frac{1}{2uu}-2 A\pi^2 -2^3 B \pi^4 u^2-2^5 C \pi^6 u^4-2^7 D \pi^8 u^6-\text{etc.},
\end{equation*}
and hence it follows that

\begin{equation*}
\frac{1}{2 uu}-\frac{\pi}{2 u}\cot \pi u= 2 A\pi^2 +2^3 B \pi^4 u^2+2^5 C \pi^6 u^4+2^7 D \pi^8 u^6+\text{etc.}
\end{equation*}
Since we just found that 

\begin{equation*}
\frac{1}{2uu}-\frac{\pi}{2u}\cot \pi u= \mathfrak{a}+\mathfrak{b}u^2+\mathfrak{c}u^4+\mathfrak{d}u^6+\text{etc.},
\end{equation*}
it is necessary that 

\begin{alignat*}{9}
& \mathfrak{a}&&= 2 &&A \pi^2 &&=\frac{~~~~~~~~2 \alpha ~~~~~~~~}{1 \cdot 2 \cdot 3}&&\pi^2 &&=\frac{~~~~~~~~2 \mathfrak{A}~~~~~~~}{1 \cdot 2}&&\pi^2 \\[3mm]
& \mathfrak{b}&&= 2^3 &&A \pi^4 &&=\frac{\quad \quad 2^3 \beta  \quad \quad}{1 \cdot 2 \cdot 3 \cdot 4 \cdot 5}&&\pi^4 &&=\frac{\quad \quad 2^3 \mathfrak{B} \quad \quad }{1 \cdot 2 \cdot 3 \cdot 4}&&\pi^4 \\[3mm]
& \mathfrak{c}&&= 2^5 &&A \pi^6 &&=\frac{\quad \quad 2^5 \gamma \quad \quad}{1 \cdot 2 \cdot 3  \cdots 7}&&\pi^6 &&=\frac{\quad \quad 2^3 \mathfrak{C} \quad \quad }{1 \cdot 2 \cdot 3 \cdots 6}&&\pi^6 \\[3mm]
& \mathfrak{d}&&= 2^7 &&A \pi^8 &&=\frac{\quad \quad 2^7 \delta \quad \quad }{1 \cdot 2 \cdot 3  \cdots 9}&&\pi^8 &&=\frac{\quad \quad 2^5 \mathfrak{D} \quad \quad}{1 \cdot 2 \cdot 3 \cdots 8}&&\pi^8 \\[3mm]
& \mathfrak{e}&&= 2^9 &&E \pi^{10} &&=\frac{2^{9} \varepsilon}{1 \cdot 2 \cdot 3  \cdots 11}&&\pi^{10} &&=\frac{2^9 \mathfrak{E}}{1 \cdot 2 \cdot 3 \cdots 10}&&\pi^{10} \\[3mm]
& \mathfrak{f}&&= 2^{11} &&F \pi^{12} &&=\frac{2^{11} \zeta}{1 \cdot 2 \cdot 3  \cdots 13}&&\pi^{12} &&=\frac{2^{11} \mathfrak{F}}{1 \cdot 2 \cdot 3 \cdots 12}&&\pi^{12} \\
& && && &&\text{etc.}
\end{alignat*}

\paragraph*{§126}

Therefore, by this simple reasoning not only all series of reciprocal powers we exhibited in the preceding paragraph are conveniently summed, but at the same time it is also understood how these sums are formed from the known values of the letters $\alpha$, $\beta$, $\gamma$, $\delta$, $\varepsilon$ etc. or even from the Bernoulli numbers $\mathfrak{A}$, $\mathfrak{B}$, $\mathfrak{C}$, $\mathfrak{D}$ etc. Therefore, since we defined fifteen of these numbers in § 122, from these one will be able to assign the sums of all even [reciprocal] powers up to the sum of this series:

\begin{equation*}
1+\frac{1}{2^{30}}+\frac{1}{3^{30}}+\frac{1}{4^{30}}+\frac{1}{5^{30}}+\text{etc.};
\end{equation*}
for, the sum of this series will be

\begin{equation*}
=\frac{2^{29}\pi}{1 \cdot 2 \cdot 3 \cdots 31}\pi^{31}=\frac{2^{29}\mathfrak{P}}{1 \cdot 2 \cdots 30}\pi^{30}.
\end{equation*}
And if one wants to determine more of these letters, this is very easily done by continuing these numbers $\alpha$, $\beta$, $\gamma$ etc. or these $\mathfrak{A}$, $\mathfrak{B}$, $\mathfrak{C}$ etc.

\paragraph*{§127}

Therefore, the origin of these numbers $\alpha$, $\beta$, $\gamma$, $\delta$ etc. or those formed from them $\mathfrak{A}$, $\mathfrak{B}$, $\mathfrak{C}$, $\mathfrak{D}$ etc. is basically the expansion of the cotangent of a certain angle into an infinite series. For, if 

\begin{equation*}
\frac{1}{2}\cot \frac{1}{2}u= \frac{1}{u}-Au-Bu^3-Cu^5-Du^7-Eu^9-\text{etc.},
\end{equation*}
it will be

\begin{equation*}
Au^2+Bu^4+Cu^6+Du^8+\text{etc.}=1-\frac{u}{2}\cot \frac{1}{2}u;
\end{equation*}
therefore, if  the respective values of the letters are substituted for coefficients $A$, $B$, $C$, $D$ etc., it will be found

\begin{equation*}
\frac{\alpha u^2}{1 \cdot 2 \cdot 3}+\frac{\beta u^4}{1 \cdot 2 \cdots 5}+\frac{\gamma u^7}{1 \cdot 2 \cdots 7}+\frac{\delta u^8}{1 \cdot 2 \cdots 9}+\text{etc.}=1-\frac{u}{2}\cot \frac{1}{2}u
\end{equation*}
and by using the Bernoulli numbers it will be

\begin{equation*}
\frac{\mathfrak{A}u^2}{1 \cdot 2}+\frac{\mathfrak{B}u^4}{1 \cdot 2 \cdot 3 \cdot 4}+\frac{\mathfrak{C}u^6}{1 \cdot 2 \cdots 6}+\frac{\mathfrak{D}u^8}{1 \cdot 2 \cdots 8}+\text{etc.}=1- \frac{u}{2}\cot \frac{1}{2}u,
\end{equation*}
from which series by differentiation innumerable others can be deduced and  infinite series these numbers enter can be summed.

\paragraph*{§128}

Let us take the first equation which we want to multiply by $u$ so that 

\begin{equation*}
\frac{\alpha u^3}{1 \cdot 2 \cdot 3}+\frac{\beta u^5}{1 \cdot 2 \cdots 5}+\frac{\gamma u^7}{1 \cdot 2 \cdots 7}+\frac{\delta u^9}{1 \cdot 2 \cdots 9}+\text{etc.}= u-\frac{uu}{2}\cot \frac{1}{2}u,
\end{equation*}
which differentiated and divided by $du$ gives

\begin{equation*}
\frac{\alpha u^2}{1 \cdot 2}+\frac{\beta u^4}{1 \cdot 2 \cdot 3 \cdot 4}+\frac{\gamma u^6}{1 \cdot 2 \cdots 6}+\frac{\delta u^8}{1 \cdot 2 \cdots 8}+\text{etc.}=1-u \cot \frac{1}{2}u+\frac{uu}{4(\sin \frac{1}{2}u)^2};
\end{equation*}
and if it is differentiated again, it will be

\begin{equation*}
\frac{\alpha u}{1}+\frac{\beta u^3}{1 \cdot 2 \cdot 3}+\frac{\gamma u^5}{1 \cdot 2 \cdot 3 \cdot 4 \cdot 5}+\text{etc.}=- \cot \frac{1}{2}u+\frac{u}{(\sin \frac{1}{2}u)^2}-\frac{uu \cos \frac{1}{2}u}{4(\sin \frac{1}{2}u)^2}.
\end{equation*}
But if the other equation is differentiated, it will be

\begin{equation*}
\frac{\mathfrak{A}u}{1}+\frac{\mathfrak{B}u^3}{1 \cdot 2 \cdot 3}+\frac{\mathfrak{C}u^5}{1 \cdot 2 \cdots 5}+\frac{\mathfrak{D}u^7}{1 \cdot 2 \cdots 7}+\text{etc.}=-\frac{1}{2}\cot \frac{1}{2}u+\frac{u}{4(\sin \frac{1}{2}u)^2}.
\end{equation*}
From these, if one puts $u=\pi$, because of $\cot \frac{1}{2}\pi =0$ and $\sin \frac{1}{2}\pi =1$, these summations follow

\begin{alignat*}{19}
& \quad  1 &&=\frac{\alpha \pi^2}{1 \cdot 2 \cdot 3}&&+\frac{\beta \pi^4}{1 \cdot 2 \cdot 3 \cdot 4 \cdot 5}&&+~\frac{\gamma \pi^6}{1 \cdot 2 \cdot 3 \cdots 7}&&+\frac{\delta \pi^8}{1 \cdot 2 \cdot 3 \cdots 9}&&+\text{etc.} \\[3mm]
&1+ \frac{\pi^2}{4}&&=~~\frac{\alpha \pi^2}{1 \cdot 2}&&+~~\frac{\beta \pi^4}{1 \cdot 2 \cdot 3 \cdot 4}&&+~\frac{\gamma \pi^6}{1 \cdot 2 \cdot 3 \cdots 6}&&+\frac{\delta \pi^8}{1 \cdot 2 \cdot 3 \cdots 8}&&+\text{etc.}\\[3mm]
&\quad \pi &&=~~~\frac{\alpha \pi}{1 }&&+~~~~\frac{\beta \pi^3}{1 \cdot 2 \cdot 3}&&+\frac{\gamma \pi^5}{1 \cdot 2 \cdot 3 \cdot 4 \cdot  5}&&+\frac{\delta \pi^7}{1 \cdot 2 \cdot 3 \cdots 7}&&+\text{etc.}
\intertext{or}
&\quad ~ 1 &&=~~~~~\alpha &&+~~~~\frac{\beta \pi^2}{1 \cdot 2 \cdot 3}&&+\frac{\gamma \pi^4}{1 \cdot 2 \cdot 3 \cdot 4 \cdot 5}&&+\frac{\delta \pi^6}{1 \cdot 2 \cdot 3 \cdots 7}&&+\text{etc.};\\
\end{alignat*}
if the first is subtracted from this one, it will remain 

\begin{equation*}
\alpha = \frac{(\alpha -\beta)\pi^2}{1 \cdot 2 \cdot 3}+\frac{(\beta -\gamma)\pi^4}{1 \cdot 2 \cdot 3 \cdot 4 \cdot 5}+\frac{(\alpha -\gamma)\pi^6}{1 \cdot 2 \cdot 3 \cdots 7}+\text{etc.}
\end{equation*}
But then it will be

\begin{alignat*}{9}
&1 &&= \frac{\mathfrak{A}\pi^2}{1 \cdot 2}&&+\frac{\mathfrak{B}\pi^4}{1 \cdot 2 \cdot 3 \cdot 4}&&+\frac{\mathfrak{C}\pi^6}{1 \cdot 2 \cdot 3 \cdots 6}&&+\frac{\mathfrak{D}\pi^8}{1 \cdot 2 \cdot 3 \cdots 8}\\[3mm]
&\frac{\pi}{4}&&= \frac{\mathfrak{A} \pi}{1}&&+~~\frac{\mathfrak{B}\pi^3}{1 \cdot 2 \cdot 3}~~&&+\frac{\mathfrak{C}\pi^5}{1 \cdot 2 \cdot 3 \cdot 4 \cdot 5}&&+\frac{\mathfrak{D}\pi^7}{1 \cdot 2 \cdot 3 \cdots 7}
\intertext{or}
& \frac{1}{4}&&=~~\frac{\mathfrak{A}}{1}~~&&+~~\frac{\mathfrak{B}\pi^2}{1 \cdot 2 \cdot 3}~~&&+\frac{\mathfrak{C}\pi^4}{1 \cdot 2 \cdot 3 \cdot 4 \cdot 5}&&+\frac{\mathfrak{D}\pi^6}{1 \cdot 2 \cdot 3 \cdots 7}&&+\text{etc.}
\end{alignat*}

\paragraph*{§129}

From the table of values of the numbers $\alpha$, $\beta$, $\gamma$, $\delta$ etc. we exhibited above (§ 121) it is plain that they decrease at first, but then increase and do so to infinity. Therefore, it will be worth one's while to investigate, how these numbers grow, after they had already been continued very far. Therefore, let $\varphi$ be any number of this series of the numbers $\alpha$, $\beta$, $\gamma$, $\delta$ etc. removed very far from the beginning and let $\psi$ be the following number. Since  the sum of the reciprocal powers are defined by means of these numbers, let $2n$ be the exponent of the power, whose sum the number $\varphi$ enters; $2n+2$ will be the exponent corresponding to number $\psi$ and the number $n$ will already be immensely large. Hence, from § 125 one will have

\begin{alignat*}{9}
&1&&+\frac{~~~~1~~~}{2^{2n}}&&+\frac{~~~~1~~~}{3^{2n}}&&+\frac{~~~~1~~~}{4^{2n}}&&+\text{etc.}&&=\frac{2^{2n-1}\varphi }{1 \cdot 2 \cdot 3 \cdots (2n+1)}\pi^{2n}, \\[3mm]
&1&&+\frac{1}{2^{2n+2}}&&+\frac{1}{3^{2n+2}}&&+\frac{1}{4^{2n+2}}&&+\text{etc.}&&=\frac{2^{2n+1}\varphi }{1 \cdot 2 \cdot 3 \cdots (2n+3)}\pi^{2n+2}.
\end{alignat*} 
Therefore, if these numbers are divided by each other, it will be

\begin{equation*}
\frac{1+\frac{1}{2^{2n+2}}+\frac{1}{3^{2n+2}}+\text{etc.}}{1+\frac{1}{2^{2n}}+\frac{1}{3^{2n}}+\text{etc.}}=\frac{4 \psi \pi^2}{(2n+2)(2n+3)\varphi}.
\end{equation*}
but since $n$ is an immensely large number and since both series are very close to $1$, it will be

\begin{equation*}
\frac{\psi}{\varphi}=\frac{(2n+2)(2n+3)}{4\pi^2}=\frac{nn}{\pi \pi}.
\end{equation*}
Therefore, because $n$ denotes, how far away the number $\varphi$   was  from the first number $\alpha$, that number $\varphi$ will have almost the same ratio to the following $\psi$ as $\pi^2$ has to to $n^2$ which ratio will be exact, if $n$ was an infinite number. Since it is almost $\pi \pi =10$, if one puts $n=100$, the hundredth term will be thousand times smaller than its following term. Therefore, the numbers $\alpha$, $\beta$, $\gamma$, $\delta$ etc. as the Bernoulli numbers $\mathfrak{A}$, $\mathfrak{B}$, $\mathfrak{C}$, $\mathfrak{D}$ etc. constitute a highly divergent series which grows even faster than a geometric series whose terms increase.

\paragraph*{§130}

Therefore, having found the values of the numbers $\alpha$, $\beta$, $\gamma$, $\delta$ etc. or $\mathfrak{A}$, $\mathfrak{B}$, $\mathfrak{C}$, $\mathfrak{D}$ etc., if a series is propounded whose general term $z$ was an arbitrary function of the index $x$, the summatory term $Sz$ of this series will be expressed as follows

\begin{equation*}
Sz = \int zdx +\frac{1}{2}z+\frac{1}{6}\cdot \frac{dz}{1 \cdot 2 }-\frac{1}{30}\cdot \frac{d^3 z}{1 \cdot 2 \cdot 3 \cdot 4 dx^3}
\end{equation*}
\begin{equation*}
+\frac{1}{42}\cdot \frac{d^5}{1 \cdot 2 \cdot 3 \cdots 6 dx^5}-\frac{1}{30}\cdot \frac{d^7z}{1 \cdot 2 \cdot 3 \cdots 8 dx^7}
\end{equation*}
\begin{equation*}
+\frac{5}{66}\cdot \frac{d^9 z}{1 \cdot 2 \cdot 3 \cdot 10 dx^9}-\frac{691}{2730}\cdot \frac{d^{11}z}{1 \cdot 2 \cdot 3 \cdots 12 dx^{11}}
\end{equation*}
\begin{equation*}
+\frac{7}{6}\cdot \frac{d^{13}z}{1 \cdot 2 \cdot 3 \cdots 14 dx^{13}}-\frac{3617}{510}\cdot \frac{d^{15}z}{1 \cdot 2 \cdot 3 \cdots 16 dx^{16}}
\end{equation*}
\begin{equation*}
\frac{854513}{138}\cdot \frac{d^{21} z}{1 \cdot 2 \cdot 3 \cdots 22 dx^{21}}-\frac{236364091}{2730}\cdot \frac{d^{23}z}{1 \cdot 2 \cdot 3 \cdots 24 dx^{23}}
\end{equation*}
\begin{equation*}
+\frac{8553103}{6}\cdot \frac{d^{25}z}{1 \cdot 2 \cdot 3 \cdots 26 dx^{25}}-\frac{23749461029}{870}\cdot \frac{d^{27}z}{1 \cdot 2 \cdot 3 \cdots 28 dx^{27}}
\end{equation*}
\begin{equation*}
+\frac{861584 127 6005}{14322}\cdot \frac{d^{29}z}{1 \cdot 2 \cdot 3 \cdots 30 dx^{29}}-\text{etc.}
\end{equation*}
Therefore, if the integral $\int zdx$ or the quantity whose differential is $=zdx$ is known, the summatory term will be found by means of iterated differentiation. But it is to be noted that always a constant of such a kind is to be added to this expression that the sum becomes $=0$, if the index $x$ is put $=0$.

\paragraph*{§131}

Therefore, if $z$ was a polynomial function of $x$, since its higher order differentials  vanish eventually, the summatory terms will be expressed by a finite expression; we will illustrate this in the following examples.

\subsection*{Example 1}

\textit{Let the summatory term of this series be in question}

\begin{alignat*}{9}
&1 \quad &&2 \quad && 3 \quad && 4 \quad && 5 \quad && \quad x \\
&1~+~ &&9~+~&&25~+~ &&49~+~&&81~+~ \cdots ~+~ && (2x-1)^2.
\end{alignat*}
Since here it is $z=(2x-1)^2=4xx-4x+1$, it will be

\begin{equation*}
\int z dx=\frac{4}{3}x^3-2x^2+x;
\end{equation*}
for, from the differentiation of this series the equation $4xxdx-4xdx+dx=zdx$ results. Furthermore, by differentiation it will be

\begin{equation*}
\frac{dz}{dx}=8x-4, \quad \frac{ddz}{dx^2}=8, \quad \frac{d^3z}{dx^3}=0 \quad \text{etc.}
\end{equation*}
Hence the summatory term in question will be

\begin{equation*}
\frac{4}{3}x^3-2x^2+x+2xx-2x+\frac{1}{2}+\frac{2}{3}x \pm \text{Const.},
\end{equation*}
 which constant has to cancel the terms $\frac{1}{2}-\frac{1}{3}$; therefore, it will be

\begin{equation*}
S(2x-1)^2=\frac{4}{3}x^3-\frac{1}{3}x=\frac{x}{3}(2x-1)(2x+1).
\end{equation*}
So, for $x=4$ the sum of the first four terms will be

\begin{equation*}
1+9+25+49=\frac{4}{3}\cdot 7 \cdot 9 =84.
\end{equation*}

\subsection*{Example 2}

\textit{Let the summatory term of this series be in question}

\begin{alignat*}{9}
&1 \quad && 2 \quad && 3 \quad && 4 \quad && \quad x \\
&1~+~ &&27~+~&&125~+~ &&343~+~ \cdots ~+~&&(2x-1)^3.
\end{alignat*}
Since  $z=(2x-1)^3=8x^3-12x^2+6x-1$, it will be

\begin{equation*}
\frac{dz}{dx}=24x^2-24x+6, \quad \frac{ddz}{dx^2}= 48x-24, \quad \frac{d^3z}{dx^3}=48;
\end{equation*}
the following differential quotients all vanish. Therefore, it will be

\begin{alignat*}{9}
&S(2x-1)^3 &&=2x^4 &&-4x^3 &&+3x^3 &&-1x \\
& && &&+4x^3 &&-6x^2 &&+3x &&-~\frac{1}{2}\\[3mm]
& && && &&+2x^2&&-2x&&+~\frac{1}{2} \\[3mm]
& && && && && &&-\frac{1}{15}\pm \text{Const.},  
\end{alignat*}
i.e.

\begin{equation*}
S(2x-1)^3=2x^4-x^2=x^2(2xx-1).
\end{equation*}
So, having put $x=4$, it will be

\begin{equation*}
1+27+125+343= 16\cdot 31=496.
\end{equation*}

\paragraph*{§132}

From this general expression found for the summatory term immediately that summatory term  we gave in the above part [§ 29 and 61] for the powers of the natural numbers and whose demonstration could not be given at that point follows. For, if we put $z=x^n$, it will be $\int zdx =\frac{1}{n+1}x^ {n+1}$; the differentials on the other hand will be 

\begin{equation*}
\frac{dz}{dx}=nx^{n-1}, \quad \frac{ddz}{dx^2}=n(n-1)x^{n-2}, \quad \frac{d^3z}{dx^3}=n(n-1)(n-2)x^{n-3},
\end{equation*}
\begin{equation*}
\frac{d^5z}{dx^5}=n(n-1)(n-2)(n-3)(n-4)x^{n-5}, \quad \frac{d^7z}{dx^7}=n(n-1)\cdots (n-6)x^{n-7} \quad \text{etc.}
\end{equation*}
From these therefore the following summatory term corresponding to the general term $x^n$ will be deduced \newpage

\begin{equation*}
Sx^n=\frac{1}{n+1}x^{n+1}+\frac{1}{2}x^n+\frac{1}{6}\cdot \frac{n}{2}x^{n-1}-\frac{1}{30}\cdot \frac{n(n-1)(n-2)}{2 \cdot 3 \cdot 4}x^{n-3}
\end{equation*}
\begin{alignat*}{9}
+\frac{1}{42}&\cdot \frac{n(n-1)(n-2)(n-3)(n-4)}{2 \cdot 3 \cdot 4 \cdot 5 \cdot 6}&&x^{n-5}\\[3mm]
-\frac{1}{30}&\cdot \frac{n(n-1)\cdots \cdots \cdots \cdots \cdot (n-6)}{2 \cdot 3 \cdots 8}&&x^{n-7}\\[3mm]
+\frac{5}{66}&\cdot \frac{n(n-1)\cdots \cdots \cdots \cdots \cdot (n-8)}{2 \cdot 3 \cdots 10}&&x^{n-9}\\[3mm]
-\frac{691}{2730}&\cdot \frac{n(n-1)\cdots \cdots \cdots \cdots (n-10)}{2 \cdot 3 \cdots 12}&&x^{n-11}\\[3mm]
+\frac{7}{6}&\cdot \frac{n(n-1)\cdots \cdots \cdots \cdots(n-12)}{2 \cdot 3 \cdots 14}&&x^{n-13}\\[3mm]
-\frac{3617}{510}&\cdot \frac{n(n-1)\cdots \cdots \cdots \cdots(n-14)}{2 \cdot 3 \cdots 16}&&x^{n-13}\\[3mm]
+\frac{43867}{798}&\cdot \frac{n(n-1)\cdots \cdots \cdots \cdots (n-16)}{2 \cdot 3 \cdots 18}&&x^{n-17}\\[3mm]
-\frac{174611}{330}&\cdot \frac{n(n-1)\cdots \cdots \cdots \cdots (n-18)}{2 \cdot 3 \cdots 20}&&x^{n-19}\\[3mm]
+\frac{854513}{138}&\cdot \frac{n(n-1)\cdots \cdots \cdots \cdots (n-20)}{2 \cdot 3 \cdots 22}&&x^{n-21}\\[3mm]
-\frac{236364091}{2730}&\cdot \frac{n(n-1)\cdots \cdots \cdots \cdots (n-22)}{2 \cdot 3 \cdots 24}&&x^{n-23}\\[3mm]
+\frac{8553103}{6}&\cdot \frac{n(n-1)\cdots \cdots \cdots \cdots (n-24)}{2 \cdot 3 \cdots 26}&&x^{n-25}\\[3mm]
-\frac{23749461029}{870}&\cdot \frac{n(n-1)\cdots \cdots \cdots \cdots (n-26)}{2 \cdot 3 \cdots 28}&&x^{n-27}\\[3mm]
+\frac{8615841276005}{14322}&\cdot \frac{n(n-1)\cdots \cdots \cdots \cdots (n-28)}{2 \cdot 3 \cdots 30}&&x^{n-29}\\
 &\text{etc.};
\end{alignat*}
this expression does not differ from the one we gave above except for the fact that  here we introduced the Bernoulli numbers $\mathfrak{A}$, $\mathfrak{B}$, $\mathfrak{C}$ etc., whereas above we used the numbers $\alpha$, $\beta$, $\gamma$, $\delta$ etc.; nevertheless, the agreement is immediately clear. Therefore, it is possible to exhibit the summatory term of all series up to the sums of the thirtieth powers; this investigation, if it would been done in another way, would have required very long and most tedious calculations.

\paragraph*{§133}

Above (§ 59) we already gave an almost identical expression to define the summatory term from the general term. But that expressions used the iterated differences of the general term; therefore, it differs from the expression we gave here mainly in that regard that it does not require the integral $\int zdx$, but each difference of the general term is multiplied by certain functions of $x$. Therefore, let us find the same expression again in the following way more accommodated to the nature of series, from which at the same time the rule  how the coefficients of the differentials proceed will be seen. Therefore, let the general term of the series be $z$, a function of the index $x$; the summatory term in question the other hand shall be $=s$; since this term, as we saw above, will be a function of $x$ vanishing for $x=0$,  applying the results we demonstrated above [§ 68] on the nature of functions of this kind it will be

\begin{equation*}
s- \frac{xds}{1 dx}+\frac{x^2 dds}{1 \cdot 2 dx^2}-\frac{x^3 d^3 s}{1 \cdot 2 \cdot 3 dx^3}+\frac{x^4 d^4s}{1 \cdot 2 \cdot 3 \cdot 4 dx^4}-\text{etc.}=0.
\end{equation*}

\paragraph*{§134}

Since $s$ denotes the sum of all terms of the series from the first to the last $z$, it is perspicuous, if in $s$  one writes $x-1$ instead of $x$, that the first sum does not contain the last term $z$; it will be 

\begin{equation*}
s-z=s-\frac{ds}{dx}+\frac{dds}{2dx^2}-\frac{d^3 s}{6dx^3}+\frac{d^4 s}{24 dx^4}-\text{etc.}
\end{equation*}
and hence

\begin{equation*}
z=\frac{ds}{dx}-\frac{dds}{2dx^2}+\frac{d^3 s}{6dx^3}-\frac{d^4s}{24dx^4}+\text{etc.},
\end{equation*}
which equation provides a way to define the general term from the given summatory term what is per se very easy. But from an  appropriate combination of this equation with that we found in the preceding paragraph, one will be able to define the value of $s$ in terms  of $x$ and $z$. For this aim, let us put 

\begin{equation*}
s-Az+\frac{Bdz}{dx}-\frac{Cddz}{dx^2}+\frac{Ddz^3}{dx^3}-\frac{Ed^4z}{dx^4}+\text{etc.}=0,
\end{equation*}
where $A$, $B$, $C$, $D$ etc. denote the necessary coefficients, either constant or variable; for, since

\begin{equation*}
z=\frac{ds}{dx}-\frac{dds}{2dx^2}+\frac{d^3 s}{6dx^3}-\frac{d^4 s}{24 dx^4}+\frac{d^5 s}{120 dx^5}-\text{etc.},
\end{equation*}
if the values for $z$, $\frac{dz}{dx}$, $\frac{ddz}{dx^2}$, $\frac{d^3 z}{dx^3}$ etc. are substituted in the above equation, it will result

\begin{equation*}
\renewcommand{\arraystretch}{2,5}
\setlength{\arraycolsep}{0.0mm}
\begin{array}{rrllcccccccccccccccccccccc}
~+~ &  s &~=~&s \\
~-~&Az &~=~&- \dfrac{Ads}{dx}&~+~&\frac{Adds}{2dx^2}&~-~&\dfrac{Ad^3s}{6dx^3}&~+~&\dfrac{Ad^4 s}{24dx^4}&~-~&\dfrac{Ad^5 s}{120 dx^5}&~+~&\text{etc.} \\
~+~&\dfrac{Bdz}{dx}&~=~& &~+~&\dfrac{Bdds}{dx^2}&~-~&\dfrac{Bd^3 s}{2 dx^3}&~+~&\dfrac{Bd^4 s}{6dx^4}&~-~&\dfrac{Bd^5 s}{24 dx^5}&~+~&\text{etc.} \\
~-~&\dfrac{Cddz}{dx^2}&~=~& && &~-~&\dfrac{Cd^3 s}{dx^3}&~+~&\dfrac{Cd^4 s}{2 dx^4}&~-~& \dfrac{Cd^5s}{6dx^5}&~+~&\text{etc.} \\
~+~& \dfrac{D d^3 z}{dx^3}&~=~& && && &~+~&\dfrac{Dd^4s}{dx^4}&~-~&\dfrac{Dd^5s}{2 dx^5}&~+~&\text{etc.} \\
~-~&\dfrac{Ed^4z}{dx^4}&~=~& && && && &~-~&\dfrac{Ed^5 s }{dx^5}&~+~&\text{etc.} \\
& && && && \text{etc.}
\end{array}
\end{equation*}
which series all added up therefore will be equal zero.

\paragraph*{§135}

Therefore, since we found before that

\begin{equation*}
0= s- \frac{xds}{dx}+\frac{x^2 dds}{2 dx^2}-\frac{x^3 d^3 s}{6dx^3}+\frac{x^4 d^4 s}{24 dx^4}-\frac{x^5 d^5 s}{120 dx^5}+\text{etc.},
\end{equation*}
if the above equation is put equal to this one, the following defining equations of the letters $A$, $B$, $C$, $D$ etc. will result

\begin{equation*}
A=x, \quad B=\frac{x^2}{2}-\frac{A}{2}, \quad C=\frac{x^3}{6}-\frac{B}{2}-\frac{A}{6},
\end{equation*}
\begin{equation*}
D=\frac{x^4}{24}-\frac{C}{2}-\frac{B}{6}-\frac{A}{24}, \quad E=\frac{x}{120}-\frac{D}{2}-\frac{C}{6}-\frac{B}{24}-\frac{A}{120} \quad \text{etc.}.
\end{equation*}
Therefore, having found the values of the letters $A$, $B$, $C$, $D$ etc.,  from the general term $z$ the summatory term $s=Sz$ will be determined in such a way that 

\begin{equation*}
Sz= Az-\frac{Bdz}{dx}+\frac{Cddz}{dx^2}-\frac{Dd^3z }{dx^3}+\frac{E d^4 z}{dx^4}-\frac{Fd^5 z}{dx^5}+\text{etc.}
\end{equation*}

\paragraph*{§136}

But since

\begin{equation*}
A=x, \quad B=\frac{1}{2}x^2-\frac{1}{2}x, \quad C=\frac{1}{6}x^3-\frac{1}{4}x^2+\frac{1}{12}x,
\end{equation*}
\begin{equation*}
D=\frac{1}{24}x^4-\frac{1}{12}x^3+\frac{1}{24}xx \quad \text{etc.},
\end{equation*}
it is clear that these coefficients are the same as those we had above (§ 59); hence that expression of the summatory term is the same as the one we found there and therefore it will be

\begin{equation*}
A=Sx^0=S1, \quad B=\frac{1}{1}Sx^1-\frac{1}{1}x, \quad C=\frac{1}{2}Sx^2-\frac{1}{2}x^2,
\end{equation*}
\begin{equation*}
D=\frac{1}{6}Sx^3-\frac{1}{6}x^3, \quad E=\frac{1}{24}Sx^4-\frac{1}{24}x^4 \quad \text{etc.}
\end{equation*}
Therefore, it will be

\begin{equation*}
Sz =xz-\frac{dz}{dx}Sx+\frac{ddz}{2dx^2}Sx^2-\frac{d^3z}{6dx^3}Sx^3+\frac{d^4z}{24dx^4}Sx^4-\text{etc.}
\end{equation*}
\begin{equation*}
+\frac{xdz}{dx}-\frac{x^2 ddz}{2 dx^2}+\frac{x^3 d^3 z}{6dx^3}-\frac{x^4d^4 z}{24 dx^4}+\text{etc.}
\end{equation*}
But if   $x=0$ in the general term $z$, the term corresponding to the index $=0$ will result; if it is put $=a$, it will be

\begin{equation*}
a= z-\frac{xdz}{dx}+\frac{x^2ddz}{2dx^2}-\frac{x^3d^3 z}{6 dx^3}+\text{etc.}
\end{equation*}
and hence

\begin{equation*}
\frac{xdz}{dx}-\frac{x^2 ddz}{2 dx^2}+\frac{x^3 d^3 z}{6dx^3}-\frac{x^4 d^4 z}{24 dx^4}+\text{etc.}=z-a, 
\end{equation*}
having substituted which value one will have

\begin{equation*}
Sz =(x+1)z-a-\frac{dz}{dx}Sx+\frac{ddz}{2 dx^2}Sx^2-\frac{d^3 z}{6dx^3}+\frac{d^4 z}{24 dx^4}Sx^4-\text{etc.}
\end{equation*}
Therefore, having found the sums of the powers  from a certain given general term, one can exhibit the summatory term corresponding to it.

\paragraph*{§137}

Therefore, since we found two expressions of the summatory term $Sz$ for the general term $z$ and  one of the formulas contains the integral $\int zdx$, if these two expressions are  equated, one will obtain the value of $\int zdx$ expressed by means of a series. For, since

\begin{equation*}
\int zdx = +\frac{1}{2}z+\frac{ \mathfrak{A}dz}{1 \cdot 2 dx}-\frac{\mathfrak{B}d^3z}{1 \cdot 2 \cdot 3 \cdot 4 dx^3}+\frac{\mathfrak{C}d^5 z}{1 \cdot 2 \cdots 6 dx^5}-\text{etc.}
\end{equation*} 
\begin{equation*}
=(x+1)z-a-\frac{dz}{1dx}Sx+\frac{ddz}{1 \cdot 2 dx^2}Sx^2-\frac{d^3 z}{1 \cdot 2 \cdot 3 }Sx^3+\text{etc.},
\end{equation*}
it will be

\begin{equation*}
\int zdx = \left(x+\frac{1}{2}\right)z-a-\frac{dz}{dx}\left(Sx+\frac{1}{2}\mathfrak{A}\right)+\frac{ddz}{2dx^2}Sx^2-\frac{d^3z}{6dx^3}\left(Sx^3-\frac{1}{4}\mathfrak{B}\right)
\end{equation*}
\begin{equation*}
+\frac{d^4 z}{24 dx^4}Sx^4-\frac{d^5 z}{120 dx^5}\left(Sx^5+\frac{1}{6}\mathfrak{C}\right)+\frac{d^6 z}{720 dx^6}Sx^6-\frac{d^7z}{5040dx^7}\left(Sx^7-\frac{1}{8}\mathfrak{D}\right)+\text{etc.},
\end{equation*}
where $\mathfrak{A}$, $\mathfrak{B}$, $\mathfrak{C}$, $\mathfrak{D}$ etc. denote the Bernoulli numbers exhibited above (§ 122).\\[2mm]
For the sake of an example, let $z=xx$; it will be $a=0$, $\frac{dz}{dx}=2x$ and $\frac{ddz}{2 dx^2}=1$; it will hence be

\begin{equation*}
\int xxdx = \left(x+\frac{1}{2}\right)xx-2x \left(\frac{1}{2}xx+\frac{1}{2}x+\frac{1}{12}\right)+1\left(\frac{1}{3}x^3+\frac{1}{2}x^2+\frac{1}{6}x\right)
\end{equation*}
or $\int xx dx = \frac{1}{3}x^3$; but $\frac{1}{3}x^3$, having differentiated it, gives $xxdx$, of course.

\paragraph*{§138}

Therefore, there is a new way to find the summatory terms of the series of powers; for, since from the coefficients  $A$, $B$, $C$, $D$ etc. assumed before these summatory terms are formed very easily, but each of these coefficients is conflated of the preceding ones, if   the values given in § 136  are substituted for these letters in the formulas given in § 135, it will be

\begin{alignat*}{9}
&Sx^1 &&-x &&=\frac{1}{2}xx&&-\frac{1}{2}x \\[3mm]
&Sx^2 &&-x^2&&=\frac{1}{3}x^3&&-\frac{1}{3}&&-\frac{2}{2}(Sx-x) \\[3mm]
&Sx^3 &&-x^3&&= \frac{1}{4}x^4&&-\frac{1}{4}x&&-\frac{3}{2}(Sx^2-x^2)&&-\frac{3 \cdot 2}{2 \cdot 3}(Sx-x)\\[3mm]
&Sx^4&&-x^4 &&=\frac{1}{5}x^5&&-\frac{1}{5}x&&-\frac{4}{2}(Sx^3-x^3)&&-\frac{4 \cdot 3}{2 \cdot 3}(Sx^2-x^2)&&-\frac{4 \cdot 3 \cdot 2}{2 \cdot 3 \cdot 4}(Sx-x)\\[3mm]
& && && && \text{etc.}
\end{alignat*}
Therefore, one will be able to form the sums of the above powers from the sums of the lower ones.

\paragraph*{§139}

But if we consider the law, which the coefficients $A$, $B$, $C$, $D$ etc. above (§ 135) were found to follow, with more attention, we will detect that they constitute a recurring series. For, if we expand this fraction

\begin{equation*}
y=
\renewcommand{\arraystretch}{1,4}
\setlength{\arraycolsep}{0mm}
\begin{array}{llllllllllllll}
  x&+&\frac{1}{2}xxu&+&\frac{1}{6}x^3 u^2&+& \frac{1}{24}x^4 u^3&+&\frac{1}{120}x^5u^4&+&\text{etc.}\\
  \hline
  1&+&\frac{1}{2}u&+&\frac{1}{6} u^2&+& \frac{1}{24} u^3&+&\frac{1}{120}u^4&+&\text{etc.}
\end{array}
\end{equation*} 
into a power series in $u$ and assume this series to result

\begin{equation*}
A+Bu+Cu^2+Du^3+Eu^4+\text{etc.},
\end{equation*}
it will be, as we found before,

\begin{equation*}
A=x, \quad B=\frac{1}{2}xx-\frac{1}{2}A \quad \text{etc.}
\end{equation*}
and so having found this series one will obtain the summatory terms of the series of powers. But that fraction, from whose expansion this series results, will go over into this form $\frac{e^{xu}-1}{e^u-1}$ which, if $x$ was a positive integer, goes over into

\begin{equation*}
1+e^{u}+e^{2u}+e^{3u}+\cdots+e^{(x-1)u};
\end{equation*}
therefore, since

\begin{equation*}
\renewcommand{\arraystretch}{2,5}
\setlength{\arraycolsep}{0.0mm}
\begin{array}{rcccccccccccccccccc}
1 &~=~&1 \\
e^u &~=~&1 &~+~&\dfrac{ \quad  u \quad }{1}&~+~&\dfrac{u^2}{1 \cdot 2}&~+~&\dfrac{u^3}{1 \cdot 2 \cdot 3}&~+~&\dfrac{u^4}{1 \cdot 2 \cdot 3 \cdot 4}&~+~&\text{etc.} \\
e^{2u} &~= ~&1 &~+~&\dfrac{ 2u }{1}&~+~&\dfrac{4u^2}{1 \cdot 2}&~+~&\dfrac{8 u^3 }{1 \cdot 2 \cdot 3}&~+~&\dfrac{16 u^4}{1 \cdot 2 \cdot 3 \cdot 4}&~+~&\text{etc.} \\
e^{3u} &~=~&1 &~+~&\dfrac{3u}{1}&~+~&\dfrac{9u^2}{1 \cdot 2}&~+~&\dfrac{27 u^3}{1 \cdot 2 \cdot 3}&~+~&\dfrac{81 u^4}{1 \cdot 2 \cdot 3 \cdot 4}&~+~&\text{etc.} \\
e^{(x-1)u} &~=~&1 &~+~&\dfrac{(x-1)u}{1}&~+~&\dfrac{(x-1)^2u^2}{1 \cdot 2}&~+~&\dfrac{(x-1)^3 u^3}{1 \cdot 2 \cdot 3}&~+~&\dfrac{(x-1)^4 u^4}{1 \cdot 2 \cdot 3 \cdot 4}&~+~&\text{etc.} 
\end{array}
\end{equation*}
it will be

\begin{alignat*}{9}
&A&&=x \\[3mm]
&B &&=S(x-1)&&=Sx &&-x \\[3mm]
&C &&=\frac{1}{2}S(x-1)^2&&=\frac{1}{2}Sx^2&&-\frac{1}{2}x^2\\[3mm]
&D &&=\frac{1}{6}S(x-1)^3&&=\frac{1}{6}Sx^3&&-\frac{1}{6}x^3\\[3mm]
& && &&\text{etc.}
\end{alignat*}
Therefore, the connection already mentioned before of these coefficients to the sums of the powers is confirmed and demonstrated completely.


\title{On the Summation of Progressions by means of infinite Series
\footnote{
Original title: "'De Summatione Progressionum per Series infinitas"', first published as part of the book \textit{Institutiones calculi differentialis cum eius usu in analysi finitorum ac doctrina serierum}, 1755, reprinted in in \textit{Opera Omnia}: Series 1, Volume 10, pp. 337 - 367 "', Eneström-Number E212, translated by: Alexander Aycock for the "`Euler-Kreis Mainz"'}}
\author{Leonhard Euler}
\date{ }
\maketitle

\paragraph*{§140}

The general expression, we found  in the preceding chapter for the summatory term of a  series, whose general term or the term corresponding to the index $x$ is $=z$, i.e.

\begin{equation*}
Sz = \int zdx +\frac{1}{2}z+\frac{\mathfrak{A}dz}{1 \cdot 2 dx}-\frac{\mathfrak{B}d^3z}{1 \cdot 2 \cdot 3 \cdot 4 dx^3}+\frac{\mathfrak{C}d^5 z}{1 \cdot 2 \cdots 6 dx^5}-\text{etc.},
\end{equation*}
especially yields the summation of series whose general terms are  polynomial functions of the index $x$, since in these cases one  eventually gets to vanishing differentials. But if $z$ was not a function of $x$ of such a kind, then its differentials proceed to infinity and so an infinite series expressing the sum of the propounded series up to the given term, whose index is $=x$, results. Therefore, the sum of the propounded series continued to infinity will result for $x=\infty$; and this way another infinite series equal to the first arises.

\paragraph*{§141}

But if one puts $x=0$, then the expression exhibiting the sum has to vanish, as we already mentioned; if this does not happen, a constant quantity has either to be added or subtracted that this condition is satisfied. If, having done this, one puts $x=1$, the found sum will yield the first term of the series; but if one puts $x=2$, the aggregate of the first and the second, if $x=3$, the aggregate of the  three initial terms of the series will result, and so forth. Therefore, in these cases, since the sum of one or two or three etc. terms is known, the value of the infinite series expressing this sum will become known, and from this source one will be able to sum innumerable series.

\paragraph*{§142}

Since, if a constant of such a kind was added to the sum that it vanishes for $x=0$,  the sum is correct in the remaining cases, i.e. for all  $x$, it is obvious, if  a constant quantity of such a kind is added to the found sum that in one special case the true sum results, that then  the true sum has to result also in all remaining cases. Hence, if  having put put $x=0$ it is not clear, a value of which kind the expression of the sum will have, and hence the constant to be added cannot be found, then any other number can be assumed for $x$ and by adding a constant  the correct sum can be forced to result; how this has to happen, will become more perspicuous in the following.

\paragraph*{§142[a]}

At first, let us consider this harmonic progression

\begin{equation*}
1+\frac{1}{2}+\frac{1}{3}+\frac{1}{4}+\cdots+\frac{1}{x}=s;
\end{equation*}
because its general term  is $=\frac{1}{x}$, it will be $z=\frac{1}{x}$ and the summatory term $s$ will be found this way. First, it will be $\int zdx= \int \frac{dx}{x}=\log x$; further, the differentials will be as follows
\begin{equation*}
\frac{dz}{dx}=-\frac{1}{x^2}, \quad \frac{ddz}{2dx^2}=\frac{1}{x^3}, \quad \frac{d^3 z}{6 dx^3}=-\frac{1}{x^4}, \quad \frac{d^4 z}{24 dx^4}=\frac{1}{x^5}, \quad \frac{d^5 z}{120 dx^5}=-\frac{1}{x^6} \quad \text{etc.}
\end{equation*}
Therefore, it will  be

\begin{equation*}
s = \log x+ \frac{1}{2x}-\frac{\mathfrak{A}}{2x^2}+\frac{\mathfrak{B}}{4x^4}-\frac{\mathfrak{C}}{6x^6}+\frac{\mathfrak{D}}{8 x^8}-\text{etc.}+\text{Constant}.
\end{equation*}
Therefore, the constant to be added here cannot be defined from the case $x=0$. So put $x=1$, since then $s=1$; it will be

\begin{equation*}
1=\frac{1}{2}-\frac{\mathfrak{A}}{2}+\frac{\mathfrak{B}}{4}-\frac{\mathfrak{C}}{6}+\frac{\mathfrak{D}}{8}-\text{etc.}+\text{Const.},
\end{equation*}
whence this constant becomes

\begin{equation*}
=\frac{1}{2}+\frac{\mathfrak{A}}{2}-\frac{\mathfrak{B}}{4}+\frac{\mathfrak{C}}{6}-\frac{\mathfrak{D}}{8}+\text{etc.},
\end{equation*}
and hence the summatory term in question will be

\begin{alignat*}{9}
&s &&= \log x &&+\frac{1}{2x}&&-\frac{\mathfrak{A}}{2x^2}&&+\frac{\mathfrak{B}}{4x^4}&&-\frac{\mathfrak{C}}{6x^6}&&+\frac{\mathfrak{D}}{8x^8}&&-\text{etc.} \\[3mm]
& && &&+~\frac{1}{2}&&+~~\frac{\mathfrak{A}}{2}&&-~~\frac{\mathfrak{B}}{4}&&+~~\frac{\mathfrak{C}}{6}&&-~~\frac{\mathfrak{C}}{8}&&+\text{etc.}
\end{alignat*}

\paragraph*{§143}

Since the Bernoulli numbers $\mathfrak{A}$, $\mathfrak{B}$, $\mathfrak{C}$, $\mathfrak{D}$ etc. constitute a divergent series, this value of the constant cannot be found. But if  a  larger number is substituted for $x$ and the sum of that many terms is actually found, the value of the constant can be  investigated conveniently. For this purpose, set $x=10$ and by adding the first ten terms, one will find their sum to be

\begin{equation*}
=2.928968253968253968,
\end{equation*}
to which the expression for the sum has to be equal, if  one puts $x=10$; that expression is

\begin{equation*}
\log 10+\frac{1}{20}-\frac{\mathfrak{A}}{200}+\frac{\mathfrak{B}}{40000}-\frac{\mathfrak{C}}{6000000}+\frac{\mathfrak{D}}{800000000}-\text{etc.}+C.
\end{equation*}
Therefore, having taken  the hyperbolic logarithm of ten for $\log 10$ and  having substituted the values found above [§ 122] for $\mathfrak{A}$, $\mathfrak{B}$, $\mathfrak{C}$ etc., one will find that constant  to be

\begin{equation*}
C=0.5772156649015325,
\end{equation*}
which number therefore is equal to the sum of the series

\begin{equation*}
\frac{1}{2}+\frac{\mathfrak{A}}{2}-\frac{\mathfrak{B}}{4}+\frac{\mathfrak{C}}{6}-\frac{\mathfrak{D}}{8}+\frac{\mathfrak{E}}{10}-\text{etc.}
\end{equation*}

\paragraph*{§144}

If  sufficiently small numbers are substituted for $x$, since the sum of the series is actually easily found, one will obtain the sum of this series

\begin{equation*}
\frac{1}{2x}-\frac{\mathfrak{A}}{2x^2}+\frac{\mathfrak{B}}{4x^4}-\frac{\mathfrak{C}}{6x^6}+\frac{\mathfrak{D}}{8x^8}-\text{etc.}=s- \log x-C.
\end{equation*}
But if $x$ denotes a very large number, since then the value of this infinite expression  is easily assigned in decimal numbers,  vice versa the sum of the series will be defined. And first it is certainly clear, if the series is continued to infinity, that its sum will be infinitely large; for, having put $x=\infty$, also $\log x$ becomes infinite, even though $\log x$ has an infinitely small ratio to $x$. But in order to assign the sum of an arbitrary number of terms of the series in a convenient manner, let us express the values of the letters $\mathfrak{A}$, $\mathfrak{B}$, $\mathfrak{C}$ etc. in decimal fractions.

\begin{alignat*}{9}
&\mathfrak{A}&&=0.1666666666666\\
&\mathfrak{B}&&=0.0333333333333\\
&\mathfrak{C}&&=0.0238095238095\\
&\mathfrak{D}&&=0.0333333333333\\
&\mathfrak{E}&&=0.0757575757575\\
&\mathfrak{F}&&=0.2531135531135\\
&\mathfrak{G}&&=1.1666666666666\\
&\mathfrak{H}&&=7.0921568627451 \quad \text{etc.}
\intertext{whence}
&\frac{\mathfrak{A}}{2} &&=0.0833333333333\\
&\frac{\mathfrak{B}}{4} &&=0.0083333333333\\
&\frac{\mathfrak{C}}{6} &&=0.0039682539682\\
&\frac{\mathfrak{D}}{8} &&=0.0041666666666\\
&\frac{\mathfrak{E}}{10}&&=0.0075757575757\\
&\frac{\mathfrak{F}}{12}&&=0.0210927960928\\
&\frac{\mathfrak{G}}{14}&&=0.0833333333333\\
&\frac{\mathfrak{H}}{16}&&=0.4432598039316\\
\end{alignat*}

\subsection*{Example 1}

\textit{To find the sum of thousand terms of the series $1+\frac{1}{2}+\frac{1}{3}+\frac{1}{4}+\frac{1}{5}+\frac{1}{6}+\text{etc.}$}\\[2mm]
Therefore, put $x=1000$, and since

\begin{alignat*}{9}
& \log 10 &&=2.3025850929940456840,
\intertext{it will be}
& \log x                                &&=6.9077553789821 \\
& \text{Const.}                        &&=0.5772156649015 \\
&\frac{1}{2x}                          &&=0.0005000000000 
\intertext{and in total}
&                                      &&=7.4854709438836 \\
&\text{subtr. }\frac{\mathfrak{A}}{2xx}&&=0.0000000833333 
\intertext{which yields}
&                                      &&=7.4854708605503 \\
&\text{add } \frac{\mathfrak{B}}{4x^4}   &&=0.0000000000000 
\intertext{Therefore}
&                                      &&=7.4854708605503
\end{alignat*}
is the  sum of thousand terms in question, which is still smaller than seven and a half units.

\subsection*{Example 2}

\textit{To find the sum of a million terms of the series $1+\frac{1}{2}+\frac{1}{3}+\frac{1}{4}+\frac{1}{5}+\frac{1}{6}+\text{etc.}$}\\[2mm]
Since $x=1000000$, it will be $\log x= 6 \cdot \log 10$, therefore

\begin{alignat*}{9}
&\log x &&= 1&&2.8155105579642 \\
&\text{Const.} &&= &&0.05772156649015 \\
&\frac{1}{2x}&&= &&0.0000005000000
\intertext{in total}
& &&= &&14.3927267228657=\text{the sum in question}
\end{alignat*}

\paragraph*{§145}

Therefore, if  one substitutes a very large number for $x$, the sum is found to a high degree of accuracy using only the first term $\log x$ increased by the constant $C$; therefore, extraordinary corollaries can  be deduced from this. So, if $x$ was a very large number and one puts

\begin{equation*}
1+\frac{1}{2}+\frac{1}{3}+\frac{1}{4}+\frac{1}{5}+\cdots +\frac{1}{x}=s
\end{equation*}
and

\begin{equation*}
1+\frac{1}{2}+\frac{1}{3}+\frac{1}{4}+\cdots + \frac{1}{x}+\cdots +\frac{1}{x+y}=t,
\end{equation*}
since  approximately $s= \log x +C$ and $t=\log (x+y)+C$, it will be

\begin{equation*}
t-s=\log(x+y)- \log x=\frac{x+y}{x}
\end{equation*}
and hence this logarithm is approximately expressed by means of a harmonic series consisting of a finite numbers of terms as follows

\begin{equation*}
\log \frac{x+y}{y}=\frac{1}{x+1}+\frac{1}{x+2}+\frac{1}{x+3}+\cdots +\frac{1}{x+y}.
\end{equation*}
But this logarithm is exhibited more accurately, if more terms of the above sums $s$ and $t$ are taken . So, since

\begin{equation*}
s=\log x+ C+\frac{1}{2x}-\frac{1}{12xx} \quad \text{and} \quad t=\log(x+y)+C+\frac{1}{2(x+y)}-\frac{1}{12(x+y)^2}
\end{equation*}
it will be

\begin{equation*}
t-s= \log \frac{x+y}{x}-\frac{1}{2x}+\frac{1}{2(x+y)}+\frac{1}{12xx}-\frac{1}{12(x+y)^2},
\end{equation*}
and hence

\begin{equation*}
\log \frac{x+y}{x}=\frac{1}{x+1}+\frac{1}{x+2}+\frac{1}{x+3}+\cdots+ \frac{1}{x+y}+\frac{1}{2x}-\frac{1}{2(x+y)}-\frac{1}{12xx}+\frac{1}{12(x+y)^2}.
\end{equation*}
But if $x$ is such a large number that the two last terms can be neglected, it will approximately be

\begin{equation*}
\log \frac{x+y}{x}=\frac{1}{x+1}+\frac{1}{x+2}+\frac{1}{x+3}+\cdots+ \frac{1}{x+y}+\frac{1}{2}\left(\frac{1}{x}-\frac{1}{x+y}\right).
\end{equation*}

\paragraph*{§ 145[a]}

Using this harmonic series, we will also be able to define the sum of this series, in which only the odd numbers occur, i.e.

\begin{equation*}
\frac{1}{1}+\frac{1}{3}+\frac{1}{5}+\frac{1}{7}+\frac{1}{9}+\cdots +\frac{1}{2x+1}.
\end{equation*}
For, since, taking all terms,

\begin{equation*}
1+\frac{1}{2}+\frac{1}{3}+\frac{1}{4}+\cdots + \frac{1}{2x}+\frac{1}{2x+1}
\end{equation*}
\begin{equation*}
= \log (2x+1)+C+\frac{1}{2(2x+1)}-\frac{\mathfrak{A}}{2(2x+1)^2}+\frac{\mathfrak{B}}{4(2x+1)^4}-\frac{\mathfrak{C}}{6(2x+1)^6}+\text{etc.},
\end{equation*}
the sum of the even terms

\begin{equation*}
\frac{1}{2}+\frac{1}{4}+\frac{1}{6}+\cdots+\frac{1}{2x}
\end{equation*}
is the half of the above series, i.e.

\begin{equation*}
\frac{1}{2}C+\frac{1}{2}\log x+\frac{1}{4x}-\frac{\mathfrak{A}}{4x^2}+\frac{\mathfrak{B}}{8x^4}-\frac{\mathfrak{C}}{12x^6}+\frac{\mathfrak{D}}{16x^8}-\text{etc.},
\end{equation*}
having subtracted this series from the latter,

\begin{equation*}
1+\frac{1}{3}+\frac{1}{5}+\frac{1}{7}+\cdots+ \frac{1}{2x+1}
\end{equation*}
\begin{equation*}
=\frac{1}{2}C+\log \frac{2x+1}{\sqrt{x}}+\frac{1}{2(2x+1)}-\frac{\mathfrak{A}}{2(2x+1)^2}+\frac{\mathfrak{B}}{4(2x+1)^4}-\text{etc.}
\end{equation*}
\begin{equation*}
-\frac{1}{4x}+\frac{\mathfrak{A}}{4x^2}-\frac{\mathfrak{B}}{8x^4}+\text{etc.}
\end{equation*}

\paragraph*{§146}

One can indeed even find the sum of any harmonic series by means of the same general expression; for, let 

\begin{equation*}
\frac{1}{m+n}+\frac{1}{2m+n}+\frac{1}{3m+n}+\frac{1}{4m+n}+\cdots+ \frac{1}{mx+n}=s;
\end{equation*}
since the general term is $z=\frac{1}{mx+n}$, it will be

\begin{equation*}
\int zdx=\frac{1}{m}\log (mx+n), \quad \frac{dz}{dx}=-\frac{m}{(mx+n)^2}, \quad \frac{ddz}{2dx^2}=\frac{mm}{(mx+n)^3},
\end{equation*}
\begin{equation*}
\frac{d^3 z}{6dx^3}=-\frac{m^3}{(mx+n)^4}, \quad \frac{d^4 z}{24dx^4}=\frac{m^4}{(mx+n)^5}, \quad \frac{d^5 z}{120 dx^5}=- \frac{m^5}{(mx+n)^6}\quad \text{etc.}
\end{equation*}
From these one hence finds

\begin{equation*}
s=D+\frac{1}{m}\log (mx+n)+\frac{1}{2(mx+n)}-\frac{\mathfrak{A}m}{2(mx+n)^2}+\frac{\mathfrak{B}m^3}{4(mx+n)^4}
\end{equation*}
\begin{equation*}
-\frac{\mathfrak{C}m^5}{6(mx+n)^6}+\frac{\mathfrak{D}m^7}{8(mx+n)^8}-\text{etc.}
\end{equation*}
Therefore, having put $x=0$, the constant to be added will be

\begin{equation*}
D=-\frac{1}{m}\log n -\frac{1}{2n}+\frac{\mathfrak{A}m}{2n^2}-\frac{\mathfrak{B}m^3}{4n^4}+\frac{\mathfrak{C}m^5}{6n^6}-\text{etc.}
\end{equation*}

\paragraph*{§147}

But if $n=0$, since the sum of the series

\begin{equation*}
\frac{1}{m}+\frac{1}{2m}+\frac{1}{3m}+\frac{1}{4m}+\cdots +\frac{1}{mx}
\end{equation*}
is

\begin{equation*}
=\frac{1}{m}C+\frac{1}{m}\log x +\frac{1}{2mx}-\frac{\mathfrak{A}}{2mx^2}+\frac{\mathfrak{B}}{4mx^4}-\text{etc.},
\end{equation*}
but the sum of this series 

\begin{equation*}
1+\frac{1}{2}+\frac{1}{3}+\frac{1}{4}+\frac{1}{5}+\cdots +\frac{1}{mx}
\end{equation*}
is

\begin{equation*}
=C+ \log mx +\frac{1}{2mx}-\frac{\mathfrak{A}}{2m^2x^2}+\frac{\mathfrak{B}}{4m^4 x^4}-\text{etc.},
\end{equation*}
if from this series the other series is subtracted $m$ times  that this series results

\begin{alignat*}{9}
&1+\frac{1}{2}+\cdots &&+\frac{1}{m}+\cdots + &&+\frac{1}{2m}+ \cdots &&+\frac{1}{3m}&&+ \cdots &&+\frac{1}{mx}\\[3mm]
&  &&-\frac{m}{m}&&-\frac{m}{2m}&&-\frac{m}{3m}&& &&-\frac{m}{mx},
\end{alignat*}
its sum will be

\begin{alignat*}{9}
&= \log m &&+\frac{1}{2mx}&&-\frac{\mathfrak{A}}{2m^2x^2}&&+\frac{\mathfrak{B}}{4m^4x^4}&&-\text{etc.} \\[3mm]
& &&-~~\frac{1}{2x}&&+~~\frac{\mathfrak{A}}{2xx}&&-~~\frac{\mathfrak{B}}{4x^4}&&+\text{etc.},
\end{alignat*}
and if one sets $x=\infty$, the sum will be $=\log m$. Hence, by taking the numbers $2$, $3$, $4$ etc. for $m$ it will be

\begin{alignat*}{19}
&\log 2 &&=1 &&-\frac{1}{2}&&+\frac{1}{3}&&-\frac{1}{4}&&+\frac{1}{5}&&-\frac{1}{6}&&+\frac{1}{7}&&-\frac{1}{8}&&+\text{etc.} \\[3mm]
& \log 3 &&=1&&+\frac{1}{2}&&-\frac{2}{3}&&+\frac{1}{4}&&+\frac{1}{5}&&-\frac{2}{6}&&+\frac{1}{7}&&+\frac{1}{8}&&-\frac{2}{9}+\text{etc.}\\[3mm]
&\log 4 &&=1 &&+\frac{1}{2}&&+\frac{1}{3}&&-\frac{3}{4}&&+\frac{1}{5}&&+\frac{1}{6}&&+\frac{1}{7}&&-\frac{3}{8}&&+\text{etc.}\\[3mm]
& \log 5 &&=1 &&+\frac{1}{2}&&+\frac{1}{3}&&+\frac{1}{4}&&-\frac{4}{5}&&+\frac{1}{6}&&+\frac{1}{7}&&+\frac{1}{8}&&+\frac{1}{9}-\frac{4}{10}+\text{etc.}\\
& && && && && && &&\text{etc.}
\end{alignat*}

\paragraph*{§148}

Having discussed the harmonic series now,  let us proceed to the reciprocal series of the squares and let

\begin{equation*}
s=1+\frac{1}{4}+\frac{1}{9}+\frac{1}{16}+\cdots + \frac{1}{xx};
\end{equation*}
since its general term is $z=\frac{1}{xx}$, it will be $\int zdx=-\frac{1}{x}$ and the differentials of $z$ will be

\begin{equation*}
\frac{dz}{2dx}=-\frac{1}{x^3}, \quad \frac{ddz}{2 \cdot 3 dx^2}=\frac{1}{x^4}, \quad \frac{d^3 z}{2 \cdot 3 \cdot 4 dx^3}=-\frac{1}{x^5} \quad \text{etc.},
\end{equation*}
whence the sum will be

\begin{equation*}
s=C-\frac{1}{x}+\frac{1}{2xx}-\frac{\mathfrak{A}}{x^3}+\frac{\mathfrak{B}}{x^5}-\frac{\mathfrak{C}}{x^7}+\frac{\mathfrak{D}}{x^9}-\frac{\mathfrak{E}}{x^{11}}+\text{etc.},
\end{equation*}
in which the constant $C$ to be added is to be defined from a special case, in which the sum is known. Therefore, let us put $x=1$; since $x=1$, it has to be

\begin{equation*}
C=1+1-\frac{1}{2}+\mathfrak{A}-\mathfrak{B}+\mathfrak{C}-\mathfrak{D}+\mathfrak{E}-\text{etc.},
\end{equation*}
which series, since highly divergent, does not show the value of the constant $C$. But since we demonstrated above [§ 125] that the sum of this series continued to infinity is $=\frac{\pi \pi}{6}$, for $x=\infty$, if one puts $s=\frac{\pi \pi}{6}$, it will be $C=\frac{\pi \pi}{6}$ because all remaining terms vanish. Therefore, it will be

\begin{equation*}
1+1-\frac{1}{2}+\mathfrak{A}-\mathfrak{B}+\mathfrak{C}-\mathfrak{D}+\mathfrak{E}-\text{etc.}=\frac{\pi \pi}{6}.
\end{equation*}

\paragraph*{§149}

But if the sum of this series would have not been known, the value of that constant $C$ would have to be determined from another case, in which the sum was actually found. For this aim, let us put $x=10$ and actually adding ten terms one will find

\begin{equation*}
\renewcommand{\arraystretch}{2,0}
\setlength{\arraycolsep}{0mm}
\begin{array}{llclll}
                          & &s&~=~&1.549767731166540690\\
  \text{then} \quad & &\\
  &\text{add } & \dfrac{1}{x}&~=~&0.1 \\[2mm]
  &\text{subtr. } & \dfrac{1}{2xx}&=&0.0005 \\[2mm]
  \cline{3-5}
  &            &              & &1.644767731166540690 \\[2mm]
 
  &\text{add } & \dfrac{\mathfrak{A}}{x^3}&=&0.000166666666666666 \\[2mm]
  &            &              & &1.644934397833207356 \\[2mm]
  \cline{3-5}
  &\text{subtr. } & \dfrac{\mathfrak{B}}{x^5}&=&0.000000333333333333 \\[2mm]
  &            &              & &1.644934064499874023 \\[2mm]
  \cline{3-5}
   &\text{add } & \dfrac{\mathfrak{C}}{x^7}&=&0.000000002380952381 \\[2mm]
  &            &              & &1.644934066880826404 \\[2mm]
  \cline{3-5}
   &\text{subtr. } & \dfrac{\mathfrak{D}}{x^9}&=&0.000000000033333333 \\[2mm]
  &            &              & &1.6444066847493071 \\[2mm]
  \cline{3-5}
\end{array}
\end{equation*}
\begin{equation*}
\renewcommand{\arraystretch}{2,0}
\setlength{\arraycolsep}{0mm}
\begin{array}{llclll}
   &\text{add } & \dfrac{\mathfrak{E}}{x^{11}}&~=~&0.000000000000757575 \\[2mm]
  &            &              & &1.644934066848250646 \\[2mm]
  \cline{3-5}
   &\text{subtr. } & \dfrac{\mathfrak{F}}{x^{13}}&~=~&0.000000000000025311 \\[2mm]
  &            &              & &1.644934066848225335 \\
  \cline{3-5}
    &\text{add } & \dfrac{\mathfrak{G}}{x^{15}}&~=~&0.000000000000001166 \\[2mm]
    &\text{subtr. } & \dfrac{\mathfrak{H}}{x^{17}}&~=~&0.000000000000000071 \\[2mm]
  &            &              & &1.644934066848226430 &~=~ C. 
\end{array}
\end{equation*}
And this value at the same time is the value of the expression $\frac{\pi \pi}{6}$, as  anyone carrying out the calculation using the known value of $\pi$ will discover. Therefore, it is understood at the same time, even though the series $\mathfrak{A}$, $\mathfrak{B}$, $\mathfrak{C}$ etc. diverges, that nevertheless the true sum results this way.

\paragraph*{§150}

Now let $z=\frac{1}{x^3}$ and

\begin{equation*}
s=1+\frac{1}{2^3}+\frac{1}{3^3}+\frac{1}{4^3}+\cdots +\frac{1}{x^3};
\end{equation*} 
since 

\begin{equation*}
\int zdx= -\frac{1}{2xx}, \quad \frac{dz}{1 \cdot 2 \cdot 3 dx}=-\frac{1}{2x^4}, \quad \frac{ddz}{1 \cdot 2 \cdot 3 \cdot 4 dx^2}=\frac{1}{2x^5},
\end{equation*}
\begin{equation*}
\frac{d^3 z}{1 \cdot 2 \cdots 5 dx^3}=-\frac{1}{2x^6}, \quad \frac{d^4z }{1 \cdot 2 \cdots 6 dx^4}=\frac{1}{2x^7}, \quad \frac{d^5 z}{1 \cdot 2 \cdots 7 dx^5}=-\frac{1}{2x^6} \quad \text{etc.},
\end{equation*}
it will be

\begin{equation*}
s=C-\frac{1}{2xx}+\frac{1}{2x^3}-\frac{3 \mathfrak{A}}{2x^4}+\frac{5 \mathfrak{B}}{2x^6}-\frac{7 \mathfrak{C}}{2x^6}+\text{etc.}
\end{equation*}
and hence, having put $x=1$, because of $s=1$, it will be

\begin{equation*}
C=1+\frac{1}{2}-\frac{1}{2}+\frac{3}{2}\mathfrak{A}-\frac{5}{2}\mathfrak{B}+\frac{7}{2}\mathfrak{C}-\frac{9}{2}\mathfrak{D}+\text{etc.}.
\end{equation*}
and this value of $C$ at the same time will show  the sum of the propounded series, if it is continued to infinity. Since the sums of the odd powers are not known as the sum of the even powers, this value of $C$ has to be defined from the known sum of some terms. Therefore, let $x=10$; it will be

\begin{equation*}
C=s+\frac{1}{2xx}-\frac{1}{2x^3}+\frac{2 \mathfrak{A}}{2x^4}-\frac{5 \mathfrak{B}}{2x^6}+\frac{7 \mathfrak{C}}{2x^8}-\text{etc.}
\end{equation*}
But, in order to perform the calculation more easily, note that 

\begin{alignat*}{9}
&~\frac{3 \mathfrak{A}}{2}&&= &&0.250000000000\\[3mm]
&~\frac{5 \mathfrak{B}}{2}&&= &&0.833333333333\\[3mm]
&~\frac{7 \mathfrak{C}}{2}&&= &&0.833333333333\\[3mm]
&~\frac{9 \mathfrak{D}}{2}&&= &&0.150000000000\\[3mm]
&\frac{11 \mathfrak{E}}{2}&&=&&0.416666666666\\[3mm]
&\frac{13 \mathfrak{F}}{2}&&=&&1.645280952380\\[3mm]
&\frac{15 \mathfrak{G}}{2}&&=&&8.750000000000\\[3mm]
&\frac{17 \mathfrak{H}}{2}&&=6&&0.283333333333\\[3mm]
& && &&\text{etc.}
\end{alignat*}
Therefore, the terms to be added to $s$ will become

\begin{equation*}
\renewcommand{\arraystretch}{2,0}
\setlength{\arraycolsep}{0mm}
\begin{array}{clclll}
 \dfrac{1}{2xx}                  &~=~&0.005000000000000000\\[2mm]
 \dfrac{3 \mathfrak{A}}{2x^4}    &~=~&0.000025000000000000\\[2mm]
 \dfrac{7 \mathfrak{C}}{2x^8}    &~=~&0.000000000833333333\\[2mm]
 \dfrac{11 \mathfrak{E}}{2x^{12}}&~=~&0.000000000000416666\\[2mm]
 \dfrac{13 \mathfrak{F}}{2x^{16}}&~=~&0.000000000000000875\\[2mm]
 \cline{2-3}
                                & &0.005025000833750875
\end{array}
\end{equation*}
but the terms to be subtracted are

\begin{equation*}
\renewcommand{\arraystretch}{2,0}
\setlength{\arraycolsep}{0mm}
\begin{array}{clclll}
 \dfrac{1}{2x^3}                    &~=~&0.005000000000000000\\[2mm]
 \dfrac{5 \mathfrak{B}}{2x^6}       &~=~&0.000000000833333333\\[2mm]
 \dfrac{9 \mathfrak{D}}{2x^{10}}    &~=~&0.000000000001500000\\[2mm]
 \dfrac{13 \mathfrak{F}}{2x^{14}}   &~=~&0.000000000000016452\\[2mm]
 \dfrac{17 \mathfrak{F}}{2x^{18}}   &=&0.000000000000000060\\[2mm]
 \cline{3-3}
                                   & &0.000500083348349845\\[2mm]
 \text{from}                       & &0.005025000833750875\\[2mm]
 \cline{2-3}
                                   & &0.004524917485401030\\[2mm]
        s                          &~=~&1.197531985674193251\\[2mm]
        \cline{3-3}
        C                          &~=~&1.202056903159594281\\ 
\end{array}
\end{equation*}

\paragraph*{§151}

If we continue this way, we will find the sum of all series of reciprocal powers expressed in decimal fractions.

\begin{alignat*}{19}
&1&&+\frac{1}{2^2}&&+\frac{1}{3^2}&&+\frac{1}{4^2}&&+\text{etc.}&&=1.6449340668482264&&=\frac{\quad ~ 2 \mathfrak{A} \quad ~~}{1 \cdot 2}&&\pi^2 \\[3mm]
&1&&+\frac{1}{2^3}&&+\frac{1}{3^3}&&+\frac{1}{4^3}&&+\text{etc.}&&=1.2020569031595942&&\\[3mm]
&1&&+\frac{1}{2^4}&&+\frac{1}{3^4}&&+\frac{1}{4^4}&&+\text{etc.}&&=1.0823232337111381&&=\frac{2^3 \mathfrak{B}}{1 \cdot 2 \cdot 3 \cdot 4}&&\pi^4 \\[3mm]
&1&&+\frac{1}{2^5}&&+\frac{1}{3^5}&&+\frac{1}{4^5}&&+\text{etc.}&&=1.0369277551433699&& \\[3mm]
&1&&+\frac{1}{2^6}&&+\frac{1}{3^6}&&+\frac{1}{4^6}&&+\text{etc.}&&=1.0173430619844491&&=\frac{ \quad ~ 2^5 \mathfrak{E} \quad ~}{1 \cdot 2 \cdots 6}&&\pi^6 \\[3mm]
&1&&+\frac{1}{2^7}&&+\frac{1}{3^7}&&+\frac{1}{4^7}&&+\text{etc.}&&=1.0083492773819288&& \\[3mm]
&1&&+\frac{1}{2^8}&&+\frac{1}{3^8}&&+\frac{1}{4^8}&&+\text{etc.}&&=1.0040773561979443&&=\frac{\quad ~ 2^7 \mathfrak{D}\quad ~}{1 \cdot 2 \cdots 8}&&\pi^8 \\[3mm]
&1&&+\frac{1}{2^9}&&+\frac{1}{3^9}&&+\frac{1}{4^9}&&+\text{etc.}&&=1.0020083928260822&&\\[3mm]
&1&&+\frac{1}{2^{10}}&&+\frac{1}{3^{10}}&&+\frac{1}{4^{10}}&&+\text{etc.}&&=1.0009945751278180&&=\frac{2^9 \mathfrak{E}}{1 \cdot 2 \cdots 10}&&\pi^{10} \\[3mm]
&1&&+\frac{1}{2^{11}}&&+\frac{1}{3^{11}}&&+\frac{1}{4^{11}}&&+\text{etc.}&&=1.0004941886041194&& \\[3mm]
&1&&+\frac{1}{2^{12}}&&+\frac{1}{3^{12}}&&+\frac{1}{4^{12}}&&+\text{etc.}&&=1.0002460865533080&&=\frac{2^{11} \mathfrak{F}}{1 \cdot 2 \cdots 12}&&\pi^{12} \\[3mm]
&1&&+\frac{1}{2^{13}}&&+\frac{1}{3^{13}}&&+\frac{1}{4^{13}}&&+\text{etc.}&&=1.0001227133475784&& \\[3mm]
&1&&+\frac{1}{2^{14}}&&+\frac{1}{3^{14}}&&+\frac{1}{4^{14}}&&+\text{etc.}&&=1.0000612481350587&&=\frac{2^{13} \mathfrak{G}}{1 \cdot 2 \cdots 14}&&\pi^{14} \\[3mm]
&1&&+\frac{1}{2^{15}}&&+\frac{1}{3^{15}}&&+\frac{1}{4^{15}}&&+\text{etc.}&&=1.0000305882363070&&\\[3mm]
&1&&+\frac{1}{2^{16}}&&+\frac{1}{3^{16}}&&+\frac{1}{4^{16}}&&+\text{etc.}&&=1.0000152822594086&&=\frac{2^{15} \mathfrak{H}}{1 \cdot 2 \cdots 16}&&\pi^{16} \\
& && && && && && \text{etc.}
\end{alignat*}

\paragraph*{§152}

From these vice versa the sums of those series consisting of the Bernoulli numbers can be exhibited. For, it will be

\begin{alignat*}{9}
&1+-\frac{1}{2}+\frac{\mathfrak{A}}{2}-\frac{\mathfrak{B}}{4}+\frac{\mathfrak{C}}{6}-\frac{\mathfrak{D}}{8}+\text{etc.}=0.57721~\text{etc.} \\[3mm]
&1+1-\frac{1}{2}+\mathfrak{A}-\mathfrak{B}+\mathfrak{C}-\mathfrak{D}+\text{etc.}=\frac{2 \mathfrak{A}}{1 \cdot 2}\pi^2 \\[3mm]
&1+\frac{1}{2}-\frac{1}{2}+\frac{2 \mathfrak{A}}{2}-\frac{5 \mathfrak{B}}{2}+\frac{7 \mathfrak{C}}{2}-\frac{9 \mathfrak{D}}{2}+\text{etc.}=1.2020~\text{etc.}\\[3mm]
&1+\frac{1}{3}-\frac{1}{2}+\frac{3 \cdot 4 \mathfrak{A}}{2 \cdot 3}-\frac{5 \cdot 5 \mathfrak{B}}{2 \cdot 3}+\frac{7 \cdot  8 \mathfrak{C}}{2 \cdot 3}-\frac{9 \cdot 10 \mathfrak{D}}{2 \cdot 3}+\text{etc.}=\frac{2^3 \mathfrak{B}}{1 \cdot 2 \cdot 3 \cdot 4}\pi^4 \\[3mm]
&1+\frac{1}{4}-\frac{1}{2}+\frac{3 \cdot 4 \cdot 5 \mathfrak{A}}{2 \cdot 3 \cdot 4}-\frac{5 \cdot 6 \cdot 7 \mathfrak{B}}{2 \cdot 3 \cdot 4}+\frac{7 \cdot 8 \cdot 9 \mathfrak{C}}{2 \cdot 3 \cdot 4}-\frac{ 9 \cdot 10 \cdot 11 \mathfrak{D}}{2 \cdot 3 \cdot 4}+\text{etc.}=1.0369~\text{etc.} \\[3mm]
&1+\frac{1}{5}-\frac{1}{2}+\frac{3 \cdot 4 \cdot 5 \cdot 6 \mathfrak{A}}{2 \cdot 3 \cdot 4 \cdot 5}-\frac{  4 \cdot 5 \cdot 6 \cdot 7 \mathfrak{B}}{2 \cdot 3 \cdot 4 \cdot 5}+\frac{5 \cdot 6 \cdot 7 \cdot 8 \mathfrak{C}}{2 \cdot 3 \cdot 4 \cdot 5}-\text{etc.}=\frac{2^5 \mathfrak{C}}{1 \cdot 2  \cdots 6}\pi^6 \\
& \qquad \qquad \qquad \qquad \qquad \text{etc.}
\end{alignat*}
Therefore, each second of these series can be summed by means of the quadrature of the circle; it is not known to this day, on which transcendental quantity the remaining series depend; for, they cannot be reduced to powers of $\pi$ with odd exponents, such that the  coefficients would be rational numbers. But that this at least becomes approximately clear, of what nature the coefficients of powers of $\pi$ will be for odd exponents, we added the following table.

\begin{alignat*}{9}
&1&&+\frac{1}{2}&&+\frac{1}{3}&&+\frac{1}{4}&&+\text{etc. to infinity}&&=\frac{\pi}{0.0000}&&=\infty \\[3mm]
&1&&+\frac{1}{2^2}&&+\frac{1}{3^2}&&+\frac{1}{4^2}&&+\text{etc. to infinity}&&=\frac{\pi^2}{6.0000}&&~\text{exactly}\\[3mm]
&1&&+\frac{1}{2^3}&&+\frac{1}{3^3}&&+\frac{1}{4^3}&&+\text{etc. to infinity}&&=\frac{\pi^3}{25.79436}&&~\text{approximately}\\[3mm]
&1&&+\frac{1}{2^4}&&+\frac{1}{3^4}&&+\frac{1}{4^4}&&+\text{etc. to infinity}&&=\frac{\pi^2}{90.00000}&&~\text{exactly}\\[3mm]
&1&&+\frac{1}{2^5}&&+\frac{1}{3^5}&&+\frac{1}{4^5}&&+\text{etc. to infinity}&&=\frac{\pi^5}{295.1215}&&~\text{approximately}\\[3mm]
&1&&+\frac{1}{2^6}&&+\frac{1}{3^6}&&+\frac{1}{4^6}&&+\text{etc. to infinity}&&=\frac{\pi^6}{945.000}&&~\text{exactly}\\[3mm]
&1&&+\frac{1}{2^7}&&+\frac{1}{3^7}&&+\frac{1}{4^7}&&+\text{etc. to infinity}&&=\frac{\pi^2}{2995.284}&&~\text{approximately}\\[3mm]
&1&&+\frac{1}{2^8}&&+\frac{1}{3^8}&&+\frac{1}{4^8}&&+\text{etc. to infinity}&&=\frac{\pi^2}{9450.0000}&&~\text{exactly}\\[3mm]
&1&&+\frac{1}{2^9}&&+\frac{1}{3^9}&&+\frac{1}{4^9}&&+\text{etc. to infinity}&&=\frac{\pi^9}{29749.35}&&~\text{approximately}\\
& && && && &&\text{etc.}
\end{alignat*}

\paragraph*{§153}

From this source the series of Bernoulli numbers

\begin{alignat*}{19}
&1 \quad && 2 \quad && 3 \quad && 4 \quad && 5 \quad && 6 \quad && 7 \quad && 8 \quad && 9 \\[3mm]
&\mathfrak{A}, \quad && \mathfrak{B}, \quad && \mathfrak{C}, \quad && \mathfrak{D}, \quad && \mathfrak{E}, \quad && \mathfrak{F}, \quad && \mathfrak{G}, \quad && \mathfrak{H}, \quad && \mathfrak{I} \quad \text{etc.},
\end{alignat*}
even though it seems to be rather irregular, can be interpolated, i.e. the terms lying in the middle between any two  consecutive ones can be assigned; for, if the term falling in the middle between the first $\mathfrak{A}$ and the second $\mathfrak{B}$ or the one corresponding to the index $1\frac{1}{2}$ was $=p$, it will be

\begin{equation*}
1+\frac{1}{2^3}+\frac{1}{3^3}+\text{etc.}=\frac{2^2 p}{1 \cdot 2 \cdot 3}\pi^3
\end{equation*} 
and hence

\begin{equation*}
p=\frac{3}{2\pi^3}\left(1+\frac{1}{2^3}+\frac{1}{3^3}+\text{etc.}\right)=0.05815227.
\end{equation*}
If in like manner the term falling in the middle between $\mathfrak{B}$ and $\mathfrak{C}$ or the one corresponding to the index $2\frac{1}{2}$ is put $=q$, since it will be

\begin{equation*}
1+\frac{1}{2^5}+\frac{1}{3^5}+\text{etc.}=\frac{2^4 q}{1 \cdot 2 \cdot 3 \cdot 4 \cdot 5}\pi^5,
\end{equation*}
it will also be

\begin{equation*}
q=\frac{15}{2\pi^5}\left(1+\frac{1}{2^5}+\frac{1}{3^5}+\text{etc.}\right)=0.02541327.
\end{equation*}
Therefore, if the sums of these series, in which the exponents of the powers are odd numbers, could be exhibited, then also the series of the Bernoulli numbers could be interpolated.

\paragraph*{§154}

Now, let us put $z=\frac{1}{nn+xx}$ and find the sum of this series

\begin{equation*}
s=\frac{1}{nn+1}+\frac{1}{nn+4}+\frac{1}{nn+9}+\cdots + \frac{1}{nn+xx}.
\end{equation*}
Since $\int zdx = \int \frac{dx}{nn+xx}$, it will be

\begin{equation*}
\int zdx = \frac{1}{n}\arctan \frac{x}{n}.
\end{equation*}
Put $\operatorname{arccot} \frac{x}{n}=u$; it will be

\begin{equation*}
\int z dx = \frac{1}{n}\left(\frac{\pi}{2}-u\right)
\end{equation*}
and

\begin{equation*}
\frac{x}{n}=\cot u= \frac{\cos u}{\sin u} \quad \text{and} \quad \frac{nn+xx}{nn}=\frac{1}{\sin^2 u} \quad \text{and} \quad z= \frac{\sin^2 u}{nn} \quad \text{and} \quad \frac{dx}{n}=-\frac{du}{\sin^2 u},
\end{equation*}
whence

\begin{equation*}
du=-\frac{dx \sin^2 u}{n}.
\end{equation*}
Therefore, the differentials of $z$ will be found this way
 \begin{equation*}
dz= \frac{2 du \sin u \cdot \cos u}{nn}=-\frac{dx \sin^2 u \sin 2u}{n^3} \quad \text{and} \quad \frac{dz}{dx}=-\frac{\sin^2 u \cdot \sin 2u}{n^3}, 
 \end{equation*}
 \begin{equation*}
 \frac{ddz}{dx^2}=-\frac{du(\sin u \cdot \cos u \cdot \sin 2u +\sin^2 u \cdot \cos 2u)}{n^3}=\frac{dx \sin^3 \cdot 3u}{n^4}
 \end{equation*}
 and
 
 \begin{equation*}
 \frac{ddz}{2dx^2}=\frac{\sin^3 u \cdot \sin 3u }{n^4}.
 \end{equation*}
In the same way, as we already found above [§ 87] for the same case, it will be

\begin{equation*}
\frac{d^3 z}{2 \cdot 3 dx^3}=-\frac{\sin^4 \cdot \sin 4u}{n^5}, \quad \frac{d^4 z}{2 \cdot 3 \cdot 4 dx^4}=\frac{\sin^5 \cdot \sin 5u}{n^6} \quad \text{etc.},
\end{equation*} 
from which the sum in question will be formed

\begin{equation*}
s= \frac{\pi }{2n}-\frac{u}{n}+\frac{\sin u \cdot \sin u}{2nn}-\frac{\mathfrak{A}}{2}\cdot \frac{\sin^2 u \cdot \sin 2u}{n^3}+\frac{\mathfrak{B}}{4}\cdot \frac{\sin^4 u \cdot \sin 4u}{n^5}
\end{equation*}
\begin{equation*}
-\frac{\mathfrak{C}}{6}\cdot \frac{\sin^6 u \cdot \sin 6u }{n^7}+\frac{\mathfrak{D}}{8}\cdot \frac{\sin^8 u \cdot \sin 8u}{n^9}-\text{etc.}+\text{Const.}
\end{equation*}
If here in order to determine the constant one sets $x=0$, in which case $s=0$, it will be $\cot u=0$ and hence $u$ the angle of $90^{\circ}$ and therefore $\sin u=1$, $\sin 2u =0$, $\sin 4u =0$, $\sin 6u =0$ etc.; therefore, it seems that 

\begin{equation*}
0=\frac{\pi}{2n}-\frac{\pi}{2n}+\frac{1}{2nn}+C \quad \text{and hence} \quad C=-\frac{1}{2nn};
\end{equation*}
but on the other hand it is to be noted, even though the remaining terms vanish, that nevertheless, since the coefficients $\mathfrak{A}$, $\mathfrak{B}$, $\mathfrak{C}$ etc. eventually grow to infinity, their sum can be finite.

\paragraph*{§155}
To determine this constant correctly, let us put that  $x=\infty$; for, we defined the sum of this series running to infinity already in the \textit{Introdctio} and showed it to be

\begin{equation*}
=-\frac{1}{2nn}+\frac{\pi }{2n}+\frac{\pi}{n(e^{2 \pi b}-1)}.
\end{equation*}
But having put $x=\infty$, it will be $u=0$ and hence $\sin u=0$ and at the same time the sines of all multiple arcs will vanish. But since in this series the powers of $\sin u$ grow, the divergence of the series is no obstruction for the vanishing of the value in this case. Therefore, it will become $s=\frac{\pi}{2n}+C$; hence, it will be

\begin{equation*}
\frac{\pi}{2n}+C=-\frac{1}{2nn}+\frac{\pi}{2n}+\frac{\pi}{n(e^{2n \pi}-1)} \quad \text{and} \quad C=-\frac{1}{2nn}+\frac{\pi}{n(e^{2n \pi}-1)}.
\end{equation*}
Therefore, the sum of the series in question will be

\begin{equation*}
s= \frac{\pi }{2n}-\frac{u}{n}-\frac{1}{2nn}+\frac{\sin^2 u}{2nn}-\frac{\mathfrak{A}}{2}\cdot \frac{\sin^2 u \cdot \sin 2u}{n^3}
\end{equation*}
\begin{equation*}
+\frac{\mathfrak{B}}{4}\cdot \frac{\sin^4 u \cdot \sin 4u}{n^5}- \frac{\mathfrak{C}}{6}\cdot \frac{\sin^6 u \cdot \sin 6u}{n^7}+\text{etc.}+\frac{\pi}{n(e^{2n \pi}-1)}.
\end{equation*}
Where it is to be noted, if $n$ was a modestly large number, that the last term $\frac{\pi}{n(e^{2 \pi n}-1)}$ will become so small that it can be neglected.

\paragraph*{§156}

Let us put $x=n$ so that
\begin{equation*}
s=\frac{1}{nn+1}+\frac{1}{nn+4}+\frac{1}{nn+9}+\cdots + \frac{1}{nn+nn}.
\end{equation*}
Then it will be $\cot u=1$ and $u=45^{\circ}=\frac{\pi}{4}$. Therefore, one will have $\sin u=\frac{1}{\sqrt{2}}$, $\sin 2u =1$, $\sin 4u=0$, $\sin 6u =-1$, $\sin 8u =0$, $\sin 10 u=1$ etc. Therefore, it will be

\begin{equation*}
s= \frac{\pi}{4n}-\frac{1}{2nn}+\frac{1}{4nn}-\frac{\mathfrak{A}}{2 \cdot 2 n^3}+\frac{\mathfrak{C}}{6 \cdot 8 n^7}-\frac{\mathfrak{E}}{10 \cdot 2^5 n^{11}}+\frac{\mathfrak{G}}{14 \cdot 2^7 n^{15}}-\text{etc.}+\frac{\pi}{n(e^{2n \pi}-1)},
\end{equation*}
in which expression only each second Bernoulli number appears. Therefore, if the value of $s$ was already found by actual calculation, hence the quantity $\pi$ can be defined; for, it will be

\begin{equation*}
\pi =4ns+\frac{1}{n}+\frac{\mathfrak{A}}{1 \cdot n^2}-\frac{\mathfrak{C}}{3 \cdot 2^2 n^6}+\frac{\mathfrak{E}}{5 \cdot 2^4 n^{10}}-\frac{\mathfrak{G}}{7 \cdot 2^6 n^{14}}+\text{etc.}-\frac{4 \pi}{e^{2n \pi}-1}.
\end{equation*}
For, even though  the last term contains $\pi$, nevertheless, since it is so small, it suffices to determine the value of $\pi$ approximately.

\subsection*{Example}

Let $n=5$; it will be

\begin{equation*}
s=\frac{1}{26}+\frac{1}{29}+\frac{1}{34}+\frac{1}{41}+\frac{1}{50};
\end{equation*} 
these terms, if they are actually added, will give

\begin{equation*}
s=0.146746306590549494;
\end{equation*}
Hence  the terms will be

\begin{equation*}
\renewcommand{\arraystretch}{2,0}
\setlength{\arraycolsep}{0mm}
\begin{array}{rcl}
4ns                    &~=~&2.93492611381098988 \\[2mm]
\dfrac{1}{n}            &~=~&0.20000000000000000 \\[2mm]
\dfrac{\mathfrak{A}}{nn}&~=~&0.00666666666666666 \\[2mm]
\cline{3-3}
                       & &3.14159278047765654 \\[2mm]
 \dfrac{\mathfrak{C}}{3 \cdot 2^2 \cdot n^6}   &=&0.00000012698412698 \\[2mm]
 \cline{3-3}        
                       & &3.14159265349352956 \\[2mm]
 \dfrac{\mathfrak{D}}{5 \cdot 2^4 \cdot n^{10}}&=&0.00000000009696969                                 \\[2mm]
 \cline{3-3}
                       & &3.14159265359049925 \\[2mm]
  \dfrac{\mathfrak{E}}{7 \cdot 2^6 \cdot n^{14}}&=&0.00000000000042666 \\
  \cline{3-3}
                       & &3.14159265359007259 \\[2mm]
   \dfrac{\mathfrak{F}}{9 \cdot 2^8 \cdot n^{18}}&=&0.00000000000000625 \\
   \cline{3-3}
                       & &3.14159265359007884.                                            
\end{array}
\end{equation*}
This value already comes so close to the true value that one has to wonder why by such a simple calculation one can get this far. This expression is indeed a little bit larger than the correct value; for, one has still to subtract $\frac{4 \pi}{e^{2 n \pi}-1}$, whose value, as long as $\pi$ is sufficiently accurately known, can be exhibited; this will be achieved by means of logarithms. \\[2mm]
Since $\pi \log e=1.3643763538$, it will be

\begin{equation*}
\log e^{2n \pi}=10 \pi \log e=13.6437635.
\end{equation*}
Since

\begin{equation*}
\frac{4 \pi}{e^{2n \pi}-1}=\frac{4\pi}{e^{2n\pi}}+\frac{4 \pi}{e^{4 n\pi}}+\text{etc.},
\end{equation*}
it is easily understood that for our calculation it suffices to have taken the first term. Therefore, let us  increase the characteristic by the number $17$, since we have the same number of decimal places; it will be

\begin{equation*}
\renewcommand{\arraystretch}{1,4}
\setlength{\arraycolsep}{0mm}
\begin{array}{lllrlr}
                         &  & \log \pi & ~=~ & &17.4971498 \\
                         &  & \log 4   & ~=~ &  &0.6020600 \\
 \cline{6-6}
                &       &   &   & &18.0992098 \\
                & \text{subtr. }& \log e^{2n \pi} &~=~& &13.6437635 \\           
 \text{Therefore} \quad     &       &   &   & &4.4554463 \\
 \cline{6-6}        
                &               & \dfrac{4 \pi}{e^{2n \pi}}&~=~& &28539 \\
  \text{subtract from} \\
                     &   &       & & &3.14159265369007884 \\
   \cline{6-6}                  
   \text{it will be} &   &   \pi &~=~& &3.14159265358979345                                    
\end{array}
\end{equation*}
which expression deviates from the true value just at the penultimate digit; this is not surprising, since we would have  to subtract the term $\frac{\mathfrak{L}}{11 \cdot 2^{10}\cdot n^{22}}$, which gives $22$,  and so not even the last figure would have been wrong. Moreover, it is understood, if  we would have taken a larger number for $n$, e.g.  $10$, that   the circumference $\pi$ could have been found up to $25$ and more digits very easily.

\paragraph*{§157}

Now let us also substitute transcendental functions of $x$ for $z$ and let be $z=\log x$  taking hyperbolic logarithms, since ordinary ones are easily reduced to them, and let 

\begin{equation*}
s= \log 1+ \log 2+ \log 3+ \log 4 +\cdots + \log x.
\end{equation*}
Therefore, since $z= \log x$, it will be

\begin{equation*}
\int zdx =x \log x -x;
\end{equation*}
for, its differential gives $dx \log x$. Furthermore, 

\begin{equation*}
\frac{dz}{dx}=\frac{1}{x}, \quad \frac{ddz}{dx^2}=-\frac{1}{x^2}, \quad \frac{d^3z}{1 \cdot 2 dx^3}=\frac{1}{x^3}, \quad \frac{d^4 z}{1 \cdot 2 \cdot 3 dx^4}=-\frac{1}{x^4}, \quad \frac{d^5 z}{1 \cdot 2 \cdot 3 \cdot 4 dx^5}=\frac{1}{x^5} \quad \text{etc.}
\end{equation*}
Therefore, one will conclude that 

\begin{equation*}
s=x \log x-x +\frac{1}{2}\log x+\frac{\mathfrak{A}}{1 \cdot 2 x}-\frac{\mathfrak{B}}{3 \cdot 4 x^3}+\frac{\mathfrak{C}}{5 \cdot 6 x^5}-\frac{\mathfrak{D}}{7 \cdot 8 x^7}+\text{etc.}+\text{Const.}
\end{equation*}
But this constant, putting $x=1$, since $s=\log 1=0$, will be defined as follows

\begin{equation*}
C=1-\frac{\mathfrak{A}}{1 \cdot 2}+\frac{\mathfrak{B}}{3 \cdot 4}-\frac{\mathfrak{C}}{5 \cdot 6}+\frac{\mathfrak{D}}{7 \cdot 8}-\text{etc.},
\end{equation*}
which series, because of the too strong divergence, is inept to find the value of $C$ at least approximately.

\paragraph*{§158}

But we will not only find an approximate value, but even the true value of $C$, if we consider Wallis's expression found for the value of $\pi$ and demonstrated in the \textit{Introductio}, which was

\begin{equation*}
\frac{\pi}{2}=
\renewcommand{\arraystretch}{1,4}
\setlength{\arraycolsep}{0.1mm}
\begin{array}{cccccccccccccccccccccccccccc}
2 &\cdot & 2 & \cdot &4 & \cdot & 4 & \cdot & 6 & \cdot & 6 &\cdot & 8 \cdot & 8 \cdot & 10 & \cdot & 10 & \cdot &12 &\cdot & \text{etc.} \\
\hline
1 &\cdot & 3 & \cdot &3 & \cdot & 5 & \cdot & 5 & \cdot & 7 &\cdot & 7 \cdot & 9 \cdot & 9 & \cdot & 11 & \cdot &11 & \cdot & \text{etc.} 
\end{array}
\end{equation*}
For, by taking logarithms, it will be

\begin{equation*}
\log \pi - \log 2 = 2 \log 2 +2 \log 4 +2 \log 6+ 2 \log 8 +2 \log 10+ \log 12 +\text{etc.} 
\end{equation*}
\begin{equation*}
- \log 1 -2 \log 3 -2 \log 5-2 \log 7 -2 \log 9 - 2 \log 11 -\text{etc.}
\end{equation*}
Therefore, in the assumed series let us put $x=\infty$, and since

\begin{alignat*}{9}
&\log 1 &&+\log 2 &&+\log 3 && +\log 4 &&+ \cdots &&+ \log x =C+\left(x+\frac{1}{2}\right)\log x -x, 
\intertext{it will be}
&\log 1 &&+\log 2 &&+\log 3 && +\log 4 &&+ \cdots &&+ \log 2x =C+\left(2x+\frac{1}{2}\right)\log 2x -2x
\intertext{and}
&\log 2 &&+\log 4 &&+\log 6 && +\log 8 &&+ \cdots &&+ \log 2x =C+\left(x+\frac{1}{2}\right)\log x -x \log 2-x, 
\intertext{hence}
&\log 1 &&+\log 3 &&+\log 5 && +\log 7 &&+ \cdots &&+ \log (2x-1) =C+\left(x+\frac{1}{2}\right)\log 2 -x. \\
\end{alignat*}
Therefore, since 

\begin{alignat*}{9}
& \log \frac{\pi}{2}&&= 2 \log 2 &&+2 \log 4 &&+2 \log 6 &&+ \cdots &&+2 \log 2x -\log 2x \\
& &&-2 \log 1 &&-2\log 3 &&-2 \log 5 &&- \cdots &&- 2 \log (2x-1),
\end{alignat*}
having put $x=\infty$, it will be

\begin{equation*}
\log \frac{\pi}{2}=2 C+(2x+1)\log x +2x \log 2 -2x - \log 2 -\log x -2x \log x -(2x+1)\log 2 +2x
\end{equation*}
and hence

\begin{equation*}
\log \frac{\pi}{2}=2 C- 2 \log 2, \quad \text{therefore} \quad 2C= \log 2 \pi \quad \text{and} \quad C=\frac{1}{2}\log 2 \pi,
\end{equation*}
whence in decimal fractions it is found to be

\begin{equation*}
C=0.9189385332046727417803297,
\end{equation*}
and at the same time the following series is summed

\begin{equation*}
1-\frac{\mathfrak{A}}{1 \cdot 2}+\frac{\mathfrak{B}}{3 \cdot 4}-\frac{\mathfrak{C}}{5 \cdot 6}+\frac{\mathfrak{D}}{7 \cdot 8}-\frac{\mathfrak{E}}{9 \cdot 10}+\text{etc.}=\frac{1}{2}\log 2 \pi.
\end{equation*}

\paragraph*{§159}

Now, having found this constant $C=\frac{1}{2}\log 2 \pi$, the sum of any number of logarithms from this series $\log 1 +\log 2 + \log 3 +\text{etc.}$ can be exhibited. For, if one puts

\begin{equation*}
s= \log 1 +\log 2 +\log 3 +\log 4 +\cdots +\log x,
\end{equation*}
it will be

\begin{equation*}
s=\frac{1}{2}\log 2 \pi +\left(x+\frac{1}{2}\right)\log x-x+\frac{\mathfrak{A}}{1 \cdot 2x}-\frac{\mathfrak{B}}{3 \cdot 4 x^3}+\frac{\mathfrak{C}}{5 \cdot 6x^5}-\frac{\mathfrak{D}}{7 \cdot 8 x^7}+\text{etc.},
\end{equation*}
if the propounded logarithms were hyperbolic; but if common logarithms are propounded, then in the terms $\frac{1}{2}\log 2 \pi+\left(x+\frac{1}{2}\right)\log x$ for $\log 2\pi$ and $\log x$ one also has to take common logarithms, but the remaining terms of the series

\begin{equation*}
-x+\frac{\mathfrak{A}}{1 \cdot 2 x}-\frac{\mathfrak{B}}{3 \cdot 4 x^3}+\text{etc.}
\end{equation*}
have to be multiplied by $0.434294481903251827=n$. Therefore, in this case for common logarithms it will be

\begin{equation*}
\renewcommand{\arraystretch}{1,2}
\setlength{\arraycolsep}{0.1mm}
\begin{array}{rrl}
\log \pi    &~=~&0.497149872694133854351268 \\[2mm]
\log 2      &~=~&0.301029995663981195213738 \\[2mm]
\cline{3-3}
\log 2 \pi            &~=~&0.798179868358115049565006 \\[2mm]
\frac{1}{2} \log 2 \pi&~=~&0.399089934179057524782503.
\end{array}
\end{equation*}

\subsection*{Example}

\textit{Let the aggregate of thousand tabulated logarithms be in question}

\begin{equation*}
s= \log 1 +\log 2+ \log 3 +\cdots +\log 1000.
\end{equation*}
Therefore, it will be $x=1000$ and

\begin{equation*}
\renewcommand{\arraystretch}{2,0}
\setlength{\arraycolsep}{0.1mm}
\begin{array}{lrrrr}
 &\log x &~=~&                           3.000000000000 \\[2mm]
\text{whence} \qquad & x \log x   &~=~&3000.000000000000 \\[2mm]
                 &\dfrac{1}{2}\log x &~=~&1.500000000000 \\[2mm]
            & \dfrac{1}{2}\log 2 \pi &~=~&0.399089341790 \\[2mm]
   \cline{4-4} 
            &                       & &3001.8990899341790 \\[2mm]
      \text{subtr. }& nx            &~=~&2567.6046080309272 \\[2mm]                   
 \text{Furthermore,}  \\
 & \dfrac{n \mathfrak{A}}{1 \cdot 2 x}&~=~&                 0.0000361912068 \\
 \text{subtr.} & \dfrac{n \mathfrak{B}}{3 \cdot 4 x^3}&~=~& 0.0000000000012 \\[2mm]
 & \text{add. }&~=~&              2567.6046080309272 \\[2mm]
 \cline{4-4}
 \text{sum in question} &s &~=~& 2567.6046442221328.                   
\end{array}
\end{equation*}
Therefore, since $s$ is the logarithm of the product of the numbers

\begin{equation*}
1 \cdot 2 \cdot 3 \cdot 4 \cdot 5 \cdot 6 \cdots 1000,
\end{equation*}
it is clear that this product, if actually multiplied, has $2568$ digits and  the first numbers will be $4023872$, which will be followed by $2561$ other numbers.

\paragraph*{§160}

Therefore, by means of this summation of logarithms the products of arbitrary many factors proceeding according to the natural numbers, can be assigned approximately. An application of this is the solution of the problem, in which the middle term or the largest term in any power of the binomial $(a+b)^m$ is in question, where it is certainly to be noted, if $m$ is an odd number, that two equal middle terms are given, which, if they are added, yield the middle term in the following even power. Because hence the largest coefficient in any even power is twice as large as the middle coefficient in the preceding odd power, it will be sufficient to have determined the largest middle term for the even powers. Therefore, let $m=2n$ and the middle coefficient will be expressed as follows

\begin{equation*}
\frac{2n(2n-1)(2n-2)(2n-3) \cdots (n+1)}{1 \cdot 2 \cdot 3 \cdot 4 \cdots n}.
\end{equation*}
Let us call this middle coefficient in question $=u$ and  one will be able to represent it this way

\begin{equation*}
u=\frac{1 \cdot 2 \cdot 3 \cdot 4 \cdot 5 \cdots 2n}{(1 \cdot 2 \cdot 3 \cdot 4 \cdots 2n)^2},
\end{equation*}
and, having taken logarithms, it will be

\begin{equation*}
\log u= \log 1+\log 2+\log 3 +\log 4 +\log 5 +\cdots + \log 2n
\end{equation*}
\begin{equation*}
-2 \log 1 -2 \log 2 - 2 \log 3 - 2 \log 4 - 2 \log 5 - \cdots - 2 \log 2n.
\end{equation*}

\paragraph*{§161}

Now, by assuming these logarithms to be hyperbolic logarithms, it will be

\begin{equation*}
\log 1+ \log 2 + \log 3 + \log 4 + \cdots + \log 2n =\frac{1}{2}\log 2 \pi +\left(2n+\frac{1}{2}\right)\log n+\left(2n +\frac{1}{2}\right)\log 2 -2n
\end{equation*}
\begin{equation*}
+\frac{\mathfrak{A}}{1 \cdot 2 \cdot 2 n}-\frac{\mathfrak{B}}{3 \cdot 4 \cdot 2^3 n^3}+\frac{\mathfrak{C}}{5 \cdot 6 \cdot 2^5 n^5}-\text{etc.}
\end{equation*}
and

\begin{equation*}
2 \log 1 + 2 \log 2 + 2 \log 3 +2 \log 4 +\cdots + 2 \log n
\end{equation*}
\begin{equation*}
= \log 2 \pi +(2n+1)\log n-2n +\frac{2 \mathfrak{A}}{1 \cdot 2 n}-\frac{2 \mathfrak{B}}{3 \cdot 4 n^3}+\frac{2 \mathfrak{C}}{5 \cdot 6 n^5}-\text{etc.},
\end{equation*}
having subtracted which expression from the other it will remain

\begin{equation*}
\log u= -\frac{1}{2}\log \pi -\frac{1}{2}\log n +2n \log 2 +\frac{\mathfrak{A}}{1 \cdot 2 \cdot 2 n}-\frac{\mathfrak{B}}{3 \cdot 4 \cdot 2^3 n^3}+\frac{\mathfrak{C}}{5 \cdot 6 \cdot 2^5 n^5}-\text{etc.}
\end{equation*}
\begin{equation*}
-\frac{2 \mathfrak{A}}{1 \cdot 2 n}+\frac{2 \mathfrak{B}}{3 \cdot 4 n^3}-\frac{2 \mathfrak{C}}{5 \cdot 6 n^5}+\text{etc.};
\end{equation*}
by collecting each two terms it will be

\begin{equation*}
\log u= \log \frac{2^{2n}}{\sqrt{n \pi}}-\frac{3 \mathfrak{A}}{1 \cdot 2 \cdot 2 n}+\frac{15 \mathfrak{B}}{3 \cdot 4 \cdot 2^3 n^3}-\frac{63 \mathfrak{C}}{5 \cdot 6 \cdot 2^5 n^5}+\frac{255 \mathfrak{D}}{7 \cdot 8 \cdot 2^7 n^7}-\text{etc.}
\end{equation*}
Let

\begin{equation*}
\frac{3 \mathfrak{A}}{1 \cdot 2 \cdot 2^2 n^2}-\frac{15 \mathfrak{B}}{3 \cdot 4 \cdot 2^4 n^4}+\frac{63 \mathfrak{C}}{5 \cdot 6 \cdot 2^6 n^6}-\frac{255 \mathfrak{D}}{7 \cdot 8 \cdot 2^8 n^8}+\text{etc.}
\end{equation*}
\begin{equation*}
= \log \left(1+\frac{A}{2^2 n^2}+\frac{B}{2^4 n^4}+\frac{C}{2^6 n^6}+\frac{D}{2^8 n^8}+\text{etc.} \right);
\end{equation*}
it will be

\begin{equation*}
\log u= \log \frac{2^{2n}}{\sqrt{n \pi}}-2n \log \left(1+\frac{A}{2^2 n^2}+\frac{B}{2^4 n^4}+\frac{C}{2^6 n^6}+\text{etc.}\right)
\end{equation*}
and hence

\begin{equation*}
\frac{u^{2n}}{\left(1+\frac{A}{2^2 n^2}+\frac{B}{2^4 n^4}+\frac{C}{2^6 n^6}+\text{etc.}\right)^{2n}\sqrt{n \pi}}.
\end{equation*}
Having put $2n=m$, it will be

\begin{alignat*}{9}
&\log \bigg(1 &&+\frac{A}{2^2n^2}&&+\frac{B}{2^4 n^4}&&+\frac{C}{2^6 n^6}&&+\frac{D}{2^8 n^8}&&+\text{etc.} \bigg) \\[3mm]
&=\frac{A}{m^2}&&+~\frac{B}{m^4}&&+~\frac{C}{m^6}&&+~\frac{D}{m^8}&&+~\frac{E}{m^{10}}&&+\text{etc.} \\[3mm]
& && -\frac{A^2}{2m^4}&&-\frac{AB}{m^6}&&-\frac{AC}{m^8}&&-\frac{AD}{m^{10}}&&-\text{etc.} \\[3mm]
& &&  && &&-\frac{BB}{2 m^8}&&-\frac{BC}{m^{10}}&&-\text{etc.} \\[3mm]
& &&  &&+\frac{A^3}{3 m^6}&&+\frac{A^2 B}{m^8}&&+\frac{A^2 C}{m^{10}}&&+\text{etc.} \\[3mm]
& &&  && && &&+\frac{AB^2}{m^{10}}&&+\text{etc.} \\[3mm]
& &&  && &&-\frac{A^4}{4 m^8}&&-\frac{A^3 B}{m^{10}}&&- \text{etc.} \\[3mm]
& &&  &&  && &&+\frac{A^5}{5 m^{10}}&&+\text{etc.};
\end{alignat*}
since this expression has to be equal to this one

\begin{equation*}
\frac{2 \mathfrak{A}}{1 \cdot 2 m^2}-\frac{15 \mathfrak{B}}{3 \cdot 4}+\frac{63 \mathfrak{C}}{5 \cdot 6 m^6}-\frac{255 \mathfrak{D}}{7 \cdot 8 m^8}+\text{etc.},
\end{equation*}
it will be

\begin{alignat*}{9}
&A&&= \frac{3 \mathfrak{A}}{1 \cdot 2} \\[3mm]
&B&&= \frac{A^2}{2}-\frac{25 \mathfrak{B}}{3 \cdot 4}\\[3mm]
&C&&=AB- \frac{1}{3}A^3 +\frac{63 \mathfrak{C}}{5 \cdot 6} \\[3mm]
&D&&=AC+\frac{1}{2}B^2-A^2B+\frac{1}{4}A^4-\frac{255 \mathfrak{D}}{7 \cdot 8}\\[3mm]
&E&&=AD+BC-A^2C-AB^2+A^3B-\frac{1}{5}A^5+\frac{1023 \mathfrak{E}}{9 \cdot 10}\\[3mm]
& && \qquad \qquad \text{etc.}
\end{alignat*}

\paragraph*{§162}

Since $\mathfrak{A}=\frac{1}{6}$, $\mathfrak{B}=\frac{1}{30}$, $\mathfrak{C}=\frac{1}{42}$, $\mathfrak{D}=\frac{1}{30}$, $\mathfrak{E}=\frac{5}{66}$, it will be

\begin{equation*}
A=\frac{1}{4}, \quad B=-\frac{1}{96}, \quad C=\frac{27}{640}, \quad D=-\frac{90031}{2^{11}\cdot 3^2 \cdot 5 \cdot 7} \quad \text{etc.}
\end{equation*}
Therefore, one finds

\begin{equation*}
u=\frac{2^{2n}}{\left(1+\frac{1}{2^4n^2}-\frac{1}{2^9 \cdot 3 n^4}+\frac{27}{2^{13}\cdot 5n^6}-\frac{90031}{2^{19}\cdot 3^2 \cdot 5 \cdot 7n^8}+\text{etc.}\right)^{2n}\sqrt{n \pi}}
\end{equation*}
or

\begin{equation*}
u=\frac{2^{2n}\left(1-\frac{1}{2^4 n^2}+\frac{7}{2^9 \cdot 3n^4}-\frac{121}{2^{13} \cdot 3 \cdot 5n^6}+\frac{107489}{2^{19}\cdot 3^2 \cdot 5 \cdot 7n^8}- \text{etc.} \right)^{2n}}{\sqrt{n \pi}},
\end{equation*}
or if this expansion of the series is actually done, it will approximately be

\begin{equation*}
u=\frac{2^{2n}}{\sqrt{n \pi}\left(1+\frac{1}{4n}+\frac{1}{32n^2}-\frac{1}{128n^5}-\frac{5}{16 \cdot 128 n^4}+\text{etc.}\right)};
\end{equation*}
therefore, the middle term in $(1+1)^{2n}$ will have the same ratio to the sum of all terms $2^{2n}$ 

\begin{equation*}
\text{as} \quad 1 \quad \text{to} \quad \sqrt{n \pi \left(1+\frac{1}{4n}+\frac{1}{32n^2}-\frac{1}{128 n^3}-\frac{5}{16 \cdot n^4}+\text{etc.}\right)};
\end{equation*}
or,  for the sake of brevity having put $4n= \nu$, this ratio will be

\begin{equation*}
\text{as} \quad \text{to} \quad \sqrt{n \pi \left(1+\frac{1}{\nu}+\frac{1}{2 \nu^2}-\frac{1}{2 \nu^3}-\frac{5}{8 \nu^4}+\frac{23 }{8 \nu^5}+\frac{53}{16 \nu^6}-\text{etc.}\right)}.
\end{equation*}

\subsection*{Example 1}

\textit{Let the middle term in the expanded binomial $(a+b)^{10}$ be in question, which is known to be}

\begin{equation*}
=\frac{10 \cdot 9 \cdot 8 \cdot 7 \cdot 6}{1 \cdot 2 \cdot 3 \cdot 4 \cdot 5}=252.
\end{equation*}
Applying the last formula found for $u$ it will be $n=5$ and hence

\begin{equation*}
\renewcommand{\arraystretch}{1,2}
\setlength{\arraycolsep}{0.1mm}
\begin{array}{lrrl}
&\dfrac{1}{4n}&~=~& 0.0500000 \\[2mm]
&\dfrac{1}{32n^2}&~=~& 0.0012500  \\[2mm]
\cline{4-4}
& & &0.0512500 \\
\text{subtract} & \dfrac{1}{128n^3}&~=~&0.0000625 \\[2mm]
\cline{4-4} 
& & & 0.0511875 \\[2mm]              
\text{Therefore} \quad &1+\dfrac{1}{4n}+\text{etc.} &~=~&1.0511836 \\[2mm]    
\cline{4-4}
\text{the log. of this} \quad &  &~=~&0.0216784\\[2mm]
   & \log n   &~=~&0.6989700 \\[2mm]
   & \log \pi &~=~&0.4971498 \\[2mm]
   \cline{4-4} 
   & &  &1.2177982 \\[2mm]
   & \log \sqrt{n \pi(1+\text{etc.})}&~=~&0.6088991 \\[2mm]
   \text{from} & \log 2^{2n} &~ = ~& 3.0102999 \\[2mm]
   \cline{4-4}
    & \log u & ~=~ & 2.4014008 \\[2mm]
    \text{whence} & u&~=~& 252.
\end{array}
\end{equation*}

\subsection*{Example 2}

\textit{Investigate the ratio which   the middle term has to sum of all terms in the hundredth power of the binomial $1+1$, which sum is $2^{100}$.}\\[2mm]
For this, let us use the formula found first

\begin{equation*}
\log u= \log \frac{2^{2n}}{\sqrt{n \pi}}-\frac{3 \mathfrak{A}}{1 \cdot 2 \cdot 2 n}+\frac{15 \mathfrak{B}}{3 \cdot 4 \cdot 2^3 n^3}-\frac{63 \mathfrak{C}}{5 \cdot 6 \cdot 2^5 n^5}+\text{etc.},
\end{equation*}
in which, having put $2n=m$, that one has this power $(1+1)^m$ and having substituted the values for $\mathfrak{A}$, $\mathfrak{B}$, $\mathfrak{C}$, $\mathfrak{D}$ etc., it will be

\begin{equation*}
\log u= \log \frac{2^m}{\sqrt{\frac{1}{2}m \pi}}-\frac{1}{4m}+\frac{1}{24m^3}-\frac{1}{20m^5}+\frac{17}{112 m^7}-\frac{31}{36 m^9}+\frac{691}{88m^{11}}-\text{etc.};
\end{equation*}
since these logarithms are hyperbolic, multiply them by

\begin{equation*}
k=0.434294481903251,
\end{equation*}
that they are transformed into tabulated ones, and it will be

\begin{equation*}
\log u= \log \frac{2^m}{\sqrt{\frac{1}{2}m \pi}}-\frac{k}{4m}+\frac{k}{24m^3}-\frac{k}{20m^5}+\frac{17k}{112m^7}-\frac{31k}{36m^9}+\text{etc.},
\end{equation*}
whence, because the middle term is $u$, the ratio in question will be $2^m :u$ and hence

\begin{equation*}
\log \frac{2^m}{u}=\log \frac{1}{2}m \pi +\frac{k}{4m}-\frac{k}{24m^3}+\frac{k}{20m^5}-\frac{17 k}{122 m^7}+\frac{31k}{36m^9}-\frac{691 k}{88m^{11}}+\text{etc.}
\end{equation*}
Therefore, since, because of the exponent $m=100$,

\begin{equation*}
\frac{k}{m}=0.0043429448, \quad \frac{k}{m^3}=0.0000004343, \quad \frac{k}{m^5}=0.0000000000,
\end{equation*}
it will be

\begin{equation*}
\renewcommand{\arraystretch}{1,8}
\setlength{\arraycolsep}{0.1mm}
\begin{array}{lrrll}
 & \dfrac{k}{4m}&~=~&0.0010857362 \\[2mm]
 & \dfrac{k}{24 m^3}&~=~& 0.0000000181 \\[2mm]
 \cline{4-4}
 & & &0.0010857181 \\[2mm]
\text{Then} \\ 
 & \log \pi &~=~&0.4971498726 \\[2mm]
 &   \log \dfrac{1}{2}m &~=~& 1.6989700043 \\[2mm]
 \cline{4-4}
 &\log \dfrac{1}{2}m \pi &~=~& 2.1961198769 \\[2mm]
 \cline{4-4}
 &\log \sqrt{\dfrac{1}{2}m\pi}&~=~&1.0980599384 \\[2mm]
 &\dfrac{k}{4m}-\dfrac{k}{24m^3}+\text{etc.}&~=~& 0.0010857181 \\[2mm]
 \cline{4-4}
 & & & 1.0991456565 & = \log \dfrac{2^{100}}{u}.
\end{array}
\end{equation*}
Therefore, it will be $\frac{2^{100}}{u}=12.56451$ and hence in the expanded power $(1+1)^{100}$ the middle term will have the same ratio to the sum of all terms $2^{100}$ as $1$ to $12.56451$.

\paragraph*{§163}

Now, let the general term $z$ denote the exponential function $a^x$ so that this geometric series has to be summed

\begin{equation*}
s=a+a^2+a^3+a^4+\cdots +a^x;
\end{equation*}
since it is a geometric series, its sum is already known; for, it will be $s=\frac{(a^x-1)a}{a-1}$. But let us investigate this sum in the way explained here. Since $z=a^x$, it will be $\int zdx =\frac{a^x}{\log a}$; for, the differential of this is $a^xdx$; but then it will be

\begin{equation*}
\frac{dz}{dx}=a^x \log a, \quad \frac{ddz}{dx^2}=a^x(\log a)^2, \quad \frac{d^3 z}{dx^3}=a^x(\log a)^3 \quad \text{etc.},
\end{equation*}
whence it follows that

\begin{equation*}
s=a^x \left(\frac{1}{\log a}+\frac{1}{2}+\frac{\mathfrak{A}}{1 \cdot 2}\log a -\frac{\mathfrak{B}}{1 \cdot 2 \cdot 3 \cdot 4}(\log a)^3+\frac{\mathfrak{C}}{1 \cdot 2 \cdot 3 \cdots 6}(\log a)^5-\text{etc.}\right)+C.
\end{equation*}
To define the constant $C$ put $x=0$ and, because of $s=0$, it will be

\begin{equation*}
C=-\frac{1}{\log a}-\frac{1}{2}-\frac{\mathfrak{A}}{1 \cdot 2}\log a+\frac{\mathfrak{B}}{1 \cdot 2 \cdot 3 \cdot 4}(\log a)^3-\text{etc.}
\end{equation*}
and hence it will be

\begin{equation*}
s=(a^x-1) \left(\frac{1}{\log a}+\frac{1}{2}+\frac{\mathfrak{A}}{1 \cdot 2}\log a-\frac{\mathfrak{B}}{1 \cdot 2 \cdot 3 \cdot 4}(\log a)^3+\frac{\mathfrak{C}}{1 \cdot 2 \cdot 3 \cdots 6}(\log a)^5-\text{etc.} \right)
\end{equation*}
Therefore, because the sum is $\frac{(a^x-1)a}{a-1}$, it will be

\begin{equation*}
\frac{a}{a-1}=\frac{1}{\log a}+\frac{1}{2}+\frac{\mathfrak{A}}{1 \cdot 2}\log a- \frac{\mathfrak{B}}{1 \cdot 2 \cdot 3 \cdot 4}(\log a)^3+\frac{\mathfrak{C}}{1 \cdot 2 \cdot 3 \cdots 6}(\log a)^5-\text{etc.},
\end{equation*}
where $\log a$ denotes the hyperbolic logarithm of $a$; therefore,

\begin{equation*}
\frac{(a+1)\log a}{2(a-1)}=1+\frac{\mathfrak{A}(\log a)^2}{1 \cdot 2}-\frac{\mathfrak{B}(\log a)^4}{1 \cdot 2 \cdot 3 \cdot 4}+\frac{\mathfrak{C}(\log a)^5}{1 \cdot 2 \cdot 3 \cdots 6}-\text{etc.}
\end{equation*}
and so one will be able to exhibit the sum of this series.

\paragraph*{§164}

Let the general term be $z=\sin ax$ and

\begin{equation*}
s=\sin a+ \sin 2a+\sin 3a +\cdots +\sin ax;
\end{equation*}
this series, since it is recurring, can also be summed; for, it will be

\begin{equation*}
s=\frac{\sin a +\sin ax-\sin (ax+a)}{1-2 \cos a +1}=\frac{\sin a+(1- \cos a)\sin ax -\sin a \cdot \cos ax}{2(1-\cos a)}.
\end{equation*}
It will be

\begin{equation*}
\int zdx = \int dx \sin ax = -\frac{1}{a}\cos ax
\end{equation*}
and

\begin{equation*}
\frac{dz}{dx}=a \cos ax, \quad \frac{d^3 z}{dx^3}=-a^3 \cos ax, \quad \frac{d^5 z}{dx^5}=a^5 \cos ax \quad \text{etc.}
\end{equation*}
Therefore,

\begin{equation*}
s=C-\frac{1}{a}\cos ax+\frac{1}{2}\sin ax +\frac{\mathfrak{A}a \cos ax}{1 \cdot 2}+\frac{\mathfrak{B}a^3 \cos ax}{1 \cdot 2 \cdot 3 \cdot 4}
\end{equation*}
\begin{equation*}
+\frac{\mathfrak{C}a^5 \cos ax}{1 \cdot 2 \cdot 3 \cdot 4 \cdot 5 \cdot 6}+\frac{\mathfrak{D}a^7 \cos ax}{1 \cdot 2 \cdots 8}+\text{etc.}
\end{equation*}
Put $x=0$ so that $s=0$, and it will be
\begin{equation*}
C=\frac{1}{a}-\frac{\mathfrak{A}a}{1 \cdot 2}-\frac{\mathfrak{B}a^3}{1 \cdot 2 \cdot 3 \cdot 4}-\frac{\mathfrak{C}a^5}{1 \cdot 2 \cdots 6}-\text{etc.},
\end{equation*}
therefore

\begin{equation*}
s=\frac{1}{2}\sin ax+(1-\cos sx) \left(\frac{1}{a}-\frac{\mathfrak{A}a}{1 \cdot 2}-\frac{\mathfrak{B}a^3}{1 \cdot 2 \cdot 3 \cdot 4}-\frac{\mathfrak{C}a^5}{1 \cdot 2 \cdots 6}-\text{etc.}\right).
\end{equation*}
But since

\begin{equation*}
s=\frac{1}{2}\sin ax +\frac{(1- \cos ax)\sin a}{2(1- \cos a)},
\end{equation*}
it will become

\begin{equation*}
\frac{\sin a}{2(1-\cos a)}=\frac{1}{2}\cot \frac{1}{2}a=\frac{1}{a}-\frac{\mathfrak{A}a}{1 \cdot 2}-\frac{\mathfrak{B}a^3}{1 \cdot 2 \cdot 3 \cdot 4}-\frac{\mathfrak{C}a^5}{1 \cdot 2 \cdots 6}-\text{etc.},
\end{equation*}
which same series we already found above (§ 127).

\paragraph*{§165}

Now let $z=\cos ax $ and the series to be summed

\begin{equation*}
s= \cos a+ \cos 2a +\cos 3a +\cdots + \cos ax;
\end{equation*}
the sum of this series, since it is recurring, will be

\begin{equation*}
s=\frac{\cos a -1 +\cos ax -\cos (ax+a)}{1-2 \cos a +1}=-\frac{1}{2}+\frac{1}{2}\cos ax +\frac{1}{2}\cot \frac{1}{2}a \cdot \sin ax.
\end{equation*}
But on the other hand to express the sum by means of our method note that it will be

\begin{equation*}
\int zdx= \int dx \cos ax =\frac{1}{a}\sin ax
\end{equation*}
and

\begin{equation*}
\frac{dz}{dx}=-a\sin ax, \quad \frac{d^3z }{dx^3}=a^3 \sin ax, \quad \frac{d^5 z}{dx^5}=-a^5 \sin ax \quad \text{etc.}
\end{equation*}
Therefore, 

\begin{equation*}
s=C+\frac{1}{a}\sin ax +\frac{1}{2}\cos ax -\frac{\mathfrak{A}a \sin ax}{1 \cdot 2}-\frac{\mathfrak{B}a^3 \sin ax}{1 \cdot 2 \cdot 3 \cdot 4}-\text{etc.}
\end{equation*}
Let $x=0$, it will be $s=0$ and $C=-\frac{1}{2}$ and hence it will be

\begin{equation*}
s=-\frac{1}{2}+\frac{1}{2}\cos ax +\frac{1}{a}\sin ax -\frac{\mathfrak{A}a \sin ax}{1 \cdot 2}-\frac{\mathfrak{B}a^3 \sin ax}{1 \cdot 2 \cdot 3 \cdot 4}-\text{etc.}
\end{equation*}
Therefore, since

\begin{equation*}
s=-\frac{1}{2}+\frac{1}{2}\cos ax + \frac{1}{2}\cot \frac{1}{2}a \cdot \sin ax,
\end{equation*}
it will be, as we  just found [§ 164],

\begin{equation*}
\frac{1}{2}\cot \frac{1}{2}a=\frac{1}{a}-\frac{\mathfrak{A}a}{1 \cdot 2}-\frac{\mathfrak{B}a^3}{1 \cdot 2 \cdot 3 \cdot 4}-\frac{\mathfrak{C}a^5}{1 \cdot 2 \cdot 3 \cdot 4 \cdot 5 \cdot 6}-\text{etc.}
\end{equation*}

\paragraph*{§166}

Since we found above [§ 92], if $a$ denotes any arc, that 

\begin{equation*}
\frac{\pi}{2}= \frac{a}{2}+\sin a+ \frac{1}{2}\sin 2a +\frac{1}{3}\sin 3a +\frac{1}{4}\sin 4a +\text{etc.},
\end{equation*}
let us consider this series and let $z=\frac{1}{x}\sin ax $ that 

\begin{equation*}
s= \sin a +\frac{1}{2}\sin 2a +\frac{1}{3} \sin 3a + \cdots +\frac{1}{x}\sin ax.
\end{equation*}
But in this case $\int zdx =\int \frac{dx}{x}\sin ax$, which integral cannot be exhibited. Therefore, it will be

\begin{equation*}
\frac{dz}{dx}=\frac{a}{x}\cos ax -\frac{1}{xx}\sin ax, \quad \frac{ddz}{dx^2}=-\frac{a^2}{x}\sin ax - \frac{2a}{xx}\cos ax +\frac{2}{x^3}\sin ax,
\end{equation*}
\begin{equation*}
\frac{d^3 z}{dx^3}=-\frac{a^3}{x}\cos ax +\frac{3a^2}{x^2}\sin ax+
\frac{6a}{x^3}\cos ax -\frac{6}{x^4}\sin ax,
\end{equation*}
\begin{equation*}
\frac{d^4z}{dx^4}=\frac{a^4}{x}\sin ax +\frac{4a^3}{xx}\cos ax -\frac{12a^2}{x^3}\sin ax -\frac{24 a}{x^4}\cos ax +\frac{24}{x^5}\sin ax.
\end{equation*}
Therefore, since neither the integral formula $\int zdx$ can be exhibited nor  it is possible to express these differentials conveniently, we are not able to define the sum of this series by means of this expression, such that anything could be concluded from this. The same inconvenience occurs in many other series, if the general term is not  simple enough that its differentials can be expressed in general. But in the following chapter we will find other general expressions for the sums of the series whose general terms are either too composite or cannot be given at all; these can be applied successfully. But the insufficiency of the method treated here is especially revealed, if the signs of the terms of the propounded series alternate; for, then,  even though the general terms are simple, the summatory terms can nevertheless not be expressed conveniently using this method.


\title{A further Generalization of the summation method treated in chapter V
\footnote{
Original title: "'Methodus summandi superior ulterius promota"', first published as part of the book \textit{Institutiones calculi differentialis cum eius usu in analysi finitorum ac doctrina serierum}, 1755, reprinted in \textit{Opera Omnia}: Series 1, Volume 10, pp. 368 - 395, Eneström-Number E212, translated by: Alexander Aycock for the "`Euler-Kreis Mainz"'}}
\author{Leonhard Euler}
\date{ }
\maketitle

\paragraph*{§167}

To cure the defect of the summation method  we presented before,  in this chapter we will consider series, whose general terms are more complex. Since the expression found before in the case of geometric progressions, even though by  means of other methods they can be summed very easily, does not yield the true sum in terms of a finite formula, series  whose terms are products of terms of a geometric series and an arbitrary other series will be contemplated here at first. Therefore, let this series be propounded

\begin{alignat*}{9}
&  &&1 \quad && 2 \quad && 3 \quad && 4 \quad && \quad && x \\
&s~=~&&ap~+~ &&bp^2~+~&&cp^3~+~&&dp^4~+~ && \cdots \cdots ~+~&&yp^x,
\end{alignat*}
which is conflated of the geometric series $p$, $p^2$, $p^3$ etc. and  another arbitrary series $a+b+c+d+\text{etc.}$, whose general term or the one corresponding to index $x$ is $=y$, and let us investigate the general term for the value of its sum $s=S.yp^x$.

\paragraph*{§168}

Let us perform the calculation the same way as we did above, and let $v$ be the term preceding $y$ in the series $a+b+c+d+\text{etc.}$ and $A$ the one preceding $a$ or the one corresponding to the index $0$, then $vp^{x-1}$ will be the general term of this series

\begin{alignat*}{9}
& 1 \quad &&2 \quad && 3 \quad && 4 \quad && \quad && x \\
&A~+~&&ap~+~ &&bp^2~+~&&cp^3~+~&& \cdots \cdots ~+~&&yp^x;
\end{alignat*}
if its sum is indicated by $S.vp^{x-1}$, it will be

\begin{equation*}
S.vp^{x-1}=\frac{1}{p}S. vp^x=S. yp^x-yp^x+A.
\end{equation*}
But since

\begin{equation*}
v=y-\frac{dy}{dx}+\frac{ddy}{2dx^2}-\frac{d^3 y}{6 dx^3}+\frac{d^4 y}{24dx^4}-\frac{d^5 y}{120 dx^5}+\text{etc.},
\end{equation*}
it will be

\begin{equation*}
S.yp^x-yp^x+A= \frac{1}{p}S.yp^x-\frac{1}{p}S. \frac{dy}{dx}p^x+\frac{1}{2p}S. \frac{ddy}{dx^2}p^x
\end{equation*}
\begin{equation*}
-\frac{1}{6p}S. \frac{d^3 y}{dx^3}p^x+\frac{1}{24p}S. \frac{d^4 y}{dx^4}p^x -\text{etc.}
\end{equation*}
From this 

\begin{equation*}
S.yp^x=\frac{1}{p-1}\left(yp^{x+1}-Ap-S. \frac{dy}{dx}p^x+S.\frac{ddy}{2dx^2}p^x -S. \frac{d^3y}{6dx^3}+\text{etc.}\right)
\end{equation*}
Therefore, if one has the summatory terms of the series whose general terms are $\frac{dy}{dx}p^x$, $\frac{ddy}{dx^2}p^x$, $\frac{d^3y}{dx^3}p^x$ etc., from them one will be able to define the summatory term $S.yp^x$.

\paragraph*{§169}

Therefore, one will be able to find the sums of the series whose general terms are contained in the form $x^np^x$. For, let $y=x^n$; it will be $A=0$, if  not $n=0$, in which case  $A=1$, and since 

\begin{equation*}
\frac{dy}{dx}=nx^{n-1}, \quad \frac{ddy}{2dx^2}=\frac{n(n-1)}{1 \cdot 2}x^{n-2}, \quad \frac{d^3 y}{6dx^3}=\frac{n(n-1)(n-2)}{1 \cdot 2 \cdot 3}x^{n-3} \quad \text{etc.},
\end{equation*}
it will be

\begin{equation*}
S.x^n p^x= \frac{1}{p-1} \left\lbrace \begin{aligned}
& \qquad \qquad x^np^{x+1}-Ap-n S.x^{n-1}p^x+\frac{n(n-1)}{1 \cdot 2}S.x^{n-2}p^x \\
&  -\frac{n(n-1)(n-2)}{1 \cdot 2 \cdot 3}S.x^{n-3}p^x + \frac{n(n-1)(n-2)(n-3)}{1 \cdot 2 \cdot 3 \cdot 4}S. x^{n-4}p^x- \text{etc.}
\end{aligned} 
\right\rbrace
\end{equation*}
From this form, successively substituting the numbers $0$, $1$, $2$, $3$ etc. for $n$, one will obtain the following summations; first, if $n=0$,  $A=1$, but in the remaining cases it will be $A=0$.

\begin{equation*}
S.x^0 p^x= S. p^x=\frac{1}{p-1}(p^{x+1}-p)=\frac{p^{x+1}-p}{p-1}=\frac{p(p^x-1)}{p-1},
\end{equation*}
which is the known sum of the geometric progression; 

\begin{equation*}
S.xp^x= \frac{1}{p-1}(xp^{x+1}-S.p^x)=\frac{x p^{x+1}}{p-1}-\frac{p^{x+1}-p}{(p-1)^2}
\end{equation*}
or

\begin{equation*}
S.xp^x =\frac{pxp^x}{p-1}-\frac{p(p^x-1)}{(p-1)^2};
\end{equation*}
\begin{equation*}
S.x^2 p^x =\frac{1}{p-1}(x^2 p^{x+1}-2S. xp^x+S.p^x)
\end{equation*}
or

\begin{equation*}
S.x^2p^x=\frac{x^2 p^{x+1}}{p-1}-\frac{2xp^{x+1}}{(p-1)^2}+\frac{p(p+1)(p^x-1)}{(p-1)^3}.
\end{equation*}
Further,

\begin{equation*}
S.x^3 p^x =\frac{1}{p-1}(x^3 p^{x+1}-3S. x^2 p^x+3S.xp^x-S.p^x)
\end{equation*}
or

\begin{equation*}
S.x^3p^x = \frac{x^3 p^{x+1}}{p-1}-\frac{3x^2 p^{x+1}}{(p-1)^2}+\frac{3(p+1)xp^{x+1}}{(p-1)^3}-\frac{p(pp+4p+1)(p^x-1)}{(p-1)^4}
\end{equation*}
and proceeding this way, one will be able to define the sums of the higher powers $x^4p^x$, $x^5 p^x$, $x^6 p^x$ etc.; but this is achieved  more conveniently by means of the general expression we will now investigate.

\paragraph*{§170}

Since we found that 

\begin{equation*}
S.yp^x = \frac{1}{p-1} \left(yp^{x+1}-Ap-S. \frac{dy}{dx}p^x +S. \frac{ddy}{2 dx^2}p^x-S. \frac{d^3 y}{6 dx^3 }p^x +\text{etc.}\right),
\end{equation*}
where $A$ is a constant of such a kind that the sum becomes $=0$, if one puts $x=0$ (for, in this case  $y=A$ and $yp^{x+1}=Ap$), we will be able to omit this constant, if we only always remember that to a certain sum always a constant of such a kind is to be added that having put $x=0$ the sum vanishes or that in another case the correct sum is obtained using our summation formula. Therefore,  let us write $z$ instead of $y$ and it will be

\begin{equation*}
S. p^x z= \frac{p^{x+1}z}{p-1}-\frac{1}{p-1}S. p^x \frac{dz}{dx}+\frac{1}{2(p-1)}S.p^x \frac{ddz}{dx^2}-\frac{1}{6(p-1)}S.p^x \frac{d^3 z}{dx^3}
\end{equation*}
\begin{equation*}
+\frac{1}{24(p-1)}S.p^x \frac{d^4 z}{dx^4}-\frac{1}{120(p-1)}S.p^x \frac{d^5 z}{dx^5}+\text{etc.}
\end{equation*}
Furthermore, let us successively write $\frac{dz}{dx}$, $\frac{ddz}{dx^2}$, $\frac{d^3 z}{dx^3}$ etc. instead of $y$ and it will be

\begin{alignat*}{19}
&S.\frac{p^x dz}{dx}&&= \frac{p^{x+1}}{p-1}&&\cdot \frac{dz}{dx}&&-\frac{1}{p-1}S. &&\frac{p^x ddz}{dx^2}&&+\frac{1}{2(p-1)}S. &&\frac{p^x d^3 z}{dx^3}&&-\text{etc.} \\[3mm]
&S.\frac{p^x ddz}{dx^2}&&= \frac{p^{x+1}}{p-1}&&\cdot \frac{ddz}{dx^2}&&-\frac{1}{p-1}S. &&\frac{p^x d^3z}{dx^3}&&+\frac{1}{2(p-1)}S. &&\frac{p^x d^4 z}{dx^4}&&-\text{etc.} \\[3mm]
&S.\frac{p^x d^3z}{dx^3}&&= \frac{p^{x+1}}{p-1}&&\cdot \frac{d^3z}{dx^3}&&-\frac{1}{p-1}S. &&\frac{p^x d^4z}{dx^4}&&+\frac{1}{2(p-1)}S. &&\frac{p^x d^5 z}{dx^5}&&-\text{etc.} \\
& && && && &&\text{etc.}
\end{alignat*}
Therefore, if these values are successively substituted, $S.p^x z$ will be expressed by a form of this kind

\begin{equation*}
S.p^x z= \frac{p^{x+1}z}{p-1}-\frac{\alpha p^{x+1}}{p-1}\cdot \frac{dz}{dx}+\frac{\beta p^{x+1}}{p-1}\cdot \frac{ddz}{dx^2}-\frac{\gamma p^{x+1}}{p-1}\cdot \frac{d^3 z}{dx^3}
\end{equation*}
\begin{equation*}
+\frac{\delta p^{x+1}}{p-1}\cdot \frac{d^4 z}{dx^4}-\frac{\varepsilon p^{x+1}}{p-1}\cdot \frac{d^5 z}{dx^5}+\text{etc.}
\end{equation*}

\paragraph*{§171}

To define the values of the letters $\alpha$, $\beta$, $\gamma$, $\delta$, $\varepsilon$ etc. for each term substitute the series found before, i.e.

\begin{alignat*}{9}
& \frac{~~~~p^{x+1}z~~~}{p-1}&&=S.p^xz &&+\frac{1}{p-1}S.\frac{p^x dz}{dx}&&-\frac{1}{2(p-1)}S.\frac{p^x ddz}{dx^2}&&+\frac{1}{6(p-1)}S.\frac{p^x d^3z}{dx^3}-\text{etc.} \\[3mm]
&\frac{~~p^{x+1}dz~~}{(p-1)dx}&&=S. \frac{p^x dz}{dx}&&+\frac{1}{p-1}S. \frac{p^x ddz}{dx^2}&&-\frac{1}{2(p-1)}S.\frac{p^x d^3 z}{dx^3}&&+\text{etc.} \\[3mm]
&\frac{p^{x+1}ddz}{(p-1)dx^2}&&=S. \frac{p^x ddz}{dx^2}&&+\frac{1}{p-1}S. \frac{p^x d^3 z}{dx^3}&&-\text{etc.} \\[3mm]
&\frac{p^{x+1}d^3 z}{(p-1)dx^3}&&=S. \frac{p^x d^3 z}{dx^3}&&+\text{etc.}
\end{alignat*}
Therefore, we will have

\begin{equation*}
S.p^x z=S.p^x z
\end{equation*}
\begin{alignat*}{9}
&+\frac{1}{p-1}S. \frac{p^x dz}{dx}&&-\frac{1}{2(p-1)}S. \frac{p^x ddz}{dx^2}&&+\frac{1}{6(p-1)}S. \frac{p^x d^3 z}{dx^3}&&-\frac{1}{24(p-1)}S.\frac{p^x d^4 z }{dx^4}+\text{etc.} \\[3mm]
&-\alpha &&-\frac{\alpha}{p-1}&&+\frac{\alpha}{2(p-1)}&&-\frac{\alpha}{6(p-1)} \\[3mm]
& &&+\beta &&+\frac{\beta}{p-1}&&-\frac{\beta }{2(p-1)} \\[3mm]
& &&       &&-\gamma &&-\frac{\gamma}{p-1} \\[3mm]
& &&       &&       &&+\delta
\end{alignat*}
whence the following values of the coefficients $\alpha$, $\beta$, $\gamma$, $\delta$ etc. will be obtained

\begin{equation*}
\alpha =\frac{1}{p-1}, \quad \beta =\frac{1}{p-1}\left(\alpha +\frac{1}{2}\right), \quad \gamma = \frac{1}{p-1}\left(\beta +\frac{\alpha }{2}+\frac{1}{6}\right),
\end{equation*}
\begin{equation*}
\delta =\frac{1}{p-1}\left(\gamma +\frac{\beta}{2}+\frac{\alpha}{6}+\frac{1}{24}\right), \quad \varepsilon= \frac{1}{p-1}\left(\delta +\frac{\gamma}{2}+\frac{\beta }{6}+\frac{\alpha}{24}+\frac{1}{120}\right)
\end{equation*}
\begin{equation*}
\text{etc.}
\end{equation*}

\paragraph*{§172}

For the sake of brevity, let $\frac{1}{p-1}=q$; it will be

\begin{alignat*}{9}
&\alpha &&=q \\[3mm]
&\beta && = \alpha q+\frac{1}{2}q=qq+\frac{1}{2}q \\[3mm]
&\gamma && = \beta +\frac{1}{2}\alpha q +\frac{1}{6}q= q^3 +qq +\frac{1}{6}q \\[3mm]
&\delta &&= \gamma q +\frac{1}{2}\beta q+\frac{1}{6}\alpha q +\frac{1}{24}q=q^4+\frac{3}{2}q^3+\frac{7}{12}q^2+\frac{1}{24}q \\[3mm]
&\varepsilon &&= \delta q+ \frac{1}{2}\gamma q +\frac{1}{6}\beta q + \frac{1}{24}\alpha q +\frac{1}{120}q=q^5+2q^4+\frac{5}{4}q^3+\frac{1}{4}q^2+\frac{1}{120}q\\[3mm]
&\zeta &&= q^6+\frac{5}{2}q^5+\frac{13}{6}q^4+\frac{3}{4}q^3+\frac{31}{360}q^2+\frac{1}{720}q\\[3mm]
& && \text{etc.}
\end{alignat*}
or express them this way

\begin{alignat*}{9}
&\alpha &&=\frac{q}{1}\\[3mm]
&\beta &&=\frac{2qq+q}{1 \cdot 2}\\[3mm]
&\gamma &&= \frac{6q^3 +6q^2+q}{1 \cdot 2 \cdot 3} \\[3mm]
&\delta &&= \frac{24 q^4+36 q^3+14q^2+q}{1 \cdot 2 \cdot 3 \cdot 4} \\[3mm]
&\varepsilon &&= \frac{120 q^5+240 q^4+150q^3 +30 q^2+q}{1 \cdot 2 \cdot 3 \cdot 4 \cdot 5} \\[3mm]
&\zeta &&=\frac{720q^6+1800q^5+1560q^4+540q^3+62q^2+q}{1 \cdot 2 \cdot 3 \cdot 4 \cdot 5 \cdot 6}\\[3mm]
&\eta &&= \frac{5040 q^7+15120q^6+16800q^5+8400q^4+1806q^3+126q^2+q}{1 \cdot 2 \cdot 3 \cdot 4 \cdot 5 \cdot 6 \cdot 7} \\
& &&\text{etc.},
\end{alignat*}
where the coefficient $16800$  results, if the sum of the two above ones, $1560+1800$, is multiplied by the exponent of $q$, which is $5$ here.

\paragraph*{§173}

Let us substitute the value $\frac{1}{p-1}$ for $q$ again:

\begin{alignat*}{9}
&\alpha &&= \frac{1}{1(p-1)}\\[3mm]
& \beta &&=\frac{p+1}{1 \cdot 2(p-1)^2} \\[3mm]
&\gamma &&=\frac{pp+4p+1}{1 \cdot 2 \cdot 3 (p-1)^3}\\[3mm]
&\delta && =\frac{p^3+11p^2 +11p+1}{1 \cdot 2 \cdot 3 \cdot 4(p-1)^4}\\[3mm]
&\varepsilon &&= \frac{p^4 +26p^3+66p^2 +26p^2+1}{1 \cdot 2 \cdot 3 \cdot 4 \cdot 5(p-1)^5}\\[3mm]
&\zeta &&= \frac{p^5 +57p^4+302p^3+302p^2 +57p+1}{1 \cdot 2 \cdot 3 \cdot 4 \cdot 5 \cdot 6 (p-1)^6}\\[3mm]
&\eta &&=\frac{p^6+120 p^5+1191p^4+2416p^3+1191p^2+120p+1}{1 \cdot 2 \cdot 3 \cdot 4 \cdot 5 \cdot 6 \cdot 7 (p-1)^7} \\
& &&\text{etc.}
\end{alignat*} 
The structure of these formulas is that, if any term is set

\begin{equation*}
=\frac{p^{n-2}+Ap^{n-3}+Bp^{n-4}+Cp^{n-5}+Dp^{n-6}+\text{etc.}}{1 \cdot 2 \cdot 3 \cdots (n-1)(p-1)^{n-1}},
\end{equation*}
it will be

\begin{alignat*}{9}
&A &&=2^{n-1}-n \\[3mm]
&B &&=3^{n-1}-n \cdot 2^{n-1}+\frac{n(n-1)}{1 \cdot 2}\\[3mm]
&C &&=4^{n-1}-n \cdot 3^{n-1}+\frac{n(n-1)}{1 \cdot 2}2^{n-1}-\frac{n(n-1)(n-2)}{1 \cdot 2 \cdot 3}\\[3mm]
&D &&=5^{n-1}-n \cdot 4^{n-1}+\frac{n(n-1)}{1 \cdot 2}3^{n-1}-
\frac{n(n-1)(n-2)}{1 \cdot 2 \cdot 3}2^{n-1}+\frac{n(n-1)(n-2)(n-3)}{1 \cdot 2 \cdot 3 \cdot 4} \\[3mm]
&  &&\text{etc.},
\end{alignat*}
whence these coefficients $\alpha$, $\beta$, $\gamma$, $\delta$ etc. can be continued arbitrarily far.

\paragraph*{§174}

But if  on the other hand we consider the law, according to which these coefficients depend on each other, it will easily become clear that they constitute a recurring series and result, if this fraction is expanded

\begin{equation*}
\frac{1}{1-\frac{u}{p-1}-\frac{u^2}{2(p-1)}-\frac{u^3}{6(p-1)}-\frac{u^4}{24(p-1)}-\text{etc.}};
\end{equation*}
for, this series will result

\begin{equation*}
1+\alpha u+\beta u^2 +\gamma u^3 +\delta u^4 +\varepsilon u^5 +\zeta u^6 +\text{etc.};
\end{equation*}
put that fraction $=V$, and since

\begin{equation*}
V=\frac{p-1}{p-1-u-\frac{u^2}{2}-\frac{u^3}{6}-\frac{u^4}{24}-\text{etc.}},
\end{equation*}
it will be

\begin{equation*}
V=\frac{p-1}{p-e^u},
\end{equation*}
where $e$ is the number whose hyperbolic logarithm is $=1$. And if the value of $V$ is expressed by means of a power series in $u$, this equation will result

\begin{equation*}
V=1+\alpha u+\beta u^2+ \gamma u^3+ \delta u^4+\varepsilon u^5+\zeta u^6+\text{etc.},
\end{equation*}
whose coefficients $\alpha$, $\beta$, $\gamma$, $\delta$ etc. will be those, we need in the present task. Therefore, having found those, it will be

\begin{equation*}
S.p^x z= \frac{p^{x+1}}{p-1}\left(z-\frac{\alpha dz}{dx}+\frac{\beta ddz}{dx^2}-\frac{\gamma d^3 z}{dx^3}+\frac{\delta d^4 z}{dx^4}-\text{etc.}\right)\pm \text{Const.},
\end{equation*}
which expression therefore is the summatory of this series

\begin{equation*}
ap+bp^2+cp^3+ \cdots +p^x z,
\end{equation*}
whose general term is $=p^xz$.

\paragraph*{§175}

Since we found that $V=\frac{p-1}{p-e^u}$, it will be

\begin{equation*}
e^u=\frac{pV-p+1}{V}
\end{equation*}
and, by taking logarithms, it will be

\begin{equation*}
u=\log (pV-p+1)- \log V
\end{equation*}
and hence by differentiating

\begin{equation*}
du= \frac{(p-1)dV}{pV^2-(p-1)V},
\end{equation*}
whence it will be

\begin{equation*}
pV^2 = (p-1)V+\frac{(p-1)dV}{du}.
\end{equation*}
Therefore, since 

\begin{equation*}
V=1+\alpha u +\beta u^2 + \gamma u^3+ \delta u^4 +\varepsilon u^5+\text{etc.},
\end{equation*}
it will be

\begin{equation*}
\renewcommand{\arraystretch}{1,4}
\setlength{\arraycolsep}{0.0mm}
\begin{array}{llllllllllllllll}
p V^2=p~+~2\alpha u &~+~2 \beta p u^2 &~+~2 \gamma pu^3 &~+~2 \delta p u^4 &~+~2 \varepsilon p u^5 &~+~\text{etc.}\\[2mm]
& ~+~\alpha^2 pu^2 &~+~2 \alpha \beta p u^3 &~+~ 2 \alpha \gamma pu^4 & ~+~2 \alpha \delta p u^5 &~+~\text{etc.} \\[2mm]
&   & &~+~\beta \beta pu^4 &~+~2 \beta \gamma p u^5 &~+~\text{etc.}\\
\hline
\end{array}
\end{equation*}
\begin{equation*}
(p-1)V =(p-1)+\alpha (p-1)u+\beta (p-1)u^2 +\gamma (p-1) u^3
\end{equation*}
\begin{equation*}
+\delta (p-1)u^4+\varepsilon (p-1)u^5+\text{etc.}
\end{equation*}
\begin{equation*}
\frac{(p-1)dV}{du}=(p-1)\alpha +2(p-1)\beta u +2(p-1) \gamma u^2+4(p-1) \delta u^3
\end{equation*}
\begin{equation*}
+5(p-1)\varepsilon u^4+6(p-1) \zeta u^5+\text{etc.},
\end{equation*}
equating which expressions one will find

\begin{alignat*}{9}
&1(p-1) \alpha &&= &&1 \\
&2(p-1) \beta  &&=\alpha &&(p+1) \\
&3(p-1) \gamma &&= \beta &&(p+1)&&+\alpha^2 p \\
&4(p-1) \delta &&= \gamma &&(p+1)&&+2 \alpha \beta p \\
&5(p-1) \varepsilon &&= \delta &&(p+1)&&+2 \alpha \gamma p &&+\beta \beta p \\
&6(p-1) \zeta && =\varepsilon &&(p+1) &&+2 \alpha \delta p && +2 \beta \gamma p \\
&7(p-1) \eta &&=\zeta&&(p+1) &&+2 \alpha \varepsilon p &&+2 \beta \delta p &&+\gamma \gamma p \\
& &&  && \text{etc.}
\end{alignat*}
from which formulas, if given the number $p$, the values of the coefficients $\alpha$, $\beta$, $\gamma$, $\delta$ etc. can be determined more easily than from the rule we  found first.

\paragraph*{§176}

Before considering special cases of the value $p$, let us put $z=x^n$, so that this series has to be summed

\begin{equation*}
s=p+2^np^2+3^np^3+4^np^4+\cdots + x^n p^x,
\end{equation*}
and by  the expression found before  it will be

\begin{equation*}
s=p^x \left\lbrace \begin{aligned}
&\frac{p}{p-1}\cdot x^n -\frac{p}{(p-1)^2}n x^{n-1}+\frac{pp+p}{(p-1)^3}\cdot \frac{n(n-1)}{1 \cdot 2}x^{n-2} \\
& \qquad -\frac{p^3 +4p^2+p}{(p-1)^4}\cdot \frac{n(n-1)(n-2)}{1 \cdot 2 \cdot 3}x^{n-3}+\text{etc.} \end{aligned}
\right\rbrace
\end{equation*}
$\pm C$, which renders $s=0$, if one puts $x=0$.\\[2mm]
Therefore,  successively substituting the numbers $0$, $1$, $2$, $3$,  $4$ etc. for $n$, it will be

\begin{alignat*}{9}
&S.x^0 p^x &&=p^x \frac{p}{p-1}-\frac{p}{p-1} \\[3mm]
&S.x^1 p^x &&=p^x \left(\frac{px}{p-1}-\frac{p}{(p-1)^2}\right)+\frac{p}{(p-1)^2}\\[3mm]
&S.x^2 p^x &&=p^x \left(\frac{px^2}{p-1}-\frac{2px}{(p-1)^2}+\frac{p(p+1)}{(p-1)^3}\right)-\frac{p(p+1)}{(p-1)^3}\\[3mm]
&S.x^3 p^x &&=p^x \left(\frac{px^3}{p-1}-\frac{3 px^2}{(p-1)^2}+\frac{3p(p+1)x}{(p-1)^3}-\frac{p(p^2+4p+1)}{(p-1)^4} \right)+\frac{p(p^2+4p+1)}{(p-1)^4}\\[3mm]
&S.x^4 p^x &&=p^x \left(\frac{px^4}{p-1}-\frac{4px^3}{(p-1)^2}+\frac{6p(p+1)x^2}{(p-1)^3}-\frac{4p(p^2+4p+1)x}{(p-1)^4}+\frac{p(p^3+11p^2+11p+1)}{(p-1)^5}\right)\\[3mm]
& &&-\frac{p(p^3+11p^2+11p+1}{(p-1)^5} \\[3mm]
&S.x^5 p^x &&= \frac{p^{x+1}x^5}{p-1}-\frac{5p^{x+1}x^4}{(p-1)^2}+\frac{10(p+1)p^{x+1}x^3}{(p-1)^3}-\frac{10(p^2+4p+1)p^{x+1}x^2}{(p-1)^4} \\[3mm]
& &&+\frac{5(p^3+11p^2+11p+1)p^{x+1}x}{(p-1)^5}-\frac{(p^4 +26p^3 +66p^2 +26 p +1)(p^{x+1}-p)}{(p-1)^6}\\
\end{alignat*}
\begin{equation*}
S.x^6 p^x =\frac{p^{x+1}x^6}{p-1}-\frac{6p^{x+1}x^5}{(p-1)^2}+\frac{15(p+1)p^{x+1}x^4}{(p-1)^3}-\frac{20(p^2+4p+1)(p^{x+1}-p)}{(p-1)^4}
\end{equation*}
\begin{equation*}
+\frac{15(p^3+11p^2+11p+1)p^{x+1}x^2}{(p-1)^5}-\frac{6(p^4+26p^3+66p^2+26p+1)p^{x+1}x}{(p-1)^6}
\end{equation*}
\begin{equation*}
+\frac{(p^5+57p^4+302p^3+302p^2+57 p+1)(p^{x+1}-p)}{(p-1)}
\end{equation*}
\begin{equation*}
\text{etc.}
\end{equation*}

\paragraph*{§177}

Therefore, it is understood, if $z$ was a polynomial function of $x$, that the sum of the series whose general term is $p^x z$ can be exhibited, since by taking the differentials of $z$ one finally gets to vanishing ones. If this series is propounded

\begin{equation*}
p+3p^2+6p^3+10p^4+\cdots +\frac{xx+x}{2}p^x,
\end{equation*}
because of

\begin{equation*}
z=\frac{xx+x}{2} \quad \text{and} \quad \frac{dz}{dx}=x+\frac{1}{2} \quad \text{and} \quad \frac{ddz}{dx^2}=1,
\end{equation*}
the summatory term will be

\begin{equation*}
s=\frac{p^{x+1}}{p-1}\left(\frac{1}{2}xx+\frac{1}{2}x-\frac{2x+1}{2(p-1)}+\frac{p+1}{2(p-1)^2}\right)-\frac{p}{p-1}\left(\frac{p+1}{2(p-1)^2}-\frac{1}{2(p-1)}\right)
\end{equation*}
or

\begin{equation*}
s=p^{x+1} \left(\frac{xx}{2(p-1)}+\frac{(p-3)}{2(p-1)^2}+\frac{1}{(p-1)^3}\right)-\frac{p}{p-1}.
\end{equation*}
But if $z$ was not a polynomial function, then this expression of the summatory term runs to infinity. So, if $z=\frac{1}{x}$, that this series is to be summed

\begin{equation*}
s=p+\frac{1}{2}p^2+\frac{1}{3}p^3+\frac{1}{4}p^4+\cdots + \frac{1}{x}p^x, 
\end{equation*}
because of 

\begin{equation*}
\frac{dz}{dx}=-\frac{1}{xx}, \quad \frac{ddz}{dx^2}=\frac{2}{x^3}, \quad \frac{d^3 z}{dx^3}=-\frac{2 \cdot 3}{x^4}, \quad \frac{d^4 z}{dx^4}=\frac{2 \cdot 3 \cdot 4}{x^5} \quad \text{etc.},
\end{equation*}
this summatory term will result

\begin{equation*}
s=\frac{p^{x+1}}{p-1}\left(\frac{1}{x}+\frac{1}{(p-1)x^2}+\frac{p+1}{(p-1)^2 x^3}+\frac{pp+4p+1}{(p-1)^3x^4}+\frac{p^3+11p^2+11p+1}{(p-1)^4 x^5}+\text{etc.}\right)+C.
\end{equation*}
In this case, the constant $C$ can therefore not be defined from the case $x=0$; to define it, put $x=1$, and since  $s=p$, it will be

\begin{equation*}
C=p-\frac{pp}{p-1}\left(1+\frac{1}{p-1}+\frac{p+1}{(p-1)^2}+\frac{pp+4p+1}{(p-1)^3}+\text{etc.}\right).
\end{equation*}

\paragraph*{§178}

From these results it is perspicuous, if $p$ does not denote a specified number, that hardly anything of use to exhibit the sums of series approximately follows from this. But  first it is plain that  one cannot write $1$ for $p$, since all coefficients $\alpha$, $\beta$, $\gamma$, $\delta$ etc. would become infinitely large. Therefore, because the series we are considering now, goes over into that one we already contemplated before, if one puts $p=1$, it is surprising that this case cannot be found  from this general expression, although it is the  simplest case. Then, on the other hand it is also notable that in the case $p=1$ the summation requires the integral $\int zdx$, although in general the sum can be exhibited without a single integral. So it happens that, while all the coefficients $\alpha$, $\beta$, $\gamma$, $\delta$ etc. grow to infinity, at the same time that integral formula is enters the expression. And this case, in which $p=1$, is the only one to which this general expression found here cannot be applied. But even in this case the general formula is not to be considered to be wrong; for, even though each term becomes infinite, nevertheless all the infinities indeed cancel each other and a finite quantity remains which is equal to the one we found by means of the first method; we will elaborate on this in more detail below.

\paragraph*{§179}

Therefore, let $p=-1$ and the signs in the series to be summed will alternate

\begin{alignat*}{19}
& &&1 && &&2 && &&3 && &&4 && \qquad &&x \\
&~-~&&a &&~+~&&b &&~-~&&c &&~+~&& d&&- \cdots ~\pm~ && z,
\end{alignat*}
where $z$ will be positive, if $x$ was an even number, but negative, if $x$ is an odd number. Therefore, having put

\begin{equation*}
-a+b-c+d-\cdots \pm z =s,
\end{equation*}
it will be

\begin{equation*}
s=\pm \frac{ 1}{2}\left(z- \frac{\alpha dz}{dx}+\frac{\beta ddz}{dx^2}-\frac{\gamma d^3 z}{dx^3}+\frac{\delta d^4 z}{dx^4}-\text{etc.}\right)+C,
\end{equation*}
where the upper of the ambiguous signs hold, if $x$ is an even number, the other, if $x$ is an odd number. Therefore, by changing the sings, it will be

\begin{equation*}
a-b+c-d+e-f + \cdots \mp z = \mp \frac{1}{2}\left(z-\frac{\alpha dz}{dx}+\frac{\beta ddz}{dx^2}-\frac{\gamma d^3 z}{dx^3}+\frac{\delta d^4 z}{dx^4}-\text{etc.}\right)+C,
\end{equation*}
where the ambiguity of the signs follows the same rule.

\paragraph*{§180}

In this case the coefficients $\alpha$, $\beta$, $\gamma$, $\delta$, $\varepsilon$, $\zeta$ etc. can be found from the values given before by putting $p=-1$ everywhere. But they will be found more easily from the general formulas given in § 175, whence at the same time it is understood that each second of these coefficients vanishes. For, having put $p=-1$, these formulas will go over into

\begin{equation*}
-2 \alpha =0, \quad -4 \beta =0, \quad - 6 \gamma =0- \alpha^2, \quad -8 \delta =0 -2 \alpha \beta,
\end{equation*}
\begin{equation*}
-10 \varepsilon =0 - 2 \alpha \gamma - \beta \beta, \quad -12 \zeta =0 -2 \alpha \delta- 2\beta \gamma \quad \text{etc.}
\end{equation*}
whence, because $\beta =0$, it also will be $\delta =0$ and further $\zeta=0$, $\theta =0$ and the remaining letters will be determined  as follows

\begin{equation*}
\alpha =- \frac{1}{2}, \quad \gamma =\frac{\alpha^2}{6}, \quad \varepsilon = \frac{2 \alpha \gamma}{10}, \quad \eta =\frac{2 \alpha \varepsilon +\gamma \gamma}{14}, \quad \iota = \frac{2 \alpha \eta +2 \gamma \varepsilon}{18} \quad \text{etc.}
\end{equation*}

\paragraph*{§181}

That this calculation can be performed in a more convenient way, let us introduce new letters and let

\begin{equation*}
\alpha =- \frac{A}{1 \cdot 2}, \quad \gamma = \frac{B}{1 \cdot 2 \cdot 3 \cdot 4}, \quad \varepsilon =- \frac{C}{1 \cdot 2 \cdot 3 \cdot 4 \cdot 5 \cdot 6},
\end{equation*}
\begin{equation*}
\eta= \frac{D}{1 \cdot 2 \cdot 3 \cdots 8}, \quad \iota = - \frac{E}{1 \cdot 2 \cdot 3 \cdots 10}
\end{equation*}
\begin{equation*}
\text{etc.}
\end{equation*}
 And the sum exhibited before will be
 
 \begin{equation*}
 \mp \frac{1}{2}\left(z+\frac{Adz}{1 \cdot 2 dx}-\frac{Bd^3 z}{1 \cdot 2 \cdot 3 \cdot 4 dx^3}+\frac{C d^5 z}{1 \cdot 2 \cdots 6 dx^5}-\frac{D d^7 z}{1 \cdot 2 \cdots 8 dx^7}+\text{etc.}\right)+C.
 \end{equation*}
 But the coefficients will be defined from the following formulas
 
 \begin{alignat*}{9}
 &1A &&=1 \\[3mm]
 &3B &&= \frac{~4 \cdot 3~}{1 \cdot 2}&&\cdot \frac{AA}{2}\\[3mm]
 &5C &&= \frac{~6 \cdot 5~}{1 \cdot 2}&&\cdot AB \\[3mm]
 &7D &&= \frac{~8 \cdot 7~}{1 \cdot 2}&&\cdot AC &&+\frac{~~~8 \cdot 7 \cdot 6 \cdot 5~~}{1 \cdot 2 \cdot 3 \cdot 4}&&\cdot\frac{BB}{2}\\[3mm]
 &9E &&=\frac{~10 \cdot 9}{1 \cdot 2}&&\cdot AD &&+\frac{~~10 \cdot 9 \cdot 8 \cdot 7~~}{1 \cdot 2 \cdot 3 \cdot 4}&&\cdot BC \\[3mm]
 &11F &&=\frac{12 \cdot 11}{1 \cdot 2}&&\cdot AE &&+\frac{12 \cdot 11 \cdot 10 \cdot 9}{1 \cdot 2 \cdot 3 \cdot 4}&&\cdot BD + \frac{12 \cdot 11 \cdot 10 \cdot 9 \cdot 8 \cdot 7}{1 \cdot 2 \cdot 3 \cdot 4 \cdot 5 \cdot 6} \cdot \frac{CC}{2}\\
 & && &&\text{etc.},
 \end{alignat*}
 which can be represented more easily and are more accommodated to calculation  purposes  this way 

\begin{alignat*}{19}
&A&&=1 \\[3mm]
&B &&=2 &&\cdot \frac{AA}{2} \\[3mm]
&C &&=3 && \cdot AB \\[3mm]
&D &&=4 && \cdot AC &&+4 \cdot ~~\frac{6 \cdot 5}{3 \cdot 4}&& \cdot \frac{BB}{2} \\[3mm]
&E &&=5 && \cdot AC &&+5 \cdot ~~\frac{8 \cdot 7}{3 \cdot 4}&& \cdot BC \\[3mm]
&F &&=6 &&\cdot AD  &&+6 \cdot ~\frac{10 \cdot 9}{3 \cdot 4}&&\cdot BD &&+6 \cdot \frac{10 \cdot 9 \cdot 8 \cdot 7}{\cdot 3 \cdot 4 \cdot 5 \cdot 6} \cdot \frac{CC}{2} \\[3mm]
&G &&=7 &&\cdot AE  &&+6 \cdot \frac{12 \cdot 11}{3 \cdot 4}&& \cdot BE &&+7 \cdot \frac{12 \cdot 11 \cdot 10 \cdot 9}{\cdot 3 \cdot 4 \cdot 5 \cdot 6} \cdot CD \\[3mm]
& && && &&\text{etc.}
\end{alignat*}
Therefore, having done the calculation, one will find

\begin{alignat*}{9}
& A&&=1 \\
&B &&=1 \\
&C &&=3 \\
&D &&=17 \\
&E &&=155 &&=5 \cdot 31 \\
&F &&=2073 &&=691 \cdot 3 \\
&G &&=38227 &&=7 \cdot 5461 &&=7 \cdot \frac{127 \cdot 129}{3} \\
&H &&=929569 &&= 3617 \cdot 257 \\
&I &&= 28820619 &&=43867 \cdot 9 \cdot 73 \\
&  && &&\text{etc.}
\end{alignat*}

\paragraph*{§182}

If we consider these numbers with more attention, from the factors $691$, $3617$, $43867$ one can easily conclude that these numbers  have a connection to the Bernoulli numbers exhibited above [§ 122] and can hence be determined. Therefore,  anyone investigating this relation will immediately see that these numbers can be formed from the Bernoulli numbers $\mathfrak{A}$, $\mathfrak{B}$, $\mathfrak{C}$, $\mathfrak{D}$, $\mathfrak{E}$ etc. the following way:
\begin{alignat*}{9}
&A &&= 2 &&\cdot ~~1&&\cdot ~~ 3 &&\mathfrak{A}&&~=~2(2^2 &&-1)\mathfrak{A}\\ 
&B &&= 2 &&\cdot ~~3&&\cdot ~~ 5 &&\mathfrak{B}&&~=~2(2^4 &&-1)\mathfrak{B} \\  
&C &&= 2 &&\cdot ~~7&&\cdot ~~ 9 &&\mathfrak{C}&&~=~2(2^6 &&-1)\mathfrak{C} \\ 
&D &&= 2 &&\cdot ~15&&\cdot ~ 17 &&\mathfrak{D}&&~=~2(2^8 &&-1)\mathfrak{D} \\
&E &&= 2 &&\cdot ~31&&\cdot ~ 33 &&\mathfrak{E}&&~=~2(2^{10} &&-1)\mathfrak{E} \\
&F &&= 2 &&\cdot ~63&&\cdot ~ 67 &&\mathfrak{F}&&~=~2(2^{12} &&-1)\mathfrak{F} \\
&G &&= 2 &&\cdot 127&&\cdot 129 &&\mathfrak{G}&&~=~2(2^{14} &&-1)\mathfrak{G} \\
&H &&= 2 &&\cdot 255&&\cdot 257 &&\mathfrak{H}&&~=~2(2^{16} &&-1)\mathfrak{H} \\
&  &&    &&         && && \text{etc.}           
\end{alignat*}
Since the Bernoulli numbers are fractions, but our coefficients on the other hand are integers, it is plain that these factors always cancel the fractions and they will therefore be

\begin{alignat*}{9}
&A&&=1 \\
&B&&=1 \\
&C &&=3 \\
&D&&=17 \\
&E&&=5 \cdot 31 =155 \\
&F&&= 3 \cdot 691 = 2073 \\
&G&&=7 \cdot 43 \cdot 127 =38277 \\
&H &&=257 \cdot 3617 =929569 \\
&I &&= 9 \cdot 73 \cdot 43867 =28820619 \\
&K &&= 5 \cdot 31 \cdot 41 \cdot 174611 = 1109652905 \\
&L &&= 89 \cdot 683 \cdot 854513 = 51943281713 \\
&M &&= 3 \cdot 4097 \cdot 236364091 = 2905151042481 \\
&N &&= 2731 \cdot 8191 \cdot 8553103 =191329672483963 \\
&  && \qquad \qquad \text{etc.}
\end{alignat*}
Therefore, using these numbers, one will  vice versa be able to find the Bernoulli numbers.

\paragraph*{§183}

Therefore, by applying these Bernoulli numbers to the propounded series

\begin{alignat*}{9}
&1   && 2  && 3   && 4   && 5             && x \\
&a~-~&&b~+~&& c~-~&& d~+~&& e- \cdots \mp && z
\end{alignat*}
the sum will be

\begin{equation*}
\mp \left(\frac{1}{2}z+\frac{(2^2-1)\mathfrak{A}dz}{1 \cdot 2 dx}-\frac{(2^4-1)\mathfrak{B}d^3z}{1 \cdot 2 \cdot 3 \cdot 4 dx^3}+\frac{(2^6-1)\mathfrak{C}d^5 z}{1 \cdot 2 \cdots 6 dx^5}-\frac{(2^8-1)\mathfrak{D}d^7 z}{1 \cdot 2 \cdots 8 dx^7}+\text{etc.}\right)+\text{Const.}
\end{equation*}
But hence it is understood that these numbers do not enter this expression  accidentally; for, as the propounded series results, if from this one

\begin{equation*}
a+b+c+d+\cdots + z,
\end{equation*}
where all terms have the sign $+$, the sum of all second terms $b+d+f+\text{etc.}$ is subtracted twice, so also the found expression can be resolved into two parts of which the one is the sum of all terms affected with the sign $+$ which will be

\begin{equation*}
\int zdx +\frac{1}{2}z+\frac{\mathfrak{A}dz}{1 \cdot 2 dx}-\frac{\mathfrak{B}d^3 z}{1 \cdot 2 \cdot 3 \cdot 4 dx^3}+\frac{\mathfrak{C}d^5z}{1 \cdot 2 \cdots 6 dx^5}-\text{etc.}
\end{equation*}
But the sum of all second terms is found the same way we did it above [chap. V]. Since the last term  $z$ corresponds to the index $x$, the term corresponding to the index $x-2$ will be

\begin{equation*}
z-\frac{2dz}{dx}+\frac{2^2 ddz}{1 \cdot 2 dx^2}-\frac{2^3 d^3 z}{1 \cdot 2 \cdot 3 dx^3}+\frac{2^4 d^4 z}{1 \cdot 2 \cdot 3 \cdot 4 dx^4}-\text{etc. },
\end{equation*}
which form results from the other expressing the preceding term, if  one writes $\frac{x}{2}$ instead of $x$. Therefore, one will have the sum of all the second terms, if in the sum of all terms  one writes $\frac{x}{2}$ instead of $x$ everywhere; therefore, this sum will be

\begin{equation*}
\frac{1}{2}\int z dx +\frac{1}{2}z+\frac{2 \mathfrak{A}dz}{1 \cdot 2 dx}-\frac{2^3 \mathfrak{B}d^3 z}{1 \cdot 2 \cdot 3 \cdot 4 dx^3}+\frac{2^5 \mathfrak{C}d^5 z }{1 \cdot 2 \cdots 6 dx^5}-\text{etc.};
\end{equation*}
if its double is subtracted from the preceding while $x$ is an even number or if the preceding sum is subtracted from the double of this one, if $x$ is an odd number, the residue will show the sum of the series

\begin{alignat*}{9}
&1   && 2  && 3   && 4   && 5             && x \\
&a~-~&&b~+~&& c~-~&& d~+~&& e- \cdots \mp && z,
\end{alignat*}
which will therefore be

\begin{equation*}
\mp \left(\frac{1}{2}z+\frac{(2^2-1)\mathfrak{A}dz}{1 \cdot 2 dx}-\frac{(2^4-1)\mathfrak{B}d^3 z}{1 \cdot 2 \cdot 3 \cdot 4}+\text{etc.}\right)+C,
\end{equation*}
which is the same expression we just found.

\paragraph*{§184}

Take a power of $x$ for $z$, namely $x^n$, so that one finds the sum of the series

\begin{equation*}
1-2^n+3^n-4^n+\cdots -\pm x^n
\end{equation*}
Because of 

\begin{equation*}
\frac{dz}{1 dx}=\frac{n}{1}x^{n-1}, \quad \frac{d^3 z}{1\cdot 2 \cdot 3 dx^3}= \frac{n(n-1)(n-2)}{1 \cdot 2 \cdot 3}x^{n-3}\quad \text{etc.},
\end{equation*}
introducing the coefficients $A$, $B$, $C$, $D$, $E$ etc., it will be

\begin{equation*}
\mp \frac{1}{2}\left\lbrace
\begin{aligned}
&x^n +\frac{A}{2}nx^{n-1}-\frac{B}{4}\cdot \frac{n(n-1)(n-2)}{1 \cdot 2 \cdot 3}x^{n-3}+\frac{C}{6}\cdot \frac{n(n-1)(n-2)(n-3)(n-4)}{1 \cdot 2 \cdot 3 \cdot 4 \cdot 5}x^{n-5} \\
& \qquad \qquad \qquad -\frac{D}{8}\cdot \frac{n(n-1)\cdots (n-6)}{1 \cdot 2 \cdots 7}x^{n-7}+\text{etc.}+\text{Const.}
\end{aligned}
 \right\rbrace
\end{equation*}
where the upper sign holds, if $x$ is an even number, the lower on the other hand, if it is an odd number. But the constant has to defined in such a way that the sum vanishes, if $x=0$, in which case the above sign holds. Successively substituting the numbers $0$, $1$, $2$, $3$ etc. for $n$, the following sums will result

\begin{equation*}
\text{I.} \quad 1-1+1-1+\cdots \mp 1=\mp \frac{1}{2}(1)+\frac{1}{2};
\end{equation*}
if the number of terms was even, the sum will be $=0$, if odd, it will be $=+1$.

\begin{equation*}
\text{II.} \quad 1-2+3-4+\cdots \mp x =\mp \frac{1}{2}\left(x+\frac{1}{2}\right)+\frac{1}{4};
\end{equation*}
if the number of terms is even, the sum will be $=-\frac{1}{2}x$ and for an odd number of terms $=+\frac{1}{2}x+\frac{1}{2}$.

\begin{equation*}
\text{III.} \quad 1-2^2+3^2-4^2+\cdots \mp x^2=\mp\frac{1}{2}\left(x^2+x\right)
\end{equation*}
for an even number $=- \frac{1}{2}xx-\frac{1}{2}x$ and for an odd number $=+\frac{1}{2}xx+\frac{1}{2}x$.

\begin{equation*}
\text{IV.} \quad 1-2^3+3^3-4^3+ \cdots \mp x^3 = \mp \frac{1}{2}\left(x^3+\frac{3}{2}xx -\frac{1}{4}\right)-\frac{1}{8};
\end{equation*}
for even number of terms $=-\frac{1}{2}x^3-\frac{3}{4}x^2$ and for an odd number $=\frac{1}{2}x^3+\frac{3}{4}x^2-\frac{1}{4}$.

\begin{equation*}
\text{V.} \quad 1-2^4+3^4-4^4+\cdots \mp x^4 = \mp \frac{1}{2}(x^4-2x^3-x);
\end{equation*}
for an even number of terms $=-\frac{1}{2}x^4-x^3+\frac{1}{2}x^2$ and for an odd number $=\frac{1}{2}x^4 +x^3-\frac{1}{2}x$ etc.

\paragraph*{§185}

Therefore, it is clear that in the sums of even powers except for the case $n=0$ the constant to be added vanishes and in these cases the sum of an even number or an odd number of terms only differs with regard to the sign. Therefore, if $x$ was an infinite number, since it is neither even nor odd, this consideration is not valid and the ambiguous signs are to be rejected; therefore, it follows that the sum of series, if they are continued infinity, of this kind are expressed by means of the constant to be added alone.\\[2mm]
Therefore, it will be

\begin{alignat*}{9}
&1 &&-1 &&+1&&-1 &&+\text{etc. to infinity}&&=+\frac{1}{2}\\[3mm]
&1 &&-2 &&+3&&-4 &&+\text{etc.}            &&=+\frac{A}{4}=+\frac{(2^2-1)\mathfrak{A}}{2}\\[3mm]
&1 &&-2^2&&+3^2 &&-4^2 &&+\text{etc.}     && =0 \\[3mm]
&1 &&-2^3&&+3^3 &&-4^3 &&+\text{etc.}     && =-\frac{B}{8}=-\frac{(2^4-1)\mathfrak{B}}{4} \\[3mm]
&1 &&-2^4&&+3^4 &&-4^4 &&+\text{etc.}     && =0 \\[3mm]
&1 &&-2^5&&+3^5 &&-4^5 &&+\text{etc.}     && =+\frac{C}{12}=+\frac{(2^6-1)\mathfrak{C}}{6} \\[3mm]
&1 &&-2^6&&+3^6 &&-4^6 &&+\text{etc.}     && =0 \\[3mm]
&1 &&-2^7&&+3^7 &&-4^7 &&+\text{etc.}     && =-\frac{D}{16}=-\frac{(2^8-1)\mathfrak{D}}{8} \\
&  &&    &&     &&     &&\text{etc.}
\end{alignat*}
These same sums are found by means of the method to sum series in which the signs $+$ and $-$ alternate treated above [§ 8].

\paragraph*{§186}

If  one takes negative numbers for $n$, the expression for the sums runs to infinity. Let $n=-1$; the sum of the series will be

\begin{equation*}
1-\frac{1}{2}+\frac{1}{3}-\frac{1}{4}+\frac{1}{5}-\frac{1}{6}+\cdots \mp \frac{1}{x}
\end{equation*}
\begin{equation*}
=\mp \frac{1}{2}\left(\frac{1}{x}-\frac{A}{2x^2}+\frac{B}{4x^4}-\frac{C}{6x^6}+\frac{D}{8x^8}-\text{etc.}\right)+\text{Const.}
\end{equation*}
Because here the constant cannot be defined from the case $x=0$, it is to be defined from another case. Put $x=1$ and, because of the sum $=1$ and the lower sign holds, it will be

\begin{equation*}
\text{Const.}=1-\frac{1}{2}\left(\frac{1}{1}-\frac{A}{2}+\frac{B}{4}-\frac{C}{6}+\text{etc.}\right)
\end{equation*}
or

\begin{equation*}
\text{Const.}=\frac{1}{2}+\frac{A}{4}-\frac{B}{8}+\frac{C}{12}-\frac{D}{16}+\text{etc.}
\end{equation*}
Or put $x=2$; because of the sum $=\frac{1}{2}$, and since the upper sign holds, one will find

\begin{equation*}
\text{Const.}=\frac{1}{2}+\frac{1}{2}\left(\frac{1}{2}-\frac{A}{2 \cdot 2^2}+\frac{B}{4 \cdot 2^4}-\frac{C}{6 \cdot 2^6}+\text{etc.}\right)
\end{equation*}
or

\begin{equation*}
\text{Const.}= \frac{3}{4}-\frac{A}{4 \cdot 2^2}+\frac{B}{8 \cdot 2^4}-\frac{C}{12 \cdot 2^6}+\frac{D}{16 \cdot 2^8}-\text{etc.}
\end{equation*}
But if one puts $x=4$, it will be

\begin{equation*}
\text{Const.}=\frac{17}{24}-\frac{A}{4 \cdot 4^2}+\frac{B}{8 \cdot 4^4}-\frac{C}{12 \cdot 4^6}+\frac{D}{16 \cdot 4^8}-\text{etc.}
\end{equation*}
But no matter how the constant is defined, the same value will result, which value at the same time will indicate the sum of the series continued to infinity which is $=\log 2$.

\paragraph*{§187}

Additionally, from these new numbers $A$, $B$, $C$, $D$, $E$ etc. the sums of the series of the reciprocal powers in which only the odd numbers occur can be summed in a convenient way. For, if one puts

\begin{alignat*}{9}
&1 &&+\frac{1}{2^{2n}}+\frac{1}{3^{2n}}&&+\frac{1}{4^{2n}}+\frac{1}{5^{2n}}&&+\frac{1}{6^{2n}}&&+\text{etc.}&&=s
\intertext{it will be}
&  &&+\frac{1}{2^{2n}}                &&+\frac{1}{4^{2n}}    &&+\frac{1}{6^{2n}}&&+\text{etc.}&&=\frac{s}{2^{2n}},
\end{alignat*}
which subtracted from the other gives

\begin{equation*}
1+\frac{1}{3^{2n}}+\frac{1}{5^{2n}}+\frac{1}{7^{2n}}+\text{etc.}= \frac{(2^{2n}-1)s}{2^{2n}},
\end{equation*}
Because we exhibited the values of $s$ for each number $n$ already above (§ 125), it will be

\begin{alignat*}{9}
&1+&&+\frac{1}{3^2}&&+\frac{1}{5^2}&&+\frac{1}{7^2}&&+\text{etc.}&&= \quad \quad \frac{A}{1 \cdot 2} &&\cdot \frac{\pi^2}{4} \\[3mm]
&1+&&+\frac{1}{3^4}&&+\frac{1}{5^4}&&+\frac{1}{7^4}&&+\text{etc.}&&= \quad \frac{B}{1 \cdot 2 \cdot 3 \cdot 4} &&\cdot \frac{\pi^4}{4} \\[3mm]
&1+&&+\frac{1}{3^6}&&+\frac{1}{5^6}&&+\frac{1}{7^6}&&+\text{etc.}&&= ~\frac{C}{1 \cdot 2 \cdot 3 \cdots 6} &&\cdot \frac{\pi^6}{4} \\[3mm]
&1+&&+\frac{1}{3^8}&&+\frac{1}{5^8}&&+\frac{1}{7^8}&&+\text{etc.}&&= ~\frac{D}{1 \cdot 2 \cdot 3 \cdots 8} &&\cdot \frac{\pi^8}{4} \\[3mm]
&1+&&+\frac{1}{3^{10}}&&+\frac{1}{5^{10}}&&+\frac{1}{7^{10}}&&+\text{etc.}&&= \frac{E}{1 \cdot 2 \cdot 3 \cdots 10} &&\cdot \frac{\pi^2}{4} \\[3mm]
& && && && \text{etc.}
\end{alignat*}
But if all numbers enter the expression and the signs alternate, since it will be

\begin{equation*}
1-\frac{1}{2^{2n}}+\frac{1}{3^{2n}}-\frac{1}{4^{2n}}+\text{etc.}=\frac{(2^{2n}-1)s-s}{2^{2n}},
\end{equation*}
one will have

\begin{alignat*}{19}
&1 &&- \frac{1}{2^2}&&+\frac{1}{3^2}&&-\frac{1}{4^2}&&+\frac{1}{5^2}&&-\text{etc.}&&=~\frac{A- 2\mathfrak{A}}{1 \cdot 2}&&\cdot \frac{\pi^2}{4}&&=~\frac{(2-1)\mathfrak{A}}{1 \cdot 2}&&\cdot \pi^2 \\[3mm]
&1 &&- \frac{1}{2^4}&&+\frac{1}{3^4}&&-\frac{1}{4^4}&&+\frac{1}{5^4}&&-\text{etc.}&&=~\frac{B- 2\mathfrak{B}}{1 \cdot 2 \cdot 3 \cdot 4}&&\cdot \frac{\pi^4}{4}&&=\frac{(2^3-1)\mathfrak{B}}{1 \cdot 2 \cdot 3 \cdot 4}&&\cdot \pi^4 \\[3mm]
&1 &&- \frac{1}{2^6}&&+\frac{1}{3^6}&&-\frac{1}{4^6}&&+\frac{1}{5^6}&&-\text{etc.}&&=~\frac{C- 2\mathfrak{C}}{1 \cdot 2 \cdots 6}&&\cdot \frac{\pi^6}{4}&&=\frac{(2^5-1)\mathfrak{C}}{1 \cdot 2 \cdots 6}&&\cdot \pi^6 \\[3mm]
&1 &&- \frac{1}{2^8}&&+\frac{1}{3^8}&&-\frac{1}{4^8}&&+\frac{1}{5^8}&&-\text{etc.}&&=~\frac{D- 2\mathfrak{D}}{1 \cdot 2 \cdots 8}&&\cdot \frac{\pi^8}{4}&&=\frac{(2^7-1)\mathfrak{D}}{1 \cdot 2 \cdots 8}&&\cdot \pi^8 \\[3mm]
&1 &&- \frac{1}{2^{10}}&&+\frac{1}{3^{10}}&&-\frac{1}{4^{10}}&&+\frac{1}{5^{10}}&&-\text{etc.}&&=\frac{E- 2\mathfrak{E}}{1 \cdot 2 \cdots 10}&&\cdot \frac{\pi^{10}}{4}&&=\frac{(2^9-1)\mathfrak{E}}{1 \cdot 2 \cdots 10}&&\cdot \pi^{10} \\
& && && && && && \text{etc.}
\end{alignat*}

\paragraph*{§188}

As up to this point we contemplated  series whose terms are the products of terms from the geometric progression $p$, $p^2$, $p^3$ etc. and terms of the arbitrary series $a$, $b$, $c$ etc., so we will be able to investigate the series whose general terms are the products of terms of two arbitrary series  one of which known. Let this series be known

\begin{alignat*}{9}
&1    && 2    && 3            &&z \\
&A~+~ && B~+~ && C~+ \cdots +~&&Z,
\end{alignat*}
the other on the other hand unknown, which reads 

\begin{equation*}
a+b+c+\cdots +z,
\end{equation*}
and let the sum of this series be in question

\begin{equation*}
Aa+Bb+Cc+\cdots +Zz,
\end{equation*}
which we want to put $=Zs$. Let the penultimate term in the known series be $=Y$  and having written $x-1$ instead of $x$ the expression of the sum $S.Zz$ will go over into

\begin{equation*}
Y \left(s- \frac{ds}{dx}+\frac{dds}{2dx^2}-\frac{d^3 s}{6dx^3}+\frac{d^4 s}{24 dx^4}-\text{etc.}\right).
\end{equation*}
Since this sum expresses the series $Zs$ decreased by the last term $Zz$, it will be

\begin{equation*}
Zs- Zz= Ys- \frac{Y ds}{dx}+\frac{Ydds}{2 dx^2}-\frac{Yd^3 s}{6dx^3}+\text{etc.},
\end{equation*}
which equation contains the relation how the sum $Zs$ depends on $Y$, $Z$ and $z$.

\paragraph*{§189}

To resolve this equation,  first neglect the differential terms and it will be

\begin{equation*}
s= \frac{Zz}{Z-Y};
\end{equation*}
put this value $\frac{Zz}{Z-Y}=P^{\text{I}}$  and let the true sum be $s= P^{\text{I}}+p$; having substituted this value in the equation, it will be

\begin{alignat*}{9}
& (Z-Y)p &&~=~&&-\frac{YdP^{\text{I}}}{dx}&&+\frac{YddP^{\text{I}}}{2 dx^2}&&-\text{etc.} \\[3mm]
&        &&   &&-\frac{Ydp}{dx}&&+\frac{Yddp}{2dx^2}&&-\text{etc.};
\end{alignat*}
add $YP^{\text{I}}$ on both sides, and since $P^{\text{I}}-\frac{d P^{\text{I}}}{dx}+\frac{Y ddP^{\text{I}}}{2dx^2}-$etc. is the value of $P^{\text{I}}$ which results, if  one writes $x-1$ instead of $x$, let this value be $P$ and it will be

\begin{equation*}
(Z-Y)p+YP^{\text{I}}=YP- \frac{Ydp}{dx}+\frac{Yddp}{2dx^2}-\text{etc.},
\end{equation*}
whence neglecting the differentials, it will be

\begin{equation*}
p=\frac{Y(P-P^{\text{I}})}{Z-Y}.
\end{equation*}
Put $\frac{Y(P-P^{\text{I}})}{Z-Y}=Q^{\text{I}}$ and let $p= Q^{\text{I}}+q$; it will be 

\begin{equation*}
(Z-Y)q =- \frac{Y(d Q^{\text{I}}+dq)}{dx}+\frac{Y(ddQ^{\text{I}}+ddq)}{2dx^2}-\text{etc.}
\end{equation*}
and having put $Q$ for the value of $Q^{\text{I}}$ it has, if  one writes $x-1$ instead of $x$, it will be

\begin{equation*}
(Z-Y)q+Y Q^{\text{I}} =YQ-\frac{Ydq}{dx}+\frac{Yddq}{2dx^2}-\text{etc.},
\end{equation*}
whence, neglecting the differentials,

\begin{equation*}
q= \frac{Y(Q-Q^{\text{I}})}{Z-Y}.
\end{equation*}
Put $\frac{Y(Q-Q^{\text{I}})}{Z-Y}=R^{\text{I}}$ and let the true value be $q= R^{\text{I}}+r$ and similarly one will find

\begin{equation*}
r=\frac{Y(R-R^{\text{I}})}{Z-Y};
\end{equation*}
and, proceeding this way, the sum in question will be

\begin{equation*}
Zs =Z(P^{\text{I}}+Q^{\text{I}}+R^{\text{I}}+\text{etc.}).
\end{equation*}

\paragraph*{§190}

Therefore, having propounded any series

\begin{equation*}
Aa+Bb+Cc+\cdots +Yy+Zz,
\end{equation*}
its sum will be defined the following way

\begin{alignat*}{9}
&  \quad  \quad \text{Put} \quad && &&\text{having put} \quad x-1 \quad \text{instead of} \quad x \\[3mm]
& \frac{\quad ~~ Zz \quad ~~}{Z-Y}&&=P^{\text{I}} \quad && \text{and let} \quad P^{\text{I}} \quad \text{go over into} \quad P \\[3mm]
& \frac{Y(P-P^{\text{I}})}{Z-Y}&&=Q^{\text{I}} \quad && \text{and let} \quad Q^{\text{I}} \quad \text{go over into} \quad Q \\[3mm]
& \frac{Y(Q-Q^{\text{I}})}{Z-Y}&&=R^{\text{I}} \quad && \text{and let} \quad R^{\text{I}} \quad \text{go over into} \quad R \\[3mm]
& \frac{Y(R-R^{\text{I}})}{Z-Y}&&=S^{\text{I}} \quad && \text{and let} \quad S^{\text{I}} \quad \text{go over into} \quad S \\
& &&\text{etc.}
\end{alignat*}
Having found these values, the sum of the series will be

\begin{equation*}
=ZP^{\text{I}}+ZQ^{\text{I}}+ZR^{\text{I}}+ZS^{\text{I}}+\text{etc.}
\end{equation*}
$+$ a constant rendering the sum $=0$ for $x=0$, or, what is  the same, which forces the sum to be correct in a certain case.

\paragraph*{§191}

This formula, since it does not contain any differentials, is applied most easily in the most cases and  will often even exhibit the true sum. So if this series is propounded

\begin{equation*}
p+4p^2+9p^3+16p^4+\cdots +x^2p^x,
\end{equation*} 
let $Z=p^x$ and $z=x^2$; it will be $Y=p^{x-1}$ and $\frac{Z}{Z-Y}=\frac{p}{p-1}$ and $\frac{Y}{Z-Y}=\frac{1}{p-1}$. Therefore, it will be

\begin{alignat*}{9}
&P^{\text{I}}&&=\frac{px^2}{p-1}  \qquad && P &&=\frac{pxx-2px+p}{p-1} \\[3mm]
&Q^{\text{I}}&&=\frac{-2px+p}{(p-1)^2} \qquad && Q &&=\frac{-2px+3p}{(p-1)^2} \\[3mm]
&R^{\text{I}}&&=\frac{2p}{(p-1)^3} \qquad && R &&=\frac{2p}{(p-1)^3} \\[3mm]
&S^{\text{I}}&&=0
\end{alignat*}
and all remaining vanish; therefore, the sum will be

\begin{equation*}
=p^x \left(\frac{px^2}{p-1}-\frac{2px-p}{(p-1)^2}+\frac{2p}{(p-1)^3}\right)-\frac{p}{(p-1)^2}-\frac{2p}{(p-1)^3}
\end{equation*}
\begin{equation*}
=p^{x+1} \left(\frac{x^2}{p-1}-\frac{2x}{(p-1)^2}+\frac{p+1}{(p-1)^3}\right)-\frac{p(p+1)}{(p-1)^3}.
\end{equation*}
as we already found above [§ 176].

\paragraph*{§192}

In the same way we got to this expression of the sum we will be able to find another expression, if the propounded series is not composed of two others; this can mainly be used in those cases, in which  one gets to vanishing denominators in the preceding expression. Therefore, let this series be propounded

\begin{alignat*}{9}
 &  &&1 &&2 && 3 && 4 &&  x \\
 &s~=~ &&a~+~ &&b~+~ &&c~+~ &&d~+~ \cdots ~+~&&z;
\end{alignat*}
since having written $x-1$ instead of $x$ the sum is truncated by the last term, it will be

\begin{equation*}
s-z =s- \frac{ds}{dx}+\frac{dds}{2dx^2}-\frac{d^3 s}{6 dx^3}+\frac{d^4 s}{24dx^4}-\text{etc.}
\end{equation*}
or

\begin{equation*}
z=\frac{ds}{dx}-\frac{dds}{2dx^2}+\frac{d^3s}{6dx^3}-\frac{d^4 s}{24 dx^4}+\text{etc.}
\end{equation*}
Since here the sum $s$  does not occur, neglect the higher differentials and it will be $s=\int zdx$; put $\int zdx = P^{\text{I}}$ whose value  goes over into $P$, if   one writes $x-1$ instead of $x$, and let the true value be $P^{\text{I}}$; it will be

\begin{equation*}
z=\frac{d P^{\text{I}}}{dx}-\frac{ddP^{\text{I}}}{2dx^2}+\text{etc.}+\frac{dp}{dx}-\frac{ddp}{2dx^2}+\text{etc.};
\end{equation*}
since

\begin{equation*}
P= P^{\text{I}}-\frac{dP^{\text{I}}}{dx}+\frac{ddP^{\text{I}}}{2dx^2}-\text{etc.},	
\end{equation*}
it will be

\begin{equation*}
z-P^{\text{I}}+P=\frac{dp}{dx}-\frac{ddp}{2dx^2}+\text{etc.},
\end{equation*}
whence

\begin{equation*}
p=\int (z-P^{\text{I}}-P)dx.
\end{equation*}
If one further puts $\int (z-P^{\text{I}}+P)dx=Q^{\text{I}}$ and  this value  goes over into $Q$, having written $x-1$ instead of $x$, let

\begin{equation*}
\int (z-P^{\text{I}}+P-Q^{\text{I}}+Q)dx=R^{\text{I}}=Q^{\text{I}}-\int(Q^{\text{I}}-Q)dx,
\end{equation*}
furthermore,

\begin{equation*}
R^{\text{I}}-\int(R^{\text{I}}-R)dx=S^{\text{I}}
\end{equation*}
etc.; the sum in question will be

\begin{equation*}
s=P^{\text{I}}+Q^{\text{I}}+R^{\text{I}}+S^{\text{I}}+\text{etc.}+\text{Const.},
\end{equation*}
by means of which constant one single case is to be forced to be correct.

\paragraph*{§193}

Having changed the letters a little bit, this summation reduces to this. Having propounded this series to be summed

\begin{alignat*}{9}
& && 1 && 2 && 3 && 4 && x \\
&s~=~ &&a~+~&& b~+~&&c~+~&& d~+~ \cdots ~+~ &&z
\end{alignat*}
\begin{alignat*}{9}
& &&\text{put} \quad  && && \text{having put} \quad x-1 \quad  \text{instead of} \quad x \\[3mm]
& &&\int zdx && =P \quad && \text{and let} \quad P \quad \text{go over into} \quad p \\[3mm]
&P~-~&&\int (P-p)dx && =Q \quad && \text{and let} \quad Q \quad \text{go over into} \quad q \\[3mm]
&Q~-~&&\int (Q-q)dx && =R \quad && \text{and let} \quad R \quad \text{go over into} \quad r \\
& && && \text{etc.}
\end{alignat*}
having found which values the sum in question will be

\begin{equation*}
s=P+Q+R+S+\text{etc.}
\end{equation*}
this expression will immediately show the sum, if these integral formulas can be exhibited.  To give an example of its application, let $z=xx+x$ and it will be

\begin{equation*}
P=\frac{1}{3}x^3+\frac{1}{2}xx, \quad p=\frac{1}{3}x^3-\frac{1}{2}xx+\frac{1}{6}, \quad P-p=xx-\frac{1}{6}
\end{equation*}
and

\begin{equation*}
\int (P-p)dx=\frac{1}{3}x^3-\frac{1}{6}x;
\end{equation*}
\begin{equation*}
Q=\frac{1}{2}xx+\frac{1}{6}x, \quad q=\frac{1}{2}xx-\frac{5}{6}x+\frac{1}{3}, \quad Q-q=x-\frac{1}{3}
\end{equation*}
and

\begin{equation*}
\int (Q-q)dx = \frac{1}{2}xx-\frac{1}{3}x;
\end{equation*}
\begin{equation*}
R=\frac{1}{2}x, \quad r=\frac{1}{2}x-\frac{1}{2}, \quad R-r=\frac{1}{2}
\end{equation*}
and

\begin{equation*}
\int (R-r)dx =\frac{1}{2}x;
\end{equation*}
$S=0$ and the remaining values vanish. Therefore, the sum in question will be

\begin{equation*}
\begin{large}
\left.
\renewcommand{\arraystretch}{1,5}
\setlength{\arraycolsep}{0.0mm}
\begin{array}{llllllllllllllll}
\frac{1}{3}x^3 &~+~& \frac{1}{2}xx \\
               &~+~& \frac{1}{2}xx &~+~&\frac{1}{6}x \\
               &   &               &~+~&\frac{1}{2}x 
\end{array}
\right\rbrace
\end{large}
 = \frac{1}{3}x^3+xx+\frac{2}{3}x=\frac{1}{3}x(x+1)(x+2).
\end{equation*}
And this way the sum of all series whose general terms are polynomial functions of $x$ can be found by means of iterated integration. From these considerations it is easily seen what a broad field the doctrine of the summation of series actually is and several volumes will not suffice to capture all methods which already exist and which can still be thought of.

\paragraph*{§194}

Up to now we investigated sums of series from the first term up to the one whose index is $x$  and if we know this sum, by putting $x=\infty$, the  sum of the series continued to infinity is also found. But in most cases this is achieved more easily, if not the sum of terms from the first to that one whose index is $x$ but the sum of the terms from the one whose index is $x$ up to infinity is found; for, in this case especially the last expressions we found become more tractable. Therefore, let a series be propounded whose general term or the one corresponding to the index $x$ we want to be $=z$, the following corresponding to the index $x+1$ then is $=z^{\text{I}}$ and the ones following after this one are $z^{\text{II}}$, $z^{\text{III}}$ etc. and let the sum of this infinite series be in question

\begin{equation*}
\renewcommand{\arraystretch}{1,2}
\setlength{\arraycolsep}{0.0mm}
\begin{array}{rcccccccl}
 x & & x+1  & & x+2& & x+3 & & \text{etc.} \\
 s=z &\quad + \quad & z^{\text{I}} & \quad + \quad & z^{\text{II}} & \quad + \quad & z^{\text{III}} & \quad + \quad & \text{etc. to infinity}
\end{array}
\end{equation*}
Therefore, the sum $s$ will be a function of $x$; if in this one writes $x+1$ instead of $x$, the first sum truncated by the term $z$ will result. Since by this change $s$ goes over into

\begin{equation*}
s+\frac{ds}{dx}+\frac{dds}{2dx^2}+\text{etc.},
\end{equation*}
it will be

\begin{equation*}
s-z =s +\frac{ds}{dx}+\frac{dds}{2 dx^2}+\frac{d^3 s}{6dx^3}+\frac{d^4s}{24 dx^4}+\frac{d^5 s}{120 dx^5}+\text{etc.}
\end{equation*}
or

\begin{equation*}
0=z+\frac{ds}{dx}+\frac{dds}{2dx^2}+\frac{d^3 s}{6dx^3}+\frac{d^4 s}{24 dx^4}+\frac{d^5 s}{120 dx^5}+\text{etc.}
\end{equation*}

\paragraph*{§195}

If we now argue as before, having neglected the higher differentials, it will be $s=C-\int z dx$. Therefore, put $\int zdx =P$ and let the true value be $s=C-P+p$; it will be

\begin{alignat*}{9}
&0 &&=z &&-\frac{dP}{dx}&&-\frac{ddP}{2dx^2}&&-\frac{d^3 P}{6dx^3}&&-\text{etc.} \\[3mm]
&  &&   &&+\frac{dp}{dx}&&+\frac{ddp}{2 dx^2}&&+\frac{d^3 p}{6 dx^3}&&+\text{etc.} 
\end{alignat*}
Let $P$ go over into $P^{\text{I}}$, if  one writes $x+1$ instead of $x$, and it will be

\begin{equation*}
0= z+P-P^{\text{I}}+\frac{dp}{dx}+\frac{ddp}{2 dx^2}+\frac{d^3 p}{6 dx^3}+\text{etc.}
\end{equation*}
Therefore, having neglected the higher differentials, it will be $p= \int (P^{\text{I}}-P)dx-P$. Set $ \int (P^{\text{I}}-P)dx-P=-Q$ and let $p=-Q+q$; it will be

\begin{equation*}
0=z+P-P^{\text{I}}-\frac{dQ}{dx}-\frac{ddQ}{2dx^2}-\text{etc.}+\frac{dq}{dx}+\frac{ddq}{2dx^2}+\text{etc.}
\end{equation*} 
Let $Q$ go over into $Q^{\text{I}}$, if one writes $x+1$ instead of $x$, and it will be

\begin{equation*}
0=z+P-P^{\text{I}}+Q-Q^{\text{I}}+\frac{dq}{dx}+\frac{ddq}{2dx^2}+\text{etc.},
\end{equation*}
whence $q=\int (Q^{\text{I}}-Q)dx-Q$ follows. Therefore, if the sign $\text{I}$ attached to the quantity denotes its value which it takes having put $x+1$ instead of $x$, and one puts

\begin{alignat*}{9}
&  && \int zdx &&=P \\
&P~-~&& \int (P^{\text{I}}-P)dx &&=Q \\
&Q~-~&& \int (Q^{\text{I}}-Q)dx &&=R \\
&R~-~&& \int (R^{\text{I}}-R)dx &&=S \\
&    && \text{etc.},
\end{alignat*}
the sum of the propounded series $z+z^{\text{I}}+z^{\text{II}}+z^{\text{III}}+z^{\text{IV}}+\text{etc.}$ will be

\begin{equation*}
=C-P-Q-R-S-\text{etc.},
\end{equation*}
where the constant $C$ has to defined in such a way that for $x=\infty$ the total sum vanishes. But since the application of this expression requires integrations, it is not possible to show its use here.

\paragraph*{§196}

But to avoid integral formulas, let us set the sum of the series $=ys$ while $y$ is a known function of $x$ whose values $y^{\text{I}}$, $y^{\text{II}}$ etc., which result by writing $x+1$, $x+2$ etc. instead of $x$, will be known. If one now writes $x+1$ instead of $x$, the above series truncated by the first term will result, whose sum will therefore be

\begin{equation*}
y^{\text{I}}\left(s+\frac{ds}{dx}+\frac{dds}{2dx^2}+\frac{d^3 s}{6dx^3}+\text{etc.}\right)=ys-z
\end{equation*}
or

\begin{equation*}
z+\frac{y^{\text{I}}ds}{dx}+\frac{y^{\text{I}}dds}{2dx^2}+\frac{y^{\text{I}}d^3s}{6dx^3}+\text{etc.}=(y-y^{\text{I}})s,
\end{equation*}
whence, having neglected the differentials, $s=\frac{z}{y-y^{\text{I}}}$ results. Set $\frac{z}{y^{\text{I}}-y}=P$ and let the true value be $s=-P+p$; it will be

\begin{equation*}
\begin{large}
\left.
\renewcommand{\arraystretch}{2,5}
\setlength{\arraycolsep}{0.2mm}
\begin{array}{llllllllllllllll}
- & \dfrac{y^{\text{I}}dP}{dx}&-& \dfrac{y^{\text{I}}ddP}{2dx^2}&-&\dfrac{y^{\text{I}}d^3P}{6dx^3}&-&\text{etc.} \\
+ & \dfrac{y^{\text{I}}dp}{dx}&+& \dfrac{y^{\text{I}}ddp}{2dx^2}&+&\dfrac{y^{\text{I}}d^3p}{6dx^3}&-&\text{etc.} 
\end{array}
\right\rbrace
\end{large}
 = (y-y^{\text{I}})p
\end{equation*}
and hence

\begin{equation*}
\frac{y^{\text{I}}dp}{dx}+\frac{y^{\text{I}}ddp}{2dx^2}+\frac{y^{\text{I}}d^3p}{6dx^3}+\text{etc.}=y^{\text{I}}(P^{\text{I}}-P)-(y^{\text{I}}-y)p.
\end{equation*}
Set $Q=\frac{y^{\text{I}}(P^{\text{I}}-P)}{y^{\text{I}}-y}$ and let $p=Q+q$; it will be 

\begin{equation*}
y^{\text{I}}(Q^{\text{I}}-Q)+y^{\text{I}}\left(\frac{dq}{dx}+\frac{ddq}{2dx^2}+\text{etc.}\right)=-(y^{\text{I}}-y)q.
\end{equation*}
Put $R=\frac{y^{\text{I}}(Q^{\text{I}}-Q)}{y^{\text{I}}-y}$ and let be $q=-R+r$.\\[2mm]
And if we proceed this way, the sum of the propounded series

\begin{equation*}
z+z^{\text{I}}+z^{\text{II}}+z^{\text{III}}+z^{\text{IV}}+\text{etc.}
\end{equation*}
will be found the following way. Having taken an arbitrary function of $x$ and called this function $=y$, set 

\begin{alignat*}{9}
&P &&=  \frac{\quad \quad  z \quad \quad}{y^{\text{I}}-y}&&=\frac{z}{\Delta y} \\[3mm]
&Q &&= \frac{y^{\text{I}}(P^{\text{I}}-P)}{y^{\text{I}}-y}&&=\frac{y \Delta P}{\Delta y}&&+\Delta P \\[3mm]
&R &&= \frac{y^{\text{I}}(Q^{\text{I}}-Q)}{y^{\text{I}}-y}&&=\frac{y \Delta Q}{\Delta y}&&+\Delta Q \\[3mm]
&S &&= \frac{y^{\text{I}}(R^{\text{I}}-R)}{y^{\text{I}}-y}&&=\frac{y \Delta R}{\Delta y}&&+\Delta R \\
& && &&\text{etc.}
\end{alignat*}
And hence the sum in question will be

\begin{equation*}
=C-Py+Qy-Ry+Sy-\text{etc.}
\end{equation*}
having taken a constant of such a kind for $C$ that for $x=\infty$ the sum vanishes. 

\paragraph*{§197}

Assume $y=a^x$; because of $y^{\text{I}}=a^{x+1}$, it will be $y^{\text{I}}-y=a^x(a-1)$, whence it will become

\begin{alignat*}{9}
&P &&= \frac{z}{a^x(a-1)} && && P^{\text{I}}&&=\frac{z^{\text{I}}}{a^{x+1}(a-1)} \\[3mm]
&Q &&= \frac{a(P^{\text{I}}-P)}{a-1}&&= \frac{z^{\text{I}}-az}{a^x(a-1)^2} \quad && Q^{\text{I}}&&=\frac{z^{\text{II}}-az^{\text{I}}}{a^{x+1}(a-1)^2} \\[3mm]
&R &&= \frac{a(Q^{\text{I}}-Q)}{a-1}&&= \frac{z^{\text{II}}-2az^{\text{I}}+aaz}{a^x(a-1)^3} \quad && R^{\text{I}}&&=\frac{z^{\text{III}}-2az^{\text{II}}+aaz^{\text{I}}}{a^{x+1}(a-1)^2} \\[3mm]
&S &&=\frac{a(R^{\text{I}}-R)}{a-1}&&=\frac{z^{\text{III}}-3az^{\text{II}}+3a^2 z^{\text{I}}-a^3}{a^x(a-1)^4} \qquad \\
& && &&\text{etc.}
\end{alignat*}
Therefore, the sum of the propounded series will be

\begin{equation*}
C-\frac{z}{a-1}+\frac{z^{\text{I}}-az}{(a-1)^2}-\frac{z^{\text{II}}-2az^{\text{I}}+a^2z}{(a-1)^3}+\frac{z^{\text{III}}-3a z^{\text{II}}+3a^2 z^{\text{I}}-a^3z}{(a-1)^4}-\text{etc.}
\end{equation*}
This same expression of the sum was indeed found already above in the first chapter. Hence, taking other values for $y$, infinitely many other expressions can be found, whence the one, which is the most convenient in each case, can be selected.


\title{On the Use of Differential Calculus in the Formation of Series
\footnote{
Original title: "'De Usu Calculi Differentialis in Formandis Seriebus"', first published as part of the book \textit{Institutiones calculi differentialis cum eius usu in analysi finitorum ac doctrina serierum}, 1755, reprinted in \textit{Opera Omnia}: Series 1, Volume 10, pp. 396 - 421, Eneström-Number E212, translated by: Alexander Aycock for the "`Euler-Kreis Mainz"'}}
\author{Leonhard Euler}
\date{ }
\maketitle

\paragraph*{§198}

Until now we only considered one single application of differential calculus in the doctrine of series which was the formation of series  and  which we mentioned above already, when there was the question how to expand  the fraction whose denominator is an arbitrary power of a certain function into a series. But this method is similar to that one we already used several times, where the fraction to be converted into a series is set equal to a certain series with  coefficients to  determined from the constituted equality. But this determination is often  simplified tremendously, if, before it is actually done, the equation  is  differentiated once  and sometimes even twice. Since this method has very broad applications in integral calculus, let us explain it here more diligently.

\paragraph*{§199}

Therefore, at first let us repeat what we discussed above on the expansion of fractions into series without the application of differential calculus. Let an arbitrary fraction be propounded, i.e.

\begin{equation*}
\frac{A+Bx+Cx^2+Dx^3+\text{etc.}}{\alpha + \beta x +\gamma x^2+ \delta x^3 +\varepsilon x^4 +\text{etc.}}=s,
\end{equation*}
which is to be converted into a powers series in $x$. Assume an undetermined series for $s$

\begin{equation*}
s= \mathfrak{A}+\mathfrak{B}x+ \mathfrak{C}x^2+\mathfrak{D}x^3+\mathfrak{E}x^4+\mathfrak{F}x^5+\mathfrak{G}x^6+\text{etc.}
\end{equation*}
Therefore, since, having removed the fraction by multiplication, 

\begin{equation*}
A+Bx+Cx^2+Dx^3+Ex^4+Fx^5+Gx^6+\text{etc.}
\end{equation*}
\begin{equation*}
=s(\alpha +\beta x +\gamma x^2 +\delta x^3 +\varepsilon x^4 +\zeta x^5 +\eta x^6 +\text{etc.}),
\end{equation*}
if  the assumed series is substituted for $s$, the following equation results

\begin{equation*}
\renewcommand{\arraystretch}{1,4}
\setlength{\arraycolsep}{0.0mm}
\begin{array}{rlllllllllllllll}
A &~+~& Bx &~+~& Cx^2 &~+~& Dx^3 &~+~& Ex^4 &~+~ & Fx^5 &~+~& \text{etc.} \\ \hline
~=~ \mathfrak{A}\alpha &~+~& \mathfrak{B}\alpha x &~+~& \mathfrak{C}\alpha x^2 &~+~& \mathfrak{D}\alpha x^3 &~+~& \mathfrak{E}\alpha x^4 &~+~& \mathfrak{F}\alpha x^5 &~+~&\text{etc.}\\
  &~+~
   & \mathfrak{A}\beta  &~+~& \mathfrak{B}\beta  &~+~& \mathfrak{C}\beta  &~+~& \mathfrak{D}\beta &~+~& \mathfrak{E}\beta &~+~&\text{etc.}\\
      & &   &~+~& \mathfrak{A}\gamma  &~+~& \mathfrak{B}\gamma  &~+~& \mathfrak{C}\gamma &~+~& \mathfrak{D}\gamma &~+~&\text{etc.}\\
            & &   & &    &~+~& \mathfrak{A}\delta  &~+~& \mathfrak{B}\delta &~+~& \mathfrak{C}\delta &~+~&\text{etc.}\\
             & &   & &    & &   &~+~& \mathfrak{A}\varepsilon &~+~& \mathfrak{B}\varepsilon &~+~&\text{etc.}\\
              & &   & &    & &   & &  &~+~& \mathfrak{A}\zeta &~+~&\text{etc.}\\
\end{array}
\end{equation*} 
Therefore, having equated each term containing the same powers of $x$, it will be

\begin{alignat*}{9}
&\mathfrak{A}\alpha &&- A &&=0 \\
&\mathfrak{B}\alpha &&+\mathfrak{A}\beta &&- B &&=0 \\
&\mathfrak{C}\alpha &&+\mathfrak{B}\beta &&+\mathfrak{A}\gamma &&- C &&=0 \\
& \mathfrak{D}\alpha &&+\mathfrak{C}\beta &&+ \mathfrak{B}\gamma &&+ \mathfrak{A}\delta &&- D &&=0 \\
&\mathfrak{E}\alpha &&+\mathfrak{D}\beta &&+\mathfrak{C}\gamma &&+\mathfrak{B}\delta &&+\mathfrak{A}\varepsilon &&- E =0 \\
& && && && &&\text{etc.}
\end{alignat*}
from which equations the assumed coefficients $\mathfrak{A}$, $\mathfrak{B}$, $\mathfrak{C}$, $\mathfrak{D}$ etc. are determined, and so the infinite series 

\begin{equation*}
\mathfrak{A}+\mathfrak{B}x+\mathfrak{C}x^2+\mathfrak{D}x^3+\mathfrak{E}x^4+\text{etc.}
\end{equation*}
equal to the propounded fraction $s$ is found. And this form, if both the numerator and the denominator of the fraction $s$ consist of a finite number of terms, contains all recurring series, which were treated in a lot greater detail above.

\paragraph*{§200}

But if either the numerator or the denominator or both were raised to an arbitrary power, then  the series is found rather difficultly this way, since the task, if not a binomial function was raised, becomes very laborious.  But by means of differential calculus this work can be simplified. At first, let the fraction consist only of a numerator and let 

\begin{equation*}
s=(A+Bx+Cxx)^n,
\end{equation*}
whence the power series in $x$ equal to this power of the trinomial is to be found; it is plain that the series will be finite, if the exponent $n$ was an positive integer. Again, assume an indefinite series for $s$, i.e.

\begin{equation*}
S =\mathfrak{A}+\mathfrak{B}x+\mathfrak{C}x^2+\mathfrak{D}x^3+\mathfrak{E}x^4+\mathfrak{F}x^5+\mathfrak{G}x^6+\text{etc.},
\end{equation*}
whose first term $\mathfrak{A}$ is known to be $=A^n$; for, if one puts $x=0$, from the first propounded form $s=A^n$, but from the assumed series $s=\mathfrak{A}$. But this determination of the first term is to be derived from the nature of the series itself, if we want to use differentials, since it is not possible to determine the first coefficient from the differential, as it will be seen soon.

\paragraph*{§201}

Since $S=(A+Bx+Cx^2)^n$,  taking logarithms, it will be

\begin{equation*}
\log s = n \log (A+Bx+Cx^2)
\end{equation*} 
and hence, having taken the differentials, one will have

\begin{equation*}
\frac{ds}{s}=\frac{nBdx+2nCxdx}{A+Bx+Cx^2} \quad \text{or} \quad (A+Bx+Cx^2)\frac{ds}{dx}= ns(B+2Cx).
\end{equation*}
But from the assumed series

\begin{equation*}
\frac{ds}{dx}=\mathfrak{B}+2 \mathfrak{C}x+3 \mathfrak{D}x^2+4\mathfrak{E}x^3+5 \mathfrak{F}x^4+\text{etc.}
\end{equation*}
Therefore, if this series  is substituted for $\frac{ds}{dx}$ and  for $s$ the assumed series is substituted, the following equation will result

\begin{equation*}
\renewcommand{\arraystretch}{1,4}
\setlength{\arraycolsep}{0.0mm}
\begin{array}{rrrrrrrrrrrrrrrrr}
A\mathfrak{B}&~+~& 2A \mathfrak{C}&x&~+~&3A\mathfrak{D}&x^2&~+~& 4A\mathfrak{E}&x^3&~+~& 5A\mathfrak{F}&x^4&~+~\text{etc.} \\
&~+~& B \mathfrak{B}& &~+~&2B\mathfrak{C}& &~+~& 3B\mathfrak{D}& &~+~& 4B\mathfrak{E}& &~+~\text{etc.} \\
& &  & &~+~&C\mathfrak{B}& &~+~& 2C\mathfrak{C}& &~+~& 3C\mathfrak{D}& &~+~\text{etc.} \\ \hline
=nB\mathfrak{A}&~+~& nB \mathfrak{B}& &~+~&nB\mathfrak{C}& &~+~& n B\mathfrak{D}& &~+~& nB\mathfrak{E}& &~+~\text{etc.} \\
&~+~& 2nC \mathfrak{A}& &~+~&2nC\mathfrak{B}& &~+~& 2n C\mathfrak{C}& &~+~& 2nC\mathfrak{D}& &~+~\text{etc.} \\
\end{array}
\end{equation*} 
Therefore, having equated the terms of the same power of $x$, it will be

\begin{alignat*}{9}
&\mathfrak{B}&&=\frac{nB\mathfrak{A}}{A}\\[3mm]
&\mathfrak{C}&&=\frac{(n-1)B \mathfrak{B}+2n C\mathfrak{A}}{2A}\\[3mm]
&\mathfrak{D}&&=\frac{(n-2)B \mathfrak{C}+(2n-1)C\mathfrak{B}}{3A}\\[3mm]
&\mathfrak{E}&&=\frac{(n-3)B \mathfrak{D}+(2n-2)C\mathfrak{C}}{4A}\\[3mm]
&\mathfrak{F}&&=\frac{(n-4)B \mathfrak{E}+(2n-3)C\mathfrak{D}}{5A} \\
&            && \text{etc.}
\end{alignat*}
Therefore, since, as we saw before,  $\mathfrak{A}=A^n$, it will be $\mathfrak{B}= n A^{n-1}B$ and hence the remaining coefficients will successively be defined. But the law they follow is obvious from these formulas which would have remained immensely obscure, if we would have wanted to actually expand the trinomial.

\paragraph*{§202}

The same method succeeds, if any polynomial function has to raised to a certain power. Let

\begin{equation*}
s=(A+Bx+Cx^2+Dx^3+Ex^4+\text{etc.})^n
\end{equation*}
and assume

\begin{equation*}
s+\mathfrak{A}+\mathfrak{B}x+\mathfrak{C}x^2+\mathfrak{D}x^3+\mathfrak{E}x^4+\text{etc.};
\end{equation*}
it will be $\mathfrak{A}=A^n$ which value is concluded, if one puts $x=0$. Now, having taken the logarithms,  their differentials as before will be found to be

\begin{equation*}
\frac{ds}{s}=\frac{nBdx+2nCxdx+3nDx^2dx+4nEx^3dx+\text{etc.}}{A+Bx+Cx^2+Dx^3+Ex^4+\text{etc.}}
\end{equation*}
or

\begin{equation*}
(A+Bx+Cx^2+Dx^3+Ex^4+\text{etc.})\frac{ds}{dx}
\end{equation*}
\begin{equation*}
=s(nB+2nCx+3nDx^2+4nEx^3+\text{etc.}).
\end{equation*}
Therefore, since 

\begin{equation*}
\frac{ds}{dx}=\mathfrak{B}+2\mathfrak{C}x+3\mathfrak{D}x^2+4\mathfrak{E}x^3+5 \mathfrak{F}x^4+\text{etc.},
\end{equation*}
having substituted these series for $s$ and $\frac{ds}{dx}$, it will be

\begin{equation*}
\renewcommand{\arraystretch}{1,4}
\setlength{\arraycolsep}{0.0mm}
\begin{array}{rrrrrrrrrrrrrrrrr}
A\mathfrak{B}&~+~& 2A \mathfrak{C}&x&~+~&3A\mathfrak{D}&x^2&~+~& 4A\mathfrak{E}&x^3&~+~& 5A\mathfrak{F}&x^4&~+~\text{etc.} \\
&~+~& B \mathfrak{B}& &~+~&2B\mathfrak{C}& &~+~& 3B\mathfrak{D}& &~+~& 4B\mathfrak{E}& &~+~\text{etc.} \\
& &  & &~+~&C\mathfrak{B}& &~+~& 2C\mathfrak{C}& &~+~& 3C\mathfrak{D}& &~+~\text{etc.} \\
& &  & & & & &~+~& D\mathfrak{B}& &~+~& 2D\mathfrak{C}& &~+~\text{etc.} \\ 
& &  & & & & & & & &~+~& E\mathfrak{B}& &~+~\text{etc.} \\ \hline
=nB\mathfrak{A}&~+~& nB \mathfrak{B}& &~+~&nB\mathfrak{C}& &~+~& n B\mathfrak{D}& &~+~& nB\mathfrak{E}& &~+~\text{etc.} \\
&~+~& 2nC \mathfrak{A}& &~+~&2nC\mathfrak{B}& &~+~& 2n C\mathfrak{C}& &~+~& 2nC\mathfrak{D}& &~+~\text{etc.} \\
& & & &~+~&3nD\mathfrak{A}& &~+~& 3n D\mathfrak{B}& &~+~& 3nD\mathfrak{C}& &~+~\text{etc.} \\
& & & & & & &~+~& 4n E\mathfrak{A}& &~+~& 4nE\mathfrak{B}& &~+~\text{etc.} \\
& & & & & & & & & &~+~& 5nF\mathfrak{A}& &~+~\text{etc.} \\
\end{array}
\end{equation*} 
Therefore, the following determinations are derived

\begin{alignat*}{9}
& &&A\mathfrak{B}&&=n B\mathfrak{A}\\[3mm]
&2&&A \mathfrak{C}&&=(n-1)B\mathfrak{B}+2n C\mathfrak{A}\\[3mm]
&3&&A \mathfrak{D}&&=(n-2)B\mathfrak{C}+(2n-1)C \mathbb{B}+3n D\mathfrak{A} \\[3mm]
&4&&A \mathfrak{E}&&=(n-3)B\mathfrak{D}+(2n-2)C\mathfrak{C}+(3n-1)D\mathfrak{B}+4n E\mathfrak{A}\\[3mm]
&5&&A \mathfrak{F}&&= (n-4) B\mathfrak{E}+(2n-3)C\mathfrak{D}+(3n-2)D\mathfrak{C}+(4n-1)E\mathfrak{B}+5nF\mathfrak{A}\\[3mm]
& && &&\text{etc.},
\end{alignat*}
whence it becomes clear, how these assumed coefficients $\mathfrak{A}$, $\mathfrak{B}$, $\mathfrak{C}$, $\mathfrak{D}$ etc. depend on each other and are hence determined, since $\mathfrak{A}=A^n$.

\paragraph*{§203}

Since, if the quantity $A+Bx+Cx^2+Dx^3+\text{etc.}$ consists of a finite number of terms and the number $n$ was a positive integer, any power also has to consist of a finite number of terms, it is obvious that in this case the formulas just found must finally vanish and, since  all finite terms must occur until the first vanishes, at the same time all following ones must vanish. Let us put that the propounded formula $A+Bx+Cx^2$ is a trinomial and its cube is in question, i.e. that  $n=3$; it will be

\begin{alignat*}{19}
& && &&\mathfrak{A}&&~=~&& &&2&&A^3&& \quad && \text{and hence} \quad && \mathfrak{A}&&=A^3 \\
& &&A&&\mathfrak{B}&&~=~&& &&3&&B\mathfrak{A}&& \quad && \quad && \mathfrak{B}&&=3A^2B \\
&2&&A&&\mathfrak{C}&&~=~&& &&2&&B\mathfrak{B}&&+6C\mathfrak{A} \quad && \quad && \mathfrak{C}&&=3AB^2+3A^2C \\
&3&&A&&\mathfrak{D}&&~=~&& &&1&&B\mathfrak{C}&&+5C\mathfrak{B} \quad && \quad && \mathfrak{D}&&=B^3+6ABC \\
&4&&A&&\mathfrak{E}&&~=~&& &&0&& &&+4C\mathfrak{C} \quad && \quad && \mathfrak{E}&&=3B^2C+3AC^2 \\
&5&&A&&\mathfrak{F}&&~=~&&~-~ && &&B\mathfrak{E}&&+3C\mathfrak{D} \quad && \quad && \mathfrak{F}&&=3BC^2 \\
&6&&A&&\mathfrak{G}&&~=~&&~-~ &&2&&B\mathfrak{F}&&+2C\mathfrak{E} \quad && \quad && \mathfrak{G}&&=C^3 \\
&7&&A&&\mathfrak{H}&&~=~&&~-~ &&3&&B\mathfrak{G}&&+2C\mathfrak{F} \quad && \quad && \mathfrak{H}&&=0 \\
&8&&A&&\mathfrak{I}&&~=~&&~-~ &&4&&B\mathfrak{H}&&+0 \quad && \quad && \mathfrak{I}&&=0 \\
\end{alignat*}
Therefore, since already two letters are $=0$ and any arbitrary of the following letters depends on the two preceding ones, it is plain that all following ones must also vanish. And for this reason the law according to which these coefficients were found to depend on each other is even more noteworthy.

\paragraph*{§204}

If $n$ was a negative number such that $s$ becomes equal to a real fraction, the series will continue to infinity. Therefore, let

\begin{equation*}
s=\frac{1}{(\alpha +\beta x +\gamma x^2+ \delta x^3 +\varepsilon x^4+\text{etc.})^n};
\end{equation*} 
for its value assumes this series

\begin{equation*}
s= \mathfrak{A}+\mathfrak{B}x+\mathfrak{C}x^2+\mathfrak{D}x^3+\mathfrak{E}x^4+\mathfrak{F}x^5+\text{etc.}
\end{equation*}
And if in the above formulas  one puts $\alpha$, $\beta$, $\gamma$, $\delta$ etc. for the letters $A$, $B$, $C$, $D$ etc. and at the same time $n$ becomes negative, the following determinations of the coefficients $\mathfrak{A}$, $\mathfrak{B}$, $\mathfrak{C}$, $\mathfrak{D}$ etc. will result

\begin{alignat*}{9}
& &&\mathfrak{A}&&= \alpha^{-n}=\frac{1}{\alpha^n} \\[3mm]
&\alpha &&\mathfrak{B}&&+n \beta \mathfrak{A}=0 \\[3mm]
&2 \alpha \mathfrak{C}&&+(n+1) \beta \mathfrak{B}&&+2n \gamma \mathfrak{A}=0 \\[3mm]
&3 \alpha \mathfrak{D}&&+(n+2) \beta \mathfrak{C}&&+(2n+1) \gamma \mathfrak{B}+3n \delta \mathfrak{A}=0 \\[3mm]
&4 \alpha \mathfrak{E}&&+(n+3) \beta \mathfrak{D}&&+(2n+2) \gamma \mathfrak{C}+(3n+1) \delta \mathfrak{B}+4n \varepsilon \mathfrak{A}=0 \\[3mm]
&5 \alpha \mathfrak{F}&&+(n+4) \beta \mathfrak{E}&&+(2n+3) \gamma \mathfrak{D}+(3n+2) \delta \mathfrak{C}+(4n+1) \varepsilon \mathfrak{B}+5n \zeta \mathfrak{A}=0 \\
& && &&\qquad \qquad \qquad \text{etc.}
\end{alignat*}
These formulas contain the same law of these coefficients of numbers we already observed above in the \textit{Introductio} and whose validity has therefore been demonstrated rigorously now.

\paragraph*{§205}

Their nature is the same, if the numerator of the fraction was $1$ or even any power of $x$, say $x^m$; for, in the second case it will only be necessary to multiply the series found first $\mathfrak{A}+\mathfrak{B}x+\mathfrak{C}x^2+\mathfrak{D}x^3+\text{etc.}$ by $x^m$. But if the denominator consists of two or more terms, then  we do not observe the law of progression found above: therefore, let us investigate it here by means of differentiation. Hence let

\begin{equation*}
s=\frac{A+Bx+Cx^2+Dx^3+\text{etc.}}{(\alpha +\beta x+ \gamma x^2+ \delta x^3+ \varepsilon x^4 +\text{etc.})^n}
\end{equation*}
and assume the following series for the value of this fraction

\begin{equation*}
s=\mathfrak{A}+\mathfrak{B}x+\mathfrak{C}x^2+\mathfrak{D}x^3+\mathfrak{E}x^4+\mathfrak{F}x^5+\text{etc.};
\end{equation*}
to define its first term $\mathfrak{A}$ put $x=0$ and from the first expression it will be $s=\frac{A}{\alpha^n}$, from the assumed series on the other hand $s=\mathfrak{A}$, whence it is necessary that $\mathfrak{A}=\frac{A}{\alpha^n}$. Having determined this term, the remaining ones will be found by means of differentiation.

\paragraph*{§206}

Having taken logarithms, it will be

\begin{equation*}
\log s = \log (A+Bx+Cx^2+Dx^3+\text{etc.})
\end{equation*}
\begin{equation*}
-n \log (\alpha +\beta x+ \gamma x^2+ \delta x^3+ \varepsilon x^4+\text{etc.})
\end{equation*}
and hence by differentiation this equation will result

\begin{equation*}
\frac{ds}{s}=\frac{Bdx+2Cdx+3Dx^2dx+\text{etc.}}{A+Bx+Cx^2+Dx^3+\text{etc.}}
\end{equation*}
\begin{equation*}
-\frac{n \beta dx +2n \gamma x dx +3n \delta x^2 dx +\text{etc.}}{\alpha +\beta x+\gamma x^2 +\delta x^3+ \text{etc.}}
\end{equation*}
and, having got rid of the fractions by multiplication, it will be

\begin{equation*}
\left\lbrace
\renewcommand{\arraystretch}{1,4}
\setlength{\arraycolsep}{0.0mm}
\begin{array}{rrrrrrrrrrrrrrrrr}
A \alpha &~+~& A \beta &x &~+~& A\gamma &x^2&~+~ & A \delta &x^3 &~+~ & \text{etc.} \\
  &~+~ & B \alpha & &~+~& B\beta & &~+~ & B \gamma & &~+~ & \text{etc.} \\
    & &  &  &~+~& C \alpha & &~+~ & C \beta & &~+~ & \text{etc.} \\
        & &  &  & &  & &~+~ & D \alpha & &~+~ & \text{etc.} \\
\end{array}
\right\rbrace \frac{ds}{dx}
\end{equation*} 
\begin{equation*}
=\left\lbrace
\renewcommand{\arraystretch}{1,4}
\setlength{\arraycolsep}{0.0mm}
\begin{array}{rrrrrrrrrrrrrrrrr}
B \alpha &~+~& B \beta &x &~+~& B\gamma &x^2&~+~ & B \delta &x^3 &~+~ & \text{etc.} \\
  &~+~ & 2C \alpha & &~+~& 2C\beta & &~+~ & 2C \gamma & &~+~ & \text{etc.} \\
    & &  &  &~+~& 3D \alpha & &~+~ & 3D \beta & &~+~ & \text{etc.} \\
        & &  &  & &  & &~+~ & 4E \alpha & &~+~ & \text{etc.} \\
\end{array}
\right\rbrace s
\end{equation*} 
\begin{equation*}
-\left\lbrace
\renewcommand{\arraystretch}{1,4}
\setlength{\arraycolsep}{0.0mm}
\begin{array}{rrrrrrrrrrrrrrrrr}
A \beta &~+~& 2A \gamma &x &~+~& 3A\delta &x^2&~+~ & 4A \varepsilon &x^3 &~+~ & \text{etc.} \\
  &~+~ & B \beta & &~+~& 2B\gamma & &~+~ & 3B \delta & &~+~ & \text{etc.} \\
    & &  &  &~+~& C \beta & &~+~ & 2C \gamma & &~+~ & \text{etc.} \\
        & &  &  & &  & &~+~ & D \beta & &~+~ & \text{etc.} \\
\end{array}
\right\rbrace ns.
\end{equation*} 
Since now  $\frac{ds}{dx}=\mathfrak{B}+2 \mathfrak{C}x+3 \mathfrak{D}x^2+4 \mathfrak{E}x^3+\text{etc.}$,  after the substitutions it will be

\begin{equation*}
\left.
\renewcommand{\arraystretch}{1,4}
\setlength{\arraycolsep}{0.0mm}
\begin{array}{rrrrrrrrrrrrrrrrr}
A \alpha \mathfrak{B} &~+~& n A \beta \mathfrak{A} \\
                      &~-~& B \alpha \mathfrak{A}
\end{array}
\right\rbrace =0
\end{equation*} 
\begin{equation*}
\left.
\renewcommand{\arraystretch}{1,4}
\setlength{\arraycolsep}{0.0mm}
\begin{array}{rrrrrrrrrrrrrrrrr}
2A \alpha \mathfrak{C} &~+~& (n+1) A \beta \mathfrak{B} &~+~& 2n A \gamma \mathfrak{A} \\
  &~+~& 0  B \alpha \mathfrak{B} &~+~& (n-1)B \beta \mathfrak{A} \\                     
   & &  &~-~& 2C \alpha \mathfrak{A} \\                     
\end{array}
\right\rbrace =0
\end{equation*} 
\begin{equation*}
\left.
\renewcommand{\arraystretch}{1,4}
\setlength{\arraycolsep}{0.0mm}
\begin{array}{rrrrrrrrrrrrrrrrr}
3A \alpha \mathfrak{D} &~+~& (n+2) A \beta \mathfrak{C} &~+~& (2n+1) A \gamma \mathfrak{B} &~+~& 3n A \delta \mathfrak{A} \\
    &~+~ &  B \alpha \mathfrak{C} &~+~& n B \beta \mathfrak{B} &~+~& (2n-1) B \gamma \mathfrak{A} \\
   & &   &~-~& C \alpha \mathfrak{B} &~+~& (n-2) C \beta \mathfrak{A} \\
   & &   & &  &~-~& 3 D \alpha \mathfrak{A} \\    
\end{array}
\right\rbrace =0
\end{equation*} 
\begin{equation*}
\left.
\renewcommand{\arraystretch}{1,4}
\setlength{\arraycolsep}{0.0mm}
\begin{array}{rrrrrrrrrrrrrrrrr}
4A\alpha \mathfrak{E}& ~+~& (n+3)A\beta \mathfrak{D}&~+~& (2n+2) A \gamma \mathfrak{C}&~+~& (3n+1) A \delta \mathfrak{B}&~+~& 4n A \varepsilon \mathfrak{A} \\
 & ~+~& B\alpha \mathfrak{D}&~+~& (n+1)B \beta\mathfrak{C}&~+~& 2n B \gamma \mathfrak{B}&~+~& (3n-1) B \delta \mathfrak{A} \\
  &  & &~+~& 0 C \alpha \mathfrak{C}&~+~& (n-1) C \beta \mathfrak{B}&~+~& (3n-2) C \gamma \mathfrak{A} \\
    &  & & & &~-~& 2 D \alpha \mathfrak{B}&~+~& (n-3) D \beta \mathfrak{A} \\
        &  & & & & & &~-~& 4 E \alpha \mathfrak{A} \\
\end{array}
\right\rbrace =0.
\end{equation*} 
\begin{equation*}
\text{etc.}
\end{equation*}
Therefore, the law according to which these formulas proceed is easily seen; for, the first line of each equation follows the same law we had in § 204. But then the coefficients of the second line result, if $n+1$ is subtracted from the above coefficients, and in like manner  the third line is formed from the second line and the following from the upper ones - always by subtracting $n+1$;  but the letters constituting each term are immediately obvious considering the structure of the formulas.

\paragraph*{§207}

But if also the numerator of a fraction was a certain power, i.e.

\begin{equation*}
s=\frac{(A+Bx+Cx^2+Dx^3+\text{etc.})^m}{(\alpha +\beta x+ \gamma x^2+\delta x^3 +\varepsilon x^4+\text{etc.})^n},
\end{equation*}
and one assumes this series

\begin{equation*}
s=\mathfrak{A}+\mathfrak{B}x+\mathfrak{C}x^2+\mathfrak{D}x^3+\mathfrak{E}x^4+\text{etc.},
\end{equation*}
it will be $\mathfrak{A}=\frac{A^m}{\alpha^n}$; but the remaining coefficients will be determined from the following formulas

\begin{equation*}
\left.
\renewcommand{\arraystretch}{1,4}
\setlength{\arraycolsep}{0.0mm}
\begin{array}{rrrrrrrrrrrrrrrrr}
A \alpha \mathfrak{B} &~+~& n A \beta \mathfrak{A} \\
                      &~-~& mB \alpha \mathfrak{A}
\end{array}
\right\rbrace =0
\end{equation*} 
\begin{equation*}
\left.
\renewcommand{\arraystretch}{1,4}
\setlength{\arraycolsep}{0.0mm}
\begin{array}{rrrrrrrrrrrrrrrrr}
2A \alpha \mathfrak{C} &~+~& (n+1) A \beta \mathfrak{B} &~+~& 2n A \gamma \mathfrak{A} \\
  &~-~& (m-1)  B \alpha\mathfrak{B} &~+~& (n-m)B \beta \mathfrak{A} \\                     
   & &  &~-~& 2mC \alpha \mathfrak{A} \\                     
\end{array}
\right\rbrace =0
\end{equation*} 
\begin{equation*}
\left.
\renewcommand{\arraystretch}{1,4}
\setlength{\arraycolsep}{0.0mm}
\begin{array}{rrrrrrrrrrrrrrrrr}
3A \alpha \mathfrak{D} &~+~& (n+2) A \beta \mathfrak{C} &~+~& (2n+1) A \gamma \mathfrak{B} &~+~& 3n A \delta \mathfrak{A} \\
    &~-~ & (m-2) B \alpha \mathfrak{C} &~+~& (n-m+1) B \beta \mathfrak{B} &~+~& (2n-m) B \gamma \mathfrak{A} \\
   & &   &~-~& (2m-1)C \alpha \mathfrak{B} &~+~& (n-2m) C \beta \mathfrak{A} \\
   & &   & &  &~-~& 3m D \alpha \mathfrak{A} \\    
\end{array}
\right\rbrace =0
\end{equation*} 
\begin{equation*}
\left.
\renewcommand{\arraystretch}{1,4}
\setlength{\arraycolsep}{0.0mm}
\begin{array}{rrrrrrrrrrrrrrrrr}
4A\alpha \mathfrak{E}& ~+~& (n+3)A\beta \mathfrak{D}&~+~& (2n+2) A \gamma \mathfrak{C}&~+~& (3n+1) A \delta \mathfrak{B}&~+~& 4n A \varepsilon \mathfrak{A} \\
 & ~-~& (m-3)B\alpha \mathfrak{D}&~+~& (n+m-2)B \beta\mathfrak{C}&~+~& (2n-m+1) B \gamma \mathfrak{B}&~+~& (3n-m) B \delta \mathfrak{A} \\
  &  & &~-~& (2m-2) C \alpha \mathfrak{C}&~+~& (n-2m+1) C \beta \mathfrak{B}&~+~& (n-3m) C \gamma \mathfrak{A} \\
    &  & & & &~-~& (3m-1) D \alpha \mathfrak{B}&~+~& (n-3m) D \beta \mathfrak{A} \\
        &  & & & & & &~-~& 4m E \alpha \mathfrak{A} \\
\end{array}
\right\rbrace =0.
\end{equation*} 
\begin{equation*}
\text{etc.}
\end{equation*}
The rule, how these formulas are continued, becomes clear from the inspection of the above equation more quickly than it can be described by words. While descending down the column the coefficients are decreased by the difference $m+n$; but while proceeding horizontally the differences will  continuously be increased by the difference $n-1$.

\paragraph*{§208}

 Therefore, this way the theory of recurring series is extended, since we discovered the previously unknown equations for the coefficients even for the cases, in which not only the denominator of the fraction was any power, but also the numerator consists of any arbitrary number of terms, to detect which equations induction alone did not suffice. But except for the many applications of recurring series  we already explained, they are very useful to find the sums of certain series approximately; we exhibited a specimen of this already in the first chapter of this book, where we transformed the series  into another one which often consists of a finite number of terms by means of the substitution $x=\frac{y}{1+ny}$ . And the method could have been extended further, if  other functions were substituted for $x$. Since then the law of progression of series,  which had to be substituted for the power of $x$, was not sufficiently clear,  it seemed advisable to mention this generalization just here, after the mentioned law had already been completely discovered. Nevertheless, considering this with more attention, we learn that the same task can be done without this law of progression only by using the method we used here to investigate the law.
 
 \paragraph*{§209}
Therefore, let an arbitrary series be propounded, i.e.

\begin{equation*}
s=A+Bx+Cx^2+Dx^3+Ex^4+Fx^5+\text{etc.}
\end{equation*}  
which we want to  transform into another one, each term of which are fractions whose denominators proceed according to powers of a formula of this kind

\begin{equation*}
\alpha +\beta x+ \gamma x^2+ \delta x^3+ \text{etc.}
\end{equation*}
To start from simpler cases, let us put that

\begin{equation*}
s=\frac{\mathfrak{A}}{\alpha + \beta x}+\frac{\mathfrak{B}x}{(\alpha +\beta x)^2}+\frac{\mathfrak{C}x^2}{(\alpha +\beta x)^3}+\frac{\mathfrak{D}x^3}{(\alpha +\beta x)^4}+\text{etc.};
\end{equation*}
Having equated the series to this expression, multiply by $\alpha +\beta x$ everywhere and it will be

\begin{equation*}
\left.
\renewcommand{\arraystretch}{1,4}
\setlength{\arraycolsep}{0.0mm}
\begin{array}{lllllllllllllllll}
A \alpha &~+~& B \alpha &x&~+~& C\alpha&x^2 &~+~& D\alpha &x^3& ~+~&+\text{etc.} \\
 & & A \beta &x&~+~& b \beta &  &~+~& C \beta & & ~+~&+\text{etc.} \\
\end{array}
\right\rbrace = \mathfrak{A}+\frac{\mathfrak{B}x}{\alpha +\beta x}+\frac{\mathfrak{C}x^2}{(\alpha +\beta x)^2}+\text{etc.}
\end{equation*} 
Put $\mathfrak{A}=A \alpha$ and let 

\begin{alignat*}{9}
&A \beta &&+ B \alpha &&= A^{\prime} \\
& B \beta &&+ C \alpha &&= B^{\prime} \\
&C \beta &&+D \alpha &&=C^{\prime} \\
&D \beta &&+E \alpha &&=D^{\prime} \\
& &&\text{etc.};
\end{alignat*}
having divided by $x$, it will be

\begin{equation*}
A^{\prime}+B^{\prime}x+C^{\prime}x^2+D^{\prime}x^3+\text{etc.}=\frac{\mathfrak{B}}{\alpha +\beta x}+\frac{\mathfrak{C}x}{(\alpha +\beta x)^2}+\frac{\mathfrak{D}x^2}{(\alpha +\beta x)^3}+\text{etc.}
\end{equation*}
Multiply by $\alpha +\beta x$ again and having put 

\begin{alignat*}{9}
&A^{\prime} \beta &&+ B^{\prime} \alpha &&= A^{\prime \prime} \\
& B^{\prime} \beta &&+ C^{\prime} \alpha &&= B^{\prime \prime} \\
&C^{\prime} \beta &&+D^{\prime} \alpha &&=C^{\prime \prime} \\
& &&\text{etc.}
\end{alignat*}
it will be

\begin{equation*}
A^{\prime} \alpha +A^{\prime \prime} x+B^{\prime \prime}x^2+C^{\prime \prime}x^3+ \text{etc.}= \mathfrak{B}+\frac{\mathfrak{C}x}{\alpha +\beta x}+\frac{\mathfrak{D}x^2}{(\alpha + \beta x)^2}+\text{etc.}
\end{equation*}
Therefore, let $\mathfrak{B}=A^{\prime}\alpha$; and arguing exactly as before, if 

\begin{alignat*}{9}
&A^{\prime \prime } \beta &&+ B^{\prime \prime} \alpha &&= A^{\prime \prime \prime} \qquad \qquad &&A^{\prime \prime \prime}\beta &&+B^{\prime \prime \prime}\alpha = A^{\prime \prime \prime \prime} \\
&B^{\prime \prime } \beta &&+ C^{\prime \prime} \alpha &&= B^{\prime \prime \prime} \qquad \qquad &&B^{\prime \prime \prime}\beta &&+C^{\prime \prime \prime}\alpha = B^{\prime \prime \prime \prime} \\
&C^{\prime \prime } \beta &&+ D^{\prime \prime} \alpha &&= C^{\prime \prime \prime} \qquad \qquad &&C^{\prime \prime \prime}\beta &&+D^{\prime \prime \prime}\alpha = C^{\prime \prime \prime \prime} \\
& &&\text{etc.} && && &&\text{etc.},
\end{alignat*}
it will be $\mathfrak{C}= A^{\prime \prime}\alpha$, $\mathfrak{D}=A^{\prime \prime \prime} \alpha$, $\mathfrak{E}=A^{\prime \prime \prime \prime} \alpha$; therefore, the sum of the propounded series will be expressed as follows

\begin{equation*}
s=\frac{A \alpha}{\alpha + \beta x}+\frac{A^{\prime}\alpha x}{(\alpha + \beta x)^2}+\frac{A^{\prime}\alpha x^2}{(\alpha +\beta x)^3}+\frac{A^{\prime \prime \prime}\alpha x^3}{(\alpha +\beta x)^4}+\text{etc.}
\end{equation*}
This same series would have resulted from the substitution

\begin{equation*}
\frac{x}{\alpha +\beta x}=y \quad \text{or} \quad x=\frac{\alpha y}{1- \beta y}.
\end{equation*}

\paragraph*{§210}

This transformation is applied with the greatest success, if the propounded series $A+Bx+Cx^2+\text{etc.}$ was of such a nature that it is finally confounded with a recurring series or, even better, a geometric series resulting from the fraction $\frac{P}{\alpha +\beta x}$. For, then the values $A^{\prime}$, $B^{\prime}$, $C^{\prime}$, $D^{\prime}$ etc. will finally vanish and hence  the letters $A^{\prime \prime}$, $A^{\prime \prime \prime}$, $A^{\prime \prime \prime \prime}$ etc. will even more constitute a highly converging series.\\[2mm]
In like manner we will be able to use trinomial and any polynomial denominators which will have an extraordinary use, if the propounded series is finally confounded with a recurring series. Therefore, having propounded the series

\begin{equation*}
s=A+Bx+Cx^2+Dx^3+Ex^4+Fx^5+\text{etc.},
\end{equation*} 
set

\begin{equation*}
s= \frac{\mathfrak{A}+\mathfrak{B}x}{\alpha +\beta x+ \gamma x^2}+ \frac{\mathfrak{A}^{\prime}x^2+\mathfrak{B}^{\prime}x^3}{(\alpha +\beta x+ \gamma x^2)^2}+\frac{\mathfrak{A}^{\prime \prime}x^4+\mathfrak{B}^{\prime \prime}x^5}{(\alpha +\beta x+ \gamma x^2)^3}+\frac{\mathfrak{A}^{\prime \prime \prime}x^6+\mathfrak{B}^{\prime \prime \prime}x^7}{(\alpha +\beta x+ \gamma x^2)^4}+\text{etc.}
\end{equation*}
Multiply by $\alpha +\beta x+ \gamma x^2$ everywhere and, having put

\begin{alignat*}{19}
& A \gamma &&+ B \beta &&+ C \alpha &&=A^{\prime} \quad && \text{and} \quad && \mathfrak{A}&&= A \alpha \\
& B \gamma &&+ C \beta &&+ D \alpha &&=B^{\prime} \quad && \text{and} \quad && \mathfrak{B}&&= A \beta +B \alpha \\
& C \gamma &&+ D \beta &&+ E \alpha &&=C^{\prime}, 
\end{alignat*}
having divided by $xx$, an equation similar to the first will result, i.e.

\begin{equation*}
A^{\prime}+B^{\prime}x+C^{\prime}x^2+D^{\prime}x^3+E^{\prime}x^4+\text{etc.}
\end{equation*}
\begin{equation*}
= \frac{\mathfrak{A}^{\prime}+\mathfrak{B}^{\prime}x}{\alpha +\beta x+ \gamma xx}+\frac{\mathfrak{A}^{\prime \prime}+\mathfrak{B}^{\prime \prime}x}{(\alpha +\beta x+ \gamma xx)^2}+\frac{\mathfrak{A}^{\prime \prime \prime}+\mathfrak{B}^{\prime \prime \prime}x}{(\alpha +\beta x+ \gamma xx)^3}+\text{etc.}
\end{equation*}
Therefore, if the operation is done as before by putting

\begin{alignat*}{19}
& A^{\prime} \gamma &&+ B^{\prime} \beta &&+ C^{\prime} \alpha &&=A^{\prime \prime} \quad && \text{and} \quad && \mathfrak{A}^{\prime}&&= A^{\prime} \alpha \\
& B^{\prime} \gamma &&+ C^{\prime} \beta &&+ D^{\prime} \alpha &&=B^{\prime \prime} \quad && \text{and} \quad && \mathfrak{B}&&= A^{\prime} \beta +B^{\prime} \alpha \\
& C^{\prime} \gamma &&+ D^{\prime} \beta &&+ E^{\prime} \alpha &&=C^{\prime \prime} \\
& && &&\text{etc.}
\end{alignat*}
and further

\begin{alignat*}{19}
& A^{\prime \prime} \gamma &&+ B^{\prime \prime} \beta &&+ C^{\prime \prime} \alpha &&=A^{\prime \prime \prime}  \quad && \text{and} \quad && \mathfrak{A}^{\prime \prime}&&= A^{\prime \prime} \alpha \\
& B^{\prime \prime} \gamma &&+ C^{\prime \prime} \beta &&+ D^{\prime \prime} \alpha &&=B^{\prime \prime \prime} \quad && \text{and} \quad && \mathfrak{B}&&= A^{\prime \prime} \beta +B^{\prime \prime} \alpha \\
& C^{\prime \prime} \gamma &&+ D^{\prime \prime} \beta &&+ E^{\prime \prime} \alpha &&=C^{\prime \prime \prime} \\
& && && \text{etc.}
\end{alignat*}
and by investigating  further values in this manner, it will be

\begin{equation*}
s= \frac{A \alpha +(A \beta +b \alpha)x}{\alpha + \beta x +\gamma xx}+\frac{(A^{\prime}\alpha +(A^{\prime}\beta +B^{\prime}\alpha))x^2}{(\alpha +\beta x +\gamma xx)^2}+\frac{(A^{\prime \prime}\alpha +(A^{\prime \prime}\beta +B^{\prime \prime}\alpha))x^4}{(\alpha +\beta x +\gamma xx)^3}+\text{etc.}
\end{equation*}

\paragraph*{§211}

If one puts $x=1$, what can be done without loss of generality, because $\alpha$, $\beta$, $\gamma$ can be taken arbitrarily, and it was

\begin{equation*}
s=A+B+C+D+E+F+G+\text{etc.},
\end{equation*}
having successively put the following values

\begin{alignat*}{19}
&A \gamma &&+B \beta && +C \alpha &&=A^{\prime} \qquad \qquad &&A^{\prime}\gamma &&+B^{\prime}\beta &&+C^{\prime} \alpha &&=A^{\prime \prime} \\
&B \gamma &&+C \beta && +D \alpha &&=B^{\prime} \qquad \qquad &&B^{\prime}\gamma &&+C^{\prime}\beta &&+D^{\prime} \alpha &&=B^{\prime \prime} \qquad \text{and so forth}\\
&C \gamma &&+D \beta && +E \alpha &&=C^{\prime} \qquad \qquad &&C^{\prime}\gamma &&+D^{\prime}\beta &&+E^{\prime} \alpha &&=C^{\prime \prime} \\
& && &&\text{etc.} && && &&\text{etc.}
\end{alignat*}
but  for the sake of brevity one puts

\begin{equation*}
\alpha + \beta +\gamma =m,
\end{equation*}
one will obtain the sum of the propounded series expressed this way

\begin{equation*}
s=\left\lbrace
\renewcommand{\arraystretch}{2,5}
\setlength{\arraycolsep}{0.0mm}
\begin{array}{rrrrrrrrrrrrrrrrrr}
(\alpha +\beta) & \bigg(\dfrac{A}{m}&~+&\dfrac{A^{\prime}}{m^2}&~+~&\dfrac{A^{\prime \prime}}{m^3}&~+~&\dfrac{A^{\prime \prime \prime}}{m^4}&~+~& \text{etc.} \bigg) \\
 +\alpha & \bigg(\dfrac{B}{m}&~+&\dfrac{B^{\prime}}{m^2}&~+~&\dfrac{B^{\prime \prime}}{m^3}&~+~&\dfrac{B^{\prime \prime \prime}}{m^4}&~+~& \text{etc.} \bigg) \\
\end{array}
\right\rbrace 
\end{equation*} 

\paragraph*{§212}

The same denominators consisting of more terms can be taken, and since the operation is easily understood from the preceding, let us only expand the case for the polynomial of degree three. Therefore, let

\begin{equation*}
s=A+B+C+D+E+F+G+\text{etc.}
\end{equation*}
Find the following values

\begin{alignat*}{19}
&A \delta &&+B \gamma && +C \beta &&+ D \alpha &&= A^{\prime} \\
&B \delta &&+C \gamma && +D \beta &&+ E \alpha &&= B^{\prime} \\
&C \delta &&+D \gamma && +E \beta &&+ F \alpha &&= C^{\prime} \\
&         &&          && \text{etc.}
\end{alignat*}
\begin{alignat*}{19}
&A^{\prime} \delta &&+B^{\prime} \gamma && +C^{\prime} \beta &&+ D^{\prime} \alpha &&= A^{\prime \prime} \\
&B^{\prime} \delta &&+C^{\prime} \gamma && +D^{\prime} \beta &&+ E^{\prime} \alpha &&= B^{\prime \prime} \\
&C^{\prime} \delta &&+D^{\prime} \gamma && +E^{\prime} \beta &&+ F^{\prime} \alpha &&= C^{\prime \prime} \\
&         &&          && \text{etc.}
\end{alignat*}
\begin{alignat*}{19}
&A^{\prime \prime} \delta &&+B^{\prime \prime} \gamma && +C^{\prime\prime} \beta &&+ D^{\prime \prime} \alpha &&= A^{\prime \prime \prime} \\
&B^{\prime \prime} \delta &&+C^{\prime \prime} \gamma && +D^{\prime \prime} \beta &&+ E^{\prime \prime} \alpha &&= B^{\prime \prime \prime} \\
&C^{\prime \prime} \delta &&+D^{\prime \prime} \gamma && +E^{\prime \prime} \beta &&+ F^{\prime \prime} \alpha &&= C^{\prime \prime \prime} \\
&         &&          && \text{etc.}
\end{alignat*}
But then let $\alpha +\beta + \gamma +\delta=m$ and it will be

\begin{equation*}
s=\left\lbrace
\renewcommand{\arraystretch}{2,5}
\setlength{\arraycolsep}{0.0mm}
\begin{array}{rrrrrrrrrrrrrrrrr}
(\alpha +\beta +\gamma) & \bigg(\dfrac{A}{m}&~+&\dfrac{A^{\prime}}{m^2}&~+~&\dfrac{A^{\prime \prime}}{m^3}&~+~&\dfrac{A^{\prime \prime \prime}}{m^4}&~+~& \text{etc.} \bigg) \\
 (\alpha +\beta) & \bigg(\dfrac{B}{m}&~+&\dfrac{B^{\prime}}{m^2}&~+~&\dfrac{B^{\prime \prime}}{m^3}&~+~&\dfrac{B^{\prime \prime \prime}}{m^4}&~+~& \text{etc.} \bigg) \\
 +\alpha & \bigg(\dfrac{C}{m}&~+&\dfrac{C^{\prime}}{m^2}&~+~&\dfrac{C^{\prime \prime}}{m^3}&~+~&\dfrac{C^{\prime \prime \prime}}{m^4}&~+~& \text{etc.} \bigg) \\
\end{array}
\right\rbrace 
\end{equation*} 
whence at the same time the progression, if even more parts are attributed to the denominator $m$, is most clearly seen.

\paragraph*{§213}

And it is not necessary at all that the denominators of the fractions, to which we reduced the sum of the series, are powers of the same formula

\begin{equation*}
\alpha +\beta x+ \gamma x^2+\text{etc.},
\end{equation*}
but this can be varied in each term. In order to clarify this, let us at first only take two terms and assume that the series

\begin{equation*}
s=A+Bx+Cx^2+Dx^3+Ex^4+Fx^5+\text{etc.}
\end{equation*}
is converted into this series of fractions

\begin{equation*}
s=\frac{\mathfrak{A}}{\alpha +\beta x}+\frac{\mathfrak{A}^{\prime}x}{(\alpha +\beta x)(\alpha^{\prime}+\beta^{\prime}x)}+\frac{\mathfrak{A}^{\prime \prime}x^2}{(\alpha +\beta x)(\alpha^{\prime}+\beta^{\prime}x)(\alpha^{\prime \prime}+\beta^{\prime \prime}x)}+\text{etc.}
\end{equation*}
At first, multiply both sides by $\alpha +\beta x$ and put

\begin{alignat*}{9}
&A \beta &&+B \alpha &&=A^{\prime} \\
& B\beta &&+ C \alpha &&=B^{\prime} \quad \text{and} \quad \mathfrak{A}= A \alpha \\
&C \beta &&+ D \alpha &&=C^{\prime} \\
& && \text{etc.}
\end{alignat*}
and, having divided by $x$, it will be

\begin{equation*}
A^{\prime}+B^{\prime}+C^{\prime}x^2+D^{\prime}x^3+\text{etc.}= \frac{\mathfrak{A}}{\alpha^{\prime}+\beta^{\prime}x}+\frac{\mathfrak{A}^{\prime \prime}x}{(\alpha^{\prime}+\beta^{\prime}x)(\alpha^{\prime \prime}+\beta^{\prime \prime}x)}+\text{etc.}
\end{equation*}
Further, in like manner by multiplying by $\alpha^{\prime}+\beta^{\prime}x$ and then by $\alpha^{\prime \prime}+\beta^{\prime \prime}x$ and so forth, if one sets

\begin{alignat*}{19}
& A^{\prime}\beta^{\prime} &&+B^{\prime}\alpha^{\prime} &&= A^{\prime \prime} \qquad && A^{\prime \prime}\beta^{\prime \prime} &&+B^{\prime \prime}\alpha^{\prime \prime} &&= A^{\prime \prime \prime} \qquad &&A^{\prime \prime \prime} \beta^{\prime \prime \prime} &&+B^{\prime \prime \prime }\alpha^{\prime \prime \prime} &&=A^{\prime \prime \prime \prime} \\
& B^{\prime}\beta^{\prime} &&+C^{\prime}\alpha^{\prime} &&= B^{\prime \prime} \qquad && B^{\prime \prime}\beta^{\prime \prime} &&+C^{\prime \prime}\alpha^{\prime \prime} &&= B^{\prime \prime \prime} \qquad &&B^{\prime \prime \prime} \beta^{\prime \prime \prime} &&+C^{\prime \prime \prime }\alpha^{\prime \prime \prime} &&=B^{\prime \prime \prime \prime} \qquad \text{etc.} \\
& C^{\prime}\beta^{\prime} &&+D^{\prime}\alpha^{\prime} &&= C^{\prime \prime} \qquad && C^{\prime \prime}\beta^{\prime \prime} &&+D^{\prime \prime}\alpha^{\prime \prime} &&= C^{\prime \prime \prime} \qquad &&C^{\prime \prime \prime} \beta^{\prime \prime \prime} &&+D^{\prime \prime \prime }\alpha^{\prime \prime \prime} &&=C^{\prime \prime \prime \prime} \qquad  \\
& &&\text{etc.} && && && \text{etc.} && && && \text{etc.}
\end{alignat*}
it will be

\begin{equation*}
\mathfrak{A}^{\prime}= A^{\prime}\alpha^{\prime}, \quad \mathfrak{A}^{\prime \prime}=A^{\prime \prime}\alpha^{\prime \prime}, \quad \mathfrak{A}^{\prime \prime \prime}= A^{\prime \prime \prime} \alpha^{\prime \prime \prime} \quad \text{etc.}
\end{equation*}
and hence the propounded series will be converted into this one

\begin{equation*}
s=\frac{A \alpha}{\alpha +\beta x}+\frac{A^{\prime}\alpha^{\prime}x}{(\alpha +\beta x)(\alpha^{\prime}+\beta^{\prime}x)}+\frac{A^{\prime \prime}\alpha^{\prime \prime}x}{(\alpha +\beta x)(\alpha^{\prime}+\beta^{\prime}x)(\alpha^{\prime \prime}+\beta^{\prime \prime}x)}+\text{etc.},
\end{equation*}
where the values $\alpha$, $\beta$, $\alpha^{\prime}$, $\beta^{\prime}$, $\alpha^{\prime \prime}$, $\beta^{\prime \prime}$ etc. are arbitrary, but can be taken in such a way for each case that this new series is highly convergent.

\paragraph*{§214}

Let us apply this also to trinomial factors and, having propounded an arbitrary series $s=A+B+C+D+E+F+G+\text{etc.}$, let

\begin{alignat*}{19}
&A \gamma &&+B \beta &&+C \alpha &&=A^{\prime} \qquad \qquad && A^{\prime}\gamma^{\prime}&&+ B^{\prime}\beta^{\prime} &&+C^{\prime}\alpha^{\prime} &&= A^{\prime \prime} \\
&B \gamma &&+C \beta &&+D \alpha &&=B^{\prime} \quad && B^{\prime}\gamma^{\prime}&&+ C^{\prime}\beta^{\prime} &&+D^{\prime}\alpha^{\prime} &&= B^{\prime \prime} \\
&C \gamma &&+D \beta &&+E \alpha &&=C^{\prime} \quad && C^{\prime}\gamma^{\prime}&&+ D^{\prime}\beta^{\prime} &&+E^{\prime}\alpha^{\prime} &&= C^{\prime \prime} \\
& && &&\text{etc.} && && && &&\text{etc.}
\end{alignat*}
\begin{alignat*}{19}
&A^{\prime \prime} \gamma^{\prime} &&+B^{\prime \prime} \beta^{\prime \prime} &&+C^{\prime \prime} \alpha^{\prime \prime} &&=A^{\prime \prime \prime} \qquad \qquad && A^{\prime \prime \prime}\gamma^{\prime \prime \prime}&&+ B^{\prime \prime \prime}\beta^{\prime \prime \prime} &&+C^{\prime \prime \prime}\alpha^{\prime \prime \prime} &&= A^{\prime \prime \prime \prime} \\
&B^{\prime \prime} \gamma^{\prime \prime} &&+C^{\prime \prime} \beta^{\prime \prime} &&+D^{\prime \prime} \alpha^{\prime \prime} &&=B^{\prime \prime \prime} \quad && B^{\prime \prime \prime}\gamma^{\prime \prime \prime}&&+ C^{\prime \prime \prime}\beta^{\prime \prime \prime} &&+D^{\prime \prime \prime}\alpha^{\prime \prime \prime} &&= B^{\prime \prime \prime \prime} \\
&C^{\prime \prime} \gamma^{\prime \prime} &&+D^{\prime \prime} \beta^{\prime \prime} &&+E^{\prime \prime} \alpha^{\prime \prime} &&=C^{\prime \prime \prime} \quad && C^{\prime \prime \prime}\gamma^{\prime \prime \prime}&&+ D^{\prime \prime \prime}\beta^{\prime \prime \prime} &&+E^{\prime \prime \prime}\alpha^{\prime \prime \prime} &&= C^{\prime \prime \prime \prime} \\
& && &&\text{etc.} && && && &&\text{etc.}
\end{alignat*}
Further, for the sake of brevity put

\begin{alignat*}{9}
&\alpha &&+\beta &&+ \gamma &&=m \\
&\alpha^{\prime} &&+\beta^{\prime} &&+ \gamma^{\prime} &&=m^{\prime} \\
&\alpha^{\prime \prime} &&+\beta^{\prime \prime } &&+ \gamma^{\prime \prime} &&=m^{\prime \prime} \\
&\alpha^{\prime \prime \prime} &&+\beta^{\prime \prime \prime} &&+ \gamma^{\prime \prime \prime} &&=m^{\prime \prime \prime} \\
& && &&\text{etc.}
\end{alignat*}
and the sum of the propounded series will be

\begin{equation*}
s=\frac{\alpha(A+B)}{m}+\frac{\alpha^{\prime}(A^{\prime}+B^{\prime})}{mm^{\prime}}+\frac{\alpha^{\prime \prime}(A^{\prime \prime}+B^{\prime \prime})}{mm^{\prime}m^{\prime \prime}}+\frac{\alpha^{\prime \prime\prime}(A^{\prime \prime \prime}+B^{\prime\prime \prime})}{mm^{\prime}m^{\prime \prime}m^{\prime \prime \prime}}+\text{etc.}
\end{equation*}
\begin{equation*}
+\frac{\beta A}{m}+\frac{\beta^{\prime}A^{\prime}}{mm^{\prime}}+\frac{\beta^{\prime \prime}A^{\prime \prime}}{mm^{\prime}m^{\prime \prime}}+\frac{\beta^{\prime \prime \prime}A^{\prime \prime \prime}}{mm^{\prime}m^{\prime \prime}m^{\prime \prime \prime}}+\text{etc.}
\end{equation*}

\paragraph*{§215}

Since these formulas extend so far that their use can be seen less clearly, let us restrict our considerations to the case of the transformation given in § 213 and let $x=-1$ that one has this series

\begin{equation*}
s=A-B+C-D+E-F+G-\text{etc.}
\end{equation*}
and set

\begin{alignat*}{19}
&B&&-A &&=A^{\prime} \qquad && B^{\prime}&&-2A^{\prime}&&=A^{\prime \prime} \qquad &&B^{\prime \prime}&&-3A^{\prime \prime}&&=A^{\prime \prime \prime} \qquad &&B^{\prime \prime \prime}&&-4A^{\prime \prime \prime}&&=A^{\prime \prime \prime \prime} \\
&C&&-B &&=B^{\prime} \qquad && C^{\prime}&&-2B^{\prime}&&=B^{\prime \prime} \qquad &&C^{\prime \prime}&&-3B^{\prime \prime}&&=B^{\prime \prime \prime} \qquad &&C^{\prime \prime \prime}&&-4B^{\prime \prime \prime}&&=B^{\prime \prime \prime \prime} \\
&D&&-C &&=C^{\prime} \qquad && D^{\prime}&&-2C^{\prime}&&=C^{\prime \prime} \qquad &&D^{\prime \prime}&&-3C^{\prime \prime}&&=C^{\prime \prime \prime} \qquad &&D^{\prime \prime \prime}&&-4C^{\prime \prime \prime}&&=C^{\prime \prime \prime \prime} \\
& && &&\text{etc.} \qquad && && && \text{etc.} && && &&\text{etc.} && && &&\text{etc.} \\
\end{alignat*}
Having found these values, the sum of the propounded series will be equal to the following series

\begin{equation*}
s=\frac{A}{2}-\frac{A^{\prime}}{2 \cdot 3}+\frac{A^{\prime \prime}}{2 \cdot 3 \cdot 4}-\frac{A^{\prime \prime \prime}}{2 \cdot 3 \cdot 4 \cdot 5}+\frac{A^{\prime \prime \prime \prime }}{2 \cdot 3 \cdot 4 \cdot 5 \cdot 6}- \text{etc.}
\end{equation*}
Therefore,  any propounded series can be transformed into innumerable others equal to it, among which without any doubt a most convergent series will be found, by means of which the propounded sum can be found approximately. 

\paragraph*{§216}

But let us return to the invention of series whose law of progression is revealed by differential calculus. Therefore, because this was already achieved for algebraic quantities, let us proceed to transcendental functions and let the series equal to this logarithm be in question

\begin{equation*}
s= \log (1+ \alpha x+ \beta x^2 +\gamma x^3 +\delta x^4 +\varepsilon x^5+ \text{etc.});
\end{equation*}
assume that the series in question reads

\begin{equation*}
s= \mathfrak{A}x+\mathfrak{B}x^2+ \mathfrak{C}x^3+\mathfrak{D}x^4+\mathfrak{E}x^5+\mathfrak{F}x^6+\text{etc.}
\end{equation*}
Therefore, because from the differentiation of the first equation it follows

\begin{equation*}
\frac{ds}{dx}= \frac{\alpha +2 \beta x +3 \gamma x^2+ 4 \delta x^3 +5 \varepsilon x^4+\text{etc.}}{1+\alpha x+ \beta x^2+ \gamma x^3+ \delta x^4+\varepsilon x^5+\text{etc.}},
\end{equation*}
it will be

\begin{equation*}
(1+ \alpha x+ \beta x^2+\gamma x^3+ \delta x^4+\text{etc.})\frac{ds}{dx}=\alpha +2 \beta x +3 \gamma x^2+ 4 \delta x^3+ \text{etc.}
\end{equation*}
But since from the assumed equation

\begin{equation*}
\frac{ds}{dx}= \mathfrak{A}+2 \mathfrak{B}x+ 3 \mathfrak{C}x^2+ 4 \mathfrak{D}x^3+5 \mathfrak{E}x^4+\text{etc.},
\end{equation*}
having done the substitution, this equation results

\begin{equation*}
\renewcommand{\arraystretch}{1,4}
\setlength{\arraycolsep}{0.0mm}
\begin{array}{rrrrrrrrrrrrrrrrr}
 \mathfrak{A}&~+~& 2\mathfrak{B}&x&~+~&3 \mathfrak{C}&x^2&~+~&4 \mathfrak{D}&x^3&~+~& 5 \mathfrak{E}&x^4&~+~& \text{etc.} \\
  &~+~& \mathfrak{A}\alpha & &~+~&2 \mathfrak{B}\alpha & &~+~&3 \mathfrak{C}\alpha & &~+~& 4 \mathfrak{D}\alpha & &~+~& \text{etc.} \\
   & &  & &~+~& \mathfrak{A}\beta & &~+~&2 \mathfrak{B}\beta & &~+~& 3\mathfrak{C}\beta & &~+~& \text{etc.} \\
      & &  & & &  & &~+~& \mathfrak{A}\gamma & &~+~& 2\mathfrak{B}\gamma & &~+~& \text{etc.} \\
          & &  & & &  & & &   & &~+~& \mathfrak{A}\delta & &~+~& \text{etc.} \\ \hline
          =\alpha &~+~& 2\beta &x&~+~&3 \gamma &x^2&~+~&4 \delta &x^3&~+~& 5 \varepsilon &x^4&~+~& \text{etc.} \\
\end{array}
\end{equation*} 
From it one obtains the following determinations

\begin{alignat*}{9}
& \mathfrak{A}&&=\alpha \\[3mm]
& \mathfrak{B}&&= -\frac{1}{2}\mathfrak{A}\alpha &&+\beta \\[3mm]
& \mathfrak{C}&&=- \frac{2}{3}\mathfrak{B} \alpha &&-\frac{1}{3}\mathfrak{A}\beta &&+\gamma \\[3mm]
& \mathfrak{D}&&=-\frac{3}{4}\mathfrak{C}\alpha &&-\frac{2}{4}\mathfrak{B}\beta && - \frac{1}{4}\mathfrak{A}\gamma &&+ \delta \\[3mm]
&\mathfrak{E}&&=-\frac{4}{5}\mathfrak{D}\alpha &&-\frac{3}{5}\mathfrak{C}\beta &&-\frac{2}{5}\mathfrak{B}\gamma &&-\frac{1}{5}\mathfrak{A}\delta &&+\text{etc.} \\
& && && &&\text{etc.}
\end{alignat*}

\paragraph*{§217}

Now, let this exponential quantity be propounded

\begin{equation*}
s= e^{\alpha x+ \beta x^2+ \gamma x^3+\delta x^4+\varepsilon x^5+\text{etc.}}
\end{equation*}
in which $e$ denotes the number whose hyperbolic logarithm is $=1$, and assume this series in question

\begin{equation*}
s= 1+ \mathfrak{A}x+\mathfrak{B}x^2+\mathfrak{C}x^3+\mathfrak{D}x^4+\mathfrak{E}x^5+\text{etc.}
\end{equation*}
For, from the case $x=0$ it is plain that the first term must be $1$. Therefore, since by taking logarithms,

\begin{equation*}
\log s= \alpha +\beta x^2+ \gamma x^3+ \delta x^4+\varepsilon x^5+\zeta x^6+ \text{etc.},
\end{equation*}
 having taken the differentials, it will be

\begin{equation*}
\frac{ds}{dx}= s(\alpha +2 \beta x+ 3 \gamma x^2+ 4 \delta x^3+ 5 \varepsilon x^4 +\text{etc.})
\end{equation*}
But from the assumed equation it will be

\begin{equation*}
\renewcommand{\arraystretch}{1,4}
\setlength{\arraycolsep}{0.0mm}
\begin{array}{rrrrrrrrrrrrrrrrr}
 \dfrac{ds}{dx}&~=~& \mathfrak{A}&~+~&2 \mathfrak{B}&x&~+~&3 \mathfrak{C}&x^2&~+~&4\mathfrak{D}&x^3&~+~&5\mathfrak{E}&x^4&~+~&\text{etc.} \\  \cline{2-17}
 &~=~& \alpha&~+~& \mathfrak{A}\alpha&x&~+~& \mathfrak{B}\alpha&x^2&~+~&\mathfrak{C}\alpha&x^3&~+~&\mathfrak{D}\alpha&x^4&~+~&\text{etc.} \\
  & & &~+~& 2 \beta& &~+~& 2 \mathfrak{A}\beta& &~+~&2\mathfrak{C}\beta& &~+~&2\mathfrak{C}\beta& &~+~&\text{etc.} \\
    & & & & & &~+~& 3 \gamma & &~+~&3\mathfrak{A}\gamma& &~+~&3\mathfrak{B}\gamma& &~+~&\text{etc.} \\
     & & & & & & &  & &~+~&4\delta & &~+~&4\mathfrak{A}\delta& &~+~&\text{etc.} \\
     & & & & & & &  & & & & &~+~&5\varepsilon& &~+~&\text{etc.} ,
\end{array}
\end{equation*} 
from which the following determinations of the letters $\mathfrak{A}$, $\mathfrak{B}$, $\mathfrak{C}$, $\mathfrak{D}$ etc. result

\begin{alignat*}{19}
&\mathfrak{A}&&= \alpha \\[3mm]
&\mathfrak{B}&&= \beta &&+\frac{1}{2}\mathfrak{A}\beta \\[3mm]
&\mathfrak{C}&&= \gamma &&+ \frac{2}{3}\mathfrak{A}\beta &&+ \frac{1}{3}\mathfrak{B}\alpha \\[3mm]
&\mathfrak{D}&&= \delta &&+\frac{3}{4}\mathfrak{A}\gamma &&+\frac{2}{4}\mathfrak{B}\beta &&+\frac{1}{4}\mathfrak{C}\alpha \\[3mm]
&\mathfrak{E}&&=\varepsilon &&+\frac{4}{5}\mathfrak{A}\delta &&+\frac{3}{5}\mathfrak{B}\gamma &&+\frac{2}{5}\mathfrak{C}\beta &&+\frac{1}{5}\mathfrak{D}\alpha \\
& && && &&\text{etc.}
\end{alignat*}

\paragraph*{§218}

If  the arc, whose sine or cosine is in question, is expressed by a binomial or polynomial or even an infinite series, this way one can even express its sine and cosine by means of an infinite series. But to do this the most convenient way, it does not suffice to consider only the first differentials, but it is necessary to use the differentials of the second order. Therefore, let

\begin{equation*}
s= \sin (\alpha x+ \beta x^2+ \gamma x^3+ \delta x^4+\varepsilon x^5+\text{etc.})
\end{equation*}
and assume the series in question to be

\begin{equation*}
s=\mathfrak{A}x+\mathfrak{B}x^2+ \mathfrak{C}x^3+\mathfrak{D}x^4+\mathfrak{E}x^5+\text{etc.}
\end{equation*}
For, it is plain that the first term vanishes; but since one has to descend to the second differentials, the coefficient $\mathfrak{A}$ also has to be determined from elsewhere, which will happen, if we put $x$ to be  infinitely small. For then, because of the arc $=\alpha x$, the sine itself will become equal to it and it will therefore be $\mathfrak{A}=\alpha$. Now, for the sake of brevity let us put

\begin{equation*}
z=\alpha x+ \beta x^2+ \gamma x^3+ \text{etc.},
\end{equation*}
that  $s= \sin z$; by differentiating it will be $ds=dz \cos z$ and by differentiating again it will be $dds=ddz \cos z -dz^2\sin z$. Therefore, since  $\sin z=s$ and $\cos z = \frac{ds}{dz}$, it will be

\begin{equation*}
dds= \frac{dsddz}{dz}-sdz^2 \quad \text{and} \quad dz dds+sdz^3=dsddz.
\end{equation*}	

\paragraph*{§219}

Let us put that the arc $z$ is only expressed by a binomial and 

\begin{equation*}
z= \alpha x+ \beta x^2;
\end{equation*}
it will be

\begin{equation*}
dz= (\alpha +2 \beta x)dx
\end{equation*}
and, having put $dx$ to be constant,

\begin{equation*}
ddz = 2 \beta dx^2
\end{equation*}
and

\begin{equation*}
dz^3=(\alpha^3+6 \alpha^2 \beta x+ 12 \alpha \beta^2 x^2+8 \beta^3 x^3)dx^3.
\end{equation*}
Further, because of $s=\mathfrak{A}x+\mathfrak{B}x^2+\mathfrak{C}x^3+\mathfrak{D}x^4+\text{etc.}$, it will be

\begin{equation*}
\frac{ds}{dx}=\mathfrak{A}+2 \mathfrak{B}x+ 3 \mathfrak{C}x^2+ 4 \mathfrak{D}x^3+\text{etc.}
\end{equation*}
and

\begin{equation*}
\frac{dds}{dx^2}= 2 \mathfrak{B}+6 \mathfrak{C}x+12 \mathfrak{D}x^2+\text{etc.}
\end{equation*}
Having substituted these values in the differential equation, it will be

\begin{equation*}
\renewcommand{\arraystretch}{2,5}
\setlength{\arraycolsep}{0.0mm}
\begin{array}{rrrrrrrrrrrrrrrrr}
\dfrac{dzdds}{dx^3}&~=~& 1 \cdot 2 \mathfrak{B}\alpha &~+~& 2 \cdot 3 \mathfrak{C}\alpha &x&~+~& 3 \cdot 4 \mathfrak{D}\alpha &x^2&~+~& 4 \cdot 5 \mathfrak{E}\alpha &x^3&~+~& 5 \cdot 6 \mathfrak{F}\alpha &x^4&~+~&\text{etc.} \\
 & &  &~+~& 2 \cdot 1 \cdot 2 \mathfrak{B}\beta & &~+~& 2 \cdot 2 \cdot 3 \mathfrak{C}\beta & &~+~& 2 \cdot 3 \cdot 4 \mathfrak{D}\beta & &~+~& 2 \cdot 4 \cdot 5 \mathfrak{E}\beta & &~+~&\text{etc.} \\
\dfrac{sdz^3}{dx^3} &~=~ &  &~+~&  \mathfrak{A}\alpha^3 & &~+~&  \mathfrak{B}\alpha^3 & &~+~&  \mathfrak{C}\alpha^3 & &~+~&  \mathfrak{D}\alpha^3 & &~+~&\text{etc.} \\
  & &  & &   & &~+~&  6\mathfrak{A}\alpha^2 \beta & &~+~&  6\mathfrak{B}\alpha^2 \beta & &~+~&  6\mathfrak{C}\alpha^2 \beta & &~+~&\text{etc.} \\
   & &  & &   & & &  & &~+~&  12\mathfrak{A}\alpha \beta^2 & &~+~&  12\mathfrak{B}\beta^2 & &~+~&\text{etc.} \\
     & &  & &   & & &  & & &   & &~+~&  8\mathfrak{A}\beta^3 & &~+~&\text{etc.} \\ \hline
     \dfrac{dsddz}{dx^3}&~=~&  2 \mathfrak{A}\beta &~+~& 4 \mathfrak{B}\beta & &~+~&  6 \mathfrak{C}\beta & &~+~& 8 \mathfrak{D}\beta & &~+~& 10 \mathfrak{E}\beta & &~+~&\text{etc.} \\
\end{array}
\end{equation*} 
Therefore, the coefficients will be defined the following way:

\begin{alignat*}{9}
&\mathfrak{B}&&~=~ && &&\frac{2 \mathfrak{A}\beta}{2 \alpha} \\[3mm]
&\mathfrak{C}&&~=~ && &&~~0 &&-\frac{\mathfrak{A}\alpha^2}{2 \cdot 3} \\[3mm]
&\mathfrak{D}&&~=~ &&~-~&& \frac{2 \mathfrak{C}\beta}{4 \alpha}&&-\frac{6 \mathfrak{A}\alpha \beta}{3 \cdot 4}&&- \frac{\mathfrak{B}\alpha^2}{3 \cdot 4} \\[3mm]
&\mathfrak{E}&&~=~&&~-~&&\frac{4 \mathfrak{D}\beta}{5 \alpha}&&-\frac{12\mathfrak{A}\beta^2}{4 \cdot 5}&&- \frac{6 \mathfrak{B}\alpha \beta}{4 \cdot 5}&&- \frac{\mathfrak{C}\alpha^2}{4 \cdot 5} \\[3mm]
&\mathfrak{F}&&~=~&&~-~&&\frac{6 \mathfrak{E}\beta}{6 \alpha}&&-\frac{8 \mathfrak{A}\beta^3}{5 \cdot 6 \alpha}&&-\frac{12 \mathfrak{B}\beta \beta}{5 \cdot 6}&&-\frac{6 \mathfrak{C}\alpha \beta}{5 \cdot 6}&&-\frac{\mathfrak{D}\alpha^2}{5 \cdot 6} \\[3mm]
&\mathfrak{G}&&~=~&&~-~&&\frac{8 \mathfrak{F}\beta}{7 \alpha}&&- \frac{8 \mathfrak{B}\beta^3}{6 \cdot 7\alpha}&&-\frac{12 \mathfrak{C}\beta \beta}{6 \cdot 7}&&-\frac{6 \mathfrak{D}\alpha \beta}{6 \cdot 7}&&-\frac{\mathfrak{E}\alpha^2}{6 \cdot 7} \\
& && && && &&&& \text{etc.}
\end{alignat*}
Having found these values, it will be

\begin{equation*}
\sin (\alpha x+ \beta x^2)=\mathfrak{A}x+ \mathfrak{B}x^2+\mathfrak{C}x^3+\mathfrak{D}x^4+\text{etc.}
\end{equation*}
while $\mathfrak{A}=\alpha$.

\paragraph*{§220}

In like manner, the cosine of any angle is converted into a series; but since an arc is very rarely expressed by a polynomial, let us show the use of differential equations for the invention of the series for the cosine of the arc $x$. Therefore, let $s= \cos x$ and assume

\begin{equation*}
s=1- \mathfrak{A}x^2+ \mathfrak{B}x^4-\mathfrak{C}x^6+\mathfrak{D}x^8-\text{etc.}
\end{equation*}
Since $ds= -dx \sin x$ and $dds=-dx^2 \cos x=-sdx^2$, it will be

\begin{equation*}
dds+sdx^2=0;
\end{equation*}
after the substitution it will be

\begin{equation*}
\renewcommand{\arraystretch}{1,4}
\setlength{\arraycolsep}{0.0mm}
\begin{array}{rrrrrrrrrrrrrrrrr}
\dfrac{dds}{dx^2}&~=~&-1 \cdot 2 \mathfrak{A}&~+~&3 \cdot 4 \mathfrak{B}x^2 &- 5 \cdot 6 \mathfrak{C}x^4&~+~&7 \cdot 8 \mathfrak{D}x^6&~-~&\text{etc.} \\
s&~=~&1&~-~&\mathfrak{A} x^2 &-\mathfrak{B} x^4&~-~&\mathfrak{C} x^6&~+~&\text{etc.} \\
\end{array}
\end{equation*} 
and by comparing the coefficients it follows

\begin{alignat*}{9}
&\mathfrak{A}&&=\frac{1}{1 \cdot 2} \\[3mm]
& \mathfrak{B}&&=\frac{\mathfrak{A}}{3 \cdot 4}&&=\frac{1}{1 \cdot 2 \cdot 3 \cdot 4} \\[3mm]
&\mathfrak{C}&&=\frac{\mathfrak{B}}{5 \cdot 6}&&=\frac{1}{1 \cdot 2 \cdot 3 \cdots 6}\\[3mm]
&\mathfrak{D}&&= \frac{\mathfrak{C}}{7 \cdot 8}&&= \frac{1}{1 \cdot 2 \cdot 3 \cdots 8} \\
& && &&\text{etc.}
\end{alignat*}
Therefore, it is plain, what we demonstrated in more detail above already, that 

\begin{equation*}
\cos x= 1- \frac{x^2}{1 \cdot 2}+\frac{x^4}{1 \cdot 2 \cdot 3 \cdot 4}-\frac{x^6}{1 \cdot 2 \cdot 3 \cdots 6}+\frac{x^8}{1 \cdot 2 \cdot 3 \cdots 8}-\text{etc.};
\end{equation*}
the first series for the sine, having put $\beta=0$ and $\alpha=1$, will give

\begin{equation*}
\sin x =\frac{x}{1}-\frac{x^3}{1 \cdot 2 \cdot 3}+\frac{x^5}{1 \cdot 2 \cdot 3 \cdot 4 \cdot 5}-\frac{x^7}{1 \cdot 2 \cdot 3 \cdots 7}+\frac{x^9}{1 \cdot 2 \cdot 3 \cdots 9}-\text{etc.}
\end{equation*}

\paragraph*{§221}

From the well-known series for the sine and cosine the series for the tangent, cotangent, secant, cosecant of a certain angle are deduced. For, the tangent results, if the sine is divided by the cosine, the cotangent, if the cosine is divided by the sine, the secant, if the radius $1$ is divided by the cosine, and the cosecant, if the radius is divided by the sine. But the series result from these divisions seem to be most irregular; but, with the exception of the series exhibiting the secant, all remaining the others  can  be reduced to a simple law by means of the Bernoulli numbers $\mathfrak{A}$, $\mathfrak{B}$, $\mathfrak{C}$, $\mathfrak{D}$ etc. For, since we found above (§ 127) that

\begin{equation*}
\frac{\mathfrak{A}u^2}{1 \cdot 2}+\frac{\mathfrak{B}u^4}{1 \cdot 2 \cdot 3 \cdot 4}+\frac{\mathfrak{C}u^6}{1 \cdot 2 \cdot 3 \cdots 6}+\frac{\mathfrak{D}u^8}{1 \cdot 2 \cdot 3 \cdots 8}+\text{etc.}=1-\frac{u}{2}\cot \frac{1}{2}u,
\end{equation*}
having put $\frac{1}{2}u=x$, it will be 

\begin{equation*}
\cot x=\frac{1}{x}-\frac{2^2 \mathfrak{A}x}{1 \cdot 2}-\frac{2^4\mathfrak{B}x^3}{1 \cdot 2 \cdot 3 \cdot 4}-\frac{2^6\mathfrak{C}x^5}{1 \cdot 2 \cdot 3 \cdots 6}
- \frac{2^8\mathfrak{D}x^7}{1 \cdot 2 \cdot 3 \cdots 8}-\text{etc.},
\end{equation*}
and if one puts $\frac{1}{2}x$ for $x$, it will be

\begin{equation*}
\cot \frac{1}{2}x=\frac{2}{x}-\frac{2 \mathfrak{A}x}{1 \cdot 2}-\frac{2\mathfrak{B}x^3}{1 \cdot 2 \cdot 3 \cdot 4}-\frac{2\mathfrak{C}x^5}{1 \cdot 2 \cdot 3 \cdots 6}
- \frac{2^\mathfrak{D}x^7}{1 \cdot 2 \cdot 3 \cdots 8}-\text{etc.},
\end{equation*}

\paragraph*{§222}

But hence the tangent of any arc will be expressed by means of an infinite series the following way. Since

\begin{equation*}
\tan 2x= \frac{2 \tan x}{1- \tan x},
\end{equation*}
it will be

\begin{equation*}
\cot 2x = \frac{1}{2 \tan x}-\frac{\tan x}{2}=\frac{1}{2}\cot x -\frac{1}{2}\tan x
\end{equation*}
and hence

\begin{equation*}
\tan x = \cot x -2 \cot 2x.
\end{equation*}
Therefore, since

\begin{equation*}
\cot x = \frac{1}{x}-\frac{2^2 \mathfrak{A}x}{1 \cdot 2}-\frac{2^4 \mathfrak{B}x^3}{1 \cdot 2 \cdot 3 \cdot 4}-\frac{2^6\mathfrak{C}x^5}{1 \cdot 2 \cdots 6}-\frac{2^8 \mathfrak{D}x^7}{1 \cdot 2 \cdots 8}-\text{etc.},
\end{equation*}
\begin{equation*}
2\cot 2x = \frac{1}{x}-\frac{2^4 \mathfrak{A}x}{1 \cdot 2}-\frac{2^8 \mathfrak{B}x^3}{1 \cdot 2 \cdot 3 \cdot 4}-\frac{2^{12}\mathfrak{C}x^5}{1 \cdot 2 \cdots 6}-\frac{2^{16} \mathfrak{D}x^7}{1 \cdot 2 \cdots 8}-\text{etc.},
\end{equation*}
subtracting this series from the first, it will be

\begin{equation*}
\tan x = \frac{2^2(2^2-1)\mathfrak{A}x}{1 \cdot 2}+\frac{2^4(2^4-1)\mathfrak{B}x^3}{1 \cdot 2 \cdot 3 \cdot 4}+\frac{2^6(2^6-1)\mathfrak{C}x^5}{1 \cdot 2  \cdots 6}+\frac{2^8(2^8-1)\mathfrak{D}x^7}{1 \cdot 2  \cdots 8}+\text{etc.}
\end{equation*}
Therefore, if the numbers $A$, $B$, $C$, $D$ etc. found in § 182 are introduced here, it will be

\begin{equation*}
\tan x = \frac{2Ax}{1 \cdot 2}+\frac{2^3 B x^3}{1 \cdot 2 \cdot 3 \cdot 4}+\frac{2^5 C x^5}{1 \cdot 2 \cdot 3 \cdot 6}+\frac{2^7 D x^7}{1 \cdot 2 \cdots 8}+\text{etc.}
\end{equation*}

\paragraph*{§223}

But the cosecant will be found the following way. Since 

\begin{equation*}
\cot x= \tan x +2 \cot 2x =\frac{1}{\cot x}+2 \cot 2x,
\end{equation*}
it will be

\begin{equation*}
\cot^2 x= 2 \cot x \cot 2x+1
\end{equation*}
and having extracted the root

\begin{equation*}
\cot x = \cot 2x +\csc 2x,
\end{equation*}
whence

\begin{equation*}
\csc 2x= \cot x- \cot 2x
\end{equation*}
and having put $x$ for $2x$

\begin{equation*}
\csc x =\cot \frac{1}{2}x- \cot x.
\end{equation*}
Therefore, because we have the cotangent, i.e.

\begin{alignat*}{9}
&\cot \frac{1}{2}x &&=\frac{2}{x}&&-\frac{2 \mathfrak{A}x}{1 \cdot 2}&&-\frac{2 \mathfrak{B}x^3}{1 \cdot 2 \cdot 3 \cdot 4}&&-\frac{2 \mathfrak{C}x^5}{1 \cdot 2 \cdots 6}&&-\text{etc.} \\[3mm]
&\cot x &&=\frac{1}{x}&&-\frac{2^2 \mathfrak{A}x}{1 \cdot 2}&&-\frac{2 ^4\mathfrak{B}x^3}{1 \cdot 2 \cdot 3 \cdot 4}&&-\frac{2^5 \mathfrak{C}x^5}{1 \cdot 2 \cdots 6}&&-\text{etc.},
\end{alignat*}
having subtracted this series from the first, it will be

\begin{equation*}
\csc x=\frac{1}{x}+\frac{2(2-1)\mathfrak{A}x}{1 \cdot 2}+\frac{2(2^3-1)\mathfrak{B}x^3}{1 \cdot 2 \cdot 3 \cdot 4}+\frac{2(2^5-1)\mathfrak{C}x^5}{1 \cdot 2 \cdots 6}+\text{etc.}
\end{equation*}

\paragraph*{§224}

But the secant cannot be expressed by means of these Bernoulli numbers, but it requires other numbers which enter the sums of the odd powers of the reciprocals. For, if one puts

\begin{alignat*}{9}
&1 &&-~~\frac{1}{3}&&+~~\frac{1}{5}&&-~~\frac{1}{7}&&+~~\frac{1}{9}&&-\text{etc.} &&= \quad \quad \alpha &&\cdot \frac{\pi}{2^2} \\[3mm]
&1 &&-~\frac{1}{3^3}&&+~\frac{1}{5^3}&&-~\frac{1}{7^3}&&+~\frac{1}{9^3}&&-\text{etc.} &&= \quad ~\frac{\beta}{1 \cdot 2} &&\cdot \frac{\pi^3}{2^4} \\[3mm]
&1 &&-~\frac{1}{3^5}&&+~\frac{1}{5^5}&&-~\frac{1}{7^5}&&+~\frac{1}{9^5}&&-\text{etc.} &&=  \frac{\gamma}{1 \cdot 2 \cdot 3 \cdot 4 } &&\cdot \frac{\pi^5}{2^6} \\[3mm]
&1 &&-~\frac{1}{3^7}&&+~\frac{1}{5^7}&&-~\frac{1}{7^7}&&+~\frac{1}{9^7}&&-\text{etc.} &&= ~\frac{\delta}{1 \cdot 2 \cdots 6} &&\cdot \frac{\pi^7}{2^8} \\[3mm]
&1 &&-~\frac{1}{3^9}&&+~\frac{1}{5^9}&&-~\frac{1}{7^9}&&+~\frac{1}{9^9}&&-\text{etc.} &&= ~\frac{\varepsilon}{1 \cdot 2 \cdots 8} &&\cdot \frac{\pi^9}{2^{10}} \\[3mm]
&1 &&-\frac{1}{3^{11}}&&+\frac{1}{5^{11}}&&-\frac{1}{7^{11}}&&+\frac{1}{9^{11}}&&-\text{etc.} &&= \frac{\zeta}{1 \cdot 2 \cdots 10} &&\cdot \frac{\pi^{11}}{2^{12}} \\
& && && && && &&\text{etc.},
\end{alignat*}
it will be 

\begin{alignat*}{9}
&\alpha &&=1 \\
&\beta &&=1 \\
&\gamma &&=5 \\
& \delta &&=61 \\
& \varepsilon &&= 1385 \\
& \zeta &&= 50521 \\
& \eta &&= 2702765 \\
& \theta && = 199360981 \\
&\iota &&= 19391512145 \\
&\varkappa &&= 20404879675441 
& && \text{etc.}
\end{alignat*}
and from these values one will obtain

\begin{equation*}
\sec x=\alpha +\frac{\beta}{1 \cdot 2}xx+\frac{\gamma}{1 \cdot 2 \cdot 3 \cdot 4}x^4+\frac{\delta}{1 \cdot 2 \cdots 6}x^6+\frac{\varepsilon}{1 \cdot 2 \cdots 8}x^8+\text{etc.}
\end{equation*}

\paragraph*{§225}

To show the connection of this series to the numbers $\alpha$, $\beta$, $\gamma$, $\delta$ etc., let us consider the series treated above [§ 33]

\begin{equation*}
\frac{\pi}{n \sin \frac{m}{n}\pi}=\frac{1}{m}+\frac{1}{n-m}-\frac{1}{m+n}-\frac{1}{2n-m}+\frac{1}{2n+m}+\frac{1}{3n-m}-\text{etc.}
\end{equation*}
Put $m= \frac{1}{2}n-k$ and it will be

\begin{equation*}
\frac{\pi}{2n \cos \frac{k}{n}}=\frac{1}{n-2k}+\frac{1}{n+2k}-\frac{1}{3n-2k}-\frac{1}{3n+2k}+\frac{1}{5n-2k}+\text{etc.}
\end{equation*}
Let $\frac{k \pi}{n}=x$ or $k \pi =nx$; it will be

\begin{equation*}
\frac{\pi}{2n}\sec x= \frac{\pi}{n \pi -2n x}+\frac{\pi}{n \pi+2n x}-\frac{\pi}{3n \pi- 2n x}-\frac{\pi}{3n \pi+2nx}+\frac{\pi}{5n \pi-2nx}+\text{etc.}
\end{equation*}
or 

\begin{equation*}
\sec x= \frac{2}{\pi -2x}+\frac{2}{\pi+2x}-\frac{2}{3 \pi-2x}-\frac{2}{3 \pi+2x}+\frac{2}{5 \pi -2x}+\text{etc.}
\end{equation*}
or

\begin{equation*}
\sec x= \frac{4 \pi}{\pi^2 -4x^2}-\frac{4 \cdot 3 \pi}{9 \pi^2- 4xx}+\frac{4 \cdot 5 \pi}{25\pi^2 -4xx}-\frac{4 \cdot 7 \pi}{49 \pi^2 -x^2}+\text{etc.}
\end{equation*}
If now each terms is converted into series, it will be

\begin{alignat*}{9}
&\sec x &&= ~\frac{4}{\pi}&& \bigg(1&&-~\frac{1}{3}&&+~\frac{1}{5}&&-~\frac{1}{7}&&+~\frac{1}{9}&&-\text{etc.} \bigg)\\[3mm]
&  &&+ \frac{2^4 x^2}{\pi^3}&& \bigg(1&&-\frac{1}{3^3}&&+\frac{1}{5^3}&&-\frac{1}{7^3}&&+\frac{1}{9^3}&&-\text{etc.} \bigg)\\[3mm]
&  &&+ \frac{2^6 x^4}{\pi^5}&& \bigg(1&&-\frac{1}{3^5}&&+\frac{1}{5^5}&&-\frac{1}{7^5}&&+\frac{1}{9^5}&&-\text{etc.} \bigg)\\
& && && && &&\text{etc.};
\end{alignat*}
if  the values assigned above are substituted for these series, the same series we gave for the secant will result.

\paragraph*{§226}

Therefore, at the same time the law is plain, according to which the numbers $\alpha$, $\beta$, $\gamma$, $\delta$ etc. appearing in the expressions of the sums of the odd powers, proceed. For, since

\begin{equation*}
\sec x= \frac{1}{\cos x}= \alpha +\frac{\beta}{1 \cdot 2}x^2+\frac{\gamma}{1 \cdot 2 \cdot 3 \cdot 4}x^4+\frac{\delta }{1 \cdot 2 \cdots 6}x^6+\text{etc.},
\end{equation*}
it is necessary that this series is equal to the fraction

\begin{equation*}
\frac{1}{1-\frac{xx}{1 \cdot 2}+\frac{x^4}{1 \cdot 2 \cdot 3 \cdot 4}-\frac{x^6}{1 \cdot 2 \cdots 6}+\frac{x^6}{1 \cdot 2 \cdots 8}-\text{etc.}};
\end{equation*}
therefore, having equated the two expressions, it will be

\begin{alignat*}{9}
&1= \alpha &&+\frac{\beta}{1 \cdot 2}x^2&&+\frac{\gamma}{1 \cdot 2 \cdot 3 \cdot 4}x^4&&+\frac{\delta}{~~1 \cdot 2 \cdots 6~~}x^6&&+\frac{\varepsilon}{~~~1 \cdot 2 \cdots 8~~~}x^8&&+\text{etc.} \\[3mm]
&  &&-\frac{\alpha}{1 \cdot 2}&&-\frac{\beta}{1 \cdot 2 \cdot 1 \cdot 2}&&-\frac{\gamma}{1 \cdot 2 \cdot 1 \cdots 4}&&-\frac{\delta}{~1 \cdot 2 \cdot 1 \cdots 6~}&&-\text{etc.} \\[3mm]
&  && &&+\frac{\alpha}{1 \cdot 2 \cdot 3 \cdot 4}&&+\frac{\beta}{1 \cdots 4 \cdot 1 \cdot 2}&&+\frac{\gamma}{1 \cdots 4 \cdot 1 \cdots 4}&&-\text{etc.} \\[3mm]
&  && && &&-\frac{\alpha}{~~1 \cdot 2 \cdots 6~~~}&&+\frac{\beta}{1 \cdots 6 \cdot 1 \cdots 2}&&-\text{etc.} \\[3mm]
&  && && && &&+\frac{\alpha}{~~~~1 \cdot 2 \cdots 8~~~~}&&+\text{etc.}, \\
\end{alignat*}
whence these equations follow

\begin{alignat*}{9}
&\alpha &&=1 \\
&\beta &&=\frac{2 \cdot 1}{1 \cdot 2}\alpha \\[3mm]
&\gamma &&= \frac{4 \cdot 3}{1 \cdot 2}\beta &&-\frac{4 \cdot 3 \cdot 2 \cdot 1}{1 \cdot 2 \cdot 3 \cdot 4}\alpha \\[3mm]
&\delta &&=\frac{6 \cdot 5}{1 \cdot 2}\gamma &&-\frac{6 \cdot 5 \cdot 4 \cdot 3}{1 \cdot 2 \cdot 3 \cdot 4}\beta &&+\frac{6 \cdots 1}{1 \cdots 6}\alpha \\[3mm]
&\varepsilon &&=\frac{8 \cdot 7}{1 \cdot 2}\delta &&-\frac{8 \cdot 7 \cdot 6 \cdot 5}{1 \cdot 2 \cdot 3 \cdot 4}\gamma &&+\frac{8 \cdots 3}{1 \cdots 6}\beta &&-\frac{8 \cdots 1}{1 \cdots 8}\alpha \\
& && && &&\text{etc.}
\end{alignat*}
And from these formulas the values of these letters were found which we exhibited in § 224 and by means of which the sums of the series contained in this form

\begin{equation*}
1-\frac{1}{3^n}+\frac{1}{5^n}-\frac{1}{7^n}+\frac{1}{9^n}-\text{etc.},
\end{equation*}
if $n$ was an odd number, can be expressed.


\title{On the use of Differential Calculus in the resolution of Equations
\footnote{
Original title: "'De Usu Calculi Differentialis in Aequationibus resolvendis"', first published as part of the book \textit{Institutiones calculi differentialis cum eius usu in analysi finitorum ac doctrina serierum}, 1755, reprinted in \textit{Opera Omnia}: Series 1, Volume 10, pp. 422 - 445, Eneström-Number E212, translated by: Alexander Aycock for the "`Euler-Kreis Mainz"'}}
\author{Leonhard Euler}
\date{ }
\maketitle

\paragraph*{§227}

That an equation can be reduced to the nature of functions has already been demonstrated  above. For, let $y$ denote any function of $x$; if one puts $y=0$,  this form contains  all finite equations, may they be algebraic or transcendental. But the equation $y=0$ is said to be solved, if that value $x$ is found which substituted in the function $y$ actually renders it equal to zero. But in most cases many of such values $x$ exist, which are called the roots of the equation $y=0$. Therefore, if we assume the numbers $f$, $g$, $h$, $i$ etc. to be roots of the equation $y=0$, the function $y$ will be of such a nature that, if in it  either $f$ or $g$ or $h$ or etc. is substituted for $x$, it indeed is $y=0$.

\paragraph*{§228} 

Therefore, since the function $f$ vanishes, if in it  one writes $f$ or $x+(f-x)$ instead of $x$, where $f$ is a root of the equation $y=0$,  by the properties of functions we demonstrated above [§ 48], it will be

\begin{equation*}
0= y+\frac{(f-x)dy}{dx}+\frac{(f-x)^2ddy}{2dx^2}+\frac{(f-x)^3 d^3y}{6dx^3}+\text{etc.},
\end{equation*}
from which equation the value of the root $f$ is determined in such a way that, whatever was substituted for $x$ and hence the value of the quantities $y$, $\frac{dy}{dx}$, $\frac{ddy}{2dx^2}$ etc. were substituted, always the equation expressing the true value of $f$ results. To see this more clearly, let

\begin{equation*}
y=x^3-2x^2+3x-4;
\end{equation*}
it will be

\begin{equation*}
\frac{dy}{dx}=3xx-4x+3, \quad \frac{ddy}{2dx^2}=3x-2 \quad \text{and} \quad \frac{d^3y}{6dx^3}=1.
\end{equation*}
Having substituted these values, this equation results

\begin{equation*}
0=x^3-2x^2+3x-4+(f-x)(3xx-4x+3)+(f-x)^2(3x-2)+(f-x)^3
\end{equation*}
or having actually performed the multiplications

\begin{equation*}
f^3-2ff+3f-4=0;
\end{equation*}
of course, the same equation as the propounded one results, which  therefore has the same roots.

\paragraph*{§229}

But although this way one does not get to a new equation, from which the value of the root $f$ can be determined in an easier way, nevertheless  extraordinary auxiliary theorems teaching how to find the roots can be deduced from this. For, if  a value already very close to a certain root was assumed for $x$ such that $f-x$ is a very small quantity, then the terms of the equation

\begin{equation*}
0= y+ \frac{(f-x)dy}{dx}+\frac{(f-x)^2ddy}{2dx^2}+\frac{(f-x)^3d^3y}{6dx^3}+\text{etc.}
\end{equation*}
will converge very rapidly and therefore this expression will not deviate much from the true value, if  only  the first two initial terms are considered. Therefore,  if  a value already close to a certain root of the equation $y=0$ was assumed for $x$, it will approximately be

\begin{equation*}
0= y+\frac{(f-x)dy}{dx} \quad \text{or} \quad f=x- \frac{ydx}{dy},
\end{equation*}
from which formula a, even though not true, but nevertheless very good approximate value of the root $f$ will be found, which, if it is substituted for $x$ again, will yield a even better value for $f$ and so one will  continuously get closer to the true value of the root $f$. 

\paragraph*{§230}

Therefore, first the root of all powers of any number can be calculated this way. For, let the number $a^n+b$ be propounded and let its root of the power $n$ to be extracted. Put $x^n=a^n+b$ or $x^n-a^n-b=0$ that  $y=x^n-a^n-b$; it will be

\begin{equation*}
\frac{dy}{dx}= nx^{n-1}, \quad \frac{ddy}{2dx^2}=\frac{n(n-1)}{1 \cdot 2}x^{n-2}, \quad \frac{d^3 y}{6 dx^3}=\frac{n(n-1)(n-2)}{1 \cdot 2 \cdot 3}x^{n-3} \quad \text{etc.}
\end{equation*}
Therefore, if the root in question is put $=f$, so that $f= \sqrt[n]{a^n+b}$, it will be

\begin{equation*}
0=x^n-a^n-b+n(f-x)x^{n-1}+\frac{n(n-1)}{1 \cdot 2}(f-x)^2x^{n-2}+\text{etc.}
\end{equation*}
Therefore, if  one takes a number already coming close to the value of the root $f$ in question for $x$, which will be achieved by putting $x=a$, if $b$ was such a small number that $a^n+b <(a+1)^n$, it will approximately be $b=na^{n-1}(f-a)$ and hence

\begin{equation*}
f= a+ \frac{b}{na^{n-1}},
\end{equation*}
whence a much better approximation of the value of the root will be found. But if we want to take also the third term, so that

\begin{equation*}
b= na^{n-1}(f-a)+\frac{n(n-1)}{1 \cdot 2}a^{n-2}(f-a)^2,
\end{equation*}
it will be

\begin{equation*}
(f-a)^2=- \frac{2a}{n-1}(f-a)+\frac{2b}{n(n-1)a^{n-2}}
\end{equation*}
and hence

\begin{equation*}
f=a-\frac{a}{n-1}\pm \sqrt{\frac{aa}{(n-1)^2}+\frac{2b}{n(n-1)a^{n-2}}}
\end{equation*}
or

\begin{equation*}
f= \frac{(n-2)a+\sqrt{aa+2(n-1)b:na^{n-2}}}{n-1}.
\end{equation*}
Therefore, by means of the extraction of the square root an even closer value of the root $f$ will be found.

\subsection*{Example}

\textit{Let us find the square root of any number $c$ or let $xx-c=y$.}\\[2mm]
Therefore, put the number very close to the root $=a$ and $b=c-aa$; because of $aa+b=c$ and since  $n=2$, the first formula will become $f= a+\frac{c-aa}{2a}=\frac{c+aa}{2a}$; the other gives $f=\sqrt{c}$ which is the root in question itself. Therefore, because the root approximately is $=\frac{c+aa}{2a}$, write that value for $a$ and  $f=\frac{cc+6aac+a^4}{4a(c+aa)}$ will be an even better approximation. For the sake of an example let $c=5$, from the first formula it will be $f=\frac{5}{2a}+\frac{a}{2}$. Therefore, put $a=2$, it will be $f=2.25$; now put $a=2.25$, it will be $f=2.236111$; further, set $a=2.236111$, it will be $f=2.2360679$ which value already hardly deviates from the value.

\paragraph*{§231}

But in like manner, the root of any equation can be found approximately by means of the equation $f=x- \frac{ydx}{dy}$, after having assumed a value differing hardly from a certain root of the equation for $x$, of course. To find a value of this kind for $x$, successively substitute various values for $x$ and from them chose the one which minimizes the function $y$, which means, which indicates the value closest to zero. So, if 

\begin{equation*}
y=x^3-2xx+3x-4
\end{equation*}
\begin{alignat*}{9}
&\text{having put} \quad &&x=0 \quad \text{then} \quad && y=-4 \\
&                        &&x=1 \quad               \quad &&y=-2 \\
&                        &&x=2 \quad               \quad &&y=+2, 
\end{alignat*}
whence we see that the root is contained within the values boundaries $1$ and $2$ of $x$. Therefore, since $\frac{dy}{dx}=3xx-4x+3$, in order to find the root $f$ of the equation $x^3-2xx+3x-4=0$ one has to use this equation
\begin{equation*}
f= x- \frac{y dx}{dy}= x-\frac{x^3-2xx+3x-4}{3xx-4x+3}.
\end{equation*}
Therefore, let $x=1$; it will be $f=1+\frac{2}{2}=2$. Now put $x=2$; it will be $f=2- \frac{2}{7}=\frac{12}{7}.$ Therefore, let $x=\frac{12}{7}$; it will be
$f=\frac{12}{7}-\frac{104}{1701}=\frac{2812}{1701}=1.653$. If we want to proceed, it is advisable to use logarithms.\\
Therefore, put $x=1.653$ and it will be

\begin{center}
\begin{tabular}{rclcccccc}
$\log x^1$ &$~=~$ & \multicolumn{1}{r|}{$0.2182729$} & $x^1$ & $~=~$ & $1.653000$ \\
$\log x^2$ &$~=~$ & \multicolumn{1}{r|}{$0.4365458$} & $x^2$ & $~=~$ & $2.732409$ \\
$\log x^3$ &$~=~$ & \multicolumn{1}{r|}{$0.6548187$} & $x^3$ & $~=~$ & $4.516673$ \\
$x^3$   &$~=~$& $4.516673$ \\
$3x$   &$~=~$& $4.959000$ \\ \cline{3-3}
\end{tabular}
\end{center}
And hence

\begin{center}
\begin{tabular}{rclrcr}
$x^3+3x$ &$~=~$& $9.475673$ & $3xx+3$ &$~=~$ & $11.197227$ \\
$2xx+4$ &$~=~$& $9.464818$ & $4x$ &$~=~$ & $6.612000$ \\ \cline{3-3} \cline{6-6}
num. &$~=~$& $0.010855$ & den. &$~=~$ & $4.585227$ \\
$\log$ num. &$~=~$& $8.0356298$  \\
$\log$ den. &$~=~$& $0.6613608$ & $x$ &$~=~$ & $1.6553000$ \\
$\log$ fract.  &$~=~$& $7.3742690$ & fraction &$~=~$ & $0.002367$ \\ \cline{6-6} \\
   & &   & $f$ &$~=~$ & $1.650633$, \\ \cline{6-6} \\
\end{tabular}
\end{center}
which value already comes very close to the true one.

\paragraph*{§232}

But we will be able to deduce faster approximations from the general expression. For, because having put any function $y=0$, if the root of this equation was $x=f$, we will find that

\begin{equation*}
0= y+ \frac{(f-x)dy}{dx}+\frac{(f-x)^2ddy}{2dx^2}+\frac{(f-x)^3d^3 y}{6dx^3}+\text{etc.},
\end{equation*}
let $f-x=z$, such that the root is $f=x+z$, and put

\begin{equation*}
\frac{dy}{dx}=p, \quad \frac{dp}{dx}=q, \quad \frac{dq}{dx}=r, \quad \frac{dr}{dx}=s \quad \text{etc.};
\end{equation*}
it will be

\begin{equation*}
0= y+zp+\frac{z^2q}{2}+\frac{z^3r}{6}+\frac{z^4s}{24}+\frac{z^5t}{120}+\text{etc.};
\end{equation*}
 having taken any value for $x$ in that equation, from which value at the same time $y$, $p$, $q$, $r$, $s$ etc. are determined, the quantity $z$ must be found, having found which one will have  the root $f=x+z$ of the propounded equation $y=0$. Therefore, we will have to focus on the task how to find the value of the unknown $z$ from this equation in the most convenient way possible. 
 
\paragraph*{§233}

Assume this convergent series for $z$

\begin{equation*}
z=A+B+C+D+E+\text{etc.}
\end{equation*}
and after the substitution it will be

\begin{equation*}
\renewcommand{\arraystretch}{2,5}
\setlength{\arraycolsep}{0.0mm}
\begin{array}{rrrrrrrrrrrrrrrrrrrrrr}
y &~=~& y \\
pz &~=~& Ap &~+~& Bp &~+~& Cp &~+~& Dp &~+~& Ep &~+~& \text{etc.} \\
\dfrac{1}{2}qz^2 &~=~&  &~+~& \dfrac{1}{2}A^2q &~+~& ABq &~+~& ACq &~+~& ADq &~+~& \text{etc.} \\
 & &  & &  & &  &~+~& \dfrac{1}{2}BBq &~+~& BCq &~+~& \text{etc.} \\
\dfrac{1}{6}rz^3 &~=~&  & &  & & \dfrac{1}{6}A^3r &~+~& \dfrac{1}{2}A^2 Br &~+~& \dfrac{1}{2}A^2Cr &~+~& \text{etc.} \\
 & &  & & & &  & &  &~+~& \dfrac{1}{2}AB^2 r &~+~& \text{etc.} \\
 \dfrac{1}{24}sz^4 &~=~&  & &  & &   & & \dfrac{1}{24}A^4s &~+~& \dfrac{1}{6}A^3Bs &~+~& \text{etc.} \\
 \dfrac{1}{120}tz^5 &~=~&  & &  & &   & &  & & \dfrac{1}{120}A^5t  &~+~& \text{etc.} \\
\end{array}
\end{equation*} 
Therefore, one will obtain the following equations

\begin{alignat*}{9}
&A &&=- \frac{y}{p} \\[3mm]
&B &&= - \frac{yyq}{2p^3}\\[3mm]
&C &&=- \frac{y^3qq}{2p^5}+\frac{y^3r}{6p^4} \\[3mm]
&D &&=-\frac{5y^4q^3}{8p^7}+\frac{5y^4qr}{12 p^6}- \frac{y^4 s}{24 p^5} \\
& && \text{etc.}
\end{alignat*}
and hence it will be

\begin{equation*}
z=- \frac{y}{p}-\frac{y^2q}{2p^3}-\frac{y^3 qq}{2p^5}+\frac{y^3 r}{6p^4}-\frac{5y^4q^3}{8p^7}+\frac{5y^4 qr}{12p^6}-\frac{y^4s}{24 p^5}-\text{etc.}
\end{equation*}

\subsection*{Example}

\textit{Let this equation be propounded $x^5+2x-2=0$.}\\[2mm]
Therefore, it will be

\begin{equation*}
y=x^5+2x-2, \quad \frac{dy}{dx}=p=5x^4+2, \quad \frac{dp}{dx}=q= 20x^3,
\end{equation*}
\begin{equation*}
\frac{dq}{dx}=r= 60 x^3, \quad \frac{dr}{dx}=s=120x \quad \text{etc.}
\end{equation*}
But now  put $x=1$, since this value hardly deviates from the true root, it will be

\begin{equation*}
y=1, \quad p=7, \quad q=20, \quad r=60, \quad s=120,
\end{equation*}
whence it will be

\begin{equation*}
z=-\frac{1}{7}-\frac{10}{7^3}-\frac{200}{7^5}+\frac{10}{7^4}-\frac{5 \cdot 1000}{7^7}+\frac{500}{7^6}-\frac{5}{7^5}+\text{etc.}
\end{equation*}
or 

\begin{equation*}
z=-\frac{1}{7}-\frac{10}{7^3}-\frac{130}{7^5}-\frac{1745}{7^7}-\text{etc.},
\end{equation*}
and therefore it will be $z=0.18$ and the root $f=0.82$; if this value is again substituted for $x$, a root very close to the true one will result.

\paragraph*{§234}

Therefore, we found an infinite series, which expresses the root of any equation; but it is inconvenient that the law of progression is not obvious, and it hence is too complex and not sufficiently useful. Therefore, let us consider the same problem  in another way and try to find a more regular series expressing any root of the propounded equation.\\[2mm]
As before, let the equation $y=0$ be propounded, $y$ being an arbitrary function of $x$ and the question reduces to the definition of the value of $x$ which substituted for $x$ renders the function $y$ equal to zero. But because $y$ is a function of $x$, vice versa $x$ can be considered as function of $y$ and considering it like this the value of the function $x$ is to be found which it obtains, if the quantity $y$ vanishes. Therefore, if it is propounded to find the value of $x$ which is a root of the equation $y=0$, since $x$ goes over into $f$, if one sets $y=0$,  by the results demonstrated above [§ 67], it will be

\begin{equation*}
f= x- \frac{ydx}{dy}+\frac{y^2ddx}{2dy^2}-\frac{y^3d^3}{6dy^3}+\frac{y^4 d^4x}{24 dy^4}-\text{etc.},
\end{equation*}  
in which equation the differential $dy$ is assumed to be constant. Therefore, if one puts

\begin{equation*}
\frac{dx}{dy}=p, \quad \frac{dp}{dy}=q, \quad \frac{dq}{dy}=r, \quad \frac{dr}{dy}=s \quad \text{etc.},
\end{equation*}
having introduced these values, so that it is not necessary to consider a certain differential to be constant, it will be

\begin{equation*}
f=x-py+\frac{1}{2}qy^2-\frac{1}{6}ry^3+\frac{1}{24}sy^4-\frac{1}{120}ty^5+\text{etc.}
\end{equation*}

\paragraph*{§235}

Therefore, having attributed any value to $x$, at the same time the values of $y$ and the quantities $p$, $q$, $r$, $s$ etc. will be determined and having found these values one will have an infinite series expressing the value of the root $f$. But if the equation $y=0$ has several roots, then these values result, if different values are assumed for $x$; for, because $y$ can have the same value, even though different values are attributed to $x$, it is not surprising that the same series   can often yield several values. To avoid this ambiguity in these cases and to render the series convergent, a value already close to the value of its root, which is in question, must be assumed for $x$. For, this way  the value of  $y$ will become very small and the terms of the series will decrease rapidly, such that by taking only a few terms one will already find a sufficiently good approximation to the value for $f$. If this value is then substituted for $x$, the quantity $y$ will become a lot smaller and the series will converge a lot more rapidly and this way the root $f$ is immediately found so accurately  that the error will be very small. And hence the advantages of this expression over the one we found before are clearly seen. 

\paragraph*{§236}

Let us assume that the root of the power $n$ of any number $N$ is to be extracted. Therefore, having taken an approximate power of the exponent $n$, the propounded  number will easily be resolved into this form $N=a^n+b$. Therefore, it will be

\begin{equation*}
x^n=a^n+b \quad \text{and} \quad y= x^n-a^n-b,
\end{equation*}
whence 

\begin{alignat*}{19}
& dy &&~=~ && &&nx^{n-1}dx \quad && \text{and} \quad &&\frac{dx}{dy}&&=p &&~=~&& &&\frac{1}{nx^{n-1}} \\[3mm]
& dp &&~=~ &&~-~ &&\frac{(n-1)dx}{nx^n} \quad && \text{and} \quad &&\frac{dp}{dy}&&=q &&~=~&&~-~ &&\frac{n-1}{nnx^{2n-1}} \\[3mm]
& dq &&~=~ &&  &&\frac{(n-1)(2n-1)dx}{nnx^{2n}} \quad && \text{and} \quad &&\frac{dq}{dy}&&=r &&~=~&&  &&\frac{(n-1)(2n-1)}{n^3x^{3n-1}} \\[3mm]
& dr &&~=~ &&~-~ &&\frac{(n-1)(2n-1)(3n-1)dx}{n^3x^{3n}} \quad && \text{and} \quad &&\frac{dr}{dy}&&=s &&~=~&&~-~ &&\frac{(n-1)(2n-1)(3n-1)}{n^4x^{4n-1}} \\[3mm]
&   &&  &&  &&  && \text{etc.} \quad && &&  && &&  &&
\end{alignat*}
Now, put $x=a$ and it will be $y=-b$ and the root in question $f=\sqrt[n]{a^n+b}$ will be expressed as follows

\begin{equation*}
f=a+\frac{b}{na^{n-1}}-\frac{(n-1)bb}{n \cdot 2n a^{2n-1}}+\frac{(n-1)(2n-1)b^3}{n \cdot 2n \cdot 3n a^{3n-1}}-\frac{(n-1)(2n-1)(3n-1)b^4}{n^4 \cdot 2n \cdot 3n \cdot 4na^{4n-1}}+\text{etc.}
\end{equation*}
and so the same series results which is usually found by expansion of the binomial $(a^n+b)^{\frac{1}{n}}$.

\paragraph*{§237}

Therefore, after  the approximate root $a$ was found in the actual extraction and at the same time the remainder $b$ was found then the value $\frac{b}{an^{n-1}}$ is to be added to the root, such that a root closer to the true one is obtained. But, because of $N=a^n+b$, it will be

\begin{equation*}
a^{n-1}=\frac{N-b}{a}.
\end{equation*}
 But this way  a root larger than the correct one will be found, since the third term must be subtracted. Therefore, to find a root a lot closer to the true one, by means of division of the remainder $b$, a suitable divisor must be investigated, which we want to  assume to be

\begin{equation*}
na^{n-1}+\alpha b + \beta bb +\gamma b^3+ \text{etc.}
\end{equation*}
Therefore, since it has to be

\begin{equation*}
\frac{b}{na^{n-1}+\alpha b + \beta b^2+ \gamma b^3+ \text{etc.}}
\end{equation*}
\begin{equation*}
=\frac{b}{na^{n-1}}-\frac{(n-1)bb}{2n^2a^{2n-1}}+\frac{(n-1)(2n-1)b^3}{6n^3 a^{3n-1}}-\frac{(n-1)(2n-1)(3n-1)b^4}{24 n^4 a^{4n-1}}+\text{etc.},
\end{equation*}
after the multiplication by $na^{n-1}+\alpha b +\beta b^2+ \gamma b^3+ \text{etc.}$ it will be

\begin{alignat*}{19}
&b &&=b &&- \frac{(n-1)bb}{2na^n}&&+\frac{(n-1)(2n-1)b^3}{6n^2 a^{2n}}&&-\frac{(n-1)(2n-1)(3n-1)b^4}{24n^3a^{3n}}+\text{etc.} \\[3mm]
&  &&   &&+\frac{\alpha b^2}{na^{n-1}}&&-\frac{(n-1)\alpha b^3}{2n^2 a^{2n-1}}&&+\frac{(n-1)(2n-1)\alpha b^4}{6n^3 a^{3n-1}} \\[3mm]
&  &&   && &&+\frac{\beta b^3}{n a^{n-1}} &&- \frac{(n-1)\beta b^4}{2n^2 a^{2n-1}}\\[3mm]
& && && && &&+\frac{\gamma b^4}{na^{n-1}}
\end{alignat*}
Therefore, the following determinations are deduced

\begin{alignat*}{9}
& \alpha &&=\frac{n-1}{2a}\\[3mm]
& \beta &&=\frac{(n-1)\alpha}{2n a^n}-\frac{(n-1)(2n-1)}{6na^{n+1}}=- \frac{(n-1)(n+1)}{12 n a^{n+1}}\\[3mm]
& \gamma &&=\frac{(n-1)\beta}{2n a^n}-\frac{(n-1)(2n-1)\alpha}{6nna^{2n}}+\frac{(n-1)(2n-1)(3n-1)}{24 n^2 a^{2n+1}}=\frac{(n-1)(n+1)}{24 na^{2n+1}}.
\end{alignat*}
Therefore, the fraction to be added to the  root $a$ already found will be

\begin{equation*}
\frac{b}{na^{n-1}+\frac{(n-1)b}{2a}-\frac{(nn-1)bb}{12na^{n+1}}+\frac{(nn-1)b^3}{24 na^{2n+1}}-\text{etc.}}
\end{equation*}

\paragraph*{§238}

Therefore, if the square root of the number $N$ is to be extracted and  the approximate root was already found to be $=a$ together with the remainder $=b$, to the found root  one  has to add the quotient, which results, if the remainder $b$ is  divided by

\begin{equation*}
2a+\frac{b}{2a}-\frac{bb}{8a^3}+\frac{b^3}{16a^5}-\text{etc.}
\end{equation*}
But if the cube root must be extracted, then the remainder must be divided by

\begin{equation*}
3a^2+\frac{b}{a}-\frac{2bb}{9a^4}+\frac{b^3}{9a^7}-\text{etc.},
\end{equation*}
the application of which formulas we will seen in these examples.

\subsection*{Example 1}

\textit{To extract the square root of the number $200$.}\\[2mm]
Put $N=200$ and because the closest square is $196$, it will be $a=14$ and the remainder $b=4$, which therefore must be divided by

\begin{equation*}
28+ \frac{1}{7}-\frac{1}{7 \cdot 196}+\frac{1}{7 \cdot 196 \cdot 98},
\end{equation*}
and therefore the divisor will be $=28.142135$; if $4$ is divided by it, one will obtain a decimal fraction to be added to $14$, which will be correct up to $10$ figures and more.

\subsection*{Example 2}

\textit{To extract the cube root of the number $N=10$.}\\[2mm]
The closest cube is $8$ and the remainder is $=2$, whence $a=2$ and $b=2$ and the divisor $=12+1-\frac{1}{18}=12.9444$. Therefore, the cube root in question will approximately be $=2 \frac{2}{12.9444}=2\frac{10000}{64722}$.

\paragraph*{§239}

The series found for the root can also be considered as a recurring series resulting from a certain fraction. For, this way many terms of the series will be reduced to a lot less, namely those which constitute the numerator and the denominator of the fraction. So, having paid a little attention, one will see that it will  approximately be

\begin{equation*}
(a+b)^n =a^n \cdot \frac{a+\frac{n+1}{2}b}{a-\frac{n-1}{2}b}
\end{equation*}
and even closer

\begin{equation*}
(a+b)^n= a^n \cdot \frac{aa+\frac{n+2}{2}ab+\frac{(n+1)(n+2)}{12}bb}{aa- \frac{n-2}{2}ab+\frac{(n-1)(n-2)}{12}bb}.
\end{equation*}
In like manner, by introducing several terms even more accurate fractions can be obtained:

\begin{equation*}
(a+b)^n= a^n \cdot \frac{a^3+\frac{n+3}{2}a^2b+\frac{(n+3)(n+2)}{10}ab^2+\frac{(n+3)(n+2)(n+1)}{120}b^3}{a^3-\frac{n-3}{2}a^2b+\frac{(n-3)(n-2)}{10}ab^2-\frac{(n-3)(n-2)(n-1)}{120}b^3}.
\end{equation*}
An even more general form of this kind be exhibited, to express which conveniently let 

\begin{alignat*}{9}
&A&&= \frac{m(n+m)}{1 \cdot 2m} \qquad && \mathfrak{A}&&=\frac{m(n-m)}{1 \cdot 2 m}\\[3mm]
&B&&= \frac{(m-1)(n+m-1)}{2(2m-1)}A \qquad && \mathfrak{B}&&=\frac{(m-1)(n-m+1)}{2(2m-1)}\mathfrak{A}\\[3mm]
&C&&= \frac{(m-2)(n+m-2)}{3(2m-2)}B \qquad && \mathfrak{C}&&=\frac{(m-2)(n-m+2)}{3(2m-2)}\mathfrak{B}\\[3mm]
&D&&= \frac{(m-3)(n+m-3)}{4(2m-3)}C \qquad && \mathfrak{D}&&=\frac{(m-3)(n-m+3)}{4(2m-3)}\mathfrak{C}\\
& &&\text{etc.} && &&\text{etc.}
\end{alignat*}
But having determined these values, it will be

\begin{equation*}
(a+b)^n= a^n \cdot \frac{a^m+Aa^{m-1}b+\mathfrak{B}a^{m-2}b^2+Ca^{m-3}b^3+\text{etc.}}{a^m- \mathfrak{A}a^{m-1}b+\mathfrak{B}a^{m-2}b^2-\mathfrak{C}a^{m-3}b^3+\text{etc.}}
\end{equation*}

\paragraph*{§240}

Therefore, if  a fractional number is substituted for  $n$ here, these formulas will be very useful to extract the roots. So if any root of power $n$  of the expression $a^n+b$ has to be extracted, the following formulas can be used

\begin{equation*}
(a^n+b)^{\frac{1}{n}}=a \cdot \frac{2na^n+(n+1)b}{2na^n+(n-1)}
\end{equation*}
\begin{equation*}
(a^n+b)^{\frac{1}{n}}= a \cdot \frac{12n^2a^{2n}+6n(2n+1)a^nb+(2n+1)(n+1)bb}{12n^2a^{2n}+6n(2n-1)a^nb+(2n-1)(n-1)bb}.
\end{equation*}
But if one puts $a^n+b=N$ that $a^n=N-b$, it will be

\begin{equation*}
(a^n+b)^{\frac{1}{n}}=a \cdot \frac{2nN- (n-1)b}{2nN-(n+1)b}
\end{equation*}
\begin{equation*}
(a^n+b)^{\frac{1}{n}}a \cdot \frac{12n^2N^2-6n(2n-1)Nb+(2n-1)(n-1)bb}{12n^2N^2-6n(2n+1)Nb+(2n+1)(n+1)bb}.
\end{equation*}

\paragraph*{§241}

Therefore, the general formula for finding the root of any equation, which consists of several terms, has the same use as the usual rule of a binomial has for the resolution of the pure equations $x^n=c$, and therefore our formula goes over into that rule in this case. But if the equation was affected or even transcendental, our general expression is always applied with the same success and yields an infinite series exhibiting the value of the root. Therefore, since the resolution of such  equations is its most important application, let us demonstrate its use a little more diligently. Therefore, let this affected equation consisting of three terms be propounded

\begin{equation*}
x^n+cx=N
\end{equation*}
while $c$ and $N$ denote any given quantities. Put $x^n+cx-N=y$; it will be $dy=(nx^{n-1}+c)dx$ and hence it will be $p= \frac{1}{nx^{n-1}+c}$; then 

\begin{equation*}
dp=- \frac{n(n-1)x^{n-2}dx}{(nx^{n-1}+c)^2} \quad \text{and} \quad q=- \frac{n(n-1)x^{n-2}}{(nx^{n-1}+c)^3}.
\end{equation*}
In like manner, because of $r= \frac{dq}{dy}$, $s= \frac{dr}{dy}$ etc., one will find

\begin{small}
\begin{equation*}
r= \frac{n^2(n-1)(2n-1)x^{2n-4}-n(n-1)(n-2)cx^{n-3}}{(nx^{n-1}+c)^5}
\end{equation*}
\begin{equation*}
s= \frac{-n^3(n-1)(2n-1)(3n-1)x^{3n-6}+4n^2(n-1)(n-2)(2n-1)cx^{2n-5}-n(n-1)(n-2)(n-3)c^2x^{n-4}}{(nx^{n-1}+c)^7}
\end{equation*}

\begin{equation*}
t=\dfrac{
\left\lbrace
\renewcommand{\arraystretch}{1,4}
\setlength{\arraycolsep}{0.0mm}
\begin{array}{rrrrrrrrrrrrrrrrrrrrrr}
n^4(n-1)(2n-1)(3n-1)(4n-1)x^{4n-8}-n^3(n-1)(n-2)(2n-1)(29n-11)cx^{3n-7} \\
+n^2(n-1)(n-2)(2n-1)(11n-29)c^2x^{2n-6}-n(n-1)(n-2)(n-3)(n-4)c^3x^{n-5}
\end{array}
\right\rbrace
}{(nx^{n-1}+c)^9}
\end{equation*} 
\end{small}
\begin{center}
\text{etc.}
\end{center}
Having found these values, the root of the propounded equation will be

\begin{equation*}
f=x-py+\frac{1}{2}qyy-\frac{1}{6}ry^3+\frac{1}{24}sy^4-\frac{1}{120}ty^5+\text{etc.};
\end{equation*}
for, whatever is substituted for $x$, whence at the same time the letters $y$, $p$, $q$, $r$ etc. have explicit values, the sum of the series will become equal to the value of one single root.

\subsection*{Example 1}

\textit{Let this equation be propounded $x^3+2x=2$.} \\[2mm]
It will be $c=2$, $N=2$ and $n=3$ and $y=x^3+2x-2$. Put $x=1$; it will be $y=1$ and 
 
\begin{equation*}
p=\frac{1}{5}, \quad q=-\frac{6}{5^3}, \quad r= \frac{78}{5^5}, \quad s=-\frac{16 \cdot 90}{5^7} \quad \text{etc.} 
\end{equation*}
and the root of the equation will be
\begin{equation*}
 f=-\frac{1}{5}-\frac{3}{5^3}-\frac{13}{5^5}-\frac{60}{5^7}-\text{etc.}=0.771072.
\end{equation*}
Now put $x=0.77$, and since $y=x^3+2x-2$,
 
\begin{equation*}
 p=\frac{1}{3xx+2}, \quad q=-6p^3 x, \quad r=90xxp^5-12p^5
\end{equation*}
and 
 
\begin{equation*}
 s=- 2160p^7x^3+720 p^7x,
\end{equation*}
 using logarithms, one will have
 
\begin{center}
\begin{tabular}{rclrclccc}
$\log x^1$ &$~=~$ & \multicolumn{1}{r|}{$9.8864907$} & $x^1$ & $~=~$ & $0.77$ \\
$\log x^2$ &$~=~$ & \multicolumn{1}{r|}{$9.7729814$} & $x^2$ & $~=~$ & $0.5929$ \\
$\log x^3$ &$~=~$ & \multicolumn{1}{r|}{$9.6594721$} & $x^3$ & $~=~$ & $0.456533$ \\
    & & & $2x$ &$~=~$& $1.54$ \\ \cline{6-6}
    & & & $x^3+ 2x$ &$~=~$& $1.996533$ \\ 
 Therefore       & & & $y$ &$~=~$& $-0.003467$ 
\end{tabular}
\end{center}
And furthermore

\begin{center}
\begin{tabular}{rclcrclccc}
$\log (-y)$ &$~=~$ & $7.5399538$ & \qquad \qquad & $3xx+2$ & $~=~$ & $3.7787$ \\
$\log p$ &$~=~$ & $9.4226575$ &\qquad \qquad & $\log(3xx+2)$ & $~=~$ & $0.5773424$ \\
$\log (-py)$ &$~=~$ & $6.9626113$ & \qquad \qquad & $-py$ & $~=~$ & $0.000917511$ \\ \cline{3-7}
$\log p^3$ &$~=~$ & $8.2679725$ &\qquad \qquad  \\
$\log x$ &$~=~$ & $9.8864907$ & \qquad \qquad  \\
$\log 3$ &$~=~$ & $0.4771213$ &\qquad \qquad  \\
$\log y^2$ &$~=~$ & $5.0799076$   \\ \cline{3-3} 
$\log \left(-\frac{1}{2}qyy\right)$ &$~=~$ & $3.7114921$ & \qquad \qquad   &$-\frac{1}{2}qyy$ &$~=~$& $0.000000514$\\
\end{tabular}
\end{center} 
Therefore, the root is $f=0.770916997$, which  hardly deviates from the true value, just the last figure is not correct.

\subsection*{Example 2}

\textit{Let the equation $x^4-2xx+4x=8$ be propounded.} \\[2mm]
Put $y=x^4-2xx+4x-8$; it will be $dy=4dx(x^3-x+1)$,

\begin{equation*}
p= \frac{1}{4(x^3-x+1)}, \quad \frac{dp}{dx}=\frac{-3xx+1}{4(x^3-x+1)^2}.
\end{equation*}   
 Therefore,
 
\begin{equation*}
 q= \frac{-3xx+1}{16(x^3-x+1)^3}, \quad \frac{dq}{dx}=\frac{21x^4-12xx-6x+3}{16(x^3-x+1)^4} \quad \text{and} \quad r=\frac{21x^4-12xx-6x+3}{64(x^3-x+1)^5} \quad \text{etc.},
\end{equation*}
from which the root of the propounded equation will be

\begin{equation*}
f=x- \frac{y}{4(x^3-x+1)}-\frac{(3xx-1)yy}{32(x^3-x+1)^3}-\frac{(7x^4-4xx-2x+1)y^3}{128(x^3-x+1)^5}-\text{etc.}
\end{equation*}
Therefore, it is necessary to attribute an appropriate value to $x$, so that series becomes convergent. But at first it is perspicuous, if  a value rendering $x^3-x+1=0$ would be attributed to $x$, that  all terms of there series except for the first would become infinite and nothing can be concluded from this. Therefore, it is convenient to assign such a value to $x$ that both $y$ becomes small and $x^3-x+1$ not very large. Let $x=1$; it will be $y=-5$ and

\begin{equation*}
f=1+\frac{5}{4}-\frac{25}{16}+\frac{125}{64}-\text{etc.};
\end{equation*}
because the three terms $\frac{5}{4}-\frac{25}{16}+\frac{125}{64}$ agree with the terms of the geometric progression, whose sum is $\frac{5}{9}$, it will approximately be $f=\frac{14}{9}$. Therefore, let us put $x=\frac{3}{2}$; it will be

\begin{equation*}
y= -\frac{23}{16} \quad \text{and} \quad x^3-x+1=\frac{23}{8},
\end{equation*}
whence 

\begin{equation*}
f=\frac{3}{2}+\frac{1}{8}-\frac{1}{64}+\frac{391}{256 \cdot 529}-\text{etc.}=1.61.
\end{equation*}
 Now put $x=1.61$; it will be
 
 \begin{small}
 \begin{center}
\begin{tabular}{lcl|crclccc}
$\log x$ &$~=~$ & $0.2068259$&  & $x$ & $~=~$ & $1.61$ & let & $x^3-x+1=z$ \\ 
$\log x^2$ &$~=~$ & $0.4136518$ &   & $x^2$ & $~=~$ & $2.5921$  \\ 
$\log x^3$ &$~=~$ & $0.6204777$ &   & $x^3$ & $~=~$ & $4.173281$  \\
$\log x^4$ &$~=~$ & $0.8273036$ &   & $x^4$ & $~=~$ & $6.718983$  \\
 & &  & \quad hence \quad   &  &  &   \\   
$\log (-y)$ &$~=~$ & $8.4016934$ &    & $y$ & $~=~$ & $~-~ 0.025217$  \\   
$\log z $ &$~=~$ & $0.5518502$ &    & $z$ & $~=~$ & $3.563281$  \\ \cline{3-3}
$\log \frac{-y}{z}$ &$~=~$ & $7.8498432$ &    &   &  &   \\
$\log 4$ &$~=~$ & $0.6020600$ &    &   &  &   \\      \cline{3-3}
$\log \frac{-y}{4z}$ &$~=~$ & $7.2477832$ &    & $\frac{-y}{4z}$   & $~=~$  & $0.0017692$   \\
$\log (3xx-1)$ &$~=~$ & $0.8309926$ &    & $3xx-1$   & $~=~$  & $6.7763$   \\
$\log y^2$ &$~=~$ & $6.8033868$ &    &    &   &    \\ \cline{3-3}
  &  & $7.6343794$ &    &    &   &    \\ 
  $\log z^3$ &$~=~$ & $1.6555506$ &    &    &   &    \\ \cline{3-3}
    &  & $5.9788288$ &    &    &   &    \\ 
    $\log 32$ &$~=~$ & $1.5051500$ &    & $\dfrac{(3xx-1)^2y^2}{32z^3}$   & $~=~$   & $0.00002976$    \\ \cline{3-3}
     &$~=~$ & $4.4736788$ &    &    &   &    \\ \cline{3-3}
\end{tabular}
\end{center} 
\end{small}
\begin{equation*}
\text{Therefore } f=1.6117662.
\end{equation*}

\paragraph*{§242}

This method to find the roots of equations approximately also extends to transcendental quantities. Let us find the number $x$, whose  logarithm  has a given ratio of $1$ to $n$ to the number $x$, and one will have this equation $x-n \log x =0$; but let $k$ be the modulus of these logarithms, such that these logarithms are obtained, if  the hyperbolic logarithms are multiplied by $k$; it will be $d. \log x= \frac{kdx}{x}$. Therefore, put $x- n \log x =y$ and let $f$ be the value of $f$ in question which renders $x= n \log x$. Therefore, since $y= x- n \log x$,  it will be

\begin{equation*}
dy= dx- \frac{kndx}{x}=\frac{dx(x-kn)}{x}
\end{equation*}
and 

\begin{alignat*}{9}
&\frac{dx}{dy}&&=p &&= \frac{x}{x-kn}, \quad && \text{whence} \quad && dp &&= -\frac{kndx}{(x-kn)^2}, 
\intertext{therefore} \\
&\frac{dp}{dy}&&=q &&= \frac{-knx}{(x-kn)^3}, \quad &&   && dq &&= \frac{2knxdx+k^2n^2dx}{(x-kn)^4} \\[3mm]
&\frac{dq}{dy}&&=r &&= \frac{knx(2x+kn)}{(x-kn)^5}.
\end{alignat*}
Therefore, it will be

\begin{equation*}
f= x- \frac{xy}{x-kn}-\frac{knxyy}{2(x-kn^3)}-\frac{knxy^3(2x+kn)}{6(x-kn)^5}-\text{etc.}
\end{equation*}
Below [§ 272] we will show that this problem only admits a solution, if $kn >e$ while $e$ is the number whose hyperbolic logarithm is $=1$, or it must be $kn >2.7182818$.

\subsection*{Example}

\textit{To find a number, other than $10$, whose tabulated logarithm becomes equal to the tenth part of the number itself.}\\[2mm]
Since  the question is on tabulated logarithms, it will be $k=0.43429448190325$ and, because of $n=10$, one will have $kn=4.3429448190325$. Now, having put $x=1$, it will be $y=1$ and it will be

\begin{equation*}
f=1+\frac{1}{3.3429}+\frac{2.1714724}{(3.3429)^3}-\text{etc.}
\end{equation*}
and so it will approximately be $f=1.37$. Therefore, set $x=1.37$; it will be $\log x= 0.136720567156406$ and, because of $y=x- 10 \log x$, it will be

\begin{equation*}
y=0.00279432843594 \quad \text{and} \quad -x+kn = 2.9729448190325.
\end{equation*}
Therefore, let

\begin{center}
\begin{tabular}{lclcrclccc}
$\log x$ &$~=~$ & $0.1367205$ \\ 
$\log y$ &$~=~$ & $7.4462773$ \\ \cline{3-3} 
         &      & $7.5829978$ \\
$\log (kn-x)$ &$~=~$ & $0.4731866$ \\ \cline{3-3} 
  &  & $7.1098112$ & \quad $\dfrac{-xy}{x-kn}=0.00128769$ \\    
\end{tabular}
\end{center}
Further, because the third term is $-\frac{knxyy}{2(x-kn)^2}=\frac{kny}{2(x-kn)^2} \cdot \frac{-xy}{x-kn}$, it will be
\begin{center}
\begin{tabular}{r clcrclccc}
$\log \dfrac{-xy}{x-kn}$ &$~=~$ & $7.1098112$ \\
$\log y$ &$~=~$ & $7.4462773$ \\
$\log kn$ &$~=~$ & $0.6377842$ \\ \cline{1-3}
&  & $5.1938727$ \\    
$\log (kn-x)^2$ &$~=~$ & $0.9463732$ \\ \cline{3-3}
  &  & $4.2474995$ \\      
  $\log 2$ &$~=~$ & $0.3010300$ \\ \cline{3-3}
$\log \text{third term}$ &$~=~$ & $3.9464695$ \\
I. term $x$ &$~=~$ & $1.37$ \\  
II. term  &$~=~$ & $0.00128769$ \\
III. term  &$~=~$ & $0.00000088$ \\ \cline{3-3}
$f$ &$~=~$ & $1.37128857$ \\
$\log f$ &$~=~$ & $0.137128857$ \\      
\end{tabular}
\end{center}

\paragraph*{§243}

If the equation was an exponential equation, it can be reduced to a logarithmic one; so, if the value of $x$ is in question, that  $x^x=a$, it will be $x \log x= \log a $. Therefore, having put $y= x \log x- \log a$, it will be

\begin{equation*}
dy= dx \log x+dx \quad \text{and} \quad \frac{dx}{dy}=p= \frac{1}{1+ \log x}
 \end{equation*}
and then

\begin{equation*}
dp= \frac{-dx}{x(1+ \log x)^2} \quad \text{and} \quad \frac{dp}{dy}=q= \frac{-1}{x(1+\log x)^3},
\end{equation*}
\begin{equation*}
dq= \frac{dx}{xx(1+\log x)^3}+\frac{3dx}{xx(1+ \log x)^4} \quad \text{and hence} \quad \frac{dq}{dy}=r= \frac{1}{xx(1+ \log x)^4}+\frac{3}{xx(1+\log x)^5};
\end{equation*}
further, it will be

\begin{equation*}
dr= \frac{-2dx}{x^3(1+ \log x)^4}-\frac{10dx}{x^3(1+ \log x)^5}-\frac{15dx}{x^3(1+ \log x)^6},
\end{equation*}
therefore,

\begin{equation*}
s=\frac{-2}{x^3(1+ \log x)^5}-\frac{10}{x^3(1+ \log x)^6}-\frac{15}{x^3(1+\log x)^7}
\end{equation*}
 and
 
 \begin{equation*}
t= \frac{6}{x^4(1+ \log x)^6}+\frac{40}{x^4(1+\log x)^7}+\frac{105}{x^4(1+ \log x)^8}+\frac{105}{x^4(1+\log x)^9}, 
 \end{equation*}
\begin{equation*}
u= \frac{-24}{x^5(1+\log x)^7}-\frac{196}{x^5(1+\log x)^6}-\frac{700}{x^5(1+\log x)^9}-\frac{1260}{x^5(1+ \log x)^{10}}-\frac{945}{x^5(1+ \log x)^{11}}. 
\end{equation*}
Therefore, if the true value of $x$ is $=f$, such that $f^f=a$, it will be

\begin{alignat*}{19}
&f=x-\frac{y}{1+ \log x}-\frac{yy}{2x(1+\log x)^3}&&-\frac{y^3}{2xx(1+\log x)^5}&&-\frac{5y^4}{8x^3(1+\log x)^7}&&-\frac{7y^5}{8x^4(1+\log x)^9} \\[3mm]
&    && -\frac{y^3}{6x^2(1+\log x)^4}&&-\frac{5y^4}{12x^3(1+\log x)^6}&&-\frac{7y^5}{8x^4(1+\log x)^8} \\[3mm]
& && &&-\frac{y^4}{12x^3(1+\log x)^5}&&-\frac{y^5}{3x^4(1+\log x)^7} \\[3mm]
& && && &&-\frac{y^5}{20x^4(1+\log x)^6}\\
& &&  \text{etc.}
\end{alignat*}
Therefore, this expression, having continued it to infinity, whatever value is substituted for $x$, and having taken $y=x \log x - \log a$, will give the true value of $f$. So, if one puts $x=1$, it will be $y=- \log a$ and

\begin{equation*}
f=1+ \log a- \frac{(\log a)^2}{2}+\frac{2(\log a)^3}{3}-\frac{9(\log a)^4}{8}+\frac{
32 (\log a)^5}{15}-\frac{625(\log a)^6}{144}-\text{etc.},
\end{equation*} 
where it is to be noted that $\log a$ is the hyperbolic logarithm of $a$.

\subsection*{Example}

\textit{To find the number $f$ that  $f^f=100$.}\\[2mm]
Because 

\begin{equation*}
a=100 \quad \text{and} \quad y=x \log x- \log a = x \log x- \log 100,
\end{equation*}
since obviously $f>3$ and $<4$, put $x=\frac{7}{2}$ and it will be

\begin{center}
\begin{tabular}{r clcrclccc}
$\log x$ &$~=~$ & $1.25276296849$ \\
   $x\log x$ &$~=~$ & $4.38467038972$ \\
   $\log 100$ &$~=~$ & $4.60517018599$ \\ \cline{3-3}
  $y$ &$~=~$ & $-0.22049979627$ \\
  $1+\log x$ &$~=~$ & $2.25276296849.$ \\ 
\end{tabular}
\end{center} 
Therefore,  using ordinary logarithms, it will be

\begin{center}
\begin{tabular}{lclcrclccc}
$\log (-y)$ &$~=~$ & $9.3434083$ \\
$\log (1+\log x)$ &$~=~$ & $0.3527156$  & \quad $\dfrac{-y}{1+ \log x}=0.0978797$ \\ \cline{3-3}
  & & $8.9906927$ \\
  $\log y^2$ &$~=~$ & $8.6868166$ \\
  $3 \log (1+ \log x)$ &$~=~$ & $1.0581468$ \\ \cline{3-3}
  $ $ &$ $ & $7.6286698$ \\
  $\log 2x= \log 7$ &$~=~$ & $0.8450980$  & \quad  $\dfrac{y^2}{2x(1+ \log x)^3}=0.0006075.$ \\ \cline{3-3}
    & & $6.7835718$
\end{tabular}
\end{center} 
\begin{equation*}
\text{Therefore, it will approximately be } f=3.5972722;
\end{equation*}
but having additionally taken the following terms, it will be $f=3.5972852.$

\paragraph*{§244}

But moreover differential calculus has an extraordinary use in the resolution of equations, if a certain relation among the roots was known. Let the equation $y=0$ be propounded, in which $y$ is an arbitrary function of $x$. If now, for the sake of an example, it is known that the roots of this equation differ by the given quantity $a$, these two roots will easily be found the following way. Let $x$ denote the smaller of these two roots; the larger will be $=x+a$; therefore, because the function $y$ vanishes, if $x$ denotes any of the roots of the equation $y=0$, it will also vanish, if  one writes $x+a$ instead of $x$. Therefore, it will be

\begin{equation*}
0=y+\frac{ady}{dx}+\frac{a^2ddy}{2dx^2}+\frac{a^3d^3y}{6dx^3}+\text{etc.}
\end{equation*} 
Therefore, since $y=0$, it will also be

\begin{equation*}
0=\frac{ady}{dx}+\frac{a^2ddy}{2dx^2}+\frac{a^3d^3y}{6dx^3}+\text{etc.}
\end{equation*}
which two equations taken at the same time, using the method of elimination, will give the value of the root $x$, which is smaller than the other root by the quantity $a$.

\subsection*{Example}

\textit{Let this equation $x^5-24x^4+49xx-36=0$ be propounded, which is known from anywhere to have two roots differing by $1$.}\\[2mm]
Having put $y=x^5-24x^3+49xx-36$ it will be

\begin{alignat*}{9}
&~~~~\frac{dy}{dx}&&= 5x^4 &&- 72x^2 &&+98x \\[3mm]
&~~\frac{ddy}{2dx^2}&&=10x^3 &&- 72x &&+49 \\[3mm]
&~~\frac{d^3y}{6 dx^3}&&= 10x^2 &&-24 \\[3mm]
&~ \frac{d^4 y}{24 dx^4}&&= 5x \\[3mm]
&\frac{d^5 y}{120 dx^5}&&=1.
\end{alignat*}
But, because of $a=1$, it will be

\begin{equation*}
A \cdots 5x^4+10x^3-62x^2+31x+26 =0.
\end{equation*}
But

\begin{equation*}
B \cdots x^5-24x^3+49xx-36=0.
\end{equation*}
Multiply the upper equation by $x$ and the lower equation by $5$ and subtract the one from the other and it will remain

\begin{equation*}
10x^4+58x^3-214x^2+26x+180=0
\end{equation*}
or 
\begin{equation*}
C \cdots 5x^4+229x^3-107x^2+13x+90=0,
\end{equation*}
having subtracted the first $A$ from which it will be

\begin{center}
\begin{tabular}{rrcrcrcrcrcrcrcrcrcr}
$D \cdots$ &$19x^3$ & $-$ & $45x^2$ & $-$ &$18x$ & $+$ & $64$ & $=$ & $0$. \\
$D \cdot 5x \cdots$ &$95x^4$ & $-$ & $225x^3$ & $-$ &$90x^2$ & $+$ & $320x$ & $=$ & $0$. \\
$A \cdot 19 \cdots $ &$95x^4$ & $+$ & $190x^3$ & $-$ &$1178x^2$ & $+$ & $589x$ & $+$ & $494$ & $=$ &$0$. \\ \cline{3-12}
$E \cdots$ & &  & $415x^3$ & $-$ &$1088x^2$ & $+$ & $269x$ & $+$ & $494$ & $=$ &$0$. \\ 
$D \cdot 415 \cdots $ & &  & $7885x^3$ & $-$ &$18675x^2$ & $-$ & $7470x$ & $+$ & $26560$ & $=$ &$0$. \\ 
$E \cdot 19 \cdots $ & &  & $7885x^3$ & $-$ &$20672x^2$ & $+$ & $5111x$ & $+$ & $8386$ & $=$ &$0$. \\ \cline{5-12}
$F \cdots $ & &  &  &  &$1997x^2$ & $-$ & $12581x$ & $+$ & $17174$ & $=$ &$0$. \\ 
\end{tabular}
\end{center} 
And
\begin{center}
\begin{tabular}{rrcrcrcrcrcrcrcrcrcr}
$D \cdot 247 \cdots$ &$4693x^3$ & $-$ & $11115x^2$ & $-$ &$4446x$ & $+$ & $15808$ & $=$ & $0$ \\
$D \cdot 32 \cdots$ &$13280x^3$ & $-$ & $34816x^2$ & $+$ &$8608x$ & $+$ & $15808$ & $=$ & $0$ \\ \cline{2-10}
  &$8587x^3$ & $-$ & $23701x^2$ & $+$ &$13054x$ & $=$ & $0$ &  &  \\ \cline{2-10}
  $G \cdots $  &  &  & $8587x^2$ & $-$ &$23701x$ & $+$ & $13054$ & $=$  & $0$  \\ 
    $F \cdot 8587 \cdots $  &  &  & $17148239x^2$ & $-$ &$108033047x$ & $+$ & $147473138$ & $=$  & $0$  \\ 
    $G \cdot 1997 \cdots $  &  &  & $17148239x^2$ & $-$ &$47330897x$ & $+$ & $26068838$ & $=$  & $0$  \\  \cline{2-10}
          &  &  &  &  &$60702150x$ & $-$ & $121404300$ & $=$  & $0$.     
\end{tabular}
\end{center} 
From this equation it follows that $x=2$ and therefore also $x=3$ will be a root of the equation, both of which values indeed satisfy the equation.

\paragraph*{§245}

But this operation can be done without use of differential calculus, because the same equation the differential calculus yielded, results, if in the propounded equation one writes $x+1$ instead of $x$. Furthermore, this method of elimination is too laborious, and if the equations would be of higher degree, the amount of work would be simply too much to handle;  and this is true even more for transcendental equations. But if we put that two roots of the propounded equation $y=0$ are equal to each other, then, because of $x=a$, the differential equation goes over into this one $\frac{dy}{dx}=0$. Therefore, if any equation $y=0$ had two equal roots, it will be $\frac{dy}{dx}=0$ and these two roots taken together will yield the value of $x$, to which these two roots are equal. Therefore, vice versa, if the two equations $y=0$ and $\frac{dy}{dx}=0$ have a common root, it will be a double root of the equation $y=0$. But this happens, if, after the quantity $x$  was completely eliminated by means of these two equations $y=0$ and $\frac{dy}{dx}=0$, one gets to an identical equation. So if the equation

\begin{equation}
x^3-2xx-4x+8
\end{equation}
was propounded, it will also be $3xx-4x-4=0$, whose double added to the first gives

\begin{equation*}
x^3+4xx-12x=0 \quad \text{or} \quad xx+4x-12=0,
\end{equation*}
whose triple is

\begin{center}
\begin{tabular}{rrrrrrrrrrrrrrrrrrrrrrrrrr}
 & $3xx$ &$+$ & $12x$ & $-$ & $36$& $=$   &$0$ \\
subtract \quad & $3xx$ &$-$ & $4x$ & $-$ & $4$& $=$   &$0$ \\ \cline{2-8}
  &   &  & $16x$ & $-$ & $32$& $=$   &$0$ \\
  &   &  & $x$ & $-$ & $2$& $=$   &$0$ \\
\end{tabular}
\end{center}  
Therefore, because $x=2$ results, substitute this value in one of the preceding $3xx-4x-4=0$ and the identical equation $12-8-4=0$ will result, whence one concludes that the propounded equation $x^3-2xx-4x+8=0$ has two equal roots, namely $x=2$.

\paragraph*{§246}

Therefore, if one has an algebraic equation of no matter how many dimensions

\begin{equation*}
x^n+Ax^{n-1}+Bx^{n-2}+Cx^{n-3}+Dx^{n-4}+\text{etc.}=0,
\end{equation*}
which has two equal roots, it will also be

\begin{equation*}
nx^{n-1}+(n-1)Ax^{n-2}+(n-2)Bx^{n-2}+(n-3)Cx^{n-4}+(n-4)Dx^{n-5}+\text{etc.}=0.
\end{equation*}
This double root of that equation will at the same time be a root of this last equation, of course. Multiply the first equation by $n$ and  subtract the second multiplied by $x$ from it and this new equation will result

\begin{equation*}
Ax^{n-1}+2Bx^{n-2}+3Cx^{n-3}+4Dx^{n-4}+\text{etc.}=0.
\end{equation*}
Now add the first multiplied by $a$ and the latter multiplied by $b$; it will be

\begin{equation*}
ax^n+(a+b)Ax^{n-1}+(a+2b)Bx^{n-2}+(a+3b)Cx^{n-3}+\text{etc.}=0,
\end{equation*}
which equation combined with the propounded itself will show equal roots, if the propounded one has some. Therefore, because the quantities $a$ and $b$ are arbitrary, the coefficients $a$, $a+b$, $a+2b$ etc. represent any arithmetic progression. Therefore, if any equation has two equal roots, they will be found, if the each term of the propounded equation is multiplied by terms of a certain arithmetic progression, respectively; for, the new equation resulting this way will also contain the root, which is contained twice in the propounded one. So, if the the terms of the equation

\begin{equation*}
x^n+Ax^{n-1}+Bx^{n-2}+Cx^{n-3}+Dx^{n-4}+\text{etc.}
\end{equation*}
are multiplied by this arithmetic progression

\begin{equation*}
a, \quad a+b, \quad a+2b, \quad a+3b, \quad a+4b \quad \text{etc.};
\end{equation*}
this new equation will result

\begin{equation*}
ax^n+(a+b)Ax^{n-1}+(a+2b)Bx^{n-2}+(a+3b)x^{n-2}+(a+3b)Cx^{n-3}+\text{etc.}=0,
\end{equation*}
which combined with the latter will show the equal roots. And this is the well-known rule to find equal roots of any equation.

\paragraph*{§247}

If the equation $y=0$ has three equal roots, it will not only be $\frac{dy}{dx}=0$, but it will also be $\frac{ddy}{dx^2}=0$, if  one substitutes the value of the root for $x$, which is a triple root of the equation $y=0$.  To show this let us put that the equation $y=0$ has three roots $x$, $x+a$, $x+b$ etc., of  which the first differs from the other ones by $a$ and $b$, respectively; and since $y$ vanishes, if  one writes $x+a$ or $x+b$ instead of $x$, it will be

\begin{alignat*}{9}
y && && && && && &&=0 \\
y &&~+~\frac{ady}{dx}&&~+~\frac{a^2ddy}{2dx^2}&&~+~\frac{a^3d^3y}{6dx^3}&&~+~\frac{a^4 d^4y}{24 dx^4}&&~+~\text{etc.}&&~=~0 \\[3mm]
y &&~+~\frac{bdy}{dx}&&~+~\frac{b^2ddy}{2dx^2}&&~+~\frac{b^3d^3y}{6dx^3}&&~+~\frac{b^4 d^4y}{24 dx^4}&&~+~\text{etc.}&&~=~0;
\end{alignat*}
if the first is subtracted from the two last ones, it will be

\begin{alignat*}{9}
\frac{dy}{dx}&&~+~\frac{addy}{2dx^2}&&~+~\frac{a^2d^3y}{6dx^3}&&~+~\frac{a^3 d^4y}{24 dx^4}&&~+~\text{etc.}&&~=~0 \\[3mm]
\frac{dy}{dx}&&~+~\frac{bddy}{2dx^2}&&~+~\frac{b^2d^3y}{6dx^3}&&~+~\frac{b^3 d^4y}{24 dx^4}&&~+~\text{etc.}&&~=~0 \\
\end{alignat*}
Also subtract these from each other and, having divided by $a-b$, it will be

\begin{equation*}
\frac{ddy}{2dx^2}+ \frac{(a+b)d^3y}{6dx^3}+\frac{(aa+ab+bb)d^4y}{24dx^4}+\text{etc.}=0.
\end{equation*}
Now put $a=0$ and $b=0$ such that these three roots are equal to each other, and, because of the vanishing terms, it will be

\begin{equation*}
y=0, \quad \frac{dy}{dx}=0 \quad \text{and} \quad \frac{ddy}{dx^2}=0.
\end{equation*}

\paragraph*{§248}

Therefore, if the equation $y=0$ has three equal roots, say $f,f,f$, then this quantity $f$ will also be a root not only of this equation $\frac{dy}{dx}=0$, but also of this one $\frac{ddy}{dx^2}=0$. Therefore, recalling the results we demonstrated before on two equal roots of equations, it is obvious, because $f$ is the common root of the equation $\frac{dy}{dx}=0$ and its differential $\frac{ddy}{dx^2}=0$, that it has to be contained twice in the equation $\frac{dy}{dx}=0$. Therefore, if the equation

\begin{equation*}
x^n+Ax^{n-1}+Bx^{n-2}+Cx^{n-3}+Dx^{n-4}+\text{etc.}=0
\end{equation*}
 contains three equal roots $f,f,f$, if its terms are multiplied by the terms of a certain arithmetic progression, then the resulting equation will have two equal roots $f$ and $f$; therefore, it can be multiplied by an arithmetic progression again so that an equation containing the root $f$ once results. Therefore, one will obtain three equations having the common root $f$, from whose combination this root will easily be found. For, if arithmetic progressions of such a kind are chosen, whose either first or last terms are  $=0$, then an equation of one degree lower will result and so the elimination will be even easier.
 
\paragraph*{§249}

In like manner it can be shown, if the equation $y=0$ has four equal roots $f,f,f,f$, that for $x=f$ it will not only be $y=0$, $\frac{dy}{dx}=0$ and $\frac{ddy}{dx^2}=0$, but it will also be $\frac{d^3y}{dx^3}=0$. As the equation $y=0$ contains the root $x=f$ four times, so the equation $\frac{dy}{dx}$ will contain the same three times, the equation $\frac{ddy}{dx^2}=0$ twice and the equation $\frac{d^3y}{dx^3}=0$ once.  This will also be seen more easily, if we consider that the function $y$ has to have a form of this kind $(x-f)^4X$ in this case, where $X$ denotes any function of $x$. Having assumed this form, it will be

\begin{equation*}
\frac{dy}{dx}=(x-f)^3 \left(4X+\frac{(x-f)dX}{dx}\right)
\end{equation*}
and hence be divisible by $(x-f)^3$. Further, $\frac{ddy}{dx^2}$ will have the factor $(x-f)^2$ and $\frac{d^3y}{dx^3}$ the factor $x-f$; from this it is perspicuous, if the root $f$ is contained in the equation $y=0$, it has to be contained in the equation $\frac{dy}{dx}=0$ three times, in the equation $\frac{ddy}{dx^2}=0$ twice and in $\frac{d^3y}{dx^3}=0$  still once.


\title{On Maxima and Minima
\footnote{
Original title: "'De Maximis et Minimis"', first published as part of the book \textit{Institutiones calculi differentialis cum eius usu in analysi finitorum ac doctrina serierum}, 1755, reprinted in \textit{Opera Omnia}: Series 1, Volume 10, pp. 446 - 474, Eneström-Number E212, translated by: Alexander Aycock for the "`Euler-Kreis Mainz"'}}
\author{Leonhard Euler}
\date{ }
\maketitle

\paragraph*{§250}

If a function of $x$ was of such a nature that, while the values of $x$ increase, the function itself increases  or decreases continuously, then this function will have no maximum or minimum value. For, whatever value of this function is considered, the following value will be larger, the preceding value on the other hand will be smaller. A function of this kind is $x^3+x$, whose values for increasing $x$ increase continuously, but for decreasing $x$ decrease continuously; therefore, this function cannot have  a maximum value, if not the maximum point, i.e. infinity, is attributed to $x$; and in like manner, it will have the minimum value, if one puts $x= - \infty$. But if the function was not of such a nature, that, while $x$ increases, it either increases or decreases continuously, it will have a maximum or minimum somewhere else, this means a value of such a kind, which is either greater or smaller than the preceding and following values. For example, this function $xx-2x+3$ has a minimum value, if one puts $x=1$; for, whatever other value is attributed to $x$, the function will always have a larger value.

\paragraph*{§251}

But to understand the nature of maxima and minima more clearly, let  $y$ be a function of $x$, which has a maximum  for $x=f$, and it is seen, if $x$  is then assumed to be either greater or smaller than $f$, that  the value of $y$ to result from this will be smaller than the value it has for $x=f$. In like manner, if  the function $y$  has a minimum value for $x=f$, it is necessary, that,  no matter whether $x$ is assumed to be larger  or smaller than $f$, always a larger value of $y$ results; and this is the definition of absolute maxima and minima. But furthermore, this function $y$ is also said to have a maximum value for $x=f$, if that value of $y$ was larger than the closest ones, either the following or the preceding, which result, if $x$ is assumed to be a little bit larger or smaller than $f$, even though by substituting other values for $x$ the function $y$ might have larger values. Similarly, the function $y$ is said to have a minimum value for $x=f$, if that value was smaller than those, which it has, if  the closest larger or smaller values than $f$ are substituted for $x$. And in that sense we will use the terms maxima and minima.

\paragraph*{§252}

But before we show how to find these maxima and minima, we note that this investigation  extends only to those functions of $x$, which we called \textit{uniform} above and which are of such a nature that for each value of $x$ they in like manner have only one value. But we called functions \textit{biform} and \textit{multiform} which for each value of $x$ have two or several values at the same time; examples of this kind of functions are the roots of quadratic equations and higher order equations. Therefore, if $y$ was a biform or multiform function of such a kind of $x$, it cannot be said to take on a maximum or minimum value for $x=f$; for, since  it has either two or more values at the same time for $x=f$ and  the same is true for the preceding and the following contiguous values, one cannot decide whether $x=f$ is a maximum or minimum, if not by coincidence all values corresponding to the respective values of $x$, except one, of the function $y$, are imaginary; in this case functions of this kind are counted to the class of uniform functions. Therefore, we will at first consider the class of uniform functions in this chapter; but then, in the following chapters, we will show, how maxima and minima of multiform functions must be considered.

\paragraph*{§253}

Therefore, let $y$ be a uniform function, which hence, no matter which value is substituted for $x$, always has precisely one real value, and let $x$ denote the value, which induces the maximal or minimal value to the function $y$. Therefore, in the first case, no matter whether one substitutes $x+\alpha$ or $x-\alpha$ for $x$, the value of $y$ will be smaller than for $\alpha =0$, in the second case on the other hand larger. Therefore, because, having written $x+\alpha$ instead of $x$, the function $y$ goes over into

\begin{equation*}
y+\frac{\alpha dy}{dx}+\frac{\alpha^2ddy}{2dx^2}+\frac{\alpha^3 d^3y}{6 dx^3}+\text{etc.},
\end{equation*}
but, having written $x- \alpha$ instead of $x$, it goes over into

\begin{equation*}
x- \frac{\alpha dy}{dx}+\frac{\alpha^2 ddy}{2dx^2}-\frac{\alpha^3d^3y}{6dx^3}+\text{etc.},
\end{equation*}
it is necessary that  in the case of a maximum

\begin{alignat*}{19}
&y &&> y &&+\frac{\alpha dy}{dx}&&+\frac{\alpha^2 ddy}{2dx^2}&&+\frac{\alpha^3 d^3y}{6 dx^3}&&+\text{etc.} 
\intertext{and}
&y &&> y &&-\frac{\alpha dy}{dx}&&+\frac{\alpha^2 ddy}{2dx^2}&&-\frac{\alpha^3 d^3y}{6 dx^3}&&+\text{etc.}\\
\end{alignat*}
But in the case, in which the value of $y$ is a minimum value, it will be

\begin{alignat*}{19}
&y &&<y &&+\frac{\alpha dy}{dx}&&+\frac{\alpha^2 ddy}{2dx^2}&&+\frac{\alpha^3 d^3y}{6 dx^3}&&+\text{etc.} \\[3mm]
&y &&< y &&-\frac{\alpha dy}{dx}&&+\frac{\alpha^2 ddy}{2dx^2}&&-\frac{\alpha^3 d^3y}{6 dx^3}&&+\text{etc.}\\
\end{alignat*}

\paragraph*{§254}

Since these equations have to hold, if $\alpha$ denotes a very small quantity, let us assume $\alpha$ to be so small that its higher powers can be omitted,  and then so for the case of the maximum as the minimum it has to be $\frac{\alpha dy}{dx}=0$. For, if $\frac{\alpha dy}{dx}$ was not $=0$, the value of $y$ could be neither be maximum nor minimum value. Therefore,  to investigate  maxima or minima one has the general rule, that the differential of propounded $y$ is to be put equal to zero, and that  value of $x$, which renders the function  either maximal or minimal, will be a root of that equation. But whether the value of $y$ found this way is a maximum or a minimum value, is  not clear at this point; it can even happen, that $y$ is neither a maximum nor a minimum value in this case; for, we only found that in both cases it will be $\frac{dy}{dx}=0$ and we  did not vice versa prove, if $\frac{dy}{dx}=0$, that also a maximum or minimum value of $y$  results.

\paragraph*{§255}

Nevertheless, to investigate cases, in which the value of $y$ is either a maximum or minimum value, first all the roots of equation $\frac{dy}{dx}=0$ are to be found. Having found these, it is  to be checked, whether for those values the function $y$ has a maximum or minimum value or  none of both is the case. For, we will show that it can happen that there is neither a maximum or minimum, even though $\frac{dy}{dx}=0$.\\
Let $f$ be one of the values of $x$ satisfying the equation

\begin{equation*}
\frac{dy}{dx}=0,
\end{equation*}
and substitute this value in the expressions $\frac{ddy}{dx^2}$, $\frac{d^3y}{dx^3}$ etc. and by this substitution let

\begin{equation*}
\frac{ddy}{dx^2}=p, \quad \frac{d^3y}{dx^3}=q, \quad \frac{d^4y}{dx^4}=r \quad \text{etc.}
\end{equation*}
But, having written $f$ instead of $x$,  let the function $y$ go over into $F$, and if one writes $f+\alpha$ instead of $x$, this function will go over into

\begin{equation*}
F+\frac{1}{2}\alpha^2 p+ \frac{1}{6}\alpha^3q+\frac{1}{24}\alpha^4 r+\text{etc.};
\end{equation*}
but if one writes $f-\alpha$ instead of $x$, this expression will result

\begin{equation*}
F+\frac{1}{2}\alpha^2 p- \frac{1}{6}\alpha^3q+\frac{1}{24}\alpha^4 r- \text{etc.};
\end{equation*}
hence it is plain, if $p$ was a positive quantity, that both values will be larger than $F$, at least if $\alpha$ denotes a very small quantity, and therefore the value $F$, which the function $y$ has for $x=f$, will be a minimum value. But if $p$ is a negative quantity, then the value $x=f$ will induce a maximal value to the function $y$.

\paragraph*{§256}

But if it was $p=0$, then the value of $q$ is to be considered; if it was not $=0$, the value of $y$ will be neither a maximum nor a minimum value; for,  having put $x=f+ \alpha$, it will be $F+\frac{1}{6}\alpha^3 q >F$ and, having put $x=f-\alpha$, it will be $F- \frac{1}{6}\alpha^3 q<F$. But if it also was $q=0$, the quantity $r$ is to be considered; if it had a positive value, the value of the function $F$, which is has for $x=f$, will be a minimum value; but if $r$ has a negative value, $F$ will be a maximum value. But if also $r$ vanishes, one has to consider the value of the following letter $s$, and has to argue as  for the letter $q$. If $s$ was not $=0$, then the value $F$ will be neither a maximum nor a minimum value; but if also $s=0$, then the following letter $t$, if it has a positive value, will indicate a minimum; but if it has a negative value, it will indicate a maximum. But if also this letter $t$ vanishes, then  one has to proceed  to the next letter and argue as  for the preceding ones. And this way  one can decide for any root of the equation $\frac{dy}{dx}=0$, whether  the function $y$ has a maximum or minimum value or none of both; and this way all maxima and minima which the function $y$ can have will be found.

\paragraph*{§257}

Therefore, if the equation $\frac{dy}{dx}=0$ has two equal roots, such that it has the quadratic factor $(x-f)^2$, then  at the same time $\frac{ddy}{dx^2}$ will vanish for $x=f$ and it will be $p=0$, but not $q=0$. In this case the function $y$ will have neither a maximum nor a minimum value. But if the equation $\frac{dy}{dx}=0$ has three equal roots or $\frac{dy}{dx}$ has the cubic factor $(x-f)^3$, then, having put $x=f$, it will be $\frac{ddy}{dx^2}=0$ and $\frac{d^3y}{dx^3}=0$, but not $\frac{d^4y}{dx^4}=0$. Therefore, if the value of this term was positive, it will indicate a minimum value, if negative, a maximum value. Therefore, the rule explained before reduces to this, that, if the fraction $\frac{dy}{dx}$ had a factor $(x-f)^n$, while $n$ is an odd number, the function $y$, if in it one puts $x=f$, will have  either a maximum or minimum value, but if the exponent $n$ was an even number, the substitution $x=f$ will produce neither a maximum nor a minimum value.

\paragraph*{§258}

Furthermore, the finding a maximum or minimum is often simplified tremendously by the following considerations. In cases in which the function $y$ has a maximum or minimum value, each multiple of it, say $a y$, if $a$ was a positive quantity, will also be a maximum or minimum value, and in the same way $y^3$, $y^5$, $y^7$ etc. and in general $y^{n}$, if $n$ was a positive odd number, since formulas of this kind are of such a nature that for increasing $y$ they also increase and for decreasing $y$ they decrease. But in these cases, in which $y$ has maximum or minimum value, $-y$, $-ay$, $b-ay$ and in general $b-ay^n$, while $n$ is an odd positive integer, in reversed order, will have either minimum or maximum values. In the same way in the cases, in which $y$ has maximum or minimum value, these formulas $\frac{a}{y}$, $\frac{a}{y^3}$, $\frac{a}{y^5}$ etc. and in general $\frac{a}{y^n}\pm b$, while $a$ denotes a positive quantity and $n$ a positive odd number, in reverse order, will have a minimum or a maximum value; but if $a$ was a negative quantity, then these formulas will have a maximum value, if $y$ was a maximum value, and a minimum value, if $y$ is a minimum value.

\paragraph*{§259}

But these rules cannot be transferred to even powers the same way; for, since, if $y$ has a negative value, its even powers $y^2$, $y^4$ etc. induce positive values, it can happen, that, if $y$ has a minimum value, a negative one of course, that its even powers have maximum values. Therefore, having taken this into account, we will be able to affirm, if $y$ has a maximum or minimum value, while its value is positive, that its even powers $y^2$, $y^4$ etc. will also have a maximum or minimum, but if a negative value of $y$ was a maximum value, that its square $yy$ will be a maximum value, and otherwise, if a negative value of $y$ was a minimum value, that  $y^2$, $y^4$ etc. will have a maximum value. But if the even exponents of $y$ were negative, then the opposite will happen. Furthermore, what we mentioned here  on the even and odd exponents, will not only hold for integer numbers, but also for fractional ones, whose denominators are odd numbers; for,  the fractions $\frac{1}{3}$, $\frac{5}{3}$, $\frac{7}{3}$, $\frac{1}{5}$, $\frac{3}{5}$ etc. are equivalent to odd numbers in this case, but $\frac{2}{3}$, $\frac{4}{3}$, $\frac{2}{5}$, $\frac{4}{5}$, $\frac{6}{7}$ etc. are equivalent to even numbers.

\paragraph*{§260}

But if the denominators were even numbers, then, because, if $y$ has a negative value, its powers $y^{\frac{1}{2}}$, $y^{\frac{3}{4}}$ etc. become imaginary, here one can say only the following about them: If a positive value of $y$ was a maximum or a minimum value, $y^{\frac{1}{2}}$, $y^{\frac{3}{2}}$, $y^{\frac{1}{4}}$ etc. will also have either a maximum or minimum value, but on the other hand $y^{-\frac{1}{2}}$, $y^{-\frac{3}{2}}$, $y^{-\frac{1}{4}}$ etc. can have minimum or maximum values. But if these irrationalities  take on two values at the same, one positive, one negative, for the negative ones the contrary of that, what we said about the positive ones here, holds. But if a negative value of $y$ is a maximum or minimum value, since all powers 
of this kind become imaginary, one will not be able to  count them to maxima or minima. Therefore, by means of these auxiliary remarks, the investigation of maxima and minima is often simplified and otherwise would be extremely difficult.

\paragraph*{§261}

Because these things extend only to rational functions, which are the only uniform functions, at first let us expand polynomial functions and find their maxima and minima. Therefore, because functions of this kind are reduced to this form

\begin{equation*}
x^n+Ax^{n-1}+Bx^{n-2}+Cx^{n-3}+\text{etc.},
\end{equation*}
at first it is plain that their values cannot be greater than if one sets $x= \infty$; then on the other hand, if $x=- \infty$, the value of these formulas become $=\infty^n$, if $n$ is an even number, but $-\infty^n$, if $n$ is an odd number, which value therefore will be the smallest of all. But furthermore often other maxima and minima, in the sense we understand those terms, are given, what we will illustrate in the following examples.

\subsection*{Example 1}

\textit{To find the values of $x$ for which the function $(x-a)^n$ takes on maximum or a minimum value.}\\[2mm]
Having put $(x-a)^n=y$ it will be

\begin{equation*}
\frac{dy}{dx}=n(x-a)^{n-1};
\end{equation*} 
having put this expression $=0$, it will be $x=a$.  Therefore, because $\frac{dy}{dx}$ contains the factor $(x-a)^{n-1}$, from § 257 it is understood that $y$ cannot have a maximum or a minimum value, if $n-1$ is not an odd number or $n$ is not even. But since

\begin{equation*}
\frac{d^ny}{dx^n}=n(n-1)(n-2)\cdots 1,
\end{equation*}
i.e. a positive number, it follows that the value of $y$ for $x=a$ will turn out to be a minimum value. This is obvious, of course; for, having written $x=a$,  $y=0$, and if $x$ is put to be either greater or smaller than $a$, because of the even number $n$, the function $y$ will be positive, i.e. greater than zero; but if $n$ was an odd number, then the function $y=(x-a)^n$ can have neither a maximum nor a minimum value. But it is perspicuous that the same holds, if $n$ was a fractional number, no matter  whether it is odd or even, $(x-a)^{\frac{\mu}{\nu}}$ will have a minimum value for $x=a$, if $\mu$ was an even number and $\nu$ was an odd number; but if both were odd, neither a maximum nor a minimum value will exist.

\subsection*{Example 2}

\textit{To find the cases in which the value of this formula $xx+3x+2$ has a maximum or a minimum value.}\\[2mm]
Put $xx+3x+2=y$; it will be

\begin{equation*}
\frac{dy}{dx}=2x+3, \quad \frac{ddy}{2dx^2}=1.
\end{equation*}
Therefore, set $2x+3=0$; it will be $x=-\frac{3}{2}$. Whether this case produces a maximum or minimum value, will be found from the value $\frac{ddy}{2dx^2}=1$; since it is positive, whatever $x$ is, it indicates a minimum value. But for $x=-\frac{3}{2}$ we have $y=-\frac{1}{4}$, and if any other values are attributed to $x$, the value of $y$ to result from this will always be larger than $-\frac{1}{4}$. From the nature of the formula $xx+3x+2$ it is also seen that it has to have a minimum value; for, because it grows to infinity, if one puts $x= + \infty$ or $x=- \infty$, it is necessary, that a certain value of $x$ leads to a smallest quantity of $y$.
 
\subsection*{Example 3}

\textit{To find the cases in which this expression $x^3-axx+bx-c$ takes on the maximum or minimum value.}\\[2mm]
Having put $y=x^3-axx+bx-c$, it will be

\begin{equation*}
\frac{dy}{dx}=3xx-2ax+b \quad \text{and} \quad \frac{ddy}{2dx^2}=3x-a, \quad \frac{d^3y}{6dx^3}=1.
\end{equation*}
Therefore, set $\frac{dy}{dx}=3xx-2ax+b=0$; it will be

\begin{equation*}
x= \frac{a \pm \sqrt{aa-3b}}{3},
\end{equation*}
from which it it clear, if it is not $aa>3b$, that the propounded formula will neither have a maximum nor a minimum value. Therefore, this equation results

\begin{equation*}
\frac{ddy}{2dx^2}= \pm \sqrt{aa-3bb},
\end{equation*}
whence it is understood, if it not $aa=3b$, that the value $x= \frac{a+\sqrt{aa-3b}}{3}$ renders the formula $y=x^3-axx+bx-c$ minimal, the other value $x= \frac{a- \sqrt{aa-3b}}{3}$ on the other hand renders it maximal. But how large will these value of $y$ be? Since $3xx-2ax+b=
0$ or $x^3-\frac{2}{3}axx+\frac{1}{3}bx=0$, it will be

\begin{equation*}
y=- \frac{1}{3}axx+\frac{2}{3}bx-c
\end{equation*} 
and, because of $\frac{1}{3}axx-\frac{2aa}{9}x+\frac{ab}{9}=0$, 

\begin{equation*}
y=\frac{2}{9}(3b-aa)x+\frac{ab}{9}-c= - \frac{2a(aa-3b)}{27}-\mp \frac{2(aa-3b)\sqrt{aa-3b}}{27}+\frac{ab}{9}-c
\end{equation*}
or

\begin{equation*}
y=- \frac{2a^3}{27}+\frac{ab}{3}-c \mp \frac{2}{27}(aa-3b)^{\frac{3}{2}},
\end{equation*}
where the upper sign holds for the minimum value, but the lower sign for the maximum value.\\[2mm]
Therefore, the case $aa=3b$ remains; because in that case $\frac{ddy}{dx^2}=0$, but the following term $\frac{d^3y}{6dx^3}=1$ is not $=0$, it follows that in this case the propounded formula has neither a maximum nor a minimum value.

\subsection*{Example 4}

\textit{To find the cases in which this function of $x$, $x^4-8x^3+22x^2-24x+12$ has a maximum or minimum value.} \\[2mm]
Having put $y=x^4-8x^3+22x^2-24x+12$, it will be

\begin{equation*}
\frac{dy}{dx}=4x^3-24x^2+44x-24, \quad \frac{ddy}{2dx^2}=6x^2-24x+22.
\end{equation*}
Now set

\begin{equation*}
\frac{dy}{dx}=4x^3-24x^2+44x-24=0 \quad \text{or} \quad x^3-6x^2+11x-6=0;
\end{equation*}
three real values for $x$ are found, i.e.

\begin{equation*}
\text{I.} ~~x=1, \quad \text{II.}~~ x=2, \quad \text{III.} \quad x=3.
\end{equation*}
From the first value $\frac{ddy}{2dx^2}=4$ and hence for $x=1$ the propounded function has a minimum value. From the second value $x=2$ we find $\frac{ddy}{2dx^2}=-2$ and hence the propounded function has a maximum value. From the third value $x=2$ we deduce $\frac{ddy}{2dx^2}=+4$ and hence the propounded function has  a minimum value again.

\subsection*{Example 5}

\textit{Let this function be propounded $y=x^5-5x^4+5x^3+1$; it is in question in which cases it has a maximum or a minimum value.}\\[2mm]
Since

\begin{equation*}
\frac{dy}{dx}=5x^4-20x^3+15xx,
\end{equation*}
form the equation $x^4-4x^3+3xx=0$, whose roots are

\begin{equation*}
\text{I.}~~\text{and} \quad \text{II.} ~~x=0, \quad \text{III.}~~x=1, \quad \text{IV.}~~ x=3.
\end{equation*}
Since the first and second root are the same, from them neither a maximum nor a minimum value follows; for,  $\frac{ddy}{dx^2}=0$, but $\frac{d^3y}{dx^3}$ does not vanish. But the third root $x=1$, because of $\frac{ddy}{2dx^2}=10x^3-30x^2+15x$, yields $\frac{ddy}{2dx^2}=-5$ and in this case the function takes on a maximum value. From the fourth root $x=3$ we have $\frac{ddy}{2dx^2}=45$ and hence the propounded function has minimum value.

\subsection*{Example 6}

\textit{To find the cases in which this formula $y=10x^6-12x^5+15x^4-20x^3+20$ has a maximum or a minimum value.}\\[2mm]
Therefore, it will be

\begin{equation*}
\frac{dy}{dx}=60x^5-60x^4+60x^3-60x^2 \quad \text{and} \quad \frac{ddy}{60dx^2}= 5x^4-4x^3++3x^2-2x.
\end{equation*}
Form the equation $x^5-x^4+x^3-xx=0$; since, having resolved it into factors,  $x^2(x-1)(xx+1)=0$, it has two equal roots $x=0$ and furthermore the root $x=1$ and additionally two imaginary ones from $xx+1=0$. Therefore, since the two equal roots $x=0$ exhibit neither a maximum nor a minimum, only the root $x=1$ is to be considered, from which  $\frac{ddy}{60dx^2}=2$, whose positive value indicates a minimum value.

\paragraph*{§262}

Therefore, the determination of maxima and minima depends on the roots of the equation $\frac{dy}{dx}=0$; because its highest power is one degree lower than the highest power in the propounded equation, it is obvious, if in general this function is propounded

\begin{equation*}
x^n+Ax^{n-1}+Bx^{n-2}+Cx^{n-3}+Dx^{n-4}+\text{etc.}=y,
\end{equation*}
that its maxima and minima are determined by means of the roots of this equation

\begin{equation*}
nx^{n-1}+(n-1)Ax^{n-2}+(n-2)Bx^{n-3}+(n-3)Cx^{n-4}+\text{etc.}=0.
\end{equation*}
Let us put that the real roots of this equation ordered according to  their magnitude are $\alpha$, $\beta$, $\gamma$, $\delta$ etc. such that $\alpha$ is the largest, $\beta <\alpha$, $\gamma < \beta$ etc. And first, if these roots are all different, each of them will lead to a maximum or minimum value of the propounded formula $y$ and hence the function $y$ will have as many maxima or minima as the  equation $\frac{dy}{dx}=0$ had real different roots. But if two or more  roots were equal to each other, it will be as follows: two equal roots will exhibit neither a maximum nor a minimum value, but three on the other hand will be equivalent to a single one; and in general, if the number of equal roots was an even number, hence no maximum nor minimum results; but if the number is odd, one maximum or minimum results from this.

\paragraph*{§263}

But which roots produce maxima and which produce minima,  can be defined without using the rule given before as follows. Since the function $y$ for $x= \infty$ also becomes infinite and  values of $x$ within the boundaries $\infty$ and $\alpha$ produce neither a maximum nor a minimum value, it is perspicuous that the values of the function $y$, if  successively all values from $\infty$ up to $\alpha$ are substituted for $x$, have to decrease continuously; and hence the value $x= \alpha +\omega$  will lead to a larger value of the function $y$ than the value $x=\alpha$; therefore, because $x= \alpha$ produces a maximum or minimum, it is necessary that in this case the function $y$ takes on a minimum value. Therefore, decreasing $x$ or putting $x=\alpha -\omega$, the value of $y$ will increase again, until  finally  $x= \beta$, which is the second root  of the equation $\frac{dy}{dx}=0$ producing a maximum or minimum; therefore,  this second root $x=\beta$ will yield a maximum  and the value $x=\beta -\omega$ will cause the function $y$ to be smaller than for $x=\beta$, until one gets to $x=\gamma$, which as a logical consequence will generate a minimum value again. From this reasoning it is understood that the first, third, fifth etc. root of the equation $\frac{dy}{dx}=0$ will exhibit minima but the first, second, fourth, sixth etc. exhibit maxima. But at the same time it is hence understood that in the case of two equal roots a maximum and minimum coalesce and so none of both actually occurs.

\paragraph*{§264}

Therefore, if in the propounded function

\begin{equation*}
y=x^n+Ax^{n-1}+Bx^{n-2}+Cx^{n-3}+\text{etc.}
\end{equation*}
the greatest exponent $n$ was an even number, the equation

\begin{equation*}
\frac{dy}{dx}=x^{n-1}+(n-1)Ax^{n-2}+\text{etc.}=0
\end{equation*}
will be of odd degree and will hence have one or three or five or any odd number of real roots. If just one root was real, it will give a minimum; if three were real, the largest will yield a minimum, the middle one a maximum point and the smallest a minimum; and if five roots were real, the function $y$ will have three minimum values and two maximum values; and so forth.\\[2mm]
But if the exponent $n$ was an odd number, the equation $\frac{dy}{dx}=0$ will have an odd degree and will have either no or two or four or six etc. real roots. In the first case the function $y$ will have neither a maximum nor a minimum value; in the other case, in which two roots are given, the greater one will indicate a minimum value, the smaller a maximum value; but the first  (which is the largest) and the third of four roots will produce a minimum value, the second and the fourth on the other hand a maximum value. But no matter how many real roots were there, they will give maxima and minima alternately.

\paragraph*{§265}

Let us proceed to rational functions which constitute the other kind of uniform functions. Therefore, let

\begin{equation*}
y=\frac{P}{Q}
\end{equation*}
where $P$ and $Q$ are any polynomial functions of $x$; and at first it is certainly clear, if  a value of such a kind is attributed to $x$ that $Q=0$, if not at the same time $P$ vanishes, that the function $y$  becomes infinite, what appears to be a maximum value. Nevertheless, this case cannot be considered as a maximum value; for, because the inverse fraction $\frac{Q}{P}$  takes on  a minimum value in the same cases the propounded $\frac{P}{Q}$ takes on a maximum value, the fraction $\frac{Q}{P}$ would have to take on a minimum value, if $Q$ vanishes; but this does not always happen, since even smaller values, i.e. negative values,  could occur. Therefore, this at the same time confirms the rule given before, that maxima and minima must be found from the equation $\frac{dy}{dx}=0$. Therefore, in the propounded cases it will be

\begin{equation*}
\frac{dy}{dx}=\frac{QdP-PdQ}{QQdx}
\end{equation*}
and hence the roots of this equation
\begin{equation*}
QdP-PdQ=0
\end{equation*}
will give the a maximum or a minimum values of the function $y$. And if there is any doubt, whether $y$ takes on a maximum or a minimum value, one has to check the value $\frac{ddy}{dx^2}$; if it was positive, it will indicate a minimum value, but if it was negative a minimum value. If also this value $\frac{ddy}{dx^2}$ vanishes, which happens, if the equation $\frac{dy}{dx}=0$ has two or more equal roots, just recall that an equal number of equal roots produce neither a maximum nor a minimum value.

\subsection*{Example 1}

\textit{To find the cases in which the function $\frac{x}{1+xx}$ takes on a maximum or a minimum value.}\\[2mm]
At first it is certainly clear that this function goes over into zero in the three cases $x=\infty$, $x=0$ and $x=-\infty$, whence it will have at least either two maximum or two minimum values. To find them put $y=\frac{x}{1+xx}$ and it will be

\begin{equation*}
\frac{dy}{dx}= \frac{1-xx}{(1+xx)^2} \quad \text{and} \quad \frac{ddy}{dx^2}=\frac{-6x+2x^3}{(1+xx)^3}.
\end{equation*}
Now set $\frac{dy}{dx}=0$; it will be $1-xx=0$ and either $x=+1$ or $x=-1$. In the first case $x=+1$ we have $\frac{ddy}{dx^2}=-\frac{4}{2^3}$ and hence $y$ has a maximum value $=\frac{1}{2}$; in the second case $x=-1$ we have $\frac{ddy}{dx^2}=+\frac{4}{2^3}$ and hence $y$ has a minimum value $=-\frac{1}{2}$.\\[2mm]
These are also found more easily, if the propounded fraction $\frac{x}{1+xx}$ is inverted or by putting $y=\frac{1+xx}{x}=x+\frac{1}{x}$, if we recall that then all values, which were found to be maximum values,  are to be turned into minimum values and vice versa. But it will be

\begin{equation*}
\frac{dy}{dx}=1-\frac{1}{xx} \quad \text{and} \quad \frac{ddy}{dx^2}=\frac{2}{x^3}.
\end{equation*}
Therefore, having set $\frac{dy}{dx}=0$, $xx-1=0$ and hence either $x=+1$ or $x=-1$ as before. And in the case $x=+1$ we have $\frac{ddy}{dx^2}=2$ and hence $y$ has a minimum value and the propounded formula $\frac{1}{y}$ a maximum value. But in the case $x=-1$ we find $\frac{ddy}{dx^2}=-2$, whence $y$ has a maximum value and $\frac{1}{y}$ a minimum value.

\subsection*{Example 2}

\textit{To find the cases in which the formula $\frac{2-3x+xx}{2+3x+xx}$ has a maximum or a minimum value.}\\[2mm]
Having put $y= \frac{xx-3x+2}{xx+3x+2}$, it will be

\begin{equation*}
\frac{dy}{dx}=\frac{6x^2-12}{(xx+3x+2)^2} \quad \text{and} \quad \frac{ddy}{dx^2}=\frac{-12x^3+72x+72}{(xx+3x+2)^3}.
\end{equation*}
Set $\frac{dy}{dx}=0$; it will be either $x=+\sqrt{2}$ or $x=-\sqrt{2}$. In the first case it will be 

\begin{equation*}
\frac{ddy}{dx^2}=\frac{48\sqrt{2}+72}{(4+3\sqrt{2})^3}
\end{equation*}
and hence will be positive, because of the positive numerator and denominator; Therefore, $y$ will take on a minimum value

\begin{equation*}
= \frac{4-\sqrt{3}2}{4+\sqrt{3}2}=12\sqrt{2}-17=-0.02943725.
\end{equation*}
In the second case $x=-\sqrt{2}$ 

\begin{equation*}
\frac{ddy}{dx^2}=\frac{-48\sqrt{2}+72}{(4-3\sqrt{2})^3}=\frac{24(3-2\sqrt{2})}{(4-3\sqrt{2})^3},
\end{equation*}
whose value, because of the affirmative numerator and negative denominator, will be negative and hence $y$ will take on a maximum value

\begin{equation*}
=\frac{4+3\sqrt{2}}{4-3\sqrt{2}}=-12\sqrt{2}-17=-33.97056274.
\end{equation*}
This value, even though it is smaller than the first minimum value, is nevertheless a maximum value, since it is larger than the values in its neighborhood, which result, if  a values little bit greater or smaller than $-\sqrt{2}$ are substituted for $x$. Therefore, because $\sqrt{2}$ is contained within the limits $\frac{4}{3}$ and  $\frac{3}{2}$, the crosscheck will easily be done this way:

\begin{center}
\renewcommand{\arraystretch}{1,8}
\begin{tabular}{llcllll}
if & $x=$ & $\dfrac{4}{3}$, & we have & $y=$ & $- \dfrac{2}{70}=-0.0285$ & \\
if & $x=$ & $\sqrt{2}$, & we have & $y=$ & $+ 12\sqrt{2}-17=0.0294$ &\quad minimum \\
if & $x=$ & $\dfrac{3}{2}$, & we have & $y=$ & $- \dfrac{1}{35}=-0.0285$ & \\ 
 &  &  &  &  &  &\\ \cline{3-6} \\
if & $x=$ & $-\dfrac{4}{3}$, & we have & $y=$ & $-35$ & \\
if & $x=$ & $-\sqrt{2}$, & we have & $y=$ & $-33.970$ &\quad maximum \\
if & $x=$ & $-\dfrac{3}{2}$, & we have & $y=$ & $- 35.$ & \\ 
\end{tabular}
\end{center}
 
\subsection*{Example 3}

\textit{To find the cases in which the formula $\frac{xx-x+1}{xx+x-1}$ has a maximum or a minimum value.}\\[2mm]
Put $y=\frac{xx-x+1}{xx+x-1}$ and it will be

\begin{equation*}
\frac{dy}{dx}=\frac{2xx-4x}{(xx+x-1)^2} \quad \text{and} \quad \frac{ddy}{dx^2}=\frac{-4x^3+12xx+4}{(xx+x-1)^3}.
\end{equation*}
Set $\frac{dy}{dx}=0$; it will be either $x=0$ or $x=1$; in the first case  $\frac{ddy}{dx^2}=\frac{4}{-1}$ and hence $y$ will have a maximum value $=-1$. In the second case  $\frac{ddy}{dx^2}=\frac{20}{5^2}$ and hence $y$ has the minimum value $=\frac{3}{5}$, even though that maximum value is smaller than this minimum value. That this is indeed the case is seen from the following calculations

\begin{center}
\renewcommand{\arraystretch}{2,0}
\begin{tabular}{lllllll}
if & $x=$ & $-\dfrac{1}{3}$, & it will be & $y=$ & $- \dfrac{13}{11}$ & \\
if & $x=$ & $+0$, & it will be & $y=$ & $-1$ &\quad maximum \\
if & $x=$ & $+\dfrac{1}{3}$, & it will be & $y=$ & $- \dfrac{7}{5}$ & \\ 
 &  &  &  &  &  &\\ \cline{3-6}
 &  &  &  &  &  &\\
 if & $x=$ & $2-\dfrac{1}{3}$, & it will be & $y=$ & $- \dfrac{19}{13}$ & \\
if & $x=$ & $2$, & it will be & $y=$ & $+\dfrac{3}{5}$ &\quad minimum \\
if & $x=$ & $2+\dfrac{1}{3}$, & it will be & $y=$ & $- \dfrac{37}{61}$. & \\ 
\end{tabular}
\end{center}
But that, if one puts $x=1$,  $y=1$ and hence $>-1$, is the reason why  between the values  $0$ and $1$ there is one value of $x$, for which $y=\infty$.

\subsection*{Example 4}

\textit{To find the the cases in which this fraction $\frac{x^3+x}{x^4-xx+1}$ has a maximum or a minimum value.}\\[2mm]
Having put $y= \frac{x^3+x}{x^4-xx+1}$, it will be

\begin{equation*}
\frac{dy}{dx}=\frac{-x^6-4x^4+4xx+1}{(x^4-xx+1)^2} \quad \text{and} \quad \frac{ddy}{dx^2}=\frac{2x^9+18x^7-30x^5-16x^3+12x}{(x^4-xx+1)^3}.
\end{equation*}
Therefore, we will have this equation

\begin{equation*}
x^6+4x^4-4xx-1=0,
\end{equation*}
which is resolved into these two

\begin{equation*}
xx-1=0 \quad \text{and} \quad x^4+5x^2+1=0;
\end{equation*}
the roots of the first  of these equations are $x=+1$ and $x=-1$, the other  gives $x=-\frac{5\pm \sqrt{21}}{2}$,  from which no real root emerges. Therefore, the first of the two roots, $x=+1$, gives $\frac{ddy}{dx^2}=-14$ and therefore $y$ has a maximum value $=2$; the other root $x=-1$ gives $\frac{ddy}{dx^2}=+14$ and therefore $y$ has a minimum value $=-2$.

\subsection*{Example 5}

\textit{To find the cases in which this fraction $\frac{x^3-x}{x^4-xx+1}$ has a maximum or a minimum value.}\\[2mm]
Having put $y= \frac{x^3-x}{x^4-xx+1}$, it will be

\begin{equation*}
\frac{dy}{dx}=\frac{-x^6+2x^4+2x^2-1}{(x^4-x^2+1)^2} \quad \text{and} \quad \frac{ddy}{dx^2}=\frac{2x^9-6x^7-18x^5+10x^3}{(x^4-x^2+1)^3}.
\end{equation*}
But having put $\frac{dy}{dx}=0$, it will be

\begin{equation*}
x^6-2x^4-2x^2+1=0,
\end{equation*}
which divided by $xx+1$ gives

\begin{equation*}
x^4-3x^2+1=0,
\end{equation*}
and this is further resolved into

\begin{equation*}
xx-x-1=0 \quad \text{and} \quad xx+x-1=0,
\end{equation*}
whence the following four real roots result

\begin{alignat*}{19}
&\text{I.} \quad && x&&~=~&& &&\frac{1+\sqrt{5}}{2} \qquad && \text{II.} \quad && x &&~=~ && &&\frac{1-\sqrt{5}}{2}, \\[3mm]
&\text{III.} \quad && x&&~=~&&- &&\frac{1+\sqrt{5}}{2} \qquad && \text{IV.} \quad && x &&~=~ &&-&&\frac{1-\sqrt{5}}{2}.
\end{alignat*}
Since all are contained in the equation $x^4-3xx+1=0$, having put $x^4=3xx-1$, the following equations will hold for all of them

\begin{equation*}
\frac{ddy}{dx^2}=\frac{2x(10-20xx)}{8x^6}=\frac{5(1-2xx)}{2x^5}=\frac{5(1-2xx)}{2x(3xx-1)} \quad \text{and} \quad y= \frac{x^3-x}{2xx}=\frac{xx-1}{2x}.
\end{equation*}
But, for the first two, resulting from the equation $xx=x+1$, it will be

\begin{equation*}
\frac{ddy}{dx^2}=-\frac{5(2x+1)}{2x(3x+2)}=-\frac{5(2x+1)}{2(5x+3)} \quad \text{and} \quad y=\frac{1}{2}.
\end{equation*}
Therefore, the first root $x=\frac{1+\sqrt{5}}{2}$ gives

\begin{equation*}
\frac{ddy}{dx^2}=-\frac{5(2+\sqrt{5})}{11+5\sqrt{5}}
\end{equation*}
and hence $y$ has a maximum value. The second root $x=\frac{1-\sqrt{5}}{2}$ gives

\begin{equation*}
\frac{ddy}{dx^2}=-\frac{5(2-\sqrt{5})}{11-5\sqrt{5}}=- \frac{5(\sqrt{5}-2)}{5\sqrt{5}-11}
\end{equation*}
and hence $y=\frac{1}{2}$ will also have a maximum value. The two remaining roots give $y=-\frac{1}{2}$, a minimum value.

\paragraph*{§266}

Therefore, in these examples the exploration, whether a certain found value produces a maximum or a minimum, can be simplified; for, because $\frac{dy}{dx}=0$, the value of the term $\frac{ddy}{dx^2}$, having taken into account its equation, can be expressed in an easier way. For, let the fraction $y=\frac{P}{Q}$ be propounded; since

\begin{equation*}
dy= \frac{QdP-PdQ}{QQ} \quad \text{and} \quad QdP-PdQ=0,
\end{equation*}
it will be

\begin{equation*}
ddy=\frac{d(QdP-PdQ)}{Q^2}-\frac{2dQ(QdP-PdQ)}{Q^3}.
\end{equation*}
But, because of $QdP-PdQ=0$, this last term vanishes and it will be

\begin{equation*}
ddy= \frac{d(QdP-PdQ)}{QQ}=\frac{QddP-PddQ}{Q^2}.
\end{equation*}
But because the decision is to be made from the either positive or negative value of this term, and the denominator $Q^2$ is always positive, this  can be done considering only the numerator  in such a way, that, if $QddP-PddQ$ or $\frac{d(QdP-PdQ)}{dx^2}$ was positive, a minimum will be indicated,  if it is negative, a maximum. Or after $\frac{dy}{dx}$ was found, whose form will be of this kind $\frac{R}{QQ}$, only find $\frac{dR}{dx}$, and from the root, which causes this expression to be to positive, a minimum value will result and otherwise a maximum value.

\paragraph*{§267}

If the denominator of the propounded fraction was a square or any higher power such that $y=\frac{P}{Q^n}$, it will be

\begin{equation*}
dy=\frac{QdP-nPdQ}{Q^{n+1}}
\end{equation*}
and, having put $\frac{QdP-nPdQ}{dx}=R$, it will be

\begin{equation*}
\frac{dy}{dx}=\frac{R}{Q^{n+1}}
\end{equation*}
and the maximum and minimum values will be determined from the roots of the equation $R=0$. Further, since

\begin{equation*}
\frac{ddy}{dx}=\frac{QdR-(n+1)RdQ}{Q^{n+2}},
\end{equation*}
because of $R=0$, it will be

\begin{equation*}
\frac{ddy}{dx}=\frac{dR}{Q^{n+1}};
\end{equation*}
its positive value will indicate a minimum value, a negative a maximum value. But it is perspicuous, if $n$ was an odd number, that, because of the always positive $Q^{n+1}$, the decision can be made  considering only $\frac{dR}{dx}$; but if $n$ is an even number, use the formula $\frac{QdR}{dx}$.\\[2mm]
But let us further put that a fraction of this kind is propounded $\frac{P^m}{Q^n}=y$; it will be

\begin{equation*}
dy=\frac{(mQdP-nPdQ)P^{m-1}}{Q^{n+1}};
\end{equation*}
therefore, if one puts $\frac{mQdP-nPdQ}{dx}=R$, the roots of the equation $R=0$ will indicate the cases, in which the function $y$ has a maximum or a minimum value. Therefore, since

\begin{equation*}
\frac{dy}{dx}=\frac{P^{m-1}R}{Q^{n+1}},
\end{equation*}
it will be

\begin{equation*}
\frac{ddy}{dx}=\frac{P^{m-2}R((m-1)QdP-(n+1)PdQ)}{Q^{n+2}}+\frac{P^{m-1}dR}{Q^{n+1}},
\end{equation*}
and, because of $R=0$, it will be

\begin{equation*}
\frac{ddy}{dx^2}=\frac{P^{m-1}dR}{Q^{n+1}dx};
\end{equation*}
this can additionally be divided by any square $\frac{P^{2 \mu}}{Q^{2 \nu}}$ to  make the decision. Furthermore, also the equation $P=0$ will give a maximum or a minimum point, if $m$ was an even number; and in like manner  considering the inverse formula $\frac{Q^n}{P^m}$ a maximum or minimum point will result by putting $Q=0$, if $n$ was an even number, as we showed above (§ 257); but here we do not consider the maximum or minimum  values to result from this, but only, in order to explain the use of the method, investigate those, which result from the equation $R=0$.

\subsection*{Example 1}

\textit{Let the fraction $\frac{(\alpha +\beta x)^m}{(\gamma +\delta x)^n}$ be propounded; in which case it takes on a minimum or a maximum value, is in question.}\\[2mm]
Having put $y=\frac{(\alpha +\beta x)^m}{(\gamma +\delta x)^n}$,  first it is certainly clear that it will be $y=0$, if $x=-\frac{\alpha}{\beta}$, and  $y= \infty$, if $x=-\frac{\gamma}{\delta}$; the latter of these cases will give a minimum value, the first a maximum value, if $m$ and $n$ were even numbers. Furthermore, it will be

\begin{equation*}
\frac{dy}{dx}=\frac{(\alpha +\beta x)^{m-1}}{(\gamma +\delta x)^{n+1}}((m-n)\beta \delta x+ m \beta \gamma - n \alpha \delta)
\end{equation*}
and hence

\begin{equation*}
R=(m-n)\beta \delta x+m \beta \gamma - n \alpha \delta.
\end{equation*}
Therefore, having put $R=0$, it will be 

\begin{equation*}
x= \frac{n \alpha \delta - m \beta \gamma}{(m-n)\beta \delta}.
\end{equation*}
Further, because of $\frac{dR}{dx}=(m-n)\beta \delta$, it is perspicuous, whether

\begin{equation*}
\frac{P^{m-1}dR}{Q^{n+1}dx}=\frac{m^{m-1}\beta^{n+1}}{n^{n+1} \delta^{m-1}}\left(\frac{\alpha \delta - \beta \gamma}{m-n}\right)^{m-n-2}\frac{dR}{dx}
\end{equation*}
is a positive or negative quantity. In the first case, the propounded formula will have minimum value, in the second a maximum value. \\[2mm]
So, if it was $y=\frac{(x+3)^3}{(x+2)^2}$, it will be $\frac{P^{m+1}dR}{Q^{n+1}dx}=\frac{9}{8}$ and hence the formula $\frac{(x+3)^3}{(x+2)^2}$
will have a minimum value for $x=0$.\\[2mm]
But if $y= \frac{(x-1)^m}{(x+1)^m}$, it will be

\begin{equation*}
\frac{P^{m-1}dR}{Q^{n+1}dx}=\frac{m^{m-1}}{n^{n+1}}\left(\frac{n-m}{2}\right)^{n-m+2}(m-n)
\end{equation*}
and $x=\frac{n+m}{n-m}$. But because $m$ and $n$ are put to be positive numbers, the decision is  to be made from the formula $(n-m)^{n-m+2}(m-n)$ or $(n-m)^{n-m}(m-n)$. Therefore, if it was $n>m$, the found value $x=\frac{n+m}{n-m}$ will always give a maximum value; but if $n<m$, the number $m-n$ will give a minimum value, but an odd a maximum value; so $\frac{(x-1)^3}{(x+1)^2}$ will have a maximum value for $x=-5$; for, $y=-\frac{6^3}{4^2}=-\frac{27}{2}$.

\subsection*{Example 2}

\textit{Let the formula $y=\frac{(1+x)^3}{(1+xx)^2}$ be propounded.}\\[2mm]
It will be

\begin{equation*}
\frac{dy}{dx}=\frac{(1+x)^2}{(1+xx)^3}(3-4x-xx) \quad \text{and} \quad \frac{P^{m-1}}{Q^{n+1}}\cdot \frac{dR}{dx}=-\frac{(1+x)^2}{(1+xx)^3}(2x+4);
\end{equation*}
because here $(1+x)^2$ and $(1+xx)^3$ always have a positive value, the decision is to be made from the formula $-x-2$;  if it was positive, it indicates a minimum value, if negative, a maximum value. But from the equation $3-4x-xx=0$ it follows that either

\begin{equation*}
x=-2+\sqrt{7} \quad \text{or} \quad x=-2-\sqrt{7}.
\end{equation*}
In the first case $-x-2=-\sqrt{7}$ and hence the propounded fraction will have a maximum value, in the other case a minimum value, because of $-x-2=+\sqrt{7}$. But having put $x=-2+\sqrt{7}$, it will be $1+x=-1+\sqrt{7}$ and $1+xx=12-4\sqrt{7}$, whence

\begin{equation*}
y= \left(\frac{-1+\sqrt{7}}{12-4\sqrt{7}}\right)^2 (\sqrt{7}-1)=\frac{(2+\sqrt{7})^2(\sqrt{7}-1)}{16}=\frac{17+7\sqrt{7}}{16}=2.220.
\end{equation*}
But having put $x=-2-\sqrt{7}$, it will be

\begin{equation*}
y= \frac{17-7\sqrt{7}}{16}=-0.0950.
\end{equation*}

\paragraph*{§268}

There are also irrational and transcendental functions  which have the property of uniform functions, and therefore their maxima and minima can be found the same way. For, the cube and all odd roots are indeed uniform, since they exhibit only one single real value; but even though square roots and roots of all even powers, if they are real,  actually denote two values, the one positive, the other negative,   each of them can nevertheless be considered separately and in this sense one can even investigate the maxima and minima. So, if $y$ was any function of $x$, even though $\sqrt{y}$ takes on two values, one can nevertheless treat each one separately. $+\sqrt{y}$  will have a maximum or minimum value, if $y$  had  one, if it was affirmative, since otherwise $\sqrt{y}$ would become imaginary. But vice versa $-\sqrt{y}$ will have a maximum or minimum value in the same cases, in which $+\sqrt{y}$ has a maximum or minimum value. But any power $y^{\frac{m}{n}}$ takes on a maximum or minimum value in the same cases, if $n$ was an odd number; but if $n$ was an even number, only those cases remain, in which $y$ has a positive value, and in these cases, because of the  two values, two maximum or minimum values will result.

\paragraph*{§269}

Since the differential equation, which results from the power of the function $y^m$, is $\frac{y^{m-1}dy}{dx}=0$, whose roots  at the same time indicate the cases, in which a surdic power $y^{\frac{m}{n}}$ has a maximum or minimum value, to investigate this value one has two equations, the one $y^{m-1}=0$, the other $\frac{dy}{dx}=0$, the latter of which goes over into  $y=0$ and exhibits maxima and minima only, if $m-1$ was an odd number or if $m$ was an even number, for the reasons mentioned in § 257. Therefore, because $n$ is an odd number, if $m$ was an even number, if we denote the even numbers by $2\mu$ and the odd numbers by $2\nu -1$, the function $y^{2 \mu : (2 \nu-1)}$ will have a maximum or minimum value, if those values are attributed to $x$ which  are found so from the equation $y=0$ as from this $\frac{dy}{dx}=0$. But if $m$ is an odd number, the function $y^{(2 \mu -1):2 \nu}$ or $y^{(2 \mu -1):(2 \nu)}$ has a maximum or minimum value only, if a value resulting from this equation $\frac{dy}{dx}=0$ is substituted for $x$. And in the second case $y^{(2\mu -1):2 \nu}$ maxima and minima only result, if $y$  has positive values for the values found from the equation $\frac{dy}{dx}=0$.

\paragraph*{§270}

So, this formula $x^{\frac{2}{3}}$  takes on a minimum value for $x=0$, because in this case $x^2$  has minimum value. But if we do not reduce the formula $x^{\frac{2}{3}}$ to the form $x^2$, the method given before would not indicate this at all, since in the case $x=0$ the terms of the series

\begin{equation*}
y+\frac{\omega dy}{dx}+\frac{\omega^2 ddy}{2dx^2}+\frac{\omega^3 d^3 y}{6dx^3}+\text{etc.},
\end{equation*}
whence the decision is to be made, except for the first all become infinite. For, having put $y=x^{\frac{2}{3}}$, it will be

\begin{equation*}
\frac{dy}{dx}=\frac{2}{3x^{\frac{1}{3}}}, \quad \frac{ddy}{dx^2}=\frac{-2}{9x^{\frac{4}{3}}}, \quad \frac{d^3y}{dx^3}=\frac{2 \cdot 4}{27 x^{\frac{7}{3}}} \quad \text{etc.}
\end{equation*}
And hence neither the equation $\frac{dy}{dx}=\frac{2}{3x^{\frac{1}{3}}}=0$  shows the value $x=0$ nor the following terms indicate whether it is a maximum or a minimum value. Therefore, since we assumed that the series

\begin{equation*}
y+\frac{\omega dy}{dx}+\frac{\omega^2 ddy}{2dx^2}+\frac{\omega^3 d^3y}{6dx^3}+\text{etc.}
\end{equation*}
becomes convergent, if $\omega$ is assumed to be a very small quantity, in those cases, in which this series becomes divergent, the general method is not applicable, what happens in the example $y=x^{\frac{2}{3}}$ mentioned here, if one puts $x=0$. Therefore, in these cases the same reduction we used before will be necessary to reduce the propounded expression to another form, which is not subjected to this inconvenience. But this only happens in very few cases  which are contained in the formula $y^{\frac{2 \mu}{2 \nu -1}}$ or are easily reduced to it. So, if the maxima and minima of the formula $y^{\frac{2 \mu}{2 \nu -1}}z$ are in question, where $z$ is any function of $x$, investigate this form $y^{2 \mu}z^{2 \nu -1}$, which  has a maximum or minimum value in the same cases as the propounded one.

\paragraph*{§271}

Having excluded this case, which is now easily handled, functions containing irrational quantities, can be treated the same way as rational functions and their maxima and minima can be determined, what we will illustrate in the following examples.

\subsection*{Example 1}

\textit{Let the formula $\sqrt{aa+xx}-x$ be propounded; in which cases it has   maximum or minimum values is in question.}\\[2mm]
Having put $y=\sqrt{aa+xx}-x$, it will be

\begin{equation*}
\frac{dy}{dx}=\frac{x}{\sqrt{aa+xx}}-1 \quad \text{and} \quad \frac{ddy}{dx^2}=\frac{aa}{(aa+xx)^{3:2}}.
\end{equation*}
Having put $\frac{dy}{dx}=0$, it will be $x=\sqrt{aa+xx}$ and hence $x=\infty$ and  $\frac{ddy}{dx^2}=0$. In a like manner the following terms $\frac{d^3y}{dx^3}$, $\frac{d^4y}{dx^4}$ etc. all become $=0$; hence one cannot decide, whether it is a maximum or a minimum. The reason for this is that it actually is so $x=-\infty$ as $x= + \infty$. Putting $x=\infty$, because of

\begin{equation*}
\sqrt{aa+xx}=x+\frac{aa}{2x},
\end{equation*} 
we have $y=0$, which value is the smallest of all.

\subsection*{Example 2}

\textit{Let the cases be in question in which this form $\sqrt{aa+2bx+mxx}-nx$ takes on a maximum or a minimum value.}\\[2mm]
Having put $y=\sqrt{aa+2bx+mxx}-nx$, it will be

\begin{equation*}
\frac{dy}{dx}=\frac{b+mx}{\sqrt{aa+2bx+mxx}}-n;
\end{equation*}
having put this $=0$, it will be

\begin{equation*}
bb+2mbx+mmxx=nnaa+2nnbx+mnnxx
\end{equation*}
or

\begin{equation*}
xx=\frac{2bx(nn-m)+nnaa-bb}{mm-mnn}
\end{equation*}
and hence

\begin{equation*}
x= \frac{(nn-m)b \pm \sqrt{mnn(m-nn)aa-nn(m-nn)bb}}{m(m-nn)}
\end{equation*}
or

\begin{equation*}
x= -\frac{b}{m}\pm \frac{n}{m}\sqrt{\frac{maa-bb}{m-nn}};
\end{equation*}
hence 

\begin{equation*}
\sqrt{aa+2bx+mxx}=\frac{b+mx}{n}= \pm \sqrt{\frac{maa-bb}{m-nn}}.
\end{equation*}
Therefore, since

\begin{equation*}
\frac{ddy}{dx^2}=\frac{maa-bb}{(aa+2bx+mxx)^{\frac{3}{2}}},
\end{equation*}
it will be

\begin{equation*}
\frac{ddy}{dx^2}=\frac{maa-bb}{\pm \left(\frac{maa-bb}{m-nn}\right)^{\frac{3}{2}}}= \frac{\pm (m-nn)\sqrt{m-nn}}{\sqrt{maa-bb}}.
\end{equation*}
Therefore, a maximum or minimum only exists, if $\frac{m-nn}{maa-bb}$ was a positive quantity. But if it is a positive quantity, the upper sign will give a minimum, if $m >nn$, a maximum on the other hand, if $m<nn$; the contrary happens, if the lower sign holds. Therefore, if $m=2$, $n=1$ and $b=0$, the formula $\sqrt{aa+2xx}-x$ has a minimum value for $x=+\frac{1}{2}\sqrt{2aa}=\frac{a}{\sqrt{2}}$, but a maximum value for $x=-\frac{a}{\sqrt{2}}$. Therefore, the minimum value will be $=a\sqrt{2}-\frac{a}{\sqrt{2}}=\frac{a}{\sqrt{2}}$ and the maximum value $=a \sqrt{2}+\frac{a}{\sqrt{2}}=\frac{3a}{\sqrt{2}}$.

\subsection*{Example 3}

\textit{To find the cases in which this expression $\sqrt[4]{1+mx^4}+\sqrt[4]{1-nx^4}$ has a maximum or a minimum value.}\\[2mm]
Because $\frac{dy}{dx}=\frac{mx^3}{(1+mx^4)^{\frac{3}{4}}}-\frac{nx^3}{(1-nx^4)^{\frac{3}{4}}}$, it will be

\begin{equation*}
mx^3(1-nx^4)^{\frac{3}{4}}=nx^3(1+mx^4)^{\frac{3}{4}} \quad \text{and hence} \quad m^4(1-nx^4)^3=n^4(1+mx^4)^3
\end{equation*}
or

\begin{equation*}
n^4-m^4+3mn(n^3+m^3)x^4+3m^2n^2(n^2-m^2)x^8+m^3n^3(n+m)x^{12}=0.
\end{equation*}
Therefore, only if this equation has a positive root for $x^4$, a maximum or minimum exists. For, this equation in general cannot be solved in a convenient manner, since it will be

\begin{equation*}
x^4= \frac{m^{\frac{4}{3}}-n^{\frac{4}{3}}}{mn(\sqrt[3]{m}+\sqrt[3]{n})} \quad \text{or} \quad x^4= \frac{m-\sqrt[3]{m^2n}+\sqrt[3]{mn^2}-n}{mn}
\end{equation*}
let us consider a special case and put $m=8n$, and then it will be

\begin{equation*}
-4095+24 \cdot 513 x^4-3 \cdot 63 \cdot 64 n^2 x^8+9 \cdot 512 x^{12}=0
\end{equation*}
or

\begin{equation*}
512 n^3 x^{12}-1344 n^2 x^8+1368 nx^4-455=0;
\end{equation*}
put $8nx^4=z$; it will be

\begin{equation*}
z^3-21z^2+171z-455=0,
\end{equation*}
which has the divisor $z-5$, and the other factor will be $zz-16zz+91=0$ containing imaginary roots. Therefore, it will only be $z=8nx^4=5$ and hence $x= \sqrt[4]{\frac{5}{8n}}$, which value is a maximum or minimum point of the expression $\sqrt[4]{1+8nx^4}+\sqrt[4]{1-nx^4}$. To find out which of both is the case, consider

\begin{equation*}
\frac{ddy}{dx^2}=\frac{3mxx}{(1+mx^4)^{\frac{7}{4}}}-\frac{3nxx}{(1-nx^4)^{\frac{7}{4}}}.
\end{equation*}
But, because of $m=8n$, having put $x^4=\frac{5}{8n}$, it will be

\begin{equation*}
\frac{ddy}{dx^2}=\left(\frac{24n}{6^{\frac{7}{4}}}-\frac{3n}{(3:8)^{\frac{7}{4}}}\right)xx=- \frac{360nxx}{6^{\frac{7}{4}}}
\end{equation*}
and hence negative; therefore, $\sqrt[4]{1+8nx^4}+\sqrt[4]{1-nx^4}$ will have a maximum value for $x= \sqrt[4]{\frac{5}{8n}}$. This maximum value will be $=\sqrt[4]{6}+\sqrt[4]{\frac{3}{8}}=\frac{3\sqrt[4]{6}}{2}$. If  we write $u$ instead of $nx^4$, it is plain that this expression $\sqrt[4]{1+8u}+\sqrt[4]{1-u}$ has a maximum value for $u=\frac{5}{8}$ and that this maximum value will be $=\frac{3 \sqrt[4]{6}}{2}=2.347627$. Therefore, whatever value except for $\frac{5}{8}$ is written for $u$, the expression will have a smaller value. 

\paragraph*{§272}

Maxima and minima will be determined the same way, even if transcendental quantities are contained in the propounded expression. For, if the propounded function was not multiform and one has not to consider several values of it at the same time, the roots of the differential equation will show maxima or minima, if those roots were not equal roots, whose total amount  is even. Therefore, we will demonstrate this investigation in several examples.

\subsection*{Example 1}

\textit{To find the number which has the smallest ratio to its logarithm.} \\[2mm]
A smallest ratio $\frac{x}{\log x}$ is obvious to exist, because this ratio becomes infinite  so for $x=1$ as for $x= \infty$. Therefore, vice versa the fraction $\frac{\log x}{x}$ will take on a maximum value somewhere; of course in the same case, in which $\frac{x}{\log x}$ has a minimum value. To find this case put $y= \frac{\log x}{x}$ and it will be

\begin{equation*}
\frac{dy}{dx}=\frac{1}{xx}-\frac{\log x}{xx}.
\end{equation*} 
Having put this equal to zero it will be $\log x=1$, and since  we assume the hyperbolic logarithm here, if $e$ is put for the number whose hyperbolic logarithm is $=1$, it will be $x=e$. Therefore, because all logarithms have a certain ratio to the hyperbolic ones,  $\frac{ e}{\log e}$ will also be a minimum point for the common logarithm or $\frac{\log e}{e}$ will be a maximum value. Since, assuming tabulated logarithms,  $\log e= 0.4342944819$, the fraction $\frac{\log x}{x}$ will always be smaller than $\frac{4342944819}{27182818284}$ or approximately  $\frac{47}{305}$ and no other number exists, which  has a smaller ratio to its logarithm than $305$ to $47$. That in this case $\frac{\log x}{x}$ has a maximum value is obvious, since, because of $\frac{dy}{dx}=\frac{1- \log x}{xx}$,

\begin{equation*}
\frac{ddy}{dx^2}=-\frac{1}{x^3}-\frac{2(1- \log x)}{x^3}=-\frac{1}{x^3}
\end{equation*} 
since $1-\log x =0$ and hence the differential is negative.

\subsection*{Example 2}

\textit{To find the number $x$ that this power $x^{1:x}$ takes on a maximum value.} \\[2mm]
That a maximal value of this formula exists is obvious, because by substituting numbers for $x$

\begin{alignat*}{9}
&1^{1:1} &&=1.000000 \\
&2^{1:2} &&=1.414213 \\
&3^{1:3} &&=1.442250 \\
&1^{1:4} &&=1.414213.
\end{alignat*}
Therefore, put $x^{1:x}=y$ and it will be

\begin{equation*}
\frac{dy}{dx}=x^{1:x}\left(\frac{1}{xx}-\frac{\log x}{xx}\right).
\end{equation*}
Having put this value equal to zero, it will be $\log x=1$ and $x=e$ where $e=2.718281828$. And since $\frac{dy}{dx}=(1- \log x)\frac{x^{1:x}}{xx}$, it will be

\begin{equation*}
\frac{ddy}{dx^2}=- \frac{x^{1:x}}{x^3}+(1- \log x) \frac{d}{dx}\cdot \frac{x^{1:x}}{xx}=-\frac{x^{1:x}}{x^3}
\end{equation*}
because $1- \log x=0$. Therefore, because $\frac{ddy}{dx^2}$ is a negative quantity, $x^{1:x}$ has a maximum value in the case $x=e$. But because $e=2.718281828$, one finds $e^{\frac{1}{e}}=1.444667861009764$, which value is easily obtained from the series

\begin{equation*}
e^{\frac{1}{e}}=1+\frac{1}{e}+\frac{1}{2e^2}+\frac{1}{6e^3}+\frac{1}{24e^4}+\text{etc.}
\end{equation*}
This example is also solved using the result of the preceding example; for, if $x^{1:x}$ takes on a maximum value, also its logarithm, which is $\frac{\log x}{x}$, will have to take on a maximum value; for  this to happen, it has to be $x=e$, as we found.

\subsection*{Example 3}

\textit{To find the arc $x$ that its sine has a maximum or minimum value.}\\[2mm]
Having put $\sin x=y$, it will be $\frac{dy}{dx}=\cos x$ and hence $\cos x=0$, whence the following values for $x$ result: $\pm \frac{\pi}{2}$, $\pm \frac{3 \pi}{2}$, $\pm \frac{5 \pi}{2}$ etc. But $\frac{ddy}{dx^2}=-\sin x$. Therefore, because these values, having  substituted them for $x$, give either $+1$ or $-1$ for $\sin x$, the latter will be maximum, the first a minimum, as it is known.

\subsection*{Example 4}

\textit{To find the arc $x$ such that the rectangle $x \sin x$ has the maximum value.}\\[2mm]
That a maximum value exists is obvious, because having put either $x=0^{\circ}$ or $180^{\circ}$ the propounded rectangle vanishes. Therefore, let $y= x\sin x$; it will be

\begin{equation*}
\frac{dy}{dx}=\sin x +x \cos x
\end{equation*}
and hence

\begin{equation*}
\tan x=-x.
\end{equation*}
Let $x=90^{\circ}+u$; it will be $\tan x= -\cot u$, therefore $\cot u= 90^{\circ}+u$. To resolve the equation in the way explained above [§ 234], put $z= 90^{\circ}+u- \cot u$ and let $f$ be the value  of the arc $u$ in question. Because $dz=du+\frac{du}{\sin^2 u}$, it will be

\begin{equation*}
p= \frac{du}{dz}=\frac{\sin^2 u}{1+\sin^2 u} \quad dp= \frac{2du \sin u \cos u}{(1+ \sin u)^2}
\end{equation*}
and hence

\begin{equation*}
\frac{dp}{dz}=q= \frac{2 \sin^3 u \cos u}{(1+\sin^2 u)^3}, \quad dq=\frac{6du \sin^2 u \cos^2 u- 2 du \sin^4 u}{(1+\sin^2 u)^3}-\frac{12 du \sin^4 u \cos^2 u}{(1+\sin^2 u)^4}.
\end{equation*}
Therefore,

\begin{equation*}
\frac{dq}{dz}=r= \frac{6 \sin^4 u \cos^2 u - 2 \sin^6 u}{(1+\sin^2 u)^4}-\frac{12 \sin^2 u^6 \cos^2 u}{(1+\sin^2 u)^5}=\frac{6 \sin^4 -14 \sin u^6 +4 \sin^8 u}{(1+\sin^2 u)^5}.
\end{equation*}
From these it will be

\begin{equation*}
f=u-pz+\frac{1}{2}qzz-\frac{1}{6}rz^3+\text{etc.}
\end{equation*}
After by trying  several values  an approximate value of $f$ was detected, put $u=26^{\circ}15^{\prime}$; it will be $90^{\circ}+u=116^{\circ}15^{\prime}$ and the arc equal to the cotangent $u$ is defined this way. From

\begin{center}
\begin{tabular}{lrllllll}
& $\log \cot u$ & $=$ & $10.3070250$ \\
subtract &  & & $4.6855749$ \\ \cline{4-4}
  &  & & $5.6214501$ \\ \cline{4-4} \\
 Therefore &  $\cot u$ & $=$ & $418263.7^{\prime \prime}$ \\ 
 or &  $\cot u$ & $=$ & $116^{\circ}11^{\prime}3\frac{7}{10}^{\prime \prime}$ \\ 
 hence &  $z$ & $=$ & $3^{\prime}56\frac{3}{10}^{\prime \prime}=256.3^{\prime \prime}$. \\ 
\end{tabular}
\end{center}
To find the value of the term $pz$ perform this calculation:

\begin{center}
\begin{tabular}{lrllllll}
& $\log \sin u$ & $=$ & $\underline{9.6457058}$ \\
& $\log \sin^2 u$ & $=$ & $\underline{9.2914116}$ \\
& $1+ \sin^2 u$ & $=$ & ${1.19561}$ \\ \cline{4-4}
& $\log (1+\sin^2 u)$ & $=$ & ${0.0775895}$ \\  \cline{4-4}
& $\log p$ & $=$ & $9.2138221$ \\
& $\log z$ & $=$ & ${2.2724637}$ \\  \cline{4-4}
 & $\log pz$ & $=$ & $1.5872858$ \\
Therefore & $pz$ & $=$ & $38.6621$ seconds\\
or & $pz$ & $=$ & $38^{\prime \prime}39^{\prime \prime \prime}43^{\prime \prime \prime \prime}$ \\
from & $u$ & $=$ & ${26^{\circ}15^{\prime}}$ \\  \cline{4-4}
it will be & $f$ & $=$ & $26^{\circ}14^{\prime}21^{\prime \prime}20^{\prime \prime \prime}17^{\prime \prime \prime \prime}$ \\
and the arc in question & $x$ & $=$ & $116^{\circ}14^{\prime}21^{\prime \prime}20^{\prime \prime \prime}17^{\prime \prime \prime \prime}$
\end{tabular}
\end{center}
But the third term $\frac{1}{2}qzz= \frac{\sin^3 u \cos u}{(1+\sin^2 u)^3}zz$ has to added. To find its value, $z$ must be expressed in parts of the radius this way:

\begin{center}
\begin{tabular}{lrllllll}
& $\log \sin z$ & $=$ & $2.3734637$ \\
add &   & & ${4.6855749}$ \\  \cline{4-4}
&  &  & $7.0590386$ \\
add & $\log \dfrac{\sin^2 u}{1+ \sin^2 u}z$ & $=$ & ${1.5872858}$ \\  \cline{4-4}
 &  &  & $8.6463244$ \\
 add & $\log \sin u$ & $=$ & ${9.6457058}$ \\
 & $\log \cos u$ & $=$ & ${9.9527308}$ \\  \cline{4-4}
 &   &  & $8.2447600$ \\
subtract & $\log (1+ \sin^2 u)^2$ & $=$ & ${0.1551790}$ \\  \cline{4-4}
& $\log \dfrac{1}{2}qzz$ & $=$ & $8.0895810$ \\
\end{tabular}
\end{center}
Therefore

\begin{equation*}
\frac{1}{2}qzz= 0.012291
\end{equation*}
or

\begin{equation*}
\frac{1}{2}qzz=44^{\prime \prime \prime \prime}15^{\prime \prime \prime \prime \prime}.
\end{equation*}
Hence, having also used this term, the arc in question will be

\begin{equation*}
x=116^{\circ}14^{\prime}21^{\prime \prime}21^{\prime \prime \prime}0^{\prime \prime \prime \prime};
\end{equation*}
but taking into account even more terms, one finds

\begin{equation*}
x=116^{\circ}14^{\prime}21^{\prime \prime}20^{\prime \prime \prime}35^{\prime \prime \prime \prime}47^{\prime \prime \prime \prime \prime}.
\end{equation*}


\title{On Maxima and Minima of multiform functions and functions of several variables
\footnote{
Original title: "'De Maximis and Minimis Functionum multiformium pluresque Variables complectentium"', first published as part of the book \textit{Institutiones calculi differentialis cum eius usu in analysi finitorum ac doctrina serierum}, 1755"', reprinted in \textit{Opera Omnia}: Series 1, Volume 10, pp. 475 - 500, Eneström-Number E212, translated by: Alexander Aycock for the Euler-Kreis Mainz}}
\author{Leonhard Euler}
\date{ }
\maketitle

\paragraph*{§273}
If $y$ was a multiform function of $x$ such that for any value of $x$ it has several real values, then, having varied $x$, those several values of $y$ are connected to each other in such a way that they represent several series of successive values. For, if we consider $y$ as ordinate of a curve, while $x$ is the abscissa,  so many  different branches of the same curved line will correspond to the same abscissa $x$ as many real different values $y$ had; and hence those successive values of $y$, which constitute the same branch, are to be considered to be connected; but the values related to the different branches will be different from each other. Therefore, we will have as many series of connected values of $y$ as  many different real values  it will have for each value of $x$; and in each arbitrary series the values of $y$, while $x$ is assumed to increase, either grow or decrease or, after they had increased, decrease  again or vice versa. From this it is perspicuous that in each series of connected values maxima and minima exist as in the case of uniform functions.

\paragraph*{§274}

To determine these maxima and minima also the same method we introduced in the preceding chapter for uniform functions can be applied. For, because, if the variable $x$ is augmented by the increment $\omega$, the function $y$ will always have this form

\begin{equation*}
y+ \frac{\omega dy}{dx}+\frac{\omega^2 ddy}{2dx^2}+\frac{\omega^3 d^3y}{6 dx^3}+\text{etc.},
\end{equation*}
it is necessary, that in the case of a maximum or minimum the term $\frac{\omega dy}{dx}$ vanishes and $\frac{dy}{dx}=0$. Therefore, the roots of this equation $\frac{dy}{dx}=0$ will indicate those values of $x$, the maxima or minima in each series of the connected values of $y$ correspond to. And it will indeed not be ambiguous, in which series of connected values a maximum or minimum is given. For, because in the equation $\frac{dy}{dx}=0$ both variables $x$ and $y$ are contained, the values of $x$ can only be defined, if by means of the equation containing the relation of the function $y$ to $x$  the variable $y$ is eliminated; but before this happens, one gets to an equation expressing the value of $y$ by means of a rational or uniform function of $x$. Hence, having found the values of $x$, the corresponding value of $y$ will be found, which will be the maximum or minimum values in the series of connected successive values one gets to.

\paragraph*{§275}

But the decision, whether these values of $y$ are maxima or minima, will be made the same way as we explained it before: Find the value of $\frac{ddy}{dx^2}$ expressed in finite terms and in it  successively substitute a value of $x$ found before for $x$; but at the same time  substitute the value for $y$  corresponding to it for the respective value of $x$; having done this, see, whether the expression $\frac{ddy}{dx^2}$ will have a positive or a negative value, and in the first case a minimum will be indicated, in the second a maximum. But if also $\frac{ddy}{dx^2}$ vanishes, then one will have to proceed to the formula $\frac{d^3y}{dx^3}$; if it does not vanish in the same case, one will have neither a maximum nor a minimum; but if also $\frac{d^3y}{dx^3}$ vanishes, the decision is to be made from the formula $\frac{d^4y}{dx^4}$ the same way as we described it for the formula $\frac{ddy}{dx^2}$. And if also $\frac{d^4y}{dx^4}$ vanishes in a certain case, one will have to proceed to the fifth differential of $y$; but,  no matter how far one had to proceed, the decisions to be made from the differentials of odd order are always  those which we explained are to be made from the formula $\frac{d^3y}{dx^3}$. In these cases in the formulas $\frac{ddy}{dx^2}$, $\frac{d^3y}{dx^3}$, $\frac{d^4y}{dx^4}$ etc. one will have to proceed so far until one gets one, which does not vanish; if this was a differential of odd order, neither a maximum nor a minimum will be indicated; but if it was of even order, its positive value will imply a minimum, its negative value a maximum.

\paragraph*{§276}

Let us put that the function $y$ is determined by $x$ through an arbitrary equation; if this equation is differentiated, it will have a form of this kind $Pdx+Qdy=0$. Therefore, having put $\frac{dy}{dx}=0$, it will be $\frac{P}{Q}=0$ and hence either $P=0$ or $Q=\infty$. The second equation, if the relation among $x$ and $y$ is expressed by means of a polynomial equation, cannot hold, because either $x$ or $y$ or both would have to become infinite. Therefore, the decision is to be made from the equation $P=0$, whose roots or values of $x$ it obtains, after by means of the propounded equation the variable $y$ was completely eliminated, will indicate the cases, in which the values of $y$ are maximum or minimum values. But to make the decision, whether a maximum or a minimum value results, examine the formula $\frac{ddy}{dx^2}$. But another differentiation of the differential equation $Pdx+Qdy=0$, if we put

\begin{equation*}
dP=Rdx+Sdy \quad \text{and} \quad dQ=Tdx+Vdy
\end{equation*}
and assume $dx$ to be constant, will give

\begin{equation*}
Rdx^2+Sdxdy+Tdxdy+Vdy^2+Qddy=0.
\end{equation*}
But since $\frac{dy}{dx}=0$, having divided by $dx^2$, it will be

\begin{equation*}
R+\frac{Qddy}{dx^2}=0 \quad \text{and hence} \quad \frac{ddy}{dx^2}=-\frac{R}{Q}.
\end{equation*}
Therefore, in the differential equation $Pdx+Qdy=0$ differentiate only the quantity $P$, while assuming $y$ to be constant, and $Rdx$ will result; then investigate the value of the fraction $\frac{R}{Q}$, which, if it was positive, will indicate a maximum, if negative, a minimum value.

\paragraph*{§277}

Let $y$ be a biform function of $x$, which is determined by this equation $yy+py+q=0$ while $p$ and $q$ denote any uniform functions of $x$. Therefore, by differentiating, it will be $2ydy+pdy+ydp+sq=0$ and hence $Pdx=ydp+dq$. Therefore, having put $P=0$, it will be $ydp+dq=0$ and $y=-\frac{dq}{dp}$ will result and so $y$ is expressed by means of an uniform function of $x$ such that, whatever value  was found for $x$, from it and $y$ it acquires one single determined value. But now the elimination of $y$ will be easy; for, if in the propounded equation $yy+py+q=0$  the value $-\frac{dq}{dp}$ is substituted for $y$, one will have $dq^2-pdpdq+qdp^2=0$, which equation, having divided it by $dx^2$ and resolved it, will yield all values of $x$ corresponding to maxima or minima; this will become more clear in the following examples.

\subsection*{Example 1}

\textit{Having propounded the equation $yy+mxy+aa+bx+nxx=0$, to define the maxima or minima of the function $y$.}\\[2mm]
Having differentiated the equation, we will have

\begin{equation*}
2ydy+mxdy+mydx+bdx+2nxdx=0,
\end{equation*}
whence 

\begin{equation*}
P= my+b+2nx \quad \text{and} \quad Q=2y+mx.
\end{equation*}
Therefore, having put $P=0$, it will be $y=-\frac{b+2nx}{m}$; this value, having  substituted it in the equation, gives

\begin{equation*}
\frac{4nn}{mm}xx+\frac{4nb}{mm}x+\frac{bb}{mm}-2nxx- bx+aa+nxx+bx=0
\end{equation*}
or

\begin{equation*}
xx=\frac{4nbx+bb+mmaa}{mmn-4nn},
\end{equation*}
whence

\begin{equation*}
x=\frac{2nb \pm \sqrt{mmnbb+mmn(mm-4n)aa}}{mmn-4nn}
\end{equation*}
or

\begin{equation*}
x=\frac{2nb \pm m \sqrt{nbb+n(mm-4n)aa}}{n(mm-4n)} \quad \text{and} \quad y= \frac{-mb \mp \sqrt{nbb+n(mm-4n)aa}}{mm-4n}.
\end{equation*}
Then having put only $x$ to be variable, $dP=2ndx$ and hence $R=2n$. But 

\begin{equation*}
Q=2y+mx= \pm \frac{\sqrt{nbb+n(mm-4n)aa}}{n},
\end{equation*}
whence

\begin{equation*}
\frac{R}{Q}=\frac{\pm 2nn}{\sqrt{nbb+n(mm-4n)aa}};
\end{equation*}
because its numerator $2nn$ is always positive, if the upper sign holds, a maximum value of $y$ will result, if the lower holds, a minimal value will result. Here, the following things have to be mentioned.\\[2mm]
I. If $m=0$, from the equation $P=0$ it immediately follows $x=-\frac{b}{2n}$ so that no elimination is necessary. And  a double value of $y$ corresponds to this value, because $y= \pm \frac{1}{2n}\sqrt{nbb-4nnaa}$, of which the positive one is a maximum value, the other negative one a minimum value.\\[2mm]
II. If  $n=0$,  $y=-\frac{b}{m}$ and $x$ grows to infinity and $y$  always has the same value such that it is neither a maximum nor a minimum value.\\[2mm]
III. If $mm=4n$, it will be $4nbx+bb+mmaa=0$ or $x= \frac{bb+mmaa}{-mmb}$ and it will be

\begin{equation*}
y=- \frac{b+2nx}{m}=- \frac{2b+mmx}{2m}=-\frac{2b}{2m}+\frac{bb+mmaa}{2mb}=\frac{mmaa-bb}{2mb}.
\end{equation*} 
Therefore,  the other value of $y$, $\frac{mmaa-bb}{2mb}$,  which will be a maximum or a minimum value, corresponds to this value of $x=\frac{-mmaa+bb}{mmb}$. But because for this value of $y$ to result in the expression

\begin{equation*}
y= \frac{-mb \mp 2\sqrt{nbb+n(mm-4n)aa}}{mm-4n}
\end{equation*}
the lower sign has to hold, the value of $y$ will be a minimum value.

\subsection*{Example 2}

\textit{Having propounded the equation $yy-xxy+x-x^3=0$, to define the maximum or minimum values of $y$.}\\[2mm]
Having differentiated the propounded equation, this equation results
\begin{equation*}
2ydy-xxdy-2xydx+dx-3xxdx=0.
\end{equation*}
And

\begin{equation*}
P=1-3xx-2xy \quad \text{and} \quad Q=2y-xx.
\end{equation*}
Therefore, having put $P=0$, it will be $y=\frac{1-3xx}{2x}$ and hence, having substituted this value,  it will be

\begin{equation*}
\frac{1}{4xx}-\frac{3}{2}+\frac{9xx}{4}-\frac{x}{2}+\frac{3}{2}x^3+x-x^3=0
\end{equation*}
or

\begin{equation*}
1-6xx+2x^3+9x^4+2x^5=0.
\end{equation*}
One of its roots is $x=-1$ and $y=1$ corresponds to it. But having assumed $y$ to be constant,  $R=-6x-2y$, therefore,

\begin{equation*}
\frac{ddy}{dx^2}=\frac{2y+6x}{2y-xx};
\end{equation*} 
in the case $x=-1$  and $y=1$ it goes over into $-4$ so that the value of $y=1$ is a maximum value. But two values of $y$ from the equation $yy-y=0$ correspond to $x=-1$; therefore, the other is $y=0$, which is neither a maximum nor a minimum value. If this equation of degree five is divided by $x+1$, an equation results whose roots cannot be exhibited in a simple manner.

\subsection*{Example 3}

\textit{Let this equation be propounded $yy+2xxy+4x-3=0$; the maximum and minimum values of $y$ are in question.} \\[2mm]
By means of differentiation this equation will result

\begin{equation*}
2ydy+2xxdy+4xydx+4dx=0.
\end{equation*}
Having put $\frac{dy}{dx}=0$, it will be $xy+1=0$ and hence $y=-\frac{1}{x}$, which value substituted in the the propounded equation gives

\begin{equation*}
\frac{1}{xx}-2x+4x-3=0= 2x^3-3xx+1,
\end{equation*}
whose roots are $x=1$, $x=1$ and $x=-\frac{1}{2}$. Since

\begin{equation*}
\frac{dy}{dx}=-\frac{4xy+4}{2y+2xx}=-\frac{2xy+2}{y+xx},
\end{equation*}
by differentiating this equation, it will be $\frac{ddy}{dx^2}=-\frac{2y}{y+xx}$, having put $y$  to be constant and so $dy=0$ and having put $xy+1=0$. Therefore, the values will be as follows:

\begin{center}
\begin{tabular}{r|r|cr}
$x$ & $y$ & $\dfrac{ddy}{dx^2}$ \\ 
$1$ & $-1$ & $\infty$ \\
$1$ & $-1$ & $\infty$ \\
$-\dfrac{1}{2}$ & $2$ & $\dfrac{-16}{9}$ & for a maximum. \\   
\end{tabular}
\end{center}
Since for the equal roots $\frac{ddy}{dx^2}=\infty$, it is not determined, whether in this case a maximum or a minimum value results. But because  at the same time $y+xx=0$, it will not even  be $\frac{dy}{dx}=0$ in this case, since $P=0$ and $Q=0$ in the fraction $\frac{dy}{dx}=-\frac{P}{Q}$; therefore, because the primary property is missing, neither a maximum nor a minimum can exist in this case. But it is indicated that in this case $x=1$ both values of $y$ become equal to each other. We will explain this nature in more detail in the third book, when we will get to the application of the differential calculus in the doctrine of curved lines. For, even though this subject also extends to this, we will nevertheless discuss it completely just in the following treatise, that it is not necessary to do it twice.

\paragraph*{§278}

There is  another kind of maxima and minima  in multiform functions; those maxima and minima,  whose nature can be explained most easily considering biform functions, are not found by the method given up to now. To see this, let $y$ be any biform function of $x$, such that, whatever value is attributed to $x$, for $y$ two values result, either both real or both imaginary. Let us put that these values of $y$ become imaginary, if one puts $x>f$, but are real, if one sets $x<f$; and having put $x=f$, both values will coalesce into one, which value we want to  be $y=g$. Therefore, because, if one takes $x>f$, the function $y$ has no real value, if it happens that for $x<f$ both values of $y$ become either greater than $g$ or smaller than $g$, in the first case the value $y=g$ will be a minimum value, in the second a maximum value, since in that case it is smaller than both preceding ones, but larger in the other case. And this maximum or minimum value cannot be found by means of the method treated up to now, because here it is not $\frac{dy}{dx}=0$. But these maxima or minima are also of a different kind, since such are not maxima or minima with respect to the preceding or following  values connected by a series, but with respect to only two distinct either preceding or following values.

\paragraph*{§279}

This happens, if the propounded equation was of this kind

\begin{equation*}
y=p \pm (f-x)\sqrt{f-x}q
\end{equation*}
 while $p$ and $q$ are functions of $x$ not divisible by $f-x$; let $q$  obtain a positive value, if one puts either $x=f$ or assumes $x$ to be little larger or a little smaller. Let $p=g$ for $x=f$ and it is obvious that in the case $x=f$ the two values of $y$ coalesce into the single one, i.e. $y=g$; but for $x>f$ both values of $y$ will become imaginary. Therefore, if we put $x$ to be a little bit smaller than $f$, say $x=f- \omega$, the function $p$ will go over into
 
 \begin{equation*}
g-\frac{\omega dp}{dx}+\frac{\omega^2 ddp}{2dx^2}-\text{etc.} 
 \end{equation*}
and $q$ into

\begin{equation*}
q- \frac{\omega dq}{dx}+\frac{\omega^2 ddq}{2dx^2}-\text{etc.},
\end{equation*} 
whence   in this case it will be

\begin{equation*}
y= g- \frac{\omega dp}{dx}+\frac{\omega^2 ddp}{2dx^2}-\text{etc.}\pm \omega \sqrt{\omega}\left(q- \frac{\omega dq}{dx}+\frac{\omega^2 ddq}{2dx^2}-\text{etc.}\right).
\end{equation*}
Let us put $\omega$ to be very small that with respect to $\omega$ its higher powers vanish, and it will be $y=g-\frac{\omega dp}{dx}\pm \omega \sqrt{\omega}q$; these two values of $y$ will both be smaller than $g$, if $\frac{dp}{dx}$ was positive, but larger, if negative. Therefore, the double value of $y=g$ in that case will be a maximum value, in the other a minimum value.

\paragraph*{§280}

Therefore, these maxima and minima result, first since  both values of $y$ become equal for $x=f$, but imaginary for $x>f$ and finally real for  $x<f$; further, since, having put $x=f- \omega$, the other irrational term yields higher powers of $\omega$ than the rational term. Therefore,  this also happens, if it was $y=p \pm (f-x)^n \sqrt{f-x}q$, as long as $n$ is an integer $>0$. But because not only the square root, but also any other root of even power introduces the same ambiguity of the sign, the same will happen, if it was $y=p \pm (f-x)^{\frac{2n+1}{2m}}q$, as long as it is $2n+1 >2m$; therefore, it will be $(y-p)^{2m}=(f-x)^{2n+1}q^{2m}$ or $(y-p)^{2m}=(f-x)^{2n+1}Q$. Therefore, if the function $y$ is expressed by means of such an equation such that  $2n+1 >2m$, having put $x=f$, the value of $y$ will be a maximum or minimum; the first, if $\frac{dp}{dx}$ was a positive quantity, the latter, if  $\frac{dp}{dx}$ is negative for $x=f$. But if in this case $\frac{dp}{dx}=0$, then it will be

\begin{equation*}
y=g+\frac{\omega^2ddp}{2dx^2}\pm \omega^{\frac{2n+1}{2m}}q.
\end{equation*}
Therefore, only if $\frac{2n+1}{2m}>2$, a maximum or a minimum can exist; but if $\frac{2n+1}{2m}>2$, then $y=g$ will have a maximum value, if $\frac{ddp}{dx^2}$ had a negative value, a minimum value, if it was positive; and like this the decision will have to be made, if also $\frac{ddp}{dx^2}$ vanishes.
 
\paragraph*{§281}

Therefore, if $y$ was a function of $x$ of such a kind, it can happen, that except for the maxima and minima, which the first method exhibits, also maxima and minima of the other species are given, which can be explored by the method explained here. We will show this in the following examples.

\subsection*{Example 1}

\textit{To determine the maxima and minima of the function $y$ defined by this equation}

\begin{equation*}
yy-2xy-2xx-1+3x+x^3=0.
\end{equation*}
To investigate the maxima and minima of the first kind, differentiate the equation and it will be

\begin{equation*}
2ydy-2xdy-2ydx-4x+3dx+3xxdx=0
\end{equation*} 
and, having put $\frac{dy}{dx}=0$, it will be

\begin{equation*}
y=\frac{3}{2}-2x+\frac{3}{2}xx,
\end{equation*}
which value substituted in the first equation gives

\begin{equation*}
9x^4-32x^3+42xx-24x+5=0,
\end{equation*}
which is resolved into

\begin{equation*}
9xx-14x+5=0 \quad \text{and} \quad xx-2x+1=0.
\end{equation*}
The first gives $x=1$ twice and the corresponding value is $y=1$, whence in this case in the fraction

\begin{equation*}
\frac{dy}{dx}= \frac{2y-3+4x-3xx}{2y-2x}
\end{equation*}
the denominator also vanishes and so a maximum or minimum of the first kind can not exist in this case; the first equation $9xx-14x+5=0$ on the other hand will give $x=1$ and $x=\frac{5}{9}$, the first of which values leads to the same inconvenience as the first. But having put $x=\frac{5}{9}$, $y=\frac{3}{2}-\frac{10}{9}+\frac{25}{54}=\frac{23}{27}$. And because $\frac{dy}{dx}=\frac{2y-3+4x-3xx}{2y-2x}$, it will be

\begin{equation*}
\frac{ddy}{dx^2}=\frac{4-6x}{2y-2x}=\frac{-3x+2}{y-x},
\end{equation*}
because of $dy=0$ and the numerator $=0$. Therefore, it will be $\frac{ddy}{dx^2}=\frac{9}{8}$, whence this value $x=\frac{5}{9}$ is a minimum  of the first kind. Further, since $(y-x)^2=(1-x)^3$, it will be

\begin{equation*}
y=x \pm (1-x)\sqrt{1-x}
\end{equation*}
and hence, having put $x=1$, a maximum value of the second kind results; for, having put $x=1- \omega$ it will be $y=1- \omega \pm \omega \sqrt{\omega}$, both of which values are smaller than $1$, if a very small $\omega$ is taken.

\subsection*{Example 2}

\textit{To find the maxima and minima of the function $y=2x-xx\pm (1-x)^2\sqrt{1-x}$.} \\[2mm]
To find the maxima and minima of the first kind, differentiate the equation and it will be

\begin{equation*}
\frac{dy}{dx}=2-2x \mp \frac{5}{2}(1-x)\sqrt{1-x},
\end{equation*}
which fraction put $=0$ at first yields $x=1$, and because

\begin{equation*}
\frac{ddy}{dx^2}=-2 \pm \frac{15}{4}\sqrt{1-x},
\end{equation*}
$y$ in this case will be a maximum value of the first kind and  $y=1$. Having divided the equation $\frac{dy}{dx}=0$ by $1-x$, it will be

\begin{equation*}
4 \mp 5\sqrt{1-x}=0 \quad \text{or} \quad 16=25-25x,
\end{equation*}
whence  $x=\frac{9}{25}$ and $\frac{ddy}{dx^2}=-2 \pm 3$. Hence, if the upper sign holds,  $y=\frac{2869}{3125}$ will be a minimum value; but if the lower sign holds, it will be $y=\frac{821}{3125}$, which might seem to be a maximum value;  but indeed only the upper sign can hold, since $4 \mp 5 \sqrt{1-x}$ can only be $=0$, if $\sqrt{1-x}=+\frac{4}{5}$. Therefore, we find a maximum value of the first kind in the case $x=1$ and $y=1$ and a minimum in the case $x=\frac{9}{25}$ and $y=\frac{2689}{3125}$.  A maximum of the other kind also results, if $x=1$, in which case $y=1$. For, having put $x=1- \omega$, it will be $y=1- \omega \omega \pm \omega^2 \sqrt{\omega}$, this is $<1$ in both cases . Therefore, here, if $x=1$, a maximum value of the first and the second kind coalesce and constitute a mixed maximum value.

\paragraph*{§282}

From these examples not only the nature of this other kind of maxima and minima is understood, but also one can also form functions of this kind arbitrarily, which admit  maxima or minima of the second kind. But how, if any function was propounded, it can be decided, whether it has a maximum or minimum value of such a kind or not,  will be shown in the following book, since the nature of curved lines is  illustrated in the best way by this investigation. Furthermore, it is easily understood, if $y$ was a function of $x$ of such a kind which has a maximum or minimum value of the second kind, that then also vice versa $x$ will be a function of such a kind of $y$. For, since from this equation $(y-x)^2=(1-x)^3$ for $x=1$ the function $y$ has a maximum value of the second kind, if the variables $x$ and $y$ are interchanged, this equation $(x-y)^2=(1-y)^3$ for $y$ will exhibit a function of such a kind of $x$, which has a maximum value of the second kind. For, having put $x=1+\omega$, it will be $(1+\omega -y)^2=(1-y)^3$; hence, if we set $y=1+\varphi$, it will be $(\omega -\varphi)^2=(-\varphi)^3=-\varphi^3$ and hence $\varphi$ must be negative. Therefore, let $y=1- \varphi$; it will be $(\omega +\varphi)^2=\varphi^3$, because having taken a very small $\varphi$,  $\varphi^3$ vanishes with respect to $\varphi^2$, $\omega$ will have to be negative; therefore, no real values of $y$ correspond to the value $x=1+\omega$ . But having put $x=1- \omega$ and $y=1- \varphi$, because of $(\varphi - \omega)^2= \varphi^3$, it will be $\varphi = \omega \pm \omega \sqrt{\omega}$ and hence $y=1- \omega \mp \omega \sqrt{\omega}$, whence both values of $y$ corresponding to $x=1- \omega$ are smaller than the value $y=1$ which corresponds to the value $x=1$; and as a logical consequence this value of $y$ will be a maximum value of the second kind.

\paragraph*{§283}

Up to now we only considered biform functions, whose maxima or minima, since  both values can easily be expressed by means of the resolution of  the quadratic equation, can be checked for their correctness. But if the function $y$ is expressed by means of an equation of higher order, the method explained before we used to investigate maxima and minima of the first kind can be applied with the same success. But let us reserve the discussion of maxima and minima of the second kind for the following treatise. Therefore, let us show how to treat triform and multiform functions in several examples.

\subsection*{Example 1}

\textit{Define the function $y$, whose maxima and minima are in question, through this equation}

\begin{equation*}
y^3+x^3=3axy.
\end{equation*}    
Having differentiated this equation, 

\begin{equation*}
3y^2dy+3xxdy=3axdy+3aydx
\end{equation*}
and hence

\begin{equation*}
\frac{dy}{dx}=\frac{ay-xx}{yy-ax}.
\end{equation*}
Therefore, a maximum or minimum will exist, if it was $ay=xx$ or $y=\frac{xx}{a}$, which value, having substituted it in the propounded equation gives

\begin{equation*}
\frac{x^6}{a^3}+x^3=3x^3 \quad \text{and} \quad x^6=2a^3x^3.
\end{equation*}
Therefore, it will be $x=0$ trice, in which case because of $y=\frac{xx}{a}=0$ also the denominator $yy-ax=0$. Whether in this case a maximum or minimum results,  will be seen, if we attribute a value  differing hardly from $0$ to $x$. Therefore, let $x=\omega$ and $y=\varphi$; because of $\varphi^3+\omega^3=3a \omega \varphi$, it will be either $\varphi = \alpha \sqrt{\omega}$ or $\varphi = \beta \omega^2$. In the first case it will be $\alpha^3 \omega \sqrt{\omega}=3 \alpha a \omega \sqrt{\omega}$ and hence $\alpha =\sqrt{3}a$. Therefore, having put $x=\omega$, it will be $y= \pm \sqrt{3}\alpha \omega$. Hence, even though $\omega$ cannot be taken negatively, nevertheless one of the two values of $y$ will be greater than $0$, the other will be smaller and hence $y=0$ will be neither a maximum nor a minimum value. But if one sets $\varphi = \beta \omega^2$, it will be $\omega^3=3a a \beta \omega^3$ and hence $\beta =\frac{1}{3a}$ and $\varphi= \frac{\omega^2}{3a}$. Therefore, in this case, no matter whether $x$ is taken $=+\omega$ or $=- \omega$, the value of $y=\varphi$ will be greater than zero and hence in this case $y=0$ will be a minimum value. Therefore, finally the third case to be examined follows from the equation $x^3=2a^3$, which gives $x=a \sqrt[3]{2}$ and $y=a \sqrt[3]{4}$. Whether this is a maximum or a minimum, is to be investigated using the second differential of the equation $\frac{dy}{dx}=\frac{ay-xx}{yy-ax}$, which differential, because of $dy=0$ and $ay-xx=0$,  will be $\frac{ddy}{dx^2}=\frac{-2x}{yy-ax}$, whose value in the present case is $-\frac{2a \sqrt[3]{2}}{2a^2 \sqrt[3]{2}-aa\sqrt[3]{2}}=-\frac{2}{a}$, which indicates that the value of $y$ is a maximum value.

\subsection*{Example 2}

\textit{If the function $y$ is defined by means of this equation $y^4+x^4+ay^3+ax^3=b^3x+b^3y$, to find its maximum and minimum values.}\\[2mm]
Because by  differentiation this equation results

\begin{equation*}
4y^3dy+3ayydy-b^3dy=b^3-3axxdx-4x^3dx,
\end{equation*}
it will be
\begin{equation*}
\frac{dy}{dx}=\frac{b^3-3axx-4x^3}{4y^3+3ayy-b^3}
\end{equation*}
and one has to put $b^3=3axx+4x^3$. Therefore, the question is reduced to this, that  the maxima and minima of the uniform function $b^3-ax^3-x^4$ are investigated, which at the same time will be the maxima and the minima of the function $y$. Let $a=2$ and $b=3$ or let this equation be propounded $y^4+x^4+2y^3+2x^3=27x+27y$; it will be $\frac{dy}{dx}=\frac{27-6xx-4x^3}{4y^3+6yy-27}$ and $4x^3+6xx-27=0$, which divided by $2x-3=0$ gives $2xx+6x+9=0$; because the roots of the second equation are imaginary, it will be $x=\frac{3}{2}$ and $y^4+2y^3-27y=\frac{459}{16}$, each root of which will be either a maximum or minimum. But since $\frac{dy}{dx}=\frac{27-6xx-4x^3}{4y^3+6yy-27}$, it will be $\frac{ddy}{dx^2}=\frac{-12x-12xx}{4y^3+6yy-27}$, which, if positive for $x=\frac{3}{2}$, will indicate a minimum, otherwise a maximum.

\subsection*{Example 3}

\textit{If it was $y^m+ax^n=by^px^q$, to define the maximum and minimum values of $y$.}\\[2mm]
By means of differentiation 

\begin{equation*}
\frac{dy}{dx}=\frac{qby^px^{q-1}-nax^{n-1}}{my^{m-1}-pby^{p-1}x^q},
\end{equation*}
having put which fraction $=0$ it will at first be $x=0$, if $n$ and $q$ are greater than  $1$, and at the same time it was $y=0$. To decide, whether in this case a maximum or a minimum value is given, the closest values are to be investigated, since also the denominator becomes $=0$; this investigation will depend mainly on the exponents. Furthermore, the equation $\frac{dy}{dx}=0$ will  give $y^p=\frac{na}{qb}x^{n-q}$, which value, having substituted it in the propounded equation, by putting $\frac{na}{qp}=g$ will give

\begin{equation*}
g^{\frac{m}{p}}x^{\frac{mn-mq}{p}}+ax^n=\frac{na}{q}x^n \quad \text{or} \quad g^{\frac{m}{p}}x^{\frac{mn-mq-np}{p}}=\frac{(n-q)a}{q},
\end{equation*}
whence

\begin{equation*}
x= \left(\frac{(n-q)a}{q}\right)^{p:(mn-mq-np)}:g^{m:(mn-mq-np)}
\end{equation*}
and at the same time the value of $y$ becomes known. Further, it is to be considered, whether the second differential

\begin{equation*}
\frac{ddy}{dx^2}=\frac{q(q-q)by^px^{q-2}-n(n-1)ax^{n-2}}{my^{m-1}-pby^{p-1}x^q}
\end{equation*}
has a positive or a negative value, since in the first  case a minimum value will be indicated, in the second a maximum value.

\subsection*{Example 4}

\textit{If it was $y^4+x^4=4xy-2$, to assign the maxima and minima of the function $y$.}\\[2mm]
By differentiation 

\begin{equation*}
\frac{dx}{dy}=\frac{y-x^3}{y^3-x}
\end{equation*}
and hence $y=x^3$ results; therefore, it will be $x^{12}=3x^4-2$ or $x^{12}-3x^4+2=0$, which equation is resolved into  $x^4-1=0$ and $x^8+x^4-2=0$ and the latter into $x^4-1=0$ and $x^4+2=0$. Therefore, it will be either $x=\pm 1$ or $x=-1$ twice; in both cases also the denominator of the fraction $\frac{dy}{dx}$ vanishes. Therefore, to investigate, whether in these cases a maximum or a minimum exists, let us put $x=1-\omega$ and $y=1- \varphi$; it will be

\begin{equation*}
1-4\varphi +6 \varphi^2-4 \varphi^3+\varphi^4+1-4 \omega +6 \omega^2-4 \omega^3+\omega^4
\end{equation*}
\begin{equation*}
=4-4 \omega -4 \varphi +4 \omega \varphi-2
\end{equation*}
and hence

\begin{equation*}
4 \omega \varphi = 6 \varphi^2 +6 \omega^2 -4 \varphi^3-4 \omega^3+ \varphi^4+\omega^4
\end{equation*}
and, because of the very small $\omega$ and $\varphi$, $4 \omega \varphi =6 \varphi^2+ 6 \omega^2$. Therefore, the value of $\varphi$ will be imaginary, no matter whether $\omega$ is taken positively or negatively. Or if $y$ and $x$ denote the coordinates of a curve,  in the case $x=1$ and $y=1$ it will have a conjugated point. Therefore, also this value can be neither a maximum nor a minimum, since the preceding and following, to which it would have been compared, become  imaginary.

\paragraph*{§284}

If the equation expressing the relation among $x$ and $y$ was of such a nature that the function of $y$ becomes equal to a function of $x$, say $Y=X$, then, to find the maxima or minima one will have to put $dX=0$; therefore, $y$ will have a maximum or a minimum value in the same cases, in which $X$ has a maximum or a minimum value. In like manner, if $x$ is considered as a function of $y$, $x$ will have a maximum or a minimum value, if $dY=0$, this means, if $Y$ had a maximum or a minimum value. And from this it does nevertheless not follow that $y$ at the same time has a maximum or a minimum value. For, if it was $2ay-yy=2bx-xx$, $y$ will have a maximum or a minimum value, if $x=b$, and it will be $y=a \pm \sqrt{aa-bb}$. On the other hand, $x$ has a maximum or a minimum value, if $y=a$, and $x=b \pm \sqrt{bb-aa}$ and therefore $y$ is neither a maximum nor a minimum, if $x=b \pm \sqrt{bb-aa}$, in which case $x$  nevertheless has a maximum or a minimum value. Furthermore, in this case, if $y$ has maximum or minimum values, $x$ will have none at all; for, $y$ cannot have a maximum or a minimum value, if it was not $a>b$, in which case the maximum or minimum of $x$ becomes imaginary.

\paragraph*{§285}

Furthermore, it can  happen that not all roots of the equation $dX=0$ yield maximum or minimum values of $y$; for, if that equation had two equal roots, from this neither a maximum nor a minimum follows; the same happens, if an even number or roots were equal to each other. So, if the equation $b(y-a)^2=(x-b)^3+c^3$ is propounded, since having taken the differentials it is $2bdy(y-a)=3dx(x-b)^2$, the function $y$ will have neither a maximum nor a minimum value for $x=b$, since here two equal roots occur. But if $x$ is considered as a function of $y$, it will have a maximum or minimum value, if one sets $y=a$, and $x=b-c$ will be a minimum. Finally, since in equations of this kind $Y=X$ the variables $x$ and $y$ are not mixed, if  a value is attributed to $x$, which is a root of the equation $dX=0$, all values of $y$, no matter how many were real, will be maxima or minima; this does not happen, if in the equation the two variables were mixed.

\paragraph*{§286}

 We will reserve the things, which are still to be explained about the nature of maxima and minima,  for the following book, since  they can be represented and explained in a more convenient way showing figures. Therefore, let us proceed to functions which  are composited of several variables, and let us investigate the values, which have to be attributed to each variable that the function has a maximum or a minimum value. And at first it is certainly clear, if the variables were not mixed, such that a function of this kind is propounded $X+Y$, where $X$ is a function only of $x$ and $Y$ is a function only of $y$, that  the propounded function $X+Y$ will have a maximum, if $X$ and $Y$ have maximum value at the same time, and a minimum, if  $X$ and $Y$  have a minimum value at the same time. Therefore, to find the maxima and minima, find the values of $x$, for which $X$ has a maximum value, and in like manner the values of $y$, for which $Y$ has a maximum value;  analogously  for minima, of course. Therefore, one must be careful not to combine two values of $x$ and $y$ of different nature, of one which renders $X$ maximal, the other $Y$ minimal, or vice versa. For, if this would happen, the function $X+Y$ would have neither a maximum nor a minimum value in this case. But a function of this kind $X-Y$ will have maximum value, if $X$ has a maximum value and at the same time $Y$ has a minimum value; on the other hand $X-Y$ will have a minimum value, if $X$ had a minimum value and $Y$ a maximum value. But if both of the functions $X$ and $Y$  would have either a maximum or minimum value, their difference $X-Y$ would have neither a maximum nor a minimum value; all these things are clear and perspicuous from the nature of maxima and minima explained before.

\paragraph*{§287}

Therefore, if the maximum or minimum values of a function of two variables are in question, the question is a lot more intricate than the same question for only one variable. For, not only for each variable the cases, in which a maximum or a minimum is produced, are to be distinguished carefully, but also from these two of such a kind are to be combined that the propounded function has a maximum or a minimum value; this will become more clear from the consideration of some examples.

\subsection*{Example 1}

\textit{Let this function of the two variables $x$ and $y$ be propounded $y^4-8y^3+18y^2-8y+x^3-3xx-3x$ and let the values for $y$ and $x$ be in question which are to be substituted that this function has a maximum or minimum value.}\\[2mm]
Since this expression is resolved into two parts of this kind $Y+X$, of which one is a function only of $y$, the other a function only of $x$, investigate the cases, in which each of  has a maximum or minimum value. Therefore, since

\begin{equation*}
Y=y^4-8y^3+18y^2-8y,
\end{equation*}
it will be

\begin{equation*}
\frac{dY}{dy}=4y^3-24y^2+36y-8;
\end{equation*}
having put this expression equal to zero and having divided by $4$, it will be

\begin{equation*}
y^3-6y^2+9y-2=0,
\end{equation*}
whose roots are $y=2$ and $y=2 \pm \sqrt{3}$. Therefore, because  $\frac{ddY}{4dy^2}=3yy-12y+9$, in the case $y=2$ a maximum will result. For the remaining two roots $y=2 \pm \sqrt{3}$, which result from the equation $yy-4y+1=0$, it will be $\frac{ddY}{12dy^2}=yy-4y+3=2$, whence both of them give a minimum. But in these cases it will be as follows:

\begin{center}
\begin{tabular}{l|ll}
$y=2$  & $Y=+8$ & maximum \\
$y=2-\sqrt{3}$  & $Y=-1$ & minimum \\
$y=2+\sqrt{3}$  & $Y=-1$ & minimum \\
\end{tabular}
\end{center}
In like manner, since

\begin{equation*}
X=x^3-3xx-3x,
\end{equation*}
it will be

\begin{equation*}
\frac{dX}{dx}=3xx-6x-3,
\end{equation*}
whence this equation results

\begin{equation*}
xx=2x+1
\end{equation*}
and $x=1 \pm \sqrt{2}$. But  $\frac{ddX}{6dx^2}=x-1=\pm \sqrt{2}$. Therefore, the root $x=1+\sqrt{2}$ gives a minimum, namely $X=-5-4 \sqrt{2}$, and $x=1-\sqrt{2}$ gives a maximum, namely $X=-5-4\sqrt{2}$. Therefore, the propounded formula

\begin{equation*}
X+Y=y^4-8y^3+18yy-8y+x^3-3xx-3x
\end{equation*}
will have a maximum value, if one puts $y=2$ and $x=1-\sqrt{2}$, and $X+Y=3+4\sqrt{2}$ will result. But the same formula $X+Y$ will have a minimum, if one takes either $y=2-\sqrt{3}$ or $y=2+\sqrt{3}$ and $x=1+\sqrt{2}$; in both cases it will be $X+Y=-6 -4\sqrt{2}$.

\subsection*{Example 2}

\textit{If this function of two variables is propounded $y^4-8y^3+18y^2-8y-x^3+3xx+3x$, to investigate in which which cases it has either a maximum or minimum value.}\\[2mm]
Having put, as we did in the preceding example,

\begin{equation*}
Y=y^4-8y^3+18y^2-8y \quad \text{and} \quad X=x^3-3xx-3x
\end{equation*}
the propounded formula will be $Y-X$ and will hence have a maximum, if $Y$ had a maximum value and $X$ a minimum value. Therefore, because  we found these cases before already, it is plain that $Y-X$ has a maximum value, if one puts $y=2$ and $x=1+\sqrt{2}$; and it will be $Y-X=13+4\sqrt{2}$. The value of $Y-X$ will be a minimum, if $Y$ has a minimum value and $X$ a maximum value, which happens by putting $y=2 \pm \sqrt{3}$ and $x=1-\sqrt{2}$;  it will be $Y-X=4-4 \sqrt{2}$. Furthermore, in each of the two examples it is plain that these values we found are neither the largest nor the smallest of all; for, if, for the sake of an example, one would put $y=100$ and $x=0$, without any doubt a value greater than the one we found would result; and in like manner, putting $y=0$ and  either $x=-100$ or $x=+100$, a smaller value than those we found  for the case of a minimum would result. Therefore, recall the definition of maxima and minima we gave above, i.e. that we called a value a maximum, which is larger than the closest preceding and the following values, but this value is a minimum, if it was smaller than both the closest preceding and following values. So, in this example  the value of $Y-X$, which results by putting $y=2$ and $x=1+\sqrt{2}$, is larger than those which result, if one puts $y=2 \pm \omega$ and $x=1+\sqrt{2}\pm \varphi$ having taken sufficiently small quantities for $\omega$ and $\varphi$.

\paragraph*{§288}

Having treated these examples, the way to find the general solution will be easier. Let $V$ denote any function of the two variables $x$ and $y$ and  let the values to be found for $x$ and $y$ which render the function maximal or minimal. Therefore, because to achieve this   a determined value is to be attributed to both variables $x$ and $y$, let us put that the one $y$ already has the value, which is required to render the function $V$ either maximal or minimal, and having put this it will only be necessary that for the other variable $x$ also an appropriate value is found, which will happen, if the function $V$ is differentiated considering only $x$ as variable and the differential is set equal to zero. In a like manner, if we assume the variable $x$ to already have the value, which is apt to render the function $V$ either maximal or minimal, the value of $y$ will be found by differentiating with respect to $y$ only and putting this differential equal to zero. Hence, if the differential of the function $V$ was $=Pdx+Qdy$, it will be necessary that $P=0$ and $Q=0$, from which two equations the value of both variables $x$ and $y$ can be found.

\paragraph*{§289}

Since this way, without any difference, the values for $x$ and $y$ are found, by which the function $V$ is rendered either maximal or minimal, the cases, in which either a maximum or a minimum value results, are to be distinguished carefully from each other. For, that a function $V$ has a maximum value, it is necessary, that both variables conspire for this; for, if the one would exhibit a maximum, the other a minimum, the function itself would become neither a maximum nor a minimum. Therefore, having found the values of $x$ and $y$ from the equations $P$ and $Q$, it is to be investigated, whether both a the same time render the function $V$ either maximal or minimal; and just then, after it is certain that the value of both variables found from this lead to a maximum value, we will be able to affirm that the function  in this case has a maximum value. The same is to be understood for a minimum, such that the function $V$ can only have a minimum value, if at the same time both variables $x$ and $y$ produce a minimum. Therefore, all those cases are to be rejected, in which the one variable is detected to indicate a maximum, and the other a minimum. Sometimes it even happens that the values resulting from the equations $P=0$ and $Q=0$ of one of the variables or of both exhibit neither a maximum nor a minimum, which cases are to be rejected as inept in exactly the same way.

\paragraph*{§290}

But whether the values found for $x$ and $y$ lead to a maximum or a minimum, will be investigated for each of them separately in the same way as above, i.e. as if only one variable was there. To make the decision for the variable $x$ consider the other $y$ as constant, and because  $dV=Pdx$ or $\frac{dV}{dx}=P$, differentiate $P$ again having put $y$ to be constant that $\frac{ddV}{dx^2}=\frac{dP}{dx}$ results, and consider, whether the value of $\frac{dP}{dx}$, after  the values found before were substituted for $x$ and $y$, becomes positive or negative; for, in the first case a minimum will be indicated, in the second  a maximum. Because in like manner for constant $x$ $dV=Qdy$ or $\frac{dV}{dy}=Q$, differentiate $Q$ again having put only $y$ to be variable and examine the value $\frac{dQ}{dy}$ after having substituted the values for $x$ and $y$, which were found from the equations $P=0$ and $Q=0$; if  it was positive,  it will indicate a minimum, otherwise a maximum. Therefore, it is concluded, if from the values found for $x$ and $y$ the formulas $\frac{dP}{dx}$ and $\frac{dQ}{dy}$ have values affected with different signs, the one positive, the other negative, then the function $V$ will have neither a maximum nor minimum value; but if both of the formulas $\frac{dP}{dx}$ and $\frac{dQ}{dy}$ are positive, a minimum value will result, and otherwise, if both of them become negative, a maximum value. 

\paragraph*{§291}

If one of the formulas $\frac{dP}{dx}$ and $\frac{dQ}{dx}$ or even both vanish, if  the found values are substituted for $x$ and $y$, then one will have to proceed to the following differentials $\frac{ddP}{dx^2}$ and $\frac{ddQ}{dy^2}$;  if they not equally vanish, neither a maximum nor a minimum will exist; but if they vanish, then the decision is to be made from the following differentials $\frac{d^3P}{dx^3}$ and $\frac{d^3Q}{dy^3}$ in the same way it was made for the formulas $\frac{dP}{dx}$ and $\frac{dQ}{dy}$. To be explain the cases, in which this happens, more clearly,  let the value $x=\alpha$ have been resulted; if it renders the formula $\frac{dP}{dx}$ vanishing, it is necessary that $\frac{dp}{dx}$ has the factor $x-\alpha$; if this factor was the only one of this kind, neither a maximum nor a minimum will be indicated; the same happens, if $\frac{dP}{dx}$ had the factor $(x-\alpha)^3$ or $(x-\alpha)^5$ etc. But if the factor was either $(x-\alpha)^2$ or $(x- \alpha)^4$ etc., then a maximum or minimum will be indicated; but furthermore, one will have to see, whether it is in agreement with the case indicated by $y$.

\paragraph*{§292}

But the work to proceed to higher differentials in these cases can be reduced tremendously; for, if, in order to consider the subject in more generality, we assume that $\alpha x+ \beta =0$ was found and the formula $\frac{dP}{dx}$ has the factor $(\alpha x+ \beta)^2$ such that  $\frac{dP}{dx}=(\alpha x+ \beta)^2 T$, since  $\alpha x+ \beta =0$, it will be $\frac{d^3P}{dx^3}=2 \alpha^2 T$ and hence, because of the positive $2  \alpha^2 $, the decision can be made from the quantity $T$; if it has a positive value, it will indicate a minimum, otherwise a maximum. And the same  auxiliary theorem can be applied in the investigation of maxima and minima, if one single variable is contained, such that it never necessary to ascend to higher differentials. It is not even necessary to ascend to the second differential; for, if from the equation $P=0$  $\alpha x+ \beta =0$, it is necessary, that $P$ has a factor $\alpha x+ \beta$; let $P=(\alpha x+ \beta)T$, and because

\begin{equation*}
\frac{dP}{dx}=\alpha T+(\alpha x+ \beta)\frac{dT}{dx},
\end{equation*}
because of $\alpha x+ \beta =0$, it will be $\frac{dP}{dx}= \alpha T$ and hence already the factor $T$, depending on whether the value of $\alpha T$ was positive or negative, will immediately indicate a minimum or a maximum.

\paragraph*{§293}

Therefore, having given these prescriptions, it will not be difficult, if an arbitrary function involving two variables was propounded, to investigate the cases, in which this function has maximum or minimum value. If more is to be taken into account, the expansion of some examples will suggest this, which is why it will be convenient to illustrate the given rules by some examples.

\subsection*{Example 1}

\textit{Let this function of two variables  be propounded $V=xx+xy+yy-ax-by$; to investigate, in which cases it has maximum or minimum values.}\\[2mm]
Since $dV=2xdx+ydx+xdy+2ydy-adx-bdy$, if this expression is compared to the general formula $dV=Pdx+Qdy$, it will be

\begin{equation*}
P=2x+y-a \quad \text{and} \quad Q=2y+x-b,
\end{equation*}
whence these equations will be formed

\begin{equation*}
2x+y-a=0 \quad \text{and} \quad 2y+x-b=0,
\end{equation*}
having combined which and thus eliminating $y$, it will be $x-b=4x-2a$ and hence

\begin{equation*}
x=\frac{2a-b}{3} \quad \text{and} \quad y=a-2x= \frac{2b-a}{3}.
\end{equation*}
Therefore, since

\begin{equation*}
\frac{dP}{dx}=2 \quad \text{and} \quad \frac{dQ}{dy}=2,
\end{equation*}
both of them indicate a minimum; from this we conclude that the formula

\begin{equation*}
xx+xy+yy-ax-by
\end{equation*}
has a minimum value, if one puts $x=\frac{2a-b}{3}$ and $y=\frac{2b-a}{3}$, and this way it will be
 
 \begin{equation*}
 V=\frac{-3aa+3ab-3bb}{9}= \frac{-aa+ab-bb}{3};
 \end{equation*}
since this is the only one, it will be the smallest of all. Therefore, it can only  be

\begin{equation*}
xx+xy+yy-ax-by= \frac{-aa+ab-bb}{3},
\end{equation*} 
and since it cannot be smaller, this equation will be impossible

\begin{equation*}
xx+xy+yy-ax-by=\frac{-aa+ab-bb}{3}-cc.
\end{equation*}

\subsection*{Example 2}

\textit{If the formula $V=x^3+y^3-3axy$ is propounded, let the cases be  in question in which $V$ has a maximum or minimum value.}\\[2mm]
Because of $dV=3xxdx+3yydy-3aydx-3axdy$, it will be

\begin{equation*}
P=3xx-3ay \quad \text{and} \quad Q=3yy-3ax,
\end{equation*}
whence 

\begin{equation*}
ay=xx \quad \text{and} \quad ax=yy.
\end{equation*}
Therefore, because of $yy=x^4:aa=ax$, it will be $x^4-a^3x=0$ and hence either $x=0$ or $x=a$. In the first case  $y=0$, in the second $y=a$. Therefore, because

\begin{equation*}
\frac{dP}{dx}=6x, \quad \frac{ddP}{dx^2}=6 \quad \text{and} \quad \frac{dQ}{dy}=6y \quad \text{and} \quad \frac{ddQ}{dy^2}=6,
\end{equation*}
in the first case, in which $x=0$ and $y=0$, neither a maximum nor a minimum results. But in the second case, in which both $x=a$ and $y=a$, a minimum results, if $a$ was a positive quantity, and it will be $V=-a^3$, which value is only smaller than the closest preceding and following ones; for, without any doubt $V$ can have  much smaller values, if   negative values are attributed to both variables $x$ and $y$.

\subsection*{Example 3}

\textit{Let this function be propounded $V=x^3+ayy-bxy+cx$, whose maximal or minimal values are to be investigated.}\\[2mm]
Since $dV=3xxdx+2aydy-bydx-bxdy+cdx$, it will be

\begin{equation*}
P=3xx-by+c \quad \text{and} \quad Q=2ay-bx,
\end{equation*}
having put which values equal to zero it will be $y=\frac{bx}{2a}$ and hence

\begin{equation*}
3xx-\frac{bbx}{2a}+c=0 \quad \text{or} \quad xx=\frac{2bbx-4ac}{12a},
\end{equation*}
whence

\begin{equation*}
x=\frac{bb \pm \sqrt{b^4-48aac}}{12ac}.
\end{equation*}
Therefore, only if $b^4-48aac>0$, a maximum or a minimum can exist. Therefore, let us put that $b^4-48aac=bbff$ that $c=\frac{bb(bb-ff)}{48aa}$; it will be

\begin{equation*}
x= \frac{bb \pm bf}{12a} \quad \text{and} \quad y=\frac{bb(b \pm f)}{24aa}.
\end{equation*}
Further, since

\begin{equation*}
\frac{dP}{dx}=6x \quad \text{and} \quad \frac{dQ}{dy}=2a,
\end{equation*}
it will be

\begin{equation*}
\frac{dP}{dx}=\frac{b(b \pm f)}{2a}.
\end{equation*}
Therefore, only if $2a$ and $\frac{b(b \pm f)}{2a}$ are quantities of the same sign, either a maximum or a minimum can exist. And if they are indeed both either positive or negative, what happens, if their product $b(b \pm f)$ was positive, then the function $V$ has a minimum value, if $a$ was a positive quantity, otherwise a maximum value, if $a$ is a negative quantity.  Hence, if it was $f=0$ or $c=\frac{b^4}{48aa}$, because of the positive quantity $bb$, the function $V$ will have a minimum value, if $a$ was a positive quantity and one sets $x=\frac{bb}{12a}$ and $y= \frac{b^3}{24aa}$; otherwise, if $a$ is negative, these substitutions produce a maximum. If $f<b$, in the two cases either a maximum or a minimum results; but if $f>b$, then only the case $x=\frac{b(b+f)}{12a}$ and $y=\frac{bb(b+f)}{24aa}$ will yield a maximum or a minimum, depending on whether $a$ was negative or positive. Let $a=1$, $b=3$ and $f=1$, that one has this formula $V=x^3+yy-3xy+\frac{3}{2}x$; this will have a minimum value, because of the positive $a$, if one puts either $x=1$ and $y=\frac{3}{2}$ or $x=\frac{1}{2}$ and $y=\frac{3}{4}$. In the first case   $V=\frac{1}{4}$, in the second $V=\frac{5}{16}$. Nevertheless, it is plain that, writing negative numbers instead of $x$, a lot smaller values for $V$ can result. Therefore, the value of $V=\frac{1}{4}$ must be understood to be smaller than if one puts $x=1+\omega$ and $y=\frac{3}{2}+\varphi$, if $\omega$ and $\varphi$ are small numbers, either positive or negative; but  $\omega$ must not be larger than $-\frac{15}{4}$; for, if $\omega < - \frac{15}{4}$, it can happen that $V$ becomes smaller than $\frac{1}{4}$.

\subsection*{Example 4}

\textit{To find the maxima or minima of this function}

\begin{equation*}
V=x^4+y^4-axxy-axyy+ccxx+ccyy.
\end{equation*}
Having taken the differential, it will be

\begin{equation*}
P=4x^3-2axy-ayy+2ccx \quad \text{and} \quad Q=4y^3-axx-2axy+2ccy,
\end{equation*}
having put which equal to zero and having subtracted them from each other it will be

\begin{equation*}
4x^3-4y^3+axx-ayy+2ccx-2ccy=0;
\end{equation*}
since that equation is divisible by $x-y$, it will be $y=x$ at first and then $4x^3-3axx+2ccx=0$ which gives

\begin{equation*}
x=0 \quad \text{and} \quad 4xx=3ax-2cc \quad \text{and} \quad x=\frac{3a \pm \sqrt{9aa-32cc}}{8}.
\end{equation*}
If we take $x=0$, it will also be $y=0$ and, because of

\begin{equation*}
\frac{dP}{dx}=12xx-2ay+2cc \quad \text{and} \quad \frac{dQ}{dy}=12yy-2ax+2cc,
\end{equation*}
the function $V$ has a minimum value $=0$. If we set $x=y=\frac{3a \pm \sqrt{9aa-32cc}}{8}$, if it was $9aa >32cc$, because of $4xx= 3ax-2cc$, it will be

\begin{equation*}
\frac{dP}{dx}=\frac{dQ}{dy}= 12xx-2ax+2cc=7ax-4cc= \frac{21aa-32cc \pm 7a \sqrt{9aa -32cc}}{8};
\end{equation*}
because this value is always positive, because of $32cc <9aa$, the value $V$ will also in this case be a minimum and it will be

\begin{equation*}
V=- \frac{27}{256}a^4+\frac{9}{16}aacc-\frac{1}{2}c^4 \mp \frac{a}{256}(9aa-32cc)^{\frac{3}{2}}.
\end{equation*}
But let us divide the equation $4x^3-4y^3+axx-ayy+2ccx-2ccy=0$ by $x-y$ and it will be $4xx+4xy+4yy+ax+ay+2cc=0$. But from the equation $P=0$ it will be $yy=-2xy+\frac{4}{a}x^3+\frac{2ccx}{a}$, having substituted which value in that equation

\begin{equation*}
y=\frac{16x^3+4axx+aax+8ccx+2acc}{4ax-aa}.
\end{equation*}
But this gives

\begin{equation*}
y=-x \pm \sqrt{\frac{4x^3+axx+2ccx}{a}},
\end{equation*}
whence 

\begin{equation*}
16x^3+8axx+8ccx+2acc=(4x-a)\sqrt{4ax^3+aaxx+2accx},
\end{equation*}
which equation, having reduced it to a rational one, gives

\begin{alignat*}{19}
&256x^6+192ax^5 &&+80aa &&x^4&&+4a^3 &&x^3&&- a^2&&x^2 &&- 2a^3ccx+4a^2c^4=0; \\
&  &&+256cc && &&+160acc && &&+ 48aacc&& &&+32ac^4 \\
&  && && && && &&+ 64c^4&& && 
\end{alignat*}
whose real roots, if it has some, will indicate the maxima or minima of the function $V$, if $\frac{dP}{dx}$ and $\frac{dQ}{dy}$ become rational quantities affected with the same sign.

\subsection*{Example 5}

\textit{To find the maxima and minima of this expression}

\begin{equation*}
x^4+mxxyy+y^4+aaxx+naaxy+aayy=V.
\end{equation*}
After the differentiation it will be

\begin{alignat*}{19}
&P &&=4x^3 &&+2mxyy &&+2aax &&+naay &&=0, \\
&Q &&=4y^3  &&+2mxxy &&+2aay &&+naax &&=0,
\end{alignat*}
which equation either subtracted from or added to each other give

\begin{alignat*}{19}
&(4xx &&+4xy &&+4yy &&-2mxy &&+2aa &&-naa)&&(x-y) &&=0,\\
&(4xx &&-4xy &&+4yy &&+2mxy &&+2aa &&+naa)&&(x+y) &&=0,\\
\end{alignat*}
which divided by $x-y$ and $x+y$ and either added or subtracted again give

\begin{equation*}
4xx+4yy+2aa=0 \quad \text{and} \quad 4xy-2mxy-naa=0.
\end{equation*}
From the latter $y=\frac{naa}{2(2-m)x}$; but the first does not admit real values. Therefore, we have three cases.\\[2mm]
I. Let $y=x$ and it will be $4x^3+2mx^3+2aax+naax=0$, whence  either $x=0$ or $2(2+m)xx+(2+n)aa=0$. Let $x=0$; it will be $y=0$ and, because of

\begin{equation*}
\frac{dP}{dx}=12xx+2myy+2aa \quad \text{and} \quad \frac{dQ}{dy}=12yy+2mxx+2aa,
\end{equation*}
in this case $V=0$ will become a minimum, if the coefficient $aa$ was positive. The other case gives $xx=-\frac{(n+2)aa}{2(m+2)}$, which can only be real, if $\frac{n+2}{m+2}$ is a negative number. Let $\frac{n+2}{m+2}=-2kk$ or $n=-2kkm-4kk-2$; it will be $x=\pm ka$ and $y= \pm ka$. But

\begin{equation*}
\frac{dP}{dx}=12kkaa+2mkkaa+2aa \quad \text{and} \quad \frac{dQ}{dy}=12kkaa+2mkkaa+2aa;
\end{equation*} 
since they are equal, $V$ will have a maximum or a minimum value, depending on whether these quantities were positive or negative. \\[2mm]
II. Let $y=-x$ and it will be $2(m+2)x^3=(n-2)aax$; therefore, either $x=0$ or $xx=\frac{(n-2)aa}{2(m+2)}$. The first root $x=0$ reduces to the preceding. The second on the other hand will be real, if $\frac{(n-2)aa}{2(m+2)}$ was a positive quantity, and because $\frac{dP}{dx}=\frac{dQ}{dy}$, either maximum or a minimum will result.\\[2mm]
III. Let $y= \frac{naa}{2(2-m)x}$; it will be

\begin{equation*}
4x^3+ \frac{mn^2a^4}{2(2-m)^2 x}+2aax+\frac{nna^4}{2(2-m)x}=0 \quad \text{or} \quad 4x^4+2aaxx+\frac{nna^4}{(2-m)^2}=0,
\end{equation*} 
no root of which equation is real, if $aa$ is not a negative quantity.

\subsection*{Example 6}

\textit{Let this determined function be propounded $V=x^4+y^4-xx+xy-yy$, whose maximal or minimal values are to be investigated.}\\[2mm]
Since therefore $P=4x^3-2x+y=0$ and $Q=4y^3-2y+x=0$,  from the first it will be $y=2x-4x^3$ which, substituted in the other, gives

\begin{equation*}
256x^9-384x^7+192x^5-40x^3+3x=0.
\end{equation*}
One of its roots is $x=0$, whence  also $y=0$. Therefore, in this case, because of

\begin{equation*}
\frac{dP}{dx}=12xx-2 \quad \text{and} \quad \frac{dQ}{dy}=12yy-2,
\end{equation*}
a maximum value $V=0$ results. But having divided the found equation by $x$, it will be

\begin{equation*}
256x^8-384x^6+192x^4-40xx+3=0,
\end{equation*}
which has the factor $4xx-1$, whence $4xx=1$ and $x= \pm \frac{1}{2}$ and $y= \pm \frac{1}{2}$; then it will be $\frac{dP}{dx}=\frac{dQ}{dy}=1$; therefore, in each of the two cases a maximum  value $V=-\frac{1}{8}$ results. Divide that equation by $4xx-1$ and one will obtain

\begin{equation*}
64x^6-80x^4+28xx-3=0,
\end{equation*}
which again  contains $4xx-1=0$ twice such that the preceding case results. Furthermore, from this  $4xx-3=0$ and $x=\frac{\pm \sqrt{3}}{2}$;  $y=\frac{\mp \sqrt{3}}{2}$ corresponds to this. Therefore, it will also be $\frac{dP}{dx}=\frac{dQ}{dy}=7$ and therefore $V$ will have a minimum value $=-\frac{9}{8}$; this is the smallest value of all, which the function $V$ can have, and therefore, this equation $V=-\frac{9}{8}-cc$ is always impossible. But hence it is obvious, how to determine maxima and minima of functions containing three or more variables.


\title{On the Use of Differentials in the Investigation of the real Roots of Equations
\footnote{
Original title: "'De Usu Differentialium in ivestigandibus Radicibus realibus Aequationum"', first published as part of the book \textit{Institutiones calculi differentialis cum eius usu in analysi finitorum ac doctrina serierum}, 1755, reprinted in \textit{Opera Omnia}: Series 1, Volume 10, pp. 501 - 523, Eneström-Number E212, translated by: Alexander Aycock for the Euler-Kreis Mainz}}
\author{Leonhard Euler}
\date{ }
\maketitle

\paragraph*{§294}

The properties of maxima and minima provide us with a new method to find the nature of roots of equations, i.e., whether they are real or imaginary. For, let  an equation of any order be propounded

\begin{equation*}
x^n-Ax^{n-1}+Bx^{n-2}-Cx^{n-3}+Dx^{n-4}-\text{etc.}=0,
\end{equation*}
whose roots we put to be $p$, $q$, $r$, $s$, $t$ etc. such that $p$ is the smallest, $q$ the second smallest and likewise for the remaining roots, i.e. they are already ordered according to their magnitude; of course, $q>p$, $r>q$, $s>r$, $t>s$ etc. But let us assume that all roots of the equations are real and the largest exponent $n$ will at the same time be the number of roots $p$, $q$, $r$ etc. Additionally, let us assume all these roots to be different; this does not exclude the case of equal roots, since  non equal roots become equal, if their difference is assumed to be  infinitely small.

\paragraph*{§295}

Since the propounded expression $x^n-Ax^{n-1}+\text{etc.}$ is zero only, if  any of the values $p$, $q$, $r$ etc. is substituted for $x$, but does not vanish in all remaining cases, let us put in general

\begin{equation*}
x^n-Ax^{n-1}+Bx^{n-2}-Cx^{n-3}+\text{etc.}=z,
\end{equation*}  
such that $z$ can be considered as a function of $x$. Now let us assume that  determined values are substituted successively for $x$ starting from the smallest $x=- \infty$ and then going over to continuously larger values;  it is perspicuous that $z$ will hence have  values either greater than zero or smaller than zero and will just  vanish, if one puts $x=p$; in this case it will be $z=0$. Increase the values of $x$ beyond $p$ and the values of $z$ will either become  positive or negative until one gets to the value $x=q$; in this case it will again be $z=0$. Therefore, it is necessary that, because the values of $z$ from $0$ get back to $0$ again, $z$  had a maximum or minimum value within these limits, a maximum of course, if the values of $z$ were positive, while $x$ lies within the limits $p$ and $q$, a minimum, if they were negative. In like manner, while $x$ is increased, and becomes larger than $q$, up to $r$, the function $z$ will take on a maximum or a minimum value, a maximum value, of course, if it was a minimum value before, and vice versa. For, above [§ 263] we saw that the maxima and minima alternate.

\paragraph*{§296}

Hence, because within the limits determined by each two roots of $x$ there is a value, for which the function $z$ becomes a minimum or a maximum, the number of maxima and minimum values the function $z$ has, will be smaller than the number of real roots by $1$; and they will alternate that the maximal values of $z$ are positive and the minimal values are negative. If vice versa the function $z$ has a maximum or at least a positive value in the case $x=f$ and a minimum value or at least a negative value in the case $x=g$, since, if the values of $x$  go over from $f$ into $g$, the function $z$ goes over from the positive  into the negative values, it is necessary that in between it passes through $0$, and therefore a root of $x$ will be contained within the limits $f$ and $g$. But if this condition is not satisfied that the maximum and minimum values of $z$ alternately become positive and negative, that conclusion is incorrect. For, if  the function $z$ has minima which are also positive, it can happen that the value of $z$ goes over from a maximum  into the following minimum, although is does not vanish in between. Furthermore, from the things we explained it is understood, even though not all roots of the propounded equation were real, that nevertheless within the limits determined by each two roots there is always a maximum or a minimum, even though the converse proposition does not hold in general that there is always a real root within the limits determined by each two maxima or minima; but it holds, if the following condition is added: If the one value of $z$ was positive, the other is negative.

\paragraph*{§297}

Therefore, since we saw above that the values of $x$, for which the function

\begin{equation*}
z=x^n-Ax^{n-1}+Bx^{n-2}-Cx^{n-3}+Dx^{n-4}-\text{etc.}
\end{equation*}
has a maximum or a minimum value, are the roots of this differential equation

\begin{equation*}
\frac{dz}{dx}=nx^{n-1}-(n-1)Ax^{n-2}+(n-2)Bx^{n-3}-(n-3)Cx^{n-4}+\text{etc.}
\end{equation*}
it is obvious, if  all roots of the equation $z=0$, whose number is $=n$, were real, that also all roots of the equation $\frac{dz}{dx}=0$ will be real. For, because the function $z$ has as many maxima or minima, as the number $n-1$ contains units, the equation $\frac{dz}{dx}=0$ necessarily has the same amount of roots; and hence all its roots will be real. From this it is at the same time understood that the function $z$ cannot have more than $n-1$ maxima or minima. Therefore, we have this very far extending rule: If all roots of the equation $z=0$ were real, then also the equation $\frac{dz}{dx}=0$ will have only real roots. Hence it vice versa follows, if not all roots of the equation $\frac{dz}{dx}=0$ were real, that then also not all roots of the equation $z=0$ will be real.

\paragraph*{§298}

Since within the limits determined by each two real roots of the equation $z=0$ there is always one value for which the function $z$ has a maximum or a minimum value, it follows, if the equation $z=0$ has two real roots, that the equation $\frac{dz}{dx}=0$ will necessarily have one real root. In like manner, if the equation $z=0$ has three real roots, then the equation $\frac{dz}{dx}=0$ will certainly have two real roots. And in general, if the equation $z=0$ has $m$ real roots, it is necessary that at least $m-1$ roots of the equation $\frac{dz}{dx}=0$ are real. Hence, if the equation $\frac{dz}{dx}=0$ has less than $m-1$ real roots, then vice versa the equation $z=0$ has less  than $m$ real roots. But the converse is not true; for, even though the differential equation $\frac{dz}{dx}=0$ has several or even only real roots, it does nevertheless not follow that the equation $z=0$ will have even one real root. For, it can happen that all the roots of the equation $\frac{dz}{dx}=0$ real, although all roots of the equation $z=0$ are imaginary.

\paragraph*{§299}

Nevertheless, if the condition mentioned above is added, the converse proposition can be stated in such a way that from the real roots of the equation $\frac{dz}{dx}=0$ the number of real roots of the equation $z=0$  is known for certain. For, let us put that $\alpha$, $\beta$, $\gamma$, $\delta$ etc. are real roots of the equation $\frac{dz}{dx}=0$, of which $\alpha$ is  the largest; but the remaining follow in the order of their magnitude. Therefore, having substituted these values for $x$, the function $z$ will have either maximum or minimum values  alternately. But because the function $z$ becomes $=\infty$, if one puts $x=\infty$, it is plain that its values have to decrease continuously, while the values of $x$  are decreased from $\infty$ to $\alpha$; therefore,  $z$ will have a minimum value in the case $x=\alpha$. Therefore, if  in this case $x=\alpha$ the function $z$ has a negative value, it is necessary that   it was $=0$ somewhere else before, and so a real root $x >\alpha$ of the equation $z=0$ must exist; but if for $x=\alpha$ the function $z$ still has a positive value, it can not be smaller before; for, otherwise also a minimum would exist, before $x$ was decreased to $\alpha$, contradicting the hypothesis; hence the equation $z=0$ can have no real root larger than $\alpha$. Therefore, if we put that for $x=\alpha$ we have $z=\mathfrak{A}$, one can decide this way: If $\mathfrak{A}$ was a positive quantity,  the equation $z=0$ will have no real root larger than $\alpha$; but if $\mathfrak{A}$ was a negative quantity, the equation $z=0$ will always have one real root larger than $\alpha$ but not more.

\paragraph*{§300}

To make a further decision

\begin{center}
\begin{tabular}{l|l}
\quad if it is put \quad & \quad let \\
 \quad $x=\alpha$ & \quad $z=\mathfrak{A}$ \\
\quad $x=\beta$ & \quad $z=\mathfrak{B}$ \\
\quad  $x=\gamma$ & \quad $z=\mathfrak{C}$ \\
\quad $x=\delta$ & \quad $z=\mathfrak{D}$ \\
\quad $x=\varepsilon$ & \quad $z=\mathfrak{E}$ \\
\quad etc. & \quad etc.
\end{tabular}
\end{center} 
Therefore, because $\mathfrak{A}$ was a minimum, $\mathfrak{B}$ will be a maximum, and if $\mathfrak{A}$ was positive, also $\mathfrak{B}$ will be positive and therefore no real root of the equation $z=0$ will exist within the limits $\alpha$ and $\beta$ . Hence, if this equation has  no real root greater than $\alpha$, it will also not have one, which is greater than $\beta$. But if $\mathfrak{A}$ was a negative, in which case one root $x >\alpha$ of the equation is given, see, whether the value of $\mathfrak{B}$ is positive or negative. In the first case a there will be a root $x >\beta$, in the second there will be no root contained within the limits $\alpha$ and $\beta$. In like manner, because $\mathfrak{B}$ was a maximum, $\mathfrak{C}$ will be a minimum; hence, if $\mathfrak{B}$ had a negative value, $\mathfrak{C}$ will be a much more negative and in this case there will be no root contained within the limits $\beta$ and $\gamma$. But if $\mathfrak{B}$ was positive, there will be a real root  within the limits $\beta$ and $\gamma$, if $\mathfrak{C}$ becomes negative; but if $\mathfrak{C}$ is also positive, then there will be no root contained within the limits $\beta$ and $\gamma$  and in like manner the decision is to be made for the other cases.

\paragraph*{§301}

That these criteria for the decision seen more clearly, I summarized them in the following table. \\
\begin{center}
\begin{tabular}{c|c}
The equation $z=0$ will have one single real root & \\
 which is contained within the limits & if it was \\
\begin{tabular}{lll}
$x= \infty$ \quad & and \quad & $x= \alpha$ \\
$x= \alpha $ \quad & and \quad & $x= \beta$ \\
$x= \beta$ \quad & and \quad & $x= \gamma$ \\
$x= \gamma$ \quad & and \quad & $x= \delta$ \\
$x= \delta$ \quad & and \quad & $x= \varepsilon$ \\
 & etc. &
\end{tabular}
&
\begin{tabular}{lll}
$ \mathfrak{A}=-$ \quad &  & \\
$ \mathfrak{A}=-$ \quad & and \quad & $\mathfrak{B}=+$ \\
$ \mathfrak{B}=+$ \quad & and \quad & $\mathfrak{C}=-$ \\
$ \mathfrak{C}=-$ \quad & and \quad & $\mathfrak{D}=+$ \\
$ \mathfrak{D}=+$ \quad & and \quad & $\mathfrak{E}=-$ \\
                        & etc. &
\end{tabular}
\end{tabular}
\end{center}
The negated converses  of these propositions will hold in like manner. \\
\begin{center}
\begin{tabular}{c|c}
The equation $z=0$ will have no  real root & \\
 which is contained within the limits & if it was not\\
\begin{tabular}{lll}
$x= \infty$ \quad & and \quad & $x= \alpha$ \\
$x= \alpha $ \quad & and \quad & $x= \beta$ \\
$x= \beta$ \quad & and \quad & $x= \gamma$ \\
$x= \gamma$ \quad & and \quad & $x= \delta$ \\
$x= \delta$ \quad & and \quad & $x= \varepsilon$ \\
 & etc. &
\end{tabular}
&
\begin{tabular}{lll}
$ \mathfrak{A}=-$ \quad &  & \\
$ \mathfrak{A}=-$ \quad & and \quad & $\mathfrak{B}=+$ \\
$ \mathfrak{B}=+$ \quad & and \quad & $\mathfrak{C}=-$ \\
$ \mathfrak{C}=-$ \quad & and \quad & $\mathfrak{D}=+$ \\
$ \mathfrak{D}=+$ \quad & and \quad & $\mathfrak{E}=-$ \\
                        & etc. &
\end{tabular}
\end{tabular}
\end{center}
Therefore, by means of these rules, from the roots of the equation $\frac{dz}{dx}=0$, if they were known, not only the number of real roots of the equation $z=0$ is concluded, but also the limits will become known within which each root is contained.

\subsection*{Example}

\textit{Let this equation be propounded $x^4-14xx+24x-12=0;$ it is in question, whether it has real roots and how many.}\\[2mm]
The differential equation will be $4x^3-28x+24=0$ or $x^3-7x+6=0$, whose roots are $1$, $2$ and $-3$, which, ordered according to their magnitude, will give

\begin{center}
\begin{tabular}{r|l}
 &  whence it will be \\
 $\alpha$  $ =+2$ \quad & $\mathfrak{A}$ $=-4$ \\
 $\beta $  $ =+1$ \quad & $\mathfrak{B}$ $=-1$ \\
 $\gamma$  $ =-3 $ \quad & $\mathfrak{C}$ $=-129$.
\end{tabular}
\end{center}
Because of the negative $\mathfrak{A}$, the propounded equation will have a real root $>2$, but, because of the negative $\mathfrak{B}$, it will neither have a real root within the limits $2$ and $1$ nor within the limits $1$ and $-3$. But because for $x=-3$, $z=\mathfrak{C}=-129$, and if one sets $x=-\infty$,  $z=+\infty$, it is necessary that there is a real root  within the limits $-3$ and $-\infty$. Therefore, the propounded equation will have two real roots, the one $x>2$, the other $x<-3$; hence two roots will be imaginary. Therefore,  it has to be decided  from the last maximum or minimum of the propounded equation as from the first. If the propounded equation was of even order, the last maximum or minimum (it will be a minimum in this case), if it was negative, indicates a real root, if positive, an imaginary root. But for the equations of odd degree, since for $x=-\infty$ also $z=-\infty$, if the last maximum was positive, a real root is indicated, if negative, an imaginary one.

\paragraph*{§302}

Therefore, the rule for discovering the real and imaginary roots can be expressed conveniently this way. Having propounded any equation $z=0$, consider its differential $\frac{dz}{dx}=0$, whose real roots, ordered according to their magnitude, we want to be $\alpha$, $\beta$, $\gamma$, $\delta$ etc.; then, having put

\begin{alignat*}{9}
&x &&~=~&&\alpha, \quad && \beta, && \gamma, \quad && \delta, \quad && \varepsilon, \quad && \zeta \quad && \text{etc.} 
\intertext{let}
&z &&~=~&& \mathfrak{A}, \quad && \mathfrak{B}, \quad && \mathfrak{C}, \quad && \mathfrak{D}, \quad && \mathfrak{E}, \quad && \mathfrak{F}, \quad && \text{etc.} 
\intertext{Now, if the signs are}
 & && &&- &&+ &&- && + && - && + && \text{etc.},
\end{alignat*} 
the equation $z=0$ will have as many real roots as one has letters $\alpha$, $\beta$, $\gamma$ etc. and additionally one more. But if one of these capital letters does not have the mentioned sign written under it, two imaginary roots will be indicated. So, if $\mathfrak{A}$ has the sign $+$, there would be no root contained within the limits $\infty$ and $\beta$. If $\mathfrak{B}$ has the sign $-$, no root  will lie within the limits $\alpha$ and $\gamma$, and if $\mathfrak{C}$ has the sign $+$, there will be no root within the limits $\beta$ and $\delta$ and so forth. But in general except for the imaginary roots indicated this way the equation $z=0$ will additionally have as many imaginary roots as the equation $\frac{dz}{dx}=0$.

\paragraph*{§303}

If it happens that one of the values $\mathfrak{A}$, $\mathfrak{B}$, $\mathfrak{C}$, $\mathfrak{D}$ vanishes, then for it the equation $z=0$ will have two equal roots. If  $\mathfrak{A}=0$, then it will have two roots equal to $\alpha$; if  $\mathfrak{B}=0$, two roots will be equal $=\beta$. For, in this case the equation $z=0$ will have a root in common with the differential equation $\frac{dz}{dx}=0$; but above [§ 245] we demonstrated that this is a sign for two equal roots. But if the equation $\frac{dz}{dx}=0$ has two or more equal roots, then, if their number was even, neither a maximum nor a minimum will be indicated; therefore, for the present task, an even number of equal roots can be neglected. But if the number of equal roots of the equation $\frac{dz}{dx}=0$ was odd, all except for one are to be rejected in the decision, if not by accident in this case also the function $z$ itself vanishes. For, if this happens, the equation $z=0$ will also have equal roots and one more than the equation $\frac{dz}{dx}=0$. So, if it was $\frac{dz}{dx}=(x- \zeta)^nR$ such that this equation has $n$ roots equal to $\zeta$, if  $z$ vanishes for $x=\zeta$, the equation $z=0$ will have $n+1$ roots equal to $\zeta$.

\paragraph*{§304}

Let us apply these prescriptions to simpler equations and let us begin with the quadratic ones. Therefore, let this equation be propounded

\begin{equation*}
z=x^2-Ax+B=0;
\end{equation*}
its differential will be

\begin{equation*}
\frac{dz}{dx}=2x-A,
\end{equation*}
having put which $=0$ it will be

\begin{equation*}
x=\frac{1}{2}A \quad \text{or} \quad \alpha =\frac{1}{2}A.
\end{equation*}
Substitute this value for $x$ and it will be

\begin{equation*}
z=-\frac{1}{4}AA+B=\mathfrak{A};
\end{equation*}
hence, we conclude, if this value of $\mathfrak{A}$ was negative, i.e. if $AA>4B$, that the equation $xx-Ax+B=0$ will have two real roots, the one greater than $\frac{1}{2}A$, the other smaller. But if the value of $\mathfrak{A}$ was positive or $AA<4B$, then  both roots of the propounded equation will be imaginary. But if it was $\mathfrak{A}=0$ or $AA=4B$,  the propounded equation will have two equal roots, both of them $=\frac{1}{2}A$. Since these things are very well-known from the nature of quadratic equations, the validity of these principles is well-illustrated by this case and at the same time their utility in this task is understood.

\paragraph*{§305}

Let us proceed to the investigation of cubic equations in the same way. Therefore, let this equation be propounded

\begin{equation*}
x^3-Ax^2+Bx-C=z=0;
\end{equation*}
because its differential quotient is

\begin{equation*}
3xx-2Ax+B=\frac{dz}{dx},
\end{equation*}
if one puts this $=0$, it will be

\begin{equation*}
xx=\frac{2Ax-B}{3},
\end{equation*}
the two roots of which equation are either both  imaginary or both equal or both real and not equal. Therefore, because hence

\begin{equation*}
x=\frac{A \pm \sqrt{A^2-3B}}{3},
\end{equation*}
the two roots will be imaginary, if it was $AA<3B$; in this case the propounded cubic equation will have one single real root, which lies within limits $+\infty$ and $-\infty$, of course. Now let the two roots be equal to each other or let $AA=3B$; it will be $x=\frac{A}{3}$. Therefore, if  not $z=0$ at the same time, these two roots are to be considered as  one and the equation will have one single real root as before; but if in the case $x=\frac{A}{3}$  $z=0$ at the same time, what happens, if it was $-\frac{2}{27}A^3+\frac{1}{3}AB-C=0$ or $C=\frac{1}{3}AB- \frac{2}{27}A^3$, this means, if it was $B=\frac{1}{3}A^2$ and $C=\frac{1}{27}A^3$, the equation will have three equal roots, each of them $=\frac{1}{3}A$. Now let us expand the third case, in which both roots of the differential equation are real and different from each other, what happens, if $AA>3B$. Therefore, let $AA=3B+ff$ or $B=\frac{1}{3}AA-\frac{1}{3}ff$; those two roots will be

\begin{equation*}
x=\frac{A \pm f}{3}.
\end{equation*}
Therefore, it will be $\alpha =\frac{1}{3}A+\frac{1}{3}f$ and $\beta =\frac{1}{3}A- \frac{1}{3}f$. Therefore, find the values of $z$ corresponding  to these, i.e. $\mathfrak{A}$ and $\mathfrak{B}$, and because both roots are contained in this equation $xx=\frac{2}{3}Ax-\frac{1}{3}B$, it will be

\begin{equation*}
z=-\frac{1}{3}Axx+\frac{2}{3}Bx-C=-\frac{2}{9}AAx+\frac{1}{9}AB+\frac{2}{3}Bx-C.
\end{equation*}
Therefore,  these equations result

\begin{alignat*}{19}
&\mathfrak{A}&&=-\frac{2}{27}A^3 &&+\frac{1}{3}AB &&-\frac{2}{27}A^2f&&+\frac{2}{9}Bf&&-C&&=\frac{1}{27}A^3&&-\frac{1}{9}Aff&&-\frac{2}{27}f^3&&-C, \\[3mm]
&\mathfrak{B}&&=-\frac{2}{27}A^3 &&+\frac{1}{3}AB &&+\frac{2}{27}A^2f&&-\frac{2}{9}Bf&&-C&&=\frac{1}{27}A^3&&-\frac{1}{9}Aff&&+\frac{2}{27}f^3&&-C \\
\end{alignat*}
because of $B=\frac{1}{3}AA-\frac{1}{3}ff$. Therefore, if $\mathfrak{A}$ was a negative quantity, what happens, if  it was $C >\frac{1}{27}A^3-\frac{1}{9}Aff-\frac{2}{27}f^3$, the equation $z=0$ will have one single real root $>\alpha$, i.e. larger than $\frac{1}{3}A+\frac{1}{3}f$. Therefore, let us put that 

\begin{equation*}
C >\frac{1}{27}A^3-\frac{1}{9}Aff-\frac{2}{27}f^3 \quad \text{or that} \quad C=\frac{1}{27}A^3-\frac{1}{9}Aff-\frac{2}{27}f^3+gg
\end{equation*}
and, as we saw, the propounded cubic equation will have a real root $>\frac{1}{3}A+\frac{1}{3}f$. But of what nature the remaining roots are, will be understood from the value $\mathfrak{B}$; but it will be $\mathfrak{B}=\frac{4}{27}f^3-gg$; if it was positive, the equation will have two additional real roots, the first contained  within the limits $\alpha$ and $\beta$, i.e. within $\frac{1}{3}A+\frac{1}{3}f$ and $\frac{1}{3}A-\frac{1}{3}f$, but the one smaller than $\frac{1}{3}A-\frac{1}{3}f$. But if it was $gg>\frac{4}{27}f^3$ or $\mathfrak{B}$ was negative, the equation will have two imaginary roots. But if it was $\mathfrak{B}=0$ or $\frac{4}{27}f^3=gg$, the two roots will become equal, both of them $=\beta = \frac{1}{3}A-\frac{1}{3}f$. Finally, if the value of $\mathfrak{A}$ is positive or $C<\frac{1}{27}A^3-\frac{1}{9}Aff-\frac{2}{27}f^3$, the equations will have two imaginary roots and the third will be real and $>\frac{1}{3}A-\frac{1}{3}f$. And if the value of $\mathfrak{A}$ is $=0$, two roots will be equal and $=\alpha$, while the third remains $<\frac{1}{3}A-\frac{1}{3}f$.

\paragraph*{§306}

Therefore, for all three roots of the cubic equation $x^3-Ax^2+Bx-C=0$ to be real, three conditions are required. First, that 

\begin{equation*}
B < \frac{1}{3}AA;
\end{equation*}
therefore,  $B=\frac{1}{3}AA-\frac{1}{3}ff$. Secondly, that

\begin{equation*}
C > \frac{1}{27}A^3- \frac{1}{9}Aff-\frac{2}{27}f^3.
\end{equation*}
Thirdly, that 

\begin{equation*}
C< \frac{1}{27}A^3-\frac{1}{9}Aff+\frac{2}{27}f^3.
\end{equation*}
These two last inequalities reduce to the single one that $C$ is contained within the limits

\begin{equation*}
\frac{1}{27}A^3-\frac{1}{9}Aff-\frac{2}{27}f^3 \quad \text{and} \quad \frac{1}{27}A^3-\frac{1}{9}Aff+\frac{2}{27}f^3
\end{equation*}
or within these limits

\begin{equation*}
\frac{1}{27}(A+f)(A-2f) \quad \text{and} \quad \frac{1}{27}(A-f)^2(A+2f).
\end{equation*}
Therefore, if one of these conditions is missing, the equation will have two imaginary roots. So, if it was $A=3$, $B=2$, it will be $\frac{1}{3}ff=\frac{1}{3}AA-B=1$ and $ff=3$; therefore, this equation $x^3-3xx+2x-C=0$ can have only real roots, if $C$ is contained within the limits $-\frac{2\sqrt{3}}{9}$ and $+\frac{2\sqrt{3}}{9}$. Hence, if it was either $C<-\frac{2\sqrt{3}}{9}$ or $C<-0.3849$ or $C>+\frac{2\sqrt{3}}{9}$ or $C> 0.3849$ or together $CC>\frac{4}{27}$, the equation will have one single real root. 

\paragraph*{§307}

Since in each equation the second term can be thrown out, let us put that  $A=0$ such that we have this cubic equation

\begin{equation*}
x^3+Bx-C=0.
\end{equation*}
Therefore, for all three roots of this equation to be real, first it is necessary that  $B<0$ or $B$ must be a negative quantity. Therefore, let $B=-kk$; it will be $ff=3kk$ and additionally it is required that the quantity $C$ is contained within the limits $-\frac{2}{27}f^3$ and $+\frac{2}{27}f^3$, this means within  $-\frac{2}{9}kk\sqrt{3kk}$ and $+\frac{2}{9}kk\sqrt{3kk}$. Therefore, it will be $CC<\frac{4}{27}k^6$ or $CC<-\frac{4}{27}B^3$. Therefore, the nature of cubic equations, which have three real roots, can be expressed as one single condition, if we say that 

\begin{equation*}
4B^3+27CC
\end{equation*}
is a negative quantity. For, it is necessary that $B$ is a negative quantity, since otherwise $4B^3+27CC$ can not become negative.  Therefore, in general we can affirm that the equation $x^3+Bx\pm C=0$  only has real roots, if $4B^3+27CC$ was a negative quantity; but if this quantity was positive,  one will be real, the other two imaginary; but if $4B^3+27CC=0$,  all roots will be real, but two will be equal to each other.

\paragraph*{§308}

Let us proceed to fourth order equations, in which we want to assume the second term to be missing. Therefore, let

\begin{equation*}
x^4+Bx^2-Cx+D=0.
\end{equation*}
let us set $x=\frac{1}{u}$ and it will be

\begin{equation*}
1+Bu^2-Cu^3+Du^4=0,
\end{equation*}
the differential of which equation is

\begin{equation*}
2Bu-3Cu^2+4Du^3=0,
\end{equation*}
which has the one single root $u=0$; but then it will be

\begin{equation*}
uu=\frac{6Cu-4B}{8B}
\end{equation*}
and

\begin{equation*}
u=\frac{3C \pm \sqrt{9CC-32BD}}{8D}.
\end{equation*}
Therefore, for all four roots to be real,  first it is required that  $9CC>32BD$. Therefore, let us put that $9CC=32BD+9ff$; it will be $u=\frac{3C\pm 3f}{8D}$. Here, $C$ can always be assumed to be a positive quantity; for, if it was not such a one,  putting $u=-v$ it will become one. But soon we will demonstrate that not all roots can be real, if $B$ is not a negative quantity. Therefore, let $B=-gg$ and it will be
\begin{equation*}
9CC=9ff-32ggD \quad \text{and} \quad u=\frac{3C\pm 3f}{8D}.
\end{equation*}
Two cases are to be considered, depending on whether $D$ is a positive or negative quantity. \\[2mm]
I. Let $D$ be a positive quantity and it will be $f>C$ and the three roots of $u$, ordered according to their magnitude, will be

\begin{equation*}
1. ~~ u=\frac{3C+3f}{8D}, \quad 2. ~~ u=0, \quad 3. ~~ u=\frac{3C-3f}{8D}.
\end{equation*}
But the equation

\begin{equation*}
u^4-\frac{Cu^3}{D}+\frac{Bu^2}{D}+\frac{1}{D}=0,
\end{equation*}
having substituted these values for $u$, will give the following three values

\begin{equation*}
\mathfrak{A}=\frac{27(C+f)^3(C-3f)}{4096D^4}+\frac{1}{D}, \quad \mathfrak{B}=\frac{1}{D}, \quad \mathfrak{C}=\frac{27(C-f)^3(C+3f)}{4096}+\frac{1}{D},
\end{equation*}
of which the first and the third must be negative; both of them, because of the positive $C$ and $C<f$,  become smaller than $\frac{1}{D}$. Therefore, it must be

\begin{equation*}
\frac{1}{D}<\frac{27(C+f)^3(3f-C)}{4096D^4} \quad \text{and} \quad \frac{1}{D}<\frac{27(f-C)^3(C+3f)}{4096D^4}
\end{equation*}
or

\begin{equation*}
4096D^3 < 27(f+C)^3(3f-C) \quad \text{and} \quad 4096D^3<27(f-C)^3(C+3f).
\end{equation*}
But the first quantity is always a lot larger than the second; hence it suffices, if it was $D^3< \frac{27}{4096}(f-C)^3(C+3f)$, while $B=\frac{9CC-9ff}{32D}$ and $f>C$ and $D>0$. Therefore, if $D$ was a positive quantity, $C$ positive, $B$ negative, that  $f>C$, and $D^3< \frac{27}{4096}(f-C)^3(C+3f)$, i.e. $D<\frac{3}{16}(f-C)\sqrt[3]{3f+C}$,  the equation will have only real roots. But if it was $D>\frac{3}{16}(f-C)\sqrt[3]{3f+C}$, but nevertheless $D<\frac{3}{16}(f+C)\sqrt[3]{3f-C}$,  two roots will be real and two imaginary. But if it was $D>\frac{3}{16}(f+C)\sqrt[3]{3f-C}$, all four roots will be imaginary.\\[2mm]
Let $D$ be a negative quantity, say $=-F$, while $C$ remains positive and $B$ negative; because of $B=\frac{9CC-9ff}{32D}=\frac{9ff-9CC}{32F}$, it will be $C>f$. Therefore, because  $u=\frac{3C \pm 3f}{8D}=-\frac{3C\pm 3f}{8F}$, the three values, ordered according  to their magnitude, will be

\begin{equation*}
1. ~~ u=0, \quad 2. ~~ u=-\frac{3C-3f}{8F}, \quad 3. ~~ u=-\frac{3C+3f}{8F},
\end{equation*}
which will give the following values

\begin{equation*}
\mathfrak{A}=-\frac{1}{F}, \quad \mathfrak{B}=\frac{27(C-f)^3(C+3f)}{4096F^4}-\frac{1}{F}, \quad \mathfrak{C}=\frac{27(C+f)^3(C-3f)}{4096F^4}-\frac{1}{F}.
\end{equation*}
Therefore, since $\mathfrak{A}$ is a negative quantity, the equation will now certainly have one and therefore also two real roots. But for all roots to be real, it is necessary that $\mathfrak{B}$ is a positive quantity and therefore $27(C-f)^3(C+3f)>4096F^3$; but then it is necessary that $\mathfrak{C}$ is a negative quantity or $27(C+f)^3(C-3f)<4096F^3$. Therefore, for all roots to become real, it is required that $F^3$ is contained within the limits

\begin{equation*}
\frac{27}{4096}(C+f)^3(C-3f) \quad \text{and} \quad \frac{27}{4096}(C-f)^3(C+3f)
\end{equation*}
or that $F$ is contained within the limits

\begin{equation*}
\frac{3}{16}(C+f)\sqrt[3]{C-3f} \quad \text{and} \quad \frac{3}{16}(C-f)\sqrt[3]{C+3f};
\end{equation*}
And if $F$ is not contained within these limits, two roots will be imaginary.\\[2mm]
III. Now let us put that $B$ is a positive quantity and $D$ is positive as well; because of $B=\frac{9CC-9ff}{32D}$, it will be $C>f$, and because  $u=\frac{3C\pm 3f}{8D}$, the roots, ordered according  to their magnitude, will be

\begin{equation*}
1. ~~ u=\frac{3(C+f)}{8D}, \quad 2. ~~ u=\frac{3(C-f)}{8D} \quad \text{and} \quad u=0, 
\end{equation*}
whence the following values result

\begin{equation*}
\mathfrak{A}=\frac{27(C+f)^3(C-3f)}{4096D^4}+\frac{1}{D}, \quad \mathfrak{B}=\frac{27(C-f)^3(C+3f)}{4096D^4}+\frac{1}{D}, \quad \mathfrak{C}=\frac{1}{D};
\end{equation*}
since here $\mathfrak{C}$ is a positive quantity, two roots will certainly be imaginary. But if $\mathfrak{A}$ was negative, what happens, if $4096D^3 <27(C+f)^3(3f-C)$, two roots will be real; but if it was $4096D^3 >27(C+f)^3(3f-C)$, then all four roots will be imaginary. \\[2mm]
IV. Let $B$ remain positive, but let $D$ be negative $=-F$; because of $B=\frac{9ff-9CC}{32F}$, it will be $f>C$ and, because of $u=-\frac{3C\pm 3f}{8F}$, the three roots of $u$, ordered according to their magnitude, will be

\begin{equation*}
1. ~~ u=\frac{3(f-C)}{8F}, \quad 2. ~~ u=0 \quad \text{and} \quad 3. ~~ u=-\frac{3(C+f)}{8F},
\end{equation*}
whence these values result

\begin{equation*}
\mathfrak{A}=- \frac{27(f-C)^3(C+3f)}{4096F}-\frac{1}{F}, \quad \mathfrak{B}=-\frac{1}{F}, \quad \mathfrak{C}=-\frac{27(C+f)^3(3f-C)}{4096F^4}-\frac{1}{F},
\end{equation*}
where, because of the negative $\mathfrak{A}$ and $\mathfrak{C}$, the equation certainly has two real roots, but because  of the negative $\mathfrak{B}$ two roots will be imaginary.

\paragraph*{§309}

Therefore, if we put the letters $B$, $C$, $D$ to denote positive quantities, the following cases to be distinguished result, which, because of $f=\sqrt{CC-\frac{32}{9}BD}$, reduce to this.\\[2mm]
I. If the equation is $x^4-Bx^2\pm Cx+D=0$, all roots will be real, if it was

\begin{equation*}
D> \frac{3}{16}\left(\sqrt{CC+\frac{32}{9}BD}-C\right)\sqrt[3]{3\sqrt{CC+\frac{32}{9}BD}+C};
\end{equation*}
two roots will be real and two imaginary, if it was

\begin{equation*}
D> \frac{3}{16}\left(\sqrt{CC+\frac{32}{9}BD}-C\right)\sqrt[3]{3\sqrt{CC+\frac{32}{9}BD}+C},
\end{equation*}
but
\begin{equation*}
D< \frac{3}{16}\left(\sqrt{CC+\frac{32}{9}BD}+C\right)\sqrt[3]{3\sqrt{CC+\frac{32}{9}BD}-C};
\end{equation*}
but all roots will be imaginary, if it was

\begin{equation*}
D> \frac{3}{16}\left(\sqrt{CC+\frac{32}{9}BD}+C\right)\sqrt[3]{3\sqrt{CC+\frac{32}{9}BD}-C}.
\end{equation*}
II. If the equation is $x^4-Bx^2 \pm Cx-D=0$, two roots are always real; and the two remaining ones will also be real, if the quantity $D$ is contained within these limits

\begin{alignat*}{9}
&D &&> \frac{3}{16}\left(\sqrt{CC-\frac{32}{9}BD}+C\right)&&\sqrt[3]{C-3\sqrt{CC+\frac{32}{9}BD}} \\[2mm]
& D && < \frac{3}{16}\left(C-\sqrt{CC-\frac{32}{9}BD}\right)&&\sqrt[3]{C+3 \sqrt{CC-\frac{32}{9}BD}};
\end{alignat*}
But if $D$ is not contained within this limits, the two remaining roots will be imaginary.\\[2mm]
III. If the equation is $x^4+Bx^2\pm Cx+D=0$, two roots will always be imaginary; the remaining two will be real, if it was

\begin{equation*}
D <\frac{3}{16} \left(\sqrt{CC-\frac{32}{9}BD}+C\right)\sqrt[3]{3\sqrt{CC-\frac{32}{9}BD}-C};
\end{equation*}
but the two remaining ones will also be imaginary, if it was

\begin{equation*}
D >\frac{3}{16} \left(\sqrt{CC-\frac{32}{9}BD}+C\right)\sqrt[3]{3\sqrt{CC-\frac{32}{9}BD}-C}.
\end{equation*}
IV. If the equation is $x^4+Bx^2\pm Cx-D=0$, the two roots of this equation will always be real, the two remaining ones on the other hand will always be imaginary.

\subsection*{Example 1}

\textit{If this equation is propounded $x^4-2xx+3x+4=0$, let the nature of the roots be in question, i.e. let it be in question, whether they are real or imaginary.}\\[2mm]
Since this example extends to the first case, $B=2$, $C=3$ and $D=4$; hence

\begin{equation*}
CC+\frac{32}{9}BD=9+\frac{32 \cdot 8}{9}=\frac{337}{9} \quad \text{and} \quad \sqrt{CC+\frac{32}{9}BD}=\frac{\sqrt{337}}{3},
\end{equation*}
whence the conditions for all roots to be real are

\begin{alignat*}{9}
&4 &&< \frac{3}{16}&&\left(3+\frac{\sqrt{337}}{3}\right)&&\sqrt[3]{\sqrt{337}-3}&&=\frac{1}{16}(9+\sqrt{337})\sqrt[3]{\sqrt{337}-3},\\[2mm]
&4 &&< \frac{3}{16}&&\left(\frac{\sqrt{337}}{3}-3\right)&&\sqrt[3]{\sqrt{337}+3}&&=\frac{1}{16}(\sqrt{337}-3)\sqrt[3]{\sqrt{337}+3}.
\end{alignat*}
Using approximations it has therefore to be examined, whether  $4 <\frac{69}{16}$ and $4<\frac{24}{16}$; hence, because only the first condition holds, the equation will have two real and two imaginary roots.

\subsection*{Example 2}

\textit{Let this equation be propounded $x^4-9xx+12x-4=0$.}\\[2mm]
Since this example extends to the second case,  it will at least have two real roots. To investigate the nature of the remaining ones note that, because of $B=9$, $C=12$ and $D=4$, it will be 

\begin{equation*}
\sqrt{CC-\frac{32}{9}BD}=\sqrt{144-32\cdot 4}=4.
\end{equation*}
And therefore one has to check, whether 

\begin{alignat*}{9}
&4 &&> \frac{3}{16}\cdot 16 \sqrt[3]{0},&& \quad \text{i.e} \quad && 4&&>0, 
\intertext{and}
&4 &&< \frac{3}{16}\cdot 8 \sqrt[3]{24},&& \quad \text{i.e.} \quad && 4&&<3\sqrt[3]{3}; 
\end{alignat*}
because both is true, the propounded equation will have four real roots.

\subsection*{Example 3}

\textit{Let this equation be propounded $x^4+xx-2x+6=0$.}\\[2mm]
Since this equation extends to the third case, two roots will certainly be imaginary. But  $B=1$, $C=2$ and $D=6$ and hence

\begin{equation*}
\sqrt{CC-\frac{32}{9}BD}=\sqrt{4- \frac{64}{3}};
\end{equation*}
Since this quantity is imaginary, also the two remaining ones will certainly be imaginary.

\subsection*{Example 4}

\textit{Let this equation be propounded $x^4-4x^3+8x^2-16x+20=0$.}\\[2mm]
At first eliminate the second term; substituting $x=y+1$ it will be

\begin{center}
\begin{tabular}{rrrrrrrrrrrrrrr}
& &$x^4$ &$=$ & $y^4$ & $+$ & $4y^3$ & $+$ & $6yy$ & $+$ & $4y$ & $+$ & $1$ \\
& $-$ &$4x^3$ &$=$ &  & $-$ & $4y^3$ & $-$ & $12y^2$ & $-$ & $12y$ & $-$ & $4$ \\
&  $+$ &$8x^2$ &$=$ &  &  &  & $+$ & $8y^2$ & $+$ & $16y$ & $+$ & $8$ \\
 &  $-$ &$16x$ &$=$ &  &  &  &  &  & $-$ & $16y$ & $-$ & $16$ \\
  & $+$ &$20$ &$=$ &  &  &  &  &  &  &  & $+$ & $20$ \\ \cline{4-12}
Therefore & & & & $y^4$ & $+$ & $2yy$ & $-$ & $8y$ & $+$ & $9$ & $=$ & $0$;  
\end{tabular}
\end{center}
since it extends to the third case, it will have two imaginary roots. Then, because of $B=2$, $C=8$, $D=9$, it will be

\begin{equation*}
\sqrt{CC-\frac{32}{9}BD}=\sqrt{64-64}=0.
\end{equation*}
Therefore, compare $D=9$ to $\frac{3}{16}\cdot 8 \sqrt[3]{-8}=-3$. Therefore, because  $D=9 >-3$, also the two remaining roots will be imaginary.

\subsection*{Example 5}

\textit{Let this equation be propounded $x^4-4x^3-7x^2+34x-24=0$, whose roots are known to be $1$, $2$, $4$ and $-3$.}\\[2mm]
But if we apply the rules, having  removed the second term by putting $x=y+1$, it will be

\begin{equation*}
y^4-13yy+12y+0=0,
\end{equation*}
which compared to the second case gives $B=13$, $C=12$ and $D=0$. Therefore, it has to be $D> \frac{3}{16}\cdot 24\sqrt[3]{-24}$ or $0 > - 9\sqrt[3]{3}$ and $D<0$; therefore, because $D$ is not larger than $0$, the equation is indicated to have four real roots. For, if  $D=0$, the other equation will go over into

\begin{equation*}
D < \frac{3}{16}\left(\frac{16BD}{9C}\right)\sqrt[3]{4C} \quad \text{and hence} \quad 1 < \frac{B}{3C}\sqrt[3]{4C}
\end{equation*}
or $27CC<4B^3$; but $27 \cdot 144 < 4 \cdot 13^3$ or $36 \cdot 27>13^3$.

\paragraph*{§310}

It would be very difficult, to transfer these things to equations of higher degree, since the roots of the differential equations  cannot be exhibited in most cases; but if it is possible to assign these roots, from the given principles it is easily concluded, how many real and imaginary roots the propounded equation has. Hence the roots of all equations, which consist only of three terms, can be determined, whether they are real or imaginary. For, let this general equation be propounded

\begin{equation*}
x^{m+n}+Ax^n+B=0=z.
\end{equation*}
Take its differential

\begin{equation}
\frac{dz}{dx}=(m+n)x^{m+n-1}+nAx^{n-1};
\end{equation}
Having put it equal to zero, it will at first be $x^{n-1}=0$; hence, if $n$ was an odd number, a root exhibiting a maximum or a minimum results; but if $n$ is an even number, the  root $x=0$ is to be taken into account. But then it will be $(m+n)x^m+nA=0$; if $m$ is an even number and $A$ a positive quantity, this equation has no real root. Hence  the following cases are to be considered.\\[2mm]
I. Let $m$ be an even number and $n$ an odd number and the root $x=0$ will not exist. Therefore, if $A$ was a positive quantity, one will have no root exhibiting a maximum or a minimum; hence, because of the odd number $m+n$, the propounded equation will have one single real root. But if $A$ was a negative quantity, say $A=-E$, it will be $x=\pm \sqrt[m]{\frac{nE}{m+n}}$, whence

\begin{equation*}
 \alpha =+\sqrt[m]{\frac{nE}{m+n}} \quad \text{and} \quad \beta =- \sqrt[m]{\frac{nE}{m+n}}.
\end{equation*}
From these values 

\begin{equation*}
\mathfrak{A}=(x^m-E)x^n+B=- \frac{nE}{m+n}\left(\frac{nE}{m+n}\right)^{n:m}+B
\end{equation*}
and

\begin{equation*}
\mathfrak{B}=+\frac{mE}{m+n}\left(\frac{nE}{m+n}\right)^{m:n}+B.
\end{equation*}
Therefore, if $\mathfrak{A}$ was a negative quantity or

\begin{equation*}
\frac{mE}{m+n}\left(\frac{nE}{m+n}\right)^{n:m}>B,
\end{equation*}
the equation will have one single real root $>\alpha$. If  additionally 

\begin{equation*}
B> -\frac{mE}{m+n}\left(\frac{nE}{m+n}\right)^{n:m},
\end{equation*}
i.e., by combining both conditions into a single one, if it was

\begin{equation*}
(m+n)^{m+n}B^m< m^m n^n E^{m+n},
\end{equation*}
the equation will have three real roots, and if this condition is not satisfied, only one root of the equation will be real. These conditions are satisfied for the equation $x^{m+n}-Ex^n+B=0$, if $m$ was an even number and $n$ an odd number; if $E$ was a negative number here, the equation will always have one single real root.\\[2mm]
II. Let both numbers $m$ and $n$ be odd that $m+n$ is an even number and  the root $x=0$ is not to be taken into account. Since  $(m+n)x^m+nA=0$, it will be $x=-\sqrt[m]{\frac{nA}{m+n}}$; if one root is $=\alpha$, it will be

\begin{equation*}
\mathfrak{A}=\frac{mA}{m+n}x^n+B=-\frac{mA}{m+n}\left(\frac{nA}{m+n}\right)^{n:m}+B.
\end{equation*}
If this value was negative, the propounded equation will have two real roots, otherwise none. Therefore, the propounded equation $x^{m+n}+Ax^n+B=0$ will have two real roots, if it was

\begin{equation*}
m^m n^n A^{m+n} > (m+n)^{m+n}B^m;
\end{equation*}
if it was

\begin{equation*}
m^m n^n A^{m+n} < (m+n)^{m+n}B^m,
\end{equation*}
no root will be real. \\[2mm]
III. Let the two numbers $m$ and $n$ be even; likewise, $m+n$ will be an even number and the root $x=0$ will yield a maximum or a minimum; it will be the only one, if $A$ was a positive quantity, whence, having put $\alpha =0$, it will be $\mathfrak{A}=B$. Hence, if $B$ also was a positive quantity, the equation will have no real root; but if $B$ is also a negative quantity, one will have two real roots and not more, if $A$ was a positive quantity. But let us put that $A$ is a negative quantity or $A=-E$; it will be $x=\pm \sqrt[m]{\frac{nE}{m+n}}$ and  we will have three maxima or minima, namely

\begin{equation*}
\alpha = \pm \sqrt[m]{\frac{nE}{m+n}}, \quad \beta =0, \quad \gamma =- \sqrt[m]{\frac{nE}{m+n}}.
\end{equation*}
The following values correspond to these values of $z=x^{m+n}-Ex^n+B$ 

\begin{equation*}
\mathfrak{A}=-\frac{mE}{m+n}\left(\frac{nE}{m+n}\right)^{n:m}+B, \quad \mathfrak{B}=B, \quad  \mathfrak{C}=-\frac{mE}{m+n}\left(\frac{nE}{m+n}\right)^{n:m}+B.
\end{equation*}
Therefore, if $B$ is a negative quantity, because of the negative $\mathfrak{A}$ and $\mathfrak{C}$, the equation will have only two real roots, since also $\mathfrak{B}=B$ becomes negative. But if $B$ was a positive quantity, the equation will have four real roots, if

\begin{equation*}
(m+n)^{m+n}B^m < m^m n^n E^{m+n}.
\end{equation*}
But it will have no real root, if it was

\begin{equation*}
(m+n)^{m+n}B^m > m^n n^n E^{m+n}.
\end{equation*}
IV. Let $m$ be an odd number and $n$ an even number and the root $x=0$ will give a maximum or a minimum. Furthermore, it will be $x=-\sqrt[m]{\frac{nA}{m+n}}$. Therefore, if $A$ is a positive number, it will be $\alpha =0$ and $\beta = -\sqrt[m]{\frac{nA}{m+n}}$ and hence

\begin{equation*}
\mathfrak{A}=B \quad \text{and} \quad \mathfrak{B}=\frac{mA}{m+n}\left(\frac{nA}{m+n}\right)^{n:m}+B.
\end{equation*}
Hence, if $B$ is a negative quantity, say $B=-F$, and additionally

\begin{equation*}
m^mn^nA^{m+n}> (m+n)^{m+n}F^m,
\end{equation*}
the equation will have three real roots; otherwise only one will be real. But if $A$ is a negative quantity, say $A=-E$, it will be $x=+\sqrt[m]{\frac{nE}{m+n}}$ and 

\begin{equation*}
\alpha =\sqrt[m]{\frac{nE}{m+n}} \quad \text{and} \quad \beta=0,
\end{equation*}
to which these correspond

\begin{equation*}
\mathfrak{A}=- \frac{mE}{m+n}\left(\frac{nE}{m+n}\right)^{n:m}+B \quad \text{and} \quad \mathfrak{B}=B.
\end{equation*}
Hence the equation will have three real roots, if $B$ was a positive quantity and

\begin{equation*}
m^mn^nE^{m+n}>(m+n)^{m+n}B^m;
\end{equation*}
if this property is not present, the propounded equation will have one single root.

\paragraph*{§311}

Let all coefficients be $=1$ and while $\mu$ and $\nu$ denote integer numbers the decision for the following equations will be made this way:

\begin{equation*}
x^{2\mu+2 \nu -1}+x^{2 \nu-1} \pm 1=0
\end{equation*}
will have one single real root.
\begin{equation*}
x^{2 \mu +2 \nu -1}-x^{2 \nu -1} \pm 1 =0
\end{equation*}
will have three real roots, if it was

\begin{equation*}
(2 \mu +2 \nu -1)^{2 \mu +2 \nu -1}<(2 \mu)^{2 \mu}(2 \nu -1)^{2 \nu -1};
\end{equation*}
since this can never happen, the equation will always have one single real root.

\begin{equation*}
x^{2 \mu +2 \nu}\pm x^{2 \nu-1}-1=0
\end{equation*}
has two real roots. 

\begin{equation*}
x^{2\mu+2 \nu}\pm x^{2\nu -1}+1=0
\end{equation*}
has no real root.

\begin{equation*}
x^{2 \mu +2 \nu}\pm x^{2 \nu}+1=0
\end{equation*}
has no real root.
\begin{equation*}
x^{2\mu+2 \nu} +x^{2 \nu} \pm 1=0
\end{equation*}
has one single real root.
\begin{equation*}
x^{2\mu+2 \nu} -x^{2 \nu} \pm 1=0
\end{equation*}
has one single real root.\\[2mm]
Furthermore, since in the third case both exponents are even,   putting $xx=y$ it can be reduced to a simpler form and hence this case could have been omitted. Having done this, one will be able to affirm that no equation consisting of three terms can have more than three real roots.

\subsection*{Example}

\textit{Let the cases be in question in which this equation $x^5 \pm Ax^2 \pm B=0$ has three real roots.}\\[2mm]
Since this equation extends to the fourth case, it is plain that the quantities $A$ and $B$ must have opposite signs. Hence, if it does not have a form of this kind, it will have only one real root; but if the propounded equation was of this kind $x^5 \pm Ax^2 \mp B=0$, for it to have  three real roots, it is necessary that  $3^32^2A^5> 5^5 B^3$ or $A^5 \frac{3125}{108}B^3$. Therefore, if it was $B=1$, it is necessary that $A^5> \frac{3125}{108}$ or $A>{1.960132}$. Therefore, if $A=2$, this equation $x^5-2x^2+1=0$ has three real roots; since one of them is $x=1$, it follows that this fourth order equation $x^4+x^3+x^2-x-1=0$ has two real roots. This can be understood both from the given prescriptions and from the results demonstrated in the first part of the book; for, there we showed that any equation if even degree whose absolute term is negative always has two real roots.

\paragraph*{§312}

Using these principles one can also consider equations which consist of four terms, if the roots of the differential equations can be exhibited in a convenient way, what happens, if the exponents of $x$ in the three consecutive  terms following are terms in an arithmetic progression. But since this consideration  in general leads to several cases, let us treat them in some examples.

\section*{Example 1}

\textit{Let this equation be propounded $x^7-2x^5+x^3-a=0$.}\\[2mm]
Having put $z=x^7-2x^5+x^3-a$, it will be

\begin{equation*}
\frac{dz}{dx}=7x^6-10x^4+3xx,
\end{equation*}
having set which value equal to zero it will at first be $xx=0$, which double value is to be treated as none. But then it will be $7x^4=10x^3-3$, whence  $x^2=\frac{5 \pm 2}{7}$, and four values for $x$ will emerge, which, ordered according to their magnitude, will yield the following values for $z$:

\begin{center}
\begin{tabular}{lll|lll}
$\alpha$ & $=$ & $1$ \quad & $\mathfrak{A}$ & $=$ & $-a$ \\
$\beta$ & $=$ & $+\sqrt{\dfrac{3}{7}}$ \quad & $\mathfrak{B}$ & $=$ & $\dfrac{48}{343}\sqrt{\dfrac{3}{7}}-a$ \\
$\gamma$ & $=$ & $-\sqrt{\dfrac{3}{7}}$ \quad & $\mathfrak{C}$ & $=$ & $\dfrac{-48}{343}\sqrt{\dfrac{3}{7}}-a$ \\
$\delta$ & $=$ & $-1$ \quad & $\mathfrak{D}$ & $=$ & $-a.$ \\
\end{tabular}
\end{center}
Therefore, if $a$ is a positive number, it will be either $a>\frac{48}{343}\sqrt{\frac{3}{7}}$ or $a< \frac{48}{343}\sqrt{\frac{3}{7}}$; in the first case, because of the all negative quantities $\mathfrak{A}$, $\mathfrak{B}$, $\mathfrak{C}$, $\mathfrak{D}$,  the propounded equation will have one single root $x>1$. In the second case, if $a< \frac{48}{343}\sqrt{\frac{3}{7}}$, the equation will have three real roots, the first $>1$, the second contained within the limits $1$ and $\sqrt{\frac{3}{7}}$ and the third within the limits $+\sqrt{\frac{3}{7}}$ and $-\sqrt{\frac{3}{7}}$.\\[2mm]
But if $a$ is a negative quantity, putting $x=-y$, the equation will be reduced to the first form. Therefore, for the propounded equation to have three real roots it is necessary that  $a<0.0916134$ or $a<\frac{1}{11}$.

\subsection*{Example 2}

\textit{Let this equation be propounded $ax^8-3x^6+10x^3-12=0$.}\\[2mm]
Since here the exponents of the three last terms are terms  of an arithmetic progression, put $x=\frac{1}{y}$ and the equation will be transformed into this one

\begin{equation*}
a-3y^2+10y^5-12y^8=0;
\end{equation*}
therefore, put

\begin{equation*}
z=12y^8-10y^5+3y^2-a=0
\end{equation*}
and  by differentiation it will be

\begin{equation*}
\frac{dz}{dy}=96y^7-50y^4+6y=0,
\end{equation*}
from which equation at first  $y=0$; then, it will be

\begin{equation*}
y^6=\frac{50y^3-6}{96} \quad \text{and} \quad y^3=\frac{25 \pm 7}{96}
\end{equation*}
and hence either $y=\sqrt[3]{\frac{1}{3}}$ or $y=\sqrt[3]{\frac{3}{16}}$. Therefore, having ordered these values according to their magnitude, the corresponding values of $z$ will be as follows:

\begin{center}
\renewcommand{\arraystretch}{2,2}
\begin{tabular}{lll|lll}
$\alpha$ & $=$ & $\sqrt[3]{\dfrac{1}{3}}$ \quad & $\mathfrak{A}$ & $=$ & $\sqrt[3]{\dfrac{1}{9}}-a$ \\
$\beta$ & $=$ & $\sqrt[3]{\dfrac{3}{16}}$ \quad & $\mathfrak{B}$ & $=$ & $\dfrac{99}{64}\sqrt[3]{\dfrac{9}{256}}-a=\dfrac{99}{256}\sqrt[3]{\dfrac{9}{4}}-1$ \\
$\gamma$ & $=$ & $0$ \quad & $\mathfrak{C}$ & $=$ & $-a$. \\
\end{tabular}
\end{center}
Therefore, if it was $a>\sqrt[3]{\frac{1}{9}}$, the propounded equation will have two real roots, the one $>\sqrt[3]{\frac{1}{3}}$, the other $<0$; but except for these it will additionally have two real roots, if at the same time $\mathfrak{B}$ was a positive quantity, i.e., if it was $a<\frac{99}{256}\sqrt[3]{\frac{9}{4}}$. Therefore, the propounded equation will have four real roots, if the quantity $a$ is contained within the limits $\sqrt[3]{\frac{1}{9}}$ and $\frac{99}{256}\sqrt[3]{\frac{9}{4}}$; these limits are approximately $0.48075$ and $0.50674$. Therefore, having put $a=\frac{1}{2}$, this equation $x^8-6x^6+20x^3-24=0$ will have four real roots within the limits $\infty$, $\sqrt[3]{\frac{3}{16}}$, $\sqrt[3]{3}$, $0$, $-\infty$; therefore, three will be positive and one negative.


\title{On Criteria for imaginary roots
\footnote{
Original title: "'De Criteriis Radicum imaginarium"', first published as part of the book \textit{Institutiones calculi differentialis cum eius usu in analysi finitorum ac doctrina serierum}, 1755, reprinted in \textit{Opera Omnia}: Series 1, Volume 10, pp. 524 - 542, Eneström-Number E212, translated by: Alexander Aycock for the Euler-Kreis Mainz}}
\author{Leonhard Euler}
\date{ }
\maketitle

\paragraph*{§313}

In the preceding chapter we exhibited a method to explore the nature of roots of any equation such that by means of it, if any arbitrary equation is propounded, one can find out, how many real and imaginray roots it has. In most cases, this investigation is certainly most difficult, if the differential equation is of such a nature that its roots cannot be exhibited. But although in these cases the same operation could be applied to the differential equation itself and the nature of its roots could be explored from its differential and hence the roots of the original equation could be assigned approximately, the work would nevertheless be too cumbersome in almost every case. Therefore, in this case it often suffices to know criteria  from  which, if the their conditions are satisfied, one can conclude  that    the propounded equation contains imaginary roots, even though, if the conditions are not satisfied, one can vice versa not infer that all roots are real. Even if the knowledge is not complete then, it will be useful very often; therefore, we dedicated the present chapter to the explanation of these criteria.

\paragraph*{§314}

In the preceding chapter we saw, if any arbitrary equation

\begin{equation*}
z=x^n-Ax^{n-1}+Bx^{n-2}-Cx^{n-3}+Dx^{n-4}-\text{etc.}=0
\end{equation*}
has only real roots, that its differential

\begin{equation*}
\frac{dz}{dx}=nx^{n-1}-(n-1)Ax^{n-2}+(n-2)Bx^{n-3}-(n-3)Cx^{n-4}+\text{etc.}=0
\end{equation*}
will have also only real roots. But at the same time we showed, even though the differential equation only has real roots, hence it nevertheless does not follow that all roots of the propounded equation are real. Nevertheless, if the differential equation has imaginary roots, we will always be able to conclude that the propounded equation  must at least have as many imaginary roots, where \textit{at least} is to be stressed; for, it can happen that the equation has more imaginary roots. Therefore, this way  from the differential equation one can not conclude more than, if it has imaginary roots, that the propounded equation must also have roots of such a kind, and at least as many.

\paragraph*{§315}
If the propounded equation is multiplied by  any power $x^m$, while $m$ denotes a positive integer, because this new equation will have only real roots, if all roots of the propounded one were real, then also the roots of its differential, after having divided by $x^{m-1}$, will all be real. Hence, if this equation

\begin{equation*}
x^n-Ax^{n-1}+Bx^{n-2}-Cx^{n-3}+Dx^{n-4}-\text{etc.}=0
\end{equation*}
has only real roots, also this equation

\begin{equation*}
(m+n)x^n-(m+n-1)Ax^{n-1}+(m+n-2)Bx^{n-2}-\text{etc.}=0
\end{equation*}
will have only real roots. For the same reason, if this equation is multiplied by $x^k$ and differentiated again, the resulting equation

\begin{equation*}
(m+n)(k+n)x^n-(m+n+-1)(k+n-1)Ax^{n-1}+(m+n-2)(k+n-2)Bx^{n-2}-\text{etc.}=0
\end{equation*}
will still have only real roots and one can continue arbitrarily far like this. But if an equation of this kind is detected to have imaginary roots, it will be certain at the same time that the propounded equation will have at least as many imaginary roots.

\paragraph*{§316}

If the propounded equation, before it is differentiated, is multiplied by no power of $x$, the decision is to be made for an equation of one degree lower. So, if the propounded equation

\begin{equation*}
x^n-Ax^{n-1}+Bx^{n-2}-Cx^{n-3}+\text{etc.}=0
\end{equation*}
has only real roots,  all roots of  its differentials of all orders will also be real. Hence also the roots of all the following equations will be real
\begin{equation*}
nx^{n-1}-(n-1)Ax^{n-2}+(n-2)Bx^{n-3}-(n-3)Cx^{n-4}+\text{etc.}=0,
\end{equation*}
\begin{equation*}
n(n-1)x^{n-2}-(n-1)(n-2)Ax^{n-3}+(n-2)(n-3)Bx^{n-4}-\text{etc.}=0,
\end{equation*}
\begin{equation*}
n(n-1)(n-2)x^{n-3}-(n-1)(n-2)(n-3)Ax^{n-4}+\text{etc.}=0,
\end{equation*}
\begin{equation*}
n(n-1)(n-2)(n-3)x^{n-4}-(n-1)(n-2)(n-3)(n-4)AX^{n-5}+\text{etc.}=0
\end{equation*}
\begin{center}
etc.,
\end{center}
which equations are reduced to the following forms

\begin{alignat*}{19}
&x^{n-1} &&-\frac{n-1}{n}Ax^{n-2}&&+\frac{(n-1)(n-2)}{n(n-1)}Bx^{n-3}&&-\frac{(n-1)(n-2)(n-3)}{n(n-1)(n-2)}Cx^{n-4}&&+\text{etc.}&&=0, \\[3mm]
&x^{n-2} &&-\frac{n-2}{n}Ax^{n-3}&&+\frac{(n-2)(n-3)}{n(n-1)}Bx^{n-4}&&-\frac{(n-2)(n-3)(n-4)}{n(n-1)(n-2)}Cx^{n-5}&&+\text{etc.}&&=0, \\[3mm]
&x^{n-3} &&-\frac{n-3}{n}Ax^{n-4}&&+\frac{(n-3)(n-4)}{n(n-1)}Bx^{n-5}&&-\frac{(n-3)(n-4)(n-5)}{n(n-1)(n-2)}Cx^{n-6}&&+\text{etc.}&&=0, \\[3mm]
&x^{n-4} &&-\frac{n-4}{n}Ax^{n-5}&&+\frac{(n-4)(n-5)}{n(n-1)}Bx^{n-6}&&-\frac{(n-4)(n-5)(n-6)}{n(n-1)(n-2)}Cx^{n-7}&&+\text{etc.}&&=0\\
& && && && \text{etc.}
\end{alignat*}

\paragraph*{§317}

Therefore, this way the decision can be reduced to an equation of given lower degree than the propounded equation. So, if $m$ was any arbitrary number smaller than $n$, then, if the propounded equation has only real roots,   all roots of this equation of degree $m$ will also be real

\begin{equation*}
x^m-\frac{m}{n}Ax^{m-1}+\frac{m(m-1)}{n(n-1)}Bx^{m-2}-\frac{m(m-1)(m-2)}{n(n-1)(n-2)}Cx^{m-3}+\text{etc.}=0.
\end{equation*}
Hence, if one puts $m=2$, this equation will result

\begin{equation*}
x^2-\frac{2}{n}Ax+\frac{2 \cdot 1}{n(n-1)}B=0,
\end{equation*}
whose roots have to be real, if the propounded equation

\begin{equation*}
x^n-Ax^{n-1}+Bx^{n-2}-Cx^{n-3}+\text{etc.}=0
\end{equation*}
has only real roots. But because this quadratic equation can only have real roots, if  $\frac{AA}{nn}>\frac{2 \cdot 1}{n(n-1)}B$, it follows that all roots of the propounded equation can only be real, if $AA>\frac{2n}{2n-1}B$. Therefore, if it was $AA< \frac{2n}{2n-1}B$, this will be a certain indication that at least two roots of the propounded equation will be imaginary.

\paragraph*{§318}

Hence we derived a necessary condition, which the coefficients of the three first terms have to satisfy, if  all roots of the propounded equation were real. And this is a criterion of such kind  we mentioned at the beginning: Even though in the case $AA>\frac{2n}{n-1}B$ nothing follows for the realness of the roots, if  $AA<\frac{2n}{n-1}B$, it will nevertheless be a certain indication for the existence of at least two imaginary roots. So for all roots be real, by successively substituting the numbers $2$, $3$, $4$, $5$ etc. for $n$, it has to be as follows:

\begin{alignat*}{19}
&x^2&&-Ax   &&+B && =0 && && &&  && A^2 && >4 B \\[3mm]
&x^3&&-Ax^2 &&+Bx && -C &&=0  && && && A^2 && >\frac{6}{2} B \\[3mm]
&x^4&&-Ax^3 &&+Bx^2 && -Cx &&+D &&=0 &&   && A^2 && >\frac{8}{3} B \\[3mm]
&x^5&&-Ax^4 &&+Bx^3 && -Cx^2 &&+Dx &&-E &&=0 \quad \quad   && A^2 && >\frac{10}{4} B.
\end{alignat*}
Hence, if the second term is missing and the coefficient $B$ of the third is positive that the equation is of this kind

\begin{equation*}
x^n+Bx^{n-1}-Cx^{n-3}+Dx^{n-4}-\text{etc.}=0,
\end{equation*}
not all roots can be real, but at least two will be imaginary.

\paragraph*{§319}

Criteria of this kind can indeed be found for the coefficients of the following terms, if we consider that this equation

\begin{equation*}
1-Ay+By^2-Cy^3+Dy^4-\text{etc.}=0
\end{equation*}
has as many real and imaginary roots as the propounded equation. For, this equation results from the given one, if one puts $x=\frac{1}{y}$ such that from the roots of this equation one at the same time knows the roots  of the latter. Hence, if the propounded equation has only real roots,  also all roots of the differential equation of the reciprocal equation, i.e. of this one

\begin{equation*}
-A+2By-3Cy^2+4Dy^3-\text{etc.}=0,
\end{equation*}
will be real.  Substitute $x$ for $\frac{1}{y}$ in  this equation again and this equation will result

\begin{equation*}
Ax^{n-1}-2Bx^{n-2}+3Cx^{n-3}-4Dx^{n-4}+\text{etc.}=0,
\end{equation*}
whose roots will therefore be real, if the roots of the propounded equation were real. Hence it is now plain, if it was $n=3$, that necessarily $BB>3AC$.

\paragraph*{§320}

But now  differentiate this equation again and these equations will result

\begin{alignat*}{19}
&Ax^{n-2}&&-\frac{2(n-2)}{n-1}Bx^{n-3}&&+\frac{3(n-2)(n-3)}{(n-1)(n-2)}Cx^{n-4}&&-\text{etc.}&&=0\\[3mm]
&Ax^{n-3}&&-\frac{2(n-3)}{n-1}Bx^{n-4}&&+\frac{3(n-3)(n-4)}{(n-1)(n-2)}Cx^{n-5}&&-\text{etc.}&&=0\\[3mm]
&Ax^{n-4}&&-\frac{2(n-4)}{n-1}Bx^{n-5}&&+\frac{3(n-4)(n-5)}{(n-1)(n-2)}Cx^{n-6}&&-\text{etc.}&&=0\\
& && &&\text{etc.}
\end{alignat*}
Therefore, in general, if the number $m$ is smaller than $n$, it will be

\begin{equation*}
Ax^m-\frac{2m}{n-1}Bx^{m-1}+\frac{3m(m-1)}{(n-1)(n-2)}Cx^{m-2}-\text{etc.}=0.
\end{equation*}
If one now puts $m=2$, one will have this equation

\begin{equation*}
Ax^2-\frac{4}{n-1}Bx+\frac{6}{(n-1)(n-2)}C=0;
\end{equation*}
for its roots to be real it is necessary that  $\frac{4BB}{(n-1)^2}>\frac{6AC}{(n-1)(n-2)}.$ Hence, if the propounded equation has only real roots, it will be

\begin{equation*}
BB > \frac{3(n-1)}{2(n-2)}AC.
\end{equation*}
And if it was $BB<\frac{3(n-1)}{2(n-2)}AC$, this is a certain sign that the propounded equation has at least two imaginary roots. Therefore, if $n=3$, the criterion will be $BB>3AC$; but if $n=4$, it will be $BB> \frac{3 \cdot 3}{2 \cdot 2}AC$; if $n=5$, it will be $BB>\frac{3 \cdot 4}{2 \cdot 3}AC$ and so forth.

\paragraph*{§321}

To transfer these criteria to the following coefficients, let us consider the differential equation expressed in $y$ again 

\begin{equation*}
-A+2By-3Cy^2+4Dy^3-5Ey^4+\text{etc.}=0
\end{equation*}
and let us differentiate  it once more that we have

\begin{equation*}
2B-6Cy+12Dy^2-20Ey^3+\text{etc.}=0,
\end{equation*}
which, having substituted $\frac{1}{x}$ for $y$ again, will give

\begin{equation*}
Bx^{n-2}-3Cx^{n-3}+6Dx^{n-4}-10Ex^{n-5}+\text{etc.}=0,
\end{equation*}
from whose further differentiation this equation follows

\begin{equation*}
Bx^{n-3}-\frac{3(n-3)}{n-2}Cx^{n-4}+\frac{6(n-3)(n-4)}{(n-2)(n-3)}Dx^{n-5}-\text{etc.}=0
\end{equation*}
and in general

\begin{equation*}
Bx^m-\frac{3m}{n-2}Cx^{m-1}+\frac{6m(m-1)}{(n-2)(n-3)}Dx^{m-3}-\text{etc.}=0.
\end{equation*}
Therefore, if we put $m=2$, this quadratic equation will result

\begin{equation*}
Bx^2-\frac{2 \cdot 3}{n-2}Cx+\frac{6 \cdot 2}{(n-2)(n-3)}D=0,
\end{equation*}
whose roots will be real, if it was $\frac{)CC}{(n-2)^2}>\frac{6 \cdot 2 BD}{(n-2)(n-3)}$ or

\begin{equation*}
CC> \frac{4(n-2)}{3(n-3)}BD.
\end{equation*}
Hence, if the propounded equation has only real roots, it will be $CC> \frac{4(n-2)}{3(n-3)}BD$, and if this condition is not satisfied, the equation will certainly have at least two imaginary roots.

\paragraph*{§322}

If we differentiate the above equation $2B-6Cy+12Dy^2-\text{etc.}=0$ once again, this equation will result

\begin{equation*}
-6C+24Dy-60Ey^2+\text{etc.}=0
\end{equation*}
or

\begin{equation*}
C-4Dy+10Ey^2-20Fy^3+\text{etc.}=0,
\end{equation*}
which having substituted $x$ for $\frac{1}{y}$ again will go over into this one

\begin{equation*}
Cx^{n-3}-4Dx^{n-4}+10Ex^{n-5}-20Fx^{n-6}+\text{etc.}=0,
\end{equation*}
from whose further differentiation these  equations follow

\begin{alignat*}{9}
&Cx^{n-4}&&-\frac{4(n-4)}{n-3}Dx^{n-5}&&+\frac{10(n-4)(n-5)}{(n-3)(n-4)}Ex^{n-6}&&-\text{etc.}&&=0,\\[3mm]
&Cx^{n-5}&&-\frac{4(n-5)}{n-3}Dx^{n-6}&&+\frac{10(n-5)(n-6)}{(n-3)(n-4)}Ex^{n-7}&&-\text{etc.}&&=0
\end{alignat*}
and in general

\begin{equation*}
Cx^m-\frac{4m}{n-3}Dx^{m-1}+\frac{10m(m-1)}{(n-3)(n-4)}Ex^{m-2}-\text{etc.}=0.
\end{equation*}
Let us put $m=2$ and it will be

\begin{equation*}
Cx^2-\frac{2 \cdot 4}{n-3}Dx+\frac{2 \cdot 10}{(n-3)(n-4)}E=0,
\end{equation*}
from which, if its roots are real, it follows that 

\begin{equation*}
\frac{4 \cdot 4}{(n-3)^2}>\frac{2 \cdot 10}{(n-3)(n-4)}CE \quad \text{or} \quad DD>\frac{5(n-3)}{4(n-4)}CE.
\end{equation*}

\paragraph*{§323}

From these results one already clearly sees the relation for all coefficients. Therefore, in general, if this equation

\begin{equation*}
x^n-Ax^{n-1}+Bx^{n-2}-Cx^{n-3}+Dx^{n-4}-Ex^{n-5}+\text{etc.}=0
\end{equation*}
has only real roots, it will be

\begin{alignat*}{9}
&AA &&> \frac{2n}{1(n-1)}&&B \\[3mm]
&BB &&> \frac{3(n-1)}{2(n-2)}&&AC \\[3mm]
&CC &&> \frac{4(n-2)}{3(n-3)}&&BD \\[3mm]
&DD &&> \frac{5(n-3)}{4(n-4)}&&CE \\[3mm]
&EE &&> \frac{6(n-4)}{5(n-5)}&&DF \\[3mm]
&   && \text{etc.}
\end{alignat*}
If one of these conditions is not satisfied, the equation will have at least two imaginary roots. And if these criteria do not depend on each other, it is easily seen that there are as many pairs of imaginary roots   as the total amount of non-satisfied conditions. But even if these conditions all hold in one single equation, it hence nevertheless does not follow that no imaginary roots are given; it can even happen, because there is no reason against it, that all roots are imaginary. Therefore, one has to be careful that not more is attributed to these criteria than it can actually attributed to them considering the principles whence they were deduced.

\paragraph*{§324}

But it is easily seen that not each condition, which is not satisfied, can indicate imaginary roots; for, in an equation of $n$ dimensions, since one has $n+1$ terms and from  each, except for the first and the last, a criterion can be derived, one will in total have $n-1$ conditions; and nevertheless, if they are not satisfied, the equation can not have $2n-2$ imaginary roots, since it in total has only $n$ roots. But one condition alone always reveals two imaginary roots, and since it can happen that two conditions of this kind do not show more roots, one has to consider, whether these two conditions indicate the same imaginary root or not; in the first case the number of imaginary roots will not increase, in the second on the other hand, since the conditions involve completely different letters, each one of them will show two imaginary roots. So, even though it was

\begin{equation*}
AA<\frac{2n}{1(n-1)}B \quad \text{and} \quad BB<\frac{3(n-1)}{2(n-2)}AC,
\end{equation*}
hence nevertheless not necessarily four imaginary roots are indicated, but both of them might indicate the same two roots. On the other hand, if it was

\begin{equation*}
AA< \frac{2n}{1(n-1)B} \quad \text{and} \quad CC<\frac{4(n-2)}{3(n-3)}BD
\end{equation*}
while $BB>\frac{3(n-1)}{2(n-2)}AC$, four imaginary roots will be indicated.

\paragraph*{§325}
 
Therefore, from several criteria involving consecutive letters for imaginary roots it does not follow more than from one; but if they proceed in an interrupted order that  between each two at least one other criterion was skipped, then from each one of them one can conclude two imaginary roots. This consideration yields the following rule. Except for the first and the last term write the coefficients found  from  the criteria before over the respective terms of the propounded equation this way

\begin{alignat*}{9}
& &&\frac{2n}{1(n-1)}\quad &&\frac{3(n-1)}{2(n-2)}\quad && \frac{4(n-2)}{3(n-3)}\quad && \frac{5(n-3)}{4(n-4)}\quad &&\text{etc.} \\[3mm]
&x^n~-~&&Ax^{n-1}~+~&&Bx^{n-2}~-~&&Cx^{n-3}~+~&&Dx^{n-4}~-~&&\text{etc.}=0 \\
&~+~&&\cdots &&\cdots &&\cdots &&\cdots &&\text{etc.} \\
\end{alignat*}
Then examine the square of each coefficient, whether it is larger or smaller than the fraction written above it multiplied by the product of the  corresponding coefficients; in the first case   attribute the sign $+$ to the term, in the second the sign $-$; but always attribute term the sign $+$ to the first and the last. Having done this, the equation will have at least as many imaginary roots as there are variations of the  signs.

\paragraph*{§326}

This is the rule found by Newton to explore the imaginary roots of each equation; but one has to pay attention, as we already mentioned, since it can happen  that the equation has more imaginary roots than one detects by means of this method. Hence others tried to find other similar rules which  would yield the number of imaginary roots more exactly, such that the true number of roots would exceed the number which the rule gives less often.  In this regard, especially the  rule of Campbell added to Newtons \textit{Universal Arithmetic} stands out, which will therefore be conveniently explained here, even if it is not perfect. It is based on this lemma: If $\alpha$, $\beta$, $\gamma$, $\delta$, $\varepsilon$ etc. were some quantities and their total amount is $m$, put the sum of these quantities

\begin{equation*}
\alpha +\beta +\gamma +\delta +\text{etc.}=S,
\end{equation*}
the sum of the squares

\begin{equation*}
\alpha^2+ \beta^2 +\gamma^2 +\delta^2+\text{etc.}=V,
\end{equation*}
and it will be $V>0$. But because the product of each two is

\begin{equation*}
\alpha \beta +\alpha \gamma +\alpha \delta +\beta \gamma +\beta \delta +\text{etc.}=\frac{SS-V}{2},
\end{equation*}
it will be $(m-1)V>SS-V$ or $mV>SS$. For, if the squares of the differences of two quantities are taken, their sum will be

\begin{equation*}
=(\alpha -\beta)^2+(\alpha -\gamma)^2+(\alpha -\delta)^2+(\beta -\gamma)^2+(\beta -\delta)^2+\text{etc.}
\end{equation*}
\begin{equation*}
=(m-1)(\alpha^2 +\beta^2+\gamma^2+\delta^2+\text{etc.})-2(\alpha \beta +\alpha \gamma +\alpha \delta +\beta \gamma +\text{etc.})
\end{equation*}
\begin{equation*}
=(m-1)V-2\frac{SS-V}{2}=mV-SS.
\end{equation*}
Therefore, since the sum of real squares is always positive, it will be

\begin{equation*}
mV-SS>0 \quad \text{and hence} \quad mV>SS.
\end{equation*}

\paragraph*{§327}

Having stated this lemma in advance, if one has this equation

\begin{equation*}
x^n-Ax^{n-1}+Bx^{n-2}-Cx^{n-3}+Dx^{n-4}-Ex^{n-5}+Fx^{n-6}-\text{etc.}=0
\end{equation*}
and all its $n$ roots were real, which we want to be $a$, $b$, $c$, $d$, $e$ etc., it will be, as it is known from the nature of equations,

\begin{center}
\renewcommand{\arraystretch}{2,2}
\begin{tabular}{l|l}
& \quad numbers of terms \\
$A=a+b+c+d+\text{etc.}$ \quad & \quad $n$ \\
$B=ab+ac+ad+bc+bd+\text{etc.}$ \quad & \quad $\dfrac{n(n-1)}{1 \cdot 2}$ \\
$C=abc+abd+abe+ acd+bcd+\text{etc.}$ \quad & \quad $\dfrac{n(n-1)(n-2)}{1 \cdot 2 \cdot 3}$ \\
$B=abcd+abce+abde+\text{etc.}$ \quad & \quad $\dfrac{n(n-1)(n-2)(n-3)}{1 \cdot 2 \cdot 3 \cdot 4}$ \\
\end{tabular}
\end{center}
\begin{center}
etc.
\end{center}
Now take the squares of each term and put

\begin{alignat*}{9}
&P &&=a^2+b^2+c^2+d^2+\text{etc.}, \\
&Q &&=a^2b^2+a^2c^2+a^2d^2+b^2c^2+\text{etc.}, \\
&R &&=a^2b^2c^2+a^2b^2d^2+a^2b^2e^2+a^2c^2d^2+\text{etc.}, \\
&S &&=a^2b^2c^2d^2+a^2b^2c^2e^2+a^2b^2d^2e^2+\text{etc.}\\
&  &&\text{etc.};
\end{alignat*}
 but from the nature of combinatorics it will be

\begin{alignat*}{9}
&P &&=A^2&&-2B, \\
&Q&&=B^2&&-2AC &&+2D, \\
&R&&=C^2 &&-2BD &&+2AE &&-2F, \\
&S &&=D^2&&-2CE &&+2BF &&-2AG &&+2H \\
& && && &&\text{etc.}
\end{alignat*}

\paragraph*{§328}

By means of the lemma stated in advance we will therefore have

\begin{equation*}
nP >AA,
\end{equation*}
\begin{equation*}
\frac{n(n-1)}{1 \cdot 2}Q > BB, 
\end{equation*}
\begin{equation*}
\frac{n(n-1)(n-2)}{1 \cdot 2 \cdot 3}R> CC,
\end{equation*}
\begin{equation*}
\frac{n(n-1)(n-2)(n-3)}{1 \cdot 2 \cdot 3 \cdot 4}S >DD
\end{equation*}
\begin{center}
etc.
\end{center}
Therefore, if  the values found before are substituted for the values $P$, $Q$, $R$ etc., we will obtain the following properties of the real roots

\begin{equation*}
nAA-2nB > AA \quad \text{or} \quad AA>\frac{2n}{n-1}B,
\end{equation*}
\begin{equation*}
\frac{n(n-1)}{1 \cdot 2}BB- \frac{2n(n-1)}{1 \cdot 2}AC+\frac{2n(n-1)}{1 \cdot 2}D>BB
\end{equation*}
or

\begin{equation*}
BB > \frac{\frac{2n(n-1)}{1 \cdot 2}}{\frac{n(n-1)}{1 \cdot 2}-1}(AC-D)
\end{equation*}
and in like manner the following equations yield

\begin{alignat*}{9}
&CC &&> \frac{\frac{2n(n-1)(n-2)}{1 \cdot 2 \cdot 3}}{\frac{n(n-1)(n-2)}{1 \cdot 2 \cdot 3}-1}(BD-AE+F), \\[3mm]
&DD &&> \frac{\frac{2n(n-1)(n-2)(n-3)}{1 \cdot 2 \cdot 3 \cdot 4}}{\frac{n(n-1)(n-2)(n-3)}{1 \cdot 2 \cdot 3 \cdot 4}-1}(CE-BF+AG-H).
\end{alignat*}
Therefore, the square of each coefficient is not only compared to the product of the closest terms, but also to the rectangles of two equally distant ones, nevertheless in such a way that the signs of these rectangles alternate.

\paragraph*{§329}

Therefore,  except for the first and the last  one has to write the fractions, whose numerators are the binomial coefficients of the same power  multiplied  by two and whose denominators are the same binomial coefficients decreased by $1$, above the respective terms of the equation. So, by considering quadratic, cubic, fourth order equations etc., if their roots  all were real, it will be

\begin{alignat*}{9}
&   && \tfrac{4}{1} \\
&x^2~-~ &&Ax~+~ &&B=0; \quad A^2>4B;
\end{alignat*}
For the cubic equation

\begin{alignat*}{9}
& && \tfrac{6}{2} && \tfrac{6}{2} \\
&x^3~-~&&Ax^2~+~&&Bx~-~&&C=0
\end{alignat*}
it will be

\begin{equation*}
A^2 >3B \quad \text{and} \quad B^2 >3AC.
\end{equation*}
For the forth order equation

\begin{alignat*}{9}
& && \tfrac{8}{3}&& \tfrac{12}{5} && \tfrac{8}{3} \\
&x^4~-~&&Ax^3~+~&&Bx^2~-~&&Cx~+~&&=0
\end{alignat*}
it will be

\begin{equation*}
A^2 >\frac{8}{3}B, \quad B^2>\frac{12}{5}(AC-D), \quad C^2 >\frac{8}{3}BD.
\end{equation*}
For the equation of the fifth order

\begin{alignat*}{9}
& && \tfrac{10}{4}&& \tfrac{20}{9}&& \tfrac{20}{9}&& \tfrac{10}{4}\\
&x^5 ~-~ &&Ax^4~+~&&Bx^3~-~&&Cx^2~+~&&Dx~-~&&E=0
\end{alignat*}
it will be

\begin{equation*}
AA>\frac{10}{4}B, \quad B^2> \frac{20}{9}(AC-D), \quad C^2>\frac{20}{9}(BD-AE) \quad \text{and} \quad D^2 >\frac{10}{4}CE.
\end{equation*}
For the equation of the sixth order

\begin{alignat*}{9}
&  &&\tfrac{12}{5}&& \tfrac{30}{14} && \tfrac{40}{19} && \tfrac{30}{14} && \tfrac{12}{5} \\
&x^6~-~&&Ax^5~+~&&Bx^4~-~&&Cx^3~+~&&Dx^2~-~&&Ex~+~&&F=0
\end{alignat*}
it will be

\begin{equation*}
A^2>\frac{12}{5}B, \quad B^2 >\frac{30}{14}(AC-D), \quad C^2>\frac{40}{19}(BD-AE+F),
\end{equation*}
\begin{equation*}
D^2 >\frac{30}{14}(CE-BF), \quad E^2>\frac{12}{5}DF.
\end{equation*}
\begin{center}
etc.
\end{center}

\paragraph*{§330}

Therefore, if a certain condition is not satisfied, it will be an indication that at  least two  imaginary roots are contained in the propounded equation. But because, if each condition is not satisfied, the equation can not have twice as many imaginary roots, one has to argue in these cases as  before in the case of Newton's rule.  If the square of a certain term was larger than the fraction written above it multiplied by the product of the closest and equally distant terms, then  attribute the sign $+$ to this term, otherwise the sign $-$; but  always attribute the sign $+$ to the first and the last term. Having done this, check the progression of the signs, and as often as a variation occurs, an imaginary root will be indicated. Therefore, if this rule indicates more imaginary roots than Newton's rule, it will be closer to the truth. Nevertheless, it can happen that the equation has more imaginary roots than each of the rules indicates.

\paragraph*{§331}

Therefore, we would make a mistake, if we wanted to use these criteria as perfect indication for real and imaginary roots, since it can happen that the equation has more imaginary roots than these criteria indicate; and the error could be the greater, the higher the degree of the propounded equation was. For, in the case of quadratic equations these criteria are true in such a way that, if they do  not indicate any imaginary roots,  the equation will also have none. But the cubic equation can have two imaginary roots, even though both rules (they coincide in this case) do not exhibit them. Therefore, to someone wanting to investigate these cases, let this general cubic equation be propounded

\begin{alignat*}{9}
& &&3 && 3 \\
&x^3~-~ &&Ax^2~+~&&Bx~-~&&C=0;
\end{alignat*}
if in this it was $AA>3B$ and $BB>3AC$, none of both rules indicates imaginary roots. But above (§ 306) we saw that for  this equation to have imaginary roots, at first it is required that  $B<\frac{1}{3}AA$, which condition also both rules require. Therefore, let $B=\frac{1}{3}AA-\frac{1}{3}FF$ and it is necessary that $C$ is contained within these limits

\begin{equation*}
\frac{1}{27}A^3-\frac{1}{9}AFF-\frac{2}{27}f^2 \quad \text{and} \quad \frac{1}{27}A^3-\frac{1}{9}Aff+\frac{2}{27}f^3.
\end{equation*}
But both rules only demand that  $C<\frac{BB}{3A}$, i.e.

\begin{equation*}
C<\frac{1}{27}A^3-\frac{2}{27}Aff+\frac{f^4}{27A}.
\end{equation*}
This condition can be satisfied, even though $C$ is not contained within the mentioned limits.

\paragraph*{§332}

For, let 

\begin{equation*}
C=\frac{1}{27}A^3-\frac{2}{27}Aff+\frac{f^4}{27A}-gg
\end{equation*}
and the rules will indicate no imaginary roots. There will nevertheless be imaginary roots, if it was either

\begin{equation*}
\frac{1}{27}A^3-\frac{2}{27}Aff+\frac{f^4}{27A}-gg < \frac{1}{27}A^3-\frac{1}{9}Aff-\frac{2}{27}f^3
\end{equation*}
or

\begin{equation*}
\frac{1}{27}A^3-\frac{2}{27}Aff+\frac{f^4}{27A}-gg > \frac{1}{27}A^3-\frac{1}{9}Aff+\frac{2}{27}f^3.
\end{equation*}
Therefore, if it was either

\begin{equation*}
gg > \frac{(ff+Af)^2}{27A} \quad \text{or} \quad gg< \frac{(Af-ff)^2}{27A},
\end{equation*}
the cubic equation will have two imaginary roots, even though none of both rules indicates them. But here we assumed that $A$ is a positive quantity; for, if it was negative, by putting $x=-y$ the equation would be transformed into a form in which $A$ would be positive. Hence  one can form infinitely many cubic equations, which have two imaginary roots, even though they are not indicated by the rule. For, let $gf=\frac{(ff+Af)^2}{27A}+hh$; it will be

\begin{equation*}
C=\frac{(ff-AA)^2}{27A}-gg= \frac{1}{27}A^3-\frac{1}{9}Aff-\frac{2}{27}f^3-hh \quad \text{and} \quad B=\frac{1}{3}AA- \frac{1}{3}ff.
\end{equation*}
Or let it be $gg=\frac{(Af-ff)^2}{27A}-hh$ with $hh< \frac{(Af-ff)^2}{27A}$; it will be

\begin{equation*}
C=\frac{1}{27}A^3-\frac{1}{9}Aff+\frac{2}{27}f^3+hh \quad \text{and} \quad B=\frac{1}{3}AA-\frac{1}{3}ff.
\end{equation*}
In both cases an equation having two imaginary roots, which none of both rules indicate, will result. For the sake of an example let us put $A=4$, $f=1$; it will be $B=5$ and, because of $gg=\frac{25}{108}+hh$, it will be

\begin{equation*}
C=\frac{225}{108}-\frac{25}{108}-hh=\frac{50}{27}-hh.
\end{equation*}
Hence, if  $C<\frac{50}{27}$, the equation $x^3-4x^2+5x-C=0$ will always have two imaginary roots. But having taken $gg=\frac{1}{12}-hh$, it must be $hh<\frac{1}{12}$ and it will be

\begin{equation*}
C=\frac{25}{12}-\frac{1}{12}+hh=2+hh.
\end{equation*}
Let $hh=\frac{1}{16}$ and the equation $x^3-4x^2+5x-\frac{33}{16}=0$ will have two imaginary roots, even though none is revealed by the rules.

\paragraph*{§333}

It is even possible to form general equations of such a kind, in which none of both rules exhibits imaginary roots, even though in most cases two or more are contained in the equation. This happens, if two equal signs alternate, as in

\begin{equation*}
x^n-Ax^{n-1}-Bx^{n-2}+Cx^{n-3}+Dx^{n-4}-Ex^{n-5}-Fx^{n-6}+\text{etc.}=0
\end{equation*}
or

\begin{equation*}
x^n+Ax^{n-1}-Bx^{n-2}-Cx^{n-3}+Dx^{n-4}+Ex^{n-5}-Fx^{n-6}-\text{etc.}=0;
\end{equation*}
here, both rules do not reveal an imaginary root. But that they most often can contain roots of this kind, is also clear from the cubic equation $x^3-Ax^2-Bx+C=0$, which for $ff=AA+3B$ always has two imaginary roots, if it was either

\begin{equation*}
-C <\frac{1}{27}A^3-\frac{1}{9}Aff-\frac{2}{27}f^3 \quad \text{or} \quad -C> \frac{1}{27}A^3-\frac{1}{9}Aff+\frac{2}{27}f^3.
\end{equation*}
Nevertheless, also these cases can be found from the rules, if the equation  is transformed into another form by means of a substitution. Put $x=y+k$ and it will be

\begin{equation*}
\left.
\setlength{\arraycolsep}{0.0mm}
\begin{array}{llllllllllllll}
y^3 &~+~&3ky^2  &~-~3kky  &~+~&k^3 \\
    &~-~&Ayy   &~-~&2Aky &~-~&Akk \\
    &   &      &~-~&By   &~-~&Bk \\
    &   &      & &       &~+~&C
\end{array}
\right\rbrace =0,
\end{equation*}
which examined according to the rule will first immediately give

\begin{equation*}
(3k-A)^2 > 3(3kk-2Ak-B);
\end{equation*}
but for it to be

\begin{equation*}
(3kk-2Ak-B)^2 >3(3k-A)(k^3-Akk-Bk+C),
\end{equation*}
which is the other criterion, it is necessary that

\begin{equation*}
BB+3AC+(AB-9C)k+(AA+3B)kk>0,
\end{equation*}
whatever value is attributed to $k$. Therefore, choose $k$ in such a way that this expression has minimum value, what will happen for $k=\frac{9C-AB}{2(AA+3B)}$, and if this expression was still $>0$, it will be probable that the propounded equation has no imaginary roots. But it will be

\begin{equation*}
BB+3AC-\frac{(AB-9C)^2}{2(AA+3B)}+\frac{(AB-9C)^2}{4(AA+3B)}>0
\end{equation*}
or

\begin{equation*}
BB+3AC >\frac{(AB-9C)^2}{4(AA+3B)}.
\end{equation*}
Therefore, since $B=\frac{1}{3}ff-\frac{1}{3}AA$, it will be

\begin{equation*}
4ff\left(\frac{1}{9}f^4-\frac{2}{9}AAff+\frac{1}{9}A^4+3AC\right)>\left(\frac{1}{3}Aff-\frac{1}{3}A^3-9C\right)^2
\end{equation*}
or

\begin{equation*}
4f^6-(A^2f^4+4A^4ff+108ACff>A^2f^4-2A^4f^2-54ACff+A^6+54A^3C+729CC
\end{equation*}
or

\begin{equation*}
4f^6 > 9A^2 f^4-6A^4 ff- 162 ACff+A^6+54A^3C+729CC,
\end{equation*}
whence, having factored the equation, it will have to be

\begin{equation*}
(2f^3+A^3-3Af^2+27C)(2f^3-A^3+54A^3C+27C)>0.
\end{equation*}
And hence the rules will show imaginary roots, if it was either

\begin{alignat*}{9}
&C &&>-\frac{1}{27}A^3+\frac{1}{9}Af^2-\frac{2}{27}f^3 \quad && \text{and} \quad && C &&>-\frac{1}{27}A^3+\frac{1}{9}Af^2+\frac{2}{27}f^3 \\[3mm]
&C &&<-\frac{1}{27}A^3+\frac{1}{9}Af^2-\frac{2}{27}f^3 \quad && \text{and} \quad && C &&<-\frac{1}{27}A^3+\frac{1}{9}Af^2+\frac{2}{27}f^3.
\end{alignat*}
These are the same conditions which we found above [§ 306]. Therefore, it is plain that by means of an appropriate transformation the rules given in this chapter can be stated in such a way that they are always true, even though they are converted.

\paragraph*{§334}

From these principles also Harriot's rule can be demonstrated, according to which any arbitrary equation is predicted to have as many positive roots as there are variations of the signs, but as many negative as there are successions of the same sign; but this rule only holds for real roots. Therefore, let us put that the equation

\begin{equation*}
x^n-Ax^{n-1}+Bx^{n-2}-Cx^{n-3}+Dx^{n-4}-\text{etc.}=0
\end{equation*}
has only real and positive roots and its differential

\begin{equation*}
nx^{n-1}-(n-1)Ax^{n-2}+(n-2)Bx^{n-3}-\text{etc.}=0
\end{equation*}
will have not only also  only real and positive roots, but  the roots of this one will also constitute the limits of the roots of the propounded equation. Furthermore, having put $x=\frac{1}{y}$, this equation

\begin{equation*}
1-Ay+By^2-Cy^3+Dy^4-\text{etc.}=0
\end{equation*}
will have only real positive roots, but they are the reciprocals of the others, such that the roots, which are the maxima in that equation, are the minima in this one. Having constituted all this, if that propounded equation is continuously differentiated until one gets to an equation of first order, which will be $x-\frac{1}{n}A=0$ (§ 317), the root of this one will still be positive and hence the  coefficient of the second term will have the sign $-$, as we assumed. But if this coefficient would have the sign $+$,  it would follow that the propounded equation has not only positive real roots, but at least one will be negative, of course the one, which  corresponds to the mentioned limits.

\paragraph*{§335}

If the propounded equation is converted into its reciprocal and is differentiated, then $x$ is substituted again and the differentiations are repeated until one gets to a simple equation, which from  § 320 will be of this kind $Ax-\frac{2}{n-1}B=0$, its roots must therefore also be positive, if the propounded equation has only real and positive roots, and hence the second and the third term will have different signs. Therefore, if these two terms have the same sign, at least one negative root will be indicated corresponding to the  limit assigned in this equation, which will be different from the limit indicated by the equation, since here they were once converted into its reciprocals; hence one concludes, if the three initial terms of the equation have equal signs, that  two negative roots will be indicated.

\paragraph*{§336}

In like manner, if the conversions and differentiations are done according to § 320 and are continued until one gets to the simple equation $Bx-\frac{3}{n-2}C=0$, also the roots of this equation must be positive, if all roots of the propounded equation were positive, of course; hence, if the third and the fourth term have equal signs, one negative root will be indicated. And so forth, if any two contiguous terms have the same sign, one negative root will be indicated; and hence, no matter how many successions of the same sign occur, the propounded equation will have at least as many negative roots, since each of these criteria refer to different limits. But if the propounded equation is put to have only negative roots, then, because the roots of all differential equations deduced from it  must also be negative, all terms need to have the same sign. Hence, if two contiguous terms have different signs, from them at least one positive root will be concluded. And in like manner, no matter how many variations of the signs of two terms occur in the propounded equation, at least as many positive roots are concluded to be contained in the equation. Therefore, since each equation has as many roots as there are combinations of two contiguous signs, and not more, it follows that each equation, whose roots are all real, has as many positive roots as there were variations of contiguous signs, but has as many negative roots as there were iterations of the same sign.


\title{On Differentials of Functions in certain cases only
\footnote{
Original title: "'De Differentialibus Functionum in certis tantum Casibus"', first published as part of the book \textit{Institutiones calculi differentialis cum eius usu in analysi finitorum ac doctrina serierum}, 1755, reprinted in \textit{Opera Omnia}: Series 1, Volume 10, pp. 542 - 563, Eneström-Number E212, translated by: Alexander Aycock for the Euler-Kreis Mainz}}
\author{Leonhard Euler}
\date{ }
\maketitle

\paragraph*{§337}

If $y$ was any function of $x$ and this variable quantity is increased by $\omega$ that $x$ goes over into $x+ \omega$, the function will have this value

\begin{equation*}
y+ \frac{\omega dy}{dx}+\frac{\omega^2 ddy}{2 dx^2}+\frac{\omega ^3 d^3 y}{6 dx^3}+\frac{\omega^4 d^4 y}{24 dx^4}+\text{etc.}
\end{equation*}
and hence will receive this increment

\begin{equation*}
\frac{\omega dy}{dx}+\frac{\omega^2 ddy}{2dx^2}+\frac{\omega^3 d^3 y}{6dx^3}+\frac{\omega^4 d^4 y}{24 dx^4}+\text{etc.},
\end{equation*}
as we showed above [§ 48]. Therefore, if  $\omega =dx$ such that $x$ grows by the amount of its differential $dx$, the function $y$ will receive the increment

\begin{equation*}
=dy +\frac{1}{2}ddy+\frac{1}{6}d^3y +\frac{1}{24}d^4 y+ \text{etc.},
\end{equation*}
which will be the true differential of $y$. But since any arbitrary term of this series has an infinite ratio to the following, with respect to the first all vanish such that  $dy$ taken in usual manner yields the true differential of $y$. In like manner the true second, third, fourth etc. differentials of $y$ will be as follows

\begin{alignat*}{19}
&dd.y &&= ddy &&+\frac{3}{3}&&d^3y &&+ \frac{7}{3 \cdot 4} &&d^4y&&+ \frac{15}{3 \cdot 4 \cdot 5}&&d^5y &&+\frac{31}{3 \cdot 4 \cdot 5 \cdot 6}&&d^6y &&+\text{etc.} \\[3mm]
&d^3. y &&= d^3 y &&+\frac{6}{4}&&d^4y &&+\frac{25}{4 \cdot 5}&&d^5y &&+\frac{90}{4 \cdot 5 \cdot 6}&&d^6 y &&+\frac{301}{4 \cdot 5 \cdot 6 \cdot 7}&&d^7 y &&+\text{etc.}\\[3mm]
&d^4. y &&=d^4y &&+\frac{10}{5}&&d^5 y &&+\frac{65}{5 \cdot 6}&&d^6 y &&+\frac{350}{5 \cdot 6 \cdot 7}&&d^7y &&+\frac{1701}{5 \cdot 6 \cdot 7 \cdot 8}&&d^8 y &&+\text{etc.} \\[3mm]
&d^5. y &&= d^5 y &&+\frac{15}{6}&&d^6 y &&+\frac{140}{6 \cdot 7}&&d^7 y &&+\frac{1050}{6 \cdot 7 \cdot 8}&&d^8 y&&+\frac{6951}{6 \cdot 7 \cdot 8 \cdot 9}&&d^9 y &&+\text{etc.} \\[3mm]
&d^6. y && =d^6 y &&+\frac{21}{7}&&d^7 y &&+\frac{266}{7 \cdot 8}&&d^8 y &&+\frac{2646}{7 \cdot 8 \cdot 9}&&d^9y &&+\frac{22827}{7 \cdot 8 \cdot 9 \cdot 10}&&d^{10}y &&+\text{etc.}\\
& && && && && &&\text{etc.}
\end{alignat*}
which follow from § 56, if one writes $dx$ instead of $\omega$. Therefore, exactly those differentials of $y$ will be complete, in which not even the terms  vanishing with respect to the first are neglected. But these  terms are found, if the function $y$ is continuously differentiated and $dx$ was put to be constant. So having put $y=ax-xx$, because of

\begin{equation*}
dy=adx- 2xdx \quad \text{and} \quad ddy=-2dx^2,
\end{equation*}
the complete differentials of $y$ will be

\begin{equation*}
dy=adx-2xdx -dx^2, \quad ddy=-2dx^2;
\end{equation*}
the following are all zero.

\paragraph*{§338}

Although in general in these expressions of the differentials the following terms with respect to the first are considered to be zero,  nevertheless in special cases, in which the first term  vanishes, this assumption is not valid and the second term can not  be neglected anymore. Therefore, even though in the preceding example  the differential of the formula $y=ax-xx$ in general is $=(a-2x)dx$, having neglected the term $-dx^2$, which certainly is infinite times smaller than the first $(a-2x)dx$, here nevertheless this condition that the first term does not vanish per se is considered to be satisfied. Therefore, if the differential of $y=ax-xx$ is in question in the case  $x=\frac{1}{2}a$,  one has to say that actually $=-dx^2$; if the variable $x$ increases by the differential $dx$,  the decrement of the function $y$ in the case $x=\frac{1}{2}a$ will be $dx^2$. But, having excluded this single case, the differential of the function $y$ will always be $=(a-2x)dx$; for, if not $x=\frac{1}{2}a$, the second term $-dx^2$, with respect to the first, is always justly neglected. And the negligence of the term $dx^2$  cannot induce an error even in the case $x=\frac{1}{2}a$; for, usually the first differentials are  compared to each other; hence, because $dy=-dx^2$ in the case $x=\frac{1}{2}a$ vanishes in comparison  to the first differentials $dx$, it does not matter, whether in this case we have $dy=0$ or $dy=-dx^2$.

\paragraph*{§339}

While $y$ denotes any function of $x$,  having taken the differentials several times, let

\begin{equation*}
dy= pdx, \quad dp=qdx, \quad dq=rdx, \quad dr=sdx \quad \text{etc.}
\end{equation*}
Hence, the complete differentials of $y$ in which nothing is neglected  will be

\begin{alignat*}{19}
&d. y &&= pdx &&+ \frac{1}{2}q&&dx^2 &&+\frac{1}{6}r &&dx^3 &&+\frac{1}{24}s &&dx^4 &&+\frac{1}{120}t&& dx^5 &&+\text{etc.} \\[3mm]
&d^2.y &&= qdx^2 &&+r &&dx^3 &&+\frac{7}{12}s &&dx^4 &&+\frac{1}{4}t&& dx^5 &&+\text{etc.} \\[3mm]
&d^3.y &&= rdx^3 &&+\frac{3}{2}s&&dx^4&&+\frac{5}{4}t &&dx^5 &&+\text{etc.} \\[3mm]
&d^4. y && = sdx^4 &&+2 t && dx^5 &&+\text{etc.} \\[3mm]
&d^5. y && =tdx^5 &&+\text{etc.} \\
&\text{etc.}
\end{alignat*}
Therefore, if the first terms of these expressions do not vanish, they alone will exhibit the differentials of $y$; but, if in a certain case the first term becomes $=0$,  the following term will express the differential in question. And if  the second term also vanishes,  the third term will yield the value of the differential in question; but if even this term vanishes, the fourth, and so forth. Therefore, it is understood that the first differential of any function of $x$ never vanishes completely; for, even though $p=0$, in which case $sy$ is usually considered to vanish, this differential will then be expressed by  a higher power of $dx$,  as, e.g., either by $\frac{1}{2}qdx^2$ or, if also $q=0$, by  $\frac{1}{6}rdx^3$, and so forth.

\paragraph*{§340}

But although in these cases the differential of $y$, with respect to higher first differentials  it is compared to, is justly neglected and considered to be zero, it is nevertheless often helpful to also know its true expression. For, from the complete form of the differential  it is immediately seen in which cases the given function has a maximum or minimum value. For, if it was

\begin{equation*}
d.y = pdx +\frac{1}{2}qdx^2+\frac{1}{6}rdx^3+\text{etc.},
\end{equation*}
for $y$ to have a maximum or minimum value it is necessary that  $p=0$; therefore, in this case it will be $dy=\frac{1}{2}qdx^2$ and the function $y$, if  one writes $x\pm dx$ instead of $x$, goes over into $y+\frac{1}{2}qdx^2$ and will therefore have a minimum value, if $q$ has a positive value, but maximum value, if $q$ has a negative value. But if at the same time $q=0$, it will be $dy=\frac{1}{6}rdx^3$ and the function  will go over into $y \pm \frac{1}{6}rdx^3$ by writing $x \pm dx$ instead of $x$ and in this case neither a maximum nor a minimum value results; but if also $r=0$, then, having written $x \pm dx$ instead of $x$, the function $y$ will become $=y+\frac{1}{24}sdx^4$, which exhibits a maximum, if $s$ was a negative quantity, a minimum on the other hand, if $s$ is a positive quantity. Other occasions in which the complete expression of the differentials have a use will occur below.

\paragraph*{§341}

Let us put that $p$ vanishes in the case $x=a$, what happens, if it was  $p=(x-a)P$. But such a value results, if it was

\begin{equation*}
y=(x-a)^2 P+C
\end{equation*}
while $C$ denotes any constant quantity. For, because 

\begin{equation*}
pdx =(x-a)^2 dP+2(x-a)Pdx,
\end{equation*}
it will certainly be $p=0$ for $x=a$. Therefore,  because of

\begin{equation*}
dp dx = qdx^2 =(x-a)^2 ddP+4(x-a) dPdx+2 Pdx^2,
\end{equation*}
having put $x=a$, it will be $qdx^2 =2 Pdx^2$ and the complete differential in this case $x=a$ will be

\begin{equation*}
dy= Pdx^2,
\end{equation*}
if not by accident also $P$ vanishes for $x=a$ which cases I will contemplate later. \\[2mm]
But the present case can be exhibited more generally this way. Let 

\begin{equation*}
z=(x-a)^2P+C
\end{equation*}
and let $y$ be any function of $z$ such that $dy=Zdz$ while $Z$ denotes any function of $z=(x-a)^2P+C$. Therefore, it will be

\begin{equation*}
dz=(x-a)^2dP+2(x-a)Pdx \quad \text{and} \quad pdx =Z(x-a)^2 dP+2Z(x-a)Pdx,
\end{equation*}
which term becomes $=0$, if  $x=a$; and in the same case, having neglected the terms containing the factor $x-a$, it will be $qdx^2= 2 PZdx^2$ and hence in the case $x=a$ it will be $dy=PZdx^2$, after in $PZ$ it was put $a$ for $x$ everywhere. Hence, if $y$ was any function of $z=(x-a)^2P+C$ such that $dy=Zdz$, in the case $x=a$ the differential will be

\begin{equation*}
dy=PZdx^2.
\end{equation*}
Therefore, this function $y$ has maximum value in the case $x=a$, if in the same case $PZ$ was a negative quantity, a minimum value on the other hand, if $PZ$ was a positive quantity.

\paragraph*{§342}

If it was $p=(x-a)^2P$, in the case $x=a$ also $q$ vanishes; but such an expression  results for $p$, if it was

\begin{equation*}
y=(x-a)^3P+C.
\end{equation*}
Therefore, it will be 

\begin{equation*}
pdx= (x-a)^2 dP+3(x-a)Pdx,
\end{equation*}
\begin{equation*}
qdx^2 =(x-a)^3 ddp +6(x-a)^2dPdx+6(x-a)Pdx^2,
\end{equation*}
 both sides of which vanish in the case $x=a$; but the following will be

\begin{equation*}
rdx^3 = (x-a)^3d^3 P+9(x-a)^2 ddp dx+18(x-a)dPdx^2+6Pdx^3=6P dx^3
\end{equation*}
having put $x=a$. Hence, because $p$ and $q$ vanish in the case $x=a$, it will be

\begin{equation*}
dy= \frac{1}{6}rdx^3 =Pdx^3.
\end{equation*}
In like manner, if one puts

\begin{equation*}
z=(x-a)^3 P+C
\end{equation*}
and $y$ was any function of $z$ such that  $dy=Zdz$, because of

\begin{equation*}
dz=(x-a)^3 dP+3(x-a)^2 Pdx,
\end{equation*}
it will also be $p=0$ and $q=0$ and it will be $rdx^3 =6PZdx^3 $; hence, in the case $x=a$ it will be

\begin{equation*}
dy=PZdx^3.
\end{equation*}
Therefore, this function $y$, even though in the case $x=a$  $p=0$, will nevertheless have neither a maximum nor minimum value.

\paragraph*{§343}

These differentials can be found more easily  from the nature of the differentials. For, since the differential of $y$ results, if $y$ is subtracted from the closest following state which results, if  one writes $x+dx$ instead of $x$, let us in the first case, in which it was

\begin{equation*}
y=(x-a)^2P+C,
\end{equation*}
write $x+dx$ instead of $x$ and it will be

\begin{equation*}
y^{\text{I}}= (x-a+dx)^2P^{\text{I}}+C,
\end{equation*}
whence it will be

\begin{equation*}
dy =(x-a+dx)^2P^{\text{I}}-(x-a)P.
\end{equation*}
Therefore, in the case  $x=a$ it will be $dy =P^{\text{I}}dx^2$, and because $P^{\text{I}}$ and $P$ have the ratio of $1$, it will be

\begin{equation*}
dy = Pdx^2.
\end{equation*}
In like manner, if it was 

\begin{equation*}
z=(x-a)^2P+C,
\end{equation*}
it will be $dz =Pdx^2$; hence, if $y$ is any function of $z$ such that $dy=Zdz$, it will be

\begin{equation*}
dy=PZdx^2
\end{equation*}
in the case in which one puts $x=a$. \\[2mm]
Further, if 

\begin{equation*}
z=(x-a)^3 P+C,
\end{equation*}
it will be $z^{\text{I}}=(x-a+dx)^3P^{\text{I}}+C$ and therefore in the case $x=a$ it will be

\begin{equation*}
z^{\text{I}}-z =dz=Pdx^3.
\end{equation*}
Hence, if $y$ was any function of $z$ and $dy=Zdz$, the differential in the case $x=a$ will also be

\begin{equation*}
dy= PZdx^3,
\end{equation*}
if in the functions $P$ and $Z$  one substitutes $a$ for $x$ everywhere. Since in this case  $z=C$ and $Z$ is a function of $z$, $Z$ will become a constant quantity, such a function of $C$, of course, as it was one of $z$ before.

\paragraph*{§344}

Therefore, if it was in general

\begin{equation*}
y=(x-a)^nP+C,
\end{equation*}
since 

\begin{equation*}
y^{\text{I}}=(x-a+dx)^n P^{\text{I}}+C,
\end{equation*}
in the case $x=a$ it will be

\begin{equation*}
dx=Pdx^n;
\end{equation*}
whence, if it was $n>1$, this differential will vanish compared to the higher first differentials which are homogeneous to $dx$. Therefore, from the preceding it is obvious that the function $y$ has a maximum or minimum value in the case $x=a$, if $n$ was an even number; for, then, if for $x=a$ $P$ becomes a positive quantity, $y$ will have a minimum; but if $P$ was a negative quantity, $y$ will have a maximum. And this way the nature of maxima and minima is found a lot more easily than by using the method explained above, since it is not necessary to  consider higher differentials. If 

\begin{equation*}
z=(x-a)^n P+C
\end{equation*}
and $y$ was any function of $z$ that $dy=Zdz$, in the case $x=a$ the differential will be

\begin{equation*}
dy= PZdx^n.
\end{equation*}
But it is to be noted that here $n$ it to be taken positively or greater than $0$; for, if $n$ was a negative number, then  $(x-a)^n$ would not vanish for $x=0$, as we assumed, but would even become infinitely large.

\paragraph*{§345}

Now, we saw that this way the differential can be found a lot more conveniently than by means of the series we used to express the complete differential before; for, if $n$ was an integer number, so many terms of that series have to be considered as $n$ contains unities. But if $n$ was a fractional number,  this series will not even ever exhibit the true differential. For, let us put that 

\begin{equation*}
y=(x-a)^{\frac{3}{2}}+a\sqrt{a};
\end{equation*} 
if we consider the series

\begin{equation*}
dy= pdx+ \frac{1}{2}qdx^2+\frac{1}{6}rdx^3+\frac{1}{24}dx^4+\text{etc.},
\end{equation*}
it will be

\begin{equation*}
p=\frac{3}{2}\sqrt{x-a}, \quad q=\frac{3}{4\sqrt{x-a}}, \quad r=\frac{-3}{8(x-a)\sqrt{x-a}},
\end{equation*}
\begin{equation*}
s= \frac{9}{16(x-a)^2\sqrt{x-a}} \quad \text{etc.}
\end{equation*}
Hence, if one puts $x=a$, it will be $p=0$, but all following terms $q$, $r$, $s$ etc. will become infinite; hence, the value of the differential $dy$ cannot be defined in this case at all. But the method deduced from the nature of differentials leaves no doubt. For, since $y=(x-a)^{\frac{3}{2}}+a \sqrt{a}$, having written $x+dx$ instead of $x$, it will be $y^{\text{I}}=(x-a+dx)^{\frac{3}{2}}+a\sqrt{a}$ and, if one puts $x=a$, it will be $dy=dx\sqrt{dx}$. Therefore, this differential vanishes with respect to $dx$; but the second differentials  homogeneous to $dx^2$ will vanish with respect to the latter.

\paragraph*{§346}

Let us expand these cases in which the exponent $n$ is a fractional number a little more accurately and let be

\begin{equation*}
y= P\sqrt{x-a}+C;
\end{equation*}
because of $y^{\text{I}}=P^{\text{I}}\sqrt{x-a+dx}+C$, it will be

\begin{equation*}
dy =P\sqrt{dx}
\end{equation*}
in the case $x=a$; therefore, this differential will have an infinite ratio to $dx$ and to the differentials homogeneous to $dx$. Hence it is plain, what is to be said about  the nature of maxima and minima in this case. For, because, having written $a+dx$ instead of $x$, $y$ goes over into

\begin{equation*}
C+P\sqrt{dx},
\end{equation*}
because of the ambiguous $\sqrt{dx}$, the function $y$ will obtain two values; the one greater than $C$, which it receives for $x=a$, the other smaller; hence in the case $x=a$ it will have neither a maximum nor a minimum value. Furthermore, if $dx$ is taken negatively, the value of $y$ will even become imaginary. The same is to said, if $z=P\sqrt{x-a}+C$ and $y$ is any function of $z$ that it is $dy=Zdz$; for, then it will be $dy= PZ\sqrt{dx}$ in the case $x=a$.

\paragraph*{§347}

If this function was propounded

\begin{equation*}
y=(x-a)^{\frac{m}{n}}P+C,
\end{equation*}
whose differential is in question in the case $x=a$, it will, as we concluded from the preceding, be

\begin{equation*}
dy= Pdx^{\frac{m}{n}}.
\end{equation*}
Therefore, if it was $m > n$, this differential will vanish with respect to $dx$; but if $m<n$, the ratio $\frac{dy}{dx}$ will be infinitely large. Furthermore, if $n$ is an even number, the differential will have two values, the one positive, the other negative; and so the function $y$ which in the case $x=a$ becomes $=C$, if one puts $x=a+dx$,  will have two values, the one greater than $C$, the other smaller; but if one would put $x=a-dx$, then $y$ would even become imaginary; hence, in this case $y$ would become neither a maximum nor a minimum. Now, let us put that the denominator $n$ is an odd number; the numerator $m$ will either be even or odd. Let $m$ be an even number at first; since $dy$ retains the same value, no matter whether $dx$ is taken positively or negatively, it is perspicuous that the function $y$ in the case $x=a$ has either a maximum or minimum value depending on whether in this case $P$ was a negative or positive quantity. But if both numbers $m$ and $n$ were odd, the differential $dy$ will go over into its negative, if $dx$ was negative; and therefore,  the function $y$ will have neither a maximum nor a minimum in this case, if one puts $x=a$.

\paragraph*{§348}

If the function $y$ consists  of several terms of this kind, each of which  is divisible by $x-a$, such that 

\begin{equation*}
y=(x-a)^mP+(x-a)^n Q+C,
\end{equation*} 
in the case $x=a$ its differential will be

\begin{equation*}
dy=Pdx^m+Qdx^n;
\end{equation*}
in this expression, if it was $n>m$, the second term vanishes in comparison to the first such that only $dy= Pdx^m$ results. But if $n$ was a fraction with an even denominator, then, even though $Qdx^n$ vanishes with respect to $Pdx^m$, it can nevertheless not be completely neglected. For, from this it is clear, if $dx$ is taken negatively, that the value of $dy$ becomes imaginary, whereas from the first term $Pdx^m$ only this is not obvious. Therefore, since, if $n$ is a fraction with an even denominator, $dx$ cannot be taken negatively, but if it is taken positively, the term $Qdx^n$ yields two values, the function $y=(x-a)^mP+(x-a)^n Q+C$ which in the case $x=a$ becomes $=C$, if $x=a+dx$, will be

\begin{equation*}
y=C+Pdx^m\pm Qdx^n;
\end{equation*}
since both of these values are either greater or smaller than $C$, depending on whether $P$ was a positive or a negative quantity,  in the case $x=a$ the function $y$ will have either a minimum or a maximum value of the second kind [§ 278].

\paragraph*{§349}

Therefore, in such cases the true differentials of functions cannot be found by means of the usual rules for differentiation; those are only applicable, if the differential of the function is homogeneous to $dx$. But if in a singular case the differential of the function is expressed in terms of its power $dx^n$, then the rule yields $0$ for this differential, if $n$ was a number greater than $1$; but  on the other hand it exhibits an infinite times larger differential, if $n$ is an exponent smaller than $1$. So if the differential of $y=\sqrt{a-x}$ is in question in the case $x=a$, since $dy=-\frac{dx}{\sqrt{a-x}}$, having put $x=a$, $dy=-\frac{dx}{0}$ results. And if we wanted to derive the following differentials, because of the denominators $=0$, all of them grow to infinity in the same way such that nothing can be concluded from this. But we saw that in this case  $dy=\sqrt{-dx}$ and hence imaginary. But if  one writes $x-dx$ instead of $x$, it will be $dy=\sqrt{dx}$ and hence it will be infinitely larger than $dx$ such that $dx$ vanishes with respect to $dy$. Hence, the usual rule even in this case does not cause any errors, since it exhibits the infinite value of $dy$.

\paragraph*{§350}

Therefore, one must not apply the usual rule, if in the series

\begin{equation*}
pdx+\frac{1}{2}qdx^2+\frac{1}{6}rdx^3+\text{etc.},
\end{equation*}
expressing the complete differential of the function $y$, the first term $p$ either becomes $=0$ or infinite, but in this case the differential must be derived from the first principles. Therefore, if the differential corresponding to a given value of $x$ of the function $y$, for which the letter $p$ becomes either infinitely small or infinitely large,   is in question, one has to go back to the first principles of differentiation. In all remaining cases in which neither $p=0$ nor $p=\infty$ the usual rule will yield the true values of the the differential. Nevertheless, the case mentioned before (§ 348) is not to be neglected, if the function $y$ contains a term of the form $(x-a)^nQ$ while $n$ is fraction with an even denominator; for, even though one has lower differentials than $Qdx^n$, with respect to which this one vanishes, nevertheless, since $Qdx^n$, if $dx$ is negative, becomes imaginary, this term $Qdx^n$ transforms all remaining ones, with respect to which it vanishes, also into imaginary ones;  this circumstance is especially to be considered in the case of curves. Therefore, I will explain some particular cases in which the true differential is not indicated by the common rule in the following examples.

\subsection*{Example 1}

\textit{Let the differential of the function $y=a+x-\sqrt{xx+ax-x\sqrt{2ax-xx}}$ be in question in the case $x=a$.}\\[2mm]
It is plain that the differential of this function $x=a$ is not found by means of the usual rule; for, 

\begin{equation*}
dy= dx + \frac{-xdx-\frac{1}{2}adx+\frac{1}{2}dx\sqrt{2ax-xx}+(axdx-xxdx):\sqrt{2ax-xx}}{\sqrt{xx+ax-x\sqrt{2ax-xx}}};
\end{equation*}
for, having put $x=a$, it will be $dy=dx-\frac{adx}{a}=0$. Therefore, let us start from the first principles of differentiation and at first, having written $x+dx$ instead of $x$, it will be

\begin{equation*}
y^{\text{I}}=a+x+dx
\end{equation*}
\begin{equation*}
-\sqrt{xx+2xdx+dx^2+ax+adx-(x+dx)\sqrt{2ax-xx+2adx-2xdx-dx^2}}
\end{equation*}
But, having put $x=a$, it will be

\begin{equation*}
y^{\text{I}}=2a+dx- \sqrt{2aa+3adx+dx^2-(a+dx)\sqrt{aa-dx^2}}.
\end{equation*}
Now, because $\sqrt{aa-dx^2}=a -\frac{dx^2}{2a}$ (for, the following terms can be neglected, since not all which are infinite times larger will be cancelled, as it will be seen soon), it will be

\begin{equation*}
y^{\text{I}}=2a+dx-\sqrt{aa+2adx+\frac{3}{2}dx^2}
\end{equation*}
and further, by extracting the root, it will be

\begin{equation*}
y^{\text{I}}=2a+dx- \left(a+dx+\frac{dx^2}{4a}\right)=a- \frac{dx^2}{4a}.
\end{equation*}
But in the case $x=a$ it will be $y=a$; hence, because  $y^{\text{I}}=y+dy$, one will obtain

\begin{equation*}
dy=-\frac{dx^2}{4a};
\end{equation*}
from this it is at the same time seen that the propounded function $y$ becomes imaginary, if one puts $x=a$.

\subsection*{Example 2}

\textit{To find the differential of this function $y=2ax-xx+a\sqrt{aa-xx}$ in the case in which one puts $x=a$.}\\[2mm]
Having differentiated it in the usual way, 

\begin{equation*}
dy= 2adx -2xdx-\frac{axdx}{\sqrt{aa-xx}},
\end{equation*}
which  for $x=a$ goes over into infinity and hence is not indicated this way. But the differentials of the following orders in a like manner will become infinite such that not even from the series $pdx+\frac{1}{2}qdx^2+\frac{1}{6}rdx^2+\text{etc.}$ the true value of the differential can be found. Therefore, let us write $x+dx$ instead of $x$ and we will have

\begin{equation*}
y^{\text{I}}=2ax-xx+2adx-2xdx-dx^2+a\sqrt{aa-xx-2xdx-dx^2}
\end{equation*}
and, having put $x=a$, it will be

\begin{equation*}
y^{\text{I}}= aa-dx^2+a\sqrt{-2adx-dx^2}.
\end{equation*}
But in the same case $y=aa$; therefore, it will be $dy=-dx^2+a\sqrt{-2adx}$, and since $dx^2$ vanishes with respect to $\sqrt{-2adx}$, it will be

\begin{equation*}
dy=a \sqrt{-2adx}.
\end{equation*}
Hence, if the differential $dx$ is taken positively, $dy$ will be  imaginary; but if one writes $x-dx$ for $x$, it will be

\begin{equation*}
dy=\sqrt{2adx};
\end{equation*}
because its value is a double value, the one positive, the other negative, the function $y$  will have neither a maximum nor minimum value in the case $x=a$.

\subsection*{Example 3}

\textit{To find the differential of the function $y=3aax-3axx+x^3+(a-x)\sqrt[3]{a^3-x^3}$ in the case  $x=a$.}\\[2mm]
Since this function is transformed into this form

\begin{equation*}
y=a^3-(a-x)^3+(a-x)^{\frac{7}{3}}\sqrt[3]{aa+ax+xx},
\end{equation*}
having put $x=a+dx$,

\begin{equation*}
y^{\text{I}}= a^3+dx^3-dx^{\frac{7}{3}}\sqrt[3]{3aa},
\end{equation*}
and in the same case  $y=a^3$. Therefore, it will be $dy=dx^3-dx^{\frac{7}{3}}\sqrt[3]{3aa}$, and because $dx^3$ vanishes with respect to $dx^{\frac{7}{3}}$, it will be

\begin{equation*}
dy= -dx^{\frac{7}{3}}\sqrt[3]{3aa};
\end{equation*}
therefore, in the case $x=a$ the function $y$ has neither a maximum nor a minimum value.

\subsection*{Example 4}

\textit{To find the differential of the function $y+\sqrt{x}+\sqrt[4]{x^3}=(1+\sqrt[4]{x})\sqrt{x}$ in the case $x=0.$} \\[2mm]
Since the case $x=0$ is propounded and in this case $y=0$,   write only $dx$ instead of $x$ and one will have

\begin{equation*}
dy= dx^{\frac{1}{2}}+dx^{\frac{3}{2}} \quad \text{or} \quad dy=(1+\sqrt[4]{dx})\sqrt{dx};
\end{equation*}
hence, at first it is plain that $dx$ can not be taken negatively. But then, even though $\sqrt{dx}$ has a double value, the one positive, the other negative, nevertheless in this case, since its root $\sqrt[4]{dx}$ occurs, only the positive value can be taken. But $\sqrt[4]{dx}$ has both meanings and it will be

\begin{equation*}
dy= \sqrt{dx}\pm \sqrt[4]{dx^3} \quad \text{and} \quad y^{\text{I}}=0+\sqrt{dx}\pm \sqrt[4]{dx^3},
\end{equation*}
because of $y=0$. Because both values of $y^{\text{I}}$ are greater than $y$, it follows  that in the case $x=0$ $y$ has a minimum value. But that the function $y= \sqrt{x}+\sqrt[4]{x^3}$ does not contain this one $-\sqrt{x}+\sqrt[4]{x^3}$,  will become plain by reducing both to rational expressions. For, the first reduced this form $y-\sqrt{x}=\sqrt[4]{x}$ and squared gives $y^2-2y\sqrt{x}+x=x \sqrt{x}$ or $y^2+x=(x+2y)\sqrt{x}$ which squared again yields

\begin{equation*}
y^4-2yyx-4xxy+xx-x^3=0
\end{equation*}
The other $y+\sqrt{x}=\sqrt[4]{x^3}$ will give $y^2+x=(x-2y)\sqrt{x}$ and further

\begin{equation*}
y^4 -2yyx+4xxy+xx-x^3=0,
\end{equation*}
which is different from the first. But on the other hand, the term $\sqrt[4]{x^3}$ retains the ambiguity of the sign. Therefore, this circumstance is to be considered in detail, that, even though in general the roots of even powers include both signs  $+$ and $-$, nevertheless this ambiguity does not occur, if in the same expression  higher roots of even powers  of the same roots occur; these would be imaginary, if the first roots would be taken negatively. And from this source maxima and minima of the second kind follow, whenever such might not seem to occur.

\subsection*{Example 5}

\textit{To find the differential of the function}

\begin{equation*}
y=a+\sqrt{x-f}+(x-f)\sqrt[4]{x-f}+(x-f)^2\sqrt[8]{x-f}
\end{equation*}
\textit{in the case $x=f$.} \\[2mm]
Let us put $x-f=t$, and since $y=a+\sqrt{t}+t\sqrt[4]{t}+tt\sqrt[8]{t}$, the differential of this expression is in question in the case $t=0$ in which  $y=a$. Therefore, having written $t+dt$ or $0+dt$ instead of $t$, it will be

\begin{equation*}
y^{\text{I}}=y+dy=a+\sqrt{dt}+dt\sqrt[4]{dt}+dt^2\sqrt[8]{dt}
\end{equation*}
and hence one will have

\begin{equation*}
dy= \sqrt{dt}+dt\sqrt[4]{dt}+dt^2\sqrt[8]{dt}.
\end{equation*}
Here, at first it is plain that the differential can not be taken negatively, for, otherwise $dy$ would become imaginary. But  not only $\sqrt{dt}$, but  even $\sqrt[4]{dt}$ can not be taken negatively; for, $\sqrt[8]{dt}$ would become imaginary; hence the differential $dy$ has only the double value

\begin{equation*}
dy= \sqrt{dt}+dt \sqrt[4]{dt}\pm dt^2 \sqrt[8]{dt};
\end{equation*}
because both values are greater than  zero, it follows that the function $y$ has a minimum value of the second kind for $t=0$ or $x=f$. Although in these cases the terms $dt \sqrt[4]{dt}$ and $dt^2\sqrt[8]{dt}$ vanish with respect to the first $\sqrt{dt}$, it is nevertheless to be taken into account, if the multiplicity of the values is considered, that the imaginary quantities are avoided.

\subsection*{Example 6}

\textit{To find the differential of the function $y=ax+bxx+(x-f)^n+(x-f)^{m+\frac{1}{2}}$ in the case $x=f$.}\\[2mm]
If one puts $x=f$, it will be $y=af+bff$, and if  one writes $x+dx$ or $f+dx$ instead of $x$, the closest value will result as

\begin{equation*}
y^{\text{I}}=af+bff+adx+2bfdx+bdx^2+dx^n+dx^{m+\frac{1}{2}n},
\end{equation*}  
such that 

\begin{equation*}
dy=adx+2bfdx+bdx^2+dx^n+dx^m\sqrt{dx^n}.
\end{equation*}
Therefore, if $n$ is not an even number, the differential $dx$ cannot be taken negatively. But the last term $dx^m\sqrt{dx^n}$ has an ambiguous sign; therefore, $y^{\text{I}}$ will have two values, both greater than the one of $y$, if $a+2bf$ was a positive quantity and the exponents  $n$ and $m+\frac{1}{2}n$ were greater than $1$. Therefore, the value of the function $y$ in the case $x=f$ will be a minimum and this happens, whether $n$ is an integer number or a fraction, as long as only the numerator was not also an even number in that case. 

\paragraph*{§351}

But this method to deduce differentials from the first principles is especially useful  in the case of transcendental functions, since in certain cases the differential found in the usual way either vanishes or seems to grow to infinity. But here species of the infinite and the infinitely small appear which are not found in algebraic functions. For, because, if $i$ denotes an infinite number, $\log i$ is also infinite, but nevertheless  has an infinitely small ratio  to the number $i$ and even to any power $i^n$, no matter how small the exponent $n$ is set, the fraction $\frac{\log i}{i^n}$ will be infinitely small and  can not be finite, unless the exponent $n$ becomes infinitely small. Therefore, $\log i$ will be homogeneous to $i^n$, if the exponent $n$ was infinitely small. Now let us put $i=\frac{1}{\omega}$ while $\omega$ is an infinitely small quantity; $-\log \omega$ will be homogeneous to $\frac{1}{\omega^n}$, if the exponent $n$ is infinitely small, and hence $-\frac{1}{\log \omega}$ will be homogeneous to $\omega^n$; and hence, $-\frac{1}{\log dx}$ will be infinitely small compared to $dx^n$, while $n$ is an infinitely small fraction. So, if it was $y=-\frac{1}{\log x}$, the differential of $y$ in the case $x=0$ will be $=-\frac{1}{\log dx}=dx^n$ and hence $dy$ will have an infinite ratio to $dx$ and to  any power of $dx$; and with respect to $-\frac{1}{\log dx}$ completely all powers of $dx$ vanish, no matter how small their exponents were.

\paragraph*{§352}

Further, we also saw, if $a$ was a number greater than $1$ and $i$ was infinite, that  $a^i$ will be infinite of such a high degree that with respect to it not only $i$ but also any power of $i$ vanishes; and $i^n$ does not become homogeneous to $a^i$ until the exponent $n$ was increased to infinity. Now, let $i=\frac{1}{\omega}$ such that $\omega$ denotes the infinitely small; $a^{\frac{1}{\omega}}$ will be homogeneous to $\frac{1}{\omega^n}$ while $n$ denotes an infinitely large number and hence $a^{\frac{-1}{\omega}}$ or $\frac{1}{a^{1:\omega}}$ will be infinitely large compared to $\omega^n$. Therefore, $\frac{1}{a^{1:dx}}$ will be infinitely small, but vanishes with respect to all powers of $dx$, because it is homogeneous to the power $dx^n$, while $n$ is an infinitely large number. Hence, if the differential of $y=\frac{1}{a^{1:x}}$ is in question in the case $x=0$, since  $y=0$, it will be $dy=\frac{1}{a^{1:dx}}$ and hence is infinite times smaller than each power of $dx$.

\paragraph*{§353}

But if $a$ was a number smaller than $1$, then, because $\frac{1}{a}$ becomes larger than $1$, the question is reduced to the preceding case. If one has the expression $a^{\frac{1}{\omega}}$,  by putting $a=1:b$ it will be transformed into $b^{-\frac{1}{\omega}}$ or $\frac{1}{b^{1:\omega}}$ which, because of $b>1$, will be homogeneous to $\omega^n$ while $n$ denotes any infinitely large number. Therefore, having mentioned these things in advance, we will be able to resolve the following examples.

\subsection*{Example 1}

\textit{To find the differential of the function $y=xx-\frac{1}{\log x}$ in the case $x=0$.}\\[2mm]
Since for $x=0$ also $y=0$, if we write $x+dx$ or $0+dx$ instead of $x$, it will be

\begin{equation*}
y^{\text{I}}=dy=dx^2-\frac{1}{\log dx}.
\end{equation*}
But because $-\frac{1}{\log dx}$ is homogeneous to $dx^n$ while $n$ denotes an infinitely small number, with respect to it $dx^2$ will vanish and it will be

\begin{equation*}
dy=-\frac{1}{\log dx}=dx^n.
\end{equation*}
But because the logarithms of negative numbers are imaginary, $dx$ cannot be taken negatively and therefore in the case $x=0$ the function $y$ will have a minimum value, but a minimum extending neither to the first nor the second kind. It certainly does not extend to the first kind, since $y$ does not have any preceding very close values, but  is only smaller than the following values, if $x$  is greater than zero. But it also does not extend to second kind, since the following values, it is compared to, are not double values; and therefore a third kind of maxima and minima results that only occurs in logarithmic and transcendental functions, but never occurs in algebraic functions; this will be treated in the following book on curves in more detail.

\subsection*{Example 2}

\textit{To find the differential of the function $y=(a-x)^n-x^n(\log a- \log x)^n$ in the case  $x=a$.}\\[2mm]
This differential, if $n$ is not an integer number, can be found from the general formula

\begin{equation*}
dy= pdx+\frac{1}{2}qdx^2+\frac{1}{6}rdx^3+\text{etc.};
\end{equation*}
for, it will be

\begin{equation*}
pdx=- n(a-x)^{n-1}dx-nx^{n-1}dx(\log a- \log x)^n+nx^{n-1}(\log a- \log x)^{n-1}dx,
\end{equation*}
which value vanishes for $x=a$; for, even if $n=1$, it will be

\begin{equation*}
pdx=-dx+dx=0.
\end{equation*}
Therefore, if we proceed further, it will be

\begin{equation*}
\frac{1}{2}qdx^2 =\frac{n(n-1)}{1 \cdot 2}(a-x)^{n-2}dx^2- \frac{n(n-1)}{1 \cdot 2}x^{n-1}dx^2(\log a- \log x)^n+\frac{n^2}{2}x^{n-2}dx^2(\log a- \log x)^{n-1}
\end{equation*}
\begin{equation*}
+\frac{n(n-1)}{1 \cdot 2}x^{n-2}dx^2(\log a -\log x)^{n-1}-\frac{n(n-1)}{1 \cdot 2}dx^2(\log a- \log x)^{n-2}.
\end{equation*}
Hence, if it was $n=1$, it will be $\frac{1}{2}qdx^2= \frac{dx^2}{2a}$ for $x=a$. In like manner, if $n=2$, one would have to proceed up to the term $\frac{1}{6}rdx^3$ and so fourth. Therefore, one will more conveniently use the principles of differentiation, and because for $x=a$  $y=0$, if we write $x+dx$ or $a+dx$ instead of $x$, it will be

\begin{equation*}
y^{\text{I}}=(-dx)^n-(a+dx)^n(\log a- \log (a+dx))^n=y+dy=dy,
\end{equation*}
because of $y=0$. But 

\begin{equation*}
\log (a+dx)=\log a+\frac{dx}{a}-\frac{dx^2}{2a^2}+\frac{dx^3}{3a^3}-\text{etc.},
\end{equation*}
whence 

\begin{equation*}
dy=(-dx)^n-\left(a^n+na^{n-1}dx+\frac{n(n-1)}{1 \cdot 2}a^{n-2}dx^2+\text{etc.}\right)\left(-\frac{dx}{a}+\frac{dx^2}{2a^2}-\frac{dx^3}{3a^3}+\text{etc.}\right)^n
\end{equation*}
\begin{equation*}
=\frac{n}{2a}(-dx)^{n+1}.
\end{equation*}
Therefore, in the case $x=a$ the differential $dy$ in question of the propounded formula will be as follows:

\begin{equation*}
\renewcommand{\arraystretch}{2,0}
\setlength{\arraycolsep}{0.0mm}
\begin{array}{l|lllllllllllll}
\text{if} \quad n=1 \quad & \quad dy &~= ~&\dfrac{dx^2}{2a}, \quad & \text{as we found before} \\ 
\text{if} \quad n=2 \quad & \quad dy &~=- ~ &\dfrac{2dx^3}{2a} \\ 
\text{if} \quad n=3 \quad & \quad dy &~= ~ &\dfrac{3dx^4}{2a} \\ 
\text{if} \quad n=4 \quad & \quad dy &~=- ~ &\dfrac{4dx^5}{2a} \\
\text{etc.}               & \quad \text{etc.} 
\end{array}
\end{equation*} 
Therefore, if $n$ was an odd number, in the case $x=a$ the function  has a minimum value, but if $n$ is an even number, neither a maximum nor a minimum value; the same holds, if $n$ was a fraction with an odd denominator. But if $n$ was a fraction with an even denominator, $dx$ has to be taken negatively, so  that we do not get to imaginary quantities; and because of the ambiguous meaning, the function will also have neither  a maximum value nor a minimum value.

\subsection*{Example 3}

\textit{To find the differential of the function $y=x^x$ in the case $x=\frac{1}{e}$, where $e$ denotes the number whose hyperbolic logarithm is $=1$.}\\[2mm]
Since in general  $dy=x^x dx (\log x+1)$, this differential vanishes in the case $x=\frac{1}{e}$ or $\log x=-1$. Therefore, compare this differential to the general form $pdx+\frac{1}{2}qdx^2+\text{etc.}$; it will be

\begin{equation*}
p=x^x (\log x+1) \quad \text{and} \quad q=x^x(\log x+1)^2+x^{x-1}
\end{equation*}
and, having put $\log x =-1$ or $x=\frac{1}{e}$, it will be

\begin{equation*}
q= \left(\frac{1}{e}\right)^{\frac{1-e}{e}}=e^{\frac{e-1}{e}}.
\end{equation*}
Therefore, the differential in question will be

\begin{equation*}
dy= \frac{1}{2}e^{(e-1):e}dx^2
\end{equation*}
and therefore the function $y=x^x$ will have a minimum value in the case $x=\frac{1}{e}$.

\subsection*{Example 4}

\textit{To find the differential of this function $y=x^n +e^{-1:x}$ in the case in  $x=0$.}\\[2mm]
Since for $x=0$ also $y=0$, if one puts $x=0+dx$, it will be

\begin{equation*}
y^{\text{I}}=dy=dx^n+\frac{1}{e^{1:dx}}.
\end{equation*}
But we saw that $\frac{1}{e^{1:x}}$ is homogeneous to the infinite power of $dx$ or to $dx^{\infty}$ and  will hence vanish with respect to $dx^n$ such that 

\begin{equation*}
dy=dx^n.
\end{equation*}

\paragraph*{§354}

 What happens in the first differentials in certain cases, i.e. that they are not found by means of the usual rules of differentiation, also happens in differentials of  the second and third and higher order in the cases in which in the complete differential form
 
 \begin{equation*}
d.y = pdx+\frac{1}{2}qdx^2+\frac{1}{6}rdx^3+\frac{1}{24}sdx^4+\text{etc.} 
 \end{equation*}
 some of the quantities $q$, $r$, $s$ etc. either vanish or become infinite. Because  [§ 339]
 
\begin{equation*}
 dd.y = qdx^2+rdx^3+\frac{7}{12}sdx^4+\text{etc.},
\end{equation*}
if in which case  $q=0$,  it will be $ddy =rdx^3$; but if in the same case also $r$ vanishes, from this series the second differential can not be found at all, but one will have to go back to the principles of differentials; writing $x+dx$ instead of $x$ find the value $y^{\text{I}}$ and by writing $x+2dx$ instead of $x$ find the value of $y^{\text{II}}$, having done which the true value of the second differential will be

\begin{equation*}
ddy= dy^{\text{I}}-dy=y^{\text{II}}-2 y^{\text{I}}+y.
\end{equation*}
If in like matter if the the third differential is in question, additionally,   write $x+3dx$ instead of $x$ in $y$ and, having found the value $y^{\text{III}}$, it will be

\begin{equation*}
d^3y =y^{\text{III}}-2 y^{\text{II}}+3y^{\text{I}}-y
\end{equation*}
and so fourth. We will illustrate these cases in the following examples.

\subsection*{Example 1}

\textit{To find the second differential of the function $y=\frac{aa-xx}{aa+xx}$ in the case  $x=\frac{a}{\sqrt{3}}$}.\\[2mm]
Investigating the complete differential of $y$ from the form

\begin{equation*}
dy= pdx+ \frac{1}{2}qdx^2+\frac{1}{6}rdx^3+\frac{1}{24}sdx^4+\text{etc.}
\end{equation*}
the following values will result for $p$, $q$, $r$, $s$ etc. 

\begin{equation*}
p=-\frac{4aax}{(aa+xx)^2}, \quad q= \frac{-4a^4+12aaxx}{(aa+xx)^3} \quad \text{and} \quad r=\frac{48a^4x-48aax^3}{(aa+xx)^4}.
\end{equation*}
Since now 

\begin{equation*}
ddy= rdx^3+\frac{7}{12}sdx^4+\text{etc.},
\end{equation*}
because of $q=0$, in the case $x=\frac{a}{\sqrt{3}}$ and in the same case $r=\frac{27 \sqrt{3}}{8a^3}$, the second differential in question will be

\begin{equation*}
ddy= \frac{27dx^2\sqrt{3}}{8a^3}.
\end{equation*}

\subsection*{Example 2}

\textit{To find the third differential of the function $y=\frac{aa-xx}{aa+xx}$ in the case $x=a$.}\\[2mm]
As before by finding the complete differential

\begin{equation*}
dy= \frac{1}{2}qdx^2+\frac{1}{6}rdx^3+\frac{1}{24}sdx^4+\text{etc.},
\end{equation*}
since the third differential is $d^3y = rdx^3+\frac{3}{2}sdx^4+\text{etc.}$, because of

\begin{equation*}
r=\frac{48a^4 x-48aax^3}{(aa+xx)^4},
\end{equation*}
it will be $x=a$ in the case $r=0$; hence, one has to proceed to the value $s$ which will be

\begin{equation*}
s=\frac{48a^4-144aaxx}{(aa+xx)^4}-\frac{8x(48a^4x-48aax^3)}{(aa+xx)^5};
\end{equation*}
therefore, having put $x=a$, it will be $s=-\frac{96a^4}{2^4 a^8}=-\frac{6}{a^4}$; therefore,  in this case it will be

\begin{equation*}
d^3y =-\frac{9dx^4}{a^4}.
\end{equation*}

\subsection*{Example 3}

\textit{To find the differential of arbitrary order of the function $y=ax^m+bx^n$ in the case $x=0$.}\\[2mm]
By successively writing $x+dx$, $x+2dx$, $x+3dx$ etc. instead of $x$ the following values of the function $y$ will be

\begin{alignat*}{9}
&y^{\text{I}}&&=a(x+dx)^m &&+b(x+dx)^n, \\
&y^{\text{II}}&&=a(x+2dx)^m &&+b(x+2dx)^n, \\
&y^{\text{III}}&&=a(x+3dx)^m &&+b(x+3dx)^n, \\
& &&\text{etc.}
\end{alignat*}
Therefore, having put $x=0$, it will be $y=0$ and its differentials will be

\begin{alignat*}{9}
 dy &=adx^m+bdx^n, \\
 ddy &=(2^m-1)adx^m+(2^n-2)bdx^n, \\
 d^3 y&=(3^m- 3 \cdot 2^m+3)adx^m+(3^n-3 \cdot 2^n+3)bdx^n, \\
 d^4 y &=(4^m-4 \cdot 3^m+6 \cdot 2^m -4)adx^m +(4^n-4 \cdot 3^n+6 \cdot 2^n-4 )bdx^n \\
 &\text{etc.} 
\end{alignat*}
Therefore, if the exponent $n$ was greater than $m$, the second terms in these expressions vanish with respect to the first. Nevertheless, it is to be considered, if $n$ was a fractional number, that the cases in which these differentials become either imaginary or ambiguous can be distinguished. It will be convenient reserve the further discussion of these cases for the doctrine of curves.


\title{On the values of functions which  seem to be undetermined in certain cases
\footnote{
Original title: "'De Valoribus Functionum qui certis casibus videntur indeterminanti"', first published as part of the book \textit{Institutiones calculi differentialis cum eius usu in analysi finitorum ac doctrina serierum}, 1755, reprinted in \textit{Opera Omnia}: Series 1, Volume 10, pp. 564 - 587, Eneström-Number E212, translated by: Alexander Aycock for the Euler-Kreis Mainz}}
\author{Leonhard Euler}
\date{ }
\maketitle

\paragraph*{§355}

If a function  $y$ of $x$ was a fraction $\frac{P}{Q}$ whose numerator and denominator   vanish  for the same value of $x$,  in that case the fraction $\frac{P}{Q}$  expressing the value of the function $y$ will become $=\frac{0}{0}$; since this expression can become equal to any either finite or infinite or infinitely small quantity, from it the value of $y$ in this case  can not be concluded at all and hence that value seems to be undetermined. Nevertheless, it is easily seen, since, except in this case, the function $y$ always has a determined value, whatever is substituted for $x$, that also in this case the value of $y$ can not be undetermined. This will become obvious from the example $y=\frac{aa-xx}{a-x}$; having put $x=a$, it will be $y=\frac{0}{0}$, of course. But because, having divided the numerator by the denominator, $y=a+x$, it is evident, if one puts $x=a$, that it will be $y=2a$ such that in this case that fraction $\frac{0}{0}$ becomes equal to the quantity $2a$.

\paragraph*{§356}

Therefore, since we showed above that $\frac{0}{0}$ can express any number, in examples of this kind the ratio the numerator has to the denominator must  be investigated. But because in absolute zeros this diversity can not be seen, instead of them infinitely small quantities have to be introduced; even if  in their nature they do not differ from zero, nevertheless from the different functions of them which constitute the numerator and the denominator the value of the fraction becomes obvious immediately. So if one has this fraction $\frac{adx}{bdx}$, even if indeed the numerator and the denominator are $=0$,  nevertheless it is plain that the value of this fraction is determined, namely $=\frac{a}{b}$. But if one has this fraction $\frac{adx^2}{bdx}$, this value will be zero, as the value of $\frac{adx}{bdx^2}$ is infinitely large. Therefore, if   we introduce infinitely small quantities instead of the zeros which often enter the calculation, we will hence be able to find the ratio the zeros have to each other and no further doubt about the meaning of expressions of this kind will remain.

\paragraph*{§357}

To make these things more clear, let us put that so the numerator as the denominator of the fraction $\frac{P}{Q}$ vanish, if one sets $x=a$. But to avoid these zeros which cannot be compared to each other, let us put $x=a+dx$, which reduces to the first $x=a$ because of $dx=0$. But because, if  one writes $x+dx$ instead of $x$, the functions $P$ and $Q$ go over into $P+dP$ and $Q+dQ$, this assumption $x=a+dx$ will be justified, if in these values one sets $x=a$ everywhere, in which case we assumed $P$ and $Q$ to vanish. Hence, if  one writes $a+dx$ instead of $x$, the fraction $\frac{P}{Q}$ will be transformed into this one $\frac{dP}{dQ}$ which therefore expresses the value of the function $y=\frac{P}{Q}$ in the case $x=a$. And this expression can not  be undetermined any longer, if the true differentials of the functions $P$ and $Q$ are taken as we taught in the preceding chapter. For, this way the differentials $dP$ and $dQ$ never go over into absolutely zero, but, if they are not expressed by means of the differential $dx$, will at least be exhibited in terms of its powers. Therefore, if one finds $dP=Rdx^m$ and $dQ=Sdx^n$, in the case $x=a$ the value of the function $y=\frac{P}{Q}$  will be $=\frac{Rdx^m}{Sdx^n}$ which will therefore be finite and $=\frac{R}{S}$, if it was $m=n$; but if  $m>n$, the value of the fraction will indeed be $=0$; but if $m<n$, this value grows  to infinity.

\paragraph*{§358}

Therefore, if a fraction $\frac{P}{Q}$  occurs whose numerator and denominator in a certain case, say $x=a$, vanish at the same time the value of this fraction in this case $x=a$ will be found applying the following rule:\\[2mm]
\textit{Find the differentials of the quantities $P$ and $Q$ in the case $x=a$ and substitute them for $P$ and $Q$, having done which the fraction $\frac{dP}{dQ}$ will exhibit the value of the fraction $\frac{P}{Q}$ in question.}\\[2mm]
If the differentials $dP$ and $dQ$ found by  the usual method become neither  infinite nor vanish in the case $x=a$,  one can use them; but if both become either $=0$ or $=\infty$,   these differentials must be investigated by the method explained in the preceding chapter in the case $x=a$. In most cases the calculation is miraculously contracted, if one puts $x-a=t$ or $x=a+t$ before, that a fraction $\frac{P}{Q}$ results whose numerator and denominator vanish in the case $t=0$; for, then one will have the differentials $dP$ and $dQ$, if one substitutes $dt$ for $t$ everywhere.

\subsection*{Example 1}

\textit{Let the value of the fraction $\frac{b-\sqrt{bb-tt}}{tt}$ be in question in the case $t=0$.}\\[2mm]
Since in this case $t=0$ both the numerator and the denominator vanish,   just write $dt$ instead of $t$ and the value in question will be expressed in terms of the fraction $\frac{b-\sqrt{bb-dt^2}}{dt^2}$. But because $\sqrt{bb-dt^2}=b-\frac{dt^2}{2b}$, this fraction goes over into  $\frac{dt^2}{2bdt^2}=\frac{1}{2b}$. Therefore, in the case $t=0$ the propounded fraction $\frac{b-\sqrt{bb-tt}}{tt}$  has the value $\frac{1}{2b}$.

\subsection*{Example 2}

\textit{Let the value of the fraction $\frac{\sqrt{aa+ax+xx}-\sqrt{aa-ax+xx}}{\sqrt{a+x}-\sqrt{a-x}}$ be in question in the case $x=0$.}\\[2mm]
Here, one can again immediately substitute $dx$ for $x$; because, after the substitution,

\begin{alignat*}{9}
&\sqrt{aa+adx+dx^2}&&=a&&+\frac{1}{2}dx&&+\frac{3dx^2}{8a}, \\[3mm]
&\sqrt{aa-adx+dx^2}&&=a&&-\frac{1}{2}dx&&+\frac{3dx^2}{8a} \\
\end{alignat*}
and

\begin{alignat*}{9}
&\sqrt{a+dx}&&=\sqrt{a}&&+\frac{dx}{2\sqrt{a}},\\[3mm]
&\sqrt{a-dx}&&=\sqrt{a}&&-\frac{dx}{2\sqrt{a}},\\
\end{alignat*}
the numerator will become $=dx$ and the denominator $=\frac{dx}{\sqrt{a}}$, whence the value in question of the propounded fraction will be $=\sqrt{a}$.

\subsection*{Example 3}

\textit{Let the value of this fraction $\frac{x^3-4ax^2+7a^2x-2a^3-2a^2\sqrt{2ax-aa}}{xx-2ax-aa+2a\sqrt{2ax-xx}}$ be in question in the case $x=a$.}\\[2mm]
If  the differentials are taken in usual manner and they are substituted for the numerator and denominator, one will have

\begin{equation*}
\frac{3xx-8ax+7a^2-2a^3:\sqrt{2ax-aa}}{2x-2a+2a(a-x):\sqrt{2ax-xx}},
\end{equation*}
the numerator and denominator of which fraction vanish again for $x=a$. Therefore, for the same reason, substitute their respective differentials again and  this fraction will result

\begin{equation*}
\frac{6x-8a+2a^4:(2ax-aa)^{\frac{3}{2}}}{2-2a^3:(2ax-xx)^{\frac{3}{2}}},
\end{equation*}
whose numerator and denominator vanish again in the case $x=a$. Therefore, we  substitute their differentials for them

\begin{equation*}
\frac{6-6a^5:(2ax-aa)^{\frac{5}{2}}}{6a^3(a-x):(2ax-xx)^{\frac{5}{2}}}=\frac{1-a^5:(2ax-aa)^{\frac{5}{2}}}{a^3(a-x):(2ax-xx)^{\frac{5}{2}}}.
\end{equation*}
But   both the numerator and the denominator vanish again here for $x=a$. Therefore,  having substituted their differentials for them, it will result

\begin{equation*}
\frac{5a^6:(2ax-aa)^{\frac{7}{2}}}{-(5a^5-8a^4x+4a^3xx):(2ax-xx)^{\frac{7}{2}}}.
\end{equation*}
Now, finally write $a$ instead of $x$ and this determined fraction will result $\frac{5:a}{-1:a^2}=-5a$ which is the value of the fraction in question.\\[2mm]
But if, one puts $x=a+t$ before this investigation is done, the propounded fraction will be transformed into this one

\begin{equation*}
\frac{2a^3+2a^2t-att+t^3-2a^2\sqrt{aa+2at}}{-2aa+tt+2a\sqrt{aa-tt}};
\end{equation*}
since this has the form $\frac{0}{0}$, if one puts $t=0$, write $dt$ instead of $t$ and it will be

\begin{equation*}
\frac{2a^3+2a^2dt-adt^2+dt^3-2a^2\sqrt{aa+2adt}}{-2aa+dt^2+2a\sqrt{aa-dt^2}}.
\end{equation*}
Now, convert the irrational formulas into series, which must be continued until the rational terms are no longer cancelled:

\begin{equation*}
\sqrt{aa+2adt}=a+dt-\frac{dt^2}{2a}+\frac{dt^3}{2aa}-\frac{5dt^4}{8a^3},
\end{equation*}
\begin{equation*}
\sqrt{aa-dt^2}=a-\frac{dt^2}{2a}-\frac{dt^4}{8a^3};
\end{equation*}
having substituted these values, this fraction will result

\begin{equation*}
\frac{5dt^4:4a}{-dt^4:4aa}=-5a
\end{equation*}
which is the value of the propounded fraction already found before.

\subsection*{Example 4}

\textit{To find the value of the fraction $\frac{a+\sqrt{2aa-2ax}-\sqrt{2ax-xx}}{a-x+\sqrt{aa-xx}}$ in the case $x=a$.}\\[2mm]
Having substituted their differentials for the numerator and the denominator, the following fraction will result, which in the case $x=a$ will be equal to the propounded one:

\begin{equation*}
\frac{-a:\sqrt{2aa-2ax}-(a-x):\sqrt{2ax-xx}}{-1-x:\sqrt{aa-xx}}
\end{equation*}
whose numerator and denominator  become infinite in the case $x=a$. But if both are multiplied by $-\sqrt{a-x}$, one will have

\begin{equation*}
\frac{a:\sqrt{2a}+(a-x)^{\frac{3}{2}}:\sqrt{2ax-xx}}{\sqrt{a-x}+x:\sqrt{a+x}}
\end{equation*}
which, having put $x=a$, will give the determined value $\frac{a:\sqrt{2a}}{a:\sqrt{2a}}=1$ which is therefore equal to the propounded fraction in the case $x=a$.

\paragraph*{§359}

Therefore, if one has a fraction $\frac{P}{Q}$ whose numerator and denominator vanish in the case $x=a$, one will be able to assign its value by means of the usual rules of differentiation and it will not be necessary to go back to the differentials which we treated in the preceding chapter. For, having taken the differentials, in the case $x=a$ the propounded fraction $\frac{P}{Q}$  will become equal to the fraction $\frac{dP}{dQ}$; if its numerator and denominator  take on finite values for $x=a$, one will find the value of the propounded fraction; but if the one becomes $=0$ while the other remains finite,  the fraction will be either $=0$ or $=\infty$, depending on whether the numerator or the denominator vanishes. But if one of the two or even both become $=\infty$, what happens, if  one divides by quantities vanishing in the case $x=a$, then by multiplying both by these divisors the inconvenience is removed as it happened in the preceding example. But if so the numerator as the denominator vanish again, then just, as it was done in the beginning, the differentials are  to be taken once more such that this fraction $\frac{ddP}{ddQ}$ results, which in the case $x=a$ will still be equal to the propounded one; and if the same happens again in this fraction, i.e. that $=\frac{0}{0}$,   substitute this one $\frac{d^3P}{d^3Q}$ for it and so forth, until one reaches a fraction exhibiting a determined value, either finite or infinitely large or infinitely small.  In the third example it was necessary to proceed to the fraction $\frac{d^4P}{d^4 Q}$, before it was possible to assign the value of the propounded fraction $\frac{P}{Q}$.

\paragraph*{§360}

The use of this investigation becomes clear in the definition of the sums of series which we found above (chap. II § 22), if one puts $x=1$. For, from the results  derived there it follows that:

\begin{alignat*}{9}
&x&&+x^2&&+x^3 &&+\cdots &&+x^n  &&=\frac{x-x^{n+1}}{1-x} \\[3mm]
&x&&+x^3&&+x^5 &&+\cdots &&+x^{2n-1}&&=\frac{x-x^{2n+1}}{1-xx} \\[3mm]
&x&&+2x^2&&+3x^3&&+\cdots &&+nx^n &&=\frac{x-(n+1)x^{n+1}+nx^{n+2}}{(1-x)^2}\\[3mm]
&x &&+3x^3 &&+5x^5 &&+\cdots &&+(2n-1)x^{2n-1}&&=\frac{x+x^3-(2n+1)x^{2n+1}+(2n-1)x^{2n+3}}{(1-xx)^2} \\[3mm]
&x &&+4x^2 &&+9x^3 &&+\cdots &&+n^2 x^n &&= \frac{x+x^2-(n+1)x^{n+1}+(2nn+2n-1)x^{n+2}-nnx^{n+3}}{(1-x)^3} \\[3mm]
& && && && && &&\text{etc.}
\end{alignat*} 
If the sums of these series are desired in the case $x=1$,  so the numerator as the denominator vanish in these expressions. Therefore, the values  of these sums in the case $x=1$ can be defined by  the method explained here. Since the same sums are known from elsewhere, from the agreement with those known values the validity of this method will become more clear.

\subsection*{Example 1}

\textit{To define the value of this fraction $\frac{x-x^{n+1}}{1-x}$ in the case $x=1$ which will exhibit the sum of the series $1+1+1+\cdots +1$ consisting of $n$ terms, which sum therefore will be $=n$.}\\[2mm]
Since in the case $x=1$ the numerator and the denominator vanish, substitute their differentials for them and one will have

\begin{equation*}
\frac{1-(n+1)x^n}{-1},
\end{equation*}
which for $x=1$ gives $n$ for the sum of the series in question.

\subsection*{Example 2}

\textit{To define the value of the fraction $\frac{x-x^{2n+1}}{1-xx}$ in the case $x=1$ which  will exhibit the sum of the series $1+1+1+\cdots 1$ consisting of $n$ terms, which sum will therefore be $=n$.}\\[2mm]
Having taken the differentials, the propounded fraction will be transformed into this one

\begin{equation*}
\frac{1-(2n+1)x^{2n}}{-2x}
\end{equation*}
whose value for $x=1$ will be $=n$.

\subsection*{Example 3}

\textit{To find the value of the fraction $\frac{x-(n+1)x^{n+1}+nx^{n+2}}{(1-x)^2}$ in the case $x=1$ which will express the sum of the series $1+2+3+\cdots +n$ which is known to be $=\frac{nn+n}{2}$.}\\[2mm]
Having taken the differentials, one gets to this fraction

\begin{equation*}
\frac{1-(n+1)^2x^n+n(n+2)x^{n+1}}{-2(1-x)},
\end{equation*}
whose numerator and denominator still vanish in the case $x=1$. Therefore,  take the differentials again so that this fraction results

\begin{equation*}
\frac{-n(n+1)^2x^{n-1}+n(n+1)(n+2)x^n}{2},
\end{equation*}
which for $x=1$ goes over into $\frac{n(n+1)}{2}=\frac{nn+n}{2}$ i.e. the sum of the propounded series.

\subsection*{Example 4}

\textit{To find the value of the fraction $\frac{x+x^3-(2n+1)x^{2n+1}+(2n-1)x^{2n+3}}{(1-xx)^2}$ in the case $x=1$ which will express the sum of the series $1+3+5+\cdots +(2n-1)$ which is known to be $=nn$.}\\[2mm]
Having substituted the differentials for the numerator and the denominator, this fraction results

\begin{equation*}
\frac{1+3xx-(2n+1)^2x^{2n}+(2n-1)(2n+3)x^{2n+2}}{-4x(1-xx)};
\end{equation*}
because this still suffers from the same inconvenience, i.e. that it goes over into $=\frac{0}{0}$ for $x=1$, take the differentials again

\begin{equation*}
\frac{6x-2n(2n+1)^2x^{2n-1}+(2n-1)(2n+2)(2n+3)x^{2n+1}}{-4+12xx},
\end{equation*}
which for $x=1$ goes over into 

\begin{equation*}
\frac{6-2n(2n+1)^2+(2n-1)(2n+2)(2n+3)}{8}=nn.
\end{equation*}

\subsection*{Example 5}

\textit{To find the value of the fraction}

\begin{equation*}
\frac{x+x^2-(n+1)^2x^{n+1}+(2nn+2n-1)x^{n+2}-nnx^{n+3}}{(1-x)^3}
\end{equation*}
\textit{in the case $x=1$ which will give the sum of the series $1+4+9+\cdots +n^2$ which is known to be $=\frac{1}{3}n^3+\frac{1}{2}n^2+\frac{1}{6}n$.}\\[2mm]
Having taken the differentials of the numerator and the denominator, it will be

\begin{equation*}
\frac{1+2x-(n+1)^3 x^n+(n+2)(2nn+2n-1)x^{n+1}-nn(n+3)x^{n+2}}{-3(1-x)^2}
\end{equation*}
since in it the numerator and the denominator for $x=1$ vanish again, take the second differentials

\begin{equation*}
\frac{2-n(n+1)^3 x^{n-1}+(n+1)(n+2)(2nn+2n-1)x^n-n^2(n+2)(n+3)x^{n+1}}{6(x-1)}.
\end{equation*}
Since the same inconvenience is still present here, proceed to the third differentials that this fraction results

\begin{equation*}
\frac{-n(n-1)(n+1)x^{n-2}+n(n+1)(n+2)(2nn+2n-1)x^{n-1}-n^2 (n+1)(n+2)(n+3)x^n}{-6}
\end{equation*}
which for $x=1$ finally goes over into this determined form

\begin{equation*}
\frac{-n(n-1)(n+1)^3+n(n+1)(n+2)(nn-n-1)}{-6}=\frac{n(n+1)(2n+1)}{6}=\frac{1}{3}n^3+\frac{1}{2}n^2+\frac{1}{6}n;
\end{equation*} 
which is the value  we found to express the mentioned series.

\subsection*{Example 6}

\textit{Let this fraction be propounded $\frac{x^m-x^{m+n}}{1-x^{2p}}$ whose value in the case $x=1$ is  to be assigned.}\\[2mm]
Since this fraction is the product of these two $\frac{x^m}{1+x^p}\cdot \frac{1-x^n}{1-x^p}$, but the value of the first factor in the case $x=1$ is $=\frac{1}{2}$, it is only necessary that in the same case the value of the other factor $\frac{1-x^n}{1-x^p}$  is found which, having taken the differentials, will be $=\frac{nx^{n-1}}{px^{p-1}}=\frac{n}{p}$; therefore, the value of the propounded fraction in the case $x=1$ will be $=\frac{n}{2p}$. The same value results, if the differentials are taken immediately  in the propounded fraction; for, it will be

\begin{equation*}
\frac{mx^{m-1}-(m+n)x^{m+n-1}}{-2px^{2p-1}},
\end{equation*}
whose value, having put $x=1$, will be $=\frac{-n}{-2p}=\frac{n}{2p}$ as before.

\paragraph*{§361}

The same method is to be used, if in the propounded fraction $\frac{P}{Q}$ either the numerator or the denominator or both were a transcendental number. To explain these operations more clearly, it seems advisable to add the following examples.

\subsection*{Example 1}

\textit{Let this fraction be propounded $\frac{a^n-x^n}{\log a- \log x}$ whose value in the case $x=a$ is in question.}\\[2mm]
Having taken the differentials, one immediately gets to this form

\begin{equation*}
\frac{-nx^{n-1}}{-1:x}=nx^n
\end{equation*}
whose value for $x=a$ will be $na^n$.

\subsection*{Example 2}

\textit{Let this fraction be propounded $\frac{\log x}{\sqrt{1-x}}$ whose value in the case $x=1$ is in question.}\\[2mm]
Having taken the differentials of the numerator and the denominator, this fraction results

\begin{equation*}
\frac{1:x}{-1:2 \sqrt{1-x}}=\frac{-2 \sqrt{1-x}}{x};
\end{equation*}
since its value for $x=1$ is $=0$, it follows that the fraction $\frac{\log x}{\sqrt{1-x}}$ vanishes in the case $x=1$.

\subsection*{Example 3}

\textit{Let this fraction be propounded $\frac{a-x-a\log a- a \log x}{a-\sqrt{2ax-xx}}$ whose value for $x=a$ is in question, in which case the numerator and the denominator vanish.}\\[2mm]
Having differentiated the numerator and denominator according to the rule, it will be

\begin{equation*}
\frac{-1+a:x}{-(a-x):\sqrt{2ax-xx}}=\frac{(a-x)\sqrt{2ax-xx}}{-x(a-x)};
\end{equation*}
even if here the numerator and the denominator still vanish in the case $x=a$, nevertheless, since both are divisible by $a-x$, one will have this fraction $-\sqrt{\frac{2a-x}{x}}$ whose value in the case $x=a$ is determined and $=-1$; and therefore, the propounded fraction goes over into $-1$, if one puts $x=a$.

\subsection*{Example 4}

\textit{Let this fraction be propounded $\frac{e^x-e^{-x}}{\log (1+x)}$ whose value for $x=0$ is in question.}\\[2mm]
Having taken the differentials, one will have this fraction

\begin{equation*}
\frac{e^x+e^{-x}}{1:(1+x)}
\end{equation*} 
which for $x=0$ gives $2$ for the value in question.

\subsection*{Example 5}

\textit{To find the value of this fraction $\frac{e^x-1-\log (1+x)}{xx}$ in the case $x=0$.}\\[2mm]
If  their differentials are substituted for the numerator and the denominator, this fraction will result

\begin{equation*}
\frac{e^x-1:(1+x)}{2x},
\end{equation*}
since which goes over into $\frac{0}{0}$, if one puts $x=0$, take the differentials again that one has

\begin{equation*}
\frac{e^x+1:(1+x)^2}{2}
\end{equation*}
which for $x=0$ yields $\frac{1+1}{2}=1$. The same is plain, if  one immediately substitutes $0+dx$ for $x$, because 

\begin{equation*}
e^{dx}=1+dx+\frac{dx^2}{2}+\text{etc.} \quad \text{and} \quad \log (1+dx)=dx-\frac{1}{2}dx^2+\text{etc.},
\end{equation*}
\begin{equation*}
\frac{e^{\log x}-1-\log (1+dx)}{dx^2}=\frac{dx^2}{dx^2}=1.
\end{equation*}

\subsection*{Example 6}

\textit{Let the value of the fraction $\frac{x^n}{\log x}$  be in question in the case  $x=\infty$.}\\[2mm]
To reduce this fraction to a form which  goes over into $\frac{0}{0}$ in this case, represent it this way

\begin{equation*}
\frac{1:\log x}{1:x^n};
\end{equation*}
for, in the case $x= \infty$ so the numerator as the denominator will vanish. Further, put $x=\frac{1}{y}$ such that in the case $x=\infty$ we have $y=0$, and  this fraction is propounded

\begin{equation*}
-\frac{1: \log y}{y^n}
\end{equation*}
whose value in the case $y=0$ must be investigated. But, having taken the differentials, it will be $\frac{1:y(\log y)^2}{ny^{n-1}}=\frac{1:(\log y)^2}{ny^n}$; since for $y=0$ it goes over into $\frac{0}{0}$, take the differentials again and it will be $\frac{-2:(\log y)^3}{n^2 y^n}$; since here the same inconvenience is present, if  the differentials are taken again, $\frac{6:(\log y)^4}{n^3 y^n}$ will result and so, no matter how far we proceed, always the same inconvenience occurs. Therefore, to find the value in question despite this obstacle, let $s$ be the value of the fraction $-\frac{1: \log y}{y^n}$ in the case  $y=0$, and because in the same case also

\begin{equation*}
s=\frac{1:(\log y)^2}{ny^n},
\end{equation*}
from that equation it will be

\begin{equation*}
ss=\frac{1:(\log y)^2}{y^{2n}}
\end{equation*}
which divided by the latter will give

\begin{equation*}
s= \frac{ny^n}{y^{2n}}=\frac{n}{y^n},
\end{equation*}
from which it is understood that in the case $y=0$ $s$ becomes infinite. Therefore, the value of the fraction $-\frac{1: \log y}{y^n}$  will be infinite in the case $y=0$ and for $y=dx$ $\frac{1}{\log dx}$ will have an infinite ratio to $dx^n$ as we mentioned already above [§ 351].

\subsection*{Example 7}

\textit{Let the value of the fraction $\frac{x^n}{e^{-1:x}}$  be in question in the case $x=0$, in which both the numerator and the denominator vanish.}\\[2mm]
Let   $\frac{x^n}{e^{-\frac{1}{x}}}=s$ in this case; therefore, having taken the differentials, it will also be

\begin{equation*}
s= \frac{nx^{n-1}}{e^{-1:x}:xx}=\frac{nx^{n+1}}{e^{-1:x}},
\end{equation*}
and since here the same inconvenience occurs and always occurs again, no matter how far the differentiations are continued, we argue exactly as before. The first equation gives

\begin{equation*}
x^n =e^{-1:x}s \quad \text{and} \quad x^{n(n+1)}=e^{-(n+1):x}s^{n+1};
\end{equation*}
the other equation gives

\begin{equation*}
x^{n+1}=e^{-1:x}s:n,
\end{equation*}
whence 

\begin{equation*}
x^{n(n+1)}=e^{-n:x}s^n:n^n,
\end{equation*}
which value equated to the latter will give

\begin{equation*}
e^{-1:x}sn^n=1
\end{equation*}
and hence

\begin{equation*}
s= \frac{1}{n^n e^{-1:x}}=\infty,
\end{equation*}
if $x=0$. Therefore, having put $x$ to be infinitely small,  $dx^n$ will have an infinitely large ratio to $e^{-1:dx}$, no matter which finite number is substituted for $n$; therefore, it follows that $e^{-1:dx}$ is infinitely small and homogeneous to $dx^m$, if $m$ was an infinitely large number.

\subsection*{Example 8}

\textit{Let the value of the fraction $\frac{1- \sin x +\cos x}{\sin x+\cos x -1}$  be in question in the case  $x=\frac{\pi}{2}$ or equal to the arc of  ~$90^{\circ}$.}\\[2mm]
Having taken the differentials, one will obtain this fraction

\begin{equation*}
\frac{-\cos x -\sin x}{\cos x -\sin x}
\end{equation*}
which for $x=\frac{\pi}{2}$, because of $\sin x=1$ and $\cos x=0$, goes over into $1$ such that $1$ is the value in question of the propounded fraction. The same is obvious without differentiation; for, since  $\cos x =\sqrt{(1+\sin x)(1-\sin x)}$, the propounded fraction goes over into 

\begin{equation*}
\frac{\sqrt{1- \sin x}+\sqrt{1+\sin x}}{\sqrt{1+\sin x}-\sqrt{1-\sin x}}
\end{equation*}
which evidently becomes $=1$, if  $\sin x=1$..

\subsection*{Example 9}

\textit{To find the value of this expression $\frac{x^x-x}{1-x+\log x}$ in the case  $x=1$.}\\[2mm]
Having substituted their differentials for the numerator and the denominator, this fraction will result

\begin{equation*}
\frac{x^x(1+\log x)-1}{-1+1:x};
\end{equation*}
since this  also becomes $=\frac{0}{0}$ for $x=1$, take the differentials again that this fraction results

\begin{equation*}
\frac{x^x(1+\log x)^2+x^x:x}{-1:xx}
\end{equation*}
which, having put $x=1$, goes over into $-2$ which is the value of the propounded fraction in the case $x=1$.

\paragraph*{§362}

Since here we decided to treat all expressions which in certain cases seem to have undetermined values,  not only those fractions $\frac{P}{Q}$  whose numerator and denominator vanish in certain cases extend to this, but also fractions  whose numerator and denominator  become infinite  in a certain case are to be included, since their values  seem to be undetermined as well. If $P$ and $Q$ were functions of $x$ of such a kind that in a certain case $x=a$ both become infinite and the fraction $\frac{P}{Q}$ takes on this form $\frac{\infty}{\infty}$, since two infinities as two zeros can have any ratio to each other, hence the value can not be known at all. This case can be reduced to the preceding by transforming the fraction $\frac{P}{Q}$ into this form $\frac{1:Q}{1:P}$, the numerator and denominator of which fraction now vanish in the case $x=a$; and hence its value can be found applying the method discussed before. But even without this transformation the value will be found, if  not $a$ but $a+dx$ is substituted for $x$; having done this, not an absolute infinity $\infty$ will result, but it will be expressed as $\frac{1}{dx}$ or $\frac{A}{dx^n}$; even if these expressions are equally infinite as $\infty$, nevertheless, having done the comparison of $dx$ to its powers, the value in question will easily be calculated.

\paragraph*{§363}

Also the products consisting of two factors  the of which  vanishes, but the other goes over into infinity in the  case $x=a$ extend to the same class; for, since every quantity can be represented by a product of this kind $0 \cdot \infty$, its value seems to be undetermined. Let $PQ$ be a product of this kind in which $P=0$ and $Q=\infty$ for $x=a$; its value will be found by means of the rules given before, if one puts $Q=\frac{1}{R}$; for, then the product $PQ$ will be transformed into the fraction $\frac{P}{R}$ whose numerator and denominator both vanish in the case $x=a$; and hence its value can be investigated by   the method explained before.\\[2mm]
So, if the value of this product

\begin{equation*}
(1-x) \tan \frac{\pi x}{2}
\end{equation*}
is in question in the case $x=1$, in which  $1-x=0$ and $\tan \frac{\pi x}{2}=\infty$, convert it into this fraction

\begin{equation*}
\frac{1-x}{\cot \frac{1}{2}\pi x}
\end{equation*}
whose numerator and denominator vanish in the case $x=1$. Therefore, since the differential of the numerator is $1-x=-dx$ and the differential of the denominator $\cot \frac{\pi x}{2}$ is $=- \frac{\pi dx:2}{(\sin \frac{1}{2}\pi x)^2}$, in the case $x=1$ the value of the propounded fraction will be

\begin{equation*}
=\frac{2 }{\pi}\sin \frac{\pi x}{2}\cdot \sin \frac{\pi x}{2}=\frac{2 }{\pi},
\end{equation*}
because of $\sin \frac{\pi }{2}=1$. 

\paragraph*{§364}

But especially those expressions  which, if  a certain value is attributed to $x$, go over into $\infty -\infty$ are included here; for, since two infinities can differ by any finite quantity, it is obvious that in this case the value of the expression is not determined, if the difference of those two infinities can not be assigned. Therefore, this case occurs, if a function $P-Q$ of this kind is propounded in which for $x=a$  so $P=\infty$ as $Q=\infty$, in which case the value in question can not be assigned that easily applying the rules given before. For, even if, having put that in this case $P-Q=f$, one sets $e^{P-Q}=e^f$ such that  $e^f= \frac{e^{-Q}}{e^{-P}}$ where in the case $x=a$ so the numerator $e^{-Q}$ as the numerator $e^{-P}$ vanishes, if the rule given before is applied here, it will be $e^f=\frac{e^{-Q}dQ}{e^{-P}dP}$, whence, because of $e^f= \frac{e^{-Q}}{e^{-P}}$, it would be $1=\frac{dQ}{dP}$ and hence the value in question of $f$ will not be found from this. If $P$ and $Q$ are algebraic quantities, since these only become infinite, if there are fractions whose denominators vanish, then $P-Q$ can be contracted into one single fraction whose denominator will vanish in the same way. If, having done this, the numerator vanishes, the value can be defined by the method explained above; but if the numerator does not vanish,  its value will indeed be infinite.\\[2mm]
So if the value of this expression

\begin{equation*}
\frac{1}{1-x}-\frac{2}{1-xx}
\end{equation*}
is desired in the case $x=1$, since it goes over into

\begin{equation*}
\frac{-1+x}{1-xx}=\frac{-1}{1+x},
\end{equation*}
it is plain that value in question is $=-\frac{1}{2}$.

\paragraph*{§365}

But if the functions $P$ and $Q$ were transcendental, this transformation  would lead to most cumbersome calculations in most cases. Therefore, it will be convenient to use  a direct method is in these cases  and to write $x=a+ \omega$, while $\omega$  is an infinitely small quantity, for which one can take $dx$, instead of $x=a$, in which case both quantities $P$ and $Q$ go over into infinity. If, having done this,  $P=\frac{A}{\omega}+B$ and $Q=\frac{A}{\omega}+C$, it is obvious that the function $P-Q$ will go over into $B-C$ which value will be finite. Therefore, we will illustrate this method to investigate values of functions of this kind in the following examples.

\subsection*{Example 1}

\textit{Let the value of this expression $\frac{x}{x-1}-\frac{1}{\log x}$ be in question in the case  $x=1$.}\\[2mm]
Since so $\frac{x}{x-1}$ as $\frac{1}{\log x}$ become infinite for $x=1$, put $x=1+ \omega$ and the propounded expression will be transformed into this one

\begin{equation*}
\frac{1+\omega}{\omega}- \frac{1}{\log (1+\omega)}.
\end{equation*}
Therefore, since 

\begin{equation*}
\log (1+\omega)= \omega -\frac{1}{2}\omega^2+\frac{1}{3}\omega^3-\text{etc.} =\omega \left(1-\frac{1}{2}\omega +\frac{1}{3}\omega^2-\text{etc.}\right),
\end{equation*}
one will have

\begin{equation*}
\frac{(1+\omega)(1-\frac{1}{2}\omega+\frac{1}{3}\omega^2-\text{etc.})}{\omega(1-\frac{1}{2}\omega +\frac{1}{3}\omega^2-\text{etc.})}=\frac{\frac{1}{2}\omega -\frac{1}{6}\omega^2+\text{etc.}}{\omega\left(1-\frac{1}{2}\omega+\frac{1}{3}\omega^2-\text{etc.}\right)}=\frac{\frac{1}{2}-\frac{1}{6}\omega +\text{etc.}}{1-\frac{1}{2}\omega+\frac{1}{3}\omega^2-\text{etc.}}.
\end{equation*}
Now, having put $\omega$ to be infinitely small or $\omega =0$, it is obvious that the value in question $=\frac{1}{2}$.

\subsection*{Example 2}

\textit{While $e$ denotes the number whose hyperbolic logarithm is $=1$ and $\pi$ the half of the circumference of the circle whose radius is $=1$, to investigate the value of the expression $\frac{\pi x-1}{2xx}+\frac{\pi}{x(e^{2\pi x}-1)}$ in the case $x=0$.}\\[2mm]
This expression exhibits the sum of this series

\begin{equation*}
\frac{1}{1+xx}+\frac{1}{4+xx}+\frac{1}{9+xx}+\frac{1}{16+xx}+\frac{1}{25+xx}+\text{etc.};
\end{equation*}
hence, if one puts $x=0$, the sum of this series must result

\begin{equation*}
\frac{1}{1}+\frac{1}{4}+\frac{1}{9}+\frac{1}{16}+\text{etc.}
\end{equation*}
which is known to be $=\frac{\pi \pi}{6}$. But for $x=0$ the value of the propounded expression 

\begin{equation*}
\frac{\pi x-1}{2xx}+\frac{\pi}{x(e^{2 \pi x}-1)}
\end{equation*}
seems to be most undetermined because of all the infinite terms. Therefore, put $x=\omega$ while $\omega$ denotes an infinitely small quantity and the first term $\frac{\pi x-1}{2xx}$ goes over into 

\begin{equation*}
-\frac{1}{2\omega}^2+\frac{\pi}{2\omega}.
\end{equation*}
Further, since 

\begin{equation*}
e^{-2 \pi \omega}-1=2 \pi \omega +2 \pi^2 \omega^2+\frac{4}{3}\pi^3 \omega^3 +\text{etc.},
\end{equation*}
the other side  $\frac{\pi}{x(e^{2 \pi x}-1)}$ goes over into 

\begin{equation*}
\frac{\pi}{\omega(2 \pi \omega +2 \pi^2 \omega^2 +\frac{4}{3}\pi^3 \omega^3)+\text{etc.}}=\frac{1}{2\omega^2(1+\pi \omega +\frac{2}{3}\pi^2 \omega^2+\text{etc.})}.
\end{equation*}
But 

\begin{equation*}
\frac{1}{1+\pi \omega +\frac{2}{3}\pi^2 \omega^2+\text{etc.}}=1- \pi \omega +\frac{1}{3}\pi^2 \omega^2-\text{etc.},
\end{equation*}
whence the second term becomes

\begin{equation*}
=\frac{1}{2\omega^2}-\frac{\pi}{2 \omega}+\frac{1}{6}\pi^2-\text{etc.};
\end{equation*}
if to this the first is added, $\frac{1}{6}\pi^2$ results which is the value in question of the propounded expression in the case $x=0$.\\[2mm]
The same can also be achieved applying the method of fractions whose numerator and denominator vanish in a certain case; for, the propounded expression  is transformed into this fraction

\begin{equation*}
\frac{\pi xe^{2 \pi x}-e^{2 \pi x}+\pi x+1}{2xx e^{2 \pi x}-2xx}
\end{equation*}
whose numerator and denominator vanish in the case $x=0$. Therefore, having taken the differentials, this fraction results

\begin{equation*}
\frac{\pi e^{2\pi x}+2 \pi \pi x e^{2 \pi x}-2 \pi e^{2 \pi x}+\pi}{4xe^{2 \pi x}+4 \pi xx e^{2 \pi x}-4x}
\end{equation*}
or this one

\begin{equation*}
\frac{\pi - \pi e^{2 \pi x}+2 \pi \pi e^{2 \pi x}}{4xe^{2 \pi x}+4 \pi xxe^{2 \pi x}-4x}
\end{equation*}
whose numerator and denominator still vanish, if one puts $x=0$. Therefore, having taken the differentials again, one will have

\begin{equation*}
\frac{-2 \pi \pi e^{2 \pi x}+2 \pi \pi e^{2 \pi x}+4 \pi^3 x e^{2 \pi x}}{4e^{2\pi x}+8 \pi x e^{2 \pi x}+ 8 \pi x e^{2 \pi x}+ 8 \pi^2 xx e^{2 \pi x}-4}
\end{equation*}
or

\begin{equation*}
\frac{\pi^3 x e^{2\pi x}}{e^{2 \pi x}+4 \pi x e^{2\pi x}+2 \pi^2 x^2 e^{2 \pi x}-1}
\end{equation*}
or

\begin{equation*}
\frac{\pi^3 x}{1+4 \pi x+2 \pi^2 x^2-e^{-2\pi x}}
\end{equation*}
whose numerator and denominator still vanish in the case $x=0$. Therefore, take the differentials again

\begin{equation*}
\frac{\pi^3}{4 \pi +4 \pi^2x+2 \pi e^{-2 \pi x}}
\end{equation*}
which fraction for $x=0$ goes over into $\frac{\pi^2}{6}$ as before.

\subsection*{Example 3}

\textit{While $e$ and $\pi$ retain the same values, let the value of the following expression be in question in the case $x=0$:}

\begin{equation*}
\frac{\pi}{4x}-\frac{\pi}{2x(e^{\pi x}+1)}.
\end{equation*}\\[2mm]
This expression is transformed into this one

\begin{equation*}
\frac{\pi e^{\pi x}-\pi}{4x e^{\pi x}+4x}
\end{equation*}
whose numerator and denominator  vanish in the case $x=0$. Therefore, put $x=\omega$, and because

\begin{equation*}
e^{\pi \omega}=1+ \pi \omega +\frac{1}{2}\pi^2 \omega^2 +\frac{1}{6}\pi^3 \omega^3+ \text{etc.},
\end{equation*}
the propounded formula is transformed into this one

\begin{equation*}
\frac{\pi^2 \omega +\frac{1}{2}\pi^3 \omega^2+\frac{1}{6}\pi^4 \omega^3 +\text{etc.}}{8 \omega +4 \pi \omega^2+2 \pi^2 \omega^3+ \text{etc.}}
\end{equation*}
which, having put $\omega$ to be infinitely small, immediately gives $\frac{1}{8}\pi^2$ which is the value in question of the propounded expression in the case $x=0$. For, the propounded expression $\frac{\pi}{4x}-\frac{\pi}{2x(e^{\pi x}+1)}$ will indeed exhibit the sum of this series

\begin{equation*}
\frac{1}{1+xx}+\frac{1}{9+xx}+\frac{1}{25+xx}+\frac{1}{49+xx}+\text{etc.}
\end{equation*}
whose sum for $x=0$ gives $=\frac{1}{8}\pi^2$.

\subsection*{Example 4}

\textit{To find the value of the expression $\frac{1}{2xx}-\frac{\pi}{2x \tan \pi x}$ in the case $x=0$.}\\[2mm]
This propounded formula $\frac{1}{2xx}-\frac{\pi}{2x \tan \pi x}$ expresses the sum of this infinite series

\begin{equation*}
\frac{1}{1-xx}+\frac{1}{4-xx}+\frac{1}{9-xx}+\frac{1}{16-xx}+\text{etc.}
\end{equation*} 
If one therefore puts $x=0$, the sum of the series

\begin{equation*}
1+\frac{1}{4}+\frac{1}{9}+\frac{1}{16}+\text{etc.}
\end{equation*}
results which is $=\frac{1}{6}\pi \pi$. Since  $\tan \pi = \frac{\sin \pi x}{\cos \pi x}$, the propounded expression will take on this form

\begin{equation*}
\frac{1}{2xx}-\frac{\pi \cos \pi x}{2x \sin \pi x}= \frac{\sin \pi x- \pi x \cos \pi x}{2xx \sin \pi x}
\end{equation*}
whose numerator and denominator vanish in the case $x=0$. Therefore, put $x= \omega$, and because 

\begin{equation*}
\sin \pi \omega = \pi \omega -\frac{1}{6}\pi^3 \omega^3 +\text{etc.}, \quad \cos \pi \omega =1- \frac{1}{2}\pi^2 \omega^2 +\text{etc.},
\end{equation*}
the propounded expression will be

\begin{equation*}
\frac{\pi \omega -\frac{1}{6}\pi^3 \pi^3+\text{etc.}- \pi \omega +\frac{1}{2}\pi^3 \omega^3 -\text{etc.}}{2 \pi \omega^3 -\frac{1}{3}\pi^3 \omega^5 +\text{etc.}}=\frac{\frac{1}{3}\pi^3 \omega^3 -\text{etc.}}{2 \pi \omega^3 -\text{etc.}}
\end{equation*}
which, because of the infinitely small $\omega$, gives $\frac{1}{6}\pi^2$.

\subsection*{Example 5}

\textit{Since the sum of the following infinite series is known, for, we have}

\begin{equation*}
\frac{1}{1-xx}+\frac{1}{9-xx}+\frac{1}{25-xx}+\frac{1}{49-xx}+\text{etc.}=\frac{\pi \sin \frac{1}{2}\pi x}{4x \cos \frac{1}{2}\pi x},
\end{equation*}
\textit{to find its sum for $x=0$.}\\[2mm]
Since 

\begin{equation*}
\sin \frac{1}{2}\pi x =\frac{1}{2}\pi x- \frac{1}{48}\pi^3 x^3+\text{etc.} \quad \text{and} \quad \cos \frac{1}{2}\pi x= 1- \frac{1}{8}\pi^2 x^2+ \text{etc.},
\end{equation*}
the propounded expression will be

\begin{equation*}
=\frac{\frac{1}{2}\pi^2 x- \frac{1}{48}\pi^4x^3 + \text{etc.}}{4x-\frac{1}{2}\pi^2 x^3 +\text{etc.}}=\frac{\frac{1}{2}\pi^2- \frac{1}{48}\pi^4 x^2+ \text{etc.}}{4- \frac{1}{2}\pi^2 x^2+\text{etc.}};
\end{equation*}
if in this $x=0$, the value will obviously be $=\frac{1}{8}\pi^2$ which is the sum of the series

\begin{equation*}
1+\frac{1}{9}+\frac{1}{25}+\frac{1}{49}+\text{etc.}
\end{equation*}
as it was demonstrated above in several ways. But if  one takes any even number for $x$, the sum of the propounded series is always $=0$.

\paragraph*{§366}

In these series we treated in the last two examples and containing the variable letter $x$ one can attribute  values of such a kind to $x$ that certain terms grow to infinity in which cases the sum of the whole series will be infinite. So, one term of the series

\begin{equation*}
\frac{1}{1-xx}+\frac{1}{4-xx}+\frac{1}{9-xx}+\frac{1}{16-xx}+\text{etc.},
\end{equation*}
if  one substitutes any integer for $x$, because of the vanishing denominator, always becomes infinite and therefore the series will become infinite. But if this infinite term is thrown out of the series,  the remaining sum will without any doubt be finite and will be expressed by  the first series decreased by this infinite term in this way $\infty -\infty$; therefore, one will be able to find a determined value, which will be  seen more clearly from the examples added below.

\subsection*{Example 1}

\textit{To find the sum of the series}

\begin{equation*}
\frac{1}{1-xx}+\frac{1}{4-xx}+\frac{1}{9-xx}+\frac{1}{16-xx}+\text{etc.},
\end{equation*}
\textit{after having subtracted the first term, in the case $x=1$, in which the first term is infinite.}\\[2mm]
Since the sum  in general is

\begin{equation*}
=\frac{1}{2xx}-\frac{\pi}{2x \tan \pi x},
\end{equation*}
the sum in question will be

\begin{equation*}
=\frac{1}{2xx}-\frac{\pi}{2x \tan \pi x}-\frac{1}{1-xx}
\end{equation*}
for $x=1$. Let $x=1+ \omega$ and  for the sum in question one will have

\begin{equation*}
\frac{1}{2(1+2 \omega + \omega \omega)}- \frac{\pi}{2(1+\omega)\tan (\pi +\omega \pi)}+\frac{1}{2 \omega + \omega \omega}.
\end{equation*}
But 

\begin{equation*}
\tan (\pi + \omega \pi)= \tan \omega \pi = \pi \omega +\frac{1}{3}\pi^3 \omega^3+ \text{etc.}
\end{equation*}
therefore, because the first term $\frac{1}{2xx}$ for $x=1$ has the determined value $\frac{1}{2}$, only the remaining terms are to be considered which will be

\begin{equation*}
\frac{1}{\omega(2+\omega)}-\frac{\pi}{2 \omega(1+ \omega)(\pi +\frac{1}{3}\pi^3 \omega^2)}=\frac{1}{\omega(2+\omega)}-\frac{1}{\omega(2+2 \omega)(1+\frac{1}{3}\pi^2 \omega^2)},
\end{equation*}
if $\omega$ is infinitely small, in which case the term $\frac{1}{3}\pi^2 \omega^2$ can be neglected. But 

\begin{equation*}
\frac{\omega}{\omega(2+ \omega)(2+2 \omega)}=\frac{1}{4}
\end{equation*}
for $\omega =0$ and  therefore  $\frac{1}{2}+\frac{1}{4}=\frac{3}{4}$ is the sum of the series

\begin{equation*}
\frac{1}{3}+\frac{1}{8}+\frac{1}{15}+\frac{1}{24}+\text{etc.}
\end{equation*}
as it is known from elsewhere.

\subsection*{Example 2}

\textit{To find the sum of the series}

\begin{equation*}
\frac{1}{1-xx}+\frac{1}{4-xx}+\frac{1}{9-xx}+\frac{1}{16-xx}+\text{etc.}
\end{equation*}
\textit{in the case, in which  one substitutes any integer $n$ for $x$ and omits the term $\frac{1}{nn-xx}$ of the series which would become infinite}.\\[2mm]
Therefore, this sum  in question will be expressed this way

\begin{equation*}
\frac{1}{2xx}-\frac{\pi}{2x \tan \pi x}-\frac{1}{nn-xx},
\end{equation*}
if one sets $x=n$, in which case the first term $\frac{1}{2xx}$ goes over into $\frac{1}{2nn}$, the two remaining ones on the other hand remain infinite. Therefore, put $x=n+\omega$, and because  $\tan (\pi n+\pi \omega)=\tan \pi \omega = \pi \omega$ for infinitely small $\omega$,  for the sum in question we will have

\begin{equation*}
\frac{1}{2nn}-\frac{\pi}{2(n+ \omega)\pi \omega}+\frac{1}{2n \omega +\omega \omega}
\end{equation*}
or 

\begin{equation*}
\frac{1}{2nn}-\frac{1}{\omega(2n+2 \omega)}+\frac{1}{\omega(2n+ \omega)}=\frac{1}{2nn}+\frac{1}{(2n+2 \omega)(2n+ \omega)},
\end{equation*}
whence, if it was $\omega =0$, the sum in question will result as

\begin{equation*}
=\frac{1}{2nn}+\frac{1}{4nn}=\frac{3}{4nn}.
\end{equation*}
Therefore, it will be

\begin{equation*}
\frac{3}{4nn}=\frac{1}{1-nn}+\frac{1}{4-nn}+\frac{1}{9-nn}+\cdots +\frac{1}{(n-1)^2-nn}
\end{equation*}
\begin{equation*}
+\frac{1}{(n+1)^2-nn}+\frac{1}{(n+2)^2-nn}+\text{etc. to infinity}
\end{equation*}
or the sum of this series will be infinite

\begin{equation*}
\frac{1}{(n+1)^2-nn}+\frac{1}{(n+2)^2-nn}+\frac{1}{(n+3)^2-nn}+\text{etc.}
\end{equation*}
\begin{equation*}
=\frac{3}{4nn}+\frac{1}{nn-1}+\frac{1}{nn-4}+\frac{1}{nn-9}+\cdots +\frac{1}{nn-(n-1)^2}.
\end{equation*}

\subsection*{Example 3}

\textit{To find the sum of this series}

\begin{equation*}
\frac{1}{1-xx}+\frac{1}{9-xx}+\frac{1}{25-xx}+\frac{1}{49-xx}+\text{etc.},
\end{equation*}
\textit{if one puts $x=1$ and the first term $\frac{1}{1-xx}$, which in this case becomes infinite, is subtracted.}\\[2mm]
Since the sum of this series in general is $=\frac{\pi \sin \frac{1}{2}\pi x}{4x \cos \frac{1}{2}\pi x}$, the sum in question will be

\begin{equation*}
=\frac{\pi \sin \frac{1}{2}\pi x}{4x \cos \frac{1}{2}\pi x}-\frac{1}{1-xx},
\end{equation*}
if one puts $x=1$. Since both of these terms become infinite, put $x=1- \omega$, and because

\begin{equation*}
\sin \left(\frac{1}{2}\pi -\frac{1}{2}\pi \omega\right)=\cos \frac{1}{2}\pi \omega =1- \frac{1}{8}\pi^2 \omega^2
\end{equation*}
and
\begin{equation*}
\cos \left(\frac{1}{2}\pi -\frac{1}{2}\pi \omega\right)= \sin \frac{1}{2}\pi \omega =\frac{1}{2}\pi \omega,
\end{equation*}
because of the infinitely small $\omega$, one will have this expression

\begin{equation*}
\frac{\pi(1-\frac{1}{8}\pi^2 \omega^2)}{4(1- \omega) \frac{1}{2}\pi \omega}-\frac{1}{2 \omega - \omega \omega}=\frac{1}{2 \omega(2- 2 \omega)}-\frac{1}{\omega(2- \omega)}
\end{equation*}
which becomes $=\frac{1}{4}$ for $\omega =0$, and therefore

\begin{equation*}
\frac{1}{4}=\frac{1}{8}+\frac{1}{24}+\frac{1}{48}+\frac{1}{80}+\frac{1}{120}+\text{etc.}
\end{equation*}

\subsection*{Example 4}

\textit{To find the sum of the series}

\begin{equation*}
\frac{1}{1-xx}+\frac{1}{9-xx}+\frac{1}{25-xx}+\frac{1}{49-xx}+\text{etc.},
\end{equation*}
\textit{if  one substitutes any odd integer $2n-1$ for $x$ and this term $\frac{1}{(2n-1)^2-xx}$, which  becomes infinite in this case, is subtracted.}\\[2mm]
Therefore, the sum  in question will be

\begin{equation*}
=\frac{\pi \sin \frac{1}{2}\pi x}{4x \cos \frac{1}{2}\pi x}- \frac{1}{(2n-1)^2-xx}
\end{equation*}
for $x=2n-1$. Therefore, let us set $x=2n-1 - \omega$, while $\omega$ is infinitely small, and it will be

\begin{equation*}
\sin \frac{1}{2}\pi x= \sin \left(\frac{2n-1}{2}\pi- \frac{1}{2}\pi \omega\right)= \pm \cos \frac{1}{2}\pi \omega,
\end{equation*}
where the upper sign holds, if $n$ is an odd number, the lower, if it is even. In like manner it will be

\begin{equation*}
\cos \frac{1}{2}\pi x= \cos \left(\frac{2n-1}{2}\pi-\frac{1}{2}\pi \omega\right)= \pm \sin \frac{1}{2}\pi \omega;
\end{equation*}
and hence, no matter whether $n$ is even or odd, it will be

\begin{equation*}
\frac{\sin \frac{1}{2}\pi x}{\cos \frac{1}{2}\pi x}= \frac{1}{\tan \frac{1}{2}\pi \omega}=\frac{1}{\frac{1}{2}\pi \omega}.
\end{equation*}
Therefore,  the sum in question will be expressed this way

\begin{equation*}
\frac{1}{2 \omega(2n-1- \omega)}-\frac{1}{\omega(2(2n-1)-\omega)}
\end{equation*}
and therefore will be $=\frac{1}{4(2n-1)^2}$. So, if $n=2$, it will be

\begin{equation*}
\frac{1}{36}=-\frac{1}{8}+\frac{1}{16}+\frac{1}{40}+\frac{1}{72}+\frac{1}{112}+\text{etc.}
\end{equation*}
the validity of which summation is known from elsewhere.


\title{On the Differentiation of inexplicable Functions
\footnote{
Original title: "'De Differentiatione Functionum inexplicabilum"', first published as part of the book \textit{Institutiones calculi differentialis cum eius usu in analysi finitorum ac doctrina serierum}, 1755, reprinted in \textit{Opera Omnia}: Series 1, Volume 10, pp. 588 - 618, Eneström-Number E212, translated by: Alexander Aycock for the Euler-Kreis Mainz}}
\author{Leonhard Euler}
\date{ }
\maketitle

\paragraph*{§367}

Here, I call those functions inexplicable which can not be defined either by determined expressions or by means of the roots of equations, such that those functions are not only not algebraic but it is also uncertain which kind of transcendental quantities they must be referred to. An inexplicable function of this kind is given by

\begin{equation*}
1+\frac{1}{2}+\frac{1}{3}+\cdots + \frac{1}{x},
\end{equation*}
which certainly depends on $x$, but, if $x$ is not an integer, can not be explained in any way. In like manner, this expression

\begin{equation*}
1 \cdot 2 \cdot 3 \cdot 4 \cdots x
\end{equation*}
will be an inexplicable function of $x$, since, if $x$ is any number, its value will not only be not algebraic, but even can not be expressed by means of a certain kind of transcendental quantities. In general, the notion of such inexplicable functions can be derived from series. For, let any series be propounded

\begin{alignat*}{9}
&1 &&2 &&3 &&4 \quad \cdots &&x \\
&A~+~ &&B~+~&&C~+~ &&D ~+~\cdots ~+~&&X,
\end{alignat*}
whose sum, if it can not be expressed by means of a finite formula, will yield an inexplicable function of $x$, namely,

\begin{equation*}
S=A+B+C+D+\cdots+X.
\end{equation*}
Similarly, infinite products of terms of series as

\begin{equation*}
P= A\cdot B \cdot C \cdot D \cdots X
\end{equation*}
will exhibit inexplicable functions of $x$, which, taking logarithms, can be reduced to the first form; for, it will be

\begin{equation*}
\log P= \log A+\log B +\log C +\log D+ \cdots +\log X.
\end{equation*}

\paragraph*{§368}

Therefore, I decided to a explain a method  to investigate the differentials of inexplicable functions of this kind in this chapter. This subject, although it seems to belong to the first part of the book, where  the rules of differential calculus where treated, nevertheless, since it requires a broader cognition of the doctrine of series one was only able to get to in this second part, let us, then forced to not follow the natural order, treat it here now. But because this investigation is completely new and has not been done by anybody until now, in order to discuss this part of differential calculus, it is only required that we rather try to sketch its first elements here. Furthermore, I will propound several questions whose answer requires the differentiation of inexplicable functions of this kind, by means of which at the same time the use of this treatment, which without any doubt will be a lot greater in  the future, is seen more clearly.

\paragraph*{§369}

To differentiate inexplicable functions of this kind it is especially necessary that we investigate their values  they have, if one substitutes $x+\omega$ for $x$. Therefore, let

\begin{alignat*}{9}
&  &&1 &&2 && 3 &&4 \quad \cdots \quad &&x \\
&S~=~&&A~+~&&B~+~&&C~+~&&D~+~ \cdots ~+~&&X
\end{alignat*}
and put $\Sigma$ for the value of $S$ it has, if  one substitutes $x+\omega$ for $x$, and let $Z$ be the term  corresponding to the index $x+\omega$ of the series. Now, denote the terms  corresponding to the indices $x+1$, $x+2$, $x+3$ etc. by $X^{\prime}$, $X^{\prime \prime}$, $X^{\prime \prime \prime}$ etc. and the one  corresponding to the infinite index $x+\infty$ by $X^{|\infty|}$. And in like manner indicate the terms corresponding to the indices $x+\omega +1$, $x+ \omega +2$, $x+\omega +3$ etc. by $Z^{\prime}$, $Z^{\prime \prime}$, $Z^{\prime \prime \prime}$ etc. and let $Z^{|\infty|}$ be the term corresponding to the index $x+\omega +\infty$. Having constituted all this, it will be

\begin{alignat*}{9}
&S^{\prime} &&=S+X^{\prime} \\[3mm]
&S^{\prime \prime}&&= S+X^{\prime}+X^{\prime \prime} \\[3mm]
&S^{\prime \prime \prime}&&=S+X^{\prime}+X^{\prime \prime}+X^{\prime \prime \prime} \\[3mm]
& && \text{etc.} \\
&S^{|\infty |}&&=S+X^{\prime}+X^{\prime \prime}+X^{\prime \prime \prime}+\cdots +X^{|\infty|} \\
\end{alignat*}
Since in like manner  $\Sigma$ is also successively increased by the terms $Z^{\prime}$, $Z^{\prime \prime}$ etc., it will be

\begin{alignat*}{9}
&\Sigma^{\prime} &&=\Sigma + Z^{\prime} \\[3mm]
&\Sigma^{\prime \prime}&&= \Sigma +Z^{\prime}+Z^{\prime \prime} \\[3mm]
&\Sigma^{\prime \prime \prime}&&=\Sigma + Z^{\prime}+Z^{\prime \prime}+Z^{\prime \prime \prime} \\[3mm]
& && \text{etc.} \\
&\Sigma^{|\infty |}&&=\Sigma +Z^{\prime}+Z^{\prime \prime}+Z^{\prime \prime \prime}+\cdots +Z^{|\infty|} \\
\end{alignat*}

\paragraph*{§370}

Now, the nature of the series $S$, $S^{\prime}$, $S^{\prime \prime}$, $S^{\prime \prime \prime}$ etc. it will have, if continued to infinity, is to be considered; if the series is confounded with an arithmetic progression at infinity, what happens, if  the terms of the series $X$, $X^{\prime}$, $X^{\prime \prime}$, $X^{\prime \prime \prime}$ etc. converge to a ratio of $1$  such that the differences of the series $S$, $S^{\prime}$, $S^{\prime \prime}$ etc. become equal eventually, in this case the quantities $S^{|\infty|}$, $S^{|\infty +1|}$, $S^{|\infty +2|}$ etc. will be terms of an arithmetic progression, and because $\Sigma^{|\infty|}= S^{|\infty +\omega|}$, because of

\begin{equation*}
S^{|\infty+ \omega|}= S^{|\infty|}+\omega(S^{|\infty +1|}-S^{|\infty|})= \omega S^{|\infty +1|}+(1- \omega)S^{|\infty|},
\end{equation*}
it will be

\begin{equation*}
\Sigma^{|\infty|}= \omega S^{|\infty +1|}+(1-\omega)S^{|\infty|}.
\end{equation*}
But  $S^{|\infty +1|}= S^{|\infty|}+X^{|\infty +1|}$, whence 

\begin{equation*}
\Sigma^{|\infty|}=S^{|\infty|}+\omega X^{|\infty +1|},
\end{equation*}
from which one will obtain this equation

\begin{equation*}
\Sigma +Z^{\prime}+Z^{\prime \prime}+Z^{\prime \prime \prime}+\cdots + Z^{|\infty|}
\end{equation*}
\begin{equation*}
=S+X^{\prime}+X^{\prime \prime}+X^{\prime \prime \prime} + \cdots + X^{|\infty|}+\omega X^{|\infty +1|},
\end{equation*}
from which the value in question $\Sigma$ the functions $S$ has, if in it $x+ \omega$ is substituted for $x$, will be

\begin{alignat*}{9}
&\Sigma =S+ \omega X^{|\infty +1|} &&+X^{\prime}&&+X^{\prime \prime} &&+X^{\prime \prime \prime}&&+\text{etc.} \quad \text{to infinity} \\
& &&-Z^{\prime} &&-Z^{\prime \prime} &&-Z^{\prime \prime \prime} &&- \text{etc.} \quad \text{to infinity} 
\end{alignat*}
Hence, if the infinitesimal terms of the series $A$, $B$, $C$, $D$ etc. vanish, the term $\omega X^{|\infty +1|}$ vanishes and can be omitted.

\paragraph*{§371}

Therefore, the value of $\Sigma$ is expressed by means of new infinite series which can be exhibited, if one knows the general term of the series $A+B+C+\text{etc.}$, from which then the values of the terms $Z^{\prime}$, $Z^{\prime \prime}$, $Z^{\prime \prime \prime}$ etc. can be defined. Therefore, having put $\omega$ to be infinitely small, since $\Sigma - S$ is the differential of the function $S$, this differential $dS$ will be expressed by  an infinite series. And if not even the higher powers of $\omega$ are neglected, one will have the complete differential of this inexplicable function $S$; to show its nature quite plainly, we will consider  the following examples.

\subsection*{Example 1}

\textit{To find the differential of this inexplicable function}

\begin{equation*}
S=1+\frac{1}{2}+\frac{1}{3}+\frac{1}{4}+\cdots +\frac{1}{x}.
\end{equation*} 
Since the general term $X$ of this series is $=\frac{1}{x}$ and therefore

\begin{equation*}
\renewcommand{\arraystretch}{2,0}
\setlength{\arraycolsep}{0.0mm}
\begin{array}{lll|lllllllllllll}
X^{\prime}&~=~& \dfrac{1}{x+1} \quad & \quad Z^{\prime} &~=~&\dfrac{1}{x+1+\omega} \\
X^{\prime  \prime}&~=~& \dfrac{1}{x+2} \quad & \quad Z^{\prime \prime} &~=~&\dfrac{1}{x+2+\omega} \\
X^{\prime \prime \prime}&~=~& \dfrac{1}{x+3} \quad & \quad Z^{\prime \prime \prime} &~=~&\dfrac{1}{x+3+\omega} \\
   &\text{etc.}&   & & \text{etc.},
\end{array}
\end{equation*} 
because of

\begin{equation*}
X^{|\infty +1|}=\frac{1}{x+\infty +1}=0,
\end{equation*}
if  one writes $x+\omega$ instead of $x$, the function $S$ will go over into $\Sigma$ that 

\begin{alignat*}{9}
&\Sigma =S &&+\frac{\quad ~~~ 1\quad ~~~}{x+1}&&+\frac{\quad ~~~1 \quad ~~~}{x+2}&&+\frac{\quad ~~~1  \quad ~~~}{x+3}&&+\text{etc.} \\[3mm]
& &&-\frac{1}{x+1+ \omega}&&-\frac{1}{x+2+ \omega}&&- \frac{1}{x+3+ \omega}&&-\text{etc.},
\end{alignat*}
or, combining each to terms into  one, it will be

\begin{equation*}
\Sigma =S +\frac{\omega}{(x+1)(x+1+ \omega)}+\frac{\omega}{(x+2)(x+2 \omega)}+\frac{\omega}{(x+3)(x+3 + \omega)}+\text{etc.},
\end{equation*}
or, because 

\begin{alignat*}{9}
&\frac{1}{x+1+\omega}&&=\frac{1}{x+1}&&-\frac{\omega}{(x+1)^2}&&+\frac{\omega^2}{(x+1)^3}&&-\frac{\omega^3}{(x+1)^4}&&+\text{etc.} \\[3mm]
&\frac{1}{x+2+\omega}&&=\frac{1}{x+2}&&-\frac{\omega}{(x+2)^2}&&+\frac{\omega^2}{(x+2)^3}&&-\frac{\omega^3}{(x+2)^4}&&+\text{etc.} \\
& && && &&\text{etc.},
\end{alignat*}
having ordered the series according to powers of $\omega$, it will be

\begin{alignat*}{9}
&\Sigma = S &&+\omega &&\bigg(\frac{1}{(x+1)^2}&&+\frac{1}{(x+2)^2}&&+\frac{1}{(x+3)^2}&&+\frac{1}{(x+4)^2}&&+\text{etc.}\bigg) \\[3mm]
&  &&-\omega^2 &&\bigg(\frac{1}{(x+1)^3}&&+\frac{1}{(x+2)^3}&&+\frac{1}{(x+3)^3}&&+\frac{1}{(x+4)^3}&&+\text{etc.}\bigg) \\[3mm]
&  &&+\omega^3 &&\bigg(\frac{1}{(x+1)^4}&&+\frac{1}{(x+2)^4}&&+\frac{1}{(x+3)^4}&&+\frac{1}{(x+4)^4}&&+\text{etc.}\bigg) \\[3mm]
&  &&-\omega^4 &&\bigg(\frac{1}{(x+1)^5}&&+\frac{1}{(x+2)^5}&&+\frac{1}{(x+3)^5}&&+\frac{1}{(x+4)^5}&&+\text{etc.}\bigg) \\
& && && && && \text{etc.}
\end{alignat*}
Having taken $dx$ for $\omega$, we will obtain the complete differential of the propounded function $S$

\begin{alignat*}{9}
&dS &&=dx&&\bigg(\frac{1}{(x+1)^2}&&+\frac{1}{(x+2)^2}&&+\frac{1}{(x+3)^2}&&+\frac{1}{(x+4)^2}&&+\text{etc.}\bigg) \\[3mm]
&  &&-dx^2 &&\bigg(\frac{1}{(x+1)^3}&&+\frac{1}{(x+2)^3}&&+\frac{1}{(x+3)^3}&&+\frac{1}{(x+4)^3}&&+\text{etc.}\bigg) \\[3mm]
&  &&+dx^3 &&\bigg(\frac{1}{(x+1)^4}&&+\frac{1}{(x+2)^4}&&+\frac{1}{(x+3)^4}&&+\frac{1}{(x+4)^4}&&+\text{etc.}\bigg) \\[3mm]
&  &&-dx^4 &&\bigg(\frac{1}{(x+1)^5}&&+\frac{1}{(x+2)^5}&&+\frac{1}{(x+3)^5}&&+\frac{1}{(x+4)^5}&&+\text{etc.}\bigg) \\
& && && && && \text{etc.}
\end{alignat*}

\subsection*{Example 2}

\textit{To find the differential of this inexplicable function of $x$}

\begin{equation*}
S=1+\frac{1}{3}+\frac{1}{5}+\frac{1}{7}+\cdots +\frac{1}{2x-1}.
\end{equation*}
Since the general term of this series is $X=\frac{1}{2x-1}$, it will be

\begin{equation*}
\renewcommand{\arraystretch}{2,0}
\setlength{\arraycolsep}{0.0mm}
\begin{array}{lll|lllllllllllll}
X^{\prime}&~=~& \dfrac{1}{2x+1} \quad & \quad Z^{\prime} &~=~&\dfrac{1}{2x+1+\omega} \\
X^{\prime  \prime}&~=~& \dfrac{1}{2x+3} \quad & \quad Z^{\prime \prime} &~=~&\dfrac{1}{2x+3+\omega} \\
X^{\prime \prime \prime}&~=~& \dfrac{1}{2x+5} \quad & \quad Z^{\prime \prime \prime} &~=~&\dfrac{1}{2x+5+\omega} \\
   &\text{etc.}&   & & \text{etc.},
\end{array}
\end{equation*} 
Because of the vanishing and equal infinitesimal terms of this series, the value of $S$, if  one writes $x+\omega$ instead of $x$, will result as

\begin{alignat*}{9}
&\Sigma  = S &&+ \frac{ \quad \quad ~ 1 \quad \quad ~}{2x+1}&&+\frac{\quad \quad ~1 \quad \quad ~}{2x+3}&&+\frac{\quad \quad ~1 \quad \quad ~}{2x+5}&&+\text{etc.} \\[3mm]
&            &&-\frac{1}{2x+1+2 \omega}&&-\frac{1}{2x+3+2 \omega}&&-\frac{1}{2x+5+2 \omega}&&-\text{etc.} 
\end{alignat*} 
or

\begin{equation*}
\Sigma = S +\frac{2 \omega}{(2x+1)(2x+1+2 \omega)}+\frac{2 \omega}{(2x+3)(2x+3+2 \omega)}+\text{etc.}
\end{equation*}
But if each term is expanded into a power series in $\omega$, it will be

\begin{alignat*}{9}
& \Sigma =S &&+~~2 \omega &&\bigg(\frac{1}{(2x+1)^2}&&+\frac{1}{(2x+3)^2}&&+\frac{1}{(2x+5)^2}&&+\text{etc.}\bigg)\\[3mm]
&   &&-~4 \omega^2 &&\bigg(\frac{1}{(2x+1)^3}&&+\frac{1}{(2x+3)^3}&&+\frac{1}{(2x+5)^3}&&+\text{etc.}\bigg)\\[3mm]
&  &&+~8 \omega^3 &&\bigg(\frac{1}{(2x+1)^4}&&+\frac{1}{(2x+3)^4}&&+\frac{1}{(2x+5)^4}&&+\text{etc.}\bigg)\\[3mm]
&  &&-16 \omega^4 &&\bigg(\frac{1}{(2x+1)^4}&&+\frac{1}{(2x+3)^4}&&+\frac{1}{(2x+5)^4}&&+\text{etc.}\bigg)\\
& && && && &&\text{etc.}
\end{alignat*}
Now write $dx$ for $\omega$ and the complete differential of the propounded inexplicable function $S$ will be

\begin{alignat*}{9}
&dS &&=~2 dx &&\bigg(\frac{1}{(2x+1)^2}&&+\frac{1}{(2x+3)^2}&&+\frac{1}{(2x+5)^2}&&+\text{etc.}\bigg)\\[3mm]
&   &&-~4 dx^2 &&\bigg(\frac{1}{(2x+1)^3}&&+\frac{1}{(2x+3)^3}&&+\frac{1}{(2x+5)^3}&&+\text{etc.}\bigg)\\[3mm]
&  &&+~8 dx^3 &&\bigg(\frac{1}{(2x+1)^4}&&+\frac{1}{(2x+3)^4}&&+\frac{1}{(2x+5)^4}&&+\text{etc.}\bigg)\\[3mm]
&  &&-16 dx^4 &&\bigg(\frac{1}{(2x+1)^4}&&+\frac{1}{(2x+3)^4}&&+\frac{1}{(2x+5)^4}&&+\text{etc.}\bigg)\\
& && && && &&\text{etc.}
\end{alignat*}

\subsection*{Example 3}

\textit{To find the complete differential of this inexplicable function of $x$}

\begin{equation*}
S=1+\frac{1}{2^n}+\frac{1}{3^n}+\frac{1}{4^n}+\cdots +\frac{1}{x^n}.
\end{equation*}
Since the general term of this series is $=\frac{1}{x^n}$, the infinitesimal terms will  vanish and be equal to each other. And hence, because of 

\begin{equation*}
\renewcommand{\arraystretch}{2,0}
\setlength{\arraycolsep}{0.0mm}
\begin{array}{lll|lllllllllllll}
X^{\prime}&~=~& \dfrac{1}{(x+1)^n} \quad & \quad Z^{\prime} &~=~&\dfrac{1}{(x+1+\omega)^n} \\
X^{\prime  \prime}&~=~& \dfrac{1}{(x+2)^n} \quad & \quad Z^{\prime \prime} &~=~&\dfrac{1}{(x+2+\omega)^n} \\
X^{\prime \prime \prime}&~=~& \dfrac{1}{(x+3)^n} \quad & \quad Z^{\prime \prime \prime} &~=~&\dfrac{1}{(x+3+\omega)^n} \\
   &\text{etc.}&   & & \text{etc.},
\end{array}
\end{equation*}
it will be

\begin{alignat*}{9}
&X^{\prime} &&-Z^{\prime}&&=\frac{n \omega}{(x+1)^{n+1}}&&-\frac{n(n+1)\omega^2}{2(x+1)^{n+2}}&&+\frac{n(n+1)(n+2)\omega^3}{6(x+1)^{n+3}}&&-\text{etc.} \\[3mm]
&X^{\prime \prime} &&-Z^{\prime \prime}&&=\frac{n \omega}{(x+2)^{n+1}}&&-\frac{n(n+1)\omega^2}{2(x+2)^{n+2}}&&+\frac{n(n+1)(n+2)\omega^3}{6(x+2)^{n+3}}&&-\text{etc.} \\
& && && && \text{etc.},
\end{alignat*}
from which one finds

\begin{alignat*}{9}
 \Sigma -S=n \omega && \left(\frac{1}{(x+1)^{n+1}}+\frac{1}{(x+2)^{n+1}}+\frac{1}{(x+3)^{n+1}}+ \text{etc.}\right)\\[3mm]
  -\frac{n(n+1)}{1 \cdot 2}\omega^2 && \left(\frac{1}{(x+1)^{n+2}}+\frac{1}{(x+2)^{n+2}}+\frac{1}{(x+3)^{n+2}}+ \text{etc.}\right)\\[3mm]
   +\frac{n(n+1)(n+2)}{1 \cdot 2 \cdot 3}\omega^3 && \left(\frac{1}{(x+1)^{n+3}}+\frac{1}{(x+2)^{n+3}}+\frac{1}{(x+3)^{n+3}}+ \text{etc.}\right)\\
   \text{etc.}&&
\end{alignat*}
Therefore, having put $\omega =dx$, the complete differential in question of the function $S$ will be

\begin{alignat*}{9}
 dS =n dx && \left(\frac{1}{(x+1)^{n+1}}+\frac{1}{(x+2)^{n+1}}+\frac{1}{(x+3)^{n+1}}+ \text{etc.}\right)\\[3mm]
  -\frac{n(n+1)}{1 \cdot 2}dx^2 && \left(\frac{1}{(x+1)^{n+2}}+\frac{1}{(x+2)^{n+2}}+\frac{1}{(x+3)^{n+2}}+ \text{etc.}\right)\\[3mm]
   +\frac{n(n+1)(n+2)}{1 \cdot 2 \cdot 3}dx^3 && \left(\frac{1}{(x+1)^{n+3}}+\frac{1}{(x+2)^{n+3}}+\frac{1}{(x+3)^{n+3}}+ \text{etc.}\right)\\
   \text{etc.}&&
\end{alignat*}

\paragraph*{§372}

From these also the sums of these series can be interpolated or the values of the summatory terms can be exhibited, if the number of terms is not an integer number. For, if one puts $x=0$, it will also be $S=0$ and $\Sigma$ will express the sum of as many terms as the number $\omega$ contains units, even though this number $\omega$ is not an integer. So, if in the first example one puts

\begin{equation*}
\Sigma =1+ \frac{1}{2}+\frac{1}{3}+\cdots +\frac{1}{\omega},
\end{equation*}
it will be

\begin{equation*}
\Sigma = \frac{\omega}{1(1+\omega)}+\frac{\omega}{2(2+\omega)}+\frac{\omega}{3(3+\omega)}+\frac{\omega}{4(4+\omega)}+\text{etc.}
\end{equation*}
or 

\begin{alignat*}{9}
&\sum &&= \omega &&\bigg(1&&+~\frac{1}{4}&&+~\frac{1}{9}&&+\frac{1}{16}&&+\frac{1}{25}&&+\text{etc.}\bigg)\\[3mm]
&  &&- \omega^2 &&\bigg(1&&+\frac{1}{2^3}&&+\frac{1}{3^3}&&+~\frac{1}{4^3}&&+~\frac{1}{5^2}&&+\text{etc.}\bigg)\\[3mm]
&  &&  -\omega^3 &&\bigg(1&&+\frac{1}{2^4}&&+\frac{1}{3^4}&&+~\frac{1}{4^4}&&+~\frac{1}{5^4}&&+\text{etc.}\bigg)\\
& && && && && && \text{etc.}
\end{alignat*}
In the third example on the other hand it will be

\begin{equation*}
\Sigma =1+\frac{1}{2^n}+\frac{1}{3^n}+\frac{1}{4^n}+\cdots +\frac{1}{\omega^n}.
\end{equation*}
The value of $\Sigma$, whether $\omega$ is an integer number or a fractional number, will be expressed by  the following series

\begin{alignat*}{9}
\Sigma =n \omega &&\left(1+\frac{1}{2^{n+1}}+\frac{1}{3^{n+1}}+\frac{1}{4^{n+1}}+\text{etc.}\right)\\[3mm]
-\frac{n(n+1)}{1 \cdot 2}\omega^2 &&\left(1+\frac{1}{2^{n+2}}+\frac{1}{3^{n+2}}+\frac{1}{4^{n+2}}+\text{etc.}\right) \\[3mm]
+\frac{n(n+1)(n+2)}{1 \cdot 2 \cdot 3}\omega^3 &&\left(1+\frac{1}{2^{n+3}}+\frac{1}{3^{n+3}}+\frac{1}{4^{n+3}}+\text{etc.}\right)\\
\text{etc.}
\end{alignat*}

\paragraph*{§373}

These same considerations can also be applied to a general series; for, because

\begin{alignat*}{9}
& &&1 &&2 && 3 && 4 \quad \cdots \cdots \quad && x \\
&S~=~ && A~+~&&B~+~&& C~+~ && D~+~\cdots ~+~ &&X 
\end{alignat*}
and, having written $x+\omega$ instead of $x$, $X$ goes over into $Z$ and $S$ into $\Sigma$, it will be

\begin{equation*}
Z=X+\frac{\omega dX}{dx}+\frac{\omega^2 ddX}{1 \cdot 2 dx^2}+\frac{\omega^3 d^3 X}{1 \cdot 2 \cdot 3 dx^3}+\text{etc.},
\end{equation*}
and since in like manner $Z^{\prime}$, $Z^{\prime \prime}$, $Z^{\prime \prime \prime}$ etc. are expressed by  $X^{\prime}$, $X^{\prime \prime}$, $X^{\prime \prime \prime}$ etc., it will be

\begin{alignat*}{9}
\Sigma =S +\omega X^{|\infty +1|}-\frac{\omega}{dx}d.&&(X^{\prime}+X^{\prime \prime}+X^{\prime \prime \prime}+X^{\prime \prime \prime \prime}+\text{etc.})\\[3mm]
-\frac{\omega^2}{1 \cdot 2 dx^2}dd.&&(X^{\prime}+X^{\prime \prime}+X^{\prime \prime \prime}+X^{\prime \prime \prime \prime}+\text{etc.})\\[3mm]
-\frac{\omega^3}{1 \cdot 2 \cdot 3 dx^3}d^3. &&(X^{\prime}+X^{\prime \prime}+X^{\prime \prime \prime}+X^{\prime \prime \prime \prime}+\text{etc.}) \\
\text{etc.},
\end{alignat*}
and if $X^{|\infty +1|}$ is not $=0$, it can be expressed in this way that the consideration of infinity is avoided

\begin{equation*}
X^{|\infty +1|}=X^{\prime}+(X^{\prime \prime}-X^{\prime})+(X^{\prime \prime \prime}-X^{\prime \prime})+(X^{\prime \prime \prime \prime}-X^{\prime \prime \prime})+\text{etc.}
\end{equation*}
and therefore it will be

\begin{alignat*}{9}
&\Sigma = S+\omega X^{\prime}&&+\omega ((X^{\prime \prime}-X^{\prime})+(X^{\prime \prime \prime}-X^{\prime \prime})+(X^{\prime \prime \prime \prime}-X^{\prime \prime \prime})+\text{etc.})\\[3mm]
& &&-\frac{\omega}{dx}d. (X^{\prime}+X^{\prime \prime}+X^{\prime \prime \prime}+X^{\prime \prime \prime \prime}+\text{etc.})\\[3mm]
& &&-\frac{\omega^2}{2 dx^2}dd. (X^{\prime}+X^{\prime \prime}+X^{\prime \prime \prime}+X^{\prime \prime \prime \prime}+\text{etc.}) \\[3mm]
& &&-\frac{\omega^3}{6dx^3}d^3. (X^{\prime}+X^{\prime \prime}+X^{\prime \prime \prime}+X^{\prime \prime \prime \prime}+\text{etc.})\\
& && \text{etc.}
\end{alignat*}
If one puts $\omega =dx$, the following differential of

\begin{equation*}
S=A+B+C+\cdots+X
\end{equation*}
will result expressed this way

\begin{alignat*}{9}
dS =X^{\prime}dx+dx&((X^{\prime \prime}-X^{\prime})+(X^{\prime \prime \prime}-X^{\prime \prime})+(X^{\prime \prime \prime \prime}-X^{\prime \prime \prime})+\text{etc.})\\[3mm]
-d. &(X^{\prime}+X^{\prime \prime}+X^{\prime \prime \prime}+X^{\prime \prime \prime \prime}+\text{etc.})\\[3mm]
-\frac{1}{2}dd. &(X^{\prime}+X^{\prime \prime}+X^{\prime \prime \prime}+X^{\prime \prime \prime \prime}+\text{etc.}) \\[3mm]
-\frac{1}{6}d^3. &(X^{\prime}+X^{\prime \prime}+X^{\prime \prime \prime}+X^{\prime \prime \prime \prime}+\text{etc.})\\
\text{etc.}&
\end{alignat*}

\paragraph*{§374}

Let us put that  $x=0$; it will be

\begin{equation*}
X^{\prime}=A, \quad X^{\prime \prime}=B \quad \text{etc.}
\end{equation*}
and hence $X^{\prime}+X^{\prime \prime}+X^{\prime \prime \prime}+\text{etc.}$ will be an infinite series whose general term is $=X$. Further, form the series from these general terms

\begin{equation*}
\frac{dX}{	dx}, \quad \frac{ddX}{2dx^2}, \quad \frac{d^3X}{6dx^3}, \quad \frac{d^4 X}{24 dx^4} \quad \text{etc.}
\end{equation*}
the sum of which series, if they are continued to infinity, we want to be

\begin{equation*}
\mathcal{S} X =\mathfrak{A}, \quad \mathcal{S}  \frac{dX}{dx}=\mathfrak{B}, \quad \mathcal{S} \frac{ddX}{2dx^2}=\mathfrak{C}, \quad \mathcal{S} \frac{d^3 X}{6dx^3}=\mathfrak{D} \quad \text{etc.};
\end{equation*}
and since for $x=0$ also $S=0$,  $\Sigma$ will be the sum of the series

\begin{equation*}
A+B+C+D+\cdots +Z
\end{equation*}
containing $\omega$ terms; for, $Z$ is the term of the index $\omega$, no matter whether $\omega$ is an integer number or a fraction. Therefore, one will have

\begin{equation*}
\Sigma = \omega A+ \omega ((B-A)+(C-B)+(D-C)+\text{etc.})
\end{equation*}
\begin{equation*}
- \omega \mathfrak{B}-\omega^2 \mathfrak{C}-\omega^3 \mathfrak{D}- \omega^4 \mathfrak{E}-\text{etc.},
\end{equation*}
where the first series can be omitted, if the terms of the propounded series vanish eventually.

\paragraph*{§375}

Now, let us write $x$ instead of $\omega$ and $\Sigma$ will go over into $S$ such that 

\begin{alignat*}{9}
& &&1 &&2 && 3 && 4 \quad \cdots \cdots \quad && x \\
&S~=~ && A~+~&&B~+~&& C~+~ && D~+~\cdots ~+~ &&X 
\end{alignat*}
and the same value of $S$ will be expressed as an infinite series this way

\begin{equation*}
S=Ax+x((B-A)+(C-B)+(D-C)+\text{etc.})
\end{equation*}
\begin{equation*}
-\mathfrak{B}x-\mathfrak{C}x^2-\mathfrak{D}x^3-\mathfrak{E}x^4-\mathfrak{F}x^5-\text{etc.};
\end{equation*}
since its value is expressed equally distinctly, no matter whether $x$ is an integer number or a fraction, one is able to express the value of $S$ of any order from this easily as:

\begin{equation*}
\frac{dS}{dx}=A+(B-A)+(C-B)+(D-C)+\text{etc.}
\end{equation*}
\begin{alignat*}{9}
& && &&-\mathfrak{B}&&-2 \mathfrak{C}x&&- 3 \mathfrak{D}x^2&&-4 \mathfrak{E}x^3 &&-\text{etc.} \\[3mm]
&\frac{ddS}{2dx^2}&&~=~&&-\mathfrak{C}&&-3 \mathfrak{D}x&&- 6 \mathfrak{E}x^2 &&-10 \mathfrak{F}x^3 &&-\text{etc.} \\[3mm]
&\frac{d^3S}{6dx^3}&&~=~&&-\mathfrak{D}&&-4\mathfrak{E}x&&-10 \mathfrak{F}x^2 &&-20 \mathfrak{G}x^3 &&-\text{etc.} \\[3mm]
&\frac{d^4 S}{24 dx^4}&&~=~&&-\mathfrak{E}&&-5 \mathfrak{F}x&&-15 \mathfrak{G}x^2 &&-35 \mathfrak{H}x^3&&-\text{etc.} 
\end{alignat*}
Therefore, since the complete differential is

\begin{equation*}
=dS+\frac{1}{2}ddS+\frac{1}{6}d^3S+\frac{1}{24}d^4 S+\text{etc.},
\end{equation*}
the complete differential of the propounded function $S$ will be

\begin{equation*}
dS= Adx+(B-A)dx+(C-B)dx+(D-C)dx+\text{etc.}
\end{equation*} 
\begin{equation*}
-\mathfrak{B}dx-\mathfrak{C}(2xdx+dx^2)-\mathfrak{D}(3x^2dx+3xdx^2+dx^3)
\end{equation*}
\begin{equation*}
-\mathfrak{E}(4x^3dx+6x^2 dx^2+4xdx^3+dx^4)-\text{etc.}
\end{equation*}

\paragraph*{§376}

Therefore, this way the complete differential of any inexplicable function $S$ can be assigned, if the infinitesimal terms of the series

\begin{equation*}
A+B+C+D+\text{etc.}
\end{equation*}
either vanish or become equal to each other. For, if the infinitesimal terms of this series were not $=0$, then the sum of the series $\mathfrak{B}$ which is formed from the general term $\frac{dX}{dx}$, will become infinite, but together with the series

\begin{equation*}
A+(B-A)+(C-B)+(D-C)+\text{etc.}
\end{equation*}
it will constitute a finite sum. But it can happen that the terms of the series $A+B+C+D+$etc. are increased to infinity in such a way that not only the sum of the series $\mathfrak{B}$, but also the sum of the series $\mathfrak{C}$ becomes infinitely large, in which case it does not suffice to have added the series $A+(B-A)+(C-B)+$etc.; but since in this case the infinitesimal values considered in § 370, namely $S^{|\infty|}$, $S^{|\infty +1|}$, $S^{|\infty +2|}$, are not any longer terms in an arithmetic progression, as we had assumed before, the nature of this progression will  have to be taken into account. As we assumed the first differences of these progressions to be equal, so we will extend the method even further, if we set that just the second or the third or the higher differences of these values become constant.

\paragraph*{§377}

Arguing exactly as before in  § 369, let us put that just the second differences of the mentioned values are constant

\begin{equation*}
S^{|\infty|}, \quad S^{|\infty +1|}, \quad S^{|\infty +2|}
\end{equation*}
\begin{equation*}
\textit{First Differences}
\end{equation*}
\begin{equation*}
X^{|\infty +1|}, \quad X^{|\infty +2|}
\end{equation*}
\begin{equation*}
\textit{Second Differences}
\end{equation*}
\begin{equation*}
X^{|\infty +2|} -X^{|\infty +1|}
\end{equation*}
Therefore, it will be

\begin{equation*}
\Sigma^{|\infty|}=S^{|\infty +\omega|}=S^{|\infty|}+\omega X^{|\infty +1|}+\frac{\omega (\omega -1)}{1 \cdot 2}(X^{|\infty +2|}-X^{|\infty +1|})
\end{equation*}
\begin{equation*}
=S^{|\infty|}-\frac{\omega(\omega -3)}{1 \cdot 2}X^{|\infty +1|}+\frac{\omega(\omega -1)}{1 \cdot 2}X^{|\infty +2|}.
\end{equation*}
Therefore, one will have this equation

\begin{equation*}
\Sigma +Z^{\prime}+Z^{\prime \prime}+Z^{\prime \prime \prime}+ \cdots +Z^{|\infty|}
\end{equation*}
\begin{equation*}
=S+X^{\prime}+X^{\prime \prime}+X^{\prime \prime \prime}+\cdots +X^{|\infty|}-\frac{\omega(\omega-3)}{1 \cdot 2}X^{|\infty +1|}+\frac{\omega (\omega -1)}{1 \cdot 2}X^{|\infty +2|},
\end{equation*}
from which one finds

\begin{alignat*}{9}
&\Sigma =S &&+X^{\prime}+X^{\prime \prime}+X^{\prime \prime \prime} +X^{\prime \prime \prime \prime}+\text{etc. to infinity} \\[3mm]
&  &&-Z^{\prime}-Z^{\prime \prime}-Z^{\prime \prime \prime} -Z^{\prime \prime \prime \prime}-\text{etc. to infinity} \\[3mm]
& &&+\omega X^{|\infty +1|}+\frac{\omega (\omega -1)}{1 \cdot 2}(X^{|\infty +2|}-X^{|\infty +1|}).
\end{alignat*}
But these infinitesimal terms can be represented in such a way that

\begin{alignat*}{9}
&\Sigma =S &&+X^{\prime}&&+X^{\prime \prime}&&+X^{\prime \prime \prime} &&+X^{\prime \prime \prime \prime}&&+\text{etc. to infinity} \\[3mm]
&  &&-Z^{\prime}&&-Z^{\prime \prime}&&-Z^{\prime \prime \prime}&& -Z^{\prime \prime \prime \prime}&&-\text{etc. to infinity} 
\end{alignat*}
\begin{equation*}
+\omega X^{\prime} +\omega \left\lbrace
\renewcommand{\arraystretch}{1,7}
\setlength{\arraycolsep}{0.0mm}
\begin{array}{llllllllllllllll}
~+~&X^{\prime \prime}&~+~&X^{\prime \prime \prime}&~+~& X^{\prime \prime \prime \prime}&~+~&X^{\prime \prime \prime \prime \prime}&~+~&\text{etc.} \\
~-~&X^{\prime }&~-~&X^{\prime \prime }&~-~& X^{\prime  \prime \prime}&~+~&X^{\prime \prime \prime \prime}&~-~&\text{etc.}
\end{array}
\right\rbrace
\end{equation*}
whence the law  describing the nature of this expression, if just the third or fourth or higher differences were constant, is obvious.

\paragraph*{§378}

Because, as we demonstrated above,

\begin{equation*}
Z=X+\frac{\omega dX}{1 dx}+\frac{\omega^2 ddX}{1 \cdot 2 dx^2}+\frac{\omega^3 d^3 X}{1 \cdot 2 \cdot 3 dx^3}+\text{etc.},
\end{equation*}
if for  $Z^{\prime}$, $Z^{\prime \prime}$, $Z^{\prime \prime \prime}$ etc. we substitute the values to result from them, the value of $S$, if  one writes $x+\omega$ instead of $x$, will be the following:

\begin{equation*}
\Sigma=S+\omega X^{\prime} +\omega \left\lbrace
\renewcommand{\arraystretch}{1,7}
\setlength{\arraycolsep}{0.0mm}
\begin{array}{llllllllllllllll}
~+~&X^{\prime \prime}&~+~&X^{\prime \prime \prime}&~+~& X^{\prime \prime \prime \prime}&~+~&X^{\prime \prime \prime \prime \prime}&~+~&\text{etc.} \\
~-~&X^{\prime }&~-~&X^{\prime \prime }&~-~& X^{\prime  \prime \prime}&~+~&X^{\prime \prime \prime \prime}&~-~&\text{etc.}
\end{array}
\right\rbrace
\end{equation*}
\begin{equation*}
\begin{aligned}
&+\frac{\omega (\omega -1)}{1 \cdot 2}X^{\prime \prime} \\
&-\frac{\omega(\omega -1)}{1 \cdot 2}X^{\prime}
\end{aligned}
+\frac{\omega(\omega -1)}{1 \cdot 2}\left\lbrace
\renewcommand{\arraystretch}{1,7}
\setlength{\arraycolsep}{0.0mm}
\begin{array}{lrlrlrlrlrllllll}
~+~&X^{\prime \prime \prime}&~+~&X^{\prime \prime \prime \prime}&~+~&X^{\prime \prime \prime \prime \prime}&~+~&\text{etc.} \\
~-~&2X^{\prime \prime}&~-~&2X^{\prime \prime \prime}&~-~&2X^{\prime \prime \prime \prime }&~-~&\text{etc.} \\
~+~&X^{\prime }&~+~&X^{\prime \prime }&~+~&X^{\prime \prime \prime }&~+~&\text{etc.} \\
\end{array}
\right\rbrace
\end{equation*}
\begin{alignat*}{19}
&-~\frac{\omega}{dx}&&d. &&(X^{\prime}+X^{\prime \prime}&&+X^{\prime \prime \prime} &&+X^{\prime \prime \prime \prime}&&+\text{etc.}) \\[3mm]
&-\frac{\omega^2}{2dx^2}&&d^2. &&(X^{\prime}+X^{\prime \prime}&&+X^{\prime \prime \prime} &&+X^{\prime \prime \prime \prime}&&+\text{etc.}) \\[3mm]
&-\frac{\omega^3}{6dx^3}&&d^3. &&(X^{\prime}+X^{\prime \prime}&&+X^{\prime \prime \prime} &&+X^{\prime \prime \prime \prime}&&+\text{etc.}) \\
&  && && \text{etc.}
\end{alignat*}
If  one writes $dx$ instead of $\omega$, the complete differential of the propounded inexplicable function $S$ will result, namely

\begin{equation*}
dS= X^{\prime}dx +dx \left\lbrace
\renewcommand{\arraystretch}{1,7}
\setlength{\arraycolsep}{0.0mm}
\begin{array}{llllllllllllllll}
~+~&X^{\prime \prime}&~+~&X^{\prime \prime \prime}&~+~& X^{\prime \prime \prime \prime}&~+~&X^{\prime \prime \prime \prime \prime}&~+~&\text{etc.} \\
~-~&X^{\prime }&~-~&X^{\prime \prime }&~-~& X^{\prime  \prime \prime}&~+~&X^{\prime \prime \prime \prime}&~-~&\text{etc.}
\end{array}
\right\rbrace
\end{equation*}
\begin{equation*}
\begin{aligned}
&-X^{\prime \prime}~\frac{dx (1-dx)}{1 \cdot 2} \\
&+X^{\prime}~\frac{dx(1-dx)}{1 \cdot 2}
\end{aligned}
-\frac{dx(1-dx)}{1 \cdot 2}\left\lbrace
\renewcommand{\arraystretch}{1,7}
\setlength{\arraycolsep}{0.0mm}
\begin{array}{lrlrlrlrlrllllll}
~+~&X^{\prime \prime \prime}&~+~&X^{\prime \prime \prime \prime}&~+~&X^{\prime \prime \prime \prime \prime}&~+~&\text{etc.} \\
~-~&2X^{\prime \prime}&~-~&2X^{\prime \prime \prime}&~-~&2X^{\prime \prime \prime \prime }&~-~&\text{etc.} \\
~+~&X^{\prime }&~+~&X^{\prime \prime }&~+~&X^{\prime \prime \prime }&~+~&\text{etc.} \\
\end{array}
\right\rbrace
\end{equation*}

\begin{equation*}
\begin{aligned}
&+ X^{\prime \prime \prime} \frac{dx(1-dx)(2-dx)}{1 \cdot 2 \cdot 3} \\
&- 2X^{\prime \prime } \frac{dx(1-dx)(2-dx)}{1 \cdot 2 \cdot 3} \\
&+ X^{\prime }~~~ \frac{dx(1-dx)(2-dx)}{1 \cdot 2 \cdot 3} \\
\end{aligned}
-\frac{dx(1-dx)(2-dx)}{1 \cdot 2 \cdot 3}\left\lbrace
\renewcommand{\arraystretch}{1,7}
\setlength{\arraycolsep}{0.0mm}
\begin{array}{lrlrlrlrlrllllll}
~+~&X^{\prime \prime \prime \prime}&~+~&X^{\prime \prime \prime \prime \prime}&~+~&\text{etc.} \\
~-~&3X^{\prime \prime \prime}&~-~&3X^{\prime \prime \prime \prime}&~-~&\text{etc.} \\
~+~&3X^{\prime \prime}&~+~&3X^{\prime \prime \prime}&~+~&\text{etc.} \\
~-~&X^{\prime \prime}&~-~&X^{\prime \prime \prime}&~-~&\text{etc.} \\
\end{array}
\right\rbrace
\end{equation*}
\begin{equation*}
\text{etc.}
\end{equation*}
\begin{alignat*}{9}
-d. &&(X^{\prime}+X^{\prime \prime}+X^{\prime \prime \prime}+X^{\prime \prime \prime \prime}+X^{\prime \prime \prime \prime \prime}+\text{etc.}) \\[3mm]
-\frac{1}{2}dd. &&(X^{\prime}+X^{\prime \prime}+X^{\prime \prime \prime}+X^{\prime \prime \prime \prime}+X^{\prime \prime \prime \prime \prime}+\text{etc.}) \\[3mm]
-\frac{1}{6}d^3. &&(X^{\prime}+X^{\prime \prime}+X^{\prime \prime \prime}+X^{\prime \prime \prime \prime}+X^{\prime \prime \prime \prime \prime}+\text{etc.}) \\[3mm]
-\frac{1}{24} d^4. &&(X^{\prime}+X^{\prime \prime}+X^{\prime \prime \prime}+X^{\prime \prime \prime \prime}+X^{\prime \prime \prime \prime \prime}+\text{etc.}) 
\end{alignat*}
\begin{equation*}
\text{etc.}
\end{equation*}
which expression extends very far and, no matter at which point the differences just became constant, will exhibit the differential in question. For, this formula is accommodated to constant differences and at the same time the law is plain, if it is necessary to proceed further.

\paragraph*{§379}

If the series $A+B+C+D+$ etc., from which the inexplicable function

\begin{alignat*}{9}
& &&1 &&2 &&3 &&4 &&x \\
&S~=~&&A~+~&&B~+~&&C~+~&&D~+~\cdots ~+~&&X
\end{alignat*}
is formed, was of such a nature that the infinitesimal terms vanish, then, as we already noted, it will be

\begin{alignat*}{9}
dS =-d. &(X^{\prime}+X^{\prime \prime}+X^{\prime \prime \prime}+X^{\prime \prime \prime \prime}+\text{etc.}) \\[3mm]
-\frac{1}{2}dd. &(X^{\prime}+X^{\prime \prime}+X^{\prime \prime \prime}+X^{\prime \prime \prime \prime}+\text{etc.}) \\[3mm]
-\frac{1}{6}d^3 &(X^{\prime}+X^{\prime \prime}+X^{\prime \prime \prime}+X^{\prime \prime \prime \prime}+\text{etc.}) \\[3mm]
-\frac{1}{24}d^4. &(X^{\prime}+X^{\prime \prime}+X^{\prime \prime \prime}+X^{\prime \prime \prime \prime}+\text{etc.}) \\
\end{alignat*}
But if the infinitesimal terms of that series were not $=0$, but nevertheless have vanishing differences, this expression is to be added

\begin{equation*}
dx \left\lbrace
\renewcommand{\arraystretch}{1,2}
\setlength{\arraycolsep}{0.0mm}
\begin{array}{llllllllllllllll}
&~+~&X^{\prime \prime}&~+~&X^{\prime \prime \prime}&~+~& X^{\prime \prime \prime \prime}&~+~&X^{\prime \prime \prime \prime \prime}&~+~&\text{etc.} \\
X^{\prime}& \\
& ~-~&X^{\prime }&~-~&X^{\prime \prime }&~-~& X^{\prime  \prime \prime}&~+~&X^{\prime \prime \prime \prime}&~-~&\text{etc.}
\end{array}
\right\rbrace
\end{equation*}
But if just the second differences of the infinitesimal terms of this series $A+B+C+D+\text{etc.}$ vanish,  one furthermore has to add

\begin{equation*}
\frac{dx(dx-1)}{1 \cdot 2}\left\lbrace
\renewcommand{\arraystretch}{1,2}
\setlength{\arraycolsep}{0.0mm}
\begin{array}{lrlrlrlrlrlrllllll}
 & &~+~&X^{\prime \prime \prime}&~+~&X^{\prime \prime \prime \prime}&~+~&X^{\prime \prime \prime \prime \prime}&~+~&\text{etc.} \\
~+~& X^{\prime \prime} \\
& &~-~&2X^{\prime \prime}&~-~&2X^{\prime \prime \prime}&~-~&2X^{\prime \prime \prime \prime }&~-~&\text{etc.} \\
~-~&X^{\prime} \\
& &~+~&X^{\prime }&~+~&X^{\prime \prime }&~+~&X^{\prime \prime \prime }&~+~&\text{etc.} \\
\end{array}
\right\rbrace
\end{equation*}
And if just the third differences of the mentioned infinitesimal terms vanish,  except for the already exhibited expressions one additionally has to add

\begin{equation*}
\frac{dx(dx-1)(dx-2)}{1 \cdot 2 \cdot 3}\left\lbrace
\renewcommand{\arraystretch}{1,2}
\setlength{\arraycolsep}{0.0mm}
\begin{array}{lrlrlrlrlrlrlrllll}
 & &~+~&X^{\prime \prime \prime \prime}&~+~&X^{\prime \prime \prime \prime \prime}&~+~&X^{\prime \prime \prime \prime \prime \prime}&~+&\text{etc.} \\
~+~ &X^{\prime \prime \prime} & \\
& &~-~&3X^{\prime \prime \prime}&~-~&3X^{\prime \prime \prime \prime}&~-~&3X^{\prime \prime \prime \prime \prime}&~-~&\text{etc.} \\
~-~&2X^{\prime \prime} & \\
& &~+~&3X^{\prime \prime}&~+~&3X^{\prime \prime \prime}&~+~&3X^{\prime \prime \prime \prime}&~+~&\text{etc.} \\
~+~&X^{\prime} & \\
& &~-~&X^{\prime \prime}&~-~&X^{\prime \prime \prime}&~-~&X^{\prime \prime \prime}&~-~&\text{etc.} \\
\end{array}
\right\rbrace
\end{equation*}
And this will be the nature of the expressions to be added, if just higher differences of the infinitesimal terms of the series $A+B+C+D+\text{etc.}$ vanish. And hence, no matter what series is assumed, as long as its infinitesimal terms are finally  reduced to vanishing differences, one will be able to define the differential of the inexplicable function formed from it.

\paragraph*{§380}

If one puts $x=0$, it will be $X^{\prime}=A$, $X^{\prime \prime}=B$, $X^{\prime \prime \prime}=C$ etc. Therefore, as  $A+B+C+D+$ etc. is the series whose general term is $X$, if from the general terms

\begin{equation*}
\frac{dX}{dx}, \quad \frac{ddX}{2dx^2}, \quad \frac{d^3X}{6dx^3}, \quad \frac{d^4 X}{24 dx^4} \quad \text{etc.}
\end{equation*}
in like manner infinite series are formed and its sums are denoted by $\mathfrak{B}$, $\mathfrak{C}$, $\mathfrak{D}$, $\mathfrak{E}$ etc., respectively, the sum of $\omega$ terms of the series
\begin{equation*}
A+B+C+D+\text{etc.}
\end{equation*}
will be expressed in such a way that it does not matter whether $\omega$ is an integer or not. Therefore, let us put $x$ for $\omega$ that 

\begin{alignat*}{9}
& &&1 &&2 &&3 &&4 &&x \\
&S~=~&&A~+~&&B~+~&&C~+~&&D~+~\cdots ~+~&&X
\end{alignat*}
and if the infinitesimal terms of this series vanish, it will be

\begin{equation*}
S=-\mathfrak{B}x-\mathfrak{C}x^2-\mathfrak{D}x^3-\mathfrak{E}x^4-\text{etc.}
\end{equation*}
But if at least the infinitesimal terms have constant first differences, one additionally has to add this

\begin{equation*}
x \left\lbrace
\renewcommand{\arraystretch}{1,2}
\setlength{\arraycolsep}{0.0mm}
\begin{array}{llllllllllllllll}
&~+~&B&~+~&C&~+~& D&~+~& E&~+~&\text{etc.} \\
A& \\
& ~-~&A&~-~&B&~-~& C&~+~&D&~-~&\text{etc.}
\end{array}
\right\rbrace
\end{equation*}
But if just the second differences of those infinitesimal terms vanish, one has to add

\begin{equation*}
\frac{x(x-1)}{1 \cdot 2}\left\lbrace
\renewcommand{\arraystretch}{1,2}
\setlength{\arraycolsep}{0.0mm}
\begin{array}{lrlrlrlrlrlrlrllll}
 & &~+~&C&~+~&D&~+~&E&~+~&F&~+~&\text{etc.} \\
~+~& B \\
& &~-~&2B&~-~&2C&~-~&2D&~-~&2E &~-~&\text{etc.} \\
~-~&C \\
& &~+~&A&~+~&B&~+~&C&~+~&D &~+~&\text{etc.} \\
\end{array}
\right\rbrace
\end{equation*}
If just the third differences vanish, additionally this infinite series has to be added
\begin{equation*}
\frac{x(x-1)(x-2)}{1 \cdot 2 \cdot 3}\left\lbrace
\renewcommand{\arraystretch}{1,2}
\setlength{\arraycolsep}{0.0mm}
\begin{array}{lrlrlrlrlrlrlrllll}
 & &~+~&D&~+~&E&~+~&F&~+~&G&~+~&\text{etc.} \\
~+~ &C & \\
& &~-~&3C&~-~&3D&~-~&3E&~-~&3F&~-~&\text{etc.} \\
~-~&2B & \\
& &~+~&3B&~+~&3C&~+~&3D&~+~&3E&~+~&\text{etc.} \\
~+~&A & \\
& &~-~& A&~-~&B&~-~&C&~-~&D&~-~&\text{etc.} \\
\end{array}
\right\rbrace
\end{equation*}
\begin{equation*}
\text{etc.}
\end{equation*}

\paragraph*{§381}

Let us also apply these things to the species of inexplicable functions  consisting of continuous products of several terms of the propounded series $A+B+C+D$ etc., and let 

\begin{alignat*}{9}
& &&1 &&2 &&3 && 4 \cdots &&x \\
&S~=~&&A~\cdot ~&&B~\cdot ~&& C~\cdot ~&& D \cdots &&X
\end{alignat*}
and at first find the value $\Sigma$ $S$ is transformed into, if  one writes $x+\omega$ instead of $x$. But let us, as before, put that $Z$ is the term corresponding to  the index  $=x+ \omega$ of the series $A+B+C+D+$etc., whereas $X$ corresponds to the index $x$. To reduce this to the preceding case, let us take logarithms and it will be

\begin{equation*}
\log S= \log A+ \log B+ \log C+ \log D +\cdots + \log X.
\end{equation*}
If now the infinitesimal terms of this series vanish by applying the same method we used before, it will be

\begin{alignat*}{9}
&\log \Sigma = \log S &&+ \log X^{\prime} &&+ \log X^{\prime \prime} &&+\log X^{\prime \prime \prime} &&+\text{etc.} \\
&    &&- \log Z^{\prime} &&- \log Z^{\prime \prime} &&- \log Z^{\prime \prime \prime}&&-\text{etc.}
\end{alignat*}
and hence by going back to numbers it will be

\begin{equation*}
\Sigma = S \cdot \frac{X^{\prime}}{Z^{\prime}}\cdot \frac{X^{\prime \prime}}{Z^{\prime \prime}}\cdot \frac{X^{\prime \prime \prime}}{Z^{\prime \prime \prime}}\cdot \frac{X^{\prime \prime \prime \prime}}{Z^{\prime \prime \prime \prime}}\cdot \text{etc.};
\end{equation*}
therefore, this expression holds, if the infinitesimal terms of the series $A$, $B$, $C$, $D$ etc. become equal to $1$. But if the logarithms of the infinitesimal terms of this series do not vanish, but nevertheless have vanishing differences, then to that series we found for $\log \Sigma$ one additionally has to add this series

\begin{equation*}
\omega \log X^{\prime} +\omega \left(\log \frac{X^{\prime \prime}}{X^{\prime}}+\log \frac{X^{\prime \prime \prime}}{X^{\prime \prime}}+\log \frac{X^{\prime \prime \prime \prime}}{X^{\prime \prime \prime}}+\text{etc.} \right)
\end{equation*}
and, so by taking exponentials again, one will have

\begin{equation*}
\Sigma = S \cdot {X^{\prime}}^{\omega} \cdot \frac{{X^{\prime}}^{\omega}{X^{\prime \prime}}^{1-\omega}}{Z^{\prime}}\cdot \frac{{X^{\prime   \prime \prime}}^{\omega}{X^{\prime \prime}}^{1-\omega}}{Z^{\prime \prime}}\cdot \frac{{X^{\prime \prime \prime \prime}}^{\omega}{X^{\prime \prime \prime}}^{1-\omega}}{Z^{\prime \prime \prime}}\cdot \text{etc.}
\end{equation*}

\paragraph*{§382}

If we put $x=0$ in which case $S=1$ and $X^{\prime}=A$, $X^{\prime \prime}=B$, $X^{\prime \prime \prime}=C$ etc., $\Sigma$ will denote the product of $\omega$ terms of this series $A$, $B$, $C$, $D$ etc. If we write $x$ for $\omega$ that $\Sigma$ obtains the value we had attributed to $S$ before such that 

\begin{alignat*}{9}
& &&1 &&2 &&3 && 4 \cdots &&x \\
&S~=~&&A~\cdot ~&&B~\cdot ~&& C~\cdot ~&& D \cdots &&X,
\end{alignat*}
since now $Z^{\prime}$, $Z^{\prime \prime}$, $Z^{\prime \prime \prime}$ etc. go over into $X^{\prime}$, $X^{\prime \prime}$, $X^{\prime \prime \prime}$ etc., if the logarithms of the infinitesimal terms of this series $A$, $B$, $C$, $D$, $E$ etc. vanish, $S$ will be expressed this way

\begin{equation*}
S=\frac{A}{X^{\prime}}\cdot \frac{B}{X^{\prime \prime}}\cdot \frac{C}{X^{\prime \prime \prime}}\cdot \frac{D}{X^{\prime \prime \prime \prime}}\cdot \frac{E}{X^{\prime \prime \prime \prime \prime}}\cdot \text{etc.}
\end{equation*}
But if just the differences of the logarithms of the infinitesimal terms of the series $A$, $B$, $C$, $D$ etc. vanish,  this function $S$ will be expressed the following way 

\begin{equation*}
S=A^x \cdot \frac{B^x A^{1-x}}{X^{\prime}}\cdot \frac{C^x B^{1-x}}{X^{\prime \prime}}\cdot \frac{D^x C^{1-x}}{X^{\prime \prime \prime}}\cdot \frac{E^x D^{1-x}}{X^{\prime \prime \prime}}\cdot \text{etc.};
\end{equation*}
If just the second differences of those logarithms vanish, it is easily concluded from the preceding, factors of which kind are to be added; we omit this case here, since it usually does not occur. Moreover, I will show the use of these expressions in the task of interpolation in the following chapter.

\paragraph*{§383}

Since here mainly the differentiation of inexplicable functions is propounded, let us investigate the differential of this function

\begin{equation*}
S=A \cdot B \cdot C \cdot D \cdots X.
\end{equation*}
For this, let us go back to the equation found before

\begin{alignat*}{9}
&\log \Sigma = \log S &&+ \log X^{\prime} &&+ \log X^{\prime \prime} &&+\log X^{\prime \prime \prime} &&+\text{etc.} \\
&    &&- \log Z^{\prime} &&- \log Z^{\prime \prime} &&- \log Z^{\prime \prime \prime}&&-\text{etc.},
\end{alignat*}
and since $\log Z$ results from $\log X$, if  one writes $x+\omega$ instead of $x$, it will be

\begin{equation*}
\log Z= \log X+ \frac{\omega}{dx}d. \log X+\frac{\omega^2}{2 dx^2}dd.\log X+\frac{\omega^3}{6dx^3}d^3. \log X+\text{etc.}; 
\end{equation*}
having substituted these values for $\log Z^{\prime}$, $\log Z^{\prime \prime \prime}$, $\log Z^{\prime \prime \prime}$ etc., one will have

\begin{alignat*}{9}
& \log \Sigma = \log S&& -~~\frac{\omega}{dx}&&d. &&(\log X^{\prime}+ \log X^{\prime \prime}+\log X^{\prime \prime \prime}+ \log X^{\prime \prime \prime \prime} +\text{etc.}) \\[3mm]
&   &&-\frac{\omega^2}{2dx^2}&&dd. &&(\log X^{\prime}+ \log X^{\prime \prime}+\log X^{\prime \prime \prime}+ \log X^{\prime \prime \prime \prime} +\text{etc.}) \\[3mm]  
&   &&-\frac{\omega^3}{6dx^3}&&d^3. &&(\log X^{\prime}+ \log X^{\prime \prime}+\log X^{\prime \prime \prime}+ \log X^{\prime \prime \prime \prime} +\text{etc.}) \\    
& && &&\text{etc.}  
\end{alignat*} 
Now put $\omega =dx$ and it will be $\log \Sigma = \log S+d. \log S$ and hence it will be

\begin{alignat*}{9}
&\frac{dS}{S}&&=-d. &&(\log X^{\prime}+ \log X^{\prime \prime}+\log X^{\prime \prime \prime}+ \log X^{\prime \prime \prime \prime} +\text{etc.}) \\[3mm]
& &&-\frac{1}{2}dd. &&(\log X^{\prime}+ \log X^{\prime \prime}+\log X^{\prime \prime \prime}+ \log X^{\prime \prime \prime \prime} +\text{etc.})  \\[3mm]
& &&-\frac{1}{6}d^3. && (\log X^{\prime}+ \log X^{\prime \prime}+\log X^{\prime \prime \prime}+ \log X^{\prime \prime \prime \prime} +\text{etc.}) \\
& && \text{etc.},
\end{alignat*}
which formula holds, if the logarithms of the infinitesimal terms of the series $A$, $B$, $C$, $D$ etc. vanish; but if they do not vanish, but nevertheless have vanishing differences, then to the preceding expression of the complete differential one additionally has to add this series

\begin{equation*}
dx \log X^{\prime}+dx \left(\log \frac{X^{\prime \prime}}{X^{\prime}}+\log \frac{X^{\prime \prime \prime}}{X^{\prime \prime}}+\log \frac{X^{\prime \prime \prime \prime}}{X^{\prime \prime \prime}}+\text{etc.}\right)
\end{equation*}
in order to obtain the complete differential.

\paragraph*{§384}

The same can also be achieved  in another way. Put $x=0$ in which case $\log S$ goes over into $0$. Then, form series whose general terms are

\begin{equation*}
\log X, \quad \frac{d. \log X}{dx}, \quad \frac{dd.\log X}{2 dx^2}, \quad \frac{d^3. \log X}{6dx^3}\quad \text{etc.},
\end{equation*}
and  we want the sums of these infinite series to be $\mathfrak{A}$, $\mathfrak{B}$, $\mathfrak{C}$, $\mathfrak{D}$ etc., respectively. Write $x$ for $\omega$ that $\Sigma = S$ and it will be

\begin{equation*}
\log S= - \mathfrak{B}x- \mathfrak{C}x^2-\mathfrak{D}x^3-\mathfrak{E}x^4-\text{etc.},
\end{equation*}
if the logarithms of the infinitesimal terms of the series $A$, $B$, $C$, $D$ etc. whose general term is $X$ vanish; but if just the differences of these logarithms vanish, it will be

\begin{equation*}
\log S = x \log A+ x \left(\log \frac{B}{A}+\log \frac{C}{B}+\log \frac{D}{C}+\log \frac{E}{D}+\text{etc.}\right)
\end{equation*}
\begin{equation*}
- \mathfrak{B}x- \mathfrak{C}x^2- \mathfrak{D}x^3- \mathfrak{D}x^4-\text{etc.}
\end{equation*}
And hence the differential of $\log S$ will be

\begin{equation*}
\frac{dS}{S}=dx \log A +dx \left(\log \frac{B}{A}+\log \frac{C}{B}+\log \frac{D}{C}+\log \frac{E}{D}+\text{etc.}\right)
\end{equation*}
\begin{equation*}
-\mathfrak{B}xdx- 2\mathfrak{C}xdx- 3\mathfrak{D}x^2dx -4 \mathfrak{E}x^3dx-\text{etc.}
\end{equation*}
But if the complete differential is desired, it will be

\begin{equation*}
\frac{dS}{S}=dx \log A +dx \left(\log \frac{B}{A}+\log \frac{C}{B}+\log \frac{D}{C}+\log \frac{E}{D}+\text{etc.}\right)
\end{equation*}
\begin{equation*}
- \mathfrak{B}dx -\mathfrak{C}(2xdx+dx^2)-\mathfrak{D}(3xxdx+3xdx^2+dx^3)-\text{etc.}
\end{equation*}
To show the use of these formulas we add the following examples which we resolve in both ways.

\subsection*{Example 1}

\textit{To find the differential of this inexplicable function}

\begin{equation*}
S=\frac{1}{2}\cdot \frac{3}{4}\cdot \frac{5}{6}\cdot \frac{7}{8}\cdots \frac{2x-1}{2x}.
\end{equation*}
Here, it is especially to be noted that the infinitesimal terms of these factors go over into $1$ and hence their logarithms vanish. Since $X=\frac{2x-1}{2x}$, it will be

\begin{equation*}
X^{\prime}=\frac{2x+1}{2x+1}, \quad X^{\prime \prime}=\frac{2x+3}{2x+4}, \quad X^{\prime \prime \prime}=\frac{2x+5}{2x+6} \quad \text{etc.}
\end{equation*}
and in general

\begin{equation*}
X^{|n|}=\frac{2x+2n-1}{2x+2n};
\end{equation*}
therefore, it will be

\begin{alignat*}{9}
&\log X^{|n|}&&=+\log (2x+2n-1) &&- \log (2x+2n) \\
&d. \log X^{|n|}&&=+\quad \frac{2dx}{2x+2n-1} &&-\quad  \frac{2dx}{2x+2n} \\[3mm]
&dd. \log X^{|n|}&&=-\frac{4dx^2}{(2x+2n-1)^2} &&+ \frac{4dx^2}{(2x+2n)^2} \\[3mm]
&d^3. \log X^{|n|}&&=+\frac{2 \cdot 2 \cdot 4dx^2}{(2x+2n-1)^3} &&- \frac{2\cdot 2 \cdot 4dx^2}{(2x+2n)^3} \\[3mm]
&d^4. \log X^{|n|}&&=-\frac{2 \cdot 2 \cdot 4 \cdot 6 dx^4}{(2x+2n-1)^4} &&+ \frac{2 \cdot 2 \cdot 4 \cdot 6 dx^4}{(2x+2n)^4} \\
& && &&\text{etc.};
\end{alignat*}
therefore, the complete differential will be

\begin{equation*}
\begin{aligned}
\frac{dS}{S}=&-2dx~~
\left\lbrace
\renewcommand{\arraystretch}{1,7}
\setlength{\arraycolsep}{0.0mm}
\begin{array}{lrlrlrlrlrlrlrllll}
 &\dfrac{1}{2x+1}&~+~& \dfrac{1}{2x+3}&~+~& \dfrac{1}{2x+5}&~+~&\text{etc.} \\
 ~-~ &\dfrac{1}{2x+2}&~-~& \dfrac{1}{2x+4}&~-~& \dfrac{1}{2x+6}&~-~&\text{etc.} \\
\end{array}
\right\rbrace
 \\
&+\frac{4}{2}dx^2
\left\lbrace
\renewcommand{\arraystretch}{1,7}
\setlength{\arraycolsep}{0.0mm}
\begin{array}{lrlrlrlrlrlrlrllll}
 &\dfrac{1}{(2x+1)^2}&~+~& \dfrac{1}{(2x+3)^2}&~+~& \dfrac{1}{(2x+5)^2}&~+~&\text{etc.} \\
 ~-~ &\dfrac{1}{(2x+2)^2}&~-~& \dfrac{1}{(2x+4)^2}&~-~& \dfrac{1}{(2x+6)^2}&~-~&\text{etc.} 
\end{array}
\right\rbrace 
 \\
&-\frac{8}{3}dx^3
\left\lbrace
\renewcommand{\arraystretch}{1,7}
\setlength{\arraycolsep}{0.0mm}
\begin{array}{lrlrlrlrlrlrlrllll}
 &\dfrac{1}{(2x+1)^3}&~+~& \dfrac{1}{(2x+3)^3}&~+~& \dfrac{1}{(2x+5)^3}&~+~&\text{etc.} \\
 ~-~ &\dfrac{1}{(2x+2)^3}&~-~& \dfrac{1}{(2x+4)^3}&~-~& \dfrac{1}{(2x+6)^3}&~-~&\text{etc.} \\
\end{array}
\right\rbrace
\end{aligned}
\end{equation*}
\begin{equation*}
\text{etc.}
\end{equation*}
But if only the first differential is in question, it will be

\begin{equation*}
\frac{dS}{S}=-2dx \cdot \left(\frac{1}{(2x+1)(2x+2)}+\frac{1}{(2x+3)(2x+4)}+\frac{1}{(2x+5)(2x+6)}+\text{etc.}\right),
\end{equation*}
which same equation is found by  the other method given in § 394. Since 

\begin{equation*}
\log X=\log\frac{2x-1}{2x},
\end{equation*}
it will be

\begin{equation*}
\frac{d. \log X}{dx}=\frac{2}{2x-1}-\frac{1}{x}, \quad \frac{dd. \log X}{2dx^2}=-\frac{2}{(2x-1)^2}+\frac{1}{2xx},
\end{equation*}
\begin{equation*}
\frac{d^3 \log X}{6dx^3}=\frac{8}{3(2x-1)^3}-\frac{1}{3x^3} \quad \text{etc.}
\end{equation*}
and hence it will be

\begin{equation*}
\begin{aligned}
\mathfrak{A} & = \log \frac{1}{2}+\log \frac{3}{4}+\log \frac{5}{6}+\log \frac{7}{8}+\text{etc.} \\
\mathfrak{B} & = 
\left\lbrace
\renewcommand{\arraystretch}{1,7}
\setlength{\arraycolsep}{0.0mm}
\begin{array}{lrlrlrlrlrlrlrllll}
 &\dfrac{2}{1}&~+~& \dfrac{2}{3}&~+~& \dfrac{2}{5}&~+~&\dfrac{2}{7}&~+~&\dfrac{2}{9}&~+~&\text{etc.} \\
- &\dfrac{2}{2}&~-~& \dfrac{2}{4}&~-~& \dfrac{2}{6}&~-~&\dfrac{2}{8}&~-~&\dfrac{2}{10}&~-~&\text{etc.} \\
\end{array}
\right\rbrace 
= 2 \log 2 \\
\mathfrak{C} &= -\frac{4}{2}
~\left\lbrace
\renewcommand{\arraystretch}{2,0}
\setlength{\arraycolsep}{0.0mm}
\begin{array}{lrlrlrlrlrlrlrllll}
 &\dfrac{1}{1}&~+~& \dfrac{1}{3^2}&~+~& \dfrac{1}{5^2}&~+~&\dfrac{1}{7^2}&~+~&\text{etc.} \\
- &\dfrac{1}{2^2}&~-~& \dfrac{1}{4^2}&~-~& \dfrac{1}{6^2}&~-~&\dfrac{1}{8^2}&~-~&\text{etc.} \\
\end{array}
\right\rbrace 
 \\
\mathfrak{D}&=+\frac{8}{3}
~\left\lbrace
\renewcommand{\arraystretch}{2,0}
\setlength{\arraycolsep}{0.0mm}
\begin{array}{lrlrlrlrlrlrlrllll}
 &\dfrac{1}{1}&~+~& \dfrac{1}{3^3}&~+~& \dfrac{1}{5^3}&~+~&\dfrac{1}{7^3}&~+~&\text{etc.} \\
- &\dfrac{1}{2^3}&~-~& \dfrac{1}{4^3}&~-~& \dfrac{1}{6^3}&~-~&\dfrac{1}{8^3}&~-~&\text{etc.} \\
\end{array}
\right\rbrace  \\
\mathfrak{E}& =-\frac{16}{4}
\left\lbrace
\renewcommand{\arraystretch}{2,0}
\setlength{\arraycolsep}{0.0mm}
\begin{array}{lrlrlrlrlrlrlrllll}
 &\dfrac{1}{1}&~+~& \dfrac{1}{3^4}&~+~& \dfrac{1}{5^4}&~+~&\dfrac{1}{7^4}&~+~&\text{etc.} \\
- &\dfrac{1}{2^4}&~-~& \dfrac{1}{4^4}&~-~& \dfrac{1}{6^4}&~-~&\dfrac{1}{8^4}&~-~&\text{etc.} \\
\end{array}
\right\rbrace 
 \\
\end{aligned}
\end{equation*}
\begin{equation*}
\text{etc.}
\end{equation*}
 or it will be
 
 \begin{alignat*}{9}
&\mathfrak{B}&&=+\frac{2}{1}&&\bigg(1&&-\frac{1}{2}&&+\frac{1}{3}&&-\frac{1}{4}&&+\frac{1}{5}&&-\text{etc.}\bigg) \\[3mm]
&\mathfrak{C}&&=-\frac{4}{2}&&\bigg(1&&-\frac{1}{2^2}&&+\frac{1}{3^2}&&-\frac{1}{4^2}&&+\frac{1}{5^2}&&-\text{etc.}\bigg) \\[3mm]
&\mathfrak{D}&&=+\frac{8}{3}&&\bigg(1&&-\frac{1}{2^3}&&+\frac{1}{3^3}&&-\frac{1}{4^3}&&+\frac{1}{5^3}&&-\text{etc.}\bigg) \\[3mm]
&\mathfrak{E}&&=-\frac{16}{4}&&\bigg(1&&-\frac{1}{2^4}&&+\frac{1}{3^4}&&-\frac{1}{4^4}&&+\frac{1}{5^4}&&-\text{etc.}\bigg) \\
& && && &&\text{etc.}   
 \end{alignat*}
 Having substituted the found values, it will be 
\begin{alignat*}{9}
&\frac{dS}{S}&&=-2dx &&\bigg(1&&-\frac{1}{2}&&+\frac{1}{3}&&-\frac{1}{4}&&+\frac{1}{5}&&-\text{etc.}\bigg)\\[3mm]
& &&+~~4xdx &&\bigg(1&&-\frac{1}{2^2}&&+\frac{1}{3^2}&&-\frac{1}{4^2}&&+\frac{1}{5^2}&&-\text{etc.}\bigg)\\[3mm]
& &&-~8x^2dx &&\bigg(1&&-\frac{1}{2^3}&&+\frac{1}{3^3}&&-\frac{1}{4^3}&&+\frac{1}{5^3}&&-\text{etc.}\bigg)\\[3mm]
& &&+16x^3dx &&\bigg(1&&-\frac{1}{2^4}&&+\frac{1}{3^4}&&-\frac{1}{4^4}&&+\frac{1}{5^4}&&-\text{etc.}\bigg)\\
& && && && &&\text{etc.}
\end{alignat*} 
If $x=0$ in which case $\log S=0$ and $S=1$, it will be $dS=-2dx \log 2.$

\subsection*{Example 2}

\textit{To find the differential of this inexplicable function}

\begin{equation*}
S=1\cdot 2 \cdot 3 \cdot 4 \cdots x.
\end{equation*}
The terms of this series $1$, $2$, $3$, $4$ etc. grow to infinity in such a way that the differences of the logarithms vanish; for,

\begin{equation*}
\log (\infty +1)- \log \infty = \log \left(1+\frac{1}{\infty}\right)=\frac{1}{\infty}=0.
\end{equation*}
Since  $X=x$, it will be

\begin{equation*}
X^{\prime}=x+1, \quad X^{\prime \prime}=x+2, \quad X^{\prime \prime \prime}=x+3 \quad \text{etc.};
\end{equation*}
but, further, because of $\log X= \log x$, it will be

\begin{equation*}
d. \log X =\frac{dx}{x}, \quad dd. \log X =-\frac{dx^2}{x^2}, \quad d^3. \log X=\frac{2dx^3}{x^3}, \quad d^4. \log X =-\frac{2 \cdot 3 dx^4}{x^4} \quad \text{etc.};
\end{equation*}
hence, if the last logarithms would vanish, it would be

\begin{alignat*}{9}
&\frac{dS}{S}&&=-dx &&\bigg(\frac{\quad ~~ 1 \quad ~~}{x+1}&&+\frac{\quad ~~1  \quad ~~}{x+2}&&+\frac{\quad ~~ 1 \quad ~~}{x+3}&&+\frac{\quad ~~ 1 \quad ~~}{x+4}&&+\text{etc.}\bigg) \\[3mm]
& &&+\frac{dx^2}{2}&&\bigg(\frac{1}{(x+1)^2}&&+\frac{1}{(x+2)^2}&&+\frac{1}{(x+3)^2}&&+\frac{1}{(x+4)^2}&&+\text{etc.}\bigg)\\[3mm]
& &&-\frac{dx^3}{3}&&\bigg(\frac{1}{(x+1)^3}&&+\frac{1}{(x+2)^3}&&+\frac{1}{(x+3)^3}&&+\frac{1}{(x+4)^3}&&+\text{etc.}\bigg)\\
& && && && &&\text{etc.}
\end{alignat*}
But because just the differences of the logarithms vanish, one additionally has to add this expression

\begin{equation*}
dx \log(x+1)+dx \left(\log \frac{x+2}{x+1}+\log \frac{x+3}{x+2}+\log \frac{x+4}{x+3}+\log \frac{x+5}{x+4}+\text{etc.}\right).
\end{equation*}
But, because 

\begin{alignat*}{9}
&\log \frac{x+2}{x+1}&&=\frac{1}{x+1}&&-\frac{1}{2(x+1)^2}&&+\frac{1}{3(x+1)^3}&&-\frac{1}{4(x+1)^4}&&+\text{etc.} \\[3mm]
&\log \frac{x+3}{x+2}&&=\frac{1}{x+2}&&-\frac{1}{2(x+2)^2}&&+\frac{1}{3(x+2)^3}&&-\frac{1}{4(x+2)^4}&&+\text{etc.} \\
&  && && &&\text{etc.},
\end{alignat*}
the complete differential will be

\begin{alignat*}{9}
& \frac{dS}{S}=dx \log (x+1) &&-\frac{1}{2}(dx-dx^2)&&\bigg(\frac{1}{(x+1)^2}&&+\frac{1}{(x+2)^2}&&+\frac{1}{(x+3)^2}&&+\text{etc.}\bigg) \\[3mm]
& && +\frac{1}{3}(dx-dx^3)&&\bigg(\frac{1}{(x+1)^3}&&+\frac{1}{(x+2)^3}&&+\frac{1}{(x+3)^3}&&+\text{etc.}\bigg)\\[3mm]
& && -\frac{1}{4}(dx-dx^4)&&\bigg(\frac{1}{(x+1)^4}&&+\frac{1}{(x+2)^4}&&+\frac{1}{(x+3)^4}&&+\text{etc.}\bigg)\\[3mm]
& && -\frac{1}{5}(dx-dx^5)&&\bigg(\frac{1}{(x+1)^5}&&+\frac{1}{(x+2)^5}&&+\frac{1}{(x+3)^5}&&+\text{etc.}\bigg)\\
& && &&\text{etc.}
\end{alignat*}
But if we want to express this differential by means of the other method, since 

\begin{equation*}
\log X = \log x, \quad \frac{d. \log X}{dx}=\frac{1}{x}, \quad \frac{dd. \log X}{2dx^2}=-\frac{1}{2x^2}, \quad \frac{d^3. \log X}{6dx^3}=\frac{1}{3x^3}, \quad \frac{d^4. \log X}{24 dx^4}=-\frac{1}{4x^4} \quad \text{etc.},
\end{equation*}
one will have the following series

\begin{alignat*}{9}
&\mathfrak{A}&&=  &&\log 1 &&+\log 2 &&+\log 3 &&+\log 4 &&+\log 5 &&+\text{etc.} \\[3mm]
&\mathfrak{B}&&=+1&&\bigg(1&&+~~\frac{1}{2}&&+~~\frac{1}{3}&&+~~\frac{1}{4}&&+~~\frac{1}{5}&&+\text{etc.}\bigg)\\[3mm]
&\mathfrak{C}&&=-\frac{1}{2}&&\bigg(1&&+~\frac{1}{2^2}&&+~\frac{1}{3^2}&&+~\frac{1}{4^2}&&+~\frac{1}{5^2}&&+\text{etc.}\bigg)\\[3mm]
&\mathfrak{C}&&=+\frac{1}{3}&&\bigg(1&&+~\frac{1}{2^3}&&+~\frac{1}{3^3}&&+~\frac{1}{4^3}&&+~\frac{1}{5^3}&&+\text{etc.}\bigg)\\[3mm]
&\mathfrak{D}&&=-\frac{1}{4}&&\bigg(1&&+~\frac{1}{2^4}&&+~\frac{1}{3^4}&&+~\frac{1}{4^4}&&+~\frac{1}{5^4}&&+\text{etc.}\bigg)\\
& && && && && &&\text{etc.}
\end{alignat*}
Therefore, because of $\log A =\log 1 =0$  from § 384 it will be

\begin{alignat*}{9}
&\log S && = x &&\bigg(\log \frac{2}{1}&&+\log\frac{3}{2}&&+\log \frac{4}{3}&&+\log \frac{5}{4}&&+\text{etc.}\bigg) \\[3mm]
&  && -x && \bigg(~~~1&&+~~\frac{1}{2}&&+~~\frac{1}{3}&&+~~\frac{1}{4}&&+\text{etc.}\bigg)\\[3mm]
&  && +\frac{1}{2} x^2 && \bigg(~~~1&&+~~\frac{1}{2^2}&&+~~\frac{1}{3^2}&&+~~\frac{1}{4^2}&&+\text{etc.}\bigg)\\[3mm]
&  && -\frac{1}{3}x^3 && \bigg(~~~1&&+~~\frac{1}{2^3}&&+~~\frac{1}{3^3}&&+~~\frac{1}{4^3}&&+\text{etc.}\bigg)\\[3mm]
&  && +\frac{1}{4}x^4 && \bigg(~~~1&&+~~\frac{1}{2^4}&&+~~\frac{1}{3^4}&&+~~\frac{1}{4^4}&&+\text{etc.}\bigg)\\
& && && && &&\text{etc.}
\end{alignat*}
But the two first series $x$ is multiplied by, even though both have an infinite sum, nevertheless, taken together, have a finite sum. For, if $n$ terms are taken of both of them, it will be

\begin{equation*}
\log (n+1)-1-\frac{1}{2}-\frac{1}{3}-\frac{1}{4}-\cdots -\frac{1}{n}.
\end{equation*}
But above (§ 142) we found that

\begin{equation*}
1+\frac{1}{2}+\frac{1}{3}+\frac{1}{4}+\cdots +\frac{1}{n}= \text{Const.}+\log n + \frac{1}{2n}-\frac{\mathfrak{A}}{2n^2}+\frac{\mathfrak{B}}{4n^4}-\text{etc.}
\end{equation*}
and this constant will be found to be $=0.5772156649015325$. If one puts $n=\infty$, it will be

\begin{equation*}
1+\frac{1}{2}+\frac{1}{3}+\frac{1}{4}+\cdots +\frac{1}{\infty}=\text{Const.}+\log \infty,
\end{equation*}
whence the value of those two series continued to infinity will be

\begin{equation*}
=\log (\infty +1)-\text{Const.}-\log \infty =- \text{Const.}
\end{equation*}
From this it will be

\begin{equation*}
\log S =- x \cdot 0.5772156649015325
\end{equation*}
\begin{alignat*}{9}
&+\frac{1}{2}xx &&\bigg(1 &&+\frac{1}{2^2}&&+\frac{1}{3^2}&&+\frac{1}{4^2}&&+\frac{1}{5^2}&&+\text{etc.}\bigg) \\[3mm]
&-\frac{1}{3}x^3 &&\bigg(1 &&+\frac{1}{2^3}&&+\frac{1}{3^3}&&+\frac{1}{4^3}&&+\frac{1}{5^3}&&+\text{etc.}\bigg) \\[3mm]
&+\frac{1}{4}x^4 &&\bigg(1 &&+\frac{1}{2^4}&&+\frac{1}{3^4}&&+\frac{1}{4^4}&&+\frac{1}{5^4}&&+\text{etc.}\bigg) \\
& && && &&\text{etc.},
\end{alignat*}
whence the differentials of any order are easily found. For, it will be 

\begin{equation*}
\frac{dS}{S} =- dx \cdot 0.5772156649015325
\end{equation*}
\begin{alignat*}{9}
&+xdx &&\bigg(1 &&+\frac{1}{2^2}&&+\frac{1}{3^2}&&+\frac{1}{4^2}&&+\frac{1}{5^2}&&+\text{etc.}\bigg) \\[3mm]
&-x^2dx &&\bigg(1 &&+\frac{1}{2^3}&&+\frac{1}{3^3}&&+\frac{1}{4^3}&&+\frac{1}{5^3}&&+\text{etc.}\bigg) \\[3mm]
&+x^3dx &&\bigg(1 &&+\frac{1}{2^4}&&+\frac{1}{3^4}&&+\frac{1}{4^4}&&+\frac{1}{5^4}&&+\text{etc.}\bigg) \\
& && && &&\text{etc.}
\end{alignat*}
But if these series are collected into one sum, it will be

\begin{equation*}
\frac{dS}{S}=- dx \cdot 0.5772156649015325 +\frac{xdx}{1(1+x)}+\frac{xdx}{2(2+x)}+\frac{xdx}{3(3+x)}+\frac{xdx}{4(4+x)}+\text{etc.}
\end{equation*}
Hence, if $x=0$, it will be

\begin{equation*}
\frac{dS}{S}=-dx \cdot 0.5772156649015325.
\end{equation*}
From the first expression on the other hand it will be in this case

\begin{alignat*}{9}
&\frac{dS}{S}&&~=~&&-\frac{1}{2}dx \bigg(1&&+\frac{1}{2^2}&&+\frac{1}{3^2}&&+\frac{1}{4^2}&&+\text{etc.}\bigg)\\[3mm]
& && &&+\frac{1}{3}dx \bigg(1&&+\frac{1}{2^3}&&+\frac{1}{3^3}&&+\frac{1}{4^3}&&+\text{etc.}\bigg)\\[3mm]
& && &&-\frac{1}{4}dx \bigg(1&&+\frac{1}{2^4}&&+\frac{1}{3^4}&&+\frac{1}{4^4}&&+\text{etc.}\bigg)\\[3mm]
& && &&+\frac{1}{5}dx \bigg(1&&+\frac{1}{2^5}&&+\frac{1}{3^5}&&+\frac{1}{4^5}&& +\text{etc.}\bigg)\\
& && && &&\text{etc.}
\end{alignat*}

\paragraph*{§385}

Hence one is also able to exhibit the differentials of inexplicable functions of this kind in any special case, since here we found the complete differentials. Therefore, if such functions enter expressions which seem to be undetermined, functions of which kind we treated in the preceding chapter, one will be able to define the values applying the same method, as it will be understood from the following examples.

\subsection*{Example 1}

\textit{To determine the value of this expression}

\begin{equation*}
\frac{1+\frac{1}{2}+\frac{1}{3}+\cdots +\frac{1}{x}}{x(x-1)}-\frac{1}{(x-1)(2x-1)}.
\end{equation*}
\textit{in the case $x=1$.}\\[2mm]
Let us put

\begin{equation*}
1+\frac{1}{2}+\frac{1}{3}+\cdots +\frac{1}{x}=S;
\end{equation*}
from § 372 it will be

\begin{alignat*}{9}
&S &&=x &&\bigg(1&&+\frac{1}{2^2}&&+\frac{1}{3^2}&&+\frac{1}{4^2}&&+\text{etc.}\bigg) \\[3mm]
& &&-x^2 &&\bigg(1&&+\frac{1}{2^3}&&+\frac{1}{3^3}&&+\frac{1}{4^3}&&+\text{etc.}\bigg) \\[3mm]
& &&+x^3 &&\bigg(1&&+\frac{1}{2^4}&&+\frac{1}{3^4}&&+\frac{1}{4^4}&&+\text{etc.}\bigg) \\
& && && && &&\text{etc.},
\end{alignat*}
or, because

\begin{alignat*}{9}
&S &&~=~ &&+ \quad 1 && +\quad \frac{1}{2}&&+\quad \frac{1}{3}&&+\quad \frac{1}{4}&&+ \quad \frac{1}{5}&&+\text{etc.} \\[3mm]
& &&     &&-\frac{1}{1+x}&&-\frac{1}{2+x}&&-\frac{1}{3+x}&&-\frac{1}{4+x}&&-\frac{1}{5+x}&&-\text{etc.},
\end{alignat*}
if each term of the upper series is combined with the preceding of the lower series, it will be

\begin{equation*}
S=1+\frac{x-1}{2(1+x)}+\frac{x-1}{3(2+x)}+\frac{x-1}{4(3+x)}+\text{etc.},
\end{equation*}
which expression, since one has to put $x=1$, is more convenient. Therefore, let $x=1+ \omega$ and it will be

\begin{equation*}
S=1+\frac{\omega}{2(2+ \omega)}+\frac{\omega}{3(3 +\omega)}+\frac{\omega}{4(4+ \omega)}+\text{etc.}
\end{equation*}
or

\begin{alignat*}{9}
&S =1 &&+ \omega &&\bigg(\frac{1}{2^2}&&+\frac{1}{3^2}&&+\frac{1}{4^2}&&+\frac{1}{5^2}&&+\text{etc.}\bigg)&&=1&&+\mathfrak{B}\omega \\[3mm]
&  &&- \omega^2 &&\bigg(\frac{1}{2^3}&&+\frac{1}{3^3}&&+\frac{1}{4^3}&&+\frac{1}{5^3}&&+\text{etc.}\bigg)&& &&+\mathfrak{C}\omega^2 \\[3mm]
&  &&- \omega^3 &&\bigg(\frac{1}{2^4}&&+\frac{1}{3^4}&&+\frac{1}{4^4}&&+\frac{1}{5^4}&&+\text{etc.}\bigg)&& &&+\mathfrak{D}\omega^3 \\
& && && && && \text{etc.} && && &&\text{etc.}
\end{alignat*}
Therefore, the total expression, having put $x=1+\omega$, will go over into
\begin{equation*}
\frac{1+\mathfrak{B}\omega -\mathfrak{C}\omega^2 +\mathfrak{D}\omega^3 -\text{etc.}}{\omega(1+\omega)}-\frac{1}{\omega (1+2 \omega)}
\end{equation*}
or

\begin{equation*}
\frac{\omega +\mathfrak{B}\omega +2 \mathfrak{B}\omega^2 - \mathfrak{C}\omega^2 -\text{etc.}}{\omega(1+\omega)(1+2 \omega)}=\frac{1+ \mathfrak{B}+2 \mathfrak{B}\omega - \mathfrak{C}\omega -\text{etc.}}{(1+\omega)(1+2 \omega)}.
\end{equation*}
Now put $\omega =0$ and the propounded value of the expression in the case $x=1$ will be

\begin{equation*}
=1+ \mathfrak{B}=1+\frac{1}{2^2}+\frac{1}{3^2}+\frac{1}{4^2}+\text{etc.};
\end{equation*}
since this series is $=\frac{1}{6}\pi^2$, it follows that the value in question is $=\frac{1}{6}\pi^2$.

\subsection*{Example 2}

\textit{To find the value of this expression}

\begin{equation*}
\frac{2x-xx}{(x-1)^2}+\frac{\pi\pi x}{6(x-1)}-\frac{(2x-1)(1+\frac{1}{2}+\frac{1}{3}+\cdots +\frac{1}{x})}{x(x-1)^2}
\end{equation*}
\textit{in the case  $x=1$.}\\[2mm]
Put $1+\frac{1}{2}+\frac{1}{3}+\cdots +\frac{1}{x}=S$ and set $x=1 +\omega$; it will, as we found in the preceding example, be

\begin{equation*}
S=1+\mathfrak{B}\omega -\mathfrak{C}\omega^2 +\mathfrak{D}\omega^3-\text{etc.}
\end{equation*}
while

\begin{alignat*}{9}
&\mathfrak{B}&&= \frac{1}{2^2}+\frac{1}{3^2}+\frac{1}{4^2}+\frac{1}{5^2}+\text{etc.}&&=\frac{1}{6}\pi \pi -1 \\[3mm]
&\mathfrak{C}&&= \frac{1}{2^3}+\frac{1}{3^3}+\frac{1}{4^3}+\frac{1}{5^3}+\text{etc.}\\[3mm]
&\mathfrak{D}&&= \frac{1}{2^4}+\frac{1}{3^4}+\frac{1}{4^4}+\frac{1}{5^4}+\text{etc.}\\
& && \text{etc.}
\end{alignat*}
Therefore, having put $x=1 + \omega$, the propounded expression will take on this form

\begin{equation*}
\frac{1- \omega \omega}{\omega \omega}+\frac{(1+\mathfrak{B})(1+ \omega)}{\omega}-\frac{(1+2 \omega)(1+\mathfrak{B}\omega - \mathfrak{C}\omega^2+ \mathfrak{D}\omega^3-\text{etc.})}{(1+\omega)\omega^2},
\end{equation*}
which, reduced to the common denominator $\omega^2(1+\omega)$, reads

\begin{small}
\begin{equation*}
\frac{1+\omega -\omega^2 -\omega^3+2 \omega^2+\omega^3 +\mathfrak{B}\omega(1+2 \omega +\omega \omega)-1-\mathfrak{B}\omega +\mathfrak{C}\omega^2 -\mathfrak{D}\omega^3 -2 \omega -2 \mathfrak{B}\omega^2 +2 \mathfrak{C}\omega^3 -\text{etc.}}{\omega^2(1+\omega)},
\end{equation*}
\end{small}
which is reduced to this form

\begin{equation*}
\frac{\omega^2 +\mathfrak{C}\omega^2 +\mathfrak{B}\omega^3 +2 \mathfrak{C}\omega^3- \mathfrak{D}\omega^3+ \text{etc.}}{\omega^2(1+\omega)}
\end{equation*}
Now let $\omega=0$ and $1+\mathfrak{C}$ will result. Therefore, the value of the propounded expression in the case $x=1$ will be $=1+\mathfrak{C}$ and hence will be expressed by  this series

\begin{equation*}
1+\frac{1}{2^3}+\frac{1}{3^3}+\frac{1}{4^3}+\frac{1}{5^3}+\text{etc.};
\end{equation*}
since its sum can be exhibited neither  by means of logarithms nor the circumference of the circle $\pi$, the value in question can still not be assigned  finitely. From these two examples the use  the differentiation of inexplicable functions can have in the doctrine of series is seen sufficiently clearly.

\paragraph*{§386}

In the method to differentiate inexplicable functions treated here we assumed that the infinitesimal terms of the series $A$, $B$, $C$, $D$, $E$ etc. are either $=0$ or have  vanishing differences eventually; if both is not the case, this method can not be used. Therefore, I will explain another method not restricted by this condition and yielding the general summation of series derived from the general term and explained in more detail above [chap. V]. Therefore, let the letters $\mathfrak{A}$, $\mathfrak{B}$, $\mathfrak{C}$, $\mathfrak{D}$, $\mathfrak{E}$ etc. denote the Bernoulli numbers exhibited in § 122 and let this inexplicable function be propounded

\begin{alignat*}{9}
&     &&1   &&2   && 3  && 4 \quad \cdots  &&x \\
&S~=~ &&A~+~&&B~+~&&C~+~&&D~+~ \cdots ~+~  &&X, 
\end{alignat*}
and since we showed above (§ 130) that it will be

\begin{equation*}
S= \int Xdx +\frac{1}{2}X+\frac{\mathfrak{A}dX}{1 \cdot 2 dx}-\frac{\mathfrak{B}d^3 X}{1 \cdot 2 \cdot 3 \cdot 4 dx^3}+\frac{\mathfrak{C}d^5X}{1 \cdot 2 \cdot 3 \cdot 4 \cdot 5 \cdot 6 dx^5}-\text{etc.},
\end{equation*}
it will therefore be easy to exhibit the differential of the function $S$; for, it will be

\begin{equation*}
dS = Xdx+\frac{1}{2}dX+\frac{\mathfrak{A}ddX}{1 \cdot 2 dx}-\frac{\mathfrak{B}d^4 X}{1 \cdot 2 \cdot 3 \cdot 4 dx^3}+\frac{\mathfrak{C}d^6 X}{1 \cdot 2 \cdot 3 \cdot 4 \cdot 5 \cdot 6 dx^5}-\text{etc.}
\end{equation*}

\paragraph*{§387}

But if the propounded progression is connected to the geometric series, in which case its infinitesimal terms are never reduced to constant differences and therefore the first method can not be applied, the method treated in § 174 provides us with the solution. For, if this function is propounded

\begin{equation*}
S=Ap+Bp^2+Cp^3+Dp^4+\cdots+Xp^x,
\end{equation*}
find the values of the letters $\alpha$, $\beta$, $\gamma$, $\delta$ etc. that

\begin{equation*}
\frac{p-1}{p-e^u}=1+\alpha u+ \beta u^2+\gamma u^3 +\delta u^4+\varepsilon u^5+\text{etc.},
\end{equation*}
having found which, as we exhibited them in § 173, it will be

\begin{equation*}
S=\frac{p}{p-1}\cdot p^x \left(X-\frac{\alpha dX}{dx}+\frac{\beta ddX}{dx^2}-\frac{\gamma d^3 X}{dx^3}+\frac{\delta d^4 X}{dx^4}-\text{etc.}\right)
\end{equation*}
$\pm$ constant which renders the sum $=0$, if one puts $x=0$, or which satisfies any other case. Therefore, having taken the differential, this constant will go out of the computation and it will be

\begin{alignat*}{9}
&dS &&= \frac{p}{p-1}\cdot p^x dx \log p \left(X-\frac{\alpha dX}{dx}+\frac{\beta ddX}{dx^2}-\frac{\gamma d^3 X}{dx^3}+\text{etc.}\right)\\[3mm]
& &&+\frac{p}{p-1}\cdot p^x \left(dX-\frac{\alpha ddX}{dx}+\frac{\beta d^3 X}{dx^2}-\frac{\gamma d^4 X}{dx^3}+\text{etc.}\right)
\end{alignat*}
or

\begin{equation*}
dS= \frac{p^{x+1}}{p-1}\left(Xdx \log p-(\alpha \log p-1)dX+(\beta \log p -\alpha)\frac{ddX}{dx}-(\gamma \log p- \beta)\frac{d^3 X}{dx^2}+\text{etc.}\right),
\end{equation*}
which is the differential in question of the propounded function $S$.

\paragraph*{§388}

But if the propounded inexplicable function consists of factors, whether their infinitesimal logarithms have constant differences or not,  applying this method the differential of the function can always be exhibited. For, let 

\begin{alignat*}{9}
&    &&1         &&2        &&3          && 4~ \cdots &&x \\
&S~=~&&A~\cdot ~ &&B~\cdot ~&& C~\cdot ~&&D~\cdots &&X.
\end{alignat*}
Since 

\begin{equation*}
\log S = \log A+ \log B+ \log C+ \log D +\cdots +\log X,
\end{equation*}
using the above method involving the Bernoulli numbers it will be

\begin{equation*}
\log S = \int dx \log X+\frac{1}{2}\log X+\frac{\mathfrak{A}d. \log X}{1 \cdot 2 dx}-\frac{\mathfrak{B}d^3. \log X}{1 \cdot 2 \cdot 3 \cdot 4 dx^3}+\text{etc.},
\end{equation*}
having differentiated which expression 

\begin{equation*}
\frac{dS}{S}= dx \log X +\frac{1}{2}d. \log X +\frac{\mathfrak{A}dd.\log X}{1 \cdot 2 dx}-\frac{\mathfrak{B}d^4. \log X}{1 \cdot 2 \cdot 3 \cdot 4 dx^3}
\end{equation*}
\begin{equation*}
+\frac{\mathfrak{C}d^6. \log X}{1 \cdot 2 \cdot 3 \cdot 4 \cdot 5 \cdot 6 dx^5}-\frac{\mathfrak{D}d^8. \log X}{1 \cdot 2 \cdot 3 \cdot 4 \cdot 5 \cdot 6 \cdot 7 \cdot 8 dx^7}+\text{etc.}
\end{equation*}
Hence, if it was $X=x$ so that

\begin{equation*}
S=1 \cdot 2 \cdot 3 \cdot 4 \cdots x,
\end{equation*}
after the application it will be

\begin{equation*}
\frac{dS}{S}=dx \log x +\frac{dx}{2x}-\frac{\mathfrak{A}dx}{2xx}+\frac{\mathfrak{B}dx}{4x^4}-\frac{\mathfrak{C}dx}{6x^6}+\frac{\mathfrak{D}dx}{8x^8}-\text{etc.},
\end{equation*}
which form, if $x$ was a very large number, is  applied more conveniently than the one we found before.


\title{On the use of Differential Calculus in the Resolution of fractions
\footnote{
Original title: "'De Usu Calculi differentialis in Resolutione Fractionum"', first published as part of the book \textit{Institutiones calculi differentialis cum eius usu in analysi finitorum ac doctrina serierum}, 1755, reprinted in\textit{Opera Omnia}: Series 1, Volume 10, pp. 648- 676, Eneström-Number E212, translated by: Alexander Aycock for the Euler-Kreis Mainz}}
\author{Leonhard Euler}
\date{ }
\maketitle

\paragraph*{§403}

The method to resolve any propounded fraction into simple fractions which we explained in the \textit{Introductio}, even though it is per se simple enough,  can nevertheless  be simplified by applying differential calculus such that often the same task is completed in a much shorter calculation. Especially, if the denominator of the fraction to be resolved was of indefinite degree, the application  of the method explained before  is impeded a lot in most cases, since one has to find all factors of the denominator explicitly. But especially in these cases the division of the denominator by an already found factor becomes too cumbersome. This operation, if differential calculus is applied, can be avoided, i.e. it is not necessary to know the other factor of the denominator, which results, if the denominator is divided by the known one. For this the method  to determine the value of the fraction whose denominator and numerator vanish in a certain case can be applied; therefore, we want to teach in this chapter, how by means of that method the resolution of fractions treated above can be rendered more convenient and tractable,  and at the same time want to finish this book in which we explained the use of differential calculus in  Analysis.

\paragraph*{§404}

Therefore, if an arbitrary fraction $\frac{P}{Q}$ was propounded whose numerator and denominator are polynomial functions of the variable quantity $x$, one at first has to check, whether $x$  has as many or more dimensions in the numerator $P$ than in the denominator $Q$. If this happens, the fraction $\frac{P}{Q}$  will contain an integer part of this form $A+Bx+Cx^2+\text{etc.}$ which can be found by  division; the remaining part will be a fraction with the same denominator $Q$ but whose numerator will be a function, say $R$, containing less dimensions of $x$ than the denominator $Q$ such that  the further resolution only concerns the function $R$. Nevertheless, it is not necessary to know this new numerator $R$, but the same simple fractions  the fraction $\frac{R}{Q}$ would have yielded can be found immediately from the propounded $\frac{P}{Q}$, as we already remarked above.

\paragraph*{§405}

Therefore, except for the integral part, if the fraction $\frac{P}{Q}$ contains one, one has to find the simple fractions whose denominators are either binomials of this form $f+gx$ or trinomials of this kind $f+2x \cos \phi \cdot \sqrt{fg}+gxx$ or squares or cubes or higher powers of formulas of this kind. And these denominators will be all factors of the denominator $Q$ such that any arbitrary factor of this denominator $Q$ yields a simple fraction. If the denominator $Q$ has the factor $f+gx$, from it a simple fraction of this kind will result 

\begin{equation*}
\frac{\mathfrak{A}}{f+gx},
\end{equation*}
but if the factor was $(f+gx)^2$, the two fractions

\begin{equation*}
\frac{\mathfrak{A}}{(f+gx)^2}+\frac{\mathfrak{B}}{f+gx}.
\end{equation*}
And from a cubic factor $(f+gx)^3$ of the denominator $Q$ three simple factors of this form will result

\begin{equation*}
\frac{\mathfrak{A}}{(f+gx)^3}+\frac{\mathfrak{B}}{(f+gx)^2}+\frac{\mathfrak{C}}{f+gx}
\end{equation*}
and so forth. But if the denominator $Q$ had a trinomial factor of this kind $ff-2fg x \cos \phi +gg xx$, from it a simple function of such a form will result

\begin{equation*}
\frac{\mathfrak{A}+\mathfrak{a}x}{ff-2fg x \cos \phi +ggxx},
\end{equation*}
and if two factors of this kind were equal so that we have the factor $(ff-2fg x \cos \phi +gg xx)^2$, hence these two fractions will result

\begin{equation*}
\frac{\mathfrak{A}+\mathfrak{a}x}{(ff-2fg x \cos \varphi +ggxx)^2}+\frac{\mathfrak{B}+\mathfrak{b}x}{ff-2fgx \cos \varphi +gg xx}.
\end{equation*}
But a cubic factor of this kind $(ff-2fg x \cos \phi +gg xx)^3$ will give three simple fractions, a fourth power four and so forth.

\paragraph*{§406}

Therefore,  in order to resolve the fraction $\frac{P}{Q}$ proceed as follows. At first, find all simple either binomial or trinomial factors of the denominator $Q$, and if they were equal to each other, consider them as a single one. Then from each factor of the denominator find the simple fractions either by the method already shown above or the one we will explain here and which  can be applied instead of the first. Having done this, the aggregate of all these fractions together with the integral part, if the propounded fraction $\frac{P}{Q}$ contains one, will be equal to the value of the fraction. We assume the factors of the denominator $Q$ to be known here, since its depends on the resolution of the equation $Q=0$, and we will describe the method here  to define the simple fraction resulting from any given factor of the denominator by means of differential calculus. This, because one already has the denominators of these simple fractions, will be achieved, if we teach to investigate the numerator of each fraction.

\paragraph*{§407}

Therefore, let us put that the denominator $Q$ of the fraction $\frac{P}{Q}$ has the factor $f+gx$ such that  $Q=(f+gx)S$ and this other factor $S$ does not contain the same factor $f+gx$. Let the simple fraction resulting from this factor be $=\frac{\mathfrak{A}}{f+gx}$ and the complement will have the form $\frac{V}{S}$ such that 

\begin{equation*}
\frac{\mathfrak{A}}{f+gx}+\frac{V}{S}=\frac{P}{Q}.
\end{equation*} 
Therefore, it will be

\begin{equation*}
\frac{V}{S}=\frac{P}{Q}-\frac{\mathfrak{A}}{f+gx}=\frac{P- \mathfrak{A}S}{(f+gx)S}
\end{equation*}
and hence

\begin{equation*}
V=\frac{P- \mathfrak{A}S}{f+gx}.
\end{equation*}
Therefore, since $V$ is a polynomial function of $x$, it is necessary that $P- \mathfrak{A}S$ is divisible by $f+gx$; and therefore, if one puts $f+gx=0$ or $x=\frac{-f}{g}$, the expression $P-\mathfrak{A}S$ will vanish. Therefore, let $x=\frac{-f}{g}$, and because $P-\mathfrak{A}S=0$, it will be $\mathfrak{A}=\frac{P}{S}$, as we found already above. But because  $S=\frac{Q}{f+gx}$, it will be

\begin{equation*}
\mathfrak{A}=\frac{(f+gx)P}{Q},
\end{equation*}
if one puts $f+gx=0$ or $x=\frac{-f}{g}$ everywhere. Because in this case so the numerator $(f+gx)P$ as the denominator $Q$ vanish, applying the results we explained on the investigation of the value of fractions of this kind, it will be

\begin{equation*}
\mathfrak{A}=\frac{(f+gx)dP+Pgdx}{dQ},
\end{equation*} 
if one puts $x=\frac{-f}{g}$. In this case, because of $(f+gx)dP=0$, it will be

\begin{equation*}
\mathfrak{A}=\frac{gPdx}{dQ}
\end{equation*}
and   the value of the numerator $\mathfrak{A}$ will be found conveniently by means of differentiation.

\paragraph*{§408}

Therefore, if the denominator $Q$ of the propounded fraction $\frac{P}{Q}$ has the simple factor $f+gx$, from it this simple fraction will result

\begin{equation*}
\frac{\mathfrak{A}}{f+gx}
\end{equation*}
while $\mathfrak{A}=\frac{gPdx}{dQ}$, after here  the value $\frac{-f}{g}$ that has to result from $f+gx=0$ was substituted for $x$ everywhere. Therefore, this way it is not necessary to find the other factor $S$ of the denominator $Q$, which factor results, if $Q$ is divided by $f+gx$. Hence, if $Q$ is not expressed in factors, we can often omit this rather cumbersome division, especially if $x$  has indefinite exponents in the denominator $Q$, since the value of $A$ is obtained from the formula $\frac{gPdx}{dQ}$. But if the denominator $Q$ was already expressed in factors such that hence the value of $S$ is immediately plain, then rather the other expression should applied, by means of which we found $\mathfrak{A}=\frac{P}{S}$ putting $x=\frac{-f}{g}$ everywhere in the same way. And hence in each case one can apply that formula to find the value of $\mathfrak{A}$ which seems  more convenient. But we will illustrate the application of the new formula in some examples.

\subsection*{Example 1}

\textit{Let this fraction be propounded $\frac{x^9}{1+x^{17}}$ whose simple fraction resulting from the factor $1+x$  is to be defined.}\\[2mm]
Since here  $Q=1+x^{17}$, even if its factor $1+x$ is known, nevertheless, if, as the first method requires, we wanted to divide by it, it would result

\begin{equation*}
S=1-x+xx-x^3+\cdots +x^{16}. 
\end{equation*} 
Therefore, we will rather use the new formula $\mathfrak{A}=\frac{gPdx}{dQ}$; therefore, since  $f=1$, $g=1$ and $P=x^9$, because of $dQ=17x^{16}dx$, it will be $\mathfrak{A}=\frac{x^9}{17x^{16}}=\frac{1}{17x^7}$ for $x=-1$, whence  $\mathfrak{A}=-\frac{1}{17}$, and the simple fraction resulting from the factor $1+x$ of the denominator will be

\begin{equation*}
\frac{-1}{17(1+x)}.
\end{equation*}

\subsection*{Example 2}

\textit{Having propounded the fraction $\frac{x^m}{1-x^{2n}}$, to investigate the simple fraction resulting from the factor $1-x$.}\\[2mm]
Because of the propounded factor $1-x$, it will be $f=1$ and $g=-1$. But the denominator $Q=1-x^{2n}$ gives $dQ=-2nx^{2n-1}dx$, whence, because of $P=x^m$, one will obtain $\mathfrak{A}=\frac{-x^m}{-2nx^{2n-1}}$. And, from the equation $1-x=0$, having put $x=1$, it will be $\mathfrak{A}=\frac{1}{2n}$ such that the simple fraction will be
\begin{equation*}
\frac{1}{2n(1-x)}.
\end{equation*}

\subsection*{Example 3}

\textit{Having propounded the fraction $\frac{x^m}{1-4x^k+3x^n}$, to determine the simple fraction resulting from the factor $1-x$.}\\[2mm]
Therefore, here  $f=1$ and $g=-1$, $P=x^m$, $Q=1-4x^k+3x^n$ and $\frac{dQ}{dx}=-4 kx^{k-1}+3nx^{n-1}$; therefore,  $\mathfrak{A}=\frac{-x^m}{-4kx^{k-1}+3nx^{n-1}}$ and for $x=1$ it will be $\mathfrak{A}=\frac{1}{4k-3n}$. Therefore, the simple fraction resulting from this simple factor $1-x$ of the denominator will be

\begin{equation*}
\frac{1}{(4k-3n)(1-x)}.
\end{equation*}

\paragraph*{§409}

Let us now put that the denominator $Q$ of the fraction $\frac{P}{Q}$ has the quadratic factor $(f+gx)^2$ and the simple fraction resulting from this is

\begin{equation*}
\frac{\mathfrak{A}}{(f+gx)^2}+\frac{\mathfrak{B}}{f+gx}.
\end{equation*}
Let $Q=(f+gx)^2S$ and the complement $=\frac{V}{S}$ such that

\begin{equation*}
\frac{V}{S}=\frac{P}{Q}-\frac{\mathfrak{A}}{(f+gx)^2}-\frac{\mathfrak{B}}{f+gx} \quad \text{and} \quad V=\frac{P- \mathfrak{A}S-\mathfrak{B}(f+gx)S}{(f+gx)^2}.
\end{equation*}
Since now $V$ is an integer function, it is necessary that  $P- \mathfrak{A}S- \mathfrak{B}S(f+gx)$ is divisible by $(f+gx)^2$; and because $S$ does contain the factor $f+gx$, the expression $\frac{P}{S}-\mathfrak{A}-\mathfrak{B}(f+gx)$ will also be divisible by $(f+gx)^2$ and hence, having put $f+gx=0$ or $x=\frac{-f}{g}$, not only itself but also its differential $d. \frac{P}{S}-\mathfrak{B}gdx$ will vanish. Therefore, let $x=\frac{-f}{g}$ and from the first equation it will be $\mathfrak{A}=\frac{P}{S}$, from the second on the other hand it will be $\mathfrak{B}=\frac{1}{gdx}d. \frac{P}{S}$; having found these values, one will have the fractions in question

\begin{equation*}
\frac{\mathfrak{A}}{(f+gx)^2}+\frac{\mathfrak{B}}{f+gx}.
\end{equation*}

\subsection*{Example}

\textit{Having propounded the fraction $\frac{x^m}{1-4x^3+3x^4}$, whose numerator has the factor $(1-x)^2$, to find the simple fraction to result from it.}\\[2mm]
Since here  $f=1$, $g=-1$, $P=x^m$ and $Q=1-4x^3+3x^4$, it will be $S=1-2x+3xx$,

\begin{equation*}
\frac{P}{S}=\frac{x^m}{1+2x+3xx} \quad \text{and} \quad d. \frac{P}{S}=\frac{mx^{m-1}dx+2(m-1)x^mdx+3(m-2)x^{m+1}dx}{(1+2x+3xx)^2}.
\end{equation*}
Therefore, having put $x=1$, it will be

\begin{equation*}
\mathfrak{A}=\frac{1}{6} \quad \text{and} \quad \mathfrak{B}=-1 \cdot \frac{6m-8}{36}=\frac{4-3m}{18};
\end{equation*}
therefore, the fractions in question will be

\begin{equation*}
\frac{1}{6(1-x)^2}+\frac{4-3m}{18(1-x)}.
\end{equation*}

\paragraph*{§410}

Let the denominator $Q$ of the fraction $\frac{P}{Q}$ have three equal simple factors or let $Q=(f+gx)^3S$ and let  the simple fractions to result from this cubic factor $(f+gx)^3$ be

\begin{equation*}
\frac{\mathfrak{A}}{(f+gx)^3}+\frac{\mathfrak{B}}{(f+gx)^2}+\frac{\mathfrak{C}}{f+gx};
\end{equation*}
but let the complement of these fractions necessary to constitute the propounded fraction $\frac{P}{Q}$ be $\frac{V}{S}$ and it will be

\begin{equation*}
V=\frac{P- \mathfrak{A}S-\mathfrak{B}S(f+gx)-\mathfrak{C}S(f+gx)^2}{(f+gx)^3}.
\end{equation*}
Hence this expression $\frac{P}{S}-\mathfrak{A}-\mathfrak{B}(f+gx)-\mathfrak{C}(f+gx)^2$ will be divisible by $(f+gx)^3$; therefore, having put $f+gx=0$ or $x=\frac{-f}{g}$, not only this expression itself but also its first and second differential will become $=0$. By putting $x=\frac{-f}{g}$ it will be

\begin{equation*}
\frac{P}{S}-\mathfrak{A}-\mathfrak{B}(f+gx)-\mathfrak{C}(f+gx)^2=0
\end{equation*}
\begin{equation*}
d. \frac{P}{S}-\mathfrak{B}gdx-2 \mathfrak{C}gdx(f+gx)=0
\end{equation*}
\begin{equation*}
dd. \frac{P}{S}-2 \mathfrak{C}g^2 dx^2=0.
\end{equation*}
From the first equation, it will therefore be

\begin{equation*}
\mathfrak{A}=\frac{P}{S}.
\end{equation*}
From the second on the other hand it will be

\begin{equation*}
\mathfrak{B}=\frac{1}{gdx}d. \frac{P}{S}.
\end{equation*}
Finally, from the third one defines

\begin{equation*}
\mathfrak{C}=\frac{1}{2g^2 dx^2}dd. \frac{P}{S}.
\end{equation*}

\paragraph*{§411}

Therefore, in general, if the denominator $Q$ of the fraction $\frac{P}{Q}$ has the factor $(f+gx)^n$ such that  $Q=(f+gx)^nS$, having put the simple fractions to result from this factor $(f+gx)^n$

\begin{equation*}
\frac{\mathfrak{A}}{(f+gx)^n}+\frac{\mathfrak{B}}{(f+gx)^{n-1}}+\frac{\mathfrak{C}}{(f+gx)^{n-2}}+\frac{\mathfrak{D}}{(f+gx)^{n-3}}+\frac{\mathfrak{E}}{(f+gx)^{n-4}}+\text{etc.},
\end{equation*}
until  the last whose denominator is $f+gx$ is reached, if one reasons exactly as before, one will find that this expression

\begin{equation*}
\frac{P}{S}-\mathfrak{A}-\mathfrak{B}(f+gx)-\mathfrak{C}(f+gx)^2-\mathfrak{D}(f+gx)^3-\mathfrak{E}(f+gx)^4-\text{etc.}
\end{equation*} 
must be divisible by $(f+gx)^n$; therefore, so itself as each differential of it up to degree $n-1$ will have to vanish in the case $x=\frac{-f}{g}$. Therefore,  putting $x=\frac{-f}{g}$ everywhere, from these equations one concludes that

\begin{alignat*}{9}
&\mathfrak{A}&&=\frac{P}{S} \\[3mm]
&\mathfrak{B}&&=\frac{1}{1g dx}d. \frac{P}{S} \\[3mm]
&\mathfrak{C}&&= \frac{1}{1 \cdot 2 g^2 dx^2}dd. \frac{P}{S}\\[3mm]
&\mathfrak{D}&&= \frac{1}{1 \cdot 2 \cdot 3 g^3 dx^3}d^3. \frac{P}{S}\\[3mm]
&\mathfrak{E}&&= \frac{1}{1 \cdot 2 \cdot 3 \cdot 4 g^4 dx^4}d^4. \frac{P}{S}\\[3mm]
& &&\text{etc.}
\end{alignat*}
Here it is to be noted that these differentials of $\frac{P}{S}$ must be taken before  one substitutes $\frac{-f}{g}$ for $x$; for, otherwise the variability of $x$ would be lost.

\paragraph*{§412}

Therefore, this way these numerators $\mathfrak{A}$, $\mathfrak{B}$, $\mathfrak{C}$, $\mathfrak{D}$ etc. will be expressed more easily than by the method given in the \textit{Introductio} and often  their values are also found a lot faster by means of this new method. In order to make this comparison more easily, let us define the values of the letters  $\mathfrak{A}$, $\mathfrak{B}$, $\mathfrak{C}$, $\mathfrak{D}$ etc. by means of the first method.

\begin{equation*}
\renewcommand{\arraystretch}{2,5}
\setlength{\arraycolsep}{0.0mm}
\begin{array}{ll|l}
 & \text{Having put} \quad x=\dfrac{-f}{g} \qquad & \quad \text{Still considering} ~ x ~ \text{as a variable, set} \\
\text{it will be} \quad & \qquad \qquad \mathfrak{A}=\dfrac{P}{S} \quad & \qquad \qquad \dfrac{P- \mathfrak{A}S}{f+gx}= \mathfrak{P}, \\ 
\text{it will be} \quad & \qquad \qquad \mathfrak{B}=\dfrac{\mathfrak{P}}{S} \quad & \qquad \qquad \dfrac{\mathfrak{P}- \mathfrak{B}S}{f+gx}= \mathfrak{Q}, \\  
\text{it will be} \quad & \qquad \qquad \mathfrak{C}=\dfrac{\mathfrak{Q}}{S} \quad & \qquad \qquad \dfrac{\mathfrak{Q}- \mathfrak{C}S}{f+gx}= \mathfrak{R}, \\
\text{it will be} \quad & \qquad \qquad \mathfrak{D}=\dfrac{\mathfrak{R}}{S} \quad & \qquad \qquad \dfrac{\mathfrak{R}- \mathfrak{D}S}{f+gx}= \mathfrak{S}, \\
\text{it will be} \quad & \qquad \qquad \mathfrak{E}=\dfrac{\mathfrak{S}}{S} \quad & \qquad \qquad \text{and so forth}. 
\end{array}
\end{equation*} 

\paragraph*{§413}

But if the denominator $Q$ of the fraction $\frac{P}{Q}$  has not only simple real factors, take two imaginary ones  whose product will be real then. Therefore, let the factor of the denominator $Q$ be $ff-2fg x \cos \varphi +ggxx$, which put $=0$, gives these two imaginary values

\begin{equation*}
x=\frac{f}{g}\cos \varphi \pm \frac{f}{g\sqrt{-1}}\sin \varphi;
\end{equation*}
therefore, it will be

\begin{equation*}
x^n = \frac{f^n}{g^n}\cos n \varphi \pm \frac{f^n}{g^n \sqrt{-1}}\sin n \varphi.
\end{equation*}
Let us put that $Q=(ff-2fgx \cos \varphi +gg xx)S$ and additionally $S$ is not divisible by $ff-2fg x \cos \varphi +gg xx$. Let the fraction to result from this factor be 

\begin{equation*}
\frac{\mathfrak{A}+\mathfrak{a}x}{ff-2fg x \cos \varphi +gg xx}
\end{equation*} 
and let the complement necessary to get to the propounded fraction $\frac{P}{Q}$ be $=\frac{V}{S}$; it will be

\begin{equation*}
V=\frac{P-(\mathfrak{A}+\mathfrak{a}x)S}{ff-2fgx \cos \varphi +ggxx},
\end{equation*}
whence $P-(\mathfrak{A}+\mathfrak{a}x)S$ and therefore also $\frac{P}{S}- \mathfrak{A}- \mathfrak{a}x$ will be divisible by $ff-2fg x \cos \varphi +ggxx =0$, i.e., if one puts either

\begin{equation*}
x=\frac{f}{g}\cos \varphi +\frac{f}{g \sqrt{-1}}\sin \varphi
\end{equation*}
or 

\begin{equation*}
x=\frac{f}{g}\cos \varphi -\frac{f}{g \sqrt{-1}}\sin \varphi.
\end{equation*}

\paragraph*{§414}

Since $P$ and $S$ are polynomial functions of $x$, make both substitutions separately in both expressions; and since for any power of $x$, say $x^n$, this binomial

\begin{equation*}
x^n=\frac{f^n}{g^n}\cos n\varphi \pm \frac{f^n}{g^n \sqrt{-1}}\sin n \varphi
\end{equation*} 
has to be substituted, first let us put $\frac{f^n}{g^n}\cos n \varphi$ for $x^n$ everywhere and, having done this, let $P$ go over into $\mathfrak{P}$ and $S$ into $\mathfrak{S}$. Further, put $\frac{f^n}{g^n}\sin n \varphi$ for $x^n$ everywhere and, having done this, let $P$ go over into $\mathfrak{p}$ and $S$ into $\mathfrak{s}$; here it is to be noted that before these substitutions both functions $P$ and $S$ have to be expanded completely such that, if  they are contained in factors, one has to get rid of these factors by actual multiplication.  Having found these values $\mathfrak{P}$, $\mathfrak{p}$, $\mathfrak{S}$, $\mathfrak{s}$, it will be obvious, if one puts $x=\frac{f}{g}\cos \varphi \pm \frac{f}{g \sqrt{-1}}\sin \varphi$, that the function $P$ will go over into $\mathfrak{P}\pm \frac{\mathfrak{p}}{\sqrt{-1}}$ and the function $S$ will go over into $\mathfrak{S}\pm \frac{s}{\sqrt{-1}}$. Therefore, because $\frac{P}{S}-\mathfrak{A}-\mathfrak{a}x$ or $P-(\mathfrak{A}+\mathfrak{a}x)S$ has to vanish in both cases, it will be

\begin{equation*}
\mathfrak{P}\pm \frac{\mathfrak{p}}{\sqrt{-1}}=\left(\mathfrak{A}+\frac{\mathfrak{a}f}{g}\cos \varphi \pm \frac{\mathfrak{a}f}{g \sqrt{-1}}\sin \varphi\right)\left(\mathfrak{S}\pm \frac{\mathfrak{s}}{\sqrt{-1}}\right),
\end{equation*} 
whence, because of the ambiguous signs, these two equations will result

\begin{alignat*}{9}
&\mathfrak{P}&&=\mathfrak{AS}&&+\frac{\mathfrak{a}f\mathfrak{S}}{g}\cos \varphi &&-\frac{\mathfrak{a}f \mathfrak{s}}{\sin \varphi},\\[3mm]
&\mathfrak{p}&&=\mathfrak{As}&&+\frac{\mathfrak{a}f\mathfrak{s}}{g}\cos \varphi &&-\frac{\mathfrak{a}f \mathfrak{S}}{\sin \varphi},\\
\end{alignat*}
from which, eliminating $\mathfrak{A}$, one finds

\begin{equation*}
\mathfrak{Sp}-\mathfrak{sP}=\frac{\mathfrak{a}f(\mathfrak{S}^2+\mathfrak{s}^2)}{g}\sin \varphi;
\end{equation*}
and hence it will be

\begin{equation*}
\mathfrak{a}=\frac{g(\mathfrak{Sp}-\mathfrak{sP})}{f(\mathfrak{S}^2+\mathfrak{s}^2)\sin \varphi}.
\end{equation*}
Further,  eliminating $\sin \varphi$, it will be

\begin{equation*}
\mathfrak{S P}+\mathfrak{sp}=(\mathfrak{S}^2+\mathfrak{s}^2)(\mathfrak{A}+\frac{\mathfrak{a}f}{g}\cos \varphi).
\end{equation*}
Therefore,

\begin{equation*}
\mathfrak{A}=\frac{\mathfrak{SP}+\mathfrak{sp}}{\mathfrak{S}^2+\mathfrak{s}^2}-\frac{(\mathfrak{Sp}-\mathfrak{sP}\cos \varphi)}{(\mathfrak{S}^2+\mathfrak{s}^2)\sin \varphi}.
\end{equation*}

\paragraph*{§415}

Since 

\begin{equation*}
S=\frac{Q}{ff-2fgx \cos \varphi +gg xx}
\end{equation*}
and since, having put 

\begin{equation*}
ff-2fgx \cos \varphi +ggxx=0,
\end{equation*}
so the numerator as the denominator will vanish,  in this case it will be

\begin{equation*}
S=\frac{dQ:dx}{2ggx-2fg \cos \varphi}.
\end{equation*}
Now, let us put, if one substitutes $x^n=\frac{f^n}{g^n}\cos n \varphi$ everywhere, that the function $\frac{dQ}{dx}$ goes over into $\mathfrak{Q}$, but if one sets $x^n=\frac{f^n}{g^n}\sin n \varphi$, that it goes over into $\mathfrak{q}$; and it is obvious, if one puts $x=\frac{f}{g}\cos \varphi \pm \frac{f}{g \sqrt{-1}}\sin \varphi$, that the function $\frac{dQ}{dx}$ goes over into $\mathfrak{Q}\pm \frac{\mathfrak{q}}{\sqrt{-1}}$. From this the function $S$ will go over into

\begin{equation*}
\frac{\mathfrak{Q}\pm \mathfrak{q}:\sqrt{-1}}{\pm 2 fg \sin \varphi :\sqrt{-1}}.
\end{equation*}
Therefore, since  $S=\mathfrak{S}\pm \frac{\mathfrak{s}}{\sqrt{-1}}$ having substituted the same value for $x$, one will have

\begin{equation*}
\mathfrak{Q}\pm \frac{\mathfrak{q}}{\sqrt{-1}}=\pm \frac{2fg \mathfrak{S}}{\sqrt{-1}}\sin \varphi - 2fg \mathfrak{s}\sin \varphi.
\end{equation*}
Therefore, it will be

\begin{equation*}
\mathfrak{s}=\frac{-\mathfrak{Q}}{2fg \sin \varphi} \quad \text{and} \quad \mathfrak{S}= \frac{\mathfrak{q}}{2fg \sin \varphi}.
\end{equation*}
And having substituted these values, it will be

\begin{equation*}
\mathfrak{a}= \frac{2gg(\mathfrak{pq}+\mathfrak{PQ})}{\mathfrak{Q}^2+\mathfrak{q}^2}
\end{equation*}
and

\begin{equation*}
\mathfrak{A}=\frac{2fg(\mathfrak{Pq}-\mathfrak{pQ})\sin \varphi}{\mathfrak{Q}^2+\mathfrak{q}^2}-\frac{2fg(\mathfrak{pq}+\mathfrak{PQ})\cos \varphi}{\mathfrak{Q}^2+\mathfrak{q}^2}.
\end{equation*}

\paragraph*{§416}

Therefore, we found an appropriate way  to form  the simple fraction from each factor of second power and  here, since the denominator  of the propounded fraction is retained in the calculation, we avoid the division, by means of which the value of the letter $S$ would have to be defined and which is often very cumbersome. Therefore, if the denominator $Q$ of the fraction $\frac{P}{Q}$  has such a factor $ff-2fgx \cos \varphi +gg xx$, the simple fraction to result from this factor and we assumed to be

\begin{equation*}
=\frac{\mathfrak{A}+\mathfrak{a}x}{ff-2fgx \cos \varphi +ggxx}
\end{equation*}
will be defined the following way. Put $x=\frac{f}{g}\cos \varphi$ and for each power $x^n$ of $x$ write $\frac{f^n}{g^n}\cos n \varphi$; having done this, let $P$ go over into $\mathfrak{P}$ and the function $\frac{dQ}{dx}$ into $\mathfrak{Q}$. Further, put $x=\frac{f}{g}\sin \varphi$ and each of its powers $x^n=\frac{f^n}{g^n}\sin n\varphi$ and let $P$ go over into $\mathfrak{p}$ and $\frac{dQ}{dx}$ into $\mathfrak{q}$. And, having found these values of the letters $\mathfrak{P}$, $\mathfrak{Q}$, $\mathfrak{p}$ and $\mathfrak{q}$ this way, the quantities $\mathfrak{A}$ and $\mathfrak{a}$ will be defined in such a way that 

\begin{equation*}
\mathfrak{A}=\frac{2fg(\mathfrak{Pq}-\mathfrak{pQ}\sin \varphi)}{\mathfrak{Q}^2+\mathfrak{q}^2}-\frac{2fg(\mathfrak{PQ}+\mathfrak{pq})\cos \varphi}{\mathfrak{Q}^2+\mathfrak{q}^2},
\end{equation*}
\begin{equation*}
\mathfrak{a}=\frac{2gg(\mathfrak{PQ}+\mathfrak{pq})}{\mathfrak{Q}^2+\mathfrak{q}^2}.
\end{equation*}
Therefore, the fraction to result from the factor $ff-2fg \cos \varphi +ggxx$ of the denominator $Q$ will be

\begin{equation*}
\frac{2fg (\mathfrak{Pq}-\mathfrak{pQ})\sin \varphi +2g(\mathfrak{PQ}+\mathfrak{pq})(gx- f \cos \varphi)}{(\mathfrak{Q}^2+\mathfrak{q}^2)(ff-2fg x \cos \varphi +ggxx)}.
\end{equation*}

 \subsection*{Example 1}
 
 \textit{If the fraction $\frac{x^m}{a+bx^n}$ was propounded whose denominator $a+bx^n$ has the factor $ff-2fg x \cos \varphi +ggxx$, to find the simple fraction corresponding to this factor.}\\[2mm]
 Since here $P=x^m$ and $Q=a+bx^n$, it will be
 
\begin{equation*}
 \frac{dQ}{dx}=n b x^{n-1},
\end{equation*}
whence it will be

\begin{alignat*}{9}
&\mathfrak{P} &&=\frac{f^m}{g^m}\cos m \varphi, \quad && \mathfrak{p}&&=\frac{f^m}{g^m} \sin m \varphi , \\[3mm]
&\mathfrak{Q} &&=\frac{nb f^{n-1}}{g^{n-1}}\cos (n-1) \varphi, \quad && \mathfrak{q}&&= \frac{n b f^{n-1}}{g^{n-1}}\sin (n-1) \varphi.
\end{alignat*}
From these it will be 

\begin{equation*}
\mathfrak{Q}^2+\mathfrak{q}^2 = \frac{n^2 b^2 f^{2(n-1)}}{g^{2(n-1)}},
\end{equation*}
\begin{equation*}
\mathfrak{Pq}- \mathfrak{pQ}= \frac{nbf^{m+n-1}}{g^{m+n-1}}\sin (n-m-1)\varphi
\end{equation*}
and
\begin{equation*}
\mathfrak{PQ}+\mathfrak{pq}=\frac{nbf^{m+n-1}}{g^{m+n-1}}\cos (n-m-1)\varphi .
\end{equation*}
Therefore, the simple fraction in question will be

\begin{equation*}
\frac{2g^{n-m}(f \sin \varphi \cdot \sin (n-m-1)\varphi+gx \cos(n-m-1)\varphi -f \cos \varphi \cdot \cos (n-m-1)\varphi)}{nbf^{n-m-1}(ff-2fg x \cos \varphi +ggxx)}
\end{equation*}
or

\begin{equation*}
\frac{2g^{n-m}(gx\cos(n-m-1)\varphi -f \cos(n-m)\varphi)}{nbf^{n-m-1}(ff-2fgx \cos \varphi +ggxx)}.
\end{equation*}

\subsection*{Example 2}

\textit{Let the fraction $\frac{1}{x^m(a+bx^n)}$ be propounded whose denominator has the factor $ff-2fg x \cos \varphi +ggxx$; to find the simple fraction to result from this.}\\[2mm]
Since $P=1$ and $Q=ax^m+bx^{m+n}$, it will be

\begin{equation*}
\frac{dQ}{dx}= ma x^{m-1}+(m+n)bx^{m+n-1}
\end{equation*}
and hence, having put $x^n=\frac{f^n}{g^n}\cos n \varphi$, because of $P=x^0$, $\mathfrak{P}=1$, it will be

\begin{equation*}
\mathfrak{Q}= \frac{maf^{m-1}}{g^{m-1}}\cos (m-1)\varphi + \frac{(m+n)bf^{m+n-1}}{g^{m+n-1}}\cos (m+n-1)\varphi
\end{equation*}
and, having put $x^n=\frac{f^n}{g^n}\sin n \varphi$, it will be $\mathfrak{p}=0$ and  

\begin{equation*}
\mathfrak{q}=\frac{ma f^{m-1}}{g^{m-1}}\sin (m-1) \varphi + \frac{(m+n)b f^{m+n-1}}{g^{m+n-1}}\sin (m+n-1)\varphi.
\end{equation*}
Therefore,

\begin{equation*}
\mathfrak{Q}^2+\mathfrak{q}^2= \frac{m^2 a^2 f^{2(m-1)}}{g^2(m-1)}+\frac{2m(m+n)abf^{2m+n-2}}{g^{2m+n-2}}\cos n \varphi +\frac{(m+n)^2b^2 f^{2(m+n-1)}}{g^{2(m+n-1)}}.
\end{equation*}
If  $ff-2fgx \cos \varphi +ggxx$ is a divisor of $a+bx^n$, it will be

\begin{equation*}
a+\frac{bf^n}{g^n}\cos n \varphi = 0 \quad \text{and} \quad \frac{bf^n}{g^n}\sin n \varphi =0, \quad \text{whence} \quad aa=\frac{bbf^{2n}}{g^{2n}}.
\end{equation*}
Therefore, it will be

\begin{equation*}
\mathfrak{Q}^2+\mathfrak{q}^2= \frac{(m+n)^2 bb f^{2(m+n-1)}}{g^{2(m+n-1)}}- \frac{m(2n+m)aaf^{2(m-1)}}{g^{2(m-1)}}= \frac{nn aa f^{2(m-1)}}{g^{2(m-1)}}= \frac{nnbb f^{2(m+n-1)}}{g^{2(m+n-1)}}.
\end{equation*}
Further, it will be

\begin{equation*}
\mathfrak{Pq}-\mathfrak{pQ}=\frac{maf^{m-1}}{g^{m-1}}\sin (m-1) \varphi +\frac{(m+n)bf^{m+n-1}}{g^{m+n-1}}\sin (m+n-1)\varphi
\end{equation*}
\begin{equation*}
=\frac{bf^{m+n-1}}{g^{m+n-1}}((m+n)\sin (m+n-1)\varphi -m \cos n \varphi \cdot \sin (m-1)\varphi)
\end{equation*}
\begin{equation*}
=\frac{bf^{m+n-1}}{g^{m+n-1}}(n \cos n \varphi \cdot \sin (m-1) \varphi +(m+n)\sin n \varphi \cdot \cos (m-1) \varphi)
\end{equation*}
and

\begin{equation*}
\mathfrak{PQ}+\mathfrak{pq}=\frac{bf^{m+n-1}}{g^{m+n-1}}((m+n)\cos (m+n-1)\varphi - m \cos n \varphi \cdot \cos (m-1)\varphi).
\end{equation*}
Or because  $ff-2fg \cos \varphi +ggxx$ is also a divisor of $ax^{m-1}+bx^{m+n-1}$, it will be

\begin{equation*}
\frac{af^{m-1}}{g^{m-1}}\cos (m-1) \varphi +\frac{b f^{m+n-1}}{g^{m+n-1}}\cos (m+n-1)\varphi =0
\end{equation*}
and 
\begin{equation*}
\frac{af^{m-1}}{g^{m-1}}\sin (m-1) \varphi +\frac{b f^{m+n-1}}{g^{m+n-1}}\sin (m+n-1)\varphi =0,
\end{equation*}
whence it will be

\begin{equation*}
\mathfrak{Q}=\frac{nbf^{m+n-1}}{g^{m+n-1}}\cos (m+n-1)\varphi \quad \text{and} \quad \mathfrak{q}=\frac{nbf^{m+n-1}}{g^{m+n-1}}\sin (m+n-1)\varphi
\end{equation*}
 or
 
 \begin{equation*}
 \mathfrak{Q}=\frac{-naf^{m-1}}{g^{m-1}}\cos (m-1)\varphi \quad \text{and} \quad \mathfrak{q}=\frac{-naf^{m-1}}{g^{m-1}}\cos (m-1)\varphi.
 \end{equation*}
 From these the fraction in question will result as
 
 \begin{equation*}
 \frac{2g^m(f \cos m \varphi-gx \cos (m-1)\varphi)}{na f^{m-1}(ff-2fg x \cos  \varphi +ggxx)}.
 \end{equation*}
 This formula formula follows from the first example, if one substitutes a negative $m$, whence it would not have been necessary to consider this case.
 
 \subsection*{Example 3}
 
 \textit{If the denominator of this fraction $\frac{x^m}{a+bx^n+cx^{2n}}$ has the factor $ff-2fgx \cos \varphi +ggxx$, to investigate the simple fraction to result from this factor.}\\[2mm]
 If $ff-2fg x \cos \varphi +ggxx$ is a factor of the denominator $a+bx^n+cx^{2n}$, it will, as we showed above, be
 
 \begin{equation*}
 a+\frac{bf^n}{g^n}\cos n \varphi +\frac{c f^{2n}}{g^{2n}}= 0 \quad \text{and} \quad \frac{bf^n}{g^n}\sin n \varphi +\frac{cf^{2n}}{g^{2n}}\sin 2n \varphi =0.
 \end{equation*}
 Therefore, since $P=x^m$ and $Q=a+bx^n+cx^{2n}$, it will be
 
\begin{equation*}
\frac{dQ}{dx}=nbx^{n-1}+2ncx^{2n-1},
\end{equation*} 
whence

\begin{equation*}
\mathfrak{P}=\frac{f^m}{g^m}\cos m \varphi \quad \text{and} \quad \mathfrak{p}=\frac{f^m}{g^m}\sin m \varphi,
\end{equation*}
\begin{alignat*}{9}
& \mathfrak{Q}&&= \frac{nb f^{n-1}}{g^{n-1}}\cos (n-1) \varphi &&+\frac{2n cf^{2n-1}}{g^{2n-1}}\cos (2n-1)\varphi, \\[3mm]
& \mathfrak{q}&&= \frac{nb f^{n-1}}{g^{n-1}}\sin (n-1) \varphi &&+\frac{2n cf^{2n-1}}{g^{2n-1}}\sin (2n-1)\varphi, \\
\end{alignat*}
Therefore, we will have

\begin{equation*}
\mathfrak{Q}^2+\mathfrak{q}^2= \frac{n^2f^{2(n-1)}}{g^{2(n-1)}}\left(bb+\frac{4bcf^n}{g^n}\cos n \varphi+\frac{4ccf^{2n}}{g^{2n}}\right).
\end{equation*}
But from the first two equations 

\begin{equation*}
\frac{f^{2n}}{g^{2n}}\left(bb+\frac{2bcf^n}{g^n}\cos n \varphi + \frac{ccf^{2n}}{g^{2n}}\right)=aa
\end{equation*}
and hence

\begin{equation*}
\frac{4bcf^n}{g^n}\cos n \varphi = \frac{2g^{2n}aa}{f^{2n}}-2bb-\frac{2ccf^{2n}}{g^{2n}};
\end{equation*}
having substituted the value there, it  will be

\begin{equation*}
\mathfrak{Q}^2+\mathfrak{q}^2= \frac{n^2f^{2n-2}}{g^{2n-2}}\left(\frac{2aag^{2n}}{f^{2n}}-bb+\frac{2ccf^{2n}}{g^{2n}}\right)
\end{equation*}
or
\begin{equation*}
\mathfrak{Q}^2+\mathfrak{q}^2=\frac{n^2(2aag^{4n}-bbf^{2n}g^{2n}+2ccf^{4n})}{ffg^{4n-2}}.
\end{equation*}
Further, it will be

\begin{equation*}
\mathfrak{Pq}-\mathfrak{pQ}=\frac{nbf^{m+n-1}}{g^{m+n-1}}\sin(n-m-1)\varphi +\frac{2ncf^{m+2n-1}}{g^{m+2n-1}}\sin (2n-m-1)\varphi,
\end{equation*}
\begin{equation*}
\mathfrak{PQ}+\mathfrak{pq}=\frac{nbf^{m+n-1}}{g^{m+n-1}}\cos(n-m-1)\varphi +\frac{2ncf^{m+2n-1}}{g^{m+2n-1}}\cos (2n-m-1)\varphi.
\end{equation*}
Having found these values, the simple fraction in question will be

\begin{equation*}
\frac{2fg(\mathfrak{Pq}-\mathfrak{pQ})\sin \varphi +2g(\mathfrak{PQ}+\mathfrak{pq})(gx-f \cos \varphi)}{(\mathfrak{Q}^2+\mathfrak{q}^2)(ff-2fgx \cos \varphi +ggxx)}.
\end{equation*}

\paragraph*{§417}

But these fractions will be expressed in an easier way, if we determine the factors of the denominators. Therefore, let the denominator of the propounded fraction be

\begin{equation*}
a+bx^n;
\end{equation*}
if the trinomial factor is put

\begin{equation*}
ff-2fgx \cos \varphi +ggxx,
\end{equation*}
it will be, as we showed in the \textit{Introductio},

\begin{equation*}
a+\frac{bf^n}{g^n}\cos n \varphi =0 \quad \text{and} \quad \frac{bf^n}{g^n}\sin n \varphi =0;
\end{equation*}
therefore, because $\sin n \varphi =0$, it will be either $n \varphi = (2k-1)\pi$ or $n \pi =2 k \pi$; in the first case it will be $\cos n \varphi =-1$, in the second $\cos n \varphi =+1$. Therefore, if $a$ and $b$ are positive quantities, only the first case will hold, in which $a=\frac{bf^n}{g^n}$ and therefore

\begin{equation*}
f=a^{\frac{1}{n}} \quad \text{and} \quad g=b^{\frac{1}{n}}.
\end{equation*}
But instead of these irrational quantities let us still use the letters $f$ and $g$ or let us put $a=f^n$ and $b=g^n$ such the the factors of this function have to be investigated

\begin{equation*}
f^n+g^nx^n.
\end{equation*}
Therefore, since $\varphi = \frac{(2k-1)\pi}{n}$, where $k$ can denote any positive integer, but on the other hand for $k$ no larger numbers than those rendering $\frac{2k-1}{n}$ greater than $1$ are to be taken; therefore, the factors of the propounded fraction $f^n+g^nx^n$ will be the following

\begin{alignat*}{9}
&&ff&&-2fg x \cos \frac{\pi}{n}&&+ggxx \\[3mm]
&&ff&&-2fg x \cos \frac{3\pi}{n}&&+ggxx \\[3mm]
&&ff&&-2fg x \cos \frac{5\pi}{n}&&+ggxx \\[3mm]
&& &&\text{etc.},
\end{alignat*} 
where it is to be noted, if $n$ is an odd number, that one has this one binomial factor

\begin{equation*}
f+gx;
\end{equation*}
but if $n$ is an even number, the product will contain no binomial factor.

\subsection*{Example 1}

\textit{To resolve the fraction $\frac{x^m}{f^n+g^nx^n}$ into its simple fractions.}\\[2mm]
Since a trinomial factor of any arbitrary denominator is contained in this form

\begin{equation*}
ff- 2fg x \cos \frac{(2k-1)\pi}{n}+ggxx,
\end{equation*}
 in  example 1 of the preceding paragraph it will be $a=f^n$, $b=g^n$ and $\varphi= \frac{(2k-1)\pi}{n}$,
whence it will be

\begin{equation*}
\sin (n-m-1)\varphi =\sin (m+1)\varphi =\sin \frac{(m+1)(2k-1)\pi}{n}
\end{equation*} 
and

\begin{equation*}
\cos (n-m-1)\varphi = - \cos (m+1) \varphi = - \cos \frac{(m+1)(2k-1)\pi}{n}.
\end{equation*}
Therefore, from this factor this simple fraction results

\begin{equation*}
\frac{2f \sin \frac{(2k-1)\pi}{n}\cdot \sin \frac{(m+1)(2k-1)\pi}{n}-2 \cos \frac{(m+1)(2k-1)\pi}{n}\left(fx-f \cos \frac{(2k-1)\pi}{n}\right)}{nf^{n-m-1}g^m \left(ff-2fgx \cos \frac{(2k-1)\pi}{n}+ggxx\right)}.
\end{equation*}
Therefore, the propounded fraction will be resolved into these simple ones

\begin{equation*}
\frac{2f \sin \frac{\pi}{n}\cdot \sin \frac{(m+1)\pi}{n}-2 \cos \frac{(m+1)\pi}{n}\left(fx-f \cos \frac{\pi}{n}\right)}{nf^{n-m-1}g^m \left(ff-2fgx \cos \frac{\pi}{n}+ggxx\right)}
\end{equation*}
\begin{equation*}
+\frac{2f \sin \frac{3\pi}{n}\cdot \sin \frac{3(m+1)\pi}{n}-2 \cos \frac{3(m+1)\pi}{n}\left(fx-f \cos \frac{3\pi}{n}\right)}{nf^{n-m-1}g^m \left(ff-2fgx \cos \frac{3\pi}{n}+ggxx\right)}
\end{equation*}
\begin{equation*}
+\frac{2f \sin \frac{5\pi}{n}\cdot \sin \frac{5(m+1)\pi}{n}-2 \cos \frac{5(m+1)\pi}{n}\left(fx-f \cos \frac{5\pi}{n}\right)}{nf^{n-m-1}g^m \left(ff-2fgx \cos \frac{5\pi}{n}+ggxx\right)}
\end{equation*}
\begin{equation*}
\text{etc.}
\end{equation*}
Therefore, if $n$ was an even number, this way all simple fractions result; but if $n$ was an odd number, because of the binomial factor $f+gx$,  to the fractions resulting this way one furthermore has to add this one

\begin{equation*}
\frac{\pm 1}{nf^{n-m-1}g^m(f+gx)},
\end{equation*}
where the sign $+$ holds, if $m$ was an even number, otherwise the sign $-$. If $m$ was a number greater than $m$, then to these fractions additionally integral parts of this kind are added

\begin{equation*}
Ax^{m-n}+Bx^{m-2n}+C^{m-3n}+Dx^{m-4n}+\text{etc.},
\end{equation*}
as long as the exponents remain positive, and it will be

\begin{alignat*}{9}
&    &&   &&Ag^n &&=1 \quad && \text{therefore}  \quad && A &&=+\frac{1}{g^{n}}\\[3mm]
&Af^n&&~+~&&Bg^n &&=0 \quad &&  \quad && B &&=-\frac{f^n}{g^{2n}}\\[3mm]
&Bf^n&&~+~&&Cg^n &&=0 \quad &&  \quad && C &&=+\frac{f^{2n}}{g^{3n}}\\[3mm]
&Cf^n&&~+~&&Dg^n &&=0 \quad &&  \quad && C &&=-\frac{f^{3n}}{g^{4n}}\\
& && &&\text{etc.} && && && &&\text{etc.}
\end{alignat*}

\subsection*{Example 2}

\textit{To resolve the fraction $\frac{1}{x^m(f^n+g^nx^n)}$ into its simple fractions.}\\[2mm]
Concerning the factors of $f^n+g^nx^n$, from them the same fractions   we found in the preceding example result, if only $m$ is taken negatively; therefore, it only remains to define the simple fractions resulting from the other factor $x^m$, what is most conveniently done this way. Set the propounded fraction

\begin{equation*}
=\frac{\mathfrak{A}}{x^m}+\frac{\mathfrak{N}x^{n-m}}{f^n+g^nx^n}
\end{equation*}
and it will be

\begin{alignat*}{9}
&\mathfrak{Af^n} &&=1 \quad \text{therefore} \quad && \mathfrak{A}&&=+\frac{1}{f^n} \\[3mm]
&\mathfrak{A}g^n+\mathfrak{N} &&=0 \quad \quad && \mathfrak{N}&&=-\frac{g^n}{f^n}. \\
\end{alignat*}
If $n-m$ was a negative number, one will have to proceed in like manner, such that, if $m$ was an arbitrary large number, simple fractions of this kind result

\begin{equation*}
\frac{\mathfrak{A}}{x^m}+\frac{\mathfrak{B}}{x^{m-n}}+\frac{\mathfrak{C}}{x^{m-2n}}+\frac{\mathfrak{D}}{x^{m-3n}}+\text{etc.},
\end{equation*}
of which series so many terms are to be taken as one has positive exponents of $x$ in the denominator. And it will be

\begin{alignat*}{9}
&    &&   &&\mathfrak{A} f^n &&=1 \quad && \text{therefore}  \quad && \mathfrak{A} &&=+\frac{1}{f^{n}}\\[3mm]
&\mathfrak{A}g^n&&~+~&&\mathfrak{B} g^n &&=0 \quad &&  \quad && \mathfrak{B} &&=-\frac{g^n}{f^{2n}}\\[3mm]
&\mathfrak{B}g^n&&~+~&&\mathfrak{C} g^n &&=0 \quad &&  \quad && \mathfrak{C} &&=-\frac{g^{2n}}{f^{3n}}\\[3mm]
&\mathfrak{C}g^n&&~+~&&\mathfrak{D} g^n &&=0 \quad &&  \quad && \mathfrak{D} &&=-\frac{g^{3n}}{f^{4n}}\\
& && &&\text{etc.} && && && &&\text{etc.}
\end{alignat*}
Therefore, the propounded fraction in total will be resolved into these simple fractions

\begin{equation*}
\frac{1}{f^nx^m}-\frac{g^n}{f^{2n}x^{m-n}}+\frac{g^{2n}}{f^{3n}x^{m-2n}}-\frac{g^{3n}}{f^{4n}x^{m-3m}}+\text{etc.}
\end{equation*}
\begin{equation*}
-\frac{2fg^m \sin \frac{\pi}{n}\cdot \sin \frac{(m-1)\pi}{n}+2g^m\cos \frac{(m-1)\pi}{n}\left(gx-f \cos \frac{\pi}{n}\right)}{nf^{n+m-1}\left(ff-2fgx \cos \frac{\pi}{n}+ggxx\right)}
\end{equation*}
\begin{equation*}
-\frac{2fg^m \sin \frac{3\pi}{n}\cdot \sin \frac{3(m-1)\pi}{n}+2g^m\cos \frac{3(m-1)\pi}{n}\left(gx-f \cos \frac{3\pi}{n}\right)}{nf^{n+m-1}\left(ff-2fgx \cos \frac{3\pi}{n}+ggxx\right)}
\end{equation*}
\begin{equation*}
-\frac{2fg^m \sin \frac{5\pi}{n}\cdot \sin \frac{5(m-1)\pi}{n}+2g^m\cos \frac{5(m-1)\pi}{n}\left(gx-f \cos \frac{5\pi}{n}\right)}{nf^{n+m-1}\left(ff-2fgx \cos \frac{5\pi}{n}+ggxx\right)}
\end{equation*}
\begin{equation*}
\text{etc.}
\end{equation*}
To these formulas, if $n$ was an odd number, because of the factor $f+gx$ of the denominator, one furthermore has to add

\begin{equation*}
\frac{\pm g^m}{nf^{n+m-1}(f+gx)},
\end{equation*}
where the upper of the two signs $\pm$ holds, if $m$ was an even number, the lower on the other hand, if $m$ was odd.

\paragraph*{§418}

Now, let us also consider the formula $a+bx^n$, if $b$ was a negative number, and let this function be propounded

\begin{equation*}
f^n-g^nx^n,
\end{equation*}
 one factor of which will always be $f-gx$; and if $n$ is an even number, also $f+gx$ will be a factor. The remaining factors on the other hand are trinomial factors; if their general form is put

\begin{equation*}
ff-2fgx \cos \varphi +ggxx,
\end{equation*}
it will be

\begin{equation*}
f^n-f^n\cos n \varphi =0 \quad \text{and} \quad f^n \sin n \varphi =0
\end{equation*}
or

\begin{equation*}
\sin n \varphi =0 \quad \text{and} \quad \cos n \varphi =1.
\end{equation*}
To satisfy those equations, it is necessary that  $n \varphi = 2 k \pi$ while $k$ is any integer number and therefore it will be $\varphi= \frac{2k \pi}{n}$. Therefore, the general factor will be

\begin{equation*}
ff-2fgx \cos \frac{2k \pi}{n}+ggxx;
\end{equation*}
therefore, by taking  even numbers smaller than the exponent $n$ for $2k$ all trinomial factors will result

\begin{alignat*}{9}
&ff&&-2fg x \cos \frac{2 \pi}{n}&&+ggxx \\[3mm]
&ff&&-2fg x \cos \frac{4 \pi}{n}&&+ggxx \\[3mm]
&ff&&-2fg x \cos \frac{6 \pi}{n}&&+ggxx \\
& && \text{etc.}
\end{alignat*} 

\subsection*{Example 1}

\textit{To resolve the fraction $\frac{x^m}{f^n-g^nx^n}$ into its simple fractions.}\\[2mm]
Since one factor of the denominator is $f-gx$, hence a fraction of this kind will result $\frac{\mathfrak{A}}{f-gx}$; to find its numerator, put $P=x^m$ and $f^n-g^nx^n=Q$; it will be

\begin{equation*}
dQ=-ng^nx^{n-1}dx
\end{equation*}
and it will be

\begin{equation*}
\mathfrak{A}=\frac{-gx^m}{-ng^n x^{n-1}}=\frac{x^m}{ng^{n-1}x^{n-1}}
\end{equation*}
having put $x=\frac{f}{g}$. Therefore, it will be $\mathfrak{A}=\frac{1}{nf^{n-m-1}g^m}$ and hence the simple fraction resulting from the factor $f-gx$ will be

\begin{equation*}
\frac{1}{nf^{n-m-1}g^m(f-gx)}.
\end{equation*}
If $n$ is an even number, since then one factor of the denominator also is $f+gx$, put the simple fraction to result from this $=\frac{\mathfrak{A}}{f+gx}$; it will be

\begin{equation*}
\mathfrak{A}=\frac{-g x^m}{ng^nx^{n-1}}=\frac{-x^m}{ng^{n-1}x^{n-1}}
\end{equation*}
having put $x=\frac{-f}{g}$. Therefore, because of the odd number $n-1$, it will be $g^{n-1}x^{n-1}=-f^{n-1}$; but it will be $x^m= \frac{\pm f^m}{g^m}$, where the upper sign holds, if $m$ was an even number, the lower, if $m$ was an odd number. Therefore, because $\mathfrak{A}=\frac{\mp 1}{nf^{n-m-1}g^m}$, the simple fraction to result from the factor $f+gx$ will be this

\begin{equation*}
\frac{\mp 1}{nf^{n-m-1}g^m(f+gx)}.
\end{equation*}
Further, because the general form of trinomial factors is

\begin{equation*}
ff-2fgx \cos \frac{2k \pi}{n}+ggxx,
\end{equation*}
if we make the comparison to example 1 § 416, it will be $a=f^n$, $b=-g^n$ and $\varphi =\frac{2k \pi}{n}$; hence

\begin{equation*}
\sin n \varphi =  0 \quad \text{and} \quad \cos n \varphi =1
\end{equation*}
and
\begin{equation*}
\sin (n-m-1)\varphi =- \sin (m+1)\varphi =- \sin \frac{2k(m+1)\pi}{n}
\end{equation*}
and

\begin{equation*}
\cos (n-m-1)\varphi =\cos (m+1)\varphi = \cos \frac{2k(m+1)\pi}{n}.
\end{equation*}
From these equations the simple fraction to result from this will be

\begin{equation*}
\frac{2f \sin \frac{2k\pi}{n}\cdot \sin \frac{2k(m+1)\pi}{n}-2 \cos \frac{2k(m+1)\pi}{n}\left(gx-f \cos \frac{2k \pi}{n}\right)}{nf^{n-m-1}g^m\left(ff-2fgx \cos \frac{2k \pi}{n}+ggxx\right)}.
\end{equation*}
Therefore, the simple fractions in question will be

\begin{equation*}
\frac{1}{nf^{n-m-1}g^m(f-gx)}
\end{equation*}
\begin{equation*}
+\frac{2f \sin \frac{2\pi}{n}\cdot \sin \frac{2(m+1)\pi}{n}-2 \cos \frac{2(m+1)\pi}{n}\left(gx-f \cos \frac{2 \pi}{n}\right)}{nf^{n-m-1}g^m\left(ff-2fgx \cos \frac{2 \pi}{n}+ggxx\right)}
\end{equation*}
\begin{equation*}
+\frac{2f \sin \frac{4\pi}{n}\cdot \sin \frac{4(m+1)\pi}{n}-2 \cos \frac{4(m+1)\pi}{n}\left(gx-f \cos \frac{4 \pi}{n}\right)}{nf^{n-m-1}g^m\left(ff-2fgx \cos \frac{4 \pi}{n}+ggxx\right)}
\end{equation*}
\begin{equation*}
+\frac{2f \sin \frac{6\pi}{n}\cdot \sin \frac{6(m+1)\pi}{n}-2 \cos \frac{6(m+1)\pi}{n}\left(gx-f \cos \frac{6 \pi}{n}\right)}{nf^{n-m-1}g^m\left(ff-2fgx \cos \frac{6 \pi}{n}+ggxx\right)}
\end{equation*}
\begin{equation*}
\text{etc.},
\end{equation*}
to which, if $n$ was an even number, one furthermore has to add this fraction

\begin{equation*}
\frac{\mp 1}{nf^{n-m-1}g^m(f+gx)},
\end{equation*}
of which the upper sign $-$ is to be taken, if $m$ was an even number, the lower, if odd. Furthermore, if $m$ is a number not smaller than $n$, this integer parts are to be added

\begin{equation*}
Ax^{m-n}+Bx^{m-2n}+Cx^{m-3n}+Dx^{m-4n}+\text{etc.},
\end{equation*}
as long the exponents were not negative, and it will be

\begin{alignat*}{9}
&    &&~-~   &&Ag^n &&=1 \quad && \text{therefore}  \quad && A &&=-\frac{1}{g^{n}}\\[3mm]
&Af^n&&~-~&&Bg^n &&=0 \quad &&  \quad && B &&=-\frac{f^n}{g^{2n}}\\[3mm]
&Bf^n&&~-~&&Cg^n &&=0 \quad &&  \quad && C &&=-\frac{f^{2n}}{g^{3n}}\\[3mm]
&Cf^n&&~-~&&Dg^n &&=0 \quad &&  \quad && C &&=-\frac{f^{3n}}{g^{4n}}\\
& && &&\text{etc.} && && && &&\text{etc.}
\end{alignat*}

\subsection*{Example 2}

\textit{To resolve the fraction $\frac{1}{x^m(f^n-g^nx^n)}$ into its simple fractions.}\\[2mm]
The fractions resulting from the factor $f^n-g^nx^n$ of the denominator will be the same as before, as long as in those formulas $m$ is taken negatively. Hence one has to consider the other factor $x^m$; if we put that from this these fractions the following expression results

\begin{equation*}
\frac{\mathfrak{A}}{x^m}+\frac{\mathfrak{B}}{x^{m-n}}+\frac{\mathfrak{C}}{x^{m-2n}}+\frac{\mathfrak{D}}{x^{m-3n}}+\text{etc.},
\end{equation*}
which  series is to be continued until the exponents of $x$ become negative, it will be

\begin{alignat*}{9}
&  &&+\mathfrak{A}f^n &&=1 \quad        \text{therefore} \quad && \mathfrak{A}&&=\frac{1}{f^n} \\[3mm]
&\mathfrak{B}f^n &&- \mathfrak{A}g^n && =0  &&\mathfrak{B}&&=\frac{g^n}{f^{2n}} \\[3mm]
&\mathfrak{C}f^n &&- \mathfrak{B}g^n && =0   &&\mathfrak{C}&&=\frac{g^{2n}}{f^{3n}} \\[3mm]
&\mathfrak{D}f^n &&- \mathfrak{C}g^n && =0  &&\mathfrak{D}&&=\frac{g^{3n}}{f^{4n}} \\
& && \text{etc.} && && && \text{etc.}
\end{alignat*}
Therefore, the propounded fraction will be resolved into these simple fractions

\begin{equation*}
\frac{1}{f^n x^m}+\frac{g^n}{f^{2n}x^{m-n}}+\frac{g^{2n}}{f^{3n}x^{m-2n}}+\frac{g^{3n}}{f^{4n}x^{m-3n}}+\text{etc.}
\end{equation*}
\begin{equation*}
+\frac{g^m}{nf^{n+m-1}(f-gx)}
\end{equation*}
\begin{equation*}
-\frac{2fg^m \sin \frac{2\pi}{n}\cdot \sin \frac{2(m-1)\pi}{n}+2g^m \cos \frac{2(m-1)\pi}{n}\left(gx-f \cos \frac{2 \pi}{n}\right)}{nf^{n+m-1}\left(ff-2fgx \cos \frac{2 \pi}{n}+ggxx\right)}
\end{equation*}
\begin{equation*}
-\frac{2fg^m \sin \frac{4\pi}{n}\cdot \sin \frac{4(m-1)\pi}{n}+2g^m \cos \frac{4(m-1)\pi}{n}\left(gx-f \cos \frac{4 \pi}{n}\right)}{nf^{n+m-1}\left(ff-2fgx \cos \frac{4 \pi}{n}+ggxx\right)}
\end{equation*}
\begin{equation*}
-\frac{2fg^m \sin \frac{6\pi}{n}\cdot \sin \frac{6(m-1)\pi}{n}+2g^m \cos \frac{6(m-1)\pi}{n}\left(gx-f \cos \frac{6 \pi}{n}\right)}{nf^{n+m-1}\left(ff-2fgx \cos \frac{6 \pi}{n}+ggxx\right)}
\end{equation*}
\begin{equation*}
\text{etc.},
\end{equation*}
to which, if $n$ was an even number, one additionally has to add this fraction

\begin{equation*}
\frac{\mp g^m}{nf^{n+m-1}(f+gx)},
\end{equation*}
which is omitted, if $n$ was an odd number. The upper sign $-$ of the ambiguous signs holds, if $m$ is an even number, the lower $+$ on the other hand, if $m$ is an odd number.

\paragraph*{§419}

Therefore, this way all fractions whose denominator consists of two terms of this kind $a+bx^n$ are resolved into simple fractions. But if the denominator consists of three terms of this kind $a+bx^n+cx^{2n}$, then one has to see at first, whether it can be resolved into two real factors of this first form. For, if this is possible, the resolution into simple fractions can be done in the way explained before. For, if a fraction of this kind is propounded

\begin{equation*}
\frac{x^m}{(f^n+g^nx^n)(f^n+k^nx^n)},
\end{equation*}
it will at first be transformed into two of this kind

\begin{equation*}
\frac{\alpha x^m}{f^n+g^nx^n}+\frac{\beta x^m}{f^n +h^n x^n}
\end{equation*}
and it will be

\begin{equation*}
\alpha f^n+ \beta f^n =1 \quad \text{and} \quad \alpha h^n+\beta g^n=0,
\end{equation*}
whence

\begin{equation*}
\alpha = \frac{1}{f^n}-\beta = - \frac{\beta g^n}{h^n}.
\end{equation*}
If the exponent $m$ was greater than $n$, the transformation into the following fractions will be more convenient

\begin{equation*}
\frac{\alpha x^{m-n}}{f^n+g^n x^n}+\frac{\beta x^{m-n}}{f^n+h^n x^n},
\end{equation*}
by means of which

\begin{equation*}
\alpha +\beta =0 \quad \text{and} \quad \alpha h^n+\beta g^n =1
\end{equation*}
and hence

\begin{equation*}
\alpha =\frac{1}{h^n-g^n} \quad \text{and} \quad \beta =\frac{1}{g^n-h^n}.
\end{equation*}
But no matter which of both transformations is used, both fractions to result this way will be resolved  into its simple fractions, which taken together will be equal to the propounded fraction using the method explained before.

\paragraph*{§420}

In like manner the method treated up to this point will suffice, if the denominator consists of several terms of this kind

\begin{equation*}
a+bx^n+cx^{2n}+dx^{3n}+ex^{4n}+\text{etc.},
\end{equation*}
if it only can be resolved into factors of this form $f^n\pm g^nx^n$. For, let us assume that this fraction  is to be resolved into its simple factors 

\begin{equation*}
\frac{x^m}{(a-x^n)(b-x^n)(c-x^n)(d-x^n)\text{etc.}}.
\end{equation*}
At first, resolve it into these

\begin{equation*}
\frac{Ax^m}{a-x^n}+\frac{Bx^m}{b-x^n}+\frac{Cx^m}{c-x^n}+\frac{dx^m}{d-x^n}+\text{etc.},
\end{equation*}
the numerators of which will be determined the following way that 

\begin{alignat*}{9}
&A &&= \frac{1}{(b-a)(c-a)(d-a)\text{etc.}}\\[3mm]
&B &&=\frac{1}{(a-b)(c-b)(d-b)\text{etc.}}\\[3mm]
&C &&=\frac{1}{(a-c)(b-c)(d-c)\text{etc.}}\\
&  &&\text{etc.}
\end{alignat*}
Therefore, after this preparation each fraction will be resolved into their simple fractions applying the method explained before.

\paragraph*{§421}

If a denominator of this kind

\begin{equation*}
a+bx^n+cx^{2n}+dx^{3n}+\text{etc.}
\end{equation*}
has not only real factors of the form $f^n\pm g^n x^n$, two imaginary ones are to be combined into a single real one. Therefore, let us put that the product of two factors of this kind is

\begin{equation*}
f^{2n}-2f^ng^n \cos \omega +g^{2n}x^{2n};
\end{equation*}
and because this expression has no simple real factors, let us put that the trinomial ones are contained in this general form

\begin{equation*}
ff-2fgx \cos \varphi +ggxx,
\end{equation*}
whose number will be $=n$. Therefore, having put $x^n=\frac{f^n}{g^n}\cos n \varphi$, this equation will result

\begin{equation*}
1-2 \cos \omega \cdot \cos n \varphi +\cos 2n \varphi =0.
\end{equation*}
Furthermore, having put $x^n=\frac{f^n}{g^n}\sin n \varphi$, it will also be

\begin{equation*}
-2 \cos \omega \cdot \sin n \varphi + \sin 2 n \varphi =0,
\end{equation*}
which divided by $\sin n \varphi$ gives $\cos n \varphi = \cos \omega$ and so at the same time the first equation will be satisfied. Therefore, it will be $n \varphi = 2 k \pi \pm \omega $ while $k$ denotes any positive integer and hence it will be $\varphi = \frac{2k \pi \pm \omega}{n}$ and all factors will be contained in this form

\begin{equation*}
ff-2fg x \cos \frac{2k \pi \pm \omega}{n}+ggxx,
\end{equation*}
whence one will have the following factors

\begin{alignat*}{9}
&ff &&- 2fg x \cos \frac{\quad ~ \omega \quad ~}{n}&&+ggxx \\[3mm]
&ff &&- 2fg x \cos \frac{2 \pi - \omega}{n}&&+ggxx \\[3mm]
&ff &&- 2fg x \cos \frac{2 \pi + \omega}{n}&&+ggxx \\[3mm]
&ff &&- 2fg x \cos \frac{4 \pi - \omega}{n}&&+ggxx \\[3mm]
&ff &&- 2fg x \cos \frac{4 \pi + \omega}{n}&&+ggxx \\[3mm]
&   && \text{etc.},
\end{alignat*}
of which so many are to be taken until their number becomes $=n$.

\paragraph*{§422}

Therefore, if this fraction is propounded to be resolved into its simple fractions

\begin{equation*}
\frac{x^{m-1}}{f^{2n}-2f^ng^n \cos \omega +g^{2n}x^{2n}},
\end{equation*}
since any trinomial factor of the denominator is contained in this form

\begin{equation*}
ff-2fgx \cos \varphi +ggxx
\end{equation*}
while $\varphi =\frac{2k \pi \pm \omega}{n}$, consider this fraction

\begin{equation*}
\frac{x^m}{f^{2n}x-2 f^n g^n x^{n+1} \cos \omega +g^{2n}x^{2n+1}}
\end{equation*}
equal to it and put $x^m=P$ and the denominator

\begin{equation*}
f^{2n}x-2f^ng^nx^{n+1} \cos \omega +g^{2n}x^{2n+1}=Q;
\end{equation*}
it will be

\begin{equation*}
\frac{dQ}{dx}=f^{2n}-2(n+1)f^n g^n x^n \cos \omega +(2n+1) g^{2n}x^{2n}.
\end{equation*}
Therefore, by putting

\begin{equation*}
x^n=\frac{f^n}{g^n}\cos n \varphi,
\end{equation*}
it will be

\begin{equation*}
\mathfrak{P}=\frac{f^m}{g^m}\cos m \varphi \quad \text{and} \quad \mathfrak{P}=\frac{f^m}{g^m}\cos \frac{m(2k \pi \pm \omega)}{n}
\end{equation*}
and

\begin{equation*}
\mathfrak{Q}=f^{2n}(1-2(n+1)\cos \omega \cdot \cos n \varphi +(2n+1) \cos 2n \varphi).
\end{equation*}
But because  $\cos n \varphi = \cos \omega$, it will be

\begin{equation*}
\cos 2n \varphi = 2 \cos^2 \omega -1
\end{equation*}
and hence

\begin{equation*}
\mathfrak{Q}= f^{2n}(-2n+2n \cos^2 \omega)=-2n f^{2n} \sin^2 \omega.
\end{equation*}
Further, having put

\begin{equation*}
x^n= \frac{f^n}{g^n}\sin n \varphi,
\end{equation*}
it will be

\begin{equation*}
\mathfrak{p}= \frac{f^m}{g^m}\sin m \varphi = \frac{f^m}{g^m}\sin \frac{m(2 k \pi \pm \omega)}{n}
\end{equation*}
and

\begin{equation*}
\mathfrak{q}=-f^{2n}(2(n+1)\cos \omega \cdot \sin n \varphi -(2n+1)\sin 2n \varphi);
\end{equation*}
because of

\begin{equation*}
\sin 2n \varphi = 2 \sin n \varphi \cdot \cos n \varphi = 2 \cos \omega \cdot \sin n \varphi,
\end{equation*}
it will be

\begin{equation*}
\mathfrak{q}= 2n f^{2n} \cos \omega \cdot \sin n \varphi.
\end{equation*}
But because  $n \varphi = 2k \pi \pm \omega$, it will be $\sin n \varphi = \pm \sin \omega$ and

\begin{equation*}
\mathfrak{q}=\pm 2n f^{2n} \sin \omega \cdot \cos \omega.
\end{equation*}
Having found these, it will be

\begin{equation*}
\mathfrak{Q}^2+\mathfrak{q}^2= 4n^2 f^{4n} \sin^2 \omega,
\end{equation*}
\begin{equation*}
\mathfrak{Pq}-\mathfrak{pQ}=\frac{2n f^{m+2n}}{g^m}(\pm \cos m \varphi \cdot \sin \omega \cdot \cos \omega + \sin m \varphi \cdot \sin^2 \omega)
\end{equation*}
or

\begin{equation*}
\mathfrak{Pq}-\mathfrak{pQ}= \pm \frac{2nf^{m+2n}}{g^m}\sin \omega \cdot \cos (m \varphi \mp \omega)
\end{equation*}
or 

\begin{alignat*}{9}
&\mathfrak{Pq}&&- \mathfrak{pQ}&&~=~&& \pm ~&&\frac{2n f^{m+2n}}{g^m}\sin \omega \cdot \cos \frac{2k m \pi \pm (m-n)\omega}{n} \\[3mm]
&\mathfrak{PQ}&&+ \mathfrak{pq}&&~ =~ && &&\frac{2n f^{m+2n}}{g^m}(- \cos m \varphi \cdot \sin^2 \omega \pm \sin m \varphi \cdot \sin \omega \cos \omega), \\[3mm]
&\mathfrak{PQ}&&+ \mathfrak{pq}&&~ =~ &&~\pm~ &&\frac{2n f^{m+2n}}{g^m}\sin \omega \cdot \sin (m \varphi \mp \omega), 
\intertext{or}
&\mathfrak{PQ}&&+ \mathfrak{pQ}&&~=~&& \pm ~&&\frac{2n f^{m+2n}}{g^m}\sin \omega \cdot \sin \frac{2k m \pi \pm (m-n)\omega}{n}.
\end{alignat*}
Therefore, from the factor of the denominator

\begin{equation*}
ff- 2fg x \cos \frac{2k \pi \pm \omega}{n}+ggxx
\end{equation*}
this simple fraction results

\begin{equation*}
\frac{\pm f \sin \frac{2k \pi \pm \omega}{n}\cdot \cos \frac{2km \pi \pm (m-n)\pi}{n}\pm \sin \frac{2k m \pi \pm (m-n)\omega}{n}\left(gx- f \cos \frac{2k \pi \pm \omega}{n}\right)}{nf^{2n-m}g^{m-1} \sin \omega \left(ff- 2fgx \cos \frac{2k \pi \pm \omega}{n}+ggxx\right)}
\end{equation*}
or

\begin{equation*}
\frac{gx \sin \frac{2k m \pi \pm(m-n)\omega}{n}\pm f \sin \frac{2k(m-1)\pi \pm (m-n-1)\omega}{n}}{nf^{2n-m}g^{m-1} \sin \omega \left(ff- 2fgx \cos \frac{2k \pi \pm \omega}{n}+ggxx\right)}.
\end{equation*}

\subsection*{Example}

\textit{To resolve the fraction $\frac{x^{m-1}}{f^{2n}-2f^ng^n x^n \cos \omega +g^{2n}x^{2n}}$ into its simple fractions.}\\[2mm]
These simple fractions in question will therefore be

\begin{equation*}
\frac{f \sin \frac{\omega}{n}\cdot \cos \frac{(m-n)\omega}{n}+\sin \frac{(m-n)\omega}{n}\left(gx- f \cos \frac{\omega}{n}\right)}{nf^{2n-m}g^{m-1} \sin \omega \left(ff- 2fgx \cos \frac{ \omega}{n}+ggxx\right)}
\end{equation*}
\begin{equation*}
- \frac{f \sin \frac{2 \pi - \omega}{n} \cos \frac{2m \pi -(m-n)\omega}{n}+\sin \frac{2m \pi -(m-n)\omega}{n}\left(gx-f \cos \frac{2 \pi -\omega}{n}\right)}{nf^{2n-m}g^{m-1} \sin \omega \left(ff- 2fgx \cos \frac{2 \pi - \omega}{n}+ggxx\right)}
\end{equation*}
\begin{equation*}
+ \frac{f \sin \frac{2 \pi + \omega}{n} \cos \frac{2m \pi +(m-n)\omega}{n}+\sin \frac{2m \pi +(m-n)\omega}{n}\left(gx-f \cos \frac{2 \pi +\omega}{n}\right)}{nf^{2n-m}g^{m-1} \sin \omega \left(ff- 2fgx \cos \frac{2 \pi + \omega}{n}+ggxx\right)}
\end{equation*}
\begin{equation*}
- \frac{f \sin \frac{4 \pi - \omega}{n} \cos \frac{4m \pi -(m-n)\omega}{n}+\sin \frac{4m \pi -(m-n)\omega}{n}\left(gx-f \cos \frac{4 \pi -\omega}{n}\right)}{nf^{2n-m}g^{m-1} \sin \omega \left(ff- 2fgx \cos \frac{4 \pi - \omega}{n}+ggxx\right)}
\end{equation*}
\begin{equation*}
+ \frac{f \sin \frac{4 \pi + \omega}{n} \cos \frac{4m \pi +(m-n)\omega}{n}+\sin \frac{4m \pi +(m-n)\omega}{n}\left(gx-f \cos \frac{4 \pi +\omega}{n}\right)}{nf^{2n-m}g^{m-1} \sin \omega \left(ff- 2fgx \cos \frac{4 \pi + \omega}{n}+ggxx\right)}
\end{equation*}
\begin{equation*}
\text{etc.}
\end{equation*}
and one has to continue this way until the number of these fractions was $n$. If $m$ was a number either greater than $2n-1$ or negative, in the first case integral parts, in the second fraction are to be added, which are easily found applying the method explained before.